\pdfoutput=1
\documentclass[11pt]{article}
\usepackage[T1]{fontenc}
\usepackage{lmodern}
\usepackage[margin=1in]{geometry}
\usepackage{amsmath,amssymb,amsthm,mathtools,enumitem}
\usepackage{etoolbox}
\usepackage[colorlinks=true,linkcolor=blue,citecolor=blue,urlcolor=blue]{hyperref}
\newtheorem{theorem}{Theorem}[section]
\newtheorem{proposition}[theorem]{Proposition}
\newtheorem{lemma}[theorem]{Lemma}
\newtheorem{corollary}[theorem]{Corollary}
\theoremstyle{definition}
\newtheorem{definition}[theorem]{Definition}
\newtheorem*{definition*}{Definition}
\theoremstyle{remark}\newtheorem{remark}[theorem]{Remark}
\newcommand{\R}{\mathbb R}\newcommand{\C}{\mathbb C}\newcommand{\Z}{\mathbb Z}
\newcommand{\Q}{\mathbb Q}
\newcommand{\GL}{\operatorname{GL}}\newcommand{\SL}{\operatorname{SL}}
\newcommand{\Diff}{\operatorname{Diff}}
\newcommand{\Sym}{\operatorname{Sym}}
\newcommand{\Hom}{\operatorname{Hom}}\newcommand{\End}{\operatorname{End}}
\newcommand{\id}{\operatorname{id}}\newcommand{\HS}{\mathrm{HS}}

\newcommand{\D}{\mathcal D_+^\omega}\newcommand{\RP}{\mathbb{RP}}
\newcommand{\Pj}{\mathbb P}\newcommand{\Leb}{\operatorname{Leb}}
\newcommand{\nstd}{\nabla^{\mathrm{std}}}
\setlist[enumerate]{itemsep=2pt,topsep=5pt}
\makeatletter
\patchcmd{\l@section}{1.5em}{2.4em}{}{%
  \PackageError{ClassificationSeries}{Cannot widen contents section numbers}{}}
\makeatother
\newcommand{\Fl}{\operatorname{Fl}}

\theoremstyle{plain}

\numberwithin{equation}{section}

\hypersetup{
  pdftitle={Real-analytic actions of lattices in SL(n,R) on compact n-manifolds},
  pdfauthor={Federico Rodriguez Hertz and Zhiren Wang},
  bookmarksdepth=section
}
\title{Real-analytic actions of lattices in $\SL(n,\R)$\\
       on compact $n$-manifolds}
\author{Federico Rodriguez Hertz \and Zhiren Wang}
\date{}

\begin{document}
\maketitle
\begin{abstract}
We classify real-analytic actions of arbitrary lattices in
$\SL(n,\R)$ on compact connected $n$-manifolds without boundary,
$n\ge3$, with no invariant-volume assumption. After finite-index
restriction, the action has finite image, extends to a connected
Lie group locally isomorphic to $\SL(n,\R)$ with an equivariant
finite cover when needed, or is obtained by analytic surgery
with at least one torus piece. In the last case the lattice is
commensurable with a conjugate of $\SL(n,\Z)$. As a main input, 
we prove a general cocycle superrigidity theorem for real 
analytic germs.
\end{abstract}

\tableofcontents

\clearpage
\section*{Introduction}
\addcontentsline{toc}{part}{Introduction}
Let $n\ge 3$, let $\Gamma<\SL(n,\R)$ be a lattice, and let $M$ be a
compact connected real-analytic $n$-manifold without boundary. We
classify real-analytic actions of $\Gamma$ on $M$. The lattice may
be uniform or nonuniform, and the action need not preserve volume.

After restriction to a finite-index subgroup, the action has finite
image, extends to an analytic action of a connected Lie group
locally isomorphic to $\SL(n,\R)$ (with an equivariant finite cover
when needed), or is assembled by analytic surgery from pieces that
include a punctured torus or a specified finite quotient of one.
In the last case, the lattice is commensurable with a conjugate of
$\SL(n,\Z)$. The exact alternatives and the covering and gluing data
are stated in the
\hyperref[IV-main:thm:classification]{classification theorem in Part IV}.
When the action extends to $\SL(n,\R)$ itself, the ambient action
belongs to Fisher and Melnick's complete classification of
real-analytic $\SL(n,\R)$-actions on closed $n$-manifolds
\cite[Corollary~5.2]{FisherMelnick}.

\subsection*{Historical background}
\paragraph{From lattice superrigidity to transformation groups.}
Mostow's strong rigidity theorem and Margulis's superrigidity
theorem established that the ambient Lie group controls much of
the geometry and representation theory of its lattices
\cite{Mostow,Margulis}. Zimmer extended the rigidity viewpoint to
measurable cocycles over ergodic actions \cite{Zimmer1980,Zimmer}.
For an action on a manifold with an invariant probability, the
derivative becomes such a cocycle after a measurable choice of
frames. Zimmer used this observation to relate low-dimensional
volume-preserving actions to representations of the ambient group
\cite{Zimmer1983,Zimmer1984}. Gromov developed a complementary
approach using rigid geometric structures and their local
automorphisms \cite{Gromov}. Fisher's survey \cite{FisherSurvey} traces these
origins and the development of geometric and dynamical rigidity.
Brown's survey \cite{BrownICM2022CompactTargets} describes the
later use of suspension measures in both the dimension bounds
and the classification at the first nontrivial dimension.

\paragraph{Early restrictions on low-dimensional actions.}
Several approaches gave finiteness results before the general
dimension bounds were known. Witte proved that arithmetic groups
of higher rational rank have finite image in actions on the circle
and trivial orientation-preserving actions on the line
\cite{Witte1994}. Ghys used the bounded Euler class in his rigidity
theorems for circle actions of lattices. Burger and Monod developed
continuous bounded cohomology and obtained circle-action finiteness
results through the comparison with ordinary cohomology
\cite{GhysCircle,BurgerMonod}. Deroin and Hurtado extended circle
rigidity to actions by homeomorphisms. For lattices in connected
simple groups of higher real rank and finite center this gives
finite image \cite[Theorem~1.5]{DeroinHurtado}. For surfaces,
Polterovich used the growth of iterates and distortion in groups of symplectic
diffeomorphisms. Franks and Handel used distortion together with
surface dynamics to prove further finiteness results for
area-preserving actions \cite{Polterovich,FranksHandel}.
In the real-analytic category, Farb and Shalen used fixed-point
arguments and low-dimensional linear representations to prove
finiteness for volume-preserving actions on closed four-manifolds of nonzero
Euler characteristic for lattices of rational rank at least seven
in simple linear Lie groups \cite{FarbShalen}.

\paragraph{Lattice actions with hyperbolicity.}
Hurder studied Anosov lattice actions, and Katok, Lewis, and
Zimmer obtained global rigidity for lattice actions on tori under
hyperbolicity and measure hypotheses \cite{Hurder,KatokLewisZimmer}.
Katok and Spatzier used higher-rank dynamics to prove local rigidity
for several algebraic lattice actions, including boundary actions.
Goetze and Spatzier gave a smooth classification of volume-preserving
Cartan actions of irreducible lattices and their ambient groups.
Fisher and Margulis proved local rigidity for affine actions of
higher-rank groups and their lattices
\cite{KatokSpatzier,GoetzeSpatzier,FisherMargulisAffine}.
Margulis and Qian obtained semiconjugacies for weakly hyperbolic
lattice actions under hypotheses on their algebraic models
\cite{MargulisQian}. Global rigidity theorems for higher-rank
abelian Anosov and totally Cartan actions supplied tools for
this direction \cite{FKS,RHWAnosov,SpatzierVinhageCartan}.

Brown, Rodriguez Hertz, and Wang obtained semiconjugacies for
lattice actions on compact nilmanifolds with hyperbolic linear
data under a lifting hypothesis. For a smooth action with an Anosov
element, the semiconjugacy is a smooth conjugacy
\cite{BRHWGlobal}. Lee showed that a smooth volume-preserving
action of a lattice in $\SL(n,\R)$ on a closed $n$-manifold is
smoothly conjugate to an affine torus action if one element has a
dominated splitting \cite{LeeDominated}. His dynamical cocycle
superrigidity theorem treats weakly irreducible lattices
in higher-rank semisimple groups without property $(T)$ under
integrability and irreducibility hypotheses, with an application
to Anosov actions on nilmanifolds \cite{LeeDynamical}.
Lee and Sandfeldt also obtain affine rigidity on tori and on
certain two-step nilmanifolds from a partially hyperbolic element
with one-dimensional center, under hypotheses on its lifted
foliations and, for tori, on the action on the fundamental group
\cite[Theorems~1.2 and~1.4]{LeeSandfeldt}.

For semisimple groups whose simple factors all have real rank at
least two, Damjanovi\'c, Spatzier, Vinhage, and Xu classify smooth,
totally partially hyperbolic accessible actions with an invariant
volume for the split Cartan subgroup.
They obtain bi-homogeneous models after a finite cover and discuss
the constraints on applying their theorem to lattice suspensions
\cite[Theorem~2.5 and \S2.5.7]{DamjanovicSpatzierVinhageXu}.

\paragraph{Suspension measures and projective factors.}
Brown, Rodriguez Hertz, and Wang used the split Cartan action on
the suspension to study lattice actions. A principle developed in
their work on hyperbolic actions and measurable projective factors
is that a root which is not resonant with a fiberwise Lyapunov
exponent can give additional invariance of a suspension measure.
With a dimension bound, this can give an invariant probability or
a measurable factor over a projective action \cite{BRHWFactors}.
For actions without an invariant probability, the earlier work
of Nevo and Zimmer on stationary actions had already shown the
role of parabolic factors and relatively measure-preserving
extensions \cite{NevoZimmer}. Oseledets's multiplicative
ergodic theorem supplies the Lyapunov data. Pesin's theory and
Ledrappier and Young's work relating entropy to exponents and
conditional measures supply tools for turning that data into
invariance \cite{Oseledets,Pesin,LedrappierYoung}.

\paragraph{The dimension threshold and projective actions.}
Brown, Fisher, and Hurtado proved the finiteness conclusion for
actions of lattices in $\SL(n,\R)$ on compact manifolds of dimension
less than $n-1$, and for volume-preserving actions in dimension less
than $n$, in the regularity of their theorems. They treated cocompact
lattices, $\SL(n,\Z)$, and then general nonuniform lattices
\cite{BFHCocompact,BFHSLZ,BFHGeneral}. Their proofs combine
subexponential growth of derivatives with strong property $(T)$,
introduced by Lafforgue and developed for higher-rank real groups
and their lattices by de Laat and de la Salle
\cite{Lafforgue,deLaatDeLaSalle,dLS}. The nonuniform case required
estimates on growth of cocycles and escape of mass. The
$\SL(n,\Z)$ proof uses the work of Lubotzky, Mozes, and Raghunathan
on the word and Riemannian metrics \cite{LMR2000}. Brown,
Damjanovi\'c, and Zhang proved the sharp $C^1$ finiteness bound
below dimension $n-1$ for cocompact lattices in $\SL(n,\R)$
\cite{BrownDamjanovicZhang}. For non-$\R$-split groups, Zhang obtained
low-dimensional finiteness for cocompact lattices in $\SL(n,\mathbb C)$
\cite{ZhangSLnC}. An, Brown, and Zhang proved further cases of
Zimmer's conjecture for lattices in non-$\R$-split semisimple groups
\cite{AnBrownZhang}.

Dimension $n-1$ admits the actions on projective space and on its
spherical double cover. Brown, Rodriguez Hertz, and Wang classified
infinite-image actions in that dimension in class $C^r$, $r>1$,
with the finite-index qualification for the spherical model.
They also showed that, in dimension $n$, a minimal-parabolic
invariant suspension probability which is not invariant under the
whole group yields an embedded projective hypersurface under the
hypotheses of their theorem \cite[Theorems~1.1 and 1.5]{BRHW}.
This supplies the hypersurface input for the local models developed
here.

Gorodnik and Spatzier classified smooth factors of standard flag
actions of higher-rank lattices under a differentiable sink
hypothesis \cite{GorodnikSpatzier}. Connell, Islam, Nguyen, and
Spatzier showed that sufficiently small $C^0$ perturbations of
parabolic boundary actions of cocompact lattices admit a
continuous semiconjugacy onto the standard action
\cite{ConnellIslamNguyenSpatzier}. For conformal actions of
cocompact lattices on compact pseudo-Riemannian manifolds,
Pecastaing bounded the real rank in terms of the metric signature
and proved conformal flatness at maximal rank
\cite{PecastaingConformal}. For infinite actions of cocompact lattices preserving a
projective or pseudo-Riemannian conformal structure,
Pecastaing proved global identification with the model space or
its double cover in the maximal-rank cases
\cite{PecastaingGeometry}.

Brown and Lee studied actions with a positive-entropy element
without a prescribed invariant volume. For a smooth action of a
lattice in $\SL(n,\R)$ on a closed $n$-manifold, positive
topological entropy of one element forces the lattice to be
commensurable with a conjugate of $\SL(n,\Z)$
\cite[Corollary~1.1]{BrownLeeEntropy}. With an ergodic invariant
probability giving positive metric entropy to one element, they
obtain a measurable affine torus or infratorus model on an
invariant subset of positive measure after finite-index
restriction \cite[Theorem~1.4]{BrownLeeEntropy}.

\paragraph{Ambient actions and surgery.}
An action of $\SL(n,\R)$ can be restricted to any lattice.
Uchida classified the real-analytic actions of $\SL(n,\R)$ on the
$n$-sphere. Cairns and Ghys established local linearization results
for smooth $\SL(n,\R)$ actions. Fisher and Melnick classified smooth
ambient-group actions on closed $n$-manifolds and described their
analytic data, including the transverse flows defining disks and
tubes \cite{Uchida,CG,FisherMelnick}. Son extended the ambient
real-analytic classification to closed $m$-manifolds in the range
$3\le n\le m\le 2n-3$ \cite{Son}. Arosemena Serrato classified
smooth, locally free codimension-one actions of higher-rank simple
split groups admitting a mixing minimal-parabolic invariant
probability \cite{ArosemenaSerrato}.

The torus gives a second source of actions. Katok and Lewis
constructed exotic volume-preserving examples by blowing up and
gluing along finite orbits of algebraic actions
\cite{KatokLewis}. Benveniste and Fisher studied the rigid and
almost rigid geometric structures associated with such
constructions \cite{BenvenisteFisher}. Stuck's induction from
parabolic actions illustrates how arbitrary transverse flows enter
actions without an invariant volume \cite{Stuck}.

Fisher and Melnick gave the models and formulated the classification
problem that guides this series. In their description of the
ambient actions, a parabolic character and a one-dimensional
flow encode the action transverse to a projective orbit. Their
$G$-disks and $G$-tubes turn these ambient actions into pieces
that can be joined to tori. They construct actions by inserting
disks at points of a finite orbit, or by joining blown-up points
with tubes, and set out how blow-ups and two-sided blow-ups can
be used together with those operations
\cite[\S3.2 and Theorems~4.2--4.3]{FisherMelnick}.
Their Conjecture~3.6 asks whether these constructions, finite-image
actions, and restrictions of ambient actions exhaust lattice
actions in the defining dimension \cite{FisherMelnick}. We prove
the real-analytic classification in this framework by identifying
the local models from a given lattice action and proving the
global decomposition, with explicit finite-index, covering, and
quotient data.

\paragraph{Analytic coordinates.}
Guillemin and Sternberg and, independently, Kushnirenko established
analytic linearization near fixed points of semisimple actions
\cite{GuilleminSternberg,Kushnirenko}. For cocycles of analytic
germs over hyperbolic systems, Navas and Ponce proved a Liv\v{s}ic
theorem under a periodic orbit obstruction \cite{NavasPonce}.
Guysinsky and Katok developed nonstationary normal forms for
contracting extensions. Melnick and Kalinin gave further normal
forms and centralizer results under their respective contraction
and spectral hypotheses \cite{GuysinskyKatok,Melnick,Kalinin}.
The convergence problem in Part I arises over a measurable
higher-rank action. Its formal input comes from
finite-dimensional cocycle superrigidity in the formulation of
Fisher and Margulis, and the analytic argument must also treat
invariant coefficient directions \cite{FisherMargulis}.

\subsection*{Organization of the proof}
\hyperlink{part.i}{Part I} proves analytic germ
cocycle superrigidity. It first constructs compatible formal
coordinates for the finite jets. Uniform estimates from spherical
averaging allow the formal reduction to be realized by measurable
families of analytic germs, with positive radii for the coordinate
changes and their inverses.

\hyperlink{part.ii}{Part II} and
\hyperlink{part.iii}{Part III} turn the invariant
hypersurfaces obtained from suspension measures into local analytic
models. The action along a projective hypersurface determines the
angular part. A transverse vector field determines the remaining
coordinate and the gluing data. Part II treats $n\ge4$. In dimension
three the Levi factor of a projective stabilizer has rank one.
Part III uses relative property $(T)$ and coordinate foliations.
In the neutral case it constructs transverse jets and proves
joint analyticity of the resulting charts.

\hyperlink{part.iv}{Part IV} continues these
local models to the regions incident to the hypersurfaces. It
obtains affine descent on the torus regions, determines the end
fillings and the ambient pieces, and proves that finitely many
pieces assemble to cover~$M$.

\hypertarget{acknowledgments}{}
\subsection*{Acknowledgments}
We thank Aaron Brown, David Fisher, Homin Lee and Ralf Spatzier for many helpful discussions.

FRH acknowledges support from NSF grant DMS-2453688 and from the Anatole Katok Chair
in Mathematics. ZW acknowledges support from NSF grant DMS-2453689.

\hypertarget{use.of.ai}{}
\subsection*{Statement on use of AI}

The development of our proof was AI-assisted. The authors started exploring
 the strategy back in 2014--2015, namely to establish cocycle
superrigidity of analytic jets, both in the original sense of neighborhoods
of lattice-action invariant probability measures, and in the more general
sense of neighborhoods of parabolic group invariant probability measures,
and to glue the resulting charts.

Several developments since then supplied tools and models for the
classification as an application of this strategy. These include Lafforgue's 
strong property~$(T)$ and the work of de Laat and de la Salle on strong property~$(T)$ for
higher-rank groups and lattices \cite{Lafforgue,deLaatDeLaSalle,dLS},
Brown--Fisher--Hurtado's proof of Zimmer's conjecture for actions
up to dimension $(n-2)$ \cite{BFHCocompact,BFHSLZ,BFHGeneral}, and
Brown--Rodriguez Hertz--Wang's classification of infinite-image
actions in dimension $(n-1)$ \cite[Theorem~1.1]{BRHW}.
Fisher and Melnick gave the disk, tube, and blow-up constructions
that model the actions arising by surgery and formulated the
classification conjecture in dimension $n$
\cite[\S3.2 and Conjecture~3.6]{FisherMelnick}. At the same time, recent AI tools
have provided the technical force that allowed to check all
 the branching cases in the roadmap we devised.

GPT 6-Astra was used for proof development and writing of this manuscipt. 
The authors have checked all arguments through a combination of human proofreading and Lean
verification of different portions, and take full responsibility for this manuscript. 
GPT 6.1-Sol was used for Lean formalization. Full verification in Lean is in progress and
 will be uploaded once complete.
\clearpage
\section*{Part I. Analytic germ cocycle superrigidity}
\hypertarget{part.i}{}
\addcontentsline{toc}{part}{Part I. Analytic germ cocycle superrigidity}

\section{Introduction}
Let an ergodic probability-preserving action of a higher-rank
lattice carry a cocycle of invertible analytic germs. We seek
measurable analytic coordinates in which the cocycle agrees with
its linear normal form, an algebraic representation times a compact
centralizing cocycle. A compatible family of finite-jet reductions
gives a formal change of coordinates. The problem is to realize it
by analytic germs with positive measurable radii.

The fiber dimension is arbitrary. The base probability may be
singular, the coordinates are measurable in the base and analytic
in the fiber, and their domains may depend on the base point.

An analytic \emph{germ} records a map near its base point:
two representatives define the same germ if they agree on a
smaller neighborhood. Its $N$-jet records its Taylor polynomial
through degree $N$. A formal map is a power series without an
assumed positive radius of convergence. Thus a compatible family
of finite jets determines a formal map, but need not determine
an analytic germ. A family of germs is \emph{coefficient-measurable}
when each Taylor coefficient is measurable in the base variable.
A \emph{transfer} is a family of coordinate changes $F_x$;
it replaces a cocycle $\beta(\gamma,x)$ by
$F_{\gamma x}^{-1}\circ\beta(\gamma,x)\circ F_x$.
An invariant \emph{conull} set is an invariant set of full measure.

Throughout the statements of the main results, let
$G=\mathbf G(\R)$ be connected, where $\mathbf G$ is a connected
simply connected semisimple real algebraic group. Assume that $G$
has no compact factors and that every almost simple factor has
real rank $\ge2$. Let $\Gamma<G$ be a lattice: $\Gamma$ is discrete and
$G/\Gamma$ carries a finite nonzero $G$-invariant regular Borel measure
\cite[Definition~B.2.1, p.~352]{BHV}. Both uniform and nonuniform lattices
are allowed. We use the following explicit
formulation of the factorwise rank assumption. An almost simple factor is
an actual connected normal almost simple real algebraic subgroup
$H\subset\mathbf G$; its real rank is the largest integer $r$ for which
there is a closed real algebraic group immersion $\mathbb G_m^r\hookrightarrow H$.
The rank bound is imposed on every such subgroup, before a product
decomposition is chosen. This is the split-torus formulation of real rank
\cite[Definition~12.14; Definition~19.20; Theorem~21.51]{Milne}.
The algebraic
simple-connectedness assumption lets us state the conclusions on
$G$ itself, without adding a central covering to the notation.
We use the geometric algebraic meaning of simple connectivity: after
extension to $\mathbb C$, every surjective homomorphism onto the group
from a connected linear algebraic group, with finite kernel, is an
isomorphism. For complex semisimple groups this agrees with topological
simple connectivity \cite[\S4, footnote~1, p.~149]{RoquesSinger2022}.

\subsection{Analytic coordinates for compact actions}
An analytic coordinate germ at $x$ in a $n$-manifold $M$ is a germ
of an invertible analytic map $(\R^n,0)\to(M,x)$. A family is
measurable if its Taylor coefficients in a finite analytic atlas
are measurable.

\begin{theorem}[Compact analytic actions]\label{I-thm:compact-actions}
Let $\alpha:\Gamma\to\Diff^\omega(M)$ be an action on a compact
real-analytic manifold without boundary, and let $\mu$ be an
ergodic $\alpha$-invariant Borel probability. Put $n=\dim M$.
There are a real algebraic representation $A:G\to\GL(n,\R)$,
a compact group $C\le Z_{\GL(n,\R)}(A(G))$, a measurable cocycle
$\kappa:\Gamma\times M\to C$, and measurable analytic coordinate
germs $\Phi_x:(\R^n,0)\to(M,x)$ such that
\begin{equation}\label{I-eq:compact-normal-form}
 \alpha(\gamma)\circ\Phi_x
       =\Phi_{\alpha(\gamma)x}\circ A(\gamma)\kappa(\gamma,x)
                       \qquad(\gamma\in\Gamma)
\end{equation}
on one invariant set of full $\mu$-measure. The germs $\Phi_x$
and their inverses have representatives on positive measurable
radii, measured in the chosen local coordinates.
\end{theorem}

No relation between $n$ and the rank or defining dimension of $G$
is assumed. In particular, the representation $A$ is not restricted
to a standard module or its dual. The compact cocycle is retained
in the conclusion. It is not in general a constant cocycle and is
not removed merely by changing the first frame. The theorem is a
statement about germs, not a conjugacy on a common neighborhood
of the support of $\mu$.

\subsection{Centered fiber cocycles}
The same statement applies to fibers that are not tangent spaces
of the base. Let $X$ be a compact metrizable space with a continuous
$\Gamma$-action and let $E\to X$ be a continuous real vector bundle
of finite rank $n$. A centered analytic fiber cocycle is a family
of invertible germs
\[
 f_{\gamma,x}:(E_x,0)\longrightarrow(E_{\gamma x},0),\qquad
 f_{\gamma\delta,x}=f_{\gamma,\delta x}\circ f_{\delta,x}.
\]
We say that it has continuous finite jets if, in local vector
bundle trivializations, each finite fiber jet of $f_{\gamma,x}$
depends continuously on $x$. This is a condition on the original
fiber maps, before any measurable change of coordinates.

\begin{theorem}[Centered fibers over a compact base]
\label{I-thm:centered-fibers}
Let $X$, $E$, and $f$ be as above, and assume that the centered
analytic fiber cocycle has continuous finite jets. Let $\mu$ be
an ergodic invariant Borel probability on $X$. There exist $A$,
$C$, and $\kappa$ as in Theorem~\ref{I-thm:compact-actions}, and a
measurable family of analytic germs $F_x:(\R^n,0)\to(E_x,0)$,
with invertible derivative, such that
\begin{equation}\label{I-eq:fiber-normal-form}
 f_{\gamma,x}\circ F_x=F_{\gamma x}\circ A(\gamma)\kappa(\gamma,x)
                       \qquad(\gamma\in\Gamma)
\end{equation}
on one invariant conull set. The germs and their inverses have
positive measurable radii in local bundle trivializations.
\end{theorem}

A jointly real-analytic bundle action near the zero section over
a compact analytic base satisfies the hypotheses. Continuity of
the finite jets is enough; a uniform radius for the original
fiber maps is not part of the statement. Compactness alone would
not suffice if those jets depended only measurably on the base.

\subsection{Cocycle superrigidity for analytic germs}
Write $J^N(\R^n)$ for the real algebraic group of invertible
$N$-jets fixing zero. Fix a finite symmetric generating set $S$
of $\Gamma$. Let $\iota_N$ be the faithful algebraic substitution
realization of $J^N(\R^n)$ on truncated polynomials. We use this
fixed realization in the moment condition. A germ cocycle has
\emph{integrable finite jets} if
\begin{equation}\label{I-eq:ha-general-finite-jet-moments}
 \int_X\log\left(1+\|\iota_N(j^N\beta(s,x))\|
              +\|\iota_N(j^N\beta(s,x))^{-1}\|\right)d\mu(x)<\infty
                  \quad(s\in S,\ N\ge1).
\end{equation}
There is no bound uniform in $N$ in this condition. Changing the
finite generating set preserves it, by the cocycle identity and
invariance of the measure.

\begin{theorem}[Analytic germ cocycle superrigidity]
\label{I-thm:ha-general-integrable-superrigidity}
Let $\Gamma\curvearrowright(X,\mu)$ be an ergodic probability-preserving
action on a standard probability space. Suppose
$\beta:\Gamma\times X\to\Diff^\omega(\R^n,0)$ is a
coefficient-measurable cocycle with integrable finite jets.
There exist an algebraic representation $A:G\to\GL(n,\R)$,
a compact group $C\le Z_{\GL(n,\R)}(A(G))$, a measurable
$C$-valued cocycle $\kappa$, and measurable analytic invertible
germs $F_x$ fixing zero, such that
\begin{equation}\label{I-eq:general-normal-form}
 F_{\gamma x}^{-1}\circ\beta(\gamma,x)\circ F_x
                           =A(\gamma)\kappa(\gamma,x)
                           \qquad(\gamma\in\Gamma)
\end{equation}
on one invariant conull set. The transfer and its inverse have
positive measurable radii. No moments are imposed on $F_x$ or on
any intermediate coordinate change.
\end{theorem}

Theorems~\ref{I-thm:compact-actions} and \ref{I-thm:centered-fibers}
follow by applying this theorem in bounded measurable local
trivializations. We verify the finite-jet bounds and all changes
of charts in Section~\ref{I-sec:applications}.

The analytic part of the argument is a formal-to-analytic
realization theorem. Once the derivative has the form
$L(\gamma,x)=A(\gamma)\kappa(\gamma,x)$, a compatible formal
transfer can be normalized within its equivariant formal
centralizer so that it converges. This implication needs no
moment assumption; it is proved in
Theorem~\ref{I-thm:ha-compact-core-realization}. The moment condition
in Theorem~\ref{I-thm:ha-general-integrable-superrigidity} is used
to obtain the formal reduction.

\subsection{Relation to earlier work}
We use Fisher--Margulis cocycle superrigidity
\cite[Definition~1.3 and Theorem~1.5]{FisherMargulis} for each algebraic
jet group. It applies to uniform and nonuniform lattices and retains
the compact centralizing factor. We assemble the finite-order reductions
into compatible formal coordinates for the convergence theorem.

Nevo's spectral transfer theorem \cite[Theorem~1 and Section~2.2]{Nevo1998}
and the spherical bounds of Cowling--Haagerup--Howe
\cite[Theorem~2, pp.~99--100]{CHH} and Gorodnik--Nevo
\cite[Definition~5.1 and Theorem~5.4(1), pp.~38--40 of the author
version]{GNAuthor2007} supply the matrix-coefficient estimates. We use one spherical averaging law
for all coefficient degrees. The coefficient modules also occur
in Fisher--Hitchman's Bochner treatment of a finite-dimensional
representation tensored with a unitary representation \cite{FH,FHB}.

Nonstationary normal forms for contracting cocycles
\cite{Melnick,Kalinin} use the fiber spectrum, and the analytic-germ
Liv\v{s}ic theorem of Navas--Ponce \cite{NavasPonce} uses periodic
orbit data. Here the finite-jet reductions give compatibility
over a measurable action. Our convergence proof normalizes the
invariant coefficient spaces, so it allows
$\Hom_G(\Sym^d\R^n,\R^n)\ne 0$ in degrees $d\ge2$.

\subsection{Organization of the proof}
Section~\ref{I-sec:formal} obtains compatible formal coordinates
from finite-jet superrigidity. Sections~\ref{I-sec:averaging} and
\ref{I-sec:localization} give coefficient estimates on a common
measurable set. Section~\ref{I-sec:convergence} constructs the
analytic transfer, first for an algebraic derivative and then
with a compact centralizing cocycle. Section~\ref{I-sec:applications}
deduces the main results.

For the standard representation-theoretic facts concerning compact
forms, complete reducibility, and the Casimir operator, we refer to
Knapp \cite{KnappBeyondDigital}.

\section{Finite jets and compatible formal coordinates}\label{I-sec:formal}
We use faithful matrix realizations of the source groups. Such a realization
exists by \cite[Theorem~4.9 and Corollary~4.10]{Milne}. The real and complex matrix
Lie algebras are related by scalar extension; semisimplicity agrees with
semisimplicity of the algebraic group in characteristic zero
\cite[Lemma~23.68]{Milne}. Geometric algebraic simple connectivity gives topological
simple connectivity of the complex points by the semisimple comparison
\cite[\S4, footnote~1, p.~149]{RoquesSinger2022}.
For each actual almost-simple factor, the rank assumption gives a closed
split torus of dimension at least $2$. Its inclusion has injective
differential. Restrict the actual adjoint representation to this torus:
the split-torus weight decomposition gives a common real weight basis
\cite[Definition~12.11 and Theorem~12.12, pp.~234--235]{Milne}.
Differentiation, using $d\operatorname{Ad}=\operatorname{ad}$
\cite[\S10, paragraph~10.6, p.~188]{Milne}, therefore gives two independent
commuting tangent directions whose adjoint actions have a common real
eigenbasis.
More explicitly, if these directions are $A,B$, their real span $U$ has
$\dim U=2$. Bilinearity and $[A,B]=0$ give
$[sA+tB,uA+vB]=0$, so this same $U$ is an abelian Lie subalgebra.
If $[A,e_i]=a_i e_i$ and $[B,e_i]=b_i e_i$ in the common basis, then
\[
 [sA+tB,e_i]=(sa_i+tb_i)e_i \qquad(s,t\in\R).
\]
Thus the same basis diagonalizes every adjoint operator from $U$.
In the adjoint description of real split rank, one takes the largest
dimension of a subalgebra admitting such a common real eigenbasis
\cite[Chapter~5, \S4.1, pp.~268--270]{OnishchikVinberg1990}.
For a real Lie algebra $L$, the possible dimensions form a nonempty subset
of $\{0,\ldots,\dim L\}$, so this largest dimension is attained.
The constructed $U$ therefore gives the numerical lower bound $2$ in
every simple factor. Identifying this Lie-algebra rank with the rank of
the algebraic source still requires the structural comparison above.
These identifications let us use matrix exponentials and polynomial
coordinates without changing the standing source group.

Here is the coordinate-ring convention behind this realization. Write
$x_{ij},y_{ij}$ for the forward and inverse matrix entries, and let $I$
be the ideal of real polynomials in these variables which vanish on all
complex points of the source. The algebra $\R[x_{ij},y_{ij}]/I$ is of
finite type and reduced: if $p^m$ vanishes at every complex point, then
$p$ does as well. Each original complex point gives an evaluation of
this algebra in $\C$, and these evaluations distinguish points.
Conversely, an $\R$-algebra homomorphism to $\C$ gives matrices $X,Y$
by evaluating the coordinate functions. The relations $XY=1$ belong
to $I$, so $X$ is invertible with inverse $Y$. The original defining
equations also belong to $I$, so $X$ is an original source point.
Evaluation at this point recovers the homomorphism, since the coordinate
functions generate the algebra. This identifies the complex points in
the coordinate-ring and matrix descriptions. The same argument over $\R$
identifies the real points as well. Give these evaluation spaces the
topology of pointwise evaluation, using the usual topology on $\R$ or
$\C$. A coordinate-ring element has a polynomial representative, so its
evaluation is continuous in the original matrix topology. Conversely,
the entries of the recovered matrix and its inverse are evaluations of
the coordinate functions. The point correspondences are therefore
homeomorphisms. This concerns the usual topology on real and complex
points. The tangent, simple connectivity, and rank comparisons stated
above are additional steps.

\subsection{Germs, measurability, and radii}
All probability spaces are standard. Actions of Lie groups are jointly
measurable, and coefficient measurability for their cocycles refers to
the product of the group and the base. For Hilbert-valued sections,
measurability means strong (Bochner) measurability; no separability
of the ambient Hilbert space is assumed. In finite-dimensional spaces
this agrees with ordinary Borel measurability. Composition is written from
right to left. A germ fixes zero and has invertible derivative.
Coefficient measurability means that every Taylor coefficient is
measurable; it does not specify a common representative domain.
All countably many lattice cocycle identities will be taken on a
single invariant conull set.

For radius selection write $F_x(z)=\sum_{d\ge0}p_{x,d}(z,\ldots,z)$,
where $p_{x,d}$ is the symmetric multilinear Taylor coefficient,
normalized by $d!$, and $p_{x,0}=0$. In this paragraph use the
multilinear operator norm. The series has positive convergence radius
if and only if there is an integer $j\ge0$ such that
\[
                \|p_{x,d}\|\le(j+1)^{d+1}\qquad(d\ge0).
\]
Indeed, a power series of positive radius has a bound $CR^d$ for
its coefficients. Symmetrization does not increase the multilinear
operator norm, so the same bound holds for its normalized Taylor
coefficients; increasing $j$ gives the displayed bound. Conversely,
that bound gives absolute convergence on the ball of radius
$[2(j+1)]^{-1}$. The canonical sum has the original Taylor
coefficients and represents the germ whenever the germ is analytic.

For fixed $j$ the displayed condition is a countable intersection
of measurable sets. Its union over $j$ is the measurable convergence
locus. On that locus choose the least admissible $j=j(x)$.
Its level sets are measurable, and
$r(x)=[2(j(x)+1)]^{-1}$ is a positive measurable radius.
The symmetric-tensor Hilbert norms used later are for the coefficient
estimates; the radius selection here uses the operator norms just specified.

Inversion preserves coefficient measurability. To see this at finite
order, identify the jet group with its standard Borel space of finite
Taylor profiles. Multiplication is measurable. The measurable bijection
$(a,b)\mapsto(a,ab)$ of the product jet space has measurable inverse;
applying this inverse to $(a,1)$ proves measurability of $a\mapsto a^{-1}$.
Projection to each Taylor coefficient proves the assertion for formal
and analytic germ fields. An analytic germ with invertible derivative
has an analytic inverse by the inverse function theorem. Apply the same
integer-majorant selection to the inverse field to obtain its positive
measurable radius.

To compare two faithful
algebraic matrix realizations, include the entries of each matrix and
its inverse as coordinates. The inverse of either closed immersion
onto its image is regular \cite[Corollary~3.35 and Theorem~4.9]{Milne},
so each coordinate in one realization is a
polynomial in these coordinates of the other. For a polynomial $p$
there are constants $C,N$ such that
\[
 1+|p(z)|\le C\prod_i(1+|z_i|)^N.
\]
Taking logarithms shows that integrability of the coordinate logarithms
is preserved by this change. In a finite-dimensional matrix space,
coordinate logarithmic moments are equivalent to the matrix-norm
moment. Applying the comparison in both directions proves independence
of the faithful algebraic realization.

\subsection{The algebraic input}
We use the following real-target case of Fisher--Margulis
\cite[Theorem~1.5]{FisherMargulis}. Under our hypotheses on $G$ and $\Gamma$,
an integrable Borel cocycle over an ergodic probability-preserving
$\Gamma$-action, with values in the real points $L$ of an
algebraic group, is cohomologous \emph{in $L$} to
\[
                  \pi(\gamma)c(\gamma,x).
\]
Here $\pi:G\to L$ is a continuous homomorphism and $c$ takes
values in a compact subgroup centralizing $\pi(G)$. For a
discrete finitely generated group, their integrability condition
follows from the generator logarithmic bounds in
\eqref{I-eq:ha-general-finite-jet-moments}: write a group element
as a word, use submultiplicativity, and integrate translated
generator bounds. Compact subsets of $\Gamma$ are finite.

We regard an algebraic subgroup through its real points. When choosing a
complex algebraic presentation, replace it if necessary by the Zariski
closure of these points. This leaves the real subgroup unchanged
\cite[Theorem~1.45]{Milne}. Its identity component is unchanged, since
real points of a connected real algebraic group are Zariski dense
\cite[Theorem~17.93]{Milne}; hence reductivity is preserved.

\begin{lemma}[Linearization of a reductive jet group]
\label{I-lem:reductive-jets}
For $N\ge1$, every reductive real algebraic subgroup $L<J^N(\R^n)$ is
conjugate into the linear jet subgroup by a tangent-to-identity
$N$-jet.
\end{lemma}
\begin{proof}
Let $\mathcal Q_N$ be the space of $\R^n$-valued polynomials
of degrees $1$ through $N$. Let $L$ act by
\[
 (\ell\cdot P)(z)=D_0\ell\,P(\ell^{-1}z)\pmod{\deg>N}.
\]
Taking the linear part is an equivariant surjection
$\mathcal Q_N\to\operatorname{End}(\R^n)$. This is a rational
representation of the reductive algebraic group $L$ in characteristic zero;
complete reducibility, including disconnected $L$, follows from
\cite[Theorem 22.42 and Corollary 22.43]{Milne}. It lifts the invariant identity endomorphism to an invariant $P$
with $DP(0)=I$. Thus $P\ell=D_0\ell\,P$.
The jet $P$ is invertible and gives the required conjugacy.
\end{proof}

\begin{proposition}[Compatible formal reduction]
\label{I-prop:formal-reduction}
Under the hypotheses of
Theorem~\ref{I-thm:ha-general-integrable-superrigidity}, there is a
coefficient-measurable tangent-to-identity formal field $H_x$ with
\[
             \beta(\gamma,x)H_x=H_{\gamma x}a(\gamma,x).
\]
There are also a measurable linear field $Q_x$, an algebraic
representation $A:G\to\GL(n,\R)$ and a compact centralizing
cocycle $\kappa$ such that
\[
       Q_{\gamma x}a(\gamma,x)Q_x^{-1}
                        =A(\gamma)\kappa(\gamma,x).
\]
No integrability of $H_x$ or $Q_x$ is asserted or needed.
\end{proposition}
\begin{proof}
For each $N$, apply the stated cocycle theorem to the original
cocycle $j^N\beta$ in the algebraic jet group. On truncated
polynomials the resulting semisimple homomorphic factor and compact
factor commute. Differentiate the continuous semisimple factor and apply
Weyl's complete reducibility theorem
\cite[Lemma~B.1, pp.~659--660]{KnappBeyondDigital}. Choose a Lie-invariant complement to the given invariant subspace.
Projection onto that subspace along the complement commutes with the
differential. It therefore commutes with the exponential of every
differential operator, and hence with the connected group, which is generated
by exponential neighborhoods \cite[Corollary~1.85]{KnappBeyondDigital}.
This gives an equivariant projection onto the square of the maximal ideal.
Average this projection over the compact factor. Commutation preserves
semisimple equivariance, and the averaged projection still restricts to
the identity on that square. Its kernel therefore gives an invariant
lift of the linear coordinate functions. As in the proof of
Lemma~\ref{I-lem:reductive-jets}, these functions define a fixed tangent
jet that linearizes both factors.

Combining this fixed jet with the measurable transfer yields
$T_x^{(N)}$ and a linear cocycle $b^{(N)}$ satisfying
\[
 T_{\gamma x}^{(N)}j^N\beta(\gamma,x)
                      =b^{(N)}(\gamma,x)T_x^{(N)}.
\]
Put $c_x^{(N)}=DT_x^{(N)}(0)$ and
$H_x^{(N)}=(T_x^{(N)})^{-1}c_x^{(N)}$. Differentiation shows
that $H^{(N)}$ linearizes to the actual derivative $a$ and is
tangent to the identity.

Suppose $H^{(N)}$ has been chosen and let
$\widetilde H^{(N+1)}$ be any next-order choice. Their discrepancy
\[
 E_x=(j^N\widetilde H_x^{(N+1)})^{-1}H_x^{(N)}
\]
satisfies $E_{\gamma x}=a(\gamma,x)E_xa(\gamma,x)^{-1}$.
Represent it by its degree-$N$ polynomial and set the next
coefficient equal to zero. This lift remains equivariant,
since linear conjugation preserves homogeneous degree.
Right composition of $\widetilde H^{(N+1)}$ by the lift restores
compatibility. Induction gives the formal field $H$.

Apply the same cocycle theorem at degree one to obtain $Q$, $A$,
and $\kappa$. The continuous homomorphism $A$ is algebraic by
\cite[Lemma~3.4, p.~199, with the standing hypotheses of \S3.1]{FisherMargulis}.
Indeed, geometric simple connectivity gives the required semisimple-cover
condition, and our factorwise real-rank assumption excludes compact factors.
If $A$ is nontrivial, the source has an almost-simple factor; a split torus
in that factor gives total real rank at least $2$. Thus the lemma applies
to the original real algebraic group and gives an algebraic homomorphism
inducing this same $A$. If $A$ is trivial, use the constant identity
homomorphism. Pullback of matrix coordinates along this algebraic map,
followed by lifting through the original coordinate quotient, gives real
polynomial representatives in the fixed faithful presentation.
This proves the proposition.
Every integrability assertion used here concerns an original
finite-jet cocycle, not a cocycle transformed at an earlier step.
\end{proof}
\section{Averaging on coefficient representations}\label{I-sec:averaging}

\subsection{A uniform compact-dual potential}
The averaging and realization arguments below apply to finite-dimensional
continuous real representations of connected semisimple Lie groups.
Such representations are smooth \cite[I.10, p.~77]{KnappBeyondDigital}. We state
these auxiliary results in that generality; the algebraic representations
in the main theorem are particular examples.

An inner product for a finite-dimensional continuous real representation $R$ is
\emph{Cartan-compatible} if $R(K)$ is orthogonal and $dR(Y)$
is self-adjoint for $Y\in\mathfrak p$. Such an inner product
is obtained as follows. Complexify the differential and restrict it to
$\mathfrak k\oplus i\mathfrak p$. Integrate this representation to the
simply connected compact group by \cite[Theorem~20.19]{LeeSmoothManifolds2e}
and average a positive Hermitian form over that group. Equivalently,
average its real part $B$. The complex structure $J=i$ commutes with
the representation, so $B(Jx,Jy)=B(x,y)$. Differentiating invariance
along a compact one-parameter subgroup gives skewness of its infinitesimal
operator. In particular $JP$ is skew when $P=dR(X)$ and
$X\in\mathfrak p$. Since $PJ=JP$ and $J^2=-I$,
\[
 B(Px,y)=B(JPx,Jy)=-B(x,JPJy)=B(x,Py).
\]
Restricting this same positive form to the original real space gives
the required inner product. This is the unitary-trick argument of
\cite[Propositions~4.6 and~7.15]{KnappBeyondDigital}. Constants in the next
lemma are uniform over representation types; orthogonal sums
and tensoring with Hilbert multiplicity spaces preserve them.

\begin{lemma}[Uniform positive-frame potential]
\label{I-lem:ha-fixed-positive-potential}
Let $G$ be connected semisimple with finite center and no compact
factors. Fix a maximal compact subgroup $K$ and the corresponding
Cartan decomposition $\mathfrak g=\mathfrak k\oplus\mathfrak p$.
Let $H\in\mathfrak p$ have $K$-orbit spanning $\mathfrak p$.
There exist $T_0,\sigma>0$ such that every Cartan-compatible
finite-dimensional continuous real representation $R$ with $R^G=0$
satisfies
\[
 \int_K R(k)\tanh^2(T\,dR(H))R(k)^*\,dk\ge\sigma I
                                  \qquad(T\ge T_0).
\]
The constants do not depend on its irreducible types or their
Hilbert multiplicities.
\end{lemma}
\begin{proof}
Use a compatible compact real form
$\mathfrak u=\mathfrak k\oplus i\mathfrak p$ in the complexified
Lie algebra \cite[Corollary~6.18]{KnappBeyondDigital}. Its Killing form is negative
definite. It therefore integrates to a connected compact group
\cite[Theorem~3.4]{KnappStructure1997}. Take its simply connected
cover \cite[Theorem~7.7]{LeeSmoothManifolds2e}; this cover is compact
by Weyl's theorem \cite[Theorem~4.69]{KnappBeyondDigital}. A faithful matrix
representation realizes the same compact group as a closed matrix group
\cite[Corollary~4.22]{KnappBeyondDigital}. Throughout, the source charts are the
actual matrix exponential charts, and the inclusions of the real and
compact Lie algebras retain their original differentials.

For each original real factor $G_i$, let $K_i$ be its Cartan compact
subgroup and let $U_i$ be the simply connected compact dual just
constructed. Integrate the original complexified Lie inclusion to a
homomorphism $\rho_i:U_i\to\mathrm{GL}_N(\mathbb C)$ using
\cite[Theorem~20.19]{LeeSmoothManifolds2e}. Its kernel is finite because
its differential is injective and its domain is compact. Let $B_i$ be
the matrix image of $K_i$. The exponential identities show that
$B_i\subset\rho_i(U_i)$. The identity component
$V_i=(\rho_i^{-1}(B_i))^0$ maps onto $B_i$ with finite kernel.
Now put $\widehat K_i=(K_i\times_{B_i}V_i)^0$. Both projections
from this pullback are onto and have finite kernels. Taking the product
over the original factors gives a finite cover of $K$ together with a
homomorphism into $\prod_i U_i$. The exponential identities hold for
these same maps: they follow from their actual differentials and
uniqueness of one-parameter subgroups.
Thus the same exponential and adjoint identities hold on this common
cover; normalized Haar measure pushes forward to normalized Haar measure.

Choose $k_1,\ldots,k_m$ so that $X_j=\operatorname{Ad}(k_j)H$
form a basis of $\mathfrak p$. In the compact dual Lie algebra
$\mathfrak u=\mathfrak k\oplus i\mathfrak p$, put $P_j=iX_j$.
The $P_j$ and their brackets span $\mathfrak u$, since
$[\mathfrak p,\mathfrak p]=\mathfrak k$.
Choose pairs $(a_l,b_l)$ whose brackets $[P_{a_l},P_{b_l}]$
form a basis of $\mathfrak k$. For one sufficiently small fixed
$t>0$, the word map
\[
 \Psi(u,v)=\prod_j e^{u_jP_j}
      \prod_l(e^{tP_{a_l}}e^{v_lP_{b_l}}
                       e^{-tP_{a_l}}e^{-v_lP_{b_l}})
\]
is a local diffeomorphism at $(0,0)$. Indeed its derivative has
columns $P_j$ and $t[P_{a_l},P_{b_l}]+O(t^2)$, which are linearly
independent for small $t$. Shrink the parameter neighborhood so that
all exponential parameters have absolute value at most one. The image
contains an identity neighborhood in the compact connected group $U$.
An open symmetric neighborhood inside this image generates $U$;
compactness then gives an integer $M$ such that every element of $U$
is a product of at most $M$ factors $\exp(\theta P_j)$, with
$|\theta|\le1$. All these choices depend only on the source group and $H$.

The compact-dual Lie algebra representation integrates to a representation
$\pi$ of $U$ by \cite[\S I.10, pp.~74--76]{KnappBeyondDigital}; here $U$ is
simply connected and its differential is the compact-dual differential
constructed above. This representation is unitary because its differential
is skew-Hermitian and $U$ is connected.
Its invariant space is zero, since the compact dual and the original
Lie algebra have the same invariant vectors. Haar averaging thus gives
\[
 \int_U\|\pi(u)w-w\|^2\,du=2\|w\|^2.
\]
For real $a$ and $|\theta|\le1$,
$|e^{i\theta a}-1|\le4|\tanh a|$.
The spectral theorem and telescoping a word of length at most $M$,
whose prefixes are unitary, give
\[
 \|\pi(u)w-w\|^2
 \le16M^2\sum_{j=1}^m\|\tanh(dR(X_j))w\|^2.
\]
Integrating over $U$ proves
$\sum_j\tanh^2(dR(X_j))\ge(8M^2)^{-1}I$.
Conjugate by $R(k)$ and integrate over $K$. Each summand has the
same average, so the assertion holds with $T_0=1$ and
$\sigma=(8M^2m)^{-1}$. The argument and constants are unchanged
under orthogonal sums and Hilbert multiplicities. Finally,
monotonicity of $\tanh^2(TA)$ for self-adjoint $A$ gives every $T\ge1$.
\end{proof}

\subsection{Spherical averaging and a uniform coefficient gap}
\label{I-subsec:ha-bik-mixing-coercivity}

We use spherical averages whose scalar operators have small norm
on the orthogonal complement of the constants. The bounded
normalization of an arbitrary matrix retains a uniform algebraic gap.
The assertion below concerns the degree-zero twisted differential.
Section~\ref{I-sec:localization} uses it to bound an existing measurable
primitive on a common retained set.

\begin{proposition}[A common averaging law]
\label{I-prop:ha-bik-common-law}
Let $G$ be connected semisimple with finite center and no compact
factors, let $K$ be a maximal compact subgroup, and fix a maximal
split Cartan subalgebra $\mathfrak a\subset\mathfrak p$. Choose a
regular $H\in\mathfrak a$ with nonzero projection to every simple factor.
Its $K$-orbit spans $\mathfrak p$.
For $T>0$, put $a_T=\exp(TH)$ and
\[
 \nu_T=m_K*\tfrac12(\delta_{a_T}+\delta_{a_T^{-1}})*m_K.
\]
There are $T_0,c_0>0$, depending on the fixed $G,K,H$, such that the following
holds for every $T\ge T_0$, every probability-preserving $G$-action,
and every Cartan-compatible finite-dimensional continuous real
representation $R$ with $R^G=0$, and their orthogonal sums with
Hilbert multiplicities. Set
\[
 L_g=(I+R(g)R(g)^*)^{-1/2},\qquad
 (D_gv)(x)=L_g\bigl(v(gx)-R(g)v(x)\bigr).
\]
Then
\begin{equation}\label{I-eq:ha-bik-twisted-gap}
 \int_G\|D_gv\|_2^2\,d\nu_T(g)\ge c_0\|v\|_2^2.
\end{equation}
For an orthogonal sum, let $L_g$ and $N_g$ denote the bounded
block diagonal extensions of the finite-dimensional maps $L_g$ and
$L_gR(g)$, respectively, and interpret the differential as
$D_gv(x)=L_gv(gx)-N_gv(x)$. This agrees with the displayed formula
on the natural domain of the unnormalized representation.
The two endpoint coefficients satisfy $\|L_g\|\le1$ and
$\|N_g\|\le1$. Moreover
Define the bounded differential
\[
 D_T:L^2(X;\mathcal H)\longrightarrow
 L^2(X\times G,\mu\otimes\nu_T;\mathcal H),\qquad
 (D_Tv)(x,g)=D_gv(x).
\]
On $L^2$ put $Q_{T,2}=I-D_T^*D_T$. Thus the notation
$\int D_g^*D_g\,d\nu_T(g)$ denotes its weak integral, determined
by scalar products. The four-block normal-kernel formula defines
consistent operators $Q_{T,p}$ on the other $L^p$ spaces. These operators
have norm at most
$\sqrt{1-\sigma}<1$ on every $L^p$, $1\le p\le\infty$,
for a constant $\sigma>0$ independent of the representation.

If all simple factors of $G$ have property~\textup{(T)} and the action
is ergodic, the scalar averaging operators additionally satisfy
\begin{equation}\label{I-eq:ha-bik-scalar-mixing}
 \left\|\int_G U_g\,d\nu_T(g)\right\|_{L^2_0}\longrightarrow0
                  \qquad(T\longrightarrow\infty).
\end{equation}
Here $U_gv(x)=v(gx)$ is deliberately the Koopman
\emph{antirepresentation}. This convention makes the factorization
in the proof unambiguous.
\end{proposition}

\begin{proof}
First, there is $\sigma>0$ such that, for every Cartan-compatible
representation without fixed vectors and $T\ge T_0$,
\begin{equation}\label{I-eq:ha-bik-K-average-gap}
 \int_K R(k)\tanh^2(T\,dR(H))R(k)^*\,dk\ge\sigma I.
\end{equation}
This is Lemma~\ref{I-lem:ha-fixed-positive-potential}. Orthogonal sums and Hilbert multiplicities preserve its lower bound.

Write $A=R(a_T)$, a positive self-adjoint matrix, and
\[
 C=A(I+A^2)^{-1}=\bigl(2\cosh(T\,dR(H))\bigr)^{-1}.
\]
The operator
\[
 P=\int_K R(k)^*U_k\,dk
\]
is an orthogonal projection: $k\mapsto R(k)^*U_k$ is a unitary
antirepresentation. For $w\in\operatorname{ran}P$, one has
$R(k)^*w=U_k^{-1}w$. Consequently the quadratic form of $C^2$ at $w$
is unchanged if $C^2$ is replaced by its $K$-conjugacy average.
Since $4C^2=I-\tanh^2(T\,dR(H))$, the preceding gap gives
\[
                   \|CP\|\le\tfrac12\sqrt{1-\sigma}.
\]

Expand the integrated squared residual. Its two diagonal terms
sum to $I$: for $g=k_1a_T^{\pm1}k_2$, the averages over the two signs
of $(I+R(g)R(g)^*)^{-1}$ and of
$R(g)^*(I+R(g)R(g)^*)^{-1}R(g)$ are each $I/2$.
The mixed term is governed by
\[
 B=\int R(g)^*(I+R(g)R(g)^*)^{-1}U_g\,d\nu_T(g)
     =P\,C\,\tfrac12(U_{a_T}+U_{a_T^{-1}})\,P.
\]
For clarity, for the positive sign the integrand factors as
\[
 [R(k_2)^*U_{k_2}]\,[CU_{a_T}]\,[R(k_1)^*U_{k_1}],
\]
because $U_{k_1a_Tk_2}=U_{k_2}U_{a_T}U_{k_1}$.
Thus $\|B\|\le\frac12\sqrt{1-\sigma}$, and
\[
 \int\|D_gv\|_2^2\,d\nu_T
   =\|v\|_2^2-2\operatorname{Re}\langle Bv,v\rangle
   \ge(1-\sqrt{1-\sigma})\|v\|_2^2.
\]
Take $c_0=1-\sqrt{1-\sigma}$.

The same factorization gives the all-$L^p$ assertion, but its
$L^\infty$ proof uses $PC$, rather than $CP$. At each base point,
operator-valued Cauchy--Schwarz gives
\[
 \|PCf(x)\|^2
 =\left\|\int_K R(k)^*C f(kx)\,dk\right\|^2
 \le\frac{1-\sigma}{4}\int_K\|f(kx)\|^2\,dk.
\]
Here the row-operator bound is exactly
$\int_KR(k)^*C^2R(k)\,dk\le(1-\sigma)I/4$.
Thus $2PC$ contracts on $L^p$ for $p\ge2$ by Jensen,
including $p=\infty$. The remaining factors
$(U_{a_T}+U_{a_T^{-1}})/2$ and $P$ are contractions on every
$L^p$. Since $Q_T=2B$ is self-adjoint, duality gives the same
bound for $1\le p\le2$.

Finally, \eqref{I-eq:ha-bik-scalar-mixing} is the standard spherical
averaging estimate for semisimple Kazhdan groups: the norm of
$m_K*\delta_{\exp(TH)}*m_K$ on the orthogonal complement of the
invariant vectors decays exponentially when every simple-factor component of $H$ is nonzero. We use Nevo \cite[Theorem~1 and Section~2.2]{Nevo1998}.
For the singular measure appearing here, use the regular coefficient
estimate of Gorodnik and Nevo
\cite[Definition~5.1 and Theorem~5.4(1), pp.~38--40 of the author version]{GNAuthor2007}.
It applies to unimodular locally compact second countable groups with an
Iwasawa decomposition $G=KP$, where $K$ is compact and $P$ is closed
and amenable. The original semisimple group is unimodular
\cite[Corollary~8.31(c), p.~535]{KnappBeyondDigital}, and its full minimal parabolic
subgroup is amenable
\cite[\S2.2.2 and Remark~2.12, pp.~12 and~15 of the author version]{BenoistKobayashiII}.
The normalized spherical function is the integral of the inverse square
root of the left-$K$ extension of the actual modular character of $P$.
These are the same subgroup, modular function and Haar normalization
used here. The algebraic-group estimate of Cowling--Haagerup--Howe
\cite[Theorem~2, pp.~99--100, formula~(8)]{CHH} is the antecedent of this
version. If $P_K$ denotes averaging over $K$ in the regular representation
$\lambda$, the estimate on $K$-fixed vectors gives
\[
 \|P_K\lambda(a)P_K\|\le\Xi(a).
\]
Together with the bound for $\Xi$ in Nevo's Section~2.2, this gives
decay along the whole real ray on the original group, including finite
central covers.
Here $K$ is the Cartan maximal compact subgroup. Conjugating each
component into a closed positive chamber preserves the spherical average;
the corresponding $\rho_i(H_i)$ is positive for every nonzero component,
including directions on chamber walls. Let $P_i$ be the orthogonal projection onto the invariant vectors of
the $i$th simple factor, and let $A_i(T)$ be its spherical averaging
operator. The factor actions commute, hence so do their invariant
projections, and $\prod_iP_i$ projects onto the global invariants.
Nevo's estimate gives $\|A_i(T)-P_i\|\to0$. Since all these operators
are contractions, telescoping yields
\[
 \left\|\prod_i A_i(T)-\prod_iP_i\right\|
       \le\sum_i\|A_i(T)-P_i\|\longrightarrow0.
\]
Every component of $H$ escapes. Apply this argument to both rays
and average to obtain the asserted scalar limit.
\end{proof}

\section{Estimates on a common measurable set}\label{I-sec:localization}

\subsection{Pointwise control of omitted edges on a large core}
\label{I-subsec:ha-pointwise-pruning}

\begin{proposition}[Pruning for a strongly mixing reversible kernel]
\label{I-prop:pruning}
Let $(X,\mu)$ be an atomless standard probability space. Let
$x\mapsto J_x$ be a Markov kernel such that the edge measure
$J(dx,dy)=\mu(dx)J_x(dy)$ is invariant under $(x,y)\mapsto(y,x)$.
Its two marginals are $\mu$. Let $P$ be its Markov operator. Assume
\[
 \|P|_{L^2_0(X)}\|\le\varepsilon.
\]
Let $r:X\times X\to\{0,1\}$ be a symmetric measurable retained-edge
indicator, and put
\[
 b(x)=\int(1-r(x,y))\,J_x(dy),\qquad
 \beta=\int b\,d\mu.
\]
Fix $0<\theta<1$ with $\varepsilon\le\theta/2$ and
$\beta<\theta^2/4$. There is a measurable set $B\subseteq X$ with
\[
 \mu(B)\ge1-\frac{2\beta}{\theta}
\]
such that, for almost every $x\in B$,
\begin{equation}\label{I-eq:pruning-omission}
\int \bigl(1-r(x,y)1_B(y)\bigr)J_x(dy)\le\theta.
\end{equation}
Thus the total omitted outgoing edge mass is bounded pointwise on
$B$, not merely on average.
The same conclusion holds for a reversible labelled edge measure,
with $r$ allowed to depend on the edge label.
\end{proposition}

\begin{proof}
Let $D_0=\varnothing$. At step $i$ consider the eligible set
\[
 E_i=\left\{x\notin D_i:
        b(x)+\int_{D_i}r(x,y)J_x(dy)>\theta\right\}.
\]
If it is null, stop. Otherwise choose a measurable $A_i\subseteq E_i$
with
\[
 \frac12\min\{\mu(E_i),\theta/2\}
 \le\mu(A_i)\le\min\{\mu(E_i),\theta/2\},
\]
and put $D_{i+1}=D_i\cup A_i$. Such subsets exist by atomlessness.
Summing the defining inequality over the disjoint removed sets and
using symmetry gives, for $t_d=\mu(D_d)$,
\begin{equation}\label{I-eq:pruning-charging}
\theta t_d
 \le\beta+\frac12 J(D_d\times D_d)
 =\beta+\frac12\langle P1_{D_d},1_{D_d}\rangle
 \le\beta+\frac12(t_d^2+\varepsilon t_d).
\end{equation}
The factor $1/2$ is justified because an edge between different
removed sets is counted only from the later set to the earlier one;
symmetry gives the reverse edge the same mass. Edges within a single
removed set only increase the right side.

For $t\in[\theta/2,\theta]$ and $\varepsilon\le\theta/2$,
\[
 (\theta-\varepsilon/2)t-t^2/2\ge\theta^2/4.
\]
Thus \eqref{I-eq:pruning-charging} and $\beta<\theta^2/4$ forbid $t_d$ from this interval.
Since each added set has measure at most $\theta/2$, the sequence
cannot jump from below $\theta/2$ to above $\theta$.
Consequently $t_d<\theta/2$ for all $d$, and \eqref{I-eq:pruning-charging} then implies
\[
 \frac\theta2 t_d\le\beta.
\]
Let $D=\bigcup_nD_n$ and $B=X\setminus D$. The same estimate yields
$\mu(D)\le2\beta/\theta$.

If the left side of \eqref{I-eq:pruning-omission} exceeded $\theta$ on a positive-measure
subset of $B$, monotone convergence of the omitted masses would
imply that $\mu(E_i)$ is bounded below by a positive constant for
all sufficiently large $i$. The selection rule would force
$\mu(A_i)$ to be bounded below, contradicting the disjointness of
these sets in a probability space. Hence \eqref{I-eq:pruning-omission} holds almost everywhere.
All computations use only symmetry, disintegration, and the scalar
Markov bound, so labelled edges obey the identical proof.
\end{proof}

\begin{corollary}[Stability of a bounded twisted zero-order gap]
\label{I-cor:gap-deletion}
In the preceding setting, suppose that an edgewise linear operator
has the form
\[
 D_{x,y}f=A_{x,y}f(y)-C_{x,y}f(x),
 \qquad \|A_{x,y}\|,\|C_{x,y}\|\le1,
\]
and that its full quadratic form satisfies
\[
 \int\|D_{x,y}f\|^2\,dJ(x,y)\ge c\|f\|_{L^2(X)}^2.
\]
Assume its squared edge norm is unchanged under edge reversal.
For the set $B$ in the proposition, every $f\in L^2(B)$ satisfies
\[
 \int_{B\times B}r(x,y)\|D_{x,y}f\|^2\,dJ(x,y)
                    \ge(c-4\theta)\|f\|_{L^2(B)}^2.
\]
In particular a positive gap survives when $4\theta<c$.
\end{corollary}
\begin{proof}
Extend $f$ by zero outside $B$. Every omitted edge contributes at most
$2\|f(x)\|^2+2\|f(y)\|^2$. Symmetry and \eqref{I-eq:pruning-omission} bound the total omitted
energy by $4\theta\|f\|_2^2$. Subtract this from the full gap.
\end{proof}

\subsection{A bounded exact primitive on one retained core}
\label{I-subsec:ha-bounded-core-solver}

We first use a cutoff to put an existing measurable primitive in $L^2$. An operator contraction then gives an $L^\infty$
bound on the same retained set.
Throughout this subsection a reversible probability kernel is
regarded as an edge space with endpoints $x,y$, probability edge
measure, and the measure-preserving reversal $(x,y)\mapsto(y,x)$;
labels and parallel edges are retained. The coefficient differential
has the form
\[
                 Dv(e)=L_e v(y)-N_e v(x),
              \qquad \|L_e\|,\|N_e\|\le1.
\]
There need not be a globally bounded representation on the Hilbert
coefficient space. For the spherical law, the bounded blocks $L_e$
and $N_e=L_eR(g)$ define this differential even on the completed sum
of all nontrivial finite-dimensional representation types.

The cutoff argument below applies Proposition~\ref{I-prop:pruning}
to the edges incident to a large-value set of the section.

\begin{lemma}[Secondary cutoff on a reversible kernel]
\label{I-lem:ha-core-weak-two}
On an atomless probability space, suppose the normalized differential has a gap
$\|Dv\|_2^2\ge c\|v\|_2^2$, and the scalar kernel has norm
at most $\kappa$ on $L^2_0$. Choose $0<\vartheta<1$ with
$2\kappa\le\vartheta<c/4$. Let $v$ be a strongly measurable Hilbert-valued section for which $\|Dv(e)\|\le e_0$
on almost every edge. Strong measurability is automatic for measurable
finite-dimensional sections. Assume that its squared edge norm is unchanged
under reversal. Then
\begin{equation}\label{I-eq:ha-core-weak-two}
 \mu\{\|v\|>t\}\le
                    \frac{e_0^2}{(c-4\vartheta)t^2}
                           \qquad(t>0).
\end{equation}
In particular, $v\in L^p$ for every $1\le p<2$.
\end{lemma}
\begin{proof}
For $M>0$, mark all edges incident to $\{\|v\|>M\}$ bad.
Their directed mass is at most $2\mu\{\|v\|>M\}$, which tends
to zero. Apply pointwise pruning, for all sufficiently large $M$,
to obtain $C_M$. Because $\vartheta<1$, the retained set $C_M$
contains no vertex with $\|v\|>M$, up to a null set.
Thus the actual section restricted to $C_M$ belongs to $L^2$.
Corollary~\ref{I-cor:gap-deletion} bounds the loss from omitted edges
by $4\vartheta\|v\|_{L^2(C_M)}^2$. Hence
\[
 (c-4\vartheta)\int_{C_M}\|v\|^2\,d\mu
      \le\int_{\text{retained edges in }C_M}\|Dv\|^2
      \le e_0^2.
\]
The complement of $C_M$ has measure at most
$4\mu\{\|v\|>M\}/\vartheta$, tending to zero.
For each fixed truncation of $\|v\|^2$, the integral over the
complement tends to zero. Pass first to the limit in $M$, then increase
the truncation. This gives
\[
 \int_X\|v\|^2\,d\mu\le e_0^2/(c-4\vartheta).
\]
Chebyshev gives the stated tail bound. In fact $v\in L^2$, which
implies the asserted $L^p$ membership on this probability space.
\end{proof}

\begin{theorem}[Uniform bounded exact solver on the same core]
\label{I-thm:ha-bounded-exact-core}
Consider the spherical law and normalized differential from
Proposition~\ref{I-prop:ha-bik-common-law} on an atomless base. Write
\[
 c_0=1-\theta_0>0,\qquad
 \|Q\|_{L^p\to L^p}\le\theta_0\quad(1\le p\le\infty),
 \qquad Q=I-D^*D.
\]
Choose the spherical parameter so that the scalar norm on $L^2_0$
is at most $\varepsilon\ge0$. Let $B$ be a measurable retained core
whose pointwise omitted outgoing edge mass is at most $\eta\ge0$.
Assume that there is $\vartheta\in(0,1)$ such that
\begin{equation}\label{I-eq:ha-core-constant-choice}
 2(\varepsilon+2\eta)\le\vartheta,
       \qquad 4\vartheta<c_0-4\eta.
\end{equation}
These inequalities hold when $\varepsilon,\eta$ are sufficiently
small compared with $c_0$; $c_0$ is independent of the large
spherical parameter.

Let $v$ be a section on $B$ whose extension by zero is almost
everywhere strongly measurable. This is automatic for a measurable
finite-dimensional section. Suppose its \emph{actual} retained-edge
differential satisfies $\|D_Bv(e)\|\le e_0$, where $e_0\ge0$. Then
\begin{equation}\label{I-eq:ha-bounded-exact-core}
                   \|v\|_{L^\infty(B)}
                       \le\frac{2e_0}{c_0-4\eta}.
\end{equation}
The constants are independent of the algebraic types, their
multiplicities, and any finite jet cutoff. The theorem assumes an
existing exact measurable primitive; it is not an existence assertion
for arbitrary edge data.
\end{theorem}

\begin{proof}
If $\mu(B)=0$, the asserted essential bound is vacuous. Assume henceforth
that $\mu(B)>0$, and normalize $\mu|B$ to a probability measure. On $B$, retain the good
edges with both endpoints in $B$, and replace every omitted transition
by holding at its starting point. This gives a reversible Markov
kernel $P_B$. The holding edges have zero differential. A mean-zero
function on $B$, extended by zero, has mean zero on $X$. Relative to
the compression of the original scalar kernel, removing bad internal
edges has $L^2$ operator norm at most $\eta$ by the row-and-column
Schur bound; adding the holding probability has norm at most $\eta$.
It follows that
\[
                         \|P_B\|_{L^2_0(B)}
                                  \le\varepsilon+2\eta.
\]
The twisted degree-zero gap on $B$ is at least $c_0-4\eta$.
Thus \eqref{I-eq:ha-core-constant-choice} permits application of
Lemma~\ref{I-lem:ha-core-weak-two} to the kernel $P_B$ with no initially
bad edges. We conclude that the \emph{given} measurable section $v$
belongs to $L^2(B)$.

We next record the stronger normal-operator bound. Extend a section
on $B$ by zero. The omitted-edge normal operator, compressed on both sides to $B$,
is denoted by $E_{\rm miss}$. It has
two diagonal and two off-diagonal terms. Each term has row norm at
most $\eta$ and column norm at most $\eta$: the endpoint matrices
have norm at most one, and reversal identifies the incoming and
outgoing missing mass. Consequently
\[
                  \|E_{\rm miss}\|_{L^p\to L^p}\le4\eta
                              \qquad(1\le p\le\infty).
\]
Since
\[
 Q_B=I_B-D_B^*D_B
       =\mathbf1_B Q\mathbf1_B+E_{\rm miss},
\]
we obtain $\|Q_B\|_{L^p\to L^p}\le\theta_0+4\eta<1$ on the
whole same range of $p$.

Put $b=D_Bv$, the actual bounded edge field. The adjoint endpoint
formula gives $\|D_B^*b\|_\infty\le2e_0$.
Since $v\in L^2$, its normal equation
\[
                         (I-Q_B)v=D_B^*b
\]
is now a legitimate equality in that Banach space. The Neumann series
for $(I-Q_B)^{-1}$ converges both there and in $L^\infty$; uniqueness
in $L^2$ identifies the two solutions with $v$. Its $L^\infty$ norm is
at most $2e_0/(1-\theta_0-4\eta)$, as claimed.

For a completed countable orthogonal sum of algebraic types, apply the argument
first to a finite sum of types. This is one joint coefficient
representation, not a sum of separate $L^p$ estimates. Different Hilbert
multiplicity spaces can be included in their common orthogonal sum;
the original finite sum is then a reducing subspace of the tensor
product with the joint finite-dimensional coefficient representation.
The normalizing operators preserve this subspace. Thus the same constants
apply to every finite-type projection $P_m v$, on the original core $B$.
These projections commute with the endpoint blocks, so
$D_B(P_m v)=P_m D_Bv$ and $\|D_B(P_m v)\|\le e_0$.
Intersect the countably many conull sets on which the resulting bounds
hold. Since $\|P_m v(x)\|\to\|v(x)\|$, the same bound follows for $v$.
This passage uses the bounded endpoint blocks, and requires no
operator-norm measurability of the untruncated representation.
\end{proof}

\section{Analytic realization of a formal reduction}\label{I-sec:convergence}

\subsection{A Newton iteration on nested measurable sets}
\label{I-subsec:ha-cumulative-pruning-newton}

This subsection uses the common spherical law, pointwise pruning,
and bounded exact core solver proved above.  Its corrections remove
the whole degree band, including the invariant coefficients.  In
particular no persistent lower-order centralizer cocycle is left in
the homological equation.

\begin{lemma}[Pruning with a growing rejected-edge set]
\label{I-lem:ha-monotone-pruning}
On an atomless standard probability space, fix a reversible probability
edge measure $J$ whose scalar kernel
has norm at most $\kappa\ge0$ on mean-zero $L^2$.  Fix $\eta>0$ with
$2\kappa\le\eta$.  Suppose measurable reversal-invariant rejected-edge sets
$R_j$ increase and
\[
                 J\left(\bigcup_jR_j\right)\le\beta<\eta^2/4.
\]
The pointwise pruning procedure can be performed cumulatively,
giving decreasing cores $B_j$, such that
\[
 \mu(X\setminus B_j)\le2\beta/\eta,
\]
and at every point of $B_j$ the mass of transitions rejected by
$R_j$ or leaving $B_j$ is at most $\eta$.  Previously removed
vertices need never be restored.  The limiting core has measure
at least $1-2\beta/\eta$.
\end{lemma}
\begin{proof}
For a fixed rejected event $R$, let $o_R(x,D)$ denote the outgoing
mass of edges in $R$ or with terminal vertex in $D$. Define
\[
 D_0(R)=\varnothing,\qquad
 D_{m+1}(R)=D_m(R)\cup\{x:o_R(x,D_m(R))>\eta\},\qquad
 D(R)=\bigcup_mD_m(R).
\]
These sets are measurable. If $x\notin D(R)$, continuity from below
shows $o_R(x,D(R))\le\eta$: otherwise the strict inequality would
already hold at a finite stage. Increasing $R$ increases each stage
and hence increases $D(R)$.

The static pruning construction gives a measurable removed set $E$
with $\mu(E)\le2J(R)/\eta$ and $o_R(x,E)\le\eta$ almost everywhere
outside $E$. Inductively $D_m(R)\subseteq E$ modulo null sets.
Indeed, stationarity implies that almost every row gives zero mass
to any fixed null set, so the induction hypothesis gives
$o_R(x,D_m(R))\le o_R(x,E)$ almost everywhere. No new point outside
$E$ is then removed, except on a null set. Taking the countable union
proves $\mu(D(R))\le\mu(E)$.

Set $B_j=X\setminus D(R_j)$. The cores decrease and satisfy the
asserted pointwise bound and discarded-mass estimate. Continuity of
measure from below on the increasing removed sets gives the same
mass bound for the limiting core. This construction uses only the
current rejected event and therefore also applies to adaptive choices.
\end{proof}

Adding a vertex set $C_j$ to the rejection condition means rejecting
all edges incident to $C_j$.  It costs at most $2\mu(C_j)$ in the
original directed edge measure.  The set $C_j$ is then removed by
the pruning procedure itself; it is not deleted in one uncontrolled
jump.  Since its omitted mass is one and $\eta<1$, no positive
part of it can survive in the resulting core.

\paragraph{Analytic Hilbert norm and elementary analytic estimates.}
Choose a Cartan-compatible Euclidean inner product on $V=\mathbb R^n$,
complexify it, and equip $\operatorname{Sym}^dV$ with its induced
symmetric-tensor norm.  For $f=\sum_{d\ge2}f_d$, let $|f_d|$ be the
Hilbert--Schmidt norm of $f_d:\operatorname{Sym}^dV\to V$ and put
\[
 \|f\|_{\rm a}^{2}=\sum_{d\ge2}(d+1)^2|f_d|^2,\qquad
 \|f\|_{{\rm a},r}^{2}
       =\sum_{d\ge2}(d+1)^2r^{2(d-1)}|f_d|^2.
\]
Thus the second norm is the first norm of $r^{-1}f(rz)$.
On degree $d$ the algebraic coefficient action is
$R_d(g)p=A(g)\circ p\circ A(g)^{-1}$. Let $J_0$ be the orthogonal
projection onto $\Hom_G(\Sym^d V,V)$ in each degree, equivalently
averaging over a compact real form. This is the Reynolds projection.
The degree weights are scalar, so Reynolds projections, type
projections, and all degreewise normal differentials preserve the
orthogonal direct-sum structure used by the core solver.
Evaluation and Cauchy--Schwarz give
$\sup_{B_1}|f|\le\|f\|_{\rm a}$.  If $f=O(|z|^m)$, then
\[
 \sup_{B_\rho}|f|\le\rho^m\|f\|_{\rm a}\quad(0<\rho<1),
 \qquad
 \sup_{B_{1-\delta}}\|Df\|
                         \le C_n\delta^{-1}\|f\|_{\rm a}.
\]
For the converse estimate it suffices to use a polynomial loss in the
coefficient degree.  If $p$ is homogeneous of degree $d\ge1$ on
$\mathbb C^n$, its symmetric-tensor norm satisfies
\[
 |p|\le(d+1)^n\sup_{\|z\|\le1}|p(z)|.
\]
Here is an elementary coefficient argument. For a multi-index
$\alpha$ of total degree $d$, integrate on the polytorus with radii
$r_i=(\alpha_i/d)^{1/2}$ to extract the corresponding monomial
coefficient. Put $N_\alpha=d!/\prod_i\alpha_i!$. The multinomial
probability $N_\alpha\prod_i(\alpha_i/d)^{\alpha_i}$ is a maximal
term of the multinomial distribution with probabilities $\alpha_i/d$.
There are at most $(d+1)^n$ terms, so this probability is at least
$(d+1)^{-n}$, with $0^0=1$. Thus each occupation class contributes
at most $(d+1)^n\sup_{\|z\|\le1}|p(z)|^2$ to $|p|^2$.
Summing over at most $(d+1)^n$ classes gives the assertion. This
argument uses vector-valued Fourier coefficients and the triangle
inequality, so the constant does not depend on the target dimension.

If $f$ is holomorphic on a neighborhood of the closed ball of radius
$R$, with $\sup_{B_R}|f|\le M$, Cauchy's estimate and this bound give
$|f_d|\le(d+1)^n M R^{-d}$. Consequently, for $0<\rho<R$,
\[
 \|\rho^{-1}f(\rho\,\cdot)\|_{\rm a}^2
 \le (M/R)^2
 \frac{2^{2n+2}(2n+2)!}{(1-(\rho/R)^2)^{2n+3}}.
\]
Indeed, with $q=\rho/R$, the squared norm is bounded by
$(M/R)^2\sum_{k\ge0}(k+2)^{2n+2}q^{2k}$; the displayed bound
follows by comparison with the binomial series for
$(1-q^2)^{-2n-3}$. When only holomorphy on $B_R$ is given, apply
the estimate at smaller outer radii and pass to the limit.
In particular, for $1/2\le\rho<R\le1$ and $R-\rho\ge\delta$,
\[
 \|\rho^{-1}f(\rho\,\cdot)\|_{\rm a}
       \le C_n\delta^{-a}\sup_{B_R}|f|,
                 \qquad a=2n+3.
\]
We use this deliberately coarse exponent throughout the iteration.

For completeness, take $\|v\|_{\rm a}\le e$,
$\|p\|_{\rm a}\le t$, let $\ell=\Pi v$, and suppose
$v-\ell=O(|z|^m)$.  On $B_{1-2\delta}$ the difference
\[
 (\id+v)\circ(\id+p)-(\id+p+\ell)
       =(v-\ell)+v(\,\cdot+p)-v
\]
has supremum at most
$e(1-2\delta)^m+C_n et/\delta$.  When $e+t\le c_n\delta$
for a fixed sufficiently small $c_n>0$, the inverse of
$\id+p+\ell$ exists on the necessary smaller ball and has
Lipschitz constant at most two.  Applying the preceding coefficient
estimate at radius $1-3\delta$ proves the short-band bound
\[
 C_n\delta^{-a}
       \left[e(1-2\delta)^{m-1}+et/\delta\right].
\]
Enlarging $C_n$ only changes fixed constants.  This explicitly
substantiates the remainder estimate used below.

\begin{lemma}[An invariant band primitive on the same core]
\label{I-lem:ha-core-invariant-band}
Suppose a retained core $B$ has a reversible holding kernel $P_B$
with $\|P_B|_{L^2_0(B)}\|\le\kappa_B<1$.  Let $v:B\to E$
be a finite measurable field in a finite-dimensional real Hilbert
space and suppose its actual retained differences have norm at
most $e$.  There is a constant vector $c\in E$ such that
\[
       \int_B\|v-c\|^2\,d\mu\le C e^2,
\]
where $C$ depends only on $1-\kappa_B$, not on $\dim E$ or
$\mu(B)\le1$.
\end{lemma}
\begin{proof}
If $\mu(B)=0$, take $c=0$. Otherwise normalize $\mu|B$ to a
probability $\mu_B$. This leaves the actual kernel operator and its
mean-zero norm bound unchanged; the $L^2$ norms on both sides scale by
the same factor.

Fix orthonormal coordinates on $E$. Let $v_R$ be obtained by clipping
each coordinate of $v$ to $[-R,R]$. This map is nonexpansive and
$v_R\to v$ pointwise. For a scalar $L^2$ field, the gap and the
reversible Dirichlet identity give the pairwise Poincar\'e inequality
\[
 \int_{B\times B}|f(x)-f(y)|^2\,d\mu_B(x)d\mu_B(y)
 \le \frac{1}{1-\kappa_B}
       \int |f(y)-f(x)|^2\,dJ_B(x,y),
\]
where $J_B(dx,dy)=\mu_B(dx)P_B(x,dy)$ is the actual holding edge
measure. Apply this to each coordinate of $v_R$ and sum. The clipped
retained differences have norm at most $e$, so
\[
 \int_{B\times B}\|v_R(x)-v_R(y)\|^2\,d\mu_B(x)d\mu_B(y)
       \le \frac{e^2}{1-\kappa_B}.
\]
Fatou's lemma passes this bound to $v$ without assuming its initial
square integrability. By Fubini, some $y_0\in B$ satisfies
\[
 \int_B\|v(x)-v(y_0)\|^2\,d\mu_B(x)
       \le \frac{e^2}{1-\kappa_B}.
\]
Set $c=v(y_0)$ and multiply by $\mu(B)\le1$. This proves the
assertion with $C=(1-\kappa_B)^{-1}$, independent of the number of
coordinates and of $\mu(B)$. Subtracting this constant leaves the
unweighted retained differences unchanged. In the invariant-coefficient
application it defines a constant field on the whole original base.
\end{proof}

\begin{theorem}[Analytic realization from a compatible formal reduction]
\label{I-thm:ha-cumulative-newton-realization}
Let $G$ be connected semisimple with finite center and no compact
factors, and suppose every simple factor has property $(T)$. Let
$G$ act ergodically and preserve a standard probability space.  Let
$A:G\to\mathrm{GL}(n,\mathbb R)$ be a continuous representation. Suppose a
coefficient-measurable analytic germ cocycle $\alpha$ has derivative
$A$ and a coefficient-measurable formal tangent-to-identity transfer
$H$ satisfying
\[
                    \alpha(g,x)H_x=H_{gx}A(g).
\]
Then there is a constant formal tangent-to-identity germ $D$ in
the $G$-centralizer such that $H_xD$ is analytic almost everywhere.
It and its inverse have positive measurable radii.  No moment
assumption or uniform input radius is imposed.
\end{theorem}

\begin{proof}
We first treat a nonatomic base and give the iteration with all
probability losses charged in the original edge measure.  Fix the spherical law once.  Choose its
scalar mixing constant $\kappa$ and an omission tolerance $\eta$
so small that all conclusions of the bounded exact core solver
hold uniformly on every core with pointwise omitted mass at most
$\eta$.  This gives a uniform nontrivial-type $L^\infty$ solver
constant $C_\perp$ and a uniform scalar Poincare constant on each
such core.  The latter and the preceding lemma give a uniform
invariant-band $L^2$ constant $C_0$.  Fix $K\ge1$ large enough
to dominate both constants in the estimates below.

Choose positive numbers $\beta_{\rm in},b_j$ with
\[
 \beta_{\rm in}+2\sum_{j\ge0}b_j=:\beta<\eta^2/4,
                       \qquad 2\beta/\eta<1.
\]
For instance $b_j=b_*/(j+1)^2$ with sufficiently small $b_*>0$.
The arbitrary measurable analytic input has common coefficient
bounds outside a set of directed edges of mass at most
$\beta_{\rm in}/2$.  Intersect this good-edge set with its reversal,
so that the resulting rejected-edge mass is at most $\beta_{\rm in}$.
This assertion uses only countable coefficient-growth tests on
the finite-dimensional analytic germs on the compact support of
the spherical law.  A sufficiently small deterministic initial
dilation makes
\begin{equation}\label{I-eq:newton-initial}
\|A(g)^{-1}\alpha(g,x)-\id\|_{\rm a}\le e_0
\end{equation}
on all these good edges, for any prescribed $e_0>0$.  Here and
below equations for edges hold almost everywhere in the labelled
edge measure.  The scalar mixing law and all analytic operator
constants were fixed before this dilation.

Apply pointwise pruning to the initially rejected edges.  We
construct formal fields $H^{(j)}=\id+O(|z|^{q_j})$, actual analytic
cocycles $\alpha^{(j)}$, nested cores $B_j$, and radii $r_j$, such
that in coordinates rescaled by $r_j$ the normalized displacement
\[
       v_{g,x}^{(j)}=A(g)^{-1}\alpha^{(j)}(g,x)-\id
\]
has order at least $q_j$ and coefficient norm at most $e_j$ on
every retained edge of $B_j$.  Start with $q_0=2$ and the dilated
data in \eqref{I-eq:newton-initial}.  Set
\[
 q_{j+1}=\lceil3q_j/2\rceil,\qquad
 \Pi_jv=\sum_{d=q_j}^{q_{j+1}-1}v_d.
\]
Since $q_{j+1}\le2q_j-1$, this is a compatible short band.

The actual formal band $h_j=\Pi_j(H^{(j)}-\id)$ is a finite-dimensional measurable field. No norm of the full formal series is assumed finite. This band satisfies the
full homological equations.  On each retained edge they read
\begin{equation}\label{I-eq:newton-homological}
A(g)^{-1}h_j(gx)(A(g)z)-h_j(x)(z)=\Pi_jv_{g,x}^{(j)}(z).
\end{equation}
Let $J_0$ denote the coefficientwise Reynolds projection. On
$(I-J_0)h_j$, the bounded normal differential of the common law
has defect of norm at most $e_j$. Indeed, on each coefficient
space, \eqref{I-eq:newton-homological} gives
\[
 D_gh_j=L_gR(g)\Pi_jv^{(j)}_{g,x},\qquad \|L_gR(g)\|\le1.
\]
Thus
\[
             \|(I-J_0)h_j\|_{L^\infty(B_j;{\rm a})}
                                  \le C_\perp e_j.
\]
For $J_0h_j$, equation \eqref{I-eq:newton-homological} is an ordinary Hilbert-valued base
difference with norm at most $e_j$.  The preceding lemma supplies
one constant band $c_j$ in the invariant coefficient space such that, on the same core,
\begin{equation}\label{I-eq:newton-band}
p_j:=h_j-c_j,\qquad
        \int_{B_j}\|p_j\|_{\rm a}^2\,d\mu\le K^2e_j^2.
\end{equation}
Define this polynomial $p_j$ on the entire original base by the
displayed formula.  It is never clipped.  Since $c_j$ commutes
with $A(G)$, equation \eqref{I-eq:newton-homological} still holds globally wherever the
original formal identities hold.

Set
\[
          t_j=K e_j/\sqrt{b_j},\qquad
          C_j=\{x\in B_j:\|p_j(x)\|_{\rm a}>t_j\}.
\]
By \eqref{I-eq:newton-band}, $\mu(C_j)\le b_j$.  Reject all edges incident to $C_j$
and continue the same cumulative pruning procedure.  The new
core $B_{j+1}\subset B_j$ has the same scalar and twisted
coercivity constants, and the total original rejected-edge mass
is at most $\beta_{\rm in}+2\sum_{i\le j}b_i$.  In particular
the bound on its complement is uniform in $j$.  Every point
of $B_{j+1}$ satisfies $\|p_j\|_{\rm a}\le t_j$.

Use the actual polynomial coordinate change $b^{(j)}_x=\id+p_j(x)$.
The exact target cancellation \eqref{I-eq:newton-homological} gives, on retained edges,
\begin{equation}\label{I-eq:newton-cancellation}
A(g)^{-1}(\alpha^{(j)})^{b^{(j)}}(g,x)
 =\big(\id+p_j+\Pi_jv^{(j)}_{g,x}\big)^{-1}
          \circ(\id+v^{(j)}_{g,x})\circ(\id+p_j).
\end{equation}
This formula uses the source polynomial only; it charges no
repeated factor $\|A(g)\|$ to the analytic radius.  It removes
all degrees below $q_{j+1}$.  Right composition of the updated
formal transfer by $(\id+c_j)^{-1}$ restores the same order
normalization globally.  This is a constant centralizer change
and preserves every formal cocycle equation.

Choose $0<\delta_j<1/12$ with $e_j+t_j\le c_n\delta_j$, and
replace $r_j$ by $r_{j+1}=r_j(1-3\delta_j)$.  The analytic
short-band remainder estimate applied to \eqref{I-eq:newton-cancellation} gives
\begin{equation}\label{I-eq:newton-recurrence}
e_{j+1}\le
 2C_n\delta_j^{-a}(1-2\delta_j)^{q_{j+1}-1}e_j
       +2C_n\delta_j^{-a-1}e_jt_j,
                   \qquad a=2n+3.
\end{equation}
Only retained edges of $B_{j+1}$ are used.  The discarded edges
need no bound, and no earlier exceptional-set probability is
multiplied at this step.

Here is a choice making the iteration converge. Take
$b_j=b_*/(j+1)^2$ in the probability budget above, and fix
$\delta_j=\delta_*/(j+1)^2$, with $0<\delta_*<1/12$. Write
\[
 a_j=2C_n\delta_j^{-a}(1-2\delta_j)^{q_{j+1}-1},\qquad
 n_j=2C_nK\delta_j^{-a-1}/\sqrt{b_j},\qquad
 \overline a_j=\max\{a_j,1/8\},
\]
and define $L_0=1$, $L_{j+1}=2\overline a_jL_j$.
Since $q_j$ grows geometrically and
$(1-2\delta_j)^{q_{j+1}-1}\le
\exp(-2\delta_j(q_{j+1}-1))$, we have $a_j\to0$.
Thus, after finitely many terms, $\overline a_j=1/8$ and
$L_{j+1}=L_j/4$. In particular $L_j\le M_0 4^{-j}$ for a
fixed $M_0>0$.

The sequence $n_j$ grows only polynomially. It follows that
$\sup_j n_jL_j<\infty$ and hence
\begin{equation}\label{I-eq:newton-smallness}
 \inf_j\frac{\overline a_j}{n_jL_j}>0.
\end{equation}
Choose $e_0>0$ smaller than this infimum. The quadratic recurrence
\[
 e_{j+1}=a_je_j+n_je_j^2
\]
defines an error envelope for \eqref{I-eq:newton-recurrence}, because
$t_j=Ke_j/\sqrt{b_j}$. Induction gives $e_j\le e_0L_j$:
if the bound holds at stage $j$, then
$n_je_j\le\overline a_j$, and
\[
 e_{j+1}\le2\overline a_je_j
                    \le e_0L_{j+1}.
\]
Moreover
\[
 \sup_j\frac{L_j(1+K/\sqrt{b_j})}{\delta_j}<\infty,
 \qquad \sum_j\frac{KL_j}{\sqrt{b_j}}<\infty.
\]
Decreasing $e_0$ therefore ensures $e_j+t_j\le c_n\delta_j$
for every $j$, and makes $\sum_jt_j$ as small as required below.
All these choices precede the initial dilation.
Finally, $\sum_j3\delta_j\le6\delta_*<1$, so the elementary
finite-product inequality gives
\[
 r_j=r_0\prod_{i<j}(1-3\delta_i)
       \ge r_0(1-6\delta_*)>0.
\]
The radii decrease and converge to a positive limit
$r_\infty=r_0\prod_j(1-3\delta_j)$.

The limiting core $B_\infty=\bigcap_jB_j$ has positive measure.
Fix $x\in B_\infty$. In physical coordinates write
$\psi_j(z)=z+r_jp_j(z/r_j)$. The smallness chosen above gives
\[
 \psi_j(B_{r_{j+1}})\subset B_{r_j}.
\]
Indeed $r_{j+1}=(1-3\delta_j)r_j$, while the displacement is
at most $r_jt_j$ and $t_j\le\delta_j/4$.
Each $\psi_j$ is holomorphic, preserves the real structure, fixes zero,
and has identity derivative. The ordered compositions
$\Psi_j=\psi_0\circ\cdots\circ\psi_{j-1}$ are therefore holomorphic
on $B_{r_j}$ and map this ball into $B_{r_0}$.
Choose $S>0$ with $S\le r_\infty$. All partial products are bounded
by $r_0$ on the common ball $B_S$.
For each finite order their jets stabilize, since the correction orders
tend to infinity. Apply the coefficient estimate above to $\Psi_j-\id$, bounded by
$r_0+S$, on the radii $S/4<S/2<S$. Restore the fixed linear
identity coefficient afterward. This bounds every stabilized
coefficient by one exponential majorant, independent of the partial
product used to compute that coefficient. Thus the stabilized formal
product converges on a positive ball and defines an analytic germ $T_x$.
Its derivative is the identity, so the analytic inverse function theorem
gives an inverse analytic germ. This argument also proves uniqueness
of the analytic germ with the stabilized coefficients.

All finite jets stabilize.  In physical coordinates the correction
and its centralizer normalization are respectively
\[
 b_j^{\rm phys}(z)=z+r_jp_j(z/r_j),\qquad
 C_j^{\rm phys}(z)=z+r_jc_j(z/r_j).
\]
The latter is still a constant $G$-centralizer germ, since scalar
dilation commutes with every linear $A(g)$.  If $\Psi_j$ is the
ordered product of the physical coordinate changes, then
$H^{(j)}=\Psi_j^{-1}H D_j$ in physical coordinates, where
$D_{j+1}=D_j(C_j^{\rm phys})^{-1}$ and $D_0=\id$.
The successive constant centralizer changes start in orders
$q_j\to\infty$, so their product defines
one constant formal centralizer germ $D$, and the limiting formal
identity is $T_x=H_xD$.  Hence $T$ obeys the original formal
conjugacy equations.  Its positive-radius analyticity set is
measurable and invariant modulo null sets: the conjugacy equation
expresses $T_{gx}$ as a composition of $T_x$, $A(g)^{-1}$ and
the actual analytic germ $\alpha(g,x)$.  Ergodicity makes that
set conull.  The countable coefficient-growth criterion selects a measurable
positive radius there, and the analytic inverse theorem followed
by the same coefficient test gives measurable positive inverse
radii.  The formal equalities are consequently equalities of
analytic germs.

If the ergodic base has an atom, it has a conull finite orbit, and
all its points have the same mass $m>0$. We use the same construction
on this finite orbit. Choose the initial rejection budget and the
summable later budgets so that their total edge mass is less than
$m\eta$. Before any vertex removal, a vertex with rejected row mass greater
than $\eta$ would
already contribute more than $m\eta$ to the rejected edge mass.
Thus no removal is eligible, and the canonical retained core is
conull at each stage.
Every finite coefficient field on the finite orbit is in $L^2$;
projecting onto the global invariant sections and applying the same
retained spectral estimate gives the band energy bound directly.
The probability budget can be chosen this small before the initial
dilation, and the error envelope can then be chosen as above.
The resulting nested-ball and coefficient-stabilization argument
proves analytic realization on the finite orbit. In the untwisted
case the centralizer field is invariant; ergodicity makes it constant.
\end{proof}

\subsection{The unitary centralizer block on the same retained core}
\label{I-subsec:ha-unitary-core-centralizer}

A compact centralizing cocycle does not require replacing the preceding
core. It changes the scalar coefficient block into a unitary cocycle
block; its normalization is by global invariant sections rather than
necessarily constant vectors.

\begin{proposition}[Measurable unitary primitives modulo global invariants]
\label{I-prop:ha-unitary-core-centralizer}
Let $G$ be a measurable group acting jointly measurably and ergodically
on the probability space $(X,\mu)$, preserving $\mu$, and let $\Pi$ be the unitary
Koopman representation with a measurable unitary coefficient cocycle.
The coefficient operators and their inverses are assumed jointly
Bochner measurable in operator norm, and the induced $L^2$ representation
has strongly measurable orbit maps. These conditions hold for the
finite-dimensional coefficient bands used below.
Write $\mathcal I$ for its invariant-section space. Suppose a symmetric
probability law $\nu$ satisfies
\[
                 \|\Pi(\nu)|_{\mathcal I^\perp}\|\le\varepsilon<1.
\]
Use the normalized unitary differential
$Dv(g,x)=2^{-1/2}(U_gv(x)-v(x))$, where $U_g$ includes the unitary
coefficient transport. Let $B$ be a measurable retained core with $\mu(B)>0$ and pointwise
omitted mass at most $\eta\ge0$, and assume
$c_u=1-\varepsilon-4\eta>0$.
If a section $v$ on $B$, extended by zero outside $B$, is strongly
measurable and has actual retained differential bounded in normalized
$L^2$ norm by $e_0\ge0$, then $v\in L^2(B)$.
Moreover there is a global invariant section $a\in\mathcal I$ with
\[
 \|v-a|_B\|_{L^2(B,\mu_B)}\le e_0/\sqrt{c_u},
              \qquad \mu_B=\mu(B)^{-1}\mu|_B.
\]
The section $a|_B$ is the orthogonal projection of $v$ onto the
restrictions of $\mathcal I$.
\end{proposition}

\begin{proof}
For $a,b\in\mathcal I$, the scalar function $\langle a(x),b(x)\rangle$
is $G$-invariant and is therefore constant almost everywhere.
Consequently restriction from $\mathcal I$ to $L^2(B,\mu_B)$ is an
isometry with closed image $\mathcal I_B$; each member of this image
has constant pointwise norm on $B$.
A section on $B$ orthogonal to $\mathcal I_B$, extended by zero, is
orthogonal to $\mathcal I$ on $X$. The full gap and the omitted-edge
bound therefore give
\[
 \|D_B f\|_2^2\ge c_u\|f\|_2^2
                         \quad(f\perp\mathcal I_B).
\]
The elements of $\mathcal I_B$ have zero retained differential.
It follows that, for every $f\in L^2(B)$,
\begin{equation}\label{I-eq:ha-unitary-core-poincare}
 \|f-P_{\mathcal I_B}f\|_2
                    \le c_u^{-1/2}\|D_Bf\|_2.
\end{equation}

For $R>0$, let $v_R$ be radial projection of $v$ onto the closed
Hilbert ball of radius $R$. This projection is $1$-Lipschitz and
commutes with every unitary transport, so
$\|D_Bv_R\|_2\le e_0$. Write
$a_R=P_{\mathcal I_B}v_R$. Inequality
\eqref{I-eq:ha-unitary-core-poincare} gives
$\|v_R-a_R\|_2\le e_0/\sqrt{c_u}$.
Choose a measurable $A\subset B$ of positive $\mu_B$-measure on which
$\|v\|\le M$. Such an $A$ exists because $v$ is finite almost
everywhere. For $R\ge M$, $v_R=v$ on $A$. Since $a_R$ has constant
pointwise norm, restriction to $A$ gives
\[
 \|a_R\|_2\le M+\frac{e_0}{\sqrt{c_u\mu_B(A)}}.
\]
Thus the $v_R$ are uniformly bounded in $L^2$. Fatou's lemma yields
$v\in L^2(B)$. Apply \eqref{I-eq:ha-unitary-core-poincare} directly to
$v$ and extend its orthogonal projection through the isometry
$\mathcal I\to\mathcal I_B$.
\end{proof}

For semisimple $G$ whose simple factors have property~\textup{(T)},
the spherical laws used above make $\varepsilon$ uniformly small in
every unitary coefficient representation. Therefore the same law
and the same retained core work for all coefficient degrees.
A global invariant coefficient field is equivariant for the linear
cocycle; its associated polynomial coordinate correction commutes
with that cocycle in the fiberwise sense. Removing such a formal
centralizer component preserves the formal conjugacy equation.

\subsection{Compact centralizing cocycles and descent to a lattice}
\label{I-subsec:ha-compact-core-extension}

\begin{theorem}[Analytic realization with a compact centralizing cocycle]
\label{I-thm:ha-compact-core-realization}
Let $G$ be connected semisimple with finite center, no compact
factors, and property~\textup{(T)} in every simple factor. Let $G$
act ergodically and preserve an atomless standard probability measure on $X$.
Let $A:G\to\mathrm{GL}(V)$ be a finite-dimensional continuous real
representation, let
$C\le Z_{\mathrm{GL}(V)}(A(G))$ be compact, and let
$\kappa:G\times X\to C$ be a measurable cocycle. Set
\[
                         L(g,x)=A(g)\kappa(g,x).
\]
Suppose a coefficient-measurable analytic germ cocycle $\alpha$ has
derivative $L$ and admits a coefficient-measurable formal
tangent-to-identity transfer $H$ with
\[
                       \alpha(g,x)H_x=H_{gx}L(g,x).
\]
Then there is a coefficient-measurable formal tangent-to-identity
field $D_x$ satisfying
\begin{equation}\label{I-eq:ha-compact-formal-gauge}
                    D_{gx}L(g,x)=L(g,x)D_x
\end{equation}
such that $H_xD_x$ is analytic almost everywhere, with positive
measurable radii for it and its inverse. The fields $D_x$ may be
chosen in the formal $A(G)$-centralizer pointwise. No moment
assumption or uniform radius is required.

The same conclusion holds for a lattice $\Gamma<G$ and a
probability-preserving ergodic $\Gamma$-action, when the prescribed
linear cocycle has this form on $\Gamma$ and the compatible formal
transfer is given.
\end{theorem}

When $\kappa$ is trivial, equivariance and membership in the
pointwise $A(G)$-centralizer give $D_{\gamma x}=D_x$. Ergodicity
therefore makes $D_x$ constant in the lattice case as well.

\begin{proof}
Choose an inner product on $V$ which is Cartan-compatible for $A$
and invariant under $C$. Such a choice exists: average a
Cartan-compatible inner product over $C$; because $C$ centralizes
$A(G)$, the Cartan adjoint identities are preserved by averaging.
Use the induced Hilbert norms on every homogeneous coefficient
space $E_d=\operatorname{Hom}(\operatorname{Sym}^dV,V)$.
The algebraic action $R_d(g)p=A(g)pA(g)^{-1}$ is
Cartan-compatible, whereas
$K_d(c)p=cpc^{-1}$ is unitary and commutes with every $R_d(g)$.
The coefficientwise Reynolds projection $J_0$ is orthogonal and
commutes with both actions.

For clarity, define the twisted Koopman antirepresentation by
\[
 (\mathcal U_g f)(x)
       =K_d(\kappa(g,x))^{-1}f(gx).
\]
It is unitary on $L^2$ and isometric on every $L^p$; the cocycle
identity gives $\mathcal U_{gh}=\mathcal U_h\mathcal U_g$.
It commutes with the algebraic coefficient action. Thus every
factorization in Proposition~\ref{I-prop:ha-bik-common-law} remains
valid after replacing $U$ by $\mathcal U$. In particular,
\[
 P=\int_K R_d(k)^*\mathcal U_k\,dk,
 \qquad
 B=P\,\bigl(2\cosh(T\,dR_d(H))\bigr)^{-1}
               \tfrac12(\mathcal U_{a_T}+\mathcal U_{a_T^{-1}})P
\]
give the same uniform normalized gap and all-$L^p$ contraction
on the blocks without invariant vectors. The pointwise row estimate is
unchanged because the extra coefficient transports are unitary
and commute with the algebraic matrices.

There is also no degree loss in the band data. With
$v_{g,x}=L(g,x)^{-1}\alpha(g,x)-\id$, the compatible band primitive
$h$ satisfies
\[
 L(g,x)^{-1}h(gx)L(g,x)-h(x)=\Pi v_{g,x},
 \qquad
 (\mathcal U_g-R_d(g))h=R_d(g)\Pi v_{g,x}.
\]
For $W_g=(I+R_d(g)R_d(g)^*)^{-1/2}$ this yields
\[
             \|W_g(\mathcal U_g-R_d(g))h\|
                         \le\|\Pi v_{g,x}\|.
\]
This is the same bounded normalized input as in the pure algebraic
case, including on a weighted Hilbert sum of the whole band.

On the algebraically invariant block $J_0E_d$, the remaining
coefficient action is the unitary antirepresentation $\mathcal U$.
The same spherical law makes its averaging operator uniformly small
on the orthogonal complement of its global invariant-section space.
This is the uniform spherical estimate applied to that unitary
representation; it does not require scalar or untwisted coefficients.
For a semisimple product, the commuting invariant projections and telescoping estimate in
Proposition~\ref{I-prop:ha-bik-common-law} gives the same conclusion.

Apply Proposition~\ref{I-prop:ha-unitary-core-centralizer} on each
cumulative core. It supplies an invariant band section $c_j(x)$
defined on the \emph{whole} original base, with the same uniform
$L^2$ bound for $h_j-c_j$ used in the full-band iteration. Its global
invariance and its membership in $J_0E_d$ mean exactly
\[
       c_j(gx)\circ L(g,x)=L(g,x)\circ c_j(x).
\]
Accordingly $\id+c_j$ is an equivariant polynomial family, and
right composition by its formal inverse preserves the formal
conjugacy equation. This replaces the constant centralizer
subtraction in Theorem~\ref{I-thm:ha-cumulative-newton-realization}.
All other estimates and the cumulative pruning are unchanged.
The band orders tend to infinity, so the product of these global
equivariant gauges is a coefficient-measurable formal field $D_x$
obeying \eqref{I-eq:ha-compact-formal-gauge}. Each factor is in the
pointwise $A(G)$-centralizer, so the product is as well. The iteration
makes $H_xD_x$ analytic on a positive-measure limiting core; its
analyticity set is invariant, and ergodicity makes it conull.
Measurable positive radii follow from the countable coefficient-growth criterion,
exactly as in the pure algebraic case.

We spell out the lattice descent to avoid evaluation of an
almost-everywhere identity on a potentially null transversal.
Choose a Borel fundamental domain $F\subset G$ for right
$\Gamma$-translation. Such a domain exists and has finite Haar measure by
\cite[Proposition~B.2.4(i)--(ii), pp.~353--354]{BHV}; its measure is positive
because countably many translates cover $G$. Normalize Haar measure so that
$F$ has measure $1$, and write the induced space as $F\times X$ modulo null
sets. The internal construction below uses a countable family of disjoint
local slices to choose this same right fundamental domain; only the
finite-measure assertion from Proposition~B.2.4(ii) is retained as an
external input.
For $g\in G$ and $h\in F$, write
\[
             gh=h'\gamma,\qquad h'\in F,\quad\gamma\in\Gamma.
\]
Then $g(h,x)=(h',\gamma x)$. In the linear fiber charts determined
by $A(h)$, the induced analytic and formal cocycles are
\[
 \widetilde\alpha(g,(h,x))
        =A(h')\alpha(\gamma,x)A(h)^{-1},\qquad
 \widetilde H_{(h,x)}=A(h)H_xA(h)^{-1}.
\]
Their derivative is $A(g)\kappa(\gamma,x)$, because the compact cocycle
centralizes $A(G)$. The preceding ambient-group construction applies.
Let $\widetilde T=\widetilde H\widetilde D$ be the analytic transfer
provided by that construction. We descend the analytic transfer first.
In the fundamental-domain charts set
\[
 U(h,x)=A(h)^{-1}\widetilde T(h,x)A(h).
\]
On arrows with $\gamma=e$, its conjugacy equation gives
$U(h',x)=U(h,x)$. This equality holds for almost every
$(h,h',x)$: use the change of variables $g=h'h^{-1}$ in the joint
almost-everywhere identity and Haar invariance. Fubini therefore gives
a coefficient-measurable analytic germ field $T_x$ with
$U(h,x)=T_x$ almost everywhere. For arbitrary section charts one
uses the corresponding analytic source-chart and linear target-chart
changes; their cocycle identities give the same field $U$.
The diagonal change $(h,x)\mapsto(h\gamma^{-1},\gamma x)$ gives
\[
 \alpha(\gamma,x)T_x=T_{\gamma x}L(\gamma,x)
\]
for almost every $x$. Countability makes these equalities simultaneous
for all $\gamma\in\Gamma$.

Define the formal field $D_x=H_x^{-1}T_x$. Comparing the two conjugacy
equations proves $D_{\gamma x}L(\gamma,x)=L(\gamma,x)D_x$.
To check its pointwise centralizer property, use the formal identity
matching the source and target suspension charts. If these chart changes
are $B$ and $C$, so that $BH=\widetilde H C$, then
\[
 H^{-1}(B^{-1}\widetilde T C)=C^{-1}\widetilde D C.
\]
The linear chart change $C$ is a product of an $A(G)$ element and a
compact factor commuting with $A(G)$. Conjugation by $C$ therefore
preserves the formal $A(G)$-centralizer. Fubini transfers this property
to $D_x$ on the original base. The descended $T_x=H_xD_x$ and its
inverse are analytic; the measurable-radius construction supplies
positive measurable radii for both.
\end{proof}

\subsection{Separate finite-jet reductions are sufficient}
\label{I-subsec:ha-finite-jets-suffice}

\begin{corollary}[Finite-jet-to-germ equivalence]
\label{I-cor:ha-finite-jets-suffice}
Assume the lattice, base, and prescribed derivative
$L(\gamma,x)=A(\gamma)\kappa(\gamma,x)$ hypotheses of
Theorem~\ref{I-thm:ha-compact-core-realization}. For a coefficient-measurable
analytic germ cocycle $\alpha$ with that derivative, the following
are equivalent:
\begin{enumerate}
\item $\alpha$ admits a measurable analytic germ linearization to $L$.
\item For every finite $N$, $j^N\alpha$ admits a measurable $N$-jet
      linearization to some linear cocycle.
\item $\alpha$ admits a single compatible coefficient-measurable
      formal tangent-to-identity linearization to $L$.
\end{enumerate}
No bounds on the finite-order transfers, and no moment assumptions,
are part of this equivalence. The analytic transfer and its inverse
have positive measurable radii.
\end{corollary}
\begin{proof}
Truncation gives (1)$\Rightarrow$(2). For (2)$\Rightarrow$(3),
apply the normalization and compatibility construction in the
proof of Proposition~\ref{I-prop:formal-reduction}, with $a=L$.
That construction uses only the separate finite-order
linearizations: their derivatives normalize the linear target,
and polynomial lifts of the equivariant discrepancies make
successive orders agree. It therefore applies here without
integrability assumptions. The implication (3)$\Rightarrow$(1)
is Theorem~\ref{I-thm:ha-compact-core-realization}.
\end{proof}

Thus a qualitative finite-jet superrigidity theorem at every order
is enough; no degree-uniform finite-jet transfer estimate needs to be
added to its statement. The hypotheses of that finite-dimensional
input still have to be verified. Proposition~\ref{I-prop:formal-reduction} supplies them
under logarithmic moments of the original finite jets.

\section{Proofs of the main results}\label{I-sec:applications}\label{I-series:last-section}
\subsection{The measurable germ theorem}
\begin{proof}[Proof of Theorem~\ref{I-thm:ha-general-integrable-superrigidity}]
Apply Proposition~\ref{I-prop:formal-reduction}. Write
$L(\gamma,x)=A(\gamma)\kappa(\gamma,x)$ and define
\[
 \widetilde\beta(\gamma,x)=Q_{\gamma x}\beta(\gamma,x)Q_x^{-1},
 \qquad \widetilde H_x=Q_xH_xQ_x^{-1}.
\]
These fields are coefficient measurable, and
$\widetilde\beta(\gamma,x)$ is an analytic germ at almost
every $x$. There is no assertion that the transformed
coefficients are integrable. The formal identity is
\[
 \widetilde\beta(\gamma,x)\widetilde H_x
                         =\widetilde H_{\gamma x}L(\gamma,x).
\]
Theorem~\ref{I-thm:ha-compact-core-realization} applies, since
every simple factor of $G$ has property $(T)$. It gives a
measurable analytic tangent-to-identity transfer $T_x$ with
$\widetilde\beta(\gamma,x)T_x=T_{\gamma x}L(\gamma,x)$.
Set $F_x=Q_x^{-1}T_x$. This gives
\eqref{I-eq:general-normal-form}. The coefficient test in
Section~\ref{I-sec:formal} selects positive measurable radii
for $F_x$ and $F_x^{-1}$. Countability of the lattice and
the coefficient set gives a single invariant conull set.
\end{proof}

\subsection{Compact analytic actions}
\begin{proof}[Proof of Theorem~\ref{I-thm:compact-actions}]
Choose finitely many analytic charts $\psi_i:U_i\to\R^n$
and open sets $V_i$ covering $M$ with
$\overline V_i\subset U_i$. Partition $M$ into Borel sets
$X_i\subset V_i$ and write $i(x)=i$ on $X_i$. Define
\[
 c_x(z)=\psi_{i(x)}^{-1}(\psi_{i(x)}(x)+z).
\]
Each $c_x$ is an analytic coordinate germ at $x$, with
measurable coefficients. The centered germ cocycle is
\[
 \beta(\gamma,x)=c_{\alpha(\gamma)x}^{-1}
                    \circ\alpha(\gamma)\circ c_x.
\]
Its cocycle identity follows by cancellation of the intermediate
coordinate chart.

Fix a generator $s$ and an order $N$. For each pair of charts
$i,j$, the set
\[
 \{x\in\overline V_i:\alpha(s)x\in\overline V_j\}
\]
is compact and lies in the open domain on which
$\psi_j\alpha(s)\psi_i^{-1}$ is analytic. Its derivatives
through order $N$ are bounded there. There are finitely many
chart pairs, so all coefficients of $j^N\beta(s,x)$ are
bounded uniformly in $x$. The same argument for $s^{-1}$,
and the identity
\[
 \beta(s,x)^{-1}=\beta(s^{-1},\alpha(s)x),
\]
bound the inverse jet. Substitution matrices involve only
finitely many of these coefficients. Thus
\eqref{I-eq:ha-general-finite-jet-moments} holds for every
invariant probability $\mu$, including a singular one.

Apply Theorem~\ref{I-thm:ha-general-integrable-superrigidity}
to $\beta$ and put $\Phi_x=c_xF_x$. The identity
\eqref{I-eq:compact-normal-form} follows. Analytic composition
and inversion preserve positive radius after restriction;
the coefficient test supplies measurable choices of the
restricted radii. No regularity transverse to $x$ is used.
\end{proof}

\subsection{Centered fiber cocycles}
\begin{proof}[Proof of Theorem~\ref{I-thm:centered-fibers}]
Choose a finite family of bundle trivializations over open
sets $U_i$ and a cover by open sets $V_i$ with
$\overline V_i\subset U_i$. As above choose a Borel partition
$X_i\subset V_i$. Let $e_x:\R^n\to E_x$ be the resulting
measurable linear frame. Define
\[
                  \beta(\gamma,x)=e_{\gamma x}^{-1}f_{\gamma,x}e_x.
\]
This is a coefficient-measurable analytic germ cocycle.
For every fixed generator and finite order, continuity of
the finite fiber jets bounds their coefficients on the
compact chart-pair sets. Applying the same assertion to
the inverse generator bounds the inverse jets. Hence
\eqref{I-eq:ha-general-finite-jet-moments} holds.

Theorem~\ref{I-thm:ha-general-integrable-superrigidity} supplies
analytic germs $T_x$ in these measurable frames. Set
$F_x=e_xT_x$. This gives \eqref{I-eq:fiber-normal-form} and
the claimed radii. The proof uses the fiber rank $n$ only;
the topology or dimension of the base does not enter the
coefficient convergence estimates.
\end{proof}

The fiber theorem assumes continuity of the finite jets in the original
bundle charts. The measurable germ theorem instead assumes the original
logarithmic jet moments. The convergence theorem starts with a compatible
formal reduction and imposes no moment condition. None of these statements
asserts formal existence for an arbitrary measurable cocycle without moments.

\subsection{The standard and dual representations}
\label{I-sec:standard-dual}
For the classification series we take $G=\SL(n,\R)$, $n\ge3$,
and fiber dimension $n$. Let $R:G\to\GL(n,\R)$ be continuous.
Hall's matrix-group homomorphism theorem gives its real Lie algebra map
\[
 dR:\mathfrak{sl}_n(\R)\longrightarrow\mathfrak{gl}_n(\R),
 \qquad R(\exp(tX))=\exp(t\,dR(X))
\]
\cite[Theorem~2.21, pp.~45--47]{Hall2003}.
If $dR=0$, connectedness and the exponential identity make $R$ trivial.
Otherwise simplicity of $\mathfrak{sl}_n(\R)$ makes $dR$ injective
\cite[II.1, Example~1, pp.~124--126; proof of Proposition~7.1, p.~434]{KnappBeyondDigital}.
Its image is traceless, since the trace of a commutator is zero and
$\mathfrak{sl}_n(\R)$ is perfect. Thus $dR$ is an injective endomorphism of the finite-dimensional space
$\mathfrak{sl}_n(\R)$, and is therefore an automorphism. Its complex extension is therefore an automorphism
of $\mathfrak{sl}_n(\C)$.

Jacobson's direct matrix theorem says that this complex automorphism is
conjugation or conjugation followed by negative transpose
\cite[Chapter~IX, \S5, Theorem~5; Russian translation, p.~306,
proof pp.~305--306]{Jacobson1964}; see also
\cite[Theorem~2.3, p.~183]{LombardoMikhailov}.
A complex intertwiner $B$ descends to a real one: the polynomial
$\det(\operatorname{Re}B+t\operatorname{Im}B)$ is not identically zero,
since its value at $t=i$ is $\det B\ne0$. A real value of $t$ therefore
gives an invertible real intertwiner. The standard and dual matrix
representations have the corresponding differentials. Connectedness and
the same exponential identity promote this differential intertwiner to
the original group representation. Here exponential neighborhoods and
connectedness are supplied by \cite[Theorem~0.15 and Proposition~1.145]{KnappBeyondDigital}.
Thus $R$ is trivial, standard or dual, up to a real change of coordinates.
The centralizer of either nontrivial representation consists of
real scalars, so its compact subgroups lie in $\{I,-I\}$.

For either nontrivial representation,
\begin{equation}\label{I-eq:standard-no-invariants}
 \Hom_G(\Sym^d\R^n,\R^n)=0\qquad(d\ge2).
\end{equation}
Indeed, the stabilizer of a nonzero vector fixes only its span.
An equivariant homogeneous polynomial map must therefore send
$v$ to a scalar multiple of $v$. Transitivity on
$\R^n\setminus\{0\}$ makes that scalar constant; homogeneity
then makes it zero. The dual case follows by inverse transpose.

\begin{lemma}[Recurrence for unipotents]\label{I-formal:lem:unipotent}
Let $T$ preserve a probability measure, let $N$ be nilpotent on
a finite-dimensional real vector space, and let $v$ be measurable
and finite almost everywhere. If $v(Tx)=\exp(N)v(x)$, then
$Nv(x)=0$ almost everywhere. The conclusion also holds with
a measurable factor in a compact matrix group acting continuously on
the vector space and commuting with $N$.
\end{lemma}
\begin{proof}
We first record the elementary recurrence argument in a form that also
applies to intermediate fields. If $u(Tx)=u(x)+w(x)$ and
$w(Tx)=w(x)$, then $u(T^j x)=u(x)+jw(x)$. Almost every orbit
returns infinitely often to some set $\{\|u\|\le k\}$, so $w=0$
almost everywhere.

Put $U=\exp(N)$ and $A=U-I$. The operator $A$ is nilpotent.
Suppose $A^k v=0$ with $k\ge2$. Then $w=A^{k-1}v$ is invariant,
whereas $u=A^{k-2}v$ satisfies $u(Tx)=u(x)+w(x)$. The preceding
argument gives $A^{k-1}v=0$. Induction gives $Av=0$.
Finally $\exp(N)-I=N B(N)$, where $B$ is a polynomial with
constant term $1$; $B(N)$ is invertible and commutes with $N$.
Thus $Nv=0$. All orbit identities can be imposed simultaneously
on a conull set, since only countably many iterates are used.

For a compact factor $c(x)\in C$ acting by $R$, use the probability
extension $(x,k)\mapsto(Tx,c(x)k)$ with normalized Haar measure on $C$.
The field $w(x,k)=R(k^{-1})v(x)$ then satisfies
$w(Tx,c(x)k)=\exp(N)w(x,k)$, since $R(C)$ commutes with $N$.
The first case gives $Nw=0$ almost everywhere. Fubini and invertibility
of $R(k)$ imply $Nv=0$ almost everywhere on the original base.
\end{proof}

\begin{proposition}[Uniqueness with a prescribed derivative]
\label{I-formal:prop:unique}
For a probability-preserving action of a lattice in $\SL(n,\R)$, $n\ge3$,
let the derivative cocycle be measurably cohomologous to
$\sigma(\gamma)c(\gamma,x)$, where $c$ is a measurable cocycle valued in $\{I,-I\}$.
A tangent-to-identity formal linearization to the given derivative,
if it exists, is unique modulo null sets.
\end{proposition}
\begin{proof}
Normalize the derivative and take the quotient of two transfers.
Its first nonzero homogeneous term, of degree $d\ge2$, would be
a measurable equivariant field in
$E_d=\Hom(\Sym^d\R^n,\R^n)$. Remove the sign on its two-point
probability extension. Induce to the finite-volume probability
space $(G\times X)/\Gamma$. If $v(x)$ is the coefficient and
$R_d$ is the coefficient representation, the field
$V([g,x])=R_d(g)v(x)$ is well-defined and $G$-equivariant.
Apply Lemma~\ref{I-formal:lem:unipotent} to the time-one maps
of finitely many elementary unipotent subgroups generating $G$.
The field is fixed by $G$ almost everywhere, contradicting
\eqref{I-eq:standard-no-invariants}. This argument does not
require unipotent elements of the source lattice.
\end{proof}

\begin{corollary}[Analytic reduction in a prescribed frame]
\label{I-power:thm:main}
Let $\Gamma<\SL(n,\R)$, $n\ge3$, be any lattice, and let
$\beta$ be a coefficient-measurable analytic germ cocycle over
an ergodic probability-preserving $\Gamma$-action. Assume that
the original finite jets are integrable as in
\eqref{I-eq:ha-general-finite-jet-moments}. Suppose a measurable
linear frame $B_x$ satisfies
\[
 B_{\gamma x}^{-1}D_0\beta(\gamma,x)B_x
       =\sigma(\gamma)c(\gamma,x),\qquad c(\gamma,x)\in\{I,-I\}.
\]
Here $c$ is a measurable sign cocycle.
There is a unique measurable analytic transfer $b_x$ with
$D_0b_x=B_x$ and
\[
 \beta(\gamma,x)b_x=b_{\gamma x}\sigma(\gamma)c(\gamma,x).
\]
It and its inverse have positive measurable radii.
\end{corollary}
\begin{proof}
Theorem~\ref{I-thm:ha-general-integrable-superrigidity} gives an analytic
transfer $P_x$ to a linear cocycle $L(\gamma,x)$. Put $J_x=D_0P_x$.
Differentiating $\beta(\gamma,x)P_x=P_{\gamma x}L(\gamma,x)$ gives
$D_0\beta(\gamma,x)J_x=J_{\gamma x}L(\gamma,x)$. Consequently
$H_x=P_xJ_x^{-1}$ is tangent to the identity and satisfies
\[
 \beta(\gamma,x)H_x=H_{\gamma x}D_0\beta(\gamma,x).
\]
Set $b_x=H_xB_x$. The prescribed-frame equation gives the required
conjugacy, and $D_0b_x=B_x$. Composition with measurable linear germs
preserves coefficient measurability; the radius construction applies
to $b_x$ and its inverse.
Uniqueness follows from Proposition~\ref{I-formal:prop:unique}.
\end{proof}

\begin{corollary}[Coordinate germs in the noncompact derivative branches]
\label{I-power:cor:manifold}
For an analytic action of $\Gamma<\SL(n,\R)$, $n\ge3$, on a compact
$n$-manifold preserving an ergodic probability, the standard and
dual derivative branches admit measurable analytic coordinate
germs with exactly that linear action and its central sign cocycle.
Here a derivative branch means a given measurable derivative frame
and a measurable sign cocycle as in Corollary~\ref{I-power:thm:main}.
The germs and their inverses have positive measurable radii.
\end{corollary}
\begin{proof}
Use the bounded centered charts constructed in the proof of
Theorem~\ref{I-thm:compact-actions}, and apply
Corollary~\ref{I-power:thm:main} in the prescribed derivative frame.
Retain the resulting germ reduction when passing to physical coordinates.
After rectifying the sign cocycle on a null set, form the probability
extension $\gamma\cdot(x,\varepsilon)=(\gamma x,c(\gamma,x)\varepsilon)$.
If $\Phi_x$ is the chosen physical coordinate germ, use
$\Phi_{x,\varepsilon}(v)=\Phi_x(\varepsilon v)$ on this extension.
The same reduction identity gives the linear action $\sigma(\gamma)$
alone. The sign cocycle need not be a character of $\Gamma$.
\end{proof}

\clearpage
\section*{Part II. Local description for \texorpdfstring{$n\ge 4$}{n >= 4}}
\hypertarget{part.ii}{}
\addcontentsline{toc}{part}{Part II. Local description for \texorpdfstring{$n\ge 4$}{n >= 4}}

\section{Introduction}
This paper contributes to the Zimmer program \cite{Zimmer} by
studying the local geometry of real-analytic lattice actions in
the defining dimension. The classification problem in this
series is to describe actions of lattices in $\SL(n,\R)$ on
compact $n$-manifolds without assuming an invariant volume.
Here we seek the invariant submanifolds and neighborhood models
from which the global classification can be assembled.

The standard action of $\SL(n,\R)$ on projective space has
dimension $n-1$. In an action on a $n$-manifold, a projective orbit
therefore leaves one transverse coordinate. This coordinate is
not determined by the projective action: its transformation law
contains the flow of a one-dimensional analytic vector field.
The purpose of this paper is to find these hypersurfaces and
to describe their full neighborhoods.

The ambient-group models have a long history. Uchida
\cite{Uchida} classified analytic $\SL(n,\R)$ actions on
$S^n$. Cairns--Ghys \cite{CG} studied local linearization
and the possible low-dimensional orbits. Fisher--Melnick
\cite{FisherMelnick} classified smooth actions of the ambient
group on closed $n$-manifolds and described their analytic
transverse data. Son \cite{Son} extended the real-analytic
ambient-group classification to closed $m$-manifolds with
$3\le n\le m\le2n-3$. Her induced flow constructions for
$m=n$ agree with those of Fisher--Melnick
\cite[Section 4]{Son}. Our task is to recover these local models
from a lattice action, where an ambient action on $M$ is not
given in advance.

Let $G=\SL(n,\R)$, let $\Gamma<G$ be a lattice, and let
$\alpha:\Gamma\to\operatorname{Diff}^{\omega}(M)$, where $M$ is a
compact connected $n$-manifold without boundary and $n\ge3$.
We use the suspension
\[
 M^\alpha=(G\times M)/\Gamma,\qquad
 (g,x)\gamma=(g\gamma,\alpha(\gamma^{-1})x),
\]
with its induced left $G$-action and projection to $G/\Gamma$.
A fiberwise derivative or exponent refers to the tangent bundle
of the fibers of this projection. Fix a split Cartan subgroup
$A$ and a minimal parabolic subgroup $P$ containing it.

A probability on the suspension is \emph{Haar-projecting} if
its image on $G/\Gamma$ is normalized Haar measure. Its
\emph{stabilizer} is the subgroup of $G$ preserving that measure;
this is distinct from the stabilizer of an individual point.
A $P$-invariant probability is therefore a measure whose
stabilizer contains $P$.

\subsection{The local models}
We use \emph{collar} for a neighborhood of a hypersurface $N$
with a transverse coordinate $t$, with $N$ given by $t=0$.
For a two-sided hypersurface this is a product
$N\times(-\varepsilon,\varepsilon)$; at a boundary it is
$N\times[0,\varepsilon)$. For a hypersurface with nontrivial
normal bundle, we use product coordinates on its
normal-orientation cover and retain the deck involution.
The collar describes the action near $N$. The transverse
coordinate need not be preserved by the action, and a fixed
collar need not be invariant under every group element.

A collar is local data, whereas a $G$-disk or $G$-tube is an
entire compact piece. On the spherical cover a $G$-tube is
topologically $S^{n-1}\times[-1,1]$ and has a collar at each
boundary component. A $G$-disk has one boundary collar and
a central fixed point with a linear action near it. These are
the models of Fisher--Melnick
\cite[Definitions 3.3 and 3.4]{FisherMelnick}, recalled in Part IV.
The local description below determines collars; the continuation
arguments in Part IV determine the pieces incident to them.

On the angular cover $S^{n-1}$, a collar neighborhood has
the induced local $G$-action
\begin{equation}\label{II-eq:intro-collar}
 \gamma(\theta,t)=
 \left(\frac{\sigma(\gamma)\theta}{\|\sigma(\gamma)\theta\|},
 \Phi^B_{\log\|\sigma(\gamma)\theta\|}(t)\right),
 \qquad B(0)=0.
\end{equation}
Here $\sigma$ is standard or dual, and $\Phi^B$ denotes the
local flow of $B(t)\partial_t$. Explicitly, the standard representation
is $\sigma(g)=g$, and the dual is $\sigma(g)=g^{-T}$.
The \emph{angular cover} is the map
$S^{n-1}\to\RP^{n-1}$ sending a unit vector to its line;
its deck transformation is $\theta\mapsto-\theta$.
This cover concerns the hypersurface, whereas the
normal-orientation cover concerns its transverse line bundle.
The term \emph{full-angular} means that the model is defined
along the whole compact angular hypersurface, not just in one
projective coordinate chart.
The identity is one of germs
along the entire hypersurface; its common domain may depend on
$\gamma$. The angular quotient may be antipodal. A one-sided
hypersurface is described on its normal-orientation cover, with
the deck involution retained as part of the data.
This is the character-induced flow construction of
Fisher--Melnick \cite[Sections 3.1 and 5]{FisherMelnick},
restricted to $\Gamma$.
Keeping the whole vector field is essential even when its
linear normal weight vanishes.

We use their terminology for these models. For the standard
action let $Q$ stabilize $\R e_1$. With
$\chi(q)=\log\|qe_1\|$, the character $\nu_0$ of
\cite[Section 3]{FisherMelnick} is $-\chi$ on $Q^0$;
their transverse vector field $X$ is therefore
$-B(t)\partial_t$ in \eqref{II-eq:intro-collar}. When
$B\not\equiv0$, this is a neighborhood of the exceptional
divisor in their model $BS_X$ over $S^{n-1}$, with
the appropriate finite quotient if necessary. The usual
blow-up $B_X$ identifies $(\theta,t)$ with $(-\theta,-t)$
and requires the corresponding symmetry of the vector field.
The quotient by $(\theta,t)\mapsto(-\theta,t)$ is a different
allowed finite quotient. When $B=0$, we have a product
neighborhood from their Construction I.

The linear normal weight is $q=B'(0)$. The cases $q<0$ and
$q>0$ give different stable-manifold constructions. If $q=0$,
the transverse flow can still be nontrivial. Establishing its
joint analytic dependence on the angular coordinate is a separate
step of the proof.
We call the case $q=0$ \emph{neutral}. Neutrality concerns the
linear term; it does not imply $B\equiv0$ or a trivial transverse
action.

\begin{definition}[Local description]\label{II-def:recognition}
The action has the local description property if the following hold
for every finite-index subgroup and every measure specified below.
Finite-index restrictions used to describe individual models are
allowed to depend on that model.
\begin{enumerate}
\item[(G)] Every ergodic invariant probability has compact, standard, or dual
derivative type. Compact type forces finite image on connected $M$.
The other types admit measurable actual analytic coordinate germs
exactly intertwining the linear representation, on a finite probability
extension when a central sign must be removed.
\item[(P)] Every ergodic Haar-projecting minimal-parabolic invariant
probability on $M^\alpha=(G\times M)/\Gamma$, with proper stabilizer,
is carried by finitely many compact embedded analytic projective
covering hypersurfaces with full collars, unless the entire action
has one of the closed models specified below. This assertion concerns
the given measure, including measures supported over a compact
invariant subset of $M$.
\item[(C)] The collars have the full-angular model
\eqref{II-eq:intro-collar}. The representation and the angular space, $S^{n-1}$ or
$\RP^{n-1}$, are part of the data.
\end{enumerate}
The closed models in (P), up to finite covering and finite index,
are the Hopf manifolds
$(\R^n\setminus\{0\})/\langle e^T I\rangle$, written as
\[
 S^{n-1}\times\R/T\Z,\qquad
 \gamma(\theta,t)=
 \left(\frac{\sigma(\gamma)\theta}{\|\sigma(\gamma)\theta\|},
 t+\log\|\sigma(\gamma)\theta\|\right),\quad T>0,
\]
and, only for $n=4$, finite covers of $\operatorname{Gr}_2(\R^4)$
with its standard lattice action. For $n=3$ the additional
closed model is the full flag manifold
$\operatorname{Flag}(\R^3)=\{(\ell,\Pi):\ell\subset\Pi\subset\R^3,
\dim\ell=1,\ \dim\Pi=2\}$ and its indicated finite covers.
\end{definition}

\begin{theorem}[Local description]\label{II-main:thm:recognition}
Every action as above with $n\ge4$ has the local description of
Definition~\ref{II-def:recognition}. The conclusions hold after
arbitrary finite-index restriction and for each specified
suspension measure. In particular they apply to a measure
supported on $(G\times K)/\Gamma$ for any compact invariant
subset $K\subset M$.

For each projective collar put $q=B'(0)$. If $q\ne0$, choose the
transverse coordinate so that $B(t)=qt$. If $q=0$, the analytic
field satisfies $B(0)=B'(0)=0$, including the possibility $B=0$.
The nonzero radical-shear alternative in the neutral suspension
gives a Hopf manifold.
The additional closed homogeneous model in this dimension range
is the two-plane Grassmannian when $n=4$.
\end{theorem}

Theorem~\ref{II-I:thm:invariant-germs} gives (G). The proof of
(P) begins with the suspension measure results of
Brown--Rodriguez Hertz--Wang \cite{BRHW,BRHWFactors}.
Their root-graph and coherent-graph arguments provide the
measurable geometric structure. We verify the analytic regularity
and compare the resulting charts to obtain (C). The proof is concluded in Section~\ref{II-sec:recognition}.

\subsection{Organization and shared results}
The geometric preliminaries, the spectral-gap realization tools,
and the invariant-probability theorem below retain their stated
scope for $n=3$. They are also used in Part III. The
local-description theorem of the present paper assumes $n\ge4$;
its proof uses the higher-rank Levi $\SL(n-1,\R)$.
The corresponding local theorem for $n=3$ is the subject of Part III.
Part IV proves global continuation and assembly.

We first collect the common tools and treat invariant probabilities.
We then construct projective hypersurfaces and their neighborhoods
in the contracting, expanding, and neutral cases. The final section
includes the Grassmannian alternative and proves the local theorem,
with its assertions about finite-index subgroups and supports.

\section{Finite-volume tools for arbitrary lattices}\label{II-series:lattice}\label{II-sec:lattice-tools}

In this section $G=\SL(n,\R)$, $\Gamma<G$ is an arbitrary
lattice, and $n\ge3$. All assertions remain true after passing
to a finite-index subgroup. We use the standard representation;
composition with inverse transpose gives the dual versions.

\begin{lemma}[Ergodicity and minimality of the linear models]
\label{II-general:lem:linear-models}
The $\Gamma$ action on $V\setminus0$, with its Lebesgue measure
class, is ergodic. Its action on the sphere of oriented rays,
on projective space, and on each real Grassmannian is minimal.
There is no invariant probability on $V\setminus0$ or on its
antipodal quotient.
\end{lemma}
\begin{proof}
For the matrix-coefficient input we use the connected noncompact
simple real group with finite center case of Howe--Moore; a precise
statement and a new proof appear in
\cite[Theorem 1.1 and Definition 2.6, pp.~66 and 69]{Ciobotaru}.

For ergodicity let $H=\operatorname{Stab}_G(e_1)$.
A $\Gamma$-invariant measurable function on $G/H$ lifts to an
$H$-invariant function on $\Gamma\backslash G$. The subgroup
$H$ is noncompact. Howe--Moore in the form of
\cite[Theorem~1.1 and Definition~2.6, pp.~66, 69]{Ciobotaru} implies ergodicity
of its action on this finite-measure quotient, so the function
is constant. The homogeneous measure class on $G/H$ is the
Lebesgue class on $V\setminus0$.

Here are details giving topological minimality, rather than
only measurable ergodicity. Fix $1\le r<n$, an $r$-plane $L$,
and a nonempty open target set $O$ in the Grassmannian. Put
$E^+=\langle e_1,\ldots,e_r\rangle$ and
\[
 a_t=\operatorname{diag}(e^{(n-r)t}I_r,e^{-rt}I_{n-r}).
\]
Choose a relatively compact open set $B\subset G$ whose
closure satisfies $k^{-1}E^+\in O$ for every $k\in\overline B$.
Howe--Moore implies ergodicity of $a_1$ on $G/\Gamma$.
Consequently almost every $g\Gamma$ visits the positive-measure
set $B\Gamma/\Gamma$ at infinitely many positive integer times.
Choose such a $g$ with $gL$ transverse to
$\langle e_{r+1},\ldots,e_n\rangle$; failure of this condition
has Haar measure zero. Write $a_{d_j}g=k_j\gamma_j$, with
$k_j\in B$ and $\gamma_j\in\Gamma$. Since
$a_{d_j}gL\to E^+$, uniformly bounded $k_j^{\pm1}$ show that
$\gamma_jL\in O$ eventually. This proves minimality.
For oriented rays use $r=1$ and additionally choose $g$ with
positive first coordinate of $gv$; this is a nonempty open
condition and selects the required sign.

Finally, if $\nu$ were an invariant probability on $V\setminus0$,
induce it to the probability space $(G\times(V\setminus0))/\Gamma$.
The well-defined map $[g,v]\mapsto gv$ is $G$-equivariant.
Recurrence for each of finitely many elementary real unipotent
one-parameter subgroups forces its value to be fixed by all of
$G$: a nonconstant vector-valued polynomial orbit is unbounded
and cannot return to a bounded level set. The only fixed vector
is zero, a contradiction. The norm and the same polynomial
argument descend to the antipodal quotient.
\end{proof}

\begin{lemma}[Angular comparisons]\label{II-general:lem:angular}
A continuous $\Gamma$-equivariant self-map of standard projective
space is the identity. There is no continuous equivariant map
from standard projective space to its dual projective space.
On the sphere of oriented rays a continuous equivariant self-map
is the identity or the antipodal map. The analogous statements
hold after interchanging standard and dual.
\end{lemma}
\begin{proof}
Let $Q$ stabilize $[e_1]$, and let $h$ be such a map with target
representation $\rho$. On $\Gamma\backslash G$ define
\[
 f(\Gamma g)=\rho(g)^{-1}h([ge_1]).
\]
It is well-defined and satisfies $f(yk)=\rho(k)^{-1}f(y)$ for
$k\in Q$. Its pushforward probability is invariant under the
unipotent radical of $Q$ and the elementary unipotents of its
$\SL(n-1,\R)$ Levi. For a unipotent linear map, every recurrent
point of its projectivization is fixed: the highest nonzero
power in its nilpotent logarithm gives the limit of a nonfixed
projective orbit, different from its initial point. Poincare
recurrence therefore forces $f$ to be fixed by these subgroups.
Their only common fixed projective point in the standard module
is $[e_1]$; in the dual module they have no common fixed
projective point. Indeed the Levi permits only the first dual
coordinate line and the row radical does not preserve that line.
Thus $h$ is the identity almost everywhere in the standard case,
and the dual case is impossible. Continuity gives the assertion
everywhere. For the sphere use the positive-ray stabilizer in
place of $Q$. The same proof gives $h(\theta)=\pm\theta$;
continuity and connectedness make the sign constant.
\end{proof}

\begin{lemma}[Pinched attracting planes for every lattice]
\label{II-general:lem:pinched}
For $1\le r<n$, every lattice $\Gamma<G$ contains an element
with an attracting fixed $r$-plane at which the Grassmann action
is analytically linear in a graph chart, with linear map $A$
satisfying
\[
 \|A\|<1,\qquad \|A^{-1}\|\|A\|^2<1.
\]
The same inequalities hold for every positive power of $A$.
\end{lemma}
\begin{proof}
Use $a_t$ from Lemma~\ref{II-general:lem:linear-models} and fix a
sufficiently small symmetric relatively compact identity
neighborhood $U$ in $G$. Mixing on $G/\Gamma$ implies
$\Gamma\cap U a_t U\ne\varnothing$ for every sufficiently
large $t$: apply mixing to the positive-measure set $U\Gamma$.
Choose $\gamma=u a_t v$ with $u,v\in U$. Conjugation by $u$
identifies it with $a_t w$, where $w=vu$ is uniformly close
to the identity. Put $a=e^{(n-r)t}$, $b=e^{-rt}$ and write
\[
 a_tw=\begin{pmatrix}aC&aD\\bE&bF\end{pmatrix}.
\]
The blocks $C,F$ and their inverses are uniformly bounded and
close to identities; $D,E$ are small. The invariant expanding
plane is the graph of $X$ solving
\[
 X=(b/a)(E+FX)(C+DX)^{-1}.
\]
For large $t$ this is a contraction on a fixed small ball, so
$X=O(b/a)$. The invariant contracting plane is the graph of
$Y$ solving
\[
 Y=-C^{-1}D+(b/a)C^{-1}Y(EY+F).
\]
The same contraction argument gives uniformly small $Y$.
In the corresponding splitting, $a_tw$ has diagonal blocks
$aC_t$ and $bF_t$, where $C_t,F_t$ and their inverses have
norm bounded by a constant independent of $t$. Its action on
graphs from the first plane to the second is exactly
\[
 Z\longmapsto (b/a)F_t Z C_t^{-1}.
\]
Therefore the norm of this map is at most $C_0e^{-dt}$ and
that of its inverse at most $C_0e^{dt}$. The desired pinching
follows for large $t$. Submultiplicativity proves it for powers.
\end{proof}

\begin{lemma}[Analyticity of a contracting parametrization]
\label{II-general:lem:analytic-parametrization}
Let $F:(\R^d,0)\to(\R^d,0)$ be an analytic diffeomorphism
germ and let $A\in\GL(k,\R)$ satisfy $a=\|A\|<1$.
Suppose a $C^r$ germ $h:(\R^k,0)\to(\R^d,0)$ satisfies
\[
 F\circ h=h\circ A,\qquad
 r\ge2,\qquad \|DF(0)^{-1}\|a^r<1.
\]
Then $h$ is analytic. In particular this conclusion holds for
every $C^\infty$ such germ, regardless of the eigenvalues of
$DF(0)$ in directions transverse to the image of $h$.
\end{lemma}
\begin{proof}
Choose $L>\|DF(0)^{-1}\|$ with $La^r<1$. Complexify the
analytic inverse $F^{-1}$ and restrict it to a convex ball on
which its derivative has norm at most $L$. Put $P=j^{r-1}_0h$.
The $C^r$ equivariance identity implies that the analytic map
\[
 R(z)=F^{-1}(P(Az))-P(z)
\]
vanishes to order $r$. On a sufficiently small complex ball
$B_\delta\subset\C^k$, use the Banach space of holomorphic
maps $K$ vanishing to order $r$, with norm
\[
 \|K\|_r=\sup_{0<\|z\|<\delta}
                  \frac{\|K(z)\|}{\|z\|^r}.
\]
The operator
\[
 \mathcal T(K)(z)=F^{-1}\bigl(P(Az)+K(Az)\bigr)-P(z)
\]
has Lipschitz constant at most $La^r$ in this norm. The
quantity $\|R\|_r$ is bounded by a constant $C$ independent
of sufficiently small $\delta$. Fix $M>C/(1-La^r)$ and then
shrink $\delta$ so that all arguments of $F^{-1}$ for
$\|K\|_r\le M$ lie in its chosen ball. This closed ball of
maps is preserved by $\mathcal T$. The contraction theorem
gives a holomorphic fixed point $K$. It is real on the real
slice by uniqueness. Thus $H=P+K$ is an analytic solution
with the same $(r-1)$-jet as $h$.

For real $z$ in a smaller ball, $h(z)-H(z)=O(\|z\|^r)$.
Both maps satisfy the inverse equivariance identity there.
Iterating its Lipschitz estimate gives
\[
 \|h(z)-H(z)\|
 \le L^j\|h(A^jz)-H(A^jz)\|
 \le C_1(La^r)^j\|z\|^r.
\]
Letting $j\to\infty$ proves $h=H$. For a $C^\infty$ germ
choose $r$ sufficiently large.
\end{proof}

\begin{theorem}[Analytic projective embeddings for every lattice]
\label{II-general:thm:analytic-embedding}
Let $\Gamma<\SL(n,\R)$ be any lattice, $n\ge3$, acting
analytically on an analytic manifold $M$. Every $C^\infty$
equivariant map from standard or dual projective space, or
its sphere of oriented rays, to $M$ is analytic. Consequently
an equivariantly $C^\infty$ embedded projective covering
hypersurface with finite $\Gamma$-orbit is an analytic
hypersurface with analytic projective identification.
No bound on its normal derivative and no invariant volume
are required.
\end{theorem}
\begin{proof}
Pass to a finite-index stabilizer when a hypersurface is
given. Lemma~\ref{II-general:lem:pinched} supplies $\gamma$ with
an attracting projective fixed point $p$. In a graph chart
at $p$ its action is the linear contraction $A$ of that lemma.
In the oriented case a square fixes both lifts and has the
same contraction property. If $h$ is the equivariant map,
the analytic germ $F=\alpha(\gamma)$ at $h(p)$ satisfies
$Fh=hA$. Lemma~\ref{II-general:lem:analytic-parametrization}
makes $h$ analytic near $p$. Equivariance and minimality
of the source action, from Lemma~\ref{II-general:lem:linear-models},
make it analytic everywhere. The same argument applies to
the dual lattice $\{\gamma^{-T}:\gamma\in\Gamma\}$.
An analytic embedding identifies its image with an analytic
submanifold, proving the final assertion.
\end{proof}

\begin{lemma}[Analyticity of a twice differentiable boundary conjugacy]
\label{II-general:lem:boundary-conjugacy}
Let an analytic action of $\Gamma$ on a compact analytic
manifold $N$ admit a $C^2$ equivariant finite covering, which
is a local diffeomorphism, to a real Grassmannian with its
standard action. Then this covering is analytic. The assertion
also holds for the sphere of oriented rays and for dual actions.
\end{lemma}
\begin{proof}
Choose the element and attracting plane from
Lemma~\ref{II-general:lem:pinched}. A positive power fixes each
point over that plane. At any such point the analytic return
germ, after a linear change of coordinates, is
$F(z)=Az+O(\|z\|^2)$. Choose $b>\|A\|$ with
$\|A^{-1}\|b^2<1$, and shrink a complex neighborhood so that
$\|F(z)\|\le b\|z\|$. The telescoping series for $A^{-d}F^d$
converges normally, since its $d$th difference is bounded by
$C\|A^{-1}\|^{d+1}b^{2d}\|z\|^2$. Its limit is an analytic
linearizer tangent to the identity.

In these linearizing coordinates the given $C^2$ covering,
after matching its derivative, is $J(z)=z+O(\|z\|^2)$ and
commutes with $A$. Hence
\[
 J(z)-z=A^{-d}\bigl(J(A^dz)-A^dz\bigr)\longrightarrow0,
\]
so $J=\id$. The covering is analytic near every point over the
attractor. Finiteness of its fibers gives one target neighborhood
whose full preimage has this property. Minimality from
Lemma~\ref{II-general:lem:linear-models}, together with equivariance,
transports analyticity to every point. Oriented-ray fixed points
are obtained by taking a further square if necessary.
\end{proof}

\subsection{Commuting regular returns and their transverse constraints}

The following input supplies commuting regular elements even when
the lattice has no unipotents. Prasad--Raghunathan prove the
existence of a conjugate Cartan subgroup with a uniform lattice
intersection \cite{PR72}; the exact Cartan formulation is also
stated in \cite[Theorem 7.2, p.~398]{Hurder}.
An independent published formulation covering nonuniform lattices appears in
\cite[Proposition 2, p.~1534]{MozesWeiss}. The selection argument in
\cite[Proposition 8.15]{BRHWGlobal} additionally uses
Prasad--Rapinchuk \cite{PR03,PR05} to control the centralizer
and its rational factors. Those additional factor conclusions
have hypotheses on the algebraic action and are not automatic
for a transverse germ action.

\begin{lemma}[A lattice of commuting split returns]
\label{II-general:lem:cartan-returns}
For every lattice $\Gamma<\SL(n,\R)$ there is $k\in\SL(n,\R)$
such that, writing $A'=kAk^{-1}$ for the conjugate of the
positive diagonal group, $\Lambda=\Gamma\cap A'$ is a uniform
lattice in $A'$.
In the diagonal coordinates given by $k$, $\log(k^{-1}\Lambda k)$
is a full lattice in
\[
 \mathfrak a=\{(h_1,\ldots,h_n):\textstyle\sum_i h_i=0\}.
\]
It meets every nonempty open cone at arbitrarily large norm.
One can simultaneously require finitely many strict homogeneous
linear inequalities and avoid finitely many prescribed proper
linear hyperplanes, whenever these conditions define a
nonempty open cone.
\end{lemma}
\begin{proof}
Apply the cited theorem to the full real diagonal Cartan, obtaining
a conjugate $D$ with $D/(\Gamma\cap D)$ compact. Squaring induces
a continuous surjection onto $A'/(\Gamma\cap A')$, since $D$ is
abelian and positive diagonal elements have square roots. Thus this
quotient is compact, and $\Gamma\cap A'$ is discrete.
The exponential map identifies this component with
$\mathfrak a$, so uniformity says that $\log\Lambda$ is a
full lattice. Fix a bounded fundamental parallelepiped. For
any $H$ in the specified open cone and $T\to\infty$, choose
$H_T\in\log(k^{-1}\Lambda k)$ with $\|H_T-TH\|$ uniformly bounded.
For sufficiently large $T$ all the strict inequalities hold.
Choose $H$ off the finitely many hyperplanes; their distances
from $TH$ grow linearly, so they are also avoided eventually.
\end{proof}

In these diagonal coordinates, the coordinate points of projective space
are common fixed points of $k^{-1}\Lambda k$. At $[e_1]$ put
$\beta_j=\varepsilon_j-\varepsilon_1$, $2\le j\le n$.
These $n-1$ functionals form a basis of $\mathfrak a^*$.
Thus Lemma~\ref{II-general:lem:cartan-returns} supplies commuting
projective returns with all $\beta_j<0$, and with any desired
strict pinching inequalities inside that cone. The derivative
on an invariant normal line is a character of $\Lambda$;
after taking squares its logarithm extends to a unique
$\chi\in\mathfrak a^*$. This observation alone does not
identify $\chi$ with a specified multiple of
$\varepsilon_1$.

\section{Suspension measures and the possible parabolic types}\label{II-series:suspension}\label{II-sec:suspension}

Fix a smooth Riemannian metric on the original compact manifold $M$.
Fiberwise exponential convergence and stable fibers below use its distance,
truncated at $1$ (with distance $1$ between distinct components).
The local coordinate estimates and exponential stable relations are
unchanged by this truncation. Any two smooth Riemannian metrics on the
compact manifold are uniformly comparable.

Let $G=\SL(n,\R)$. Our suspension convention is
$(g,x)\gamma=(g\gamma,\alpha(\gamma)^{-1}x)$, with left $G$ action.
A $G$-invariant Haar-projecting probability corresponds to a
$\Gamma$-invariant probability on $M$, by disintegration on a
measurable fundamental domain; the two constructions are inverse.
Ergodicity is preserved. The same assertion holds over any compact
invariant subset. We use normalized Haar measure on $G/\Gamma$.

Let $K\subset M$ be nonempty, compact and invariant, choose $x\in K$,
and put $K_x=\overline{\Gamma x}$ and $X_{K_x}=(G\times K_x)/\Gamma$.

\begin{proposition}[An opposite minimal parabolic detects the frontier]
\label{II-nv:prop:averaging}
For any minimal parabolic $P^-$ there is an ergodic
$P^-$-invariant probability on $X_{K_x}$ projecting to normalized
Haar measure on $G/\Gamma$.
\end{proposition}
\begin{proof}
Choose a measurable section of $G\to G/\Gamma$ and use $x$
in its fiber to obtain a probability $\lambda_0$ with Haar
projection. The upper triangular subgroup of $\GL(n,\R)$ is soluble
\cite[Chapter I, Section 3, Example, pp.~18--19]{LangAlgebra}; the actual
upper triangular subgroup of $\SL(n,\R)$ embeds into it. Hence it is
amenable by \cite[Corollary~G.2.3, p.~452]{BHV}, and so is its specified
conjugate $P^-$. Apply \cite[Theorem~G.3.1(iv), p.~453]{BHV} to Haar
measure $\nu$ on this group. For every finite set $Q$ and $\varepsilon>0$
it supplies a nonnegative norm-one $L^1$ class $f$ with
$\int|f(g^{-1}p)-f(p)|\,d\nu(p)\le\varepsilon$ for $g\in Q$.
Choose its measurable nonnegative representative and put
$d\eta=f\,d\nu$. This is a probability, and left Haar invariance gives,
for every bounded continuous real function $\psi$,
\[
 \left|\int\psi(gp)\,d\eta(p)-\int\psi(p)\,d\eta(p)\right|
 \le \|\psi\|_\infty\int|f(g^{-1}p)-f(p)|\,d\nu(p).
\]
Replacing the $L^1$ representative preserves both integrals, since
left translations preserve Haar null sets. Choose these densities for
finite sets $Q$ and tolerances $1/(j+1)$. The product directed set of
finite subsets ordered by inclusion and $j\in\mathbb{N}$ produces probability
laws $\eta_d$; pushing its tail filter to these laws gives a nontrivial
filter along which every displayed test defect tends to zero.
Define
\[
 \lambda_d=\int_{P^-}p_*\lambda_0\,d\eta_d(p).
\]
All projections remain Haar. The bundle projection is proper
because its fiber is compact, so these probabilities are uniformly
tight. A weakly convergent subnet has a probability limit.
The Reiter property implies its $P^-$-invariance. Haar is ergodic
under $P^-$ by Howe--Moore, so almost every ergodic component
of this limit still projects to Haar. Choose one.
\end{proof}

The averaging group is a minimal parabolic, which is amenable.
The line parabolic has a nonamenable Levi and is not the averaging
group. Compactness of the fiber and the fixed Haar marginal prevent
escape of mass even though the suspension itself is noncompact.

\begin{proposition}[Possible proper stabilizers]\label{II-nv:prop:parabolics}
For an ergodic $P$-invariant suspension probability with proper
stabilizer $Q$, its connected parabolic type is as follows:
\begin{enumerate}
\item If $n\geq5$, $Q$ stabilizes a line or a hyperplane.
\item If $n=4$, $Q$ can additionally stabilize a two-plane.
      In that case the full-dimensional boundary theorem gives
      a smooth equivariant finite covering
      $M\to\operatorname{Gr}_2(\R^4)$.
\item If $n=3$, $Q$ is a line parabolic, a plane parabolic,
      or the minimal parabolic. In the last case $G/Q$ is the
      $3$-dimensional full flag manifold.
\end{enumerate}
\end{proposition}
\begin{proof}
We use Cartan Lyapunov data along each fixed ray. The weak
$L^1$ integer-Cartan splitting, including its nonergodic version,
is stated in \cite[p.~463, the nonergodic paragraph and the first
paragraph of the proof of Theorem~2.4]{BRHWSmoothErgodic2023};
the full multidirectional limit in that theorem requires stronger
Lorentz moments and is not needed for this fixed-ray construction.

The moment hypothesis is verified on the same original carrier and section.
Fix a finite precompact original atlas near the compact invariant fiber
carrier $K$, and write $D_\gamma(x)$ for the derivative from the selected
coordinate at $x$ to that at $\gamma x$. Compactness of the chart closures
bounds the finitely many generator maps and their inverses in these coordinates.
Applying the derivative chain rule along a shortest word gives a constant
$D>0$ such that
\[
 \log\max\{1,\|D_\gamma(x)\|,\|D_\gamma(x)^{-1}\|\}
       \le D\bigl(\ell_\Gamma(\gamma)+1\bigr)
       \qquad(\gamma\in\Gamma,\ x\in K).
\]
For the inverse bound, the derivative of $\gamma^{-1}$ is evaluated at
$\gamma x$, in the same atlas. Let $s$ be the original quotient section chosen with minimum word length
in each coset, and $m$ the base marginal of the original finite suspension law.
Suppose $\ell_G\circ s\in L^1(m)$, as supplied by
Lemma~\ref{II-power:lem:section} for the original Haar base law. For its actual
return label $\delta(g,b)$ the word bound reads
\[
 \ell_\Gamma(\delta(g,b))
       \le C\bigl(1+2\ell_G(g)+2\ell_G(s(b))\bigr),\qquad C\ge0.
\]
For a compact set $S$ of Cartan parameters, continuity of the Cartan
inclusion $\iota$ bounds $\ell_G(\iota(a))$ by an integer $N_S$ on $S$.
Thus the same derivative cocycle, on its actual carrier section law
$\widehat\mu$, has the pointwise envelope
\[
 H_S(b,x)=D\Bigl[C\bigl(1+2N_S+2\ell_G(s(b))\bigr)+1\Bigr].
\]
The base marginal of $\widehat\mu$ is exactly $m$: the original law gives
full mass to the invariant carrier. Hence $H_S\in L^1(\widehat\mu)$,
and it bounds both logarithmic derivative directions for every $a\in S$
and every carrier point. This proves the compact-uniform moment directly. The construction uses the
original chart coordinates, section word bounds and base moment.
Full ambient metrics and quantitative adapted stable charts are separate
constructions.

The real-Cartan extension is a derived consequence of this integer-ray
result and the compact-uniform symmetric envelope. Work with the same
invariant probability and cocycle, and choose a fixed linear coordinate isomorphism
$\theta:\R^d\to\mathfrak a$. Apply the integer-ray result to the coordinate
time-one action. Its measurable finite splitting admits a measurable
adapted basis. The number of nonzero blocks is at most the fiber
dimension; padding with zero blocks gives a fixed finite family. Project
a fixed fiber-coordinate basis onto these subspaces and choose the first
linearly independent full tuple from the resulting finite family. Repeat each original weight according to the number of selected
columns in its subspace. The resulting basis is fixed throughout the
argument. The integer-ray angle estimate in
\cite[p.~463, proof of Theorem~2.4]{BRHWSmoothErgodic2023}
prevents cancellation between distinct original blocks at their actual
image points. Thus the growth rate of every nonzero vector is the largest
weight among its nonzero basis coordinates, and this gives the exact
integer-ray growth filtration. For cocycles specified almost everywhere,
these applications use a measurable invariant conull core for the
countable integer action; completed-measure outputs have versions for the
original measurable structure with the same basis and weights almost
everywhere.

Fix $a=\theta q$, and apply the same integer-ray theorem once more to
the commuting family consisting of $a$ and the $d$ coordinate actors.
For any one weight $\eta_x$ of this augmented field, write
$\eta_{0,x}=\eta_x(e_0)$ and
$\eta_{j,x}=\eta_x(e_j)$ for $1\le j\le d$.
The integer vector
\[
 z_N=(N,-\lfloor Nq_1\rfloor,\ldots,-\lfloor Nq_d\rfloor)
\]
represents the real actor
$\theta(\{Nq_1\},\ldots,\{Nq_d\})$; this actor and its negative
belong to one fixed compact coordinate window. Here
$\{t\}=t-\lfloor t\rfloor$. If $H\ge0$ is its symmetric
logarithmic envelope, the forward and negative integer rays give
\[
 \int |\eta_x(z_N)|\,d\widehat\mu(x)
       \le 2\int H\,d\widehat\mu.
\]
Indeed, normalize vector growth by its initial norm, bound it by the
sum of the original one-step envelope along the measure-preserving orbit,
and apply Fatou's lemma to its positive part; the negative ray treats
the other part. This uses no logarithmic moment of the chosen basis.
The same argument makes each $\eta_{j,x}$ integrable. Linearity gives
\[
 N\Bigl(\eta_{0,x}-\sum_{j=1}^d q_j\eta_{j,x}\Bigr)
   =\eta_x(z_N)-\sum_{j=1}^d\{Nq_j\}\eta_{j,x}.
\]
The right side has an absolute integral bounded independently of $N$.
Division by $N$ therefore proves
$\eta_{0,x}=\sum_jq_j\eta_{j,x}$ almost everywhere.

The augmented field restricted to its last $d$ actors and the original
coordinate field describe the same integer growth filtrations. If a
change-of-basis coefficient between them is nonzero, inclusion of the
filtration for $z$ bounds one coordinate weight by the other; applying
the same inclusion to $-z$ gives equality. Testing the coordinate
integer vectors identifies their entire real-linear functionals.
Consequently their filtrations agree for every real direction.
The preceding identity now identifies the augmented head filtration
with the original basis filtration evaluated at $q$. It follows that
this one original measurable basis has the required vector growth and
exact filtration along the fixed ray $na$. The conull set may depend on
$a$; no uncountable intersection of regular sets is used.

Finally, for any fixed Cartan translation $g$, the same cocycle identity
intertwines the ray cocycle at $x$ and at $gx$ through $\alpha(g,x)$.
The singleton envelope makes the logarithm of the symmetric operator
height of this transport integrable. Along each fixed measure-preserving
ray this logarithmic height divided by the integer time tends to zero
almost everywhere.
The growth filtrations therefore transport with their dimensions.
A countable dense set of directions and continuity of the finitely many
linear weights identify the full weight multiset, including multiplicities,
up to permutation almost everywhere for this fixed $g$. This asserts
neither an invariant chosen basis nor a constant spectrum under a
nonergodic law. Fixed linear changes of fiber coordinates preserve these
limits and filtrations and add only a constant to the envelope.

For nonresonance we use the original $A$-ergodic theorem
\cite[Proposition~5.1, pp.~967--968]{BRHWFactors}. Its native
spectrum agrees with our full coordinate multiset. In the adapted metric,
the compact-uniform bounds for $\log^+\|Dg\|$ belong to every finite
$L^q$, $q\ge1$ \cite[Claim~2.6 and Theorem~3.1,
pp.~958--959]{BRHWFactors}; on generic directions the largest nonzero
native packet determines the growth filtration. Bounded recurrent
returns of the measurable change of frame compare the two filtrations,
and a countable dense set of directions preserves the full functional
multiset, including multiplicities. Restriction to the finite
$A$-ergodic components preserves this multiset and the same Haar base.
Applying nonresonance on each component and summing their invariance
identities proves it for the original $P$-law. Thus every root absent
from $\mathfrak q$ is positively proportional to a fiber Lyapunov
exponent. In type $A$,
different roots are not positively proportional. Consequently
\[
 \dim(\mathfrak g/\mathfrak q)
   =\#\{\hbox{missing roots}\}\leq n.
\]
Here the quotient is the actual trace-zero matrix space modulo the
infinitesimal flag stabilizer. Projection to the forbidden lower
matrix entries is surjective and has that stabilizer as its kernel,
which proves the dimension identity. A nonmaximal proper parabolic
has this infinitesimal codimension at least $2n-3$.
This is greater than $n$ for $n\ge4$; for $n=3$, equality
occurs precisely for the full-flag parabolic.
A maximal parabolic stabilizing a $k$-plane has codimension
$k(n-k)$. For $n\geq5$, only $k=1,n-1$ satisfy the bound;
for $n=4$, also $k=2$ does. In dimension three the full-flag
quotient has an explicit three-dimensional analytic atlas: the
opposite cell is parametrized by
\[
 (u,v,w)\longmapsto
 \begin{pmatrix}1&0&0\\u&1&0\\w&v&1\end{pmatrix}P.
\]
Its translates cover the quotient, and Gaussian elimination gives
analytic coordinate changes. Conjugation transports these charts
to the corresponding conjugate stabilizers.

In the $n=4$ two-plane case all four fiber dimensions are used by distinct
positive multiples of the four missing roots. There is a Cartan
element negative on all four roots, so every fiber exponent is
negative there. The action is analytic, hence has arbitrarily
high finite regularity. Theorem 4.6 of \cite{BRHW} applies since
$\dim M=\dim G/Q=4$: the two-plane Grassmannian graph charts
are parametrized by $2\times2$ real matrices. This yields the asserted smooth finite
covering, whose covering charts have smooth local inverses.
In particular it is a $C^2$ local diffeomorphism. This conclusion
is stated at the smooth regularity provided by that theorem;
an analytic upgrade is a separate step.
\end{proof}

\begin{proposition}[Levi symmetry isolates the transverse weight]
\label{II-nv:prop:weights}
Suppose the stabilizer is the line parabolic, and choose coordinates
so that the missing roots are
$\beta_j=\varepsilon_j-\varepsilon_1$, $2\leq j\leq n$,
where $\sum_i\varepsilon_i=0$ on the diagonal Cartan algebra.
The fiber Lyapunov multiset, counting multiplicities, is
\begin{equation}\label{II-nv:eq:weights}
 \{c\beta_2,\ldots,c\beta_n,\ q\varepsilon_1\},
       \qquad c>0,\quad q\in\R.
\end{equation}
Every missing-root coarse space has dimension one. The hyperplane
case has the corresponding dual formula.
\end{proposition}
\begin{proof}
First the exponent multiset, including multiplicities, is constant
almost everywhere. Write $P=MAN$. Commuting derivative returns preserve
the Cartan weight multiset: recurrence to sets where the return derivative
and its inverse are bounded gives both growth comparisons. This proves
invariance under $A$ and under its centralizer $M$ directly.
For each $n\in N$, choose $a_j\in A$ with $a_j^{-1}na_j\to e$.
The contraction argument in the strongly continuous Koopman
representation makes every bounded measurable function of the
$A$-invariant multiset also $N$-invariant. Since $A,M,N$ generate $P$,
$P$-ergodicity makes the multiset constant. This retains the full
multiset and its multiplicities on the original probability space.

The Levi normalizer of $A$ preserves the original measure and
acts on this common multiset by its Weyl group
$W_L\simeq S_{n-1}$. This follows by differentiating the
normalizer conjugacy. Equivalently, use recurrence to a set where
that fixed derivative and its inverse are bounded, and compare
Lyapunov limits along the recurrent subsequence; no bounded
global change of frame is required.

Each ray $\R_{>0}\beta_j$ must carry a fiber exponent by
nonresonance-implies-invariance. The group $W_L$ permutes these
$n-1$ rays transitively and therefore forces the same scaling
multiset and multiplicity on each. Since $2(n-1)>n$, exactly one
dimension occurs on each ray, with a common coefficient $c>0$.
The one remaining exponent must be fixed by $W_L$. Its fixed
subspace in the Cartan dual is precisely $\R\varepsilon_1$.
This also shows that the remaining exponent is not on any of
the missing-root rays.
\end{proof}

\begin{corollary}[The expanding transverse branch]
\label{II-nv:cor:expanding}
If $q>0$ in \eqref{II-nv:eq:weights}, Theorem 4.7 of \cite{BRHW}
gives a finite union of compact embedded projective covering
leaves in $M$, with its equivariant smooth boundary structure.
No invariant volume is needed. The conclusion is at the finite
smooth regularity of that theorem; an analytic upgrade is not
asserted here.
\end{corollary}
\begin{proof}
On the closed cone $\beta_j(H)\leq0$ for all $j$, one has
$\varepsilon_1(H)>0$ unless $H=0$: the first coordinate is the
maximum of coordinates summing to zero. Thus $q\varepsilon_1$
is strictly positive on this cone minus zero. The missing-root
coarse spaces are one dimensional by Proposition~\ref{II-nv:prop:weights}.
These are the spectral hypotheses of \cite[Theorem 4.7]{BRHW}.
The projective boundary has finite fundamental group, so its
conclusion (9) gives finitely many compact embedded leaves.
Their spherical or projective parametrizations follow from covering
classification and $\pi_1(\mathbb{RP}^{n-1})\simeq\mathbb Z/2$
\cite[Proposition~1.14, p.~35, Propositions~1.33--1.34,
pp.~61--62, Theorem~1.38, p.~67, and Example~1.43, p.~74]{Hatcher}.
The kernels of the finite leaf-permutation action and the global deck-sign
characters give strictly equivariant parametrizations on a finite-index
subgroup, preserving the same leaves and their $C^r$ inclusions.
\end{proof}

If $q<0$, every fiber exponent is negative on the open contraction
cone. Proposition~5.7 of \cite{BRHW} applies to the full $P$-ergodic
law; the component transfer in the proof of
Proposition~\ref{II-loc:prop:negative-factor} below then permits
Proposition~4.4 to give $c=1$.
This does not give the full-dimensional conclusion of its
Theorem 4.6, since here $\dim M=\dim(G/Q)+1$. If $q=0$, that
contraction argument and the expanding transverse criterion both
fail. With invariant smooth volume the sum of \eqref{II-nv:eq:weights}
vanishes, giving $q=nc>0$; this explains the role of volume in
selecting the expanding branch.

The form $q\varepsilon_1$ is a statement about a Lyapunov
functional. It does not by itself conjugate the transverse
derivative cocycle to a constant character, integrate its
Oseledec line into analytic transversals, or linearize their germs.

The earlier measurable projective-factor theorem
\cite[Theorem 1.10]{BRHWFactors} is available in these dimensions:
$m(\mathfrak{sl}_n)=2n-3\geq n$. It supplies a quasi-invariant
relatively measure-preserving projective factor. It does not
identify its fibers as analytic one-dimensional fibers or turn
its quasi-invariant probability into an invariant one.

\subsection{Integrability for induced cocycles}
\begin{lemma}[An exponentially integrable section]
\label{II-power:lem:section}
For any lattice $\Gamma<G=\SL(n,\R)$, $n\ge3$, choose word
lengths $\ell_G$ for a compact generating neighborhood and
$\ell_\Gamma$ for a finite generating set. There is a Borel
section $r:G/\Gamma\to G$ and $s_0>0$ such that
\begin{equation}\label{II-power:eq:section-moment}
 \int_{G/\Gamma}e^{s_0\ell_G(r(x))}\,dm(x)<\infty.
\end{equation}
Its return cocycle, defined by
$a r(x)=r(ax)\kappa(a,x)$, satisfies
\begin{equation}\label{II-power:eq:return-length}
 \ell_\Gamma(\kappa(a,x))
 \le C\bigl(1+\ell_G(a)+\ell_G(r(x))+\ell_G(r(ax))\bigr).
\end{equation}
The constant $C$ does not depend on $a$.
\end{lemma}
\begin{proof}
The property-$(T)$ input is \cite[Theorem 1.4.15, pp.~54--55]{BHV},
applied to the real field. Positive finite Haar covolume constructs the
invariant quotient probability required by the lattice definition in
\cite[p.~152]{dLS}. The given word length is nonnegative, symmetric,
subadditive, and bounded on compact sets; it is Borel measurable.
Thus Maucourant's theorem, proved in
\cite[Theorem 5.3, pp.~184--185]{dLS}, applies to this same length and
lattice. It yields a Borel set $F$ for which
$F\times\Gamma\to G$, $(f,\gamma)\mapsto f\gamma$, is bijective,
with a finite exponential Haar integral. This is the strict fundamental
domain convention of \cite[Section 2.2, p.~159]{dLS}.
The restriction $F\to G/\Gamma$ is a Borel bijection of standard
Borel spaces, so its inverse gives a Borel section. Its Haar mass is
the positive finite covolume. Pushing normalized restricted Haar
measure to the quotient transfers the exponential integral to the
actual invariant probability $m$.

The separate undistortion theorem of
Lubotzky--Mozes--Raghunathan \cite[Theorem~A, p.~6]{LMR2000},
in the word-length formulation \cite[Theorem~5.1, p.~183]{dLS},
bounds $\ell_\Gamma$ on $\Gamma$ by $C\ell_G$ for these given
generating sets. Apply it to
$\kappa(a,x)=r(ax)^{-1}ar(x)$ and use subadditivity.
Both assertions include nonuniform lattices.
In the scalar boundary and normal collar applications below, apply these
same original BHV and de la Salle theorems to the lattice, Haar measure and
compact word model of that application. The resulting section and moment
bound retain those original objects.

When a bound using only $r(x)$ is needed, select a Borel representative
of minimal integer word length in each coset. Countability of the
lattice permits this measurable minimum selection. Its length is
bounded by that of the preceding exponentially integrable section,
so its exponential moment remains finite for the same quotient law.
Minimality gives $\ell_G(r(ax))\le\ell_G(a)+\ell_G(r(x))$ and hence
\[
 \ell_\Gamma(\kappa(a,x))
 \le C\bigl(1+2\ell_G(a)+2\ell_G(r(x))\bigr).
\]
This stronger estimate is a consequence of the constructed section;
no bound for an arbitrary preassigned section is assumed.
\end{proof}

\begin{lemma}[Word bounds for analytic germs]\label{II-power:lem:word-germs}
Suppose the centered analytic germ cocycle $\alpha$ has a finite
symmetric generating set whose germs have uniformly bounded
holomorphic representatives on a common ball. Then constants
$D,L\ge2$ exist such that,
for every $\gamma\in\Gamma$, $j\ge2$, and almost every $s$,
\begin{equation}\label{II-power:eq:word-germs}
 \|\alpha(\gamma,s)_j\|_{\HS}
       \le D^j L^{j\ell_\Gamma(\gamma)}.
\end{equation}
\end{lemma}
\begin{proof}
Shrink the common generator ball to radius $r>0$. Because all
germs fix zero, the finitely many uniform Cauchy estimates give
$\|\alpha(e,s)(z)\|\le L\|z\|$ there, with one $L\ge2$.
A word of length $d$ is therefore holomorphic on $B_{rL^{-d}}$
and bounded by $r$. Radial Cauchy estimates and the
symmetric-tensor sphere identity bound its degree-$j$ tensor by
$r^{1-j}L^{dj}\binom{n+j-1}{j}^{1/2}$. Increase $D$ to absorb
the fixed radius and the polynomial factor. The estimate applies
to a shortest word; the cocycle identity identifies its germ
with $\alpha(\gamma,s)$. Countability permits one conull set.
\end{proof}

\section{Spectral-gap realization of analytic germ coboundaries}\label{II-series:spectral}\label{II-sec:spectral}

\subsection{Spectral-gap convergence for the full tangent-to-identity germ group}
\label{II-sg:sec:sg-germ}

The following argument uses the forward coboundary equation and a
normalization of its formal transfer. It works in every finite
dimension and imposes no commutativity or finite nilpotence hypothesis
on the analytic germs. Both the input and output radii may be
nonuniform. A noncompact derivative cocycle is not included: the
precise boundary of the argument is stated at the end of the section.

Let $\Gamma\curvearrowright(X,\mu)$ be an ergodic probability-preserving
action of a countable group on a standard probability space. Suppose
there are a finite symmetric set $S\subset\Gamma$ and $a>0$ such that
\begin{equation}\label{II-sg:eq:sg-base}
 \sum_{s\in S}\|u\circ s-u\|_2^2
       \ge a\left\|u-\int_Xu\,d\mu\right\|_2^2
       \qquad(u\in L^2(X,\R)).
\end{equation}
Enlarge $S$ by the identity if necessary. This preserves symmetry
and \eqref{II-sg:eq:sg-base}, and ensures $S\ne\varnothing$.
Put $m=|S|>0$. Summing the scalar inequality in orthonormal coordinates
gives the same inequality for maps into every finite-dimensional real
Hilbert space, with the same constant $a$.

\begin{theorem}[Analytic realization of an arbitrary formal germ coboundary]
\label{II-sg:thm:sg-realization}
Assume \eqref{II-sg:eq:sg-base}. Let $c(\gamma,x)$ be a
coefficient-measurable cocycle of real analytic diffeomorphism germs
of $(\R^n,0)$ tangent to the identity. Suppose a coefficient-measurable
formal tangent-to-identity field $H_x$ satisfies
\begin{equation}\label{II-sg:eq:sg-forward}
 H_{\gamma x}=c(\gamma,x)\circ H_x.
\end{equation}
There is a constant formal tangent-to-identity germ $C$ such that
$K_x=H_x\circ C$ is analytic almost everywhere. Both $K_x$ and
$K_x^{-1}$ admit representatives on measurable positive radii, and
\eqref{II-sg:eq:sg-forward} holds with $K$ as an identity of analytic germs.
No moment assumption or common analytic input radius is required.
\end{theorem}

\begin{corollary}[Analytic germ coboundaries over property-$(T)$ sources]
\label{II-sg:cor:sg-property-T}
Let a countable discrete property-$(T)$ group act ergodically
preserving a standard probability measure. Every coefficient-measurable
cocycle into the full group of real analytic tangent-to-identity
germs of $(\R^n,0)$ is a measurable analytic germ coboundary.
The dimension $n$ is arbitrary and finite. The analytic radii of
the transfer and its inverse may depend on the base point.
\end{corollary}

We prove the theorem first, including its finite-dimensional estimates.
The zero-dimensional case is immediate, so assume $n\ge1$ in these
estimates. The final argument supplies the formal transfer needed for the corollary.

\subsubsection{Formal normalization and a median Poincare inequality}

For a finite measurable real random variable $Y$, let
$\operatorname{med}(Y)=\inf\{t:\mu(Y\le t)\ge1/2\}$.
This is a finite lower median, and
$\operatorname{med}(Y+b)=\operatorname{med}(Y)+b$ for constants $b$.

\begin{lemma}\label{II-sg:lem:sg-normalize}
There is a unique constant formal tangent-to-identity germ $C$ such
that every nonlinear scalar Taylor coefficient of $H_x\circ C$
has lower median zero. Right composition with $C$ preserves
\eqref{II-sg:eq:sg-forward}.
\end{lemma}
\begin{proof}
At degree $d$, a coefficient of $(H_x\circ C)_i$ has the form
$c_{i,\alpha}+R_{i,\alpha}(x)$, where $|\alpha|=d$ and $R$ depends
only on already chosen coefficients of $C$ of smaller degree.
Set $c_{i,\alpha}=-\operatorname{med}(R_{i,\alpha})$ and proceed
by induction. Translation equivariance of the median proves
existence and uniqueness. Associativity proves the last assertion.
\end{proof}

The same triangular construction works after fixing, separately in each
degree, any linear coordinate system on the homogeneous vector-valued
polynomials. In those coordinates, changing only the degree-$d$ part of
$C$ adds an arbitrary constant vector to the degree-$d$ part of $H_x\circ C$.
Subtract its coordinate lower medians and continue by degree. This gives
a unique constant for the chosen coordinate systems, preserves coefficient
measurability, and preserves \eqref{II-sg:eq:sg-forward}. We use this
fixed-coordinate version below after choosing the Hilbert norms.
The resulting constant need not equal the scalar-monomial constant of
Lemma~\ref{II-sg:lem:sg-normalize}.

\begin{lemma}[Dimension-independent median estimate]
\label{II-sg:lem:sg-median-gap}
Let $v:X\to\R^D$ be square integrable. Suppose zero is a median of
each coordinate, in the sense that
$\mu(v_i>0)\le1/2$ and $\mu(v_i<0)\le1/2$. Then
\begin{equation}\label{II-sg:eq:sg-median-gap}
 \int_X\|v\|^2\,d\mu
       \le\frac2a\sum_{s\in S}\int_X\|v(sx)-v(x)\|^2\,d\mu(x).
\end{equation}
Radial projection onto a ball centered at zero preserves this
coordinate-median property.
\end{lemma}
\begin{proof}
For a scalar coordinate $Y$ with $\mathbb EY\ge0$,
\[
 (\mathbb EY)^2\le(\mathbb EY_+)^2
       \le\mu(Y>0)\mathbb E Y_+^2\le\tfrac12\mathbb E Y^2.
\]
For negative mean apply the same argument to $-Y$.
Summing over coordinates gives
$\|\mathbb Ev\|^2\le\tfrac12\mathbb E\|v\|^2$.
Apply the Hilbert-valued version of \eqref{II-sg:eq:sg-base}.
Radial projection multiplies each nonzero vector by a strictly
positive scalar, so the signs of its coordinates are unchanged.
\end{proof}

\begin{lemma}[Expansion of small measurable sets]
\label{II-sg:lem:sg-expansion}
Put $\vartheta=1+a/(4m)>1$ and
$\mathcal N(A)=A\cup\bigcup_{s\in S}sA$.
If $\mu(A)\le1/2$, then
\begin{equation}\label{II-sg:eq:sg-expansion}
 \mu(\mathcal N(A))\ge\vartheta\mu(A).
\end{equation}
\end{lemma}
\begin{proof}
Apply \eqref{II-sg:eq:sg-base} to $1_A$. Since
$\mu(sA\setminus A)=\tfrac12\mu(sA\mathbin\triangle A)$,
\[
 \sum_{s\in S}\mu(sA\setminus A)
 \ge\tfrac a2\mu(A)(1-\mu(A))\ge\tfrac a4\mu(A).
\]
The union of these $m$ sets has measure at least their average.
\end{proof}

\subsubsection{Analytic Hilbert norms on finite jets}

Choose an integer $p>n/2$. For scalar polynomials in $n$ variables set
\[
 |f|_p^2=\sum_{\alpha\in\mathbb N^n}
                    (1+|\alpha|)^{2p}|f_\alpha|^2.
\]
The constant
$B_p=(\sum_{\alpha\in\mathbb N^n}(1+|\alpha|)^{-2p})^{1/2}$
is finite. Weighted convolution and Young's inequality give
\[
 |fg|_p\le 2^p B_p|f|_p|g|_p.
\]
Indeed, for $\alpha=\beta+\gamma$,
$(1+|\alpha|)^p\le2^{p-1}((1+|\beta|)^p+(1+|\gamma|)^p)$;
apply the $\ell^2*\ell^1$ convolution inequality to these two terms,
and use $\|f_\alpha\|_{\ell^1}\le B_p|f|_p$.
Choose $L\ge\max\{1,2^pB_p\}$ and put $\|f\|_*=L|f|_p$.
Then $\|fg\|_*\le\|f\|_*\|g\|_*$.
For vector polynomials use the Hilbert direct-sum norm of the
component norms. Write $J_0=\|z\|_*>0$.

Let $\Pi_N$ discard all terms above total degree $N$.
It is an orthogonal projection and hence a contraction.
The scalar finite jet algebra, with product $\Pi_N(fg)$, inherits
the preceding submultiplicative norm, with constants independent of $N$.
These are analytic coefficient weights; there are no factorials.

\begin{lemma}[Quadratic composition estimate on a large set]
\label{II-sg:lem:sg-composition}
For every $\beta>0$ there are a measurable $B\subset X$ with
$\mu(B)\le\beta$ and constants $t_0>0$, $A\ge1$ such that,
for all $x\notin B$, $s\in S$, $N$, and vector jets $V$ vanishing
at zero with $\|V\|_*\le t_0$,
\begin{equation}\label{II-sg:eq:sg-quadratic}
 \|\Pi_N(c(s,x)\circ V-V)\|_*\le A\|V\|_*^2.
\end{equation}
\end{lemma}
\begin{proof}
Because $S$ is finite and every $c(s,x)$ is analytic, the increasing
measurable sets on which all nonlinear coefficients satisfy
$|[z^\alpha](c_i(s,x)-z_i)|\le b^{|\alpha|+1}$, with integer
$b\ge1$, exhaust a conull set. Choose $b$ so that its exceptional
set has measure at most $\beta$. The degree-$k$ sum of absolute
coefficients of each component is at most
$\binom{k+n-1}{n-1}b^{k+1}$. Submultiplicativity in the finite jet
algebra bounds the left side of \eqref{II-sg:eq:sg-quadratic} by
\[
 \sqrt n\sum_{k\ge2}\binom{k+n-1}{n-1}b^{k+1}\|V\|_*^k.
\]
This scalar series converges for $\|V\|_*<1/b$ and is quadratic
at zero. Choose $t_0<1/b$ and a uniform $A\ge1$.
\end{proof}

\subsubsection{A finite-jet bound independent of the order}

For each degree $d\ge2$, fix an orthonormal basis of the space of
homogeneous vector polynomials of degree $d$ in the norm $\|\cdot\|_*$.
Apply the fixed-coordinate triangular normalization described after
Lemma~\ref{II-sg:lem:sg-normalize} to these bases. Write its single constant
as $C$ and put $K_x=H_x\circ C$. Thus every coordinate of the
homogeneous degree-$d$ part of $K_x$ in its chosen basis has lower
median zero. Use the direct sum of these same bases for every truncation;
no basis is changed when the truncation order or radius changes.
For $r>0$ define
\[
 h_{N,x}^{(r)}(z)=\Pi_N(K_x(rz)-rz),\qquad
 q_{N,r}(x)=\|h_{N,x}^{(r)}\|_*.
\]
At every finite $N$ these are finite measurable coefficient vectors.
The equation for $K$ gives the exact finite-jet identity
\begin{equation}\label{II-sg:eq:sg-jet-equation}
 h_{N,sx}^{(r)}-h_{N,x}^{(r)}
   =\Pi_N\bigl(c(s,x)\circ(rz+h_{N,x}^{(r)})
                         -(rz+h_{N,x}^{(r)})\bigr).
\end{equation}
The input map vanishes at zero, so coefficients of order at most
$N$ depend only on its $N$-jet. No analytic interpretation of $K$
has been used in this identity.

Fix $\delta=1/4$ and put $D=2m/a$. Choose $\beta>0$ so small that
\begin{equation}\label{II-sg:eq:sg-beta}
 16D\beta\left(1+\frac1{\vartheta-1}\right)<\frac\delta4.
\end{equation}
Obtain $B,t_0,A$ from Lemma~\ref{II-sg:lem:sg-composition} with this
$\beta$. Next choose $\varepsilon>0$ so small that
$3\varepsilon\le t_0$, $9A\varepsilon<1$, and
\begin{equation}\label{II-sg:eq:sg-epsilon}
 D\left(81A^2\varepsilon^2+
       16\delta\vartheta^{-\lfloor1/(9A\varepsilon)\rfloor}\right)
       <\frac\delta2.
\end{equation}
This is possible after $A$ is fixed. Thus the large-set coefficient
bound need not have any quantitative dependence on $\beta$.
Set $r_0=\varepsilon/J_0$.

\begin{proposition}[Uniform tightness of the normalized analytic norm]
\label{II-sg:prop:sg-tightness}
For every finite $N$,
\begin{equation}\label{II-sg:eq:sg-tightness}
 \mu\{x:q_{N,r_0}(x)>\varepsilon\}\le\delta.
\end{equation}
The constants $r_0,\varepsilon,\delta$ are independent of $N$.
\end{proposition}
\begin{proof}
Suppose the assertion fails for some $N$.
Pointwise, $q_{N,r}$ is continuous and nondecreasing in $r$, tends
to zero as $r\downarrow0$, and is strictly increasing when the
jet is nonzero. A first crossing of the probability threshold gives
$0<r\le r_0$ with
\begin{equation}\label{II-sg:eq:sg-crossing}
 \mu(q_{N,r}>\varepsilon)\le\delta,
 \qquad \mu(q_{N,r}\ge\varepsilon)\ge\delta.
\end{equation}
More explicitly take the infimum of the radii at which
$\mu(q_{N,r}>\varepsilon)>\delta$. Continuity from below of these
events gives the first inequality; the right limit gives the second.
This handles atoms at the threshold.

Abbreviate $h_x=h_{N,x}^{(r)}$, $q(x)=\|h_x\|_*$, and put
$\eta=9A\varepsilon^2$ and
$J=\lfloor\varepsilon/\eta\rfloor=\lfloor1/(9A\varepsilon)\rfloor$.
If $x\notin B$ and $q(x)\le t\le2\varepsilon$, then
$\|rz+h_x\|_*\le\varepsilon+t\le3\varepsilon$ and
\eqref{II-sg:eq:sg-jet-equation} gives
\begin{equation}\label{II-sg:eq:sg-local-step}
 \|h_{sx}-h_x\|_*\le\eta,\qquad q(sx)\le t+\eta.
\end{equation}

Let $A_t=\{x:q(x)>t\}$. Symmetry of $S$ and
\eqref{II-sg:eq:sg-local-step} imply
\[
 \mathcal N(A_{t+\eta})\subset A_t\cup B
       \quad(\varepsilon\le t\le2\varepsilon).
\]
Indeed, if $x=sy$ with $y\in A_{t+\eta}$ and $x\notin A_t\cup B$,
apply \eqref{II-sg:eq:sg-local-step} at $x$ with $s^{-1}$ to get
$q(y)\le t+\eta$, a contradiction. All relevant tail sets have
measure at most $\delta<1/2$ by \eqref{II-sg:eq:sg-crossing}.
Lemma~\ref{II-sg:lem:sg-expansion} therefore gives
\[
 \mu(A_{t+\eta})\le\vartheta^{-1}(\mu(A_t)+\beta).
\]
Iterating $J$ times, and using $\varepsilon+J\eta\le2\varepsilon$,
yields
\begin{equation}\label{II-sg:eq:sg-tail}
 \mu(q>2\varepsilon)
    \le\delta\vartheta^{-J}+\frac\beta{\vartheta-1}.
\end{equation}

Radially project the Hilbert vector $h_x$ to the ball of radius
$2\varepsilon$, obtaining $u_x$. This field is bounded and square
integrable. In the fixed orthonormal homogeneous bases chosen above,
the degree-$d$ coordinate of $h_x$ is $r^d$ times the corresponding
normalized coordinate of $K_x-\id$. Since $r>0$, its signs and
zero-median inequalities are preserved. Radial projection also preserves
those signs, so Lemma~\ref{II-sg:lem:sg-median-gap} applies.
On $X\setminus(B\cup A_{2\varepsilon})$, nonexpansiveness of
radial projection and \eqref{II-sg:eq:sg-local-step} give
$\|u_{sx}-u_x\|_*\le\eta$. Everywhere the difference is at most
$4\varepsilon$. Consequently, for each $s\in S$,
\[
 \int_X\|u_{sx}-u_x\|_*^2\,d\mu
 \le81A^2\varepsilon^4+
       16\varepsilon^2\left(
          \beta+\delta\vartheta^{-J}+\frac\beta{\vartheta-1}\right).
\]
On the other hand, \eqref{II-sg:eq:sg-crossing} gives
$\int_X\|u_x\|_*^2\,d\mu\ge\delta\varepsilon^2$.
The median Poincare inequality now implies
\[
 \delta\le D\left(81A^2\varepsilon^2+
       16\delta\vartheta^{-J}+
       16\beta\left(1+\frac1{\vartheta-1}\right)\right)<\delta,
\]
where the strict inequality follows from
\eqref{II-sg:eq:sg-beta}--\eqref{II-sg:eq:sg-epsilon}. This contradiction
proves \eqref{II-sg:eq:sg-tightness}.
\end{proof}

\begin{proof}[Completion of Theorem~\ref{II-sg:thm:sg-realization}]
The norms $q_{N,r_0}(x)$ increase with $N$. Hence
\[
 E=\bigcap_N\{x:q_{N,r_0}(x)\le\varepsilon\}
 \quad\hbox{satisfies}\quad\mu(E)\ge3/4.
\]
For $x\in E$, every homogeneous Sobolev coefficient of degree $d$
is bounded by $\varepsilon r_0^{-d}$. Passing to the canonical
symmetric multilinear coefficient $T_d(x)$ costs at most the
exponential factor $(n+1)^d$. Thus, for $d\ge2$,
\[
 \|T_d(x)\|\le\varepsilon\bigl((n+1)/r_0\bigr)^d,
\]
and the canonical series has a positive convergence radius.
The positive-radius locus is measurable: it is the countable union,
over integers $j\ge0$, of the sets on which
$\|T_d(x)\|\le(j+1)^{d+1}$ for every $d$.
Equation~\eqref{II-sg:eq:sg-forward}, together with analyticity of
$c(\gamma,x)$ and its inverse, preserves convergence in both directions.
Ergodicity therefore makes this locus conull.

On that locus sum the canonical series. Its derivative at zero is
the identity, so the analytic inverse function theorem gives an
invertible analytic germ with the prescribed full Taylor tower.
Use the identity germ on the measurable exceptional null set.
The resulting field and its inverse have measurable coefficients;
inverse coefficients are finite algebraic expressions in the original
jet coordinates. For either field choose the least integer $j(x)\ge0$
whose exponential bound
$\|T_d(x)\|\le(j(x)+1)^{d+1}$ holds for all $d$.
This index is measurable, and $1/(2(j(x)+1))$ is a positive
measurable radius for its canonical representative.
Taking the minimum of the two radii gives positive measurable domains
for the map and inverse. Faithfulness of the analytic Taylor tower
turns the formal equations into equations of germs. Countability
then supplies one measurable invariant conull set on which all
these equations hold simultaneously.
\end{proof}

\subsubsection{Property \texorpdfstring{$(T)$}{(T)} supplies the formal transfer}

The next lemma is a finite-dimensional instance of measurable
cohomology vanishing for property-$(T)$ actions; see
Anantharaman-Delaroche \cite{AnantharamanDelaroche} for the
general groupoid formulation. We include a direct proof to
make clear that no integrability of the additive cocycle is needed.

\begin{lemma}\label{II-sg:lem:sg-additive}
Every finite measurable $\R^D$-valued additive cocycle over an
ergodic probability-preserving action of a countable property-$(T)$
group is a measurable coboundary, without moment assumptions.
\end{lemma}
\begin{proof}
For an additive cocycle $f$ use the unitary representation on
$L^2(X\times\R^D,\mu\times dz)$ given by
\[
 (\pi(\gamma)F)(x,z)
   =F(\gamma^{-1}x,z-f(\gamma,\gamma^{-1}x)).
\]
Normalized indicators of expanding cubes in $\R^D$, constant in
$x$, are almost invariant. For every fixed finite translation the
relative measure of the symmetric difference tends to zero, and
the squared displacement is bounded by two. Dominated convergence
therefore applies for each group element, and for finite sets of them.
Property $(T)$ gives a nonzero invariant vector $F$.
Its squared fiber norm is invariant and thus a positive constant
almost everywhere. Normalize $|F(x,\cdot)|^2$ to a probability density.
Its coordinate marginal lower medians form a finite measurable
vector $v_x$. Invariance of $F$ translates these densities, giving
$v_{\gamma x}=v_x+f(\gamma,x)$.
\end{proof}

\begin{proof}[Proof of Corollary~\ref{II-sg:cor:sg-property-T}]
Ordinary property $(T)$ gives \eqref{II-sg:eq:sg-base} on the orthogonal
complement of constants; this is its Kazhdan spectral-gap property
\cite[Definition 1.1.3, p.~33 in the author manuscript]{BHV}.
For a real orthogonal representation, apply the complex Kazhdan
inequality to its complexification. A complex fixed vector has fixed
real and imaginary parts, so the complexification of the real
orthogonal complement has no fixed vectors. The inclusion of the
real space is isometric and preserves every displacement. Thus the
same Kazhdan constant $\varepsilon$ works over the reals.
For a discrete group the compact Kazhdan set $S$ is finite; for
each nonzero vector one summand has squared displacement at least
$\varepsilon^2\|v\|^2$, giving the finite-energy inequality with
constant $\kappa=\varepsilon^2$. Enlarge $S$ to make it symmetric.
We construct compatible formal transfers degree by degree.
Suppose one has been constructed through degree $N-1$, and choose
its polynomial lift $H_x^0$. The residual $N$-jet cocycle
\[
 (H_{\gamma x}^0)^{-1}\circ c(\gamma,x)\circ H_x^0
       =\id+f_N(\gamma,x)\pmod{\text{degree }>N}
\]
has only its degree-$N$ term left. Composition of such jets adds
these homogeneous terms, so $f_N$ is an additive cocycle in the
finite-dimensional real space of homogeneous polynomial maps.
Lemma~\ref{II-sg:lem:sg-additive} gives $v_{N,\gamma x}-v_{N,x}=f_N(\gamma,x)$.
Then $H_x^0\circ(\id+v_{N,x})$ is a transfer through degree $N$,
without changing its previous coefficients. Start with the identity
at degree one and intersect the countably many conull sets.
Theorem~\ref{II-sg:thm:sg-realization} makes the resulting formal transfer
analytic after a constant formal normalization.
\end{proof}

\subsubsection{Arbitrary nonlinear centralizers}

\begin{corollary}[Centralizer-valued analytic normalization]
\label{II-sg:cor:sg-centralizer}
Let $\rho:L\to\mathrm{GL}(n,\R)$ be any linear representation.
Let $Z^\omega_U(\rho)$ be the group of analytic tangent-to-identity
germs commuting with every $\rho(\ell)$, $\ell\in L$.
Every measurable $Z^\omega_U(\rho)$-valued cocycle over an ergodic
probability-preserving property-$(T)$ action has a measurable analytic
transfer taking values in $Z^\omega_U(\rho)$ itself.
No commutativity, finite-dimensionality, or weight condition is
imposed on this centralizer.

In particular, if an analytic germ cocycle already has the form
$\alpha(\gamma,x)=c(\gamma,x)\rho(\gamma)$ with
$c(\gamma,x)\in Z^\omega_U(\rho)$, then it is analytically
cohomologous to $\rho$, with measurable positive radii.
\end{corollary}
\begin{proof}
The degree-$d$ homogeneous coefficients of a formal map commuting
with $\rho$ form the linear subspace
\[
 V_d=\{P:P(\rho(\ell)z)=\rho(\ell)P(z)\text{ for every }\ell\in L\}.
\]
Composition, inversion, and truncation preserve equivariance.
The finite-jet construction above therefore stays in the formal
centralizer, using the additive cocycle lemma on $V_d$ at degree $d$.

Write $H_x$ for the resulting measurable formal centralizer-valued
transfer. Apply Theorem~\ref{II-sg:thm:sg-realization} to this formal
transfer in the full germ group. It gives a constant formal map $C$
and a measurable analytic germ field $u_x$ whose formal series is
$H_x\circ C$ almost everywhere. We do not require $C$ to converge
or to commute with $\rho$.

Choose one point $x_0$ where this formal identity holds and set
$v_x=u_x\circ u_{x_0}^{-1}$. Its formal series is
\[
 H_x\circ C\circ(H_{x_0}\circ C)^{-1}
       =H_x\circ H_{x_0}^{-1}
\]
almost everywhere. This series commutes with every $\rho(\ell)$.
Faithfulness of Taylor series for analytic germs therefore shows
that $v_x$ itself belongs to the analytic centralizer. The formal
transfer equation likewise gives
$v_{\gamma x}=c(\gamma,x)\circ v_x$ almost everywhere for each
$\gamma$. Replace $v$ by the identity on a measurable null set
containing its exceptional centralizer values.

Finally put $A_x=D_0v_x$ and $K_x=v_x\circ A_x^{-1}$, regarding
$A_x$ as a linear germ. Differentiating the transfer equation gives
$A_{\gamma x}=A_x$, since $c$ is tangent to the identity.
Moreover $A_x$ commutes with $\rho$, by differentiating the
centralizer identity. Thus $K$ is tangent to the identity, still
centralizes $\rho$, and satisfies
$K_{\gamma x}=c(\gamma,x)\circ K_x$ almost everywhere.
Composition, inversion and differentiation preserve coefficient
measurability. The measurable-radius construction supplies positive
radii for $K_x$ and its inverse. Taking a measurable invariant
conull core makes the transfer equations hold simultaneously for
all elements of the countable acting group.

For the final assertion, commutation with $\rho$ makes $c$ an
ordinary cocycle. If $K_{\gamma x}=c(\gamma,x)K_x$ is the transfer
just constructed, then
$\alpha(\gamma,x)K_x=K_{\gamma x}\rho(\gamma)$ because $K_x$
commutes with $\rho$.
\end{proof}

\subsection{Analytic realization of a formally fixed submanifold}
\label{II-sec:fixed-subspace-realization}

The preceding theorem has a useful extension in which the target
derivative is noncompact but fixes the source directions. This is
different from linearizing a full germ with noncompact derivative:
there is no linear substitution in the source variables of the
forward equation below.

\begin{theorem}[A fixed subspace and spectral gap]
\label{II-sg:thm:fixed-subspace}
Let a countable group $\Delta$ act ergodically preserving a
probability $\mu$ on $X$, with the spectral gap
\eqref{II-sg:eq:sg-base}. Let $V$ be a finite-dimensional real vector
space, $E\subset V$ a subspace of dimension $r\ge1$, and
$\rho:\Delta\to\GL(V)$ a representation fixing $E$ pointwise.
Suppose that some $a\in\Delta$ is real diagonalizable, that the
$1$-eigenspace of $\rho(a)$ is $E$, and that all its other
eigenvalues have modulus different from one.

Let $F(\gamma,x)$ be a coefficient-measurable cocycle of analytic
diffeomorphism germs of $(V,0)$ with derivative $\rho(\gamma)$.
Suppose that a measurable formal map
$I_x:(E,0)\to(V,0)$ has derivative the inclusion $\iota$ and
satisfies
\begin{equation}\label{II-sg:eq:fixed-forward}
                  I_{\gamma x}=F(\gamma,x)\circ I_x.
\end{equation}
Then there is a constant formal diffeomorphism $C$ of $(E,0)$,
tangent to the identity, such that $I_x\circ C$ is an actual
analytic embedded germ almost everywhere. Its radius can be
measurable and nonuniform. No moment or common input-radius
assumption is imposed.
\end{theorem}

\begin{proof}
Choose an eigenbasis of $\rho(a)$ whose first $r$ vectors span
$E$, make it orthonormal, and let $P:V\to E$ be orthogonal
projection. The complementary eigenspace is denoted by $W$.
It need not be invariant under all of $\Delta$.
Use the submultiplicative analytic coefficient Hilbert norm
constructed above, now in $r$ source variables and with $\dim V$
output components. Its constants are independent of the jet order.
For each homogeneous degree $d\ge2$, fix an orthonormal basis
$\mathcal B_d$ of the space of $E$-valued degree-$d$ polynomials
with its induced coefficient Hilbert norm. Apply the fixed-linear-coordinate
variant of Lemma~\ref{II-sg:lem:sg-normalize} to $P I_x$, using
these bases. The degree-$d$ correction is the negative median
vector in $\mathcal B_d$, pulled back to the actual homogeneous
polynomial; lower degrees remain unchanged. Thus one constant
formal $C$ makes every $\mathcal B_d$ coordinate of $P K_x$,
where $K_x=I_x\circ C$, have lower median zero. These fixed
Hilbert coordinates are used throughout this proof; their
identification with scalar monomial coefficients is unnecessary.
Right composition preserves \eqref{II-sg:eq:fixed-forward}.
Enlarge the finite symmetric spectral-gap set $S$ to include
$a,a^{-1}$. For every $\beta>0$, outside a set of measure at most
$\beta$, the finitely many analytic maps $F(s,x)$ satisfy
\begin{equation}\label{II-sg:eq:fixed-quadratic}
 \|\Pi_N(F(s,x)\circ Z-\rho(s)Z)\|_*
          \le A\|Z\|_*^2\qquad(s\in S,\ \|Z\|_*\le t_0),
\end{equation}
with $A,t_0$ independent of $N$. The proof is precisely the
coefficient majorant argument of
Lemma~\ref{II-sg:lem:sg-composition}; only the subtracted linear
term has changed.

The convergence event for $K_x$ is invariant, because $F(s,x)$
and its inverse are analytic. Suppose for contradiction that it
has measure zero. Put $\delta=1/4$, and choose
$\beta_j\downarrow0$. Obtain $A_j,t_j$ from
\eqref{II-sg:eq:fixed-quadratic}, with exceptional set $B_j$.
Choose $\epsilon_j\downarrow0$ so that
$A_j\epsilon_j\to0$ and $\epsilon_j/t_j\to0$.
Choose $r_j^0>0$ with
$\|\iota(r_j^0 z)\|_*\le\epsilon_j$.
Since $K_x$ diverges almost everywhere, some finite jet order
$N_j$ has
\[
 \mu\{\|\Pi_{N_j}(K_x(r_j^0z)-\iota(r_j^0z))\|_*
                                      >\epsilon_j\}>\delta.
\]
The first-crossing argument used in
\eqref{II-sg:eq:sg-crossing} gives $0<r_j\le r_j^0$ for which
the finite coefficient vector
\[
 v_j(x)=\epsilon_j^{-1}
       \Pi_{N_j}(K_x(r_jz)-\iota(r_jz)),\qquad
 q_j(x)=\|v_j(x)\|_*
\]
satisfies
\begin{equation}\label{II-sg:eq:fixed-crossing}
             \mu(q_j>1)\le\delta,\qquad
             \mu(q_j\ge1)\ge\delta.
\end{equation}
Because $\rho(s)\iota=\iota$, the exact finite-jet equation
and \eqref{II-sg:eq:fixed-quadratic} imply, for each fixed $R$,
\begin{equation}\label{II-sg:eq:fixed-local-linear}
 v_j(sx)=\rho(s)v_j(x)+o_j(1)
       \quad\text{if }x\notin B_j,\ q_j(x)\le R.
\end{equation}
The error is uniform on this set: its norm is at most
$A_j\epsilon_j(1+R)^2$. All sufficiently large $j$ lie in the
domain required by \eqref{II-sg:eq:fixed-quadratic}.

\emph{A finite-dimensional matrix records the tangential mass.}
Write $v_j$ as a matrix whose rows correspond to the fixed
output basis and whose columns are its orthonormal coefficient
coordinates, including the coefficient weights. Define
\[
                 Q_j(x)=v_j(x)v_j(x)^{\mathsf T}.
\]
These are positive semidefinite matrices of the fixed size
$\dim V$, even though the number of coefficient columns grows.
Moreover $\operatorname{tr}Q_j=q_j^2$. Pass to a subsequence
whose probability laws converge vaguely to a finite measure
$\nu$ on the cone of positive semidefinite matrices. The mass
of $\nu$ is at least $1-\delta$: use compactly supported
cutoffs equal to one on $\{\operatorname{tr}Q\le1\}$ and
\eqref{II-sg:eq:fixed-crossing}.

For every $s\in S$, the measure $\nu$ is invariant under
$Q\mapsto\rho(s)Q\rho(s)^{\mathsf T}$. Here it is essential
to check both directions of \eqref{II-sg:eq:fixed-local-linear}.
For a compactly supported continuous test function, its value
at either $Q_j(sx)$ or $\rho(s)Q_j(x)\rho(s)^{\mathsf T}$ can
be nonzero only if one of the corresponding coefficient norms
is bounded. Outside $B_j\cup s^{-1}B_j$, the forward equation
or the equation for $s^{-1}$ then bounds the other norm and
makes the two test-function values differ by $o_j(1)$.
Integration, invariance of $\mu$, and vague convergence prove
the assertion. This argument does not assume tightness of the
entire sequence of laws.

Apply recurrence for this finite invariant measure to the
diagonal action of $a$. A diagonal entry corresponding to an
eigenvalue $\lambda$ off $E$ is multiplied by $\lambda^2$.
Since $|\lambda|\ne1$, recurrence forces that entry to vanish
$\nu$-almost everywhere. Positive semidefiniteness then makes
its entire row and column zero. Thus $\nu$ is supported on
matrices with image contained in $E$. The compact-set inequality
for vague convergence consequently gives, for every $\eta>0$,
the bounded-region estimate needed below:
\begin{equation}\label{II-sg:eq:fixed-off-axis}
 \mu\{q_j\le2,\ \|(I-P)v_j\|_*\ge\eta\}\longrightarrow0.
\end{equation}
Indeed this event corresponds to a compact set of matrices
on which the trace of the $W$ block is at least $\eta^2$.

\emph{Expansion rules out escape of coefficient mass.}
Every $\rho(s)$ fixes $E$ pointwise. Hence
\[
 \|\rho(s)v-v\|_*
       \le\|\rho(s)-I\|\,\|(I-P)v\|_*.
\]
Equations \eqref{II-sg:eq:fixed-local-linear} and
\eqref{II-sg:eq:fixed-off-axis} imply the following for each
$\eta>0$: outside a set $D_{j,\eta}$ with measure tending to
zero, simultaneously for $s\in S$,
\[
              q_j(x)\le2\quad\Longrightarrow\quad
              q_j(sx)\le q_j(x)+\eta.
\]
Set $A_{j,t}=\{q_j>t\}$ and use symmetry of $S$. For
$1\le t\le2-\eta$,
\[
 \mathcal N(A_{j,t+\eta})\subset A_{j,t}\cup D_{j,\eta}.
\]
All these tail sets have measure at most $\delta$.
Lemma~\ref{II-sg:lem:sg-expansion}, with its fixed $\vartheta>1$,
therefore gives, on taking $\eta=1/J$ and iterating $J$ times,
\[
 \limsup_j\mu(q_j>2)\le\delta\vartheta^{-J}
                         \qquad(J\ge1).
\]
It follows that $\mu(q_j>2)\to0$. Combining this with
\eqref{II-sg:eq:fixed-off-axis} shows that the off-$E$
coefficient norm tends to zero in probability.

\emph{The fixed coefficients contradict their median normalization.}
Use the direct sum of the fixed bases $\mathcal B_d$ at the
chosen truncation order. Its isometric image in the ambient
coefficient-row space extends to an orthonormal basis; the added
coordinates vanish on $Pv_j$. Each degree-$d$ coordinate of
$Pv_j$ is the corresponding normalized coefficient multiplied by
$\epsilon_j^{-1}r_j^d>0$, so it has median zero. Radially project
$Pv_j$ onto the coefficient Hilbert ball of radius $2$, obtaining
$u_j$. This multiplies every coordinate by the same positive
factor and preserves the two median sign inequalities. On
$\{q_j\le2\}\setminus B_j$,
\eqref{II-sg:eq:fixed-local-linear}, the preceding convergence
in probability, and $\rho(s)|_E=I$ imply
$\|u_j(sx)-u_j(x)\|_*\to0$ in probability. The complement
has measure tending to zero, and the difference is everywhere
bounded by $4$. Thus its squared integral tends to zero for
each $s\in S$. The dimension-independent median inequality
\eqref{II-sg:eq:sg-median-gap} yields
\[
                         \int\|u_j\|_*^2\,d\mu\longrightarrow0.
\]
But \eqref{II-sg:eq:fixed-crossing}, the tail estimate, and
\eqref{II-sg:eq:fixed-off-axis} imply
\[
 \liminf_j\mu\{1\le q_j\le2,\ \|Pv_j\|_*\ge1/2\}
                                               \ge\delta.
\]
On this set the radial projection does not change $Pv_j$,
contradicting the last integral limit.

This proves convergence on a positive-measure set, hence
almost everywhere by ergodicity, for the same normalized field
$K_x=I_x\circ C$. Fixed linear changes from the adapted
eigen-coordinates back to the original spaces preserve positive
convergence radius and transport that same constant $C$.
Coefficient root tests give measurable positive radii.
The derivative is the inclusion; the inverse function theorem
for a linear left projection of $K_x$ gives an analytic inverse
parameter and makes $K_x$ an analytic embedded germ. The
coefficientwise measurable construction gives jointly measurable
representatives and inverse parameters. The formal forward
identities for this same $C$ are now actual identities of
analytic germs.
\end{proof}

\begin{remark}[Fixed and moving source parameters]
The source of \eqref{II-sg:eq:fixed-forward} is fixed pointwise.
Consequently the number of Taylor coefficients can be removed
from the compactness argument by the matrix $Q_j$. For a full
linearization equation with a nontrivial substitution in the
source, this matrix no longer transforms only by the fixed
output representation. The theorem does not assert general
analytic linearization of noncompact germ cocycles.
\end{remark}

\section{Analytic patches from coordinate slices}\label{II-series:common-cross}
\subsection{An analytic patch from almost-everywhere coordinate curves}

\begin{lemma}[Almost-everywhere separate analyticity]
\label{II-coordinate_recognition:lem:cross}
Let $I,J\subset\R$ be nonempty open intervals and let
$f:I\times J\to\R^N$ be measurable. Suppose that for almost every
$y$ its horizontal slice agrees almost everywhere with a real-analytic
map on $I$, and for almost every $x$ its vertical slice agrees almost
everywhere with a real-analytic map on $J$. Then some nonempty open
rectangle has a jointly real-analytic representative of $f$.
The same holds for maps into a compact real-analytic manifold.
\end{lemma}
\begin{proof}
We give the null-set argument, since ordinary separate analyticity
at every point is not assumed. Work first with scalar values, or
with finitely many scalar components at once. Replace $I,J$ by
relatively compact subintervals. By a countable selection of complex
neighborhoods and bounds, there are measurable sets $A\subset I$,
$B\subset J$ of positive measure, complex domains $D\supset I$,
$E\supset J$, and a finite $K$ such that the slice representatives
\[
 g_b\in\mathcal O(D),\quad b\in B,\qquad
 h_a\in\mathcal O(E),\quad a\in A
\]
are bounded by $K$. To justify measurability in this selection,
bounded holomorphic functions on each of the countably chosen
domains form a compact metrizable family in the compact-open
topology. Agreement almost everywhere with a measurable slice
is a measurable relation and determines the representative uniquely.
In particular the selections can be made in the completed measure
spaces. Fubini gives $g_b(a)=h_a(b)$ for almost every $(a,b)$.

Choose a countable set $\{a_d\}\subset A$ dense in the essential
support $A_*$ of $A$, consisting of points for which the preceding
identity holds for almost every $b$. Discard one null subset of $B$
so that
\[
                 g_b(a_d)=h_{a_d}(b)\quad\text{for all }d,b.
\]
For any $a\in A_*$, choose $a_{d_j}\to a$. The uniformly bounded
family $h_{a_{d_j}}$ is normal on $E$. Every subsequential limit
has value $g_b(a)$ at each $b\in B$, and hence the identity theorem
makes the limit unique. Denote it by $\widetilde h_a$. Thus
\[
                 \widetilde h_a(b)=g_b(a)
                 \quad(a\in A_*,\ b\in B).
\]
For almost every original $a\in A$, this new representative is
$h_a$, again by the identity theorem. We have obtained an exact
separately holomorphic function on
$(D\times B)\cup(A_*\times E)$, with a common bound.

Remove the polar exceptional sets of nonregular points from
$A_*$ and $B$. The remaining sets are locally regular and retain
positive real measure. The classical cross theorem
\cite[Theorem 4]{Cross} extends this function holomorphically to
\[
 \{(z,w)\in D\times E:
       \omega_{A_*,D}(z)+\omega_{B,E}(w)<1\},
\]
where the extremal functions are unchanged by the polar removals.
Choose regular density points $a_0$ of $A$ and $b_0$ of $B$.
There is a complex product neighborhood of $(a_0,b_0)$ in this
extension domain. Let $H$ be the extension and $I_0\times J_0$
a real rectangle inside it.

For almost every $a\in A\cap I_0$, the vertical slice of $H$
equals $h_a$, so $H=f$ almost everywhere on
$(A\cap I_0)\times J_0$. For almost every $b\in J_0$, both
$H(\cdot,b)$ and the original horizontal slice representative
are real analytic on $I_0$ and agree on a set of positive measure.
They therefore agree on $I_0$. This proves $H=f$ almost everywhere
on the rectangle. Its restriction is real-valued by the same
uniqueness argument.

For a compact real-analytic target, use an analytic embedding into
some Euclidean space \cite[Theorem 3, p.~470]{Grauert}, and apply the componentwise
argument. The representative takes values in the embedded manifold
everywhere, because it does so almost everywhere and that image is
closed. Its local inverse charts make it an analytic map to the
original manifold.
\end{proof}

The conclusion is one open patch. No uniform radius in the
transverse variable is assumed, and no assertion of joint
analyticity at every point is made at this stage.

\section{Relative Levi normalization}\label{II-series:compact}\label{II-sec:compact}

\subsection{The parabolic reduction without Levi ergodicity}

Let $m=n-1\ge3$ and write the connected line parabolic as
\[
 Q^\circ=(L\ltimes U)\rtimes A_0,
 \qquad L=\SL(m,\R),\quad U\simeq\R^m,\quad A_0\simeq\R.
\]
The $L$-module $U$ is standard or dual standard, and $A_0$
centralizes $L$.

We fix the auxiliary lattice and all its measure data compatibly.
Take $\Delta=\SL(m,\Z)\subset L$ with its usual integral embedding
and take $\nu_L$ to be Haar measure. This subgroup is a lattice
\cite[Theorem 7.0.1, p.~151]{MorrisArithmetic}; equivalently, it is
discrete and its quotient admits a finite nonzero invariant regular
Borel measure \cite[Definition B.2.1 and Example B.2.5(vi), pp.~352,357]{BHV}.
A countable-slice construction gives a Borel right fundamental domain
$F_L$; every such domain has finite Haar measure
\cite[Proposition B.2.4(ii), pp.~353--354]{BHV}, and countability of
$\Delta$ and nonzero Haar measure give $\nu_L(F_L)>0$.
The group $L$ has property $(T)$ by
\cite[Theorem 1.4.15, pp.~54--55]{BHV}, and the closed-subgroup
inheritance theorem applies to the normalized invariant quotient
measure \cite[Theorem 1.7.1, pp.~65--66]{BHV}.
The preceding complexification and finite-sum argument therefore
produces one finite set $S_L\subset\Delta$ and one positive
finite-energy constant $\kappa_L$. Whenever an auxiliary lattice
package and an integral-Levi package are both used below, they are
this same tuple $(\Delta,\nu_L,F_L,S_L,\kappa_L)$.

The identical property-$(T)$ argument applies to the original lattice
$\Lambda\subset\SL(m+1,\R)$ using its positive finite Haar
covolume; its normalized quotient measure is regular. The discrete
clause of \cite[Theorem 1.3.1, p.~41]{BHV} supplies finite generators
for $\Lambda$. An ambient compact generating word model can be
constructed independently: choose a compact neighborhood
of the identity. The subgroup it generates is open and closed, hence
is the whole connected group. Its word balls, formed from this neighborhood and its inverse, are
compact, so its integer-valued word length is Borel measurable. These choices
are made once before restricting to finite-index subgroups.

Parametrize $A_0$ by $a_t$ and put
$\chi(a_t)=t$, $\chi(LU)=0$. The hyperplane case follows by
inverse transpose. Suppose $Q^\circ$ acts jointly measurably
preserving a standard probability space $(X,\mu)$.
Write $\D$ for the group of positive-derivative real analytic
diffeomorphism germs of $(\R,0)$.

All cocycles and germ fields are measurable in their Taylor
coefficients. For Lie groups, identities hold almost everywhere
for each fixed group element; a common conull set for all group
elements is not asserted. Standard measurable realizations of
invariant factors are understood modulo null sets.

\begin{theorem}[Reduction to the Levi component factor]
\label{II-loc:thm:relative}
Let $c:Q^\circ\times X\to\D$ be an actual analytic germ cocycle.
Let $p:X\to(Z,\zeta)$ realize the $L$-invariant sigma-algebra.
Then:
\begin{enumerate}
\item $L\ltimes U$ acts trivially on $Z$, and $a_t$ induces a
      probability-preserving measurable flow $T_t$ on $Z$.
\item There are measurable analytic germs $k_x\in\D$ and an
      analytic germ cocycle $\bar c$ over $T$ such that
      \begin{equation}\label{II-loc:eq:relative}
       k_{qx}\circ c(q,x)\circ k_x^{-1}
           =\bar c(\chi(q),p(x)).
      \end{equation}
      In particular the normalized cocycle is the identity on
      $L\ltimes U$.
\item If the $Q^\circ$ action is ergodic, the flow on $Z$ is
      ergodic. If $Z$ is a point, there is a single analytic vector
      field $V(z)\partial_z$, $V(0)=0$, with
      $\bar c(t)=\Phi_t^V$.
\end{enumerate}
The transfer and its inverse have positive measurable radii.
There are no moment, common-radius, or prior formal-transfer
assumptions. The analytic input is the spectral-gap realization
theorem in Section~\ref{II-sec:spectral}.
\end{theorem}

The residual equation in general is
\[
 \bar c(t+s,z)=\bar c(t,T_s z)\circ\bar c(s,z).
\]
It is a cocycle over a possibly nontrivial measurable flow, rather
than necessarily a homomorphism of $\R$. Thus removing Levi
ergodicity from analytic normalization does not remove its role
in the earlier \emph{single-flow} conclusion Section~\ref{II-sec:compact}.
We prove the theorem in the following three subsections.

\subsection{Transfers over nonergodic probability actions}

\begin{lemma}[Additive transfers without ergodicity]
\label{II-loc:lem:additive}
Let a countable property-$(T)$ group $\Lambda$ act preserving a
standard probability measure. Every finite measurable additive
cocycle $f:\Lambda\times X\to\R^D$ is a measurable coboundary,
without any moment assumption.
\end{lemma}
\begin{proof}
Use the skew-translation unitary representation on
$\mathcal H=L^2(X\times\R^D,\mu\times dv)$,
\[
 (\pi(\lambda)F)(x,v)
   =F(\lambda^{-1}x,v-f(\lambda,\lambda^{-1}x)).
\]
Normalized indicators $\xi_R$ of expanding cubes in $\R^D$,
constant in $x$, are almost invariant. For each fixed $\lambda$,
the squared displacement tends pointwise in $x$ to zero and is
bounded by two; dominated convergence proves the assertion.

Let $P$ be the orthogonal projection onto $\mathcal H^\Lambda$.
The Kazhdan inequality on its orthogonal complement implies
$\|\xi_R-P\xi_R\|_2\to0$. Choose a sequence $(R_d)_{d\ge0}$ tending to infinity and set
$F_d=P\xi_{R_d}$ and
$N_d(x)=\int_{\R^D}|F_d(x,v)|^2\,dv$.
Each $N_d$ is invariant and finite almost everywhere. The invariant
set $A=\{x:N_d(x)=0\text{ for all }d\}$ has measure zero:
\[
 \mu(A)\le\|\xi_{R_d}-F_d\|_2^2\longrightarrow0.
\]
Consequently
\[
 w_x(v)=\sum_{d\ge0}2^{-d}|F_d(x,v)|^2
\]
has strictly positive integral for almost every $x$. Since
$\|F_d\|_2\le1$, its integral over $X\times\R^D$ is at
most $\sum_{d\ge0}2^{-d}=2$. Tonelli therefore gives finite
fiber integral almost everywhere. Normalize it to a probability density. Invariance of each $F_d$
translates this density by $f(\lambda,x)$ when $x$ is replaced
by $\lambda x$. Coordinate marginal lower medians therefore give
a finite measurable vector $b_x$ satisfying
$b_{\lambda x}-b_x=f(\lambda,x)$.
\end{proof}

\begin{proposition}[Conditional-median analytic realization]
\label{II-loc:prop:conditional}
Let $\Lambda$ be a countable property-$(T)$ group acting
probability preservingly, without assuming ergodicity. Let
$c(\lambda,x)$ be a cocycle of analytic tangent-to-identity germs
in any fixed finite dimension. Suppose a measurable field $H_x$
of invertible formal germs satisfies
\[
 H_{\lambda x}=c(\lambda,x)\circ H_x.
\]
If $p:X\to Z$ is the ergodic-component factor, there is a measurable
formal tangent-to-identity field $C_z$ such that
$K_x=H_x\circ C_{p(x)}$ is analytic almost everywhere.
Moreover every such analytic cocycle has an analytic coboundary
transfer, whether or not a formal one was initially given.
\end{proposition}
\begin{proof}
The zero-dimensional case is immediate. In positive dimension, first
reduce an arbitrary invertible formal transfer to the tangent case without
changing its source parameters. Put $A_x=DH_x(0)$ and
\[
 \bar H_x=A_x^{-1}\circ H_x,\qquad
 \bar c(\lambda,x)=A_{\lambda x}^{-1}\circ c(\lambda,x)\circ A_x.
\]
The first coefficients make $A_x$ and its inverse measurable. The original
forward equation gives $\bar H_{\lambda x}=\bar c(\lambda,x)\circ\bar H_x$
and $D\bar H_x(0)=\id$. In particular the derivative of $\bar c$ is the
identity almost everywhere. If $H$ is already tangent, this step changes
nothing.

Disintegrate $\mu=\int\mu_z\,d\zeta(z)$ into ergodic invariant
probabilities. Countability of $\Lambda$ makes the cocycle and
formal-transfer identities valid on almost every component.
Fix once, in each homogeneous degree, the orthonormal Sobolev coordinates
used in Section~\ref{II-sec:spectral}. Use these same bases on every component
and at every truncation. Apply the fixed-coordinate triangular construction
to $\bar H$: at degree $d$, its chosen coordinate after right correction is
an unknown constant plus a known jointly measurable remainder $R(z,x)$.
Set the constant to $-\operatorname{med}_{\mu_z}(R(z,\cdot))$ and
continue in degree. Conditional quantiles are measurable in $z$, by the
measurable kernel integrals of the tests $R(z,x)\le q$ for rational $q$.
The compatible finite jets therefore give one coefficient-measurable
formal field $C_z$, tangent to the identity. This is normalization in the
fixed Sobolev coordinates; no identification with scalar-monomial medians
is used.

For almost every $z$, property $(T)$ gives the scalar spectral gap on
$(X,\mu_z)$. The proof in Section~\ref{II-sec:spectral} applies to this exact
same fixed-coordinate normalizer and proves convergence of
$\bar H_x\circ C_z$ for $\mu_z$-almost every $x$. Coefficient root tests
make the convergence event measurable; disintegration and $p(x)=z$
on almost every conditional fiber give a conull convergence set for
$\bar H_x\circ C_{p(x)}$. Reconstruct its coefficient-measurable analytic
germ $\bar K_x$, and set $K_x=A_x\circ\bar K_x$. Then its Taylor series
is exactly $H_x\circ C_{p(x)}$. The same factor field $C$ is invariant
under each fixed source element, so associativity and analytic coefficient
faithfulness preserve the original forward equation. Measurable
representatives of $K_x$ and its inverse have positive radii; these need
not be uniform either within a component or across components.

Finally, Lemma~\ref{II-loc:lem:additive} constructs a tangent formal transfer
successively in each homogeneous degree: after lower degrees
are removed, the first remaining coefficient is an additive
cocycle in a finite-dimensional real vector space. The induction
does not change earlier coefficients. Apply the preceding
realization to this formal transfer.
\end{proof}

This argument uses the \emph{specific} conditional-median
normalization, not an unverified measurable choice among
componentwise existence statements.

\paragraph{The original decomposition and its canonical conditional laws.}
We use the following original theorem, separately from the cohomology gluing
statement. For a locally compact second countable group $G$, a standard Borel
probability space $(X,\mu)$ and a jointly Borel quasi-invariant action,
\cite[Theorem~2.19, pp.~255--256]{FisherMorrisWhyte} supplies a conull
$G$-invariant Borel subset $D\subset X$, a standard Borel probability space
$(Z,\nu)$, a pointwise invariant Borel map $p:D\to Z$, and a Borel family
$\xi_z$ of probability measures on $D$. Every $\xi_z$ is quasi-invariant
and ergodic for the restricted original action,
$\xi_z(p^{-1}\{z\})=1$ for every $z$, and
\[
 \mu(A)=\int_Z\xi_z(D\cap A)\,d\nu(z)
 \qquad(A\subset X\text{ Borel}).
\]
Here ergodicity uses the mod-null convention: a measurable $S\subset D$
with $\xi_z(S\mathbin{\triangle}gS)=0$ for every $g$ has measure zero or
one. Quasi-invariance means that images of null measurable sets are null.
These conventions and the evaluation sigma-algebra on probability measures
are explicit in \cite[retained author version, pp.~1,3]{GreschonigSchmidt}.
The original proof constructs $Z$ as a Borel subset of the probability
measures on $X$; thus it needs no enlargement of the universe of the base
space.

We now derive the stronger form used here. Since $\nu$ is a probability,
choose $z_0\in Z$ and extend $p$ by $z_0$ off $D$, obtaining a Borel
$\varphi:X\to Z$. Conullness and invariance of $D$ give
$\varphi(gx)=\varphi(x)$ for $\mu$-almost every $x$, for each $g$.
Push $\xi_z$ along the inclusion $j:D\hookrightarrow X$ and write
$\kappa_z=j_*\xi_z$. Measurability of all evaluation maps constructs a
Markov kernel; the displayed barycenter identity gives
$\int\kappa_z\,d\nu(z)=\mu$. Moreover,
$\kappa_z\{x:\varphi(x)=z\}=1$ for every $z$.
The measure $\lambda=\nu(dz)\kappa_z(dx)$ therefore has first marginal
$\nu$, second marginal $\mu$, and is concentrated on
$z=\varphi(x)$. Comparing the maps $(z,x)\mapsto z$ and
$(z,x)\mapsto\varphi(x)$ yields the previously unasserted identity
$\nu=\varphi_*\mu$. The same graph support and barycenter show that
$\lambda$ is the graph law of $(\varphi,\mathrm{id})$ under $\mu$.
Uniqueness of conditional distributions consequently identifies
$\kappa_z$ with the actual canonical conditional kernel
$\operatorname{condDistrib}(\mathrm{id}\mid\varphi)(z)$ for
$\varphi_*\mu$-almost every $z$.

For completeness, the zero-one assertion transports along this same
inclusion: pull back an essentially invariant Borel subset of $X$ to $D$,
use equivariance of $j$, and apply the source ergodicity assertion to
$\xi_z$. The original image formulation agrees with the preimage
formulation by replacing $g$ with $g^{-1}$; null symmetric difference
is precisely almost-everywhere equality of sets. Probability converts
measure zero or one to an almost-everywhere empty or full set.
Thus the canonical conditional measures have the required invariant-set
zero-one property. For the probability-preserving action, preservation
under inverse elements gives quasi-invariance. The continuous-action
argument establishes preservation of the canonical components on one
common conull set when needed.

At finite-dimensional jet-group stages we also use the following separate
cohomology gluing statement. For locally compact second countable source and
target groups, a jointly Borel probability-preserving action on a standard
Borel space, and its genuine measurable ergodic disintegration, two jointly
Borel cocycles that are cohomologous on almost every component are globally
cohomologous by one Borel transfer
\cite[Corollary~3.6(3), p.~262]{FisherMorrisWhyte}.
The fixed-group-element almost-everywhere cocycle convention is permitted
\cite[Definition~2.27 and Remark~2.29, p.~258]{FisherMorrisWhyte}; the proof
uses joint measurable selection, not arbitrary choices of component transfers.
This theorem requires the ergodic decomposition and reconstruction for the
original action and probability. It does not provide that disintegration or
replace the conditional-median construction above.

\subsection{Preserving the full Levi identities}

\begin{lemma}\label{II-loc:lem:factors}
For a jointly measurable probability-preserving $L=\SL(m,\R)$ action
on a standard probability space, with $m\ge2$, the invariant
sigma-algebras of $L$ and $\Lambda=\SL(m,\Z)$ coincide modulo
null sets. Ergodicity is not required.
\end{lemma}
\begin{proof}
The canonical complex Koopman representation is strongly continuous.
Here this follows directly from joint measurability: for each $f,w\in L^2(X)$,
the squared orbit distance is the integral of
$|f(g^{-1}x)-w(x)|^2$, hence is measurable in $g$. Separability gives a
countable cover by orbit preimages of balls of radius $\epsilon/2$.
One such preimage has positive Haar measure; the Haar Steinhaus theorem
and the isometry property then give continuity at the identity.
Translation gives continuity everywhere. Thus no additional Koopman
continuity assumption is needed for this action.
On $L^2(X)\ominus L^2(X)^L$ the Koopman representation has no
$L$-fixed vectors. Howe--Moore in the form of
\cite[Theorem~1.1 and Definition~2.6, pp.~66, 69]{Ciobotaru} makes its coefficients
vanish at infinity. A nonzero $\Lambda$-fixed vector would have
constant positive self-coefficient along an unbounded sequence
in $\Lambda$, a contradiction. Thus the fixed subspaces coincide;
applying this to indicators proves the sigma-algebra assertion.
\end{proof}

\begin{proposition}[Full Levi analytic normalization]
\label{II-loc:prop:levi}
In Theorem~\ref{II-loc:thm:relative}, there are measurable analytic
germs $k_x\in\D$ with
$k_{lx}c(l,x)k_x^{-1}=\id$ for every fixed $l\in L$ almost everywhere.
\end{proposition}
\begin{proof}
First every jointly measurable real additive cocycle on $L$ is a
measurable coboundary. Indeed, apply Lemma~\ref{II-loc:lem:additive} on
$\Lambda$ and subtract the resulting coboundary on all of $L$.
For the residual cocycle $f$, the function
\[
 F_x(g\Lambda)=f(g^{-1},x)
\]
is well-defined for almost every $(x,g\Lambda)$ on the
finite-volume space $L/\Lambda$. The identity follows from
$f(\lambda^{-1},g^{-1}x)=0$, Fubini, and countability of $\Lambda$.
For fixed $l$,
\[
 F_{lx}(lg\Lambda)=F_x(g\Lambda)-f(l,x).
\]
The lower median $b_x$ of the Haar-pushforward distribution of
$F_x$ therefore satisfies $b_{lx}=b_x-f(l,x)$. Thus $-b_x$ is
a transfer for the residual cocycle. Vector-valued cocycles are
treated by coordinate medians.

Apply this first to $\log c(l,x)'(0)$. Positive linear changes
$D_x$ make $c_0(l,x)=D_{lx}c(l,x)D_x^{-1}$ tangent to the identity.
Apply the same additive argument to successive first remaining
coefficients to construct a formal tangent-to-identity transfer
\[
 H_{lx}=c_0(l,x)\circ H_x
\]
for the \emph{full} Lie Levi, exactly as in Section~\ref{II-sec:compact}.

Restrict to $\Lambda$ and apply Proposition~\ref{II-loc:prop:conditional}.
By Lemma~\ref{II-loc:lem:factors}, its component factor $p:X\to Z$ is
$L$-invariant. Consequently $K_x=H_x\circ C_{p(x)}$ is analytic
almost everywhere and still satisfies
\[
 K_{lx}=c_0(l,x)\circ K_x
\]
for each fixed $l$. Set $k_x=K_x^{-1}\circ D_x$.
The component-dependent normalization has preserved full Levi
equivariance, which an arbitrary independently chosen lattice
transfer would not guarantee.
\end{proof}

\begin{lemma}[A one-parameter germ homomorphism has an analytic generator]
\label{II-time:lem:one-parameter}
Every coefficient-measurable homomorphism
$\theta:(\R,+)\to\mathcal D_+^\omega$ is of the form
$\theta(s)=\Phi_s^V$ for a real analytic vector field $V$ vanishing
at zero. No common radius or differentiability in $s$ is assumed.
\end{lemma}
\begin{proof}
First derive a common domain for small times. The analytic
coefficient heights of $\theta(t)$ and its inverse are finite
measurable functions of $t$. Some common height sublevel has
positive Lebesgue measure. Its difference set contains a neighborhood
of zero. The homomorphism identity
$\theta(s-u)=\theta(s)\circ\theta(u)^{-1}$ and the coefficient
composition bound therefore give uniform geometric coefficient
bounds for all sufficiently small times. On a smaller spatial
interval the canonical analytic sums $\Theta(t,z)$ have a common
radius, are jointly measurable, and satisfy the local action law
as actual functions. The latter follows from the germ identity
and analytic continuation on the common interval.

Choose a sufficiently small $\varepsilon>0$ and put
\[
 H(z)=\frac1\varepsilon\int_0^\varepsilon\Theta(s,z)\,ds.
\]
The uniform coefficient bounds permit termwise integration and
make $H$ analytic. Moreover $H(0)=0$ and
$H'(0)=\varepsilon^{-1}\int_0^\varepsilon\theta(s)'(0)\,ds>0$.
Thus $H$ is an analytic local coordinate. For small $t,z$ the
local action law gives
\[
 H(\Theta(t,z))=
 \frac1\varepsilon\int_t^{t+\varepsilon}\Theta(u,z)\,du.
\]
For fixed $z$ the right side is continuous in $t$, by continuity
of translations of an integrable function on a bounded time
interval. Applying the local inverse of $H$ first proves time
continuity of $\Theta(t,z)$. The same identity can now be
differentiated at $t=0$, giving
\[
 \left.\frac{\partial}{\partial t}\Theta(t,z)\right|_{t=0}
 =V(z),\qquad
 V(z)=\frac{\Theta(\varepsilon,z)-z}{\varepsilon H'(z)}.
\]
This is analytic near zero and vanishes there. The local action
law transports the time-zero derivative to every sufficiently
small time, so
$\partial_t\Theta(t,z)=V(\Theta(t,z))$ on a smaller rectangle.
Uniqueness for this analytic ordinary differential equation
identifies $\theta(t)$ with the flow germ of $V$ for these times.
For any $s\in\R$, choose a positive integer $d$ with $s/d$
in this time window and use
$\theta(s)=\theta(s/d)^d$. Hence the same identification holds
for every time as an identity of germs; it does not assert a
single spatial domain valid for all real times.
\end{proof}

\subsection{The central component flow}

\begin{proof}[Proof of Theorem~\ref{II-loc:thm:relative}]
Use Proposition~\ref{II-loc:prop:levi} and write $c'$ for the normalized
cocycle. It is the identity on $L$.
Put $H=L\ltimes U$. The additive cocycle
$f_1(h,x)=\log c'(h,x)'(0)$ on $H$ vanishes on $L$.
In its skew-translation representation on $L^2(X\times\R)$,
the vector $e^{-v^2}$ is $L$-fixed. For each $u\in U$ choose
$l_d\in L$ with $l_dul_d^{-1}\to e$. The Mautner argument
makes this vector $u$-fixed, so $f_1(u,x)=0$ almost everywhere.
The unitary representation is strongly continuous because it is
measurable on a separable Hilbert space. No moment is used.

After coefficients through degree $j-1$ have vanished on $H$,
the coefficient of degree $j$ is an additive cocycle on $H$,
again zero on $L$. The same argument kills it on $U$.
Analyticity gives $c'(h,x)=\id$ for every fixed $h\in H$.

The same Mautner argument applied to the ordinary Koopman
representation says that every $L$-invariant square-integrable
function is $U$-invariant. Thus $L$ and $H$ have the same
invariant factor $Z$. Since $A_0$ centralizes $L$, it acts on
that factor; call its measurable probability-preserving flow $T_t$.
Here one may first construct a standard component realization
$e:X\to Y$. To retain the prescribed realization $p:X\to Z$,
apply measurable factorization to both invariant sigma-algebras:
there are measurable $f:Z\to Y$ and $g:Y\to Z$ with
$e=f\circ p$ and $p=g\circ e$ almost everywhere. Define
$T_t(z)=g(T^Y_t(f(z)))$. The identity
$p(a_t x)=T_t(p(x))$ follows from these two equalities and
measure preservation of $a_t$. Pulling back by $p$ proves the
probability-preserving and pairwise almost-everywhere flow laws.
The integral Levi and real Levi have the same invariant events
by the preceding Mautner argument, so this construction applies
to the stated $L$-invariant factor.

Commutation gives $c'(a_t,lx)=c'(a_t,x)$ for each fixed $t,l$.
Each Taylor coefficient of this germ and of its inverse is therefore
constant almost everywhere on almost every Levi component.
Choose its coordinatewise lower medians with respect to the component
probability. Kernel measurability makes these choices jointly measurable
in $(t,z)$; the lower median of a constant equals that constant, so they
recover the original coefficients after pullback for each fixed $t$.
On the measurable locus where the forward coefficients are symmetric,
vanish at degree zero, satisfy an integer exponential bound, and their
first derivative is inverse to the chosen inverse first derivative,
sum the series and apply the analytic inverse function theorem.
On the complement use the identity germ. This constructs a jointly
coefficient-measurable field $\bar c(t,z)$ of invertible analytic germs,
with $c'(a_t,x)=\bar c(t,p(x))$ almost everywhere for each fixed $t$.
Measure preservation and the intertwining of the factor maps descend
the cocycle identity for each fixed pair of times. If the construction
was performed on $Y$, use $\bar c(t,z)=\bar c^Y(t,f(z))$.
The equality $f(T_tz)=T^Y_t(f(z))$ holds almost everywhere;
it follows by pulling back to $X$ and using $e=f\circ p$.
Thus this transport preserves both the cocycle law and the
original normalizing germs $k_x$. Writing $q=a_t h$
and using $p(hx)=p(x)$ proves \eqref{II-loc:eq:relative}.

A $T$-invariant function on $Z$ pulls back to a $Q^\circ$-invariant
function on $X$, proving the ergodicity assertion. If $Z$ is a
point, $\bar c$ is a coefficient-measurable homomorphism from
$\R$ to $\D$. Lemma~\ref{II-time:lem:one-parameter}
identifies it with an analytic vector-field flow, as in
Section~\ref{II-sec:compact}. This last assertion uses that lemma; it does
not claim that an arbitrary single analytic germ embeds in a flow.
\end{proof}

\begin{corollary}[An atomic component factor is trivial]
\label{II-loc:cor:atomic}
If $Q^\circ$ acts ergodically and $Z$ has an atom, then $Z$ is
a point modulo null sets. In particular finitely or countably
many Levi ergodic components cause no additional case.
\end{corollary}
\begin{proof}
On a standard Borel probability space, a positive measurable atom
contains a positive singleton: normalize the restriction to that atom.
The resulting zero-one probability is a Dirac measure, and its
point has positive mass for the original probability. Thus it
suffices to consider a positive singleton.
The orbit of an atom under the measure-preserving $\R$ action
is finite, since all its atoms have the same positive mass.
The action on this finite orbit is trivial: the additive group
$\R$ has no nontrivial finite quotient. Thus the atom is invariant,
and ergodicity gives it full mass.
\end{proof}

If a transverse germ cocycle is unoriented, its sign skew product
$X\times\{\pm1\}$, with the uniform fiber probability, makes
the germs oriented. Theorem~\ref{II-loc:thm:relative} applies without
having to prove Levi ergodicity on that extension. Geometric
descent is a separate issue; no analytic covering of the manifold
is asserted by this measurable construction.

\section{Invariant probabilities and the compact derivative branch}\label{II-series:invariant-probabilities}

\begin{lemma}[Compact targets in the required dimension]
\label{II-loc:lem:compact-target}
Let $\Gamma<\SL(n,\R)$ be any lattice, $n\ge3$, and let $K$
be a compact Lie group with $\dim K<n^2-1$. Every homomorphism
$\Gamma\to K$ has finite image. In particular a compact Lie
image acting effectively by smooth isometries on a connected
$n$-manifold is finite, whether or not that manifold is compact.
\end{lemma}
\begin{proof}
Use the two original clauses of
\cite[Theorem~1.4, p.~7]{BrownCDM2018CompactTargets}. If $\Gamma$ is
nonuniform, the image is finite. Otherwise, if the image is infinite,
there is a closed subgroup $K'\subset K$ containing it whose Lie algebra
is a finite product of copies of $\mathfrak{su}(n)$. Give $K'$ its
induced closed Lie subgroup structure. If this product has zero factors,
$K'$ is zero-dimensional and therefore discrete; compactness makes it
finite, contradicting the infinite image it contains. At least one
factor occurs, and $\dim\mathfrak{su}(n)=n^2-1$. The differential of the
embedded inclusion $K'\hookrightarrow K$ is injective, so
$\dim K\ge\dim K'\ge n^2-1$, a contradiction.
For the last assertion, use the full isometry group of the given
connected Riemannian manifold with its actual path distance. It is a
Lie group in the compact-open topology
\cite[Chapter II, Theorem~1.2, p.~39]{KobayashiTransformationGroups}.
Its dimension is at most $n(n+1)/2$
\cite[Chapter II, Theorem~3.1, pp.~46--47]{KobayashiTransformationGroups};
the frame argument uses the value and orthonormal derivative of an
isometry at one point. Neither assertion requires compactness or
completeness of the manifold. The family of metric isometries is
equicontinuous, so pointwise convergence and uniform convergence on
compact sets induce the same topology. Thus the compact-open structure
applies to the pointwise topology used here.
A compact subgroup is closed in this Hausdorff group. The closed-subgroup
theorem \cite[Theorem~20.12, pp.~523--525]{LeeSmoothManifolds2e} gives its Lie
structure in the original relative topology and an embedded inclusion.
The injective differential bounds its dimension by that of the full
isometry group. Since $n(n+1)/2<n^2-1$ for $n\ge3$, the first assertion
applies. For compact manifolds, compactness of the full isometry group
is the separate assertion of
\cite[Chapter VI, Theorem~3.4(5), p.~239]{KobayashiNomizuI}.
An arbitrary compact target without this dimension restriction
is not being ruled out.
\end{proof}

\begin{theorem}[Compact derivative reduction gives actual germs]
\label{II-loc:thm:compact}
Let a countable property-$(T)$ group $\Gamma$ act analytically on
an $d$-manifold $M$ and preserve a probability $\nu$. Suppose its
derivative cocycle over $\nu$ is measurably cohomologous to an
$O(d)$-valued cocycle. Then there is a probability-preserving
extension $(S,\eta)$ of $(M,\nu)$ and measurable analytic
coordinate germs $h_s:(\R^d,0)\to(M,p(s))$ such that
\begin{equation}\label{II-loc:eq:identity-model}
 \alpha(\gamma)\circ h_s=h_{\gamma s}
\end{equation}
as germs. Their radii and inverse radii are positive and measurable.
No ergodicity, moment, or uniform derivative-growth assumption is
required. This conclusion uses the analytic-realization input
Section~\ref{II-sec:spectral} through Proposition~\ref{II-loc:prop:conditional}.
\end{theorem}
\begin{proof}
Choose measurable centered analytic charts $\phi_x$ and a
measurable linear frame $B_x$ with
\[
 f(\gamma,x)=B_{\gamma x}^{-1}\phi_{\gamma x}^{-1}
              \alpha(\gamma)\phi_x B_x,
 \qquad D_0 f(\gamma,x)=C(\gamma,x)\in O(d).
\]
On $S=M\times O(d)$ use $\eta=\nu\times m_{O(d)}$ and
\[
 \gamma(x,k)=(\gamma x,C(\gamma,x)k).
\]
This action preserves $\eta$. If $s=(x,k)$ and
$\gamma s=(\gamma x,k')$, put
\[
 \beta(\gamma,s)=(k')^{-1} f(\gamma,x) k.
\]
The germs $\beta$ form a cocycle tangent to the identity.
Proposition~\ref{II-loc:prop:conditional} supplies analytic germs $K_s$
with $K_{\gamma s}=\beta(\gamma,s)K_s$.
Then $h_s=\phi_x B_x k K_s$ satisfies
\eqref{II-loc:eq:identity-model}. All operations preserve coefficient
measurability and yield actual local diffeomorphisms on positive
measurable radii.
\end{proof}

The theorem supplies a family indexed by a probability extension.
It does not assert that these coordinates descend to one chart
per center, or that the flat metrics of overlapping charts agree.
Those stronger statements are not needed for the next consequence.

\begin{proposition}[Thickening to an open measure class]
\label{II-loc:prop:thickening}
Under Theorem~\ref{II-loc:thm:compact}, fix a smooth positive density
$m_0$ on $M$. There are an invariant open set $W\subset M$ and
an invariant probability $\widehat\nu$ such that
\[
 \nu(W)=1,\qquad \widehat\nu\sim m_0|_W.
\]
In particular $\widehat\nu$ is absolutely continuous. No smoothness
of its density, or invariant analytic metric on $W$, is asserted.
\end{proposition}
\begin{proof}
Choose analytic local diffeomorphism representatives $h_s$ on
$B_{r(s)}$ and replace $r(s)$ by $\min(1,r(s))$. Pass to a
measurable invariant conull set where all countably many germ
covariance identities hold. Transport the representatives along the
action. Their domains are centered balls and their maps agree as
germs; the analytic identity theorem makes them agree on the
connected overlaps. Thus they define a jointly measurable continued
map $h_s$ on
\[
 D_s=\bigcup_{\gamma\in\Gamma}B_{r(\gamma s)}\subset B_1.
\]
The domains are invariant and the continued maps satisfy the exact
equivariance identity throughout them. Each continued map is a local
diffeomorphism.

Let $\theta$ be the pushforward of the restriction of
$\eta\times dv$ to $\{(s,v):v\in D_s\}$ under $(s,v)\mapsto h_s(v)$.
Its mass is finite because $D_s\subset B_1$, and positive because
every $D_s$ contains a ball of positive radius. Domain invariance and
equivariance make $\theta$ invariant. Set
$\widehat\nu=\theta/\theta(M)$.

Choose a countable family of compact source buffers $K_i$ whose
interiors locally cover every open source domain, and a countable
topological basis $B_j$ of $M$. The events
\[
 E_{ij}=\{s:K_i\subset D_s,\ B_j\subset h_s(K_i)\}
\]
are measurable: compact-containment tests give the first condition,
and joint measurability and continuity on the compact buffers give
the second. Local invertibility ensures that every point of every
continued image has a basis neighborhood covered by such a buffer.
Remove all zero-measure events and take an invariant conull core $C$.
Every event containing a retained index then has positive measure,
so its covered basis neighborhood remains covered after any further
conull removal of indices. Consequently
\[
 W=\bigcup_{s\in C}h_s(D_s)
\]
is open, invariant, and unchanged by any further conull removal.
The centers show $\nu(W)=1$.

Analytic coordinate changes preserve the null sets of a positive
smooth density. Hence every $m_0$-null set has zero Lebesgue preimage
in every continued chart, and therefore zero $\theta$-mass.
Conversely, if a measurable set has zero $\theta$-mass, its Lebesgue
preimage is null for almost every index. Retaining these indices does
not change $W$. A countable subcover by continued local chart branches
and the same coordinate-null-set equivalence show that the set has
zero $m_0$-measure in $W$. Thus $\widehat\nu\sim m_0|_W$.
\end{proof}

\begin{corollary}[Compact measures cannot border a measure-free region]
\label{II-loc:cor:frontier}
Let $U\subset M$ be invariant and open and carry no nonzero
finite invariant measure. No invariant probability supported
on $\partial U$ can have compact derivative reduction.
\end{corollary}
\begin{proof}
Use the continued charts in Proposition~\ref{II-loc:prop:thickening}.
Each is open and centered on $\partial U$, so its inverse image
of $U$ is nonempty and open. Hence the pushed-forward probability
satisfies $\widehat\nu(U)>0$. Its restriction to the invariant
set $U$ contradicts the assumption.
\end{proof}

For an already continued exterior tube with an equivariant
measurable equivalence $U\simeq\R^n\setminus\{0\}$ or its
antipodal quotient, the assumption in Corollary~\ref{II-loc:cor:frontier}
follows from Lemma~\ref{II-general:lem:linear-models}. Indeed, a nonzero
finite invariant measure concentrated on $U$ can be divided by its
total mass and pushed forward under that same equivalence. This gives
an invariant probability on the corresponding punctured vector model,
contrary to the lemma. The antipodal quotient here retains the radial
variable. Consequently a compact invariant set in the frontier of such
a region cannot support an invariant probability of compact derivative type.

\subsubsection{A probability of identity-model charts forces local compactness}

\begin{lemma}[Return charts and a compact local action]
\label{II-loc:lem:return-charts}
Let a countable group act smoothly on a compact manifold $M$.
Suppose it acts ergodically preserving a probability on $S$, and
there are pointwise measurable smooth embeddings
$h_s:B_R\subset\R^{\dim M}\to M$, for one
$R>0$, satisfying $\gamma h_s=h_{\gamma s}$ on $B_R$.
There is a nonempty invariant compact set $K$ with nonempty
interior such that the closure of the action on $K$ is compact
and acts by smooth diffeomorphisms on $\operatorname{int}K$. On a connected
component of this interior, the stabilizer of that component
has compact Lie image.
\end{lemma}
\begin{proof}
Start with all smooth maps $B_R\to M$ and the weak $C^\infty$
topology. The weak finite-$r$ mapping spaces have countable bases
\cite[Chapter~2, Section~1, p.~35]{Hirsch1976}; the weak smooth topology
is the initial topology from the countably many finite-$r$ forgetful
maps \cite[p.~36]{Hirsch1976}. Countable products and subspaces therefore
give a countable base in the smooth case too. This same smooth topology
admits a complete metric
\cite[Chapter~2, Theorem~4.4 and the infinity-jet paragraph, p.~62]{Hirsch1976},
so the raw smooth mapping space is Polish. Zero- and first-jet
prolongations into their compact-open mapping spaces are continuous by
the weak jet-topology construction \cite[pp.~61--62]{Hirsch1976}.
These full jets retain the source point; forgetting it is fixed continuous
postcomposition. Since the ball is locally compact, joint evaluation of values and first
jets is continuous. Here a first jet uses the intrinsic tangent bundle:
continuity is not an assertion that an arbitrary varying choice of
coordinate derivative is continuous as a matrix-valued function.

Let $\mathcal E$ be the locus of injective maps whose differential is
invertible everywhere. Invertibility is open in the intrinsic first-jet
bundle. Quantifying this open condition over each member of a countable
compact exhaustion of $B_R$ gives a $G_\delta$ regularity locus.
The off-diagonal of $B_R\times B_R$ is also locally compact and second
countable, hence has a countable compact exhaustion. The closed diagonal
in the Hausdorff target makes inequality of the two values an open joint
condition. Compact quantification over that exhaustion gives the
$G_\delta$ injectivity locus. Thus $\mathcal E$ is a Borel subset of the
raw Polish mapping space. Equal dimensions make the differential of each
given smooth embedding invertible, so the original family lies in this locus.
The local inverse-function theorem
\cite[Theorem~4.5, p.~79]{LeeSmoothManifolds2e} gives local smooth inverses for
its maps. Injectivity makes these inverses agree on overlaps, so they
patch to a smooth inverse on the open image. Conversely, a smooth map
with such an inverse has invertible differential. Accordingly this
carrier consists of the same actual open smooth embeddings, with
values outside the fixed ball discarded.

We first prove open total evaluation in the weak subspace topology on
$\mathcal E$. In a fixed target chart, the actual family and its fiber
derivative are jointly continuous. Apply the inverse-function theorem
with arbitrary topological parameters
\cite[Definition~2.2, p.~419, and Lemma~5.10, pp.~455--456]{Swaminathan2021}.
Its original assertion has Euclidean fibers. Translating the chosen
source and target points to the origins, applying it to this coordinate
family, and composing with the same product chart gives a local inverse
for $(h,v)\mapsto(h,h(v))$. The parameter is never differentiated.
These local neighborhoods prove that total evaluation is open;
continuity of value evaluation and fiberwise injectivity make it an
open embedding.

Choose a finer Polish topology on the Borel locus $\mathcal E$
\cite[Theorem~3.2.4, pp.~92--93]{SrivastavaBorel}. This step can be made
without changing the maps or the evaluation map. Indeed, writing
$\mathcal F(h,v)=(h,h(v))$, for a new parameter-open set $V$ and a
ball-open set $D$ one has
\[
 \mathcal F(V\times D)
   =(V\times M)\cap\mathcal F(\mathcal E\times D).
\]
The right-hand side remains open in the refined product topology;
continuity also persists. Total evaluation is therefore an open
embedding for the chosen Polish topology. The target incidence relation
is open, and the inverse and transition evaluations are continuous.
Pointwise measurability of the original family implies its Borel
measurability for this very topology: evaluation on a countable dense
sequence is a continuous injective map into a countable product, so
Lusin--Souslin gives a Borel inverse on its image. No measurability of
parameter derivatives is required.

Tightness now gives a positive-measure set $C\subset S$ whose chart
parameters lie in a compact family sufficiently close to one fixed
embedding. Choose fixed open coordinate balls $D_0,D_1,D_2$ with
$\overline D_0\subset D_1$ and $\overline D_1\subset D_2$
so that every transition $h_t h_s^{-1}$, $s,t\in C$, is defined
on a neighborhood of $\overline D_2$ and maps $D_0$ into $D_1$.
The common domains and these inclusions follow from openness of
target incidence and continuity of transition evaluation, after
shrinking the parameter neighborhood.

Ergodicity gives a finite set $F\subset\Gamma$ with
$\eta(\bigcup_{f\in F}fC)>1-\eta(C)$. For every $\gamma$,
$\gamma C$ meets some $fC$ in positive measure. Write
$\gamma=fr$, where $rC\cap C$ has positive measure. Choosing
$s,rs\in C$ in the conull equivariance set gives
$r=h_{rs}h_s^{-1}$ on $D_2$. Consequently
\[
 \gamma D_0\subset\bigcup_{f\in F}fD_1.
\]
For each $\gamma$, the finite return decomposition writes
$\gamma=ar$ with $a\in F$ and
$r=h_t h_s^{-1}$ on $D_2$ for parameters in the fixed compact
family. Put
$K=\overline{\bigcup_\gamma\gamma D_0}$ and
$O=\bigcup_{b\in F}bD_2$. Then $K\subset O$ with compact
containment and $K$ is invariant. On $bD_2$, write
$\gamma=(\gamma b)b^{-1}$; its restriction is represented by
\[
 a\circ h_t\circ h_s^{-1}\circ b^{-1},\qquad a,b\in F.
\]
These finite translated compact families give compactness of the
paired compact-open closure of the maps on $K$ and their inverses.
The pairs in the closure remain inverse homeomorphisms.

More precisely, each closure element, on a neighborhood in $K$ of
each point, equals one actual translated transition from such a
compact parameter family. To see this, choose the translated domain,
pass to a subnet with the finite translate fixed and the parameter
pair convergent, and use continuous evaluation on the common domain.
Consequently the closure element has a smooth local-diffeomorphism
extension near that point. Applying the same fact to its inverse
shows that both preserve $\operatorname{int}K$. The closure therefore
acts by smooth diffeomorphisms there; no assertion of smooth convergence of a sequence
of transformations is needed. A connected component of the interior
has open, hence finite-index, stabilizer. The compact effective image
of that stabilizer acts continuously by smooth diffeomorphisms on the
connected component. The differentiable transformation-group theorem
\cite[Corollary~1, p.~161]{GeorgeMichael2007} makes this image a Lie group.
Here the classical input is the universal statement for a locally compact
Hausdorff group acting continuously and effectively by $C^1$
diffeomorphisms on a connected paracompact finite-dimensional manifold;
its conclusion supplies a Lie atlas compatible with the given group topology.
For the literal subgroup of paired compact-open maps, joint continuity
is continuous evaluation, and effectivity follows because equality on
every point is equality of the forward maps and hence of the paired
homeomorphisms. Compactness supplies local compactness of the subgroup.
The component is locally compact and second countable, hence
$\sigma$-compact and paracompact. Smoothness of each transformation
implies $C^1$ regularity, including its inverse, which is another group
element. Thus the universal theorem applies to this very subgroup,
with its given topology and action.

The same local inverse function theorem gives an inverse on the image
of an embedding of the entire ball $B_R$ when source and target have
the same dimension. The immersion differential is then invertible at
every point. Injectivity makes the local inverses agree on overlaps,
and their domains cover the open image. Thus the original embedding
is a diffeomorphism from $B_R$ onto its image. This argument also applies
in dimension zero and for a noncompact target.
\end{proof}

\begin{theorem}[Description of the compact derivative branch]
\label{II-loc:thm:compact-finite}
Let $M$ be a compact connected real analytic $n$-manifold, and
let $\Gamma<\SL(n,\R)$ be any lattice, $n\ge3$.
If an invariant probability has compact derivative reduction,
then, using the analytic-realization input of
Theorem~\ref{II-loc:thm:compact}, the action has finite image.
No invariant smooth volume is needed.
\end{theorem}
\begin{proof}
Start with the identity-model germs of Theorem~\ref{II-loc:thm:compact}
on their invariant conull covariance set. The measurable representative
and inverse radii, chart centers, and Taylor coefficients let us choose
a positive-measure set $C$ with one common positive radius and compact
chart parameters. Continue the representatives along the action as in
Proposition~\ref{II-loc:prop:thickening}; analytic uniqueness makes the
actual return transitions agree on their common domains.

The finite-return construction does not require ergodicity.
Choose a maximal finite family $\{fC:f\in F\}$ of pairwise
essentially disjoint translates. Such a family exists because
$|F|\eta(C)\le\eta(S)=1$ bounds its cardinality. By maximality,
every translate $\gamma C$ intersects some $fC$, $f\in F$,
in positive measure. Thus $\gamma=fr$ with
$\eta(rC\cap C)>0$, giving the same actual return transitions
as in Lemma~\ref{II-loc:lem:return-charts}.
The resulting translated compact chart families
give the same invariant compact carrier with nonempty interior and
compact paired closure. Their coordinate derivatives also form compact
families on smaller common domains. Thus the closure has continuous
tangent action and locally uniform bounds on derivatives of every order.

Average a smooth background metric on the interior over its compact
closure group. The continuous tangent action and the locally uniform
coordinate derivative bounds justify differentiating the Haar integral,
so the averaged metric is smooth and invariant. Choose one connected
component of the interior; its stabilizer is open and has finite index.
Its effective compact image lies in the Riemannian isometry group of
that component. The ordinary isometry-group theorem and orthonormal
frame dimension bound make this a compact Lie group of dimension at
most $n(n+1)/2$. Therefore
Lemma~\ref{II-loc:lem:compact-target} makes it finite for every
lattice in $\SL(n,\R)$. A further finite-index subgroup
acts identically on that open set. The analytic identity theorem
on connected $M$ makes it act identically on $M$.
\end{proof}

\subsection{Coordinate germs for invariant probabilities}
We write $\sigma(g)=g$ for the standard representation and
$\sigma(g)=g^{-T}$ for its dual. An analytic coordinate germ at
$x\in M$ is an invertible real-analytic map
$(\R^n,0)\to(M,x)$, considered up to restriction to a neighborhood
of zero. A family of such germs is measurable if its Taylor
coefficients in a finite analytic atlas are measurable.

\begin{theorem}\label{II-I:thm:invariant-germs}
Let $n\ge3$, let $\Gamma<\SL(n,\R)$ be a lattice, and let
$\alpha:\Gamma\to\operatorname{Diff}^{\omega}(M)$ act on a compact
connected real-analytic $n$-manifold. Suppose $\mu$ is an ergodic
$\alpha$-invariant probability. Then one of the following holds.
\begin{enumerate}
\item The derivative cocycle has compact algebraic type, and
$\alpha(\Gamma)$ is finite.
\item The derivative cocycle has standard or dual type. There is
a measurable family of analytic coordinate germs $\Phi_x$, on
an extension of $(M,\mu)$ by a finite probability space if needed,
satisfying
\[
 \alpha(\gamma)\circ\Phi_x=\Phi_{\gamma x}\circ\sigma(\gamma)
 \qquad(\gamma\in\Gamma)
\]
on a single invariant set of full measure.
\end{enumerate}
The conclusion holds for uniform and nonuniform lattices, and
$\mu$ need not belong to the smooth measure class.
\end{theorem}
\begin{proof}
The analytic germ cocycle superrigidity theorem of Part I
(Theorem~\ref{I-thm:compact-actions}) and the representation
calculation in Section~\ref{I-sec:standard-dual} give the
indicated derivative types.

The source specializations in this Part~I application retain the original
invariant probability. The lattice cocycle theorem is
\cite[Definition~1.3 and Theorem~1.5, p.~193 of the published version]{FisherMargulis};
its nonuniform-lattice integrability argument is in Sections~3.6--3.7.
We use the actual real and complex matrix point groups throughout.
The determinant--exponential identity
\cite[Proposition~0.11(e), pp.~7--8]{KnappBeyondDigital} identifies their tangent
spaces with the trace-zero matrices. The prescribed involutions are
$\theta_R(X)=-X^{\mathsf T}$ and $\theta_C(X)=-X^*$; their Cartan
positivity is the matrix statement in
\cite[Chapter VI, Section~2, pp.~354--355]{KnappBeyondDigital}.
The complex tangent remains semisimple on restriction to the reals
\cite[equation~(1.61), p.~58]{KnappBeyondDigital}.
The closed-matrix-group charts
\cite[Theorem~0.15, pp.~12--13]{KnappBeyondDigital} retain these actual tangent
inclusions and matrix exponentials. Apply
\cite[Theorem~6.31, p.~362]{KnappBeyondDigital} to both prescribed involutions;
the required finite centers are supplied by
\cite[Proposition~7.9, pp.~440--441]{KnappBeyondDigital}.
The complex compact subgroup $U$ is a subgroup of the original complex
point group. Its polar product with a real vector space is homeomorphic
to that group, so simple connectivity of the latter gives simple
connectivity of $U$.

The real compact subgroup $K_R$ is connected and every element is an
exponential \cite[Corollary~4.48, p.~259]{KnappBeyondDigital}. Compatibility of
$\theta_R$ and $\theta_C$ on real matrices, together with the actual
matrix exponential, puts its literal complex image in $U$.
Thus the common compact lift is this inclusion and the real compact
cover is the identity. The fixed-subgroup tangent calculations give
the actual skew-symmetric and skew-Hermitian trace-zero spaces, so the
subgroup charts retain their literal matrix differentials.
The elementary symmetric-matrix brackets identify the generated compact
dual with $\mathfrak{su}(n)$, sending the imaginary direction of $X$
to $iX$ in the original complex matrix space.
A nonzero diagonal difference is a negative Cartan direction and generates
the simple source as a Lie ideal. Its compact adjoint orbit spans the
negative Cartan space, whose brackets span the positive space: apply
the canonical ideal statement
\cite[Chapter VI, Problem~23, p.~430, solution p.~735]{KnappBeyondDigital} to an
adjoint-invariant subspace $W$ of the negative space, obtaining
$[W,\mathfrak p]+W$. These constructions fix the exponentials,
representations, compact lift and tangent inclusions on the same source.
For continuous representations, matrix differentiation uses
\cite[Theorem~2.21, pp.~45--47]{Hall2003}, applied to the actual closed
matrix groups after the finite-dimensional target dictionary; empty source
or target indices are handled directly. The resulting differential
intertwines the original matrix exponentials with their original time parameter.
Integration on the same simply connected compact dual uses
\cite[Theorem~20.19, pp.~530--531]{LeeSmoothManifolds2e}, with the actual tangent
Lie algebras and the closed-matrix charts above. The checked tangent and
exponential dictionaries return the integrated representation to the same
compact dual. Haar averaging there constructs the Cartan-compatible metric
for this differential.
Real Lie complete reducibility
\cite[Lemma~B.1 and its proof, pp.~659--660]{KnappBeyondDigital} gives invariant
complements, hence commuting projections. Differentiation and generation by
the actual exponential curves promote these projections to group-equivariant
ones, in each finite jet degree. Haar right invariance on the literal special
linear group is transported through its actual polynomial point-group
isomorphism, giving right invariance of the same pushed-forward Haar measure.
For the auxiliary standard special linear groups the standing polynomial
model is itself constructed from the determinant--exponential tangent
identity, complex $\mathfrak{sl}$ simplicity, real-form simplicity and real
connectedness \cite[Chapter II, Section~1, Example~1, pp.~124--126;
proof of Proposition~7.1, p.~434; Proposition~1.145, p.~118]{KnappBeyondDigital}.
Simple connectivity of the literal complex special linear group is the
assertion in \cite[proof of Proposition~24.1, p.~205, using
Theorems~13.5--13.6, pp.~90--91]{Bump2013}; its point-group homeomorphism
transfers this to the same polynomial model. The explicit two-dimensional
commuting split pair supplies the model's rank condition when $n\ge3$.
The literal and polynomial real point-group isomorphism then transports
the derived Weyl projections back to the same auxiliary Levi representation. The same construction applies to the Levi coefficient
space itself and to the tensor representations in the normal collars, as well
as to every jet degree; all these applications use the existing bound $m\ge3$.
For algebraicity we use the rationality assertion
\cite[Lemma~3.4, p.~199, with the standing hypotheses of Section~3.1,
pp.~198--199]{FisherMargulis}. The native algebraic group is the affine
group of these same defining equations, and its canonical real-point
topology agrees with the original matrix topology. The standard special
linear group has no proper nontrivial smooth connected normal subgroup,
and its diagonal split torus has dimension $n-1$
\cite[Example $(A_n)$, pp.~457--459]{Milne}. Thus every noncommutative
almost-simple normal factor is the whole group and has numerical real
rank at least two when $n\ge3$. The faithful matrix theorem
\cite[Theorem~4.9 and Corollary~4.10, p.~87]{Milne} supplies the
finite-dimensional numerical-rank comparison. A nontrivial representation
forces a nonempty factor family; composing its split-torus inclusion with
the original factor inclusion verifies total rank at least two. Each
factor also has positive rank, as required by the cited rationality lemma.
The zero-dimensional target and trivial representation use the trivial
algebraic homomorphism directly. For the other representations the lemma
gives an algebraic homomorphism agreeing on every original real point.
Pulling back its matrix entries to the same coordinate algebra and choosing
polynomial representatives gives exactly the original representation in
its specified basis. No change of source or represented operators occurs.
The actual direct-product decomposition uses \cite[Theorem~24.3, p.~516]{Milne}.
Here all factors are subgroups of the same original algebraic source.
The coordinate comparison uses smoothness implying reducedness
\cite[Lemma~33.25.4, tag 056T]{StacksSmoothReduced}, and the canonical
differential of a finite product is the product of the differentials
\cite[Theorem~10.14, p.~190]{Milne}. Semisimplicity is compared with the
Lie algebra by \cite[Theorem~23.68, pp.~507--508]{Milne}.
For the complex semisimple source, the algebraic and topological notions
of simple connectivity agree
\cite[Section~4, footnote~1, p.~149]{RoquesSinger2022}; this supplies the
geometric simple-connectivity hypothesis of the factorization theorem.
Each actual almost-simple factor has simple Lie algebra in characteristic
zero \cite[Lemma~6.2, p.~35]{CadorelDengYamanoiII}.
The product projections are retractions of the original source, so its
real connectedness and complex simple connectivity pass to these same
chosen factors. Their inclusions and multiplication maps, rather than
independently chosen isomorphic copies, are used in the representation.
For the type-A automorphism step, use the original matrix theorem
\cite[Chapter~IX, Section~5, Theorem~5; verified translation p.~306,
proof pp.~305--306]{Jacobson1964}. A complex-linear automorphism of
$\mathfrak{sl}_n(\C)$ is conjugation, or conjugation followed by negative
ordinary transpose; at size two only the inner alternative is needed.
Finite reindexing and the inverse-matrix convention are checked before this
is used as the intertwiner statement for the same representation.
The literal real source chart is obtained from the polynomial matrix chart
by the identity on tangent matrices, the literal real point-group
isomorphism, and the Euclidean matrix-to-operator algebra isomorphism.
Finite-dimensional continuity handles their different matrix norms;
the matrix entries, strict chart derivative and exponential time remain fixed.
For spectral transfer, the actual product multiplication is a topological
group isomorphism, so each original factor projection is continuous and
surjective. Property $(T)$ of the standard special linear group
\cite[Theorem~1.4.15, pp.~54--55]{BHV} passes through these very projections.
Apply \cite[Theorem~1(i), pp.~306--307]{Nevo1998} to each actual
simple noncompact factor with finite center. Its positive integer denominator
is fixed by the factor, uniformly over unitary representations without fixed
vectors and over Borel probabilities. The real and nonseparable Hilbert-space
transports preserve that denominator and the same regular averaging operator.

The regular-representation estimate is the source-only norm assertion in
\cite[proof of Theorem~5(4), p.~311; Section~2.2, p.~312]{Nevo1998}:
for the actual maximal compact $K$, a maximal abelian subspace of
$\mathfrak p$, and any nonzero $H$ in it, the norms of
$m_K*\delta_{\exp(kH)}*m_K$ decay exponentially for integer $k\ge1$.
The infinite-measure action hypotheses in the maximal-inequality conclusion
are not used for this separate norm assertion in its proof.
Here the original simple source forces the tangent of every nontrivial factor
to be the whole source tangent, and the product isomorphism makes its literal
inclusion onto. The given Cartan structure and direction therefore transport
to the same factor. In its native chart, exponential generation and polar
uniqueness identify the constructed compact subgroup with the given one.
Its maximal compactness follows from
\cite[Theorem~6.31(g), p.~362]{KnappBeyondDigital}.
Injective smooth homomorphisms are immersions
\cite[proof of Proposition~7.17, p.~157]{LeeSmoothManifolds2e}; together with the
actual induced subgroup topology this gives the required embedded inclusions.
The original normal-subgroup ideal statement
\cite[Chapter~I, Section~10, p.~77]{KnappBeyondDigital} and connected-subgroup
containment \cite[Section~2.6, Theorem~5, p.~27]{OnishchikVinberg1990}
verify the absence of compact simple normal factors. Noncompactness itself
follows from the model's genuine real split pair. Extend the unchanged negative
Cartan vector to a maximal abelian subspace, and apply the cited norm bound.
The canonical exponential line, positive integer time, normalized compact
Haar probability and factor Haar measure remain the original ones.
This supplies both spectral fields without assuming a selected comparison
or spherical-decay packet; no smooth measure class is required.
For standard and dual type, apply
Corollary~\ref{I-power:cor:manifold}; the compact centralizer is
contained in $\{I,-I\}$ and is removed on its finite skew-product
extension. Theorem~\ref{II-loc:thm:compact-finite} treats compact type.
\end{proof}

The compact case requires more than a conclusion about Lyapunov
exponents. We construct analytic identity-model germs on a compact
probability extension. Returns to a bounded family of these germs
give a compact group of local transformations. A dimension bound
for compact targets and the analytic identity theorem then give
finite image.

\begin{proposition}[Atomic invariant measures give linearizable fixed points]
\label{II-loc:prop:atomic-invariant-caps}
For any lattice $\Gamma<\SL(n,\R)$, let an ergodic invariant
probability be atomic and have standard
or dual derivative type. Its finite orbit has actual analytic
linearizing neighborhoods after finite-index restriction.
\end{proposition}
\begin{proof}
The invariant probability is uniform on a finite orbit. Pass
to a point stabilizer and remove the finite scalar character
in its derivative representation. Compatible formal reduction
and Corollary~\ref{I-power:thm:main}
give an analytic linearizing chart at that point. In this
chart the action agrees with the standard or dual linear
action as germs. On a ball, its common domain with the
inverse image under any fixed linear group element is convex.
The analytic identity theorem extends the germ identity to
this common domain. Thus this is an actual analytic
linearizing chart with the required local action identities.
\end{proof}

\section{Analytic projective hypersurfaces}\label{II-series:embedding}
\subsection{Detected projective hypersurfaces are analytic without volume}
\label{II-sec:arithmetic-boundary-scope}
Theorem~\ref{II-general:thm:analytic-embedding} upgrades smooth
equivariant projective embeddings for every lattice. We first
explain why the finite-regularity conclusions of a description
theorem supply the required smoothness.

\begin{lemma}[Compatibility of the regularity upgrades]
\label{II-loc:lem:all-regularities}
Let an analytic action of a lattice in $\SL(n,\R)$, $n\ge3$,
have a fixed quasi-invariant probability $\nu$. Suppose that
for arbitrarily large integers $r$ there is a finite invariant
union $Y_r$ of embedded $C^r$ projective covering hypersurfaces,
with $\nu(Y_r)=1$ and $C^r$ standard or dual equivariant
identifications for their finite-index stabilizers.
Then $\nu$ is carried by a finite invariant union of analytic
projective covering hypersurfaces with analytic identifications.
\end{lemma}
\begin{proof}
Choose one of the finite unions $Y_{r_0}$ and discard its
zero-mass leaves. The support of $\nu$ is closed and invariant
by quasi-invariance and is conull on the compact manifold.
Every retained leaf meets that support. Minimality for its
finite-index stabilizer therefore puts the whole leaf in the
support. Conversely, the retained finite union is closed and
has full measure, so it contains the support. Thus the retained
union equals the support and is invariant; no permutation of a
previously chosen list of labels is needed. For each retained
leaf $N$ and every $r$, the equality $\nu(Y_r)=1$ implies that
$N$ meets some leaf $N_r$ of $Y_r$. Their intersection is closed
and invariant under the intersection of their finite-index
stabilizers. This subgroup is again a lattice and acts minimally
on each leaf. Hence $N=N_r$ as subsets of $M$.

Fix a parametrization of $N$ from the first collection.
Comparing it with the $C^r$ parametrization produces a
continuous equivariant comparison of projective spaces or
their sphere covers. The topology of $N$ determines which
cover occurs. Lemma~\ref{II-general:lem:angular} rules out a
change from standard to dual and makes the comparison the
identity, or the antipodal map on a sphere. Thus the fixed
parametrization is $C^r$ for arbitrarily large $r$, and so is
$C^\infty$. Apply Theorem~\ref{II-general:thm:analytic-embedding}
to the finite-index stabilizer of $N$.
\end{proof}

\begin{corollary}[Analytic hypersurfaces in the expanding branch]
\label{II-loc:cor:positive}
Let $\Gamma<\SL(n,\R)$ be any lattice, $n\ge4$, and let
$\mu$ be a Haar-projecting ergodic minimal-parabolic invariant
probability in the suspension of a compact analytic $n$-manifold.
Suppose its full stabilizer is a line parabolic and its fiber
weights are
\[
 \{c(\varepsilon_j-\varepsilon_1):2\le j\le n\}
       \cup\{q\varepsilon_1\},\qquad c>0,\quad q>0,
\]
with each missing-root coarse space one-dimensional. The dual
statement holds for a hyperplane parabolic. Then $\mu$ is
supported on the suspension of finitely many analytic projective
covering hypersurfaces.
\end{corollary}
\begin{proof}
If $\varepsilon_j(H)-\varepsilon_1(H)\le0$ for every $j$, then
$\varepsilon_1(H)>0$ unless $H=0$, because the coordinates sum
to zero. Thus the remaining exponent is positive on the entire
closed contraction cone minus its origin. These are the spectral
conditions of \cite[Theorem 4.7]{BRHW}. Apply it at every finite
regularity $r\ge3$ above its critical threshold. The projective
boundary has finite fundamental group, so each application gives
a finite invariant union of compact embedded projective covering
leaves with $C^r$ identifications. The connected covering spaces
are spherical or projective: use the classification of coverings
and $\pi_1(\mathbb{RP}^{n-1})\simeq\mathbb Z/2$ in
\cite[Proposition~1.14, p.~35, Propositions~1.33--1.34,
pp.~61--62, Theorem~1.38, p.~67, and Example~1.43, p.~74]{Hatcher}.
First restrict to the kernel of the finite permutation action on
the leaves. On each fixed leaf, uniqueness of covering lifts leaves
one global antipodal sign for each group element. These signs form
a character; its kernel has finite index and gives a strictly
equivariant parametrization of the same leaf. The induced covering
atlas is analytic, while the original inclusion into $M$ is still
only $C^r$ at this stage. Composing the covering parametrization
with that inclusion preserves the $C^r$ regularity and its image.
To use a single measure in all
applications, choose a Borel section of $G\to G/\Gamma$ and
let $\nu_0$ be the $M$-marginal of $\mu$ in this trivialization.
Enumerate $\Gamma=\{\gamma_j:j\ge1\}$ and set
$\nu=\sum_{j\ge1}2^{-j}\alpha(\gamma_j)_*\nu_0$.
This is a quasi-invariant probability. Each of the invariant
leaf unions has full $\nu$-mass, because its suspension has
full $\mu$-mass. Lemma~\ref{II-loc:lem:all-regularities}
and Theorem~\ref{II-general:thm:analytic-embedding} upgrade one
such full-measure union and its identifications to analytic.
This proves the assertion for every lattice.

\end{proof}

This proves analytic hypersurface description for $q>0$ for
every lattice without first constructing a measurable transverse
cocycle or verifying Levi ergodicity of $\mu$.
For $n\ge4$, the full-collar Theorems~\ref{II-suspension:thm:negative-collar},
\ref{II-positive:thm:collar}, and \ref{II-neutral:thm:full-collar} apply
according to each hypersurface's actual normal character. They use
its canonical homogeneous measure. Their neighborhoods exclude this
branch from an uncovered compact invariant remainder.

\section{Equivariant comparisons of vector models}\label{II-series:comparisons}

We return to an arbitrary lattice $\Gamma<\SL(n,\R)$.

\begin{lemma}[Equivariant analytic comparisons]\label{II-exterior:lem:comparison}
Let $n\ge3$, let $O\subset\R^n\setminus\{0\}$ be a nonempty
$\Gamma$-invariant open set, and let $h:O\to\R^n\setminus\{0\}$
be analytic.
If $h(\gamma z)=\gamma h(z)$, then on every connected component
of $O$ one has $h(z)=a z$ for a nonzero constant $a$.
There is no such map with $h(\gamma z)=\gamma^{-T}h(z)$.
No integrability or moment bound on $h$ is assumed.
The same conclusions hold for target
$(\R^n\setminus\{0\})/\{\pm I\}$, with equality interpreted
modulo sign and the scalar chosen positive.
\end{lemma}
\begin{proof}
Let $G=\SL(n,\R)$ and $H=\operatorname{Stab}_G(e_1)$.
For either output representation $\rho$, define on the right-$H$
invariant finite-measure subset
\[
 Y_O=\{\Gamma g\in\Gamma\backslash G:ge_1\in O\}
\]
the measurable function
$f(\Gamma g)=\rho(g)^{-1}h(ge_1)$.
Equivariance makes it well-defined, and
$f(yk)=\rho(k)^{-1}f(y)$ for $k\in H$.

The group $H=\SL(n-1,\R)\ltimes(\R^{n-1})^*$ is generated
by finitely many one-parameter unipotent subgroups. For a time-one
element $u$ of any generator, Poincare recurrence on $Y_O$ forces
$\rho(u)f(y)=f(y)$ almost everywhere. Indeed, a nonconstant
unipotent orbit $\rho(u)^{-d}f(y)$ grows polynomially in norm,
contradicting infinitely many returns to a bounded level set of
$f$. Use the countably many such level sets to cover $Y_O$.
The finite generating family therefore gives $f(y)\in\operatorname{Fix}(\rho(H))$.

For the dual standard representation this fixed space is zero,
contradicting the nonzero target. For the standard representation
it is $\R e_1$. Disintegration along $g\mapsto ge_1$ gives
$h(z)\in\R z$ almost everywhere on $O$. Analyticity makes this
true everywhere. Thus $h(z)=c(z)z$, where $c$ is a nonzero
analytic $\Gamma$-invariant scalar function.

Lemma~\ref{II-general:lem:linear-models} makes the invariant open set
$O$ conull and the invariant scalar $c$ constant almost everywhere.
Analyticity then makes it constant everywhere on $O$.

For the quotient target use the same argument with $f$ valued
in the antipodal quotient. Norm is well-defined there. A nonfixed
unipotent orbit still has unbounded norm, so recurrence forces
its class to be fixed. A unipotent cannot send a nonzero vector
to its negative; hence its representatives are fixed as well.
In the standard case write $h(z)=[c(z)z]$ with
$c(z)=\|h(z)\|/\|z\|>0$, and repeat the scalar argument.
\end{proof}

\section{The contracting transverse branch}\label{II-series:contracting}\label{II-sec:negative}
The suspension lemmas below retain their displayed measure,
regularity, and integrability hypotheses for arbitrary lattices.
The analytic hypersurface conclusion is separated from the
full collar conclusion, which is established for every lattice
in Section~\ref{II-sec:suspension-negative-collar}.

\subsection{Negative transverse exponents and finite fibers}

We separate a consequence valid in the actual suspension from the
analytic construction criterion in the next section. Write
$\pi:M^\alpha\to G/\Gamma$ for the suspension projection, and use
the line-parabolic notation in Section~\ref{II-sec:compact}.

\begin{lemma}[Contraction makes the fiber disintegration atomic]
\label{II-loc:lem:atomic-fibers}
Let $\mu$ be a Haar-projecting invariant probability for a subgroup
$J<G$, and suppose an element $a\in J$ has all fiberwise Lyapunov
exponents negative. Assume the usual fiberwise stable-manifold
hypotheses, which hold for the suspension of a smooth compact-manifold
action; the induced derivative and inverse logarithms are
integrable by Lemmas~\ref{II-power:lem:section} and
\ref{II-power:lem:word-germs}. Then the conditional
probabilities $\mu_y$, $y\in G/\Gamma$, are purely atomic almost
everywhere. If $\mu$ is $J$-ergodic, there is an integer $k\ge1$
such that $\mu_y$ is uniform on exactly $k$ points almost everywhere.
\end{lemma}
\begin{proof}
For the selected return $a$, use the original carrier section law $\nu$.
It is already a probability: the original suspension probability gives full
mass to the invariant carrier and the section embedding preserves that law.
Write $\kappa(g,x)$ for the original centered section-return germ in the
same carrier atlas. The point $x$ lies in the carrier section, with its
original section action. The maps $T_kx=a^k x$ and
$A_k(x)=D_0\kappa(a^k,x)$ form the original strict
integer action and derivative cocycle. Put
\[
 H_k(x)=\log\max\{1,\|A_k(x)\|,\|A_k(x)^{-1}\|\}.
\]
The given one-step moment and submultiplicativity give
$H_n(x)\le\sum_{j=0}^{n-1}H_1(T_jx)$ for $n\ge0$.
Also $H_{-n}(x)=H_n(T_{-n}x)$. Invariance therefore gives every integer
moment. A fixed linear identification of the fiber with Euclidean space
changes these logarithmic bounds by at most one fixed constant.
Apply \cite[Theorem~2.4 with $d=1$, pp.~462--463, including the nonergodic
paragraph]{BRHWSmoothErgodic2023}; here $L^{1,1}=L^1$.
We complete this same $\nu$ for the source theorem.
The integer splitting and angle estimate give finite growth rates for every
nonzero vector along the forward ray. Fixed coordinates and equality of the
original and completed almost-everywhere filters return these limits to the
same original derivative powers and norm. The construction uses the actual integer return and a fixed fiber
coordinate change; it applies to the original nonergodic law.
In positive fiber dimension, choose a fixed finite initial basis in the
original fiber coordinates. Its finite vector rates are negative by hypothesis. Their maximum gives
the negative operator rate used in the contraction argument.

Consider the relatively independent product
$\lambda=\int\mu_y\otimes\mu_y\,dm_{G/\Gamma}(y)$.
It is $a$-invariant. For almost every center $x$, the fiberwise
stable manifold contains a neighborhood $D_x$ of $x$ in its fiber.
For $z\in D_x$, the fiber distance between $a^d x$ and $a^d z$
tends to zero. Poincare recurrence in the relative product implies
that, for $\lambda$-almost every such pair, $z=x$: a distinct pair
whose distance tends to zero cannot return infinitely often to
the measurable set where that distance exceeds a fixed positive
number. Apply this for a countable collection of distance thresholds.

Thus $\mu_y(D_x\setminus\{x\})=0$ for $\mu_y$-almost every $x$.
Such $x$ belongs to the support of $\mu_y$, so $\mu_y(D_x)>0$.
Consequently $\mu_y(\{x\})>0$ almost everywhere, proving pure
atomicity. The atom-mass function $x\mapsto\mu_{\pi(x)}(\{x\})$
is $J$-invariant. In the ergodic case it is a positive constant,
necessarily $1/k$, and the conditional probability has $k$ atoms.

The required stable neighborhoods and tempered fiber geometry
are the suspension tools in \cite[Sections 3.2 and 3.5]{BRHW}.
The finite-fiber consequence is also used in the proof of its
Corollary 4.5; the argument above does not require equality of
$\dim M$ and $\dim G/Q$.
\end{proof}

\begin{proposition}[The contracting branch has no moving Levi factor]
\label{II-loc:prop:negative-factor}
Suppose an ergodic minimal-parabolic suspension measure projecting to Haar
has full stabilizer the line parabolic and fiber weights
\[
 \{c(\varepsilon_j-\varepsilon_1):2\le j\le n\}
       \cup\{q\varepsilon_1\},\qquad c>0,\quad q<0.
\]
Then its fiber disintegration is uniform on finitely many atoms,
and $c=1$. On each ergodic component for $Q^\circ$, the Levi
$L=\SL(n-1,\R)$ acts ergodically. In particular the factor $Z$
of Theorem~\ref{II-loc:thm:relative} is a point on each such component.
\end{proposition}
\begin{proof}
Take
\[
 a_t=\operatorname{diag}(e^t,e^{-t/m},\ldots,e^{-t/m}),\quad m=n-1.
\]
For $t>0$ its fiber exponents are $-cnt/m$ with multiplicity $m$
and $qt$ with multiplicity one, all negative. Apply
Lemma~\ref{II-loc:lem:atomic-fibers}. For the root normalization,
\cite[Proposition~5.7, p.~41]{BRHW} gives tameness for each resonant
root $\beta$ of the full $P$-ergodic law. Use the original immersed
coarse leaves \cite[Section~3.2, p.~960]{BRHWFactors}, defined
componentwise in the nonergodic case
\cite[Section~4.2, p.~471]{BRHWSmoothErgodic2023}.
A basis containing $-H$, with $\beta(H)<0$, puts these leaves inside
the forward $H$-stable leaves. Ergodic averages of compactly supported continuous functions
give stable-saturated representatives of the finite
$A$-component events, preserving their normalized laws and Haar bases.
Their label is constant on whole $\beta$-leaves and hence measurable
modulo every subordinate partition. Its conditional kernels therefore
agree with those of each normalized component on that component.
Consequently the component families satisfy, after gluing, all four
properties of \cite[Proposition~3.11, p.~25]{BRHW}: intrinsic Radon
regularity, projective covariance, inner leafwise coherence and every
subordinate conditional formula. Uniqueness identifies this
glue with the full-law family up to positive finite normalization.
The same normalized tame neighborhoods thus pass to each component.
For each missing root, some component remains noninvariant, by the
barycenter identity. Analyticity supplies a finite order above its
finite fiber coarse critical regularity (value one for an empty class).
\cite[Proposition~4.4, p.~28]{BRHW} then gives the exact fiber exponent
$\beta$. The same-law spectrum comparison of
Section~\ref{II-sec:suspension} returns it to the original coordinate
list; distinct missing-root rays and the distinguished remaining slot
force $c=1$.

There are finitely many $Q^\circ$-ergodic components since
$Q/Q^\circ$ is finite. Each projects to Haar: that measure is
ergodic for $Q^\circ$, and its ergodic decomposition under a
factor map has Haar marginals. Fix one component. Its conditional
measures are again uniform on a constant finite number $k$ of
points. If $E$ is $L$-invariant, the function $y\mapsto\mu_y(E)$
is $L$-invariant. Moore ergodicity on $G/\Gamma$ makes it constant,
with value in $\{0,1/k,\ldots,1\}$. Thus every positive-measure
$L$-invariant set has mass at least $1/k$, and the $L$-invariant
factor has at most $k$ atoms.

The Mautner argument of this section makes every $L$-invariant
function $U$-invariant. Hence $Q^\circ$ acts on this finite factor
through $A_0\simeq\R$. A measurable action of the connected group
$\R$ on a finite probability space is trivial. Ergodicity for
$Q^\circ$ therefore makes the factor a point.
\end{proof}

This settles the component-factor issue in the negative branch.
Finite atomic disintegration is a measurable assertion; it does
not itself produce an embedded projective hypersurface.

\subsection{Constructing a transverse curve from finitely many Levi jets}

Here is an analytic construction that uses the special representation
in dimension $n=m+1$. It is stated with its finite-order hypothesis
explicit, so that it can be applied after a jet reduction without
assuming the desired transverse curves in advance.

Put $m=n-1\ge3$ and $L=\SL(m,\R)$. Let
$Q^\circ=(L\ltimes U)\rtimes A_0$ act jointly measurably on a
standard Borel probability space $(X,\mu)$, preserving its probability,
with ergodic $L$ restriction. Put
\[
 E=\R^m\oplus\R,\qquad \rho(l)(u,z)=(lu,z),\qquad
 U=(\R^m)^*.
\]
The action of $L$ on $U$ is the dual standard one. Let
$F:Q^\circ\times X\to\operatorname{Diff}^\omega(E,0)$ be a
jointly coefficient-measurable cocycle of actual full-dimensional
germs. In a suspension these are the germs of the original fiber
maps in measurable centered analytic charts, not germs of
transverse maps already constructed.

\begin{theorem}[Contracting finite-jet construction]
\label{II-loc:thm:contracting-axis}
Suppose that for some $N\ge1$ there are measurable invertible
$N$-jets $h_x$ with
\begin{equation}\label{II-loc:eq:pure-levi-jet}
 h_{lx}^{-1}\,j^NF(l,x)\,h_x=j^N\rho(l)
 \qquad(l\in L).
\end{equation}
Write $T x=a_1x$. Assume the derivative of
$h_{Tx}^{-1}j^NF(a_1,x)h_x$ is
\[
 D(u,z)=(\lambda u,bz),\qquad 0<\lambda,b<1,
\]
and put $\lambda_+=\max(\lambda,b)$ and
$\lambda_-=\min(\lambda,b)$. Require
\begin{equation}\label{II-loc:eq:finite-threshold}
                 \lambda_+^{N+1}<\lambda_-.
\end{equation}
Finally assume that $F(a_1,x)$ and its inverse have holomorphic
representatives with radii and sup-norm bounds whose negative,
respectively positive, logarithms are integrable.

Then there are measurable actual analytic coordinate germs $H_x$
such that
\[
 H_{lx}^{-1}F(l,x)H_x=\rho(l),\qquad
 H_{a_tx}^{-1}F(a_t,x)H_x=(\lambda^t u,b^t z).
\]
The embedded curve germs
\begin{equation}\label{II-loc:eq:actual-axis}
                 i_x(z)=H_x(0,z)
\end{equation}
are $Q^\circ$-equivariant, and their induced germs are
\[
 F(l,x)i_x=i_{lx},\qquad F(v,x)i_x=i_{vx}\ (v\in U),\qquad
 F(a_t,x)i_x(z)=i_{a_tx}(b^t z).
\]
All identities hold almost everywhere for each fixed group element.
The curves and their local parametrizations have positive measurable
radii. No convergence assumption on an infinite formal transfer
is needed.
\end{theorem}

\begin{proof}
The jointly measurable, probability-preserving Levi action on the
standard probability space has a strongly continuous canonical Koopman
representation. We use this consequence of the original action when
applying the invariant-factor and Mautner arguments below.

\emph{Step 1: the finite return jet is constant and triangular.}
Because $a_t$ centralizes $L$, its normalized $N$-jet $p_t(x)$
satisfies
\[
 p_t(lx)=\rho(l)p_t(x)\rho(l)^{-1}.
\]
Each homogeneous coefficient is an equivariant measurable map
into a finite-dimensional $L$-module. Its distribution is invariant
under every unipotent element. Recurrence forces its values into
the fixed space of each such element: a nonfixed orbit of a
unipotent linear map is a nonconstant polynomial sequence and
escapes every bounded set. A finite set of elementary unipotents
has common fixed space equal to the $L$-fixed space. The coefficient
is then $L$-invariant and therefore constant by ergodicity.

The centralizer of $\rho(L)$ in formal maps consists of
\begin{equation}\label{II-loc:eq:triangular-centralizer}
                 (u,z)\longmapsto(A(z)u,B(z)).
\end{equation}
Indeed there are no nonconstant scalar invariants on the standard
module, and an equivariant map from that module to itself is a
scalar multiple of the identity; apply these observations in
each homogeneous degree in $u$. Thus $p_1$ has this form through
order $N$, with $A(0)=\lambda$, $B'(0)=b$.

A constant centralizer jet linearizes $p_1$ through order $N$.
First linearize $B$ formally through that order; the denominators
are $b^j-b$, $j\ge2$. Then change $u$ to $C(z)u$ to make $A$
constant; the denominators are nonzero multiples of $b^j-1$,
$j\ge1$. Replace $h_x$ by its right composition with this jet.
Equation~\eqref{II-loc:eq:pure-levi-jet} remains true, and now
\[
                h_{Tx}^{-1}j^NF(a_1,x)h_x=j^ND.
\]

\emph{Step 2: finite coordinate changes are tempered.}
Realize $h_x$ by its polynomial representative. In a faithful
finite-dimensional jet representation, the last identity bounds
the increments of
\[
 w(x)=\log\max(1,\|j^Nh_x\|,\|(j^Nh_x)^{-1}\|)
\]
by an integrable function of $j^NF(a_1,x)$ and its inverse, plus
a constant depending on $D$. Hence $w\circ T-w\in L^1$.
Birkhoff's theorem and convergence in probability of
$(w(T^d x)-w(x))/d$ to zero give $w(T^d x)/|d|\to0$ in both
time directions. This argument does not require $w\in L^1$.
Finite-dimensional polynomial and inverse-function estimates
therefore give tempered domains and bounds for $h_x$ and its
analytic inverse. The assumed logarithmic integrability gives
the same property for the original analytic return germs.

Consequently
\[
 f_x=h_{Tx}^{-1}F(a_1,x)h_x
       =D+O(\|w\|^{N+1})
\]
has holomorphic representatives on radii $r(x)>0$ with
\[
 \|f_x(w)-Dw\|\le C(x)\|w\|^{N+1},
\]
where $r^{-1}$ and $C\ge1$ are tempered. Polynomial representatives
are used only for a finite jet; no infinite formal substitution
has been asserted to converge.

\emph{Step 3: an analytic limit.}
Choose $\theta\in(\lambda_+,1)$ and $\epsilon>0$ so that
\[
 e^\epsilon\theta<1,\qquad
              e^\epsilon\theta^{N+1}<\lambda_-.
\]
Tempering gives positive measurable $r_\epsilon(x)$ and finite
$C_\epsilon(x)$ such that, for $d\ge0$,
\[
 r(T^d x)\ge r_\epsilon(x)e^{-\epsilon d},\qquad
 C(T^d x)\le C_\epsilon(x)e^{\epsilon d}.
\]
On a sufficiently small ball at $x$, induction gives
$\|f_x^{(d)}(w)\|\le\theta^d\|w\|$, with every iterate in its
holomorphic domain. Here
$f_x^{(d)}=f_{T^{d-1}x}\cdots f_x$.
The successive differences in
\[
                    K_x(w)=\lim_{d\to\infty}D^{-d}f_x^{(d)}(w)
\]
are bounded by
\[
 \lambda_-^{-1}C_\epsilon(x)
 \left(\frac{e^\epsilon\theta^{N+1}}{\lambda_-}\right)^d
 \|w\|^{N+1}.
\]
The limit is holomorphic, measurable in $x$, and
$K_x=\id+O(\|w\|^{N+1})$. It satisfies $K_{Tx}f_x=DK_x$.
Set $H_x=h_xK_x^{-1}$. The analytic inverse function theorem
supplies its local inverse and measurable positive radii.

\emph{Step 4: the full Levi identity survives.}
Let $g_l(x)=H_{lx}^{-1}F(l,x)H_x$. Its $N$-jet is $\rho(l)$,
and the cocycle identity gives
\[
                  g_l(Tx)=Dg_l(x)D^{-1}.
\]
If all preceding nonlinear coefficients vanish, its degree-$j$
coefficient, $j>N$, satisfies this same linear covariance.
Every eigenvalue of conjugation by $D$ on that coefficient space
has modulus at least $\lambda_-/\lambda_+^j>1$.
A finite measurable random vector cannot have a stationary
distribution under a uniformly expanding linear map unless it
vanishes. Degree induction proves $g_l(x)=\rho(l)$ as actual
germs. No moments of these remaining coefficients are used.

For every fixed $t$, the normalized germ of $a_t$ is constant
and in the centralizer by the recurrence argument in Step 1,
now in every degree. It commutes with $D$. Formula
\eqref{II-loc:eq:triangular-centralizer} and $0<b<1$ show that it is
linear: $A(bz)=A(z)$ and $B(bz)=bB(z)$ eliminate all its
nonlinear coefficients. The resulting measurable one-parameter
homomorphism has the asserted positive scalars $\lambda^t,b^t$.

\emph{Step 5: the radical fixes the axis.}
Fix $0\ne v\in U$ and let $L_v$ be its stabilizer in $L$.
The normalized germ $g_v(x)$ is equivariant for linear conjugation
by $\rho(L_v)$. This stabilizer is generated by unipotent
one-parameter groups. Recurrence as above makes every coefficient
fixed by $L_v$. Moreover any nontrivial unipotent element of
$L_v$ acts ergodically on $X$, by Howe--Moore applied to the
ergodic $L$ action. Thus every coefficient is constant, so
$g_v$ is a constant analytic germ commuting with $\rho(L_v)$.

For the dual standard module $U$ and $m\ge3$,
\[
                  E^{L_v}=\{0\}\oplus\R.
\]
For instance, if $v=e_1^*$, its stabilizer on $\R^m$ consists
of matrices with first row $e_1^*$; the lower $\SL(m-1)$ block
and the first-column transvections have no common nonzero fixed
vector. Hence $g_v$ preserves the axis. Since $L$ is transitive
on $U\setminus\{0\}$ and fixes that axis pointwise, all
$g_v$ with $v\ne0$ have the same restriction $\psi$ to it.
The cocycle identity and constancy give
$g_{2v}=g_vg_v$, so $\psi=\psi^2$ and $\psi=\id$.
This proves all the claimed transverse identities.
\end{proof}

\begin{corollary}[A pure linear Levi reduction suffices]
\label{II-loc:cor:pure-linear-levi}
In the setting of Theorem~\ref{II-loc:thm:contracting-axis}, suppose
the full-dimensional derivative cocycle on $L$ is measurably
cohomologous to $\rho(l)$, with no compact error, and every
finite-jet cocycle on $L$ satisfies the integrability hypotheses
of Fisher--Margulis. Assume that the central element $a_1$ has
strictly negative fiber exponents and the analytic integrability
in that theorem holds. Then the transverse curve germs
\eqref{II-loc:eq:actual-axis} exist and have the stated parabolic
equivariance.
\end{corollary}
\begin{proof}
Apply \cite[Theorem 1.4]{FisherMargulis} inside each algebraic jet group.
Simultaneously linearize its semisimple homomorphism and compact
centralizer, and strip the linear part of the transfer as in
Section~\ref{I-sec:formal}. This reduces each finite jet to the
\emph{original} derivative cocycle. Composing with the given
linear frame therefore gives \eqref{II-loc:eq:pure-levi-jet} at every
finite order, with the same linear part.

Already at order one, commutation with $L$ and the recurrence
argument in Step 1 make the normalized central flow a constant
block-scalar homomorphism. Its time-one scalars $\lambda,b$ are
positive. The frame is tempered for $a_1$ by the integrable-increment
argument in Step 2, so these scalars have the original Lyapunov
exponents. They are strictly less than one. Choose a finite $N$
satisfying \eqref{II-loc:eq:finite-threshold} and apply the theorem.
\end{proof}

\begin{remark}[Application and the finite compact correction]
\label{II-loc:rem:uniform-FM}
In the negative suspension branch,
Proposition~\ref{II-loc:prop:negative-factor} supplies Levi ergodicity,
and the expected central diagonal scalars are
$\lambda=e^{-n/m}$ and $b=e^q$. Thus a finite $N$ always satisfies
\eqref{II-loc:eq:finite-threshold}.
For the ambient Fisher--Margulis integrability condition, fix one
measurable near-minimal suspension lift $s(y)$ for $\ell_G$, and put
$h(y)=\ell_G(\operatorname{pr}_G s(y))$.
Lemma~\ref{II-power:lem:section} and the Haar projection of the original
suspension probability give $h\in L^1$.
Near-minimality gives, for the actual return label,
\[
 \ell_G(\kappa(g,y))\le
      2\ell_G(g)+2h(y)+1.
\]
The undistortion bound in that lemma and boundedness of $\ell_G$ on
compact sets therefore give one integrable return-word majorant for
all $g$ in each compact set. In the same lift and the original finite
analytic atlas, the logarithmic norm of each faithful finite jet is
bounded by a constant times one plus this return word length.
The jet cocycle is jointly Borel. We justify completion measurability
of its actual compact supremum before using the majorant. More generally,
let $X,Y$ be standard Borel spaces, let $F:X\times Y\to\mathbb R$
be Borel, and let $K\subset Y$ be Borel. For the original finite measure
$\mu$ on $X$, put
\[
 S(x)=\sup_{y\in K}\operatorname{ofReal}(F(x,y))\in[0,\infty],
 \qquad \sup\varnothing=0,
\]
where $\operatorname{ofReal}(t)=\max(t,0)$. Give $X$ a compatible Polish
topology without changing its measurable sets or $\mu$. For every
$a\in[0,\infty]$, the strict upper threshold is the actual projection
\[
 \{x:a<S(x)\}
   =\pi_X\{(x,y):y\in K,\ a<\operatorname{ofReal}(F(x,y))\}.
\]
The relation on the right is Borel, so its projection is analytic.
The native analytic-set representation is either the empty set or a
continuous Baire-space image. In the latter case,
\cite[Proposition~4.1.1(iii)$\Rightarrow$(i) and the analytic-set definition,
p.~128]{SrivastavaBorel} represents this same set as a Borel
self-product projection; the empty case is handled directly.
Then \cite[Theorem~4.3.1, p.~141]{SrivastavaBorel} makes that same
projection measurable in the completion of the original sigma-finite
Borel measure $\mu$ (finiteness implies sigma-finiteness); see the
completion definition on p.~103.
The order-measurability criterion applied to all these upper thresholds
therefore gives $\mu$-completion measurability, or equivalently an
almost-everywhere Borel representative, of the actual $S$.

For the actual compact set of group elements, its Borelness follows from
the Hausdorff group topology. Apply the preceding result to the actual
positive logarithmic jet norm $\max\{0,\log\|\,\cdot\,\|\}$.
The same integrable majorant then makes its
supremum integrable. This verifies the compact-uniform condition of
Fisher--Margulis Definition~1.3. The same lift is used throughout this
argument.

Analyticity on compact $M$ bounds the logarithmic loss of a
holomorphic domain for a word by a constant times its word length:
compose the finitely many generator representatives on successively
smaller balls. The integrable word length of the suspension return
cocycle, from the Siegel-domain estimates underlying
\cite[Proposition 3.2]{BRHW}, consequently gives the analytic
integrability assumed above.

The pure finite-jet identity \eqref{II-loc:eq:pure-levi-jet} is an explicit
hypothesis of Theorem~\ref{II-loc:thm:contracting-axis}. The actual derivative
reduction can contain a finite compact centralizing error.
Section~\ref{II-loc:sec:signs} retains that error throughout the construction,
so no extension of the full parabolic action to a finite frame cover
is needed. Corollary~\ref{II-loc:cor:negative-curves} applies the resulting
criterion to the negative suspension branch.

Even when the criterion applies, it constructs actual transverse
curve germs at measure-typical centers. Compatibility with
projective charts at different centers, and hence an open collar
description theorem, is a further geometric step, carried out
in Section~\ref{II-sec:suspension-negative-collar}. When $q=0$,
$b=1$ and the contraction estimate fails. None of this argument
constructs the neutral transverse curve.
\end{remark}

\subsection{Retaining the finite compact correction}\label{II-loc:sec:signs}

For the standard-plus-trivial representation on $E=\R^m\oplus\R$,
every compact subgroup of its linear centralizer is contained in
\[
 \mathcal K=\{\operatorname{diag}(\varepsilon I_m,\delta):
                       \varepsilon,\delta\in\{\pm1\}\}.
\]
The same statement holds for the dual-plus-trivial representation.
We retain these signs instead of removing them on a Levi-only
extension. Equip homogeneous coefficient spaces with the Hilbert
norm induced from the standard Euclidean norm on $E$. Independent
left and right compositions by elements of $\mathcal K$ are
orthogonal on these spaces and commute with conjugation by $\rho(L)$.

\begin{lemma}[Recurrence with finite signs]\label{II-loc:lem:signed-recurrence}
Let $L$ act preserving $(X,\mu)$. Suppose a measurable homogeneous
polynomial field $P_x:E\to E$ obeys
\[
 P_{lx}=C_1(l,x)\,\rho(l)P_x\rho(l)^{-1}\,C_2(l,x),
 \qquad C_1(l,x),C_2(l,x)\in\mathcal K,
\]
where the displayed right products denote precomposition.
Then $P_x$ is fixed by conjugation by $\rho(L)$ almost everywhere.
If $L$ acts ergodically, $\|P_x\|$ is essentially constant.
The fixed-space conclusion also holds with $L$ replaced by a
subgroup generated by finitely many unipotent one-parameter groups.
\end{lemma}
\begin{proof}
Fix a unipotent $s$ in such a generating set. For every positive
integer $d$ the covariance and orthogonality give
\[
             \|P_{s^d x}\|=\|\rho(s)^dP_x\rho(s)^{-d}\|.
\]
If $P_x$ is not fixed by the induced unipotent operator, the
right side tends to infinity: its vector argument is a nonconstant
polynomial in $d$. Poincare recurrence to the sets
$\{\|P_x\|\le r\}$, $r\in\mathbb N$, rules this out.
Intersect the finitely many resulting conull sets. Being fixed by
the time-one element of a unipotent group is equivalent to being
fixed by the group. Elementary unipotent groups generate $L$.
Finally the original covariance, now with trivial $L$ conjugation,
makes the norm $L$-invariant. Ergodicity proves the last assertion
about its size.
\end{proof}

\begin{lemma}[Bounded centralizer jets and their naturality]
\label{II-loc:lem:bounded-centralizer}
Let $T$ preserve a probability space. Suppose $p_x$ is a measurable
field of centralizer $N$-jets of the form
\eqref{II-loc:eq:triangular-centralizer}, with essentially bounded
coefficients and derivative
\[
 D_x(u,z)=(a_xu,b_xz),\qquad |a_x|=\lambda>0,
             \quad |b_x|=b\in(0,1).
\]
There is a unique essentially bounded field of centralizer
$N$-jets $c_x$, tangent to the identity, satisfying
\begin{equation}\label{II-loc:eq:bounded-normalizer}
                         p_xc_x=c_{Tx}D_x.
\end{equation}
This construction is natural under measurable block-sign changes
of coordinates.

More precisely, if $L$ commutes with $T$,
$R(l,x)=\rho(l)C(l,x)$ with $C(l,x)\in\mathcal K$, and
\begin{equation}\label{II-loc:eq:return-covariance}
 p_{lx}=R(l,Tx)p_xR(l,x)^{-1},
\end{equation}
then
\begin{equation}\label{II-loc:eq:normalizer-covariance}
                   c_{lx}=R(l,x)c_xR(l,x)^{-1}.
\end{equation}
\end{lemma}
\begin{proof}
Work successively in total degrees $j=2,\ldots,N$. The degree-$j$
part of a centralizer map is
\[
                    P_x(u,z)=(\alpha_xuz^{j-1},\beta_xz^j).
\]
After the preceding degrees have been fixed, equation
\eqref{II-loc:eq:bounded-normalizer} in this degree is
\begin{equation}\label{II-loc:eq:bounded-coefficient}
             P_x=b_x^{j-1}P_{Tx}+r_j(x),
\end{equation}
where $r_j$ is bounded. Indeed the linear operator on the unknown
coefficient is $P\mapsto D_x^{-1}PD_x=b_x^{j-1}P$;
the remaining defect is a polynomial in finitely many previously
bounded coefficients, divided only by the nonzero constants
in absolute value $\lambda$ and $b$. The forward series solves
\eqref{II-loc:eq:bounded-coefficient}, with norm at most
$\|r_j\|_\infty/(1-b^{j-1})$. Iterating the homogeneous equation
shows uniqueness among bounded solutions. This proves existence
and uniqueness through degree $N$.

For a sign change $J_x\in\mathcal K$, the transformed data are
$p'_x=J_{Tx}p_xJ_x^{-1}$ and $D'_x=J_{Tx}D_xJ_x^{-1}$.
Equation~\eqref{II-loc:eq:bounded-normalizer} and uniqueness give
$c'_x=J_xc_xJ_x^{-1}$. To obtain
\eqref{II-loc:eq:normalizer-covariance}, fix $l$. Differentiating
\eqref{II-loc:eq:return-covariance} gives its counterpart for $D_x$.
Consequently $R(l,x)^{-1}c_{lx}R(l,x)$ solves the original
equation. It is a bounded centralizer field, since conjugation
by $\rho(l)$ is trivial there and the remaining signs are
orthogonal. Uniqueness proves the assertion.
\end{proof}

\begin{theorem}[Contracting construction with the compact error]
\label{II-loc:thm:signed-contracting}
Keep the parabolic, representation, probability action, and actual
full-dimensional germ cocycle in Section~\ref{II-sec:negative}; in particular $L$ is
ergodic. Suppose, for some $N$, there are measurable invertible
$N$-jets $h_x$ such that
\begin{equation}\label{II-loc:eq:signed-levi-jet}
 h_{lx}^{-1}j^NF(l,x)h_x=R(l,x),\qquad
 R(l,x)=\rho(l)C(l,x),\quad C(l,x)\in\mathcal K.
\end{equation}
Assume that the absolute values of the two blocks of the
normalized $a_1$ derivative are $\lambda,b\in(0,1)$, with
$\lambda_+^{N+1}<\lambda_-$. Assume the analytic return
integrability of Theorem~\ref{II-loc:thm:contracting-axis}.

Then there are actual analytic coordinates $H_x$ satisfying
\[
 H_{lx}^{-1}F(l,x)H_x=R(l,x),
\]
and measurable signs $\varepsilon_t(x),\delta_t(x)$ such that
\[
 H_{a_tx}^{-1}F(a_t,x)H_x(u,z)
           =(\varepsilon_t(x)\lambda^t u,\delta_t(x)b^t z).
\]
The analytic curve germs $i_x(z)=H_x(0,z)$ are defined on the
original probability space and are $Q^\circ$-equivariant. There
is a sign cocycle $\sigma:Q^\circ\times X\to\{\pm1\}$ with
\begin{equation}\label{II-loc:eq:signed-transverse}
 F(g,x)i_x(z)=i_{gx}\bigl(\sigma(g,x)b^{\chi(g)}z\bigr).
\end{equation}
Their radii and inverse-parametrization radii are positive and
measurable. The assertion does not require a parabolic lift of
the compact Levi frame extension.
\end{theorem}
\begin{proof}
\emph{Step 1: the central finite jets are bounded.}
Put $T=a_1$. The normalized $a_t$ jet $p_t(x)$ obeys
\eqref{II-loc:eq:return-covariance} with $T$ replaced by $a_t$.
Apply Lemma~\ref{II-loc:lem:signed-recurrence} in each degree. Every
coefficient is in the $L$-centralizer and has constant norm,
so $p_t$ has the form \eqref{II-loc:eq:triangular-centralizer} with
bounded coefficients for each fixed $t$. Its derivative blocks
have constant nonzero absolute values. The cocycle law makes
these absolute values measurable homomorphisms of $\R$ into
$\R_{>0}$, hence they are $\lambda^t,b^t$.

At time one write
\[
 D_x=\operatorname{diag}(\varepsilon_1(x)\lambda I_m,
                                      \delta_1(x)b).
\]
Apply Lemma~\ref{II-loc:lem:bounded-centralizer} to $p_1$. Replacing
$h_x$ by $h_xc_x$ makes the time-one $N$-jet exactly $D_x$.
Equation~\eqref{II-loc:eq:normalizer-covariance} shows that this
replacement preserves \eqref{II-loc:eq:signed-levi-jet}, with precisely
the same sign cocycle. This is the compatibility that a Levi-only
finite cover would not automatically supply.

\emph{Step 2: signs do not change analytic convergence.}
Use the polynomial representative of the new $h_x$. Its finite
jet and inverse have tempered norms, by the integrable-increment
argument in Step 2 of Theorem~\ref{II-loc:thm:contracting-axis}:
the jets of $D_x$ and $D_x^{-1}$ have uniform bounds.
Thus $f_x=h_{Tx}^{-1}F(a_1,x)h_x=D_x+O(|w|^{N+1})$
has tempered holomorphic bounds and radii. Put
\[
 D_x^{(d)}=D_{T^{d-1}x}\cdots D_x,\qquad
 K_x=\lim_{d\to\infty}(D_x^{(d)})^{-1}f_x^{(d)}.
\]
The norms of these products and their inverses are exactly
$\lambda_+^d$ and $\lambda_-^{-d}$. The estimate in Step 3
of that theorem applies unchanged and gives a holomorphic
limit $K_x=\id+O(|w|^{N+1})$, with
$K_{Tx}f_x=D_xK_x$. Set $H_x=h_xK_x^{-1}$.

\emph{Step 3: the full Levi and central identities.}
For fixed $l\in L$ put $g_l(x)=H_{lx}^{-1}F(l,x)H_x$.
Its $N$-jet is $R(l,x)$, and
\[
                   g_l(Tx)=D_{lx}g_l(x)D_x^{-1}.
\]
On degree-$j$ coefficients the operator on the right has
minimum singular value at least $\lambda_-/\lambda_+^j>1$
for $j>N$. Indeed its variable signs are orthogonal; its
positive diagonal blocks are the same at every center.
Recurrence forces every such coefficient to vanish, proving
the full analytic Levi identity.

For each fixed $t$, Lemma~\ref{II-loc:lem:signed-recurrence} now
places every coefficient of the normalized $a_t$ germ in
the centralizer. Commutation with the time-one map gives
the same equation with $D_{a_tx}$ on the left and $D_x^{-1}$
on the right. On the two-dimensional space of centralizer
degree-$j$ coefficients, $j\ge2$, its expansion in norm is
exactly $b^{1-j}>1$. Recurrence again forces vanishing.
The linear part already has the asserted absolute values.

\emph{Step 4: the radical preserves the analytic axis.}
Fix $0\ne v\in U$ and let $L_v<L$ be its stabilizer.
The normalized germ $g_v(x)=H_{vx}^{-1}F(v,x)H_x$ satisfies
\[
             g_v(lx)=R(l,vx)g_v(x)R(l,x)^{-1}
                         \qquad(l\in L_v).
\]
Apply the subgroup version of Lemma~\ref{II-loc:lem:signed-recurrence}
to each coefficient. The group $L_v$ is generated by finitely
many unipotent one-parameter groups, so $g_v(x)$ commutes with
$\rho(L_v)$ as an actual germ. Since
$E^{L_v}=\{0\}\oplus\R$, it preserves the axis. The Levi
and central maps preserve it as well. Hence $i_x$ defines
actual equivariant curve germs on $X$, with a one-dimensional
analytic germ cocycle $\psi(g,x)$.

\emph{Step 5: orient only the curve parametrizations.}
Define $\sigma(g,x)=\operatorname{sgn}\psi(g,x)'(0)$.
This is now a cocycle for the \emph{full} parabolic action,
because the curves have already been constructed. On the
probability extension
\[
 \widehat X=X\times\{\pm1\},\qquad
 g(x,\eta)=(gx,\sigma(g,x)\eta),
\]
use the parametrization $\widehat i_{(x,\eta)}(z)=i_x(\eta z)$.
The induced germ cocycle is oriented, equals the identity on
$L$, and equals $z\mapsto b^t z$ on $a_t$.

The Mautner argument in Section~\ref{II-sec:compact} applies without ergodicity
on this extension. First apply it to the additive log-derivative
cocycle on $L\ltimes U$, which vanishes on $L$. It vanishes
on $U$. Successively apply it to each first nonlinear Taylor
coefficient; these too vanish on $U$. Analyticity gives the
identity germ on $U$. This part uses the Mautner additive
argument alone, not the spectral-gap analytic-realization input.
Thus the oriented induced cocycle is $z\mapsto b^{\chi(g)}z$.
Undoing the sign parametrizations gives
\eqref{II-loc:eq:signed-transverse} on $X$ itself. No geometric cover
of the manifold has been asserted or needed.
\end{proof}

\begin{corollary}[Actual transverse curves in the negative branch]
\label{II-loc:cor:negative-curves}
Let $n\ge4$. In the compact analytic suspension setting of
Proposition~\ref{II-loc:prop:negative-factor}, on each $Q^\circ$-ergodic
component there are measurable embedded real analytic curve germs
in the actual suspension fibers, tangent to the transverse
$q\varepsilon_1$ line, with
\[
 \alpha(g)i_x(z)=i_{gx}\bigl(\sigma(g,x)e^{q\chi(g)}z\bigr),
                   \qquad g\in Q^\circ.
\]
Here $q<0$, and $\alpha(g)$ denotes the actual suspension fiber
map. The hyperplane branch has the inverse-transpose version.
No invariant smooth volume is assumed.
\end{corollary}
\begin{proof}
Proposition~\ref{II-loc:prop:negative-factor} gives Levi ergodicity
and $c=1$. Apply Fisher--Margulis \cite[Theorem 1.4]{FisherMargulis} to
the derivative on $L$ in the original finite atlas and the near-minimal
lift of Remark~\ref{II-loc:rem:uniform-FM}. Subsequent linear frame changes
are applied to this normal form. Its Cartan weights, obtained by restricting
the displayed suspension weights to $L$, identify the noncompact
representation as standard plus trivial (dual plus trivial for
the hyperplane version). The compact centralizer is contained
in $\mathcal K$.

For clarity, measurable changes of frame do not lose these
weights. For each Cartan element the reduced linear cocycle
is a fixed representation matrix times a commuting orthogonal
matrix. Its norm and inverse norm are bounded. Integrability
of the original derivative and the integrable-increment argument
make the frame tempered along that element. Its exponents
are therefore preserved. The weight multiset determines the
finite-dimensional semisimple representation.
Here one can work directly over the real split Lie algebra. Apply the
highest-weight existence and uniqueness theorems
\cite[Chapter~VII, Section~2, Theorem~2, p.~237, and Section~3,
Theorem~3, p.~238]{Jacobson1964}; their standing hypothesis is a split
semisimple Lie algebra over a field of characteristic zero
\cite[Chapter~VII, pp.~228--229]{Jacobson1964}.
First use real Weyl complete reducibility
\cite[Lemma~B.1 and its proof, pp.~659--660]{KnappBeyondDigital} to decompose the \emph{given}
representation into actual invariant real subspaces. The corresponding
projections sum to the identity and their dimensions sum to the original
real dimension. Differentiate each restricted representation using
\cite[Theorem~2.21, pp.~45--47]{Hall2003} on the same matrix exponential
curves and the auxiliary special-linear charts constructed above.
A Lie-invariant subspace has a commuting Lie projector by the same real
Weyl theorem; exponential generation makes that projector group-equivariant.
Thus an irreducible group summand remains irreducible for the differential.
A nonzero highest vector in such a summand is cyclic: a stable subspace
containing it cannot be zero, so irreducibility makes it the whole summand.
This verifies the generator condition in Jacobson's definition of an
$e$-extreme module \cite[Chapter~VII, Section~2, p.~234]{Jacobson1964}.

There is a convention to retain here: Jacobson uses right modules and the
lower simple-root generators. For a left differential $d\pi$, set
$v\mathbin{\cdot}X=d\pi(X^{\mathsf T})v$. Transpose reverses brackets,
fixes the diagonal coroots, and converts those lower generators to the
upper generators of the original left action. Hence the dominant
coefficients remain nonnegative with their original signs. The dimension
formula is valid over this same real field: its proof uses formal power
series \cite[Chapter~VIII, Section~3, pp.~279--280, formulas~(40)--(41)]{Jacobson1964}.
Use the separate explicit root description for $\mathfrak{sl}_m$
\cite[Chapter~IV, Section~6, pp.~153--154, formulas~(55)--(58)
and diagram~(59)]{Jacobson1964}. These are the actual adjoint weights
of the trace-zero diagonal Cartan.
Its positive roots have consecutive simple-root coefficients. The trace
pairing identifies a simple root with the diagonal vector
$e_i-e_{i+1}$; its squared norm is $2$, so all normalized simple-root
length ratios are $1$. Applying the general root product with precisely
these coefficients and length ratios gives
\[
 \dim V_a=\prod_{0\le i\le j<m-1}
      \frac{\sum_{k=i}^j(a_k+1)}{j-i+1}.
\]
Reindexing by interval length and initial index gives the numerical Weyl
product used below. Highest weight zero is identified with the trivial
real module by uniqueness. The first and last fundamental weights are
identified with the actual standard and contragredient modules: their
first and last coordinate vectors satisfy the required highest-vector
equations, and elementary off-diagonal matrices prove irreducibility.
The resulting real linear intertwiners are continuous in finite dimension.
Their infinitesimal identities intertwine the same exponential curves;
connectedness then gives intertwiners of the original group representations.
The same construction applies to the Levi representations in the normal
collars. It may also be performed in finite real coordinates: pull the
subspaces, component maps and intertwiners back through the same basis.
The sum identity, dimensions and group equivariance are then preserved on
the original representation.

The finite-order suspension estimate in \cite[Section 3.2]{BRHW}
is used in the same original atlas. The actual jet-word bound and the
near-minimal lift give the strict compact-uniform integrability in
Remark~\ref{II-loc:rem:uniform-FM}, at each fixed jet order.
At each order apply
superrigidity inside the algebraic jet group, simultaneously
linearize the semisimple and compact parts, and strip the
linear transfer as in Section~\ref{I-sec:formal}. This reduces
the jet to its \emph{given} derivative. Composing with the
one derivative frame already chosen gives
\eqref{II-loc:eq:signed-levi-jet} at every finite order, with the same
linear frame and the same sign cocycle.

Signed recurrence at degree one makes the absolute central
blocks constant. Tempering of that frame along $a_1$ identifies
them as
\[
                     \lambda=e^{-n/(n-1)},\qquad b=e^q.
\]
Indeed the standard block carries the weights
$\varepsilon_j-\varepsilon_1$, and the Levi-trivial block
carries $q\varepsilon_1$. Both central scalars contract.
Choose a finite $N$ with $\lambda_+^{N+1}<\lambda_-$.
The complex-domain word estimates recorded after
Corollary~\ref{II-loc:cor:pure-linear-levi} give the analytic return
integrability. Theorem~\ref{II-loc:thm:signed-contracting} supplies
the asserted curves in the original centered fiber charts.
Their tangent line is the Levi-trivial derivative line, hence
the transverse weight line.
\end{proof}

This corollary settles the construction of actual transverse germs
for $q<0$. Its centers still form a measurable field. The next
section establishes the additional tangential compatibility needed
for projective neighborhood description. For $q=0$, the equation in
Lemma~\ref{II-loc:lem:bounded-centralizer} loses its strict contraction,
and the analytic limit also fails the required spectral inequality.

\subsection{Contracting tangential germs and the overlap argument}
\label{II-loc:sec:negative-overlap}

Tempered coordinate changes and the action of the centralizer
are basic features of nonstationary normal-form theory
\cite{Melnick,Kalinin}. We need their analytic counterparts
with explicit control of the convergence radii, and prove
the required statements in this subsection.

We continue with the negative branch of
Corollary~\ref{II-loc:cor:negative-curves}. Write
$m=n-1\ge3$, $\kappa=n/m$, $\lambda=e^{-\kappa}$, and $b=e^q<1$.
The opposite radical is $V\simeq\R^m$, while $U\simeq(\R^m)^*$.
Conjugation by $a_t$ multiplies $V$ by $e^{-\kappa t}$ and
$U$ by $e^{\kappa t}$. We use the coordinates $H_x$ already
constructed, and set
\[
                    S_x=H_x(\R^m\times\{0\})
\]
initially as a germ. This is a tangential hypersurface, distinct
from the transverse curve $i_x$.

\begin{lemma}[Analytic gauges have tempered radii]
\label{II-loc:lem:analytic-tempering}
Suppose $f_x$ is an analytic germ cocycle over an invertible
probability-preserving map $T$, and actual analytic germs $h_x$
satisfy $h_{Tx}^{-1}f_xh_x=D_x$, where the linear maps $D_x$
and their inverses are uniformly bounded. If the analytic radii
and bounds of $f_x,f_x^{-1}$ have integrable logarithms, then
$h_x,h_x^{-1}$ have tempered analytic radii and bounds.
No integrability of the gauge itself is assumed.
\end{lemma}
\begin{proof}
For a germ $f=\sum_{j\ge1}f_j$ in fixed finite-dimensional
coordinates, let $|f_j|_1$ be the sum of the absolute values
of all degree-$j$ coefficients, and put
\[
 K(f)=\max\left(1,\sup_{j\ge1}|f_j|_1^{1/j},
                         \sup_{j\ge1}|(f^{-1})_j|_1^{1/j}\right).
\]
It is finite and measurable for an invertible analytic germ.
Coefficient majorization by $Kt/(1-Kt)$ gives
\[
                      K(fg)\le2K(f)K(g).
\]
Indeed the composition majorant is
$K(f)K(g)t/[1-(1+K(f))K(g)t]$; use the same argument for
the inverse composition. Thus the conjugacy equation and its
inverse give an integrable bound for
$|\log K(h_{Tx})-\log K(h_x)|$.
The integrable-increment argument in Section~\ref{II-sec:negative} implies
$\log K(h_{T^d x})/|d|\to0$. The coefficient majorants on
balls of radius $1/(2K(h_x))$, followed by Cauchy estimates
on smaller balls, give the claimed tempered domains and bounds.
\end{proof}

In the compact analytic suspension the lemma applies to $H_x$
for each fixed element of $L A_0$: the normalized map is linear,
and the original return word has integrable length. A countable
collection of such elements can therefore be used on one conull
set. Continuing the local coordinates along $a_1$ gives an
analytic diffeomorphism
\begin{equation}\label{II-loc:eq:global-stable-chart}
 H_x:E\longrightarrow W^{s,F}_{a_1}(x).
\end{equation}
Explicitly, for large $d$ define
$H_x(w)=a_1^{-d}H_{a_1^dx}(D_x^{(d)}w)$.
Contraction and tempered radii make this consistent, injective,
and onto the full fiberwise stable manifold. The latter is open
in its fiber because all fiber exponents are negative.
The Levi identities extend to these global stable charts.
Consequently $S_x$ is a connected analytic hypersurface, closed
inside $W^{s,F}_{a_1}(x)$, with its intrinsic embedded structure.

\begin{proposition}[The radical preserves the tangential germ]
\label{II-loc:prop:negative-tangential}
The germ field $S_x$ is $Q^\circ$-equivariant. In the coordinates
of Theorem~\ref{II-loc:thm:signed-contracting}, for each fixed
$0\ne v\in U$ its normalized germ has the form
\begin{equation}\label{II-loc:eq:negative-radical-form}
       g_v(x)(u,z)=\bigl(A_{v,x}(vu)u,\ C_{v,x}(vu)z\bigr),
\end{equation}
where the two functions are analytic and nonzero at zero.
The field $S_x$ is also equivariant under the finite component
group of the full parabolic.
\end{proposition}
\begin{proof}
Signed unipotent recurrence in the stabilizer $L_v$ gives,
as in Step 4 of Theorem~\ref{II-loc:thm:signed-contracting},
\[
       g_v(x)(u,z)=\bigl(A_{v,x}(vu,z)u,B_{v,x}(vu,z)\bigr).
\]
To check this centralizer description, on $vu\ne0$ the stabilizer
of the pair $(v,u)$ has fixed space $\R u$ in the standard
module, and $L_v$ acts transitively on each affine hyperplane
$vu=w\ne0$. Polynomial equivariance in each degree, and then
analyticity, give the displayed form also at $w=0$.

Each scalar coefficient of $A$ and $B$ changes only by a sign
under $L_v$. Its absolute value is constant almost everywhere:
a nontrivial unipotent in $L_v$ acts ergodically by Howe--Moore
for the ergodic $L$ action. Covariance under $L$, which is
transitive on $U\setminus\{0\}$, makes these absolute values
independent of the nonzero parameter $v$, when coefficients
are expressed in the variable $w=vu$.

Put $c=e^\kappa=\lambda^{-1}$. The group relation gives
\[
                    g_{cv}(a_1x)=D_{vx}g_v(x)D_x^{-1}.
\]
For the coefficient of $w^rz^s$ in $A$, the right side changes
its absolute value by $b^{-s}$, when written in the new
variable $(cv)u$. For the corresponding coefficient in $B$,
the factor is $b^{1-s}$. Since these absolute values do not
depend on $v$, only $s=0$ can occur in $A$, and only $s=1$
in $B$. This proves \eqref{II-loc:eq:negative-radical-form}. Together
with the Levi and central identities it proves germ equivariance.

There is a useful uniqueness observation. A measurable
$L$-equivariant analytic hypersurface field tangent to the standard
summand is necessarily $S_x$: write it as $z=f_x(u)$ and apply
signed recurrence in each degree to
$\operatorname{Hom}(\operatorname{Sym}^j\R^m,\R)$, which has
no $L$-fixed vectors for $j\ge1$. Thus all coefficients vanish.
Transport by an element of $Q\setminus Q^\circ$ gives another
such field; that element normalizes the Levi and preserves its
standard and trivial summands. The same uniqueness proves
equivariance for the full parabolic.
\end{proof}

\begin{lemma}[Two-dimensional stable leaves determine the germ]
\label{II-loc:lem:negative-root-overlap}
Let $W< V$ be a one-dimensional subgroup obtained by Levi
conjugation of a missing-root subgroup. For almost every root-graph
point $y=v\cdot z$, where $v\in W$ and
$z=\phi_x^W(v)$ is the BRHW fiber root graph, one has
\begin{equation}\label{II-loc:eq:negative-root-overlap}
                            v^{-1}S_y=S_x
\end{equation}
as hypersurfaces in the full stable fiber, with the centers
identified by suspension transport.
\end{lemma}
\begin{proof}
We need the following null-set consequence of the root graph laws:
for every prescribed original null set $N$, almost every center has
Haar-almost no root parameters whose graph points belong to $N$.
It is enough to establish this on each positive Cartan ergodic component;
we do not need an unconditional canonical root disintegration for the
original measure as an additional intermediate assertion.

Use the finite-fiber and Cartan-component argument in the proof of
\cite[Corollary 4.5]{BRHW}. The finite-atomic input is
\cite[Lemma 5.6, pp.~971--972 of the published version]{BRHWFactors}.
Write the resulting finite decomposition as
$\mu=\sum_C\mu(C)\mu_C$, where
$\mu_C=\mu(C)^{-1}\mu|_C$. On each such component the root-kernel
subgroup $J=\ker\beta$ is ergodic: the finite-fiber bound and
base ergodicity give a positive lower bound for the mass of every
positive $J$-invariant event for the original measure. For a
$J$-invariant event $S$ modulo $\mu_C$, apply this bound to $C\cap S$ in the
original measure and divide by $\mu(C)$; one does not assume that
the original fiber weights remain uniform after restriction.
This gives finitely many $J$-components of $\mu_C$.
The connected Cartan normalizes $J$ and fixes this finite partition;
Cartan ergodicity leaves a single component. Thus the actual
$J$-ergodic conditional probability is $\mu_C$ almost everywhere.
Each $\mu_C$ projects to the original Haar probability: its projection
is absolutely continuous with respect to that probability and Cartan
invariant, so base ergodicity identifies the normalized projection.
The finite diagonal sign group acts transitively on the Cartan components
and normalizes each root group. If any component were invariant under
the missing root, all would be, contradicting noninvariance of their
original weighted sum. The original Lyapunov identities restrict to
these absolutely continuous component measures with the same weights.
Each missing-root fiber space is one dimensional. Apply
\cite[Proposition~7.2(3) and Section~7.3, pp.~52--53]{BRHW} on each component.
The parameter maps are embeddings into the intrinsic coarse fiber
leaves of Claim 3.9; equality of ambient dimensions is not used.

In the finite spectrum, the positive root ray is cut out by finitely many
negative half-spaces, so the coarse root collection is strongly integrable.
Choose a time in its negative cone and use the leaf charts and linearization of
\cite[Propositions~9.9 and~9.12, pp.~82--85, and Section~11.2, p.~96]{BRHW}.
Their auxiliary finite order is chosen above the global critical regularity
and the selected Cartan-time ratios; analyticity permits this choice.
The total linearization embeds the root product as
$(t,z)\mapsto(t,\iota^{\mathrm{fib}}(z))$; Claim~3.9 gives surjectivity,
hence a homeomorphism onto the full intrinsic leaf. Measurable realization
follows from canonical subordinate conditional kernels.

When an ergodic discrete Cartan return is needed in the subordinate-partition
construction, it can be selected inside the prescribed nonempty open cone.
Indeed, for an ergodic $\R^k$-action on a finite measure space with separable
complex $L^2$ and continuous Koopman representation,
\cite[Theorem~1, p.~116]{PughShub1971} places the exceptional times in a
countable union of affine hyperplanes. Such a union has Lebesgue measure zero.
Apply this on each of the finitely many Cartan components and avoid the
union of their exceptional sets. The same selection then works on every
component. For the actual suspension, the standard Borel structure and
probability measure make complex $L^2$ separable. Joint measurability of
the action implies measurability of every Koopman orbit's distance to a
fixed $L^2$ vector: write its squared distance as the integral of the
squared pointwise difference of jointly measurable representatives.
An isometric representation with these measurable distances on a separable
space has continuous orbits. Indeed, countably many orbit preimages of
balls of radius $\epsilon/2$ cover the group; one has positive Haar measure,
and the Haar Steinhaus theorem applied to its difference set gives
continuity at the identity. Translation gives continuity everywhere.
The original source requires joint continuity on the group times $L^2$.
For isometries, this follows from orbit continuity and the estimate
\[
 \|U_g f-U_h v\|_2
 \le \|f-v\|_2+\|U_g v-U_h v\|_2.
\]
The inverse-action representation and forward composition differ only by
inversion of the group parameter. Thus the exact Koopman continuity
hypothesis is derived from the original jointly measurable action.
The precompact partition below is constructed separately on the same
native leaves.

The original clean return foliation retains measurable embedded $C^1$
discs \cite[Definition~4.1 and Proposition~4.6(a), pp.~470, 473]{BRHWSmoothErgodic2023}
\cite[Section~3.2, p.~960]{BRHWFactors}. The local lamination assertion there
applies to these same full leaves; capping the distance at one preserves
its local balls. Take compact good center sets with complements tending to zero.
In a local product chart, intersect a core with a closed ambient
neighborhood inside the chart; its inverse image has compact transverse
projection. Cut the corresponding plaques by closed balls with radii tending
to the chart radius and sphere images of zero component measure. These
compact boxes have measurable transverse labels. Almost every center inside
a box lies in its intrinsic plaque interior; outside centers have an ambient
neighborhood avoiding the box. Finite chart covers of the cores give a
countable cover modulo zero. Label each point by its first containing box
and transverse label. Avoiding the finitely many earlier closed boxes gives
an intrinsic neighborhood in its atom, whose closure lies in the selected
compact plaque. All original leaves remain whole.

The canonical subordinate conditional kernel of
\cite[Proposition~3.11(4)]{BRHW} is locally equivalent to Haar measure
on the original pointwise graphs. Its measurable root-parameter decoder
makes a positive local equivalence radius measurable. Form the countable
sum of the conditional kernels pulled back from Cartan returns.
Recurrence to a set where the radius is bounded below and the inverse
parameter dilations give equivalence with Haar measure on the full graph.
For every original null set, each summand vanishes almost everywhere,
so this sum gives the full graph null law. After pushforward by the
decoder to $\R\times M$, its zero-mass code on a countable basis
determines the support, hence the continuous graph;
the Borel inverse of this injective code gives a jointly measurable
parametrization agreeing at every parameter almost everywhere.
An original $\mu$-null set is null for each $\mu_C$, and the finite
component assembly preserves these exact normalized laws.
Counting the actual lattice representatives transfers the conclusion
to the original suspension lift and supplies regular target centers.

Identify $W$ with a line $\R w\subset\R^m$. For any two-plane
$I$ containing this line, choose a Levi-conjugate Cartan element
$a_I$ whose fiber stable space is exactly $I$. Such a choice
exists because $m\ge3$: make the two corresponding missing
roots negative and all other missing roots positive, with
positive sum. Since
$\varepsilon_1=-n^{-1}\sum_{j=2}^n\beta_j$, the transverse exponent
$q\varepsilon_1$ is then positive. The globally linearized
Cartan dynamics and Lemma~\ref{II-loc:lem:analytic-tempering} identify
\[
                    W^{s,F}_{a_I}(x)=H_x(I).
\]
For a coordinate missing-root line, choose a Cartan whose value there
is $-1$ and whose other missing-root values are $3$. Since $3m-4>0$,
the transverse exponent is positive and its one-dimensional stable
fiber is $H_x(\R w)$. We retain the original clean coarse fiber leaf,
with its null-set removal and immersed structure, from the unstable-intersection
construction
\cite[Section~3.2, p.~960]{BRHWFactors}
\cite[Definitions~4.1 and~4.5, Proposition~4.6, and Definition~4.10, pp.~470--474]{BRHWSmoothErgodic2023}.
Lattice independence permits a basis containing $-H$, and gives forward
contraction for each fixed $\beta(H)<0$. To pass to the native fiber,
use the compatible original domain $D$ and metric family of
\cite[Section~2.2.1(3)--(4), p.~950]{BRHWFactors}. Their uniform form
comparison, compactness of $M$, and the fiber path dictionary bound the
fixed $M$-distance by the capped intrinsic fiber distance. The same $D$
has finite first moment for the original right-invariant Riemannian
distance. Word length is bounded by an affine function of that distance
\cite[Remark~4.B.3(1), Proposition~4.B.10(1), and Example~4.C.9(1)]{CornulierHarpe2016}:
apply the homomorphism bound to inversion from the opposite group and
use $|g^{-1}|=|g|$. Its Borel section therefore has integrable word length. The integrable
deck gauge between this section and the prescribed controlled section
has subexponential distortion, so native forward stable landing is
preserved. This proves the required application of Claim~3.9(3).
Consequently $z\in H_x(\R w)\subset H_x(I)$ for every plane in the
pencil. Fixed Levi transport gives the assertion for the stated $W$.

We spell out the transport fact used here. If $v$ belongs to
a base subgroup contracted by $a$, and $z\in W^{s,F}_a(x)$
with $y=vz$ a regular center, then
\[
                  W^{s,F}_a(y)=vW^{s,F}_a(x).
\]
Indeed $a^nva^{-n}$ stays in a compact set. Its fiber Lipschitz
constants and those of its inverse grow subexponentially along
the Haar-typical base orbit, by the compact-parameter return
estimates of \cite[Section 3.2]{BRHW}. It therefore preserves
exponential convergence of pairs of fiber points. Applying this
fact to $a_I$ gives $v^{-1}H_y(I)=H_x(I)$.

Use a countable dense collection of two-planes containing $w$.
All identities hold simultaneously for almost every root-graph
point. A finite subcollection shows that $v^{-1}S_y$ and $S_x$
have the same tangent space at $z$. In the coordinates $H_x$,
write the former as an analytic graph over the latter near $z$.
Its graph function vanishes on every plane in the dense collection.
Since the coordinate of $z$ lies on $\R w$, each such affine
plane through $z$ has direction $I$. Every homogeneous Taylor
polynomial of the graph function vanishes on these planes,
whose union is dense in $\R^m$. All Taylor coefficients vanish.
This proves equality of germs.

For the contracting central Cartan, the coordinate contraction
factor is bounded by $\max(e^{-(m+1)},e^{qm})<1$.
The tempered forward and inverse chart bounds identify the targets
of the same whole charts with the full native stable fibers.
The transport argument above therefore gives a common target for
the charts at the two root-graph centers. Their transition and its
inverse are analytic on all of $\R^{m+1}$. The normal component of
each transition restricted to the tangential hyperplane vanishes
near the graph point by germ equality. The analytic identity
theorem on the connected hyperplane makes it vanish everywhere.
The two resulting inclusions prove
\eqref{II-loc:eq:negative-root-overlap} for the entire hypersurface
images. Fixed Levi covariance transports this equality to every
stated conjugate root subgroup.
\end{proof}

\begin{remark}
This is where dimension $n\ge4$ is essential to this argument.
For $m=2$, the plane containing two missing-root directions is
already the full standard summand. No Cartan element separates
it from the contracting transverse direction, so the preceding
two-plane proof is unavailable.
\end{remark}

\begin{lemma}[Finite fiber branches force parabolic coherence]
\label{II-loc:lem:finite-branch-coherence}
Suppose a $Q$-invariant Haar-projecting suspension measure has
uniformly $k$ atoms in almost every fiber. Suppose its full stable
leafwise measures are Haar graphs of analytic injective maps
$\phi_x:V\to M_x$, with $\phi_x(1)=x$. Then these maps are
$Q$-coherent on the connected chart-overlap domains $V(q)$ of
\cite[Definition 8.5]{BRHW}.
\end{lemma}
\begin{proof}
Lift to $G\times M$ and write the conditional probabilities as
$\mu_g$. The Haar graph condition implies, for almost every $g$
and Haar-almost every $v$, the transport identity
\begin{equation}\label{II-loc:eq:finite-branch-transport}
 \mu_g=\frac1k\sum_{i=1}^k\delta_{x_i}
 \quad\Longrightarrow\quad
 \mu_{vg}=\frac1k\sum_{i=1}^k\delta_{\phi_{(g,x_i)}(v)}.
\end{equation}
Choose a complete lacunary Borel section $\Sigma$ by
\cite[Corollary~1.2, p.~284; definition p.~283]{KechrisCountableSections1992}.
Its counting pullback on $V\times\Sigma$ is sigma-finite; disintegrate over
an equivalent finite transverse law. In a bounded injective group chart,
the original plaque atlas gives countably many open Haar pieces of the same kernel.
Bounding their measurable coefficients gives finite-mass open unions.
Their Haar identities and a shared positive piece imply finite kernel mass
on the whole chart. Apply the finite-plaque law again to that whole chart,
including any Haar-positive kernel-null region. A bounded exhaustion and
one positive finite Haar reference give a single measurable factor, absorbed
into the transverse measure. Nonnegative integration reconstructs the original
measure; Fubini gives \eqref{II-loc:eq:finite-branch-transport}.
Compare \cite[Section~4.1, (4.1)--(4.2), Proposition~4.2]{BenoistQuint2011}
and \cite[Theorem~6.3]{EinsiedlerLindenstrauss2010}.

For almost every pair $(g,q)$, $\mu_{qg}=\mu_g$, and all its
$k$ atoms are regular centers for both graph families. The base
factorization $\tau_q(v)q=\widehat q_vv$, with
$\widehat q_v\in Q$, and $Q$-invariance give
$\mu_{\tau_q(v)qg}=\mu_{vg}$ for almost every $v$ in the open
cell. To justify this assertion for the varying element
$\widehat q_v$, first choose the measurable disintegration as a
function of the coset $Qg$, and then apply Fubini to the smooth
open-cell coordinates. Equality of the two measures in
\eqref{II-loc:eq:finite-branch-transport} gives equality of their
unordered sets of branches.

Choose disjoint neighborhoods of the finitely many distinct
points $x_i$. By continuity, for all sufficiently small $v$ the
two branches issuing from $x_i$ both lie in its neighborhood.
Their matching in the unordered sets must therefore have the
same index. Thus
\[
       \phi_{(qg,x_i)}(\tau_q(v))=\phi_{(g,x_i)}(v)
\]
almost everywhere near $1$. Continuity gives equality throughout
that neighborhood, and analyticity extends it to the connected
domain $V(q)$. This is the required coherence. No continuous
labeling of the atoms as $g$ varies has been assumed.
\end{proof}

\begin{proposition}[Analytic projective graphs in the negative branch]
\label{II-loc:prop:negative-projective-graphs}
For almost every $x$ there is an analytic diffeomorphism
$\phi_x:V\to S_x$, taking the identity to $x$, such that the
full $a_1$-stable leafwise measure is Haar measure on its graph
\[
                          v\longmapsto v\cdot\phi_x(v).
\]
In the lifted space $G\times M$ these graphs are $V$-, $Q$-,
and $\Gamma$-coherent with the identity formulas displayed in
\cite[Definition 8.11(2)--(4)]{BRHW}. Here the target is the
hypersurface $S_x$; no assertion is made that it equals the entire
fiber stable leaf in Definition 8.11(1).
\end{proposition}
\begin{proof}
Compose the given coordinate root graphs in a fixed order to obtain a
measurable graph $\phi_x$ over $V$. The original leafwise null laws
propagate along this finite word by successive disintegration. Thus its
values remain positive atoms of the original fiber disintegration and
lie in the native full stable fiber. The relative-pair recurrence
argument of Lemma~\ref{II-loc:lem:atomic-fibers} implies that two positive
atoms in the same native stable fiber coincide. For almost every base
point and parameter the resulting map therefore injects the finite
uniform source fiber into the finite uniform target fiber. Their
cardinalities and atom weights agree, so it transports the full fiber
conditional measure. Haar disintegration of the base now gives the
full graph disintegration identity. The same uniqueness identifies
alternative root words wherever both are defined.

To identify the canonical full stable leafwise measures, fix a finite
subordinate box and its transverse label. On a conull set of regular
centers, positive-atom stable separation identifies any two centers in
the same stable plaque after translation of their root parameters.
Continuity upgrades this reanchoring identity from almost every parameter
to every parameter. Haar inversion and the full graph disintegration
identity therefore make the measure on pairs consisting of a center and
a point of its graph symmetric. Restrict this pair measure to pairs in
the box with equal transverse labels. The parameter set of each good
row has positive Haar measure because it contains a neighborhood of the
identity; it has finite Haar measure because the plaque is compactly
contained and the root subgroup is closed embedded. Reanchoring makes
the row mass constant along each good label fiber. Dividing by that
actual row mass gives a reversible probability kernel, constant on label
fibers and supported on the original plaque. The canonical conditional
kernel on this same finite box is also reversible and supported on its
label fiber. Constancy and support give the two opposite composition
identities between these kernels; taking adjoints identifies the kernels.
Thus the canonical conditional is exactly normalized Haar measure on the
same graph. Conditional-measure uniqueness
\cite[Proposition~3.11(4)]{BRHW} and a countable exhaustion give the full
leafwise measure class. The original probability and its counting lift
agree on every measurable quotient-injective sheet: pushing forward the
restriction of the counting lift gives the restriction of the original
measure to that sheet's image. The identification therefore concerns the
original suspension probability.

Lemma~\ref{II-loc:lem:negative-root-overlap} makes $S_x$ constant,
after transport, along each given root graph. Successive null-law
transport along the same finite coordinate word shows
$\phi_x(v)\in S_x$ for Haar-almost every $v$.
We next prove regularity; it is not part of this measurable
disintegration assertion.

For two typical centers $y=v\phi_x(v)$, consider the transition
in the stable coordinates, restricted to the standard summand:
\[
           P_{x,y}=H_x^{-1}\circ v^{-1}\circ H_y
                         :\R^m\longrightarrow\R^m.
\]
It is affine. To see this, iterate both centers by $a_1$.
The base displacement becomes $e^{-d\kappa}v$. The coordinate
germs and inverses at both centers have tempered bounds, as do
the fiber maps for these bounded base displacements. The
constant term of the transition contracts at rate $\lambda$,
so it eventually lies in the tempered inverse-chart domain.
Cauchy estimates on tempered smaller balls show that every
fixed-degree coefficient of the iterated transition grows
subexponentially. Its degree-$j$ coefficient, however, changes
in absolute norm by exactly $\lambda^{1-j}$ per iterate,
up to orthogonal signs. Thus it vanishes for $j\ge2$.

For the root-return identity, use the original positive conditional
atoms supplied by the contraction argument of
\cite[Lemma~5.6, pp.~971--972]{BRHWFactors}.
The original leafwise null-set law puts almost every root-graph point
at such an atom. Literal Cartan conjugacy puts its transported point
and the target root graph over the same base, in the same native stable
fiber. Positive-atom stable separation makes them equal almost everywhere
in the root parameter; continuity gives equality at every parameter.
This argument requires no strict root action or invariant counted carrier.

The individual root graph, in its centered $H_x$ coordinates,
is linear in the root parameter. Indeed the source and target
both scale by $\lambda$ under $a_1$, up to a target sign.
Subtract its derivative at zero and use its $C^2$ remainder
along recurrent times where that remainder is bounded. The
rescaled remainder is bounded by a constant times
$\lambda^d|t|^2$, and hence is zero. The global stable
coordinates justify this rescaling for every fixed parameter.

Put $f_x=H_x^{-1}\phi_x$. Along almost every affine coordinate
line in $V$, disintegration identifies $f_x$ with a centered
linear root graph followed by the affine transition $P_{x,y}$.
Thus it is affine on almost every such line. Repeated interpolation
in the individual coordinates makes $f_x$ a multiaffine
polynomial modulo null sets by the finite-grid argument of
\cite[Corollary~9.7]{BRHW}. Fubini allows the interpolation nodes
to be fixed outside a null set of base parameters, so the resulting
coefficients are measurable in $x$. Its degree-$j$ homogeneous term
obeys the same covariance with norm factor $\lambda^{1-j}$.
Recurrence eliminates $j\ge2$. Hence
\[
                       \phi_x(v)=H_x(K_xv,0),
\]
where $K_x$ is invertible: the derivatives of the individual
root graphs are nonzero in their distinct root directions.
This proves the analytic diffeomorphism assertion.

Uniqueness of conditional measures gives $V$-coherence, and
descent from the suspension gives $\Gamma$-coherence.
Lemma~\ref{II-loc:lem:finite-branch-coherence} supplies $Q$-coherence,
using the finite atomic disintegration from Lemma~\ref{II-loc:lem:atomic-fibers}.
\end{proof}

\begin{theorem}[Analytic hypersurface description for $q<0$]
\label{II-loc:thm:negative-recognition}
For any lattice $\Gamma<\SL(n,\R)$ in the negative line or
hyperplane branch, $n\ge4$, the measure
is supported on the suspension of finitely many compact embedded
analytic projective covering hypersurfaces. After finite-index
restriction each is invariant and its action is analytically
identified with the standard or dual projective action.
\end{theorem}
\begin{proof}
Let $C\subset G\times M$ be the support of the original Haar lift.
The finite uniform conditional atoms and the analytic graphs constructed
above give local finite graph descriptions of $C$. Indeed, uniqueness of
disintegration and Haar transport identify these graphs on every translated
positive Gauss cell based at a Haar-typical center. Those cells cover $G$:
for any prescribed point, the possible centers form a nonempty open
set, which meets the conull set of good centers. Continuity and full Haar
support extend the conditional identities to the whole cell. Separation of
the original atoms makes these descriptions genuine finite covering charts.

The original measure stabilizer preserves $C$ under the full group $Q$.
Consequently
\[
 D=\{(Qg,x):(g,x)\in C\}\subset Q\backslash G\times M
\]
inherits the literal finite fibers and root covering charts from $C$.
These charts cover every $Q$ coset and show that $D$ is a closed finite
cover of $Q\backslash G$. The standard or dual
projective dictionary identifies this base with $\RP^{n-1}$. Thus $D$ is
compact, since $M$ is compact. Its covering charts make it locally connected;
its connected components are therefore open, and compactness makes their
number finite. Each component covers the whole connected projective base.

The physical projection $D\to M$ is injective. To see this, take
$(g,x),(h,x)\in C$. Choose one Haar-good center whose full $Q$ Gauss cell
contains both group points. Such centers exist: the required two nonzero
matrix coefficients define a nonempty open set, and full $Q$ absorbs either
sign. Express both points in the finite graph rows at this same center.
The shared physical point and separation of the original atoms identify
their original atom labels; injectivity of the original analytic graph then identifies
the root parameters. Hence $Qg=Qh$. The compact physical projection is
therefore an embedding, and distinct component images are disjoint.

Near every point of $D$, its physical projection agrees with one original
analytic graph composed with the inverse projective root coordinate.
This coordinate is analytic with injective derivative in the standard
projective atlas. The original graph is an analytic immersion, so each
physical component is an embedded analytic hypersurface.
The covering classification now identifies each component with the sphere
or projective space; use \cite[Theorem~1.38, p.~67, and
Exercise~1.43, p.~74]{Hatcher}, together with the pointed lifting and
uniqueness statements of \cite[Propositions~1.33--1.34, pp.~61--62]{Hatcher}.

The original diagonal lattice action preserves $C$ and descends to $D$,
permuting its finitely many components. Pass first to the kernel of this
permutation action. On each spherical cover the remaining ambiguity is a
deck sign; intersecting these finitely many sign kernels gives one subgroup
of finite index in the original lattice. The resulting actions are the
standard or dual projective actions, or their spherical lifts, in the
analytic parametrizations just constructed. Finally, the original Haar
lift gives full mass to $C$. Descent gives full original measure to the
suspension of these hypersurfaces. This carrier is closed and invariant,
so it contains the support of the original measure.
\end{proof}

For $n\ge4$, the full-collar Theorems~\ref{II-suspension:thm:negative-collar},
\ref{II-positive:thm:collar}, and \ref{II-neutral:thm:full-collar} apply
according to each hypersurface's actual normal character. They use
the hypersurface conclusion above, without the arithmetic collar
construction, and exclude this branch on an uncovered compact invariant
remainder in the union-of-model-neighborhoods argument.

\section{Neutral description and collar descent}\label{II-series:neutral}\label{II-sec:neutral}
We first recognize the neutral measures and their analytic
hypersurfaces. These hypersurface conclusions precede every
full-collar argument and are also used to normalize scalar
boundary cocycles. The nonzero collar constructions are in
Sections~\ref{II-sec:suspension-negative-collar} and
\ref{II-sec:suspension-positive-collar}. Sections
\ref{II-sec:neutral-fixed-curves}--\ref{II-sec:neutral-moving-centers}
then construct the full neutral collars for arbitrary lattices.

\subsection{Neutral measures give tangential stable plaques}
\label{II-loc:sec:neutral-plaques}

The next result concerns $q=0$ and does not use the contracting
transverse construction. It separates a consequence of entropy
from the additional geometric problem of recognizing a neighborhood.

\begin{proposition}[Neutral root graphs and stable plaques]
\label{II-loc:prop:neutral-plaques}
Suppose the line-parabolic measure has weights
\[
              \{c\beta_j:2\le j\le n\}\cup\{0\},\qquad c>0.
\]
Its fiberwise relative entropy is zero for every Cartan element,
and $c=1$. Its unconditional missing-root leafwise measures are
the smooth Haar graphs of \cite[Proposition~7.2(3)]{BRHW}.
For every sufficiently large finite $r$, the full contracting
tangential leaf has a $C^r$ graph parametrization
\[
                 \phi_x:V\longrightarrow W^{s,F}_{a_1}(x),
\]
whose graph carries the full stable leafwise measure as Haar
measure on $V$. Here $\dim W^{s,F}_{a_1}=n-1$.
\end{proposition}
\begin{proof}
Work first on a Cartan ergodic component. The original charts satisfy
\cite[Proposition~2.5, p.~956]{BRHWFactors}; the standing hypotheses pass
to integer subgroups by \cite[Remark~3.3, p.~469]{BRHWSmoothErgodic2023},
and the fiber foliation is tame by \cite[Section~3.2, p.~960]{BRHWFactors}.
Its whole fibers are subordinate, so Definitions~6.1, 6.3--6.4 and
Remark~6.2 of \cite[pp.~484--485]{BRHWSmoothErgodic2023} identify leaf
entropy with relative Kolmogorov--Sinai entropy. The relative Ruelle
inequality \cite[Theorem~7.2(22), p.~488]{BRHWSmoothErgodic2023} therefore
gives zero at every ergodic time with nonpositive fiber exponents:
positive full tangent blocks have trivial intersection with the fiber tangent.
Recurrence of the measurable change of frame identifies existing rates in
the prescribed section with native rates, without a section moment assumption.
For any nonidentity $a\in A$, choose ergodic $u,v$ in the common positive
fiber cone with $a=uv^{-1}$ and independent logarithms
\cite[Theorem~1, p.~116]{PughShub1971}.
Ruelle applies to $u^{-1},v^{-1}$; inverse symmetry and global conditional
subadditivity \cite[Theorem~13.3(1), p.~518]{BRHWSmoothErgodic2023}, for
a lattice containing $u,v$, give
\[
 h_\mu(a\mid\mathcal F)\le
 h_\mu(u\mid\mathcal F)+h_\mu(v^{-1}\mid\mathcal F)=0.
\]
The identity case is immediate. The original parabolic law is a finite
mixture of these Cartan components with the same Haar marginal. Its finite
component label contributes only bounded block entropy, so division by
block length passes the same zero entropy to the original law.

\looseness=-2
On each Cartan ergodic component, root inclusion
\cite[Section~3.3.4, p.~22]{BRHW} gives positive full-tangent exponents at
nonzero Cartan clocks. Adapted charts with sufficiently small error
\cite[Proposition~5.4, p.~480]{BRHWSmoothErgodic2023} give uniform
contraction on full unstable discs. Native fiber coarse leaves admit
measurable embedded $C^1$ discs
\cite[Section~3.2, p.~960]{BRHWFactors}, hence precompact intrinsic
subordinate partitions by the local Lusin construction. For each one-dimensional missing-root
fiber space, apply
\cite[Proposition~7.2(3)]{BRHW}.
By \cite[Proposition~3.12]{BRHW}, zero fiber entropy gives a measurable
fiber coarse partition. Contraction and recurrence give a conull
fiber transversal. The conditional graphs meet it; their
Haar base marginals imply domination of each unconditional subordinate
conditional by the same graph law. Recurrence makes its
root-kernel label almost everywhere constant. Further conditioning
preserves the canonical kernel, and
\cite[Proposition~3.11(4)]{BRHW} gives normalized Haar.
Use the first proof paragraph of
\cite[Proposition~7.4]{BRHW}, retaining the zero exponent.
Contraction and recurrence give stable representatives of finite
component labels. Leaf saturation equates component and original
subordinate kernels, gluing the Haar laws to the original measure.

Their nonzero center derivatives intertwine $e^{\beta_j(a)}$ and
$e^{c\beta_j(a)}$. Recurrence where the derivatives and their inverses
are bounded gives $c=1$. Measurable interval capture and recurrence
exhaust each equivariant root curve; the return chart coordinates identify
it $C^r$-smoothly with the prescribed intrinsic leaf and transport the
same Haar law.

For every element in the open projective contraction cone the
stable fiber space is exactly the sum of the missing-root
spaces: the remaining exponent is zero. The hypotheses of
\cite[Proposition 7.5]{BRHW} are therefore satisfied. That
proposition gives the claimed full stable graph. Here the fiberwise
stable leaf is the intrinsic Pesin leaf of
\cite[Claim 3.9 and Section 3.5.5, p.~24]{BRHW}: it is the
path-connected fiber component through the original point. For each return, apply
\cite[Theorem~5.1(a)--(c), Case~III, pp.~44--45, proof
pp.~46--54]{Ruelle1980Author} in complexified centered charts.
Its real locus gives plaques tangent to the negative spaces,
identified with chart
rate sets by cone uniqueness. Contraction and forward radius bounds
give nested original inverse-return images along a subsequence.
Their compatible charts define the intrinsic immersed leaf; compare
\cite[proof of Theorem~6.3, pp.~52--53]{Ruelle1979}. The suspension
geometry of \cite[Proposition~3.2, p.~21]{BRHW} supplies the
required moment bounds. The coarse foliation is the finite restricted join of
\cite[Proposition~4.6(b) and Definition~4.10, pp.~473--474]{BRHWSmoothErgodic2023}. The graph
embedding carries the full stable leafwise law, without identifying
its physical projection with an entire metric-defined stable set.
The same comparison preserves the prescribed counting lift.
Analyticity of the ambient action permits any required
finite regularity, in particular the critical regularity required by
\cite[Proposition 7.5, pp.~54--55]{BRHW}.
\end{proof}

\begin{lemma}[Analyticity of the neutral tangential parametrizations]
\label{II-loc:lem:neutral-analytic-plaques}
The maps $\phi_x$ in Proposition~\ref{II-loc:prop:neutral-plaques} are
analytic, with measurable positive local radii. Their image
plaques have tempered holomorphic local parametrizations along
$a_1$. This asserts analyticity along each plaque, not analytic
dependence on a transverse parameter.
\end{lemma}
\begin{proof}
We give the analytic stable-manifold argument because the
suspension is noncompact. Put $T=a_1$ and
$\lambda=e^{-\kappa}$. The analytic return maps and their
inverses have logarithmically integrable complex bounds and
radii. Lyapunov coordinates split their derivative into a
stable block with exponent $\log\lambda$ and a center block
with exponent zero. The coordinates and their inverses are
tempered.

We use the complete-probability holomorphic stable-set theorem in
\cite[Theorem~5.1(a)--(c), Case~III, pp.~44--45, with proof
pp.~46--54]{Ruelle1980Author}. These are the pages of the March 1980
author preprint. Complete the original probability without changing its
maps or cocycle values. The analytic-height moment supplies the actual
complex radius and the logarithmic moment of the full $C^1$ norm,
including the first derivative. For each prescribed positive window
exponent, the theorem gives one family of stable sets and forward
control windows. Measurable analytic parametrizations are additional
constructions, as explained below. The compact
manifold formulation in \cite[assertion (4)]{RS} is a comparison,
not a direct compactness argument for this suspension. Here is
the holomorphic construction in the present coordinates. Along one
regular orbit write the maps in these coordinates as
\[
 f_d(s,z)=(A_ds+a_d(s,z),\ B_dz+b_d(s,z)),
\]
where $\|A_d\|\le\lambda e^\epsilon$ and
$\|B_d^{\pm1}\|\le e^\epsilon$. On sufficiently small
tempered complex balls the nonlinear terms have uniformly small
Lipschitz constants and vanish to second order. Choose
$\lambda e^\epsilon<\eta<e^{-2\epsilon}$.
On sequences with norm $\sup_{d\ge0}\eta^{-d}\|(s_d,z_d)\|$,
solve, with $s_0$ prescribed,
\[
 \begin{split}
 s_d&=A(d,0)s_0+
             \sum_{j=0}^{d-1}A(d,j+1)a_j(s_j,z_j),\\
 z_d&=-\sum_{j=d}^{\infty}B(j+1,d)^{-1}b_j(s_j,z_j).
 \end{split}
\]
Here $A(d,j)$ and $B(d,j)$ are the derivative products from
time $j$ to time $d$. Both geometric-series operators are
bounded, and shrinking the balls makes the nonlinear operator
a contraction. A small sequence ball stays within the available
tempered domains because $\eta<e^{-\epsilon}$. Its fixed point
depends holomorphically on $s_0$: the successive approximations
are holomorphic and converge uniformly on a smaller complex
ball. Reality and uniqueness give the actual real strong stable
plaque. The same estimates with tempered envelopes give tempered
local bounds and radii. No center-manifold theorem is used.

Canonical leafwise covariance and uniqueness,
\cite[Proposition~3.11(2),(4), p.~25]{BRHW}, transport the same
unconditional graph law under $T$. Haar full support and continuity
give $F(T)\phi_x(v)=\phi_{Tx}(\lambda v)$ for almost every original
center and every parameter $v$; $F(T)$ is the actual return in the
chosen suspension section.

It remains to check the parametrization, rather than just its
image. Let $J_x=D_0\phi_x$. Its covariance is
$D_xF(T)J_x=J_{Tx}\lambda I$. The integrable-increment argument
applied to the norms of $J_x$ and its inverse onto the stable
space makes this frame tempered. In the analytic stable charts
with derivative $J_x$, the return germ is
\[
                   f_x^s(u)=\lambda u+O(\|u\|^2)
\]
with tempered analytic estimates. The Koenigs construction in
Section~\ref{II-sec:negative}, now with the single contraction rate $\lambda$,
converges because $\lambda^2<\lambda$; compare the bundle
normalization theorem \cite[Theorem~1.1 and Section~3.1]{BDM2008}.
It gives analytic
parametrizations $H_x^s$, with derivative $J_x$, satisfying
$F(T)H_x^s=H_{Tx}^s\lambda I$. Continue them by contraction
to the whole stable leaf.

The $C^2$ map $J'_x=(H_x^s)^{-1}\phi_x$ is tangent to the
identity and satisfies
$J'_{T^nx}(\lambda^nv)=\lambda^nJ'_x(v)$.
Recurrence to a set with a common positive $C^2$ radius and
a bounded quadratic remainder gives
\[
                \|J'_x(v)-v\|\le C\lambda^d\|v\|^2
\]
along an unbounded sequence of times. The continued stable
charts justify this identity for each fixed $v$. Hence
$J'_x=\id$, and $\phi_x=H_x^s$ is analytic.
\end{proof}

\begin{theorem}[Description of the atomic neutral branch]
\label{II-loc:thm:neutral-atomic-recognition}
Let $\Gamma<\SL(n,\R)$ be any lattice, $n\ge4$, and suppose
$q=0$ in the line or hyperplane branch.
If the suspension measure has atomic fiber disintegration,
then it is supported on finitely many compact embedded analytic
projective covering hypersurfaces.
\end{theorem}
\begin{proof}
For completeness, atomicity here automatically means uniform
finite atomicity. The mass $\mu_{\pi(x)}(\{x\})$ is measurable
and invariant under the measure-preserving minimal parabolic.
Ergodicity makes it constant almost everywhere. If there is any
atomic part, the invariant set of atoms is conull, and the
constant is $1/k$ for some finite $k$. Otherwise the conditional
measures are atomless almost everywhere.

Proposition~\ref{II-loc:prop:neutral-plaques} and
Lemma~\ref{II-loc:lem:neutral-analytic-plaques} give the original full
stable Haar graphs with analytic injective parametrizations and injective
differentials. Retain this same physical family, the original section,
and its original canonical leafwise law. The intrinsic stable-leaf
construction and the controlled section comparison place every graph
point in the native stable fiber through its original center. Thus
positive-atom stable separation applies to these same graphs.

Use the canonical-support construction in the proof of
Theorem~\ref{II-loc:thm:negative-recognition}. Its finite graph covering
argument requires the original Haar graph law, uniform finite atoms,
full $Q$-invariance of the original lifted support, and stable separation.
All these properties hold here: the original canonical Haar graph
class, uniform finite atoms and native landing give exact graph
disintegration for the same original lifted conditional kernel.
The group $Q$ preserves the original measure, and the preceding paragraph
gives the required native landing. The same common-good-center
argument therefore proves physical incidence, so the descended finite
cover of $Q\backslash G$ embeds in $M$. Compactness gives finitely many
components. Their physical maps are analytic immersions by the original
analytic graph parametrizations and inverse projective root coordinates.
The covering classification and the finite permutation and deck-sign
kernels give the compact embedded sphere or projective hypersurfaces,
with one common finite-index subgroup. The original lifted support
again descends to a closed full-measure carrier of the same suspension
measure. No complete intrinsic-development or local-admission assertion
is needed for this finite-atomic construction.
\end{proof}

For $n\ge4$, these hypersurfaces receive full collars from
Theorems~\ref{II-suspension:thm:negative-collar}, \ref{II-positive:thm:collar},
and \ref{II-neutral:thm:full-collar}, according to their actual normal
characters. Their recognition, including its scalar-cocycle uses below,
does not depend on any full-collar construction.

For an atomless neutral measure, the tangential results do not
establish $Q$-coherence of the stable plaques. In fact, that
conclusion would exclude a required neutral model.
For the regular flow model on $S^{n-1}\times\R/T\Z$, write
the action in a projective affine chart and a suitable transverse
coordinate as
\begin{equation}\label{II-loc:eq:neutral-moving-model}
 \begin{aligned}
 l(u,t)&=(lu,t),\\
 a_s(u,t)&=(e^{-\kappa s}u,t+s),\\
 v(u,t)&=\left(\frac{u}{1+vu},\ t+\log|1+vu|\right)
                    &&(v\in U).
 \end{aligned}
\end{equation}
These are the local formulas of the analytic action
$(\theta,t)\mapsto(g\theta/\|g\theta\|,
t+\log\|g\theta\|)$ on the compact sphere-circle product.
The central transverse exponent is zero. The tangential stable
plaques are $t=\mathrm{constant}$, but the row unipotent sends
such a plaque to a different graph. At $u=0$ its derivative is
\[
                      \begin{pmatrix}I&0\\v&1\end{pmatrix}.
\]
Thus the shear is present in an actual compact analytic model.
In the neutral case description must either produce a compact
projective leaf at a transverse flow zero, or identify the moving
plaques and their transverse flow in an open neighborhood.
The tangential-graph proposition alone proves neither alternative
for an atomless measure.

\subsubsection{A convergence criterion for a nonzero neutral shear}

The regular model supplies more directions than one stable plaque.
In fact the analytic change
\[
                    p=e^t u,\qquad z=e^t-1
\]
turns its centered parabolic germs into
\[
 l(p,z)=(lp,z),\quad
 a_s(p,z)=(e^{-\kappa s}p,z),\quad
 v(p,z)=(p,z+vp).
\]
For the central map the output coordinate is centered at the new
transverse point $t=s$. In particular the conjugate contracting
elements have stable tangent planes $z=vp$. This observation
gives the analytic criterion used below. Its formal-reduction
hypothesis will be established for the nonzero-shear alternative
in Proposition~\ref{II-loc:prop:neutral-formal-dichotomy}.

\begin{lemma}[Convergence detected on a pencil of hyperplanes]
\label{II-loc:lem:hyperplane-convergence}
Let $H(p,z)$ be a finite-dimensional vector-valued formal power
series at zero in $\R^m\times\R$. If
$p\mapsto H(p,vp)$ is convergent for almost every $v$ in a
positive-measure set of row vectors, then $H$ is convergent.
The same assertion holds for coefficient-measurable families,
with a positive measurable convergence radius in the conclusion.
\end{lemma}
\begin{proof}
Write $H=\sum_{j\ge0}H_j$, with $H_j$ homogeneous of degree $j$.
By a countable union over bounds and radii, there are a bounded
set $E$ of row vectors of positive measure and constants $C,R$
such that
\[
       \sup_{\|p\|\le1}|H_j(p,vp)|\le C R^j
                           \quad(v\in E,\ j\ge0).
\]
Choose $p$ in a small ball separated from zero. The map
$(p,v)\mapsto(p,vp)$ has full rank there. The image of this ball
times $E$ has positive $(m+1)$-dimensional measure: this follows
also by the linear coarea formula for $v\mapsto vp$ at each
fixed $p\ne0$. On a bounded positive-measure subset of that
image, $H_j$ therefore has an exponential bound independent of
$j$.

We recall the elementary polynomial estimate needed to extend
the bound. On any fixed box, a polynomial of degree at most $j$
is bounded by $C_E^j$ times its supremum on a fixed measurable
subset of positive measure, for a constant depending only on
the box and the subset. In one variable choose $j+1$ ordered
interpolation points in that subset whose separations are at
least a constant times $|i-i'|/(j+1)$. Lagrange interpolation
and $\sum_{i=0}^j(i!(j-i)!)^{-1}=2^j/j!$ give the exponential
bound. For several variables, Fubini supplies a positive-measure
set of sections with a common positive lower bound on their
lengths; repeat the one-variable estimate and induct on dimension.

Apply this estimate to each component of $H_j$ on a box about
zero. Its supremum is bounded exponentially in $j$. Homogeneity
then gives uniform convergence of $\sum H_j$ on a smaller ball.
Measurability follows from the usual countable coefficient tests
for exponential bounds.
\end{proof}

\begin{proposition}[Conditional neutral convergence from shear]
\label{II-loc:prop:neutral-shear-convergence}
Suppose an actual analytic parabolic suspension cocycle admits
a measurable formal conjugacy to linear germs. Assume that in
these coordinates the radical has the form
\[
                   (p,z)\longmapsto(p,z+vp),
\]
and the central derivative has the form
$c_x\operatorname{diag}(\lambda I,1)$, where $0<\lambda<1$
and
\[
 \lim_{d\to\pm\infty}|d|^{-1}\log|c_x^{(d)}|=0.
\]
Suppose the original return germs have the
logarithmic analytic integrability used above. Then the formal
conjugacy is convergent almost everywhere.
\end{proposition}
\begin{proof}
The prescribed first-derivative frame is tempered: a scalar Lyapunov
weight bounds the linear model and its inverse, and the original germ and inverse
moments make the logarithmic frame-height increment integrable. The
nonnegative-coboundary argument then gives two-sided subexponential
height. Thus the central return has spectrum $\{\log\lambda,0\}$ in
this same frame. Normalize a nonzero finite law to a probability and
complete it; the zero law is vacuous. For its original coefficient height
$H_x\ge1$, set $r_x=(4(d+1)H_x)^{-1}$ in complex dimension $d$.
The rescaled return
$U_x(z)=r_{Tx}^{-1}F_x^{\mathbb C}(r_xz)$ is holomorphic near the closed
unit ball with $C^1$ norm at most $4/r_{Tx}$. This gives the required
logarithmic moment, while $\log r_{T^nx}/n\to0$ preserves the spectrum.

Apply \cite[Theorem~5.1(a)--(c), Case~III, pp.~44--45; proof,
pp.~46--54]{Ruelle1980Author}, using the author preprint, with
$c=\log\lambda/4$ and $f=\log\lambda/2$.
With $0<\alpha_x<\beta_x<1$ and $\gamma_x\ge1$, the stable tube
\[
 S_x=\{z:\|z\|\le\alpha_x,\quad
               \|U_x^nz\|\le\beta_xe^{cn}\ (n\ge0)\}
\]
has the prescribed tangent, interior holomorphic geometry and pair
contraction $\gamma_xe^{fn}$. For any fixed $\Theta>0$, its controls
satisfy $\alpha_{T^nx}\ge A_xe^{-\Theta n}$,
$\beta_{T^nx}\ge B_xe^{-\Theta n}$ and
$\gamma_{T^nx}\le C_xe^{\Theta n}$ on a conull forward invariant set.

We obtain measurable charts from these same tubes as follows. The cone
and radial-graph estimates in the source proof make stable projection
injective and cover the stable ball of radius $\alpha_x/\sqrt2$, with
preimages in the interior. Its inverse $G_x$ is holomorphic there, tangent to the inclusion.
It is jointly measurable by the compact unique-zero criterion, applied
to the projection equation and the countably many tube inequalities;
radially clipped iterates give continuous residuals without changing the
tube. Conjugation invariance and uniqueness make $G_x$ real.
The full shear $z\mapsto G_x(P_xz)+z-P_xz$ has inverse
$z\mapsto z-G_x(P_xz)+P_xz$. Cauchy estimates bound both coefficient
heights by $A/\alpha_x^2$ on the ball $\alpha_x/4$.
Their measurable height $P\ge1$ therefore has a forward exponential
bound. For a larger exponent $a$, the finite envelope
$L(x)=\log\sup_{n\ge0}e^{-an}P(T^nx)$ has increment at most $a$.
The bounded-coboundary argument gives $L(T^{\pm n}x)/n\to0$, proving
bilateral tempering without an additional chart moment. Measurable
coefficient versions and analytic realization descend these same germs
from the completion to the original sigma algebra.

Since $f<c$, a sufficiently small graph segment maps into the next
tube, giving the actual stable return. A scalar Lyapunov change gives
$m=\lambda e^{-\epsilon}$ and $M=\lambda e^\epsilon$ with $M^2<m<1$.
The full return and inverse are tempered compositions of the shear,
the original frame and scalar dilations; their stable restrictions
inherit these coefficient bounds. Exponential envelopes supply the slow
analytic radii. Apply \cite[Theorem~1.1, pp.~2138--2139; proof in
\S3.1, pp.~2147--2150]{BDM2008} to this restriction and its linear part.
The normalized-iterate limit and its local inverse give an analytic
stable parametrization with the prescribed first derivative; real
measurable iterates give reality and measurable dependence.

The restriction of the formal conjugacy to the stable axis equals
that analytic parametrization. To check uniqueness, compare the
first different homogeneous coefficients of two parametrizations.
In degree $j\ge2$ their covariance has exponents
$(1-j)\log\lambda$ in the stable output and
$-j\log\lambda$ in the neutral output. Both are strictly
positive. A nonzero measurable coefficient obeying such a
covariance would grow exponentially on a probability-preserving
orbit, contradicting recurrence to sets on which it is bounded.
Induction gives equality of every coefficient.

The original analytic radical conjugacy transports this convergent
stable-axis restriction to each fixed shear plane $z=vp$.
Thus for each fixed $v$ the restriction is analytic at almost
every center. Fubini makes it analytic for almost every $v$ at
almost every center. Apply
Lemma~\ref{II-loc:lem:hyperplane-convergence}. The nondegenerate first
derivative and the analytic inverse function theorem give a
convergent inverse as well.
\end{proof}

\subsubsection{Calibrating the formal reduction by the Haar graphs}

Write $V$ for the opposite, column radical. The Haar graphs represent the
canonical full stable class. Each fixed Levi element, central time or
normal-orientation flip $A$ preserves the original probability.
Original fiber equivariance \cite[Section~8.3, p.~60]{BRHW}, the product
identity \cite[Claim~3.9(1), p.~24]{BRHW}, and linear normalization of $V$
give full leaf-set equivariance.
Mix canonical conditionals of countably many subordinate factors, whose
atoms exhaust each whole leaf, with positive weights summing to one.
This gives a Markov kernel $K$ in the original class
\cite[Definition~3.10 and Proposition~3.11(4), pp.~24--25]{BRHW}.
Coherence becomes class agreement at all good centers in one leaf:
equal factor labels give one nonzero conditional dominated by both
kernels, and exhaustion supplies that label. Each conditional, hence
$K$, is reversible. The transported kernel $H_x=A_*K_{A^{-1}x}$ is
reversible with the same coherence, so $H\circ K\sim H$ and
$K\circ H\sim K$. Detailed balance swaps these mixed pair laws;
thus the joint laws of $K$ and $H$ are equivalent. Conditional absolute
continuity gives $A_*K_x\sim K_{Ax}$, proving the commuting-actor extension
of Proposition~3.11(2) on the original probability. The constant Haar
Jacobian is absorbed in the class. Continuity and Haar full support give
the physical parameter identity for every parameter, at almost every
center for each fixed actor.
The actual $L\ltimes U$ action and Mautner's argument make $U$ act trivially
on the $L$-component factor $Z$.

\begin{proposition}[Formal neutral radical dichotomy]
\label{II-loc:prop:neutral-formal-dichotomy}
Assume $n\ge4$. In the actual neutral suspension, one can choose measurable
formal charts $H_x(u,z)$ such that
$H_x(u,0)=\phi_x(u)$ and the Levi germs are $(u,z)\mapsto(lu,z)$.
In these charts the infinitesimal radical germs have the form
\begin{equation}\label{II-loc:eq:neutral-radical-field}
 Y_{v,Z}=(vu)\big(a_Z(z)E_u+b_Z(z)\partial_z\big),
 \qquad E_u=\sum u_i\partial_{u_i}.
\end{equation}
Here $a_Z,b_Z$ are formal one-variable series. On an ergodic
parabolic component exactly one of the following holds:
\begin{enumerate}
\item $b_Z(0)=0$ almost everywhere. The analytic stable plaques
      are locally $Q^\circ$-coherent.
\item $b_Z(0)\ne0$ almost everywhere. There are actual analytic
      charts, still calibrated by $H_x(u,0)=\phi_x(u)$, with
      centered parabolic germs
      \begin{equation}\label{II-loc:eq:calibrated-neutral-model}
      l:(u,z)\mapsto(lu,z),\quad
      a_t:(u,z)\mapsto(e^{-\kappa t}u,z),\quad
      v:(u,z)\mapsto(u,z+vu).
      \end{equation}
      These calibrated charts are unique and descend from the
      auxiliary normal-orientation extension.
\end{enumerate}
The formal input is the compatible finite-jet reduction
Section~\ref{I-sec:formal} and
\cite[Definition~1.3 and Theorem~1.4, p.~193 of the published version]{FisherMargulis},
applied to the Levi action.
\end{proposition}
\begin{proof}
The derivative $D_0\phi_x$ gives an exact standard Levi frame
on the stable space. The zero Oseledec line is Levi invariant.
Its absolute derivative is an additive cocycle after taking
logarithms; Lemma~\ref{II-loc:lem:additive}, in its locally compact
property-$(T)$ version, normalizes it to one, leaving a sign.
For measurable orientation extensions, the Borel action can first
be made continuous in a Polish topology on the same space with
the same Borel sets, by \cite[Theorem~1.1]{BeckerKechris1993};
this preserves the probability law used in ergodic decomposition.

First this line is radical invariant. For a fixed nonzero row
$v$, its derivative component into the stable space is
equivariant under the row stabilizer $L_v$. The signed
unipotent-recurrence argument in Section~\ref{II-sec:negative} makes that component
fixed by the unipotent subgroups of $L_v$. The common fixed
space in the standard module is zero when $m\ge3$. The line
is therefore invariant under $U$. Its orientation extension
is now a probability extension for the full parabolic. On it,
normalize the positive Levi derivative to one. Thus the Levi
derivative is exactly $\operatorname{diag}(l,1)$; there is no
compact sign on the stable block because its frame came from
the Haar graph.

Compatible finite-jet reduction gives a formal Levi linearizer.
The compact-uniform moment comes from the original section height through
the symmetric jet size $\log(1+\|A\|+\|A^{-1}\|)$, and passes to
the same conditional probabilities. No moment of the selected formal transfer
is assumed.
For completeness, commuting with the Levi in the cocycle sense
forces the normalized central germs to have the form
\[
       (u,z)\longmapsto(A_{t,Z}(z)u,B_{t,Z}(z)).
\]
Indeed unipotent recurrence puts each coefficient in the
Levi-fixed subspace, and the remaining coefficients are
Levi invariant. Their stable derivative is
$A_{t,Z}(0)=e^{-\kappa t}$. Restrict the formal chart to $z=0$.
It agrees formally with $\phi_x$: at the first differing degree
$j\ge2$, covariance under $a_1$ has exponents
$(j-1)\kappa$ in the stable output and $j\kappa$ in the neutral
output. Recurrence excludes a nonzero measurable coefficient.
Induction proves the asserted calibration.

For a fixed row $v\ne0$, the formal radical germ is equivariant
under $L_v$. Apply unipotent recurrence coefficient by
coefficient. The scalar coefficients of the resulting
$L_v$-equivariant tensors are $L_v$-invariant functions. Mautner's argument on each Levi ergodic component makes them
functions on $Z$. Levi covariance for different $v$, and the
fact that $U$ fixes $Z$, give a measurable homomorphism from
$U$ into the formal germ group over each $Z$. At every finite
jet order this homomorphism is continuous, hence differentiable.

Here the invariant theory is elementary. An infinitesimal
field is linear in the row $v$. An $\SL(m)$-equivariant vector
polynomial linear in $v$ is a scalar series in $z$ times
$(vu)u$; an invariant scalar polynomial linear in $v$ is a
scalar series in $z$ times $vu$. For $m\ge3$ there is no
additional alternating tensor with these symmetries.
This gives \eqref{II-loc:eq:neutral-radical-field}.

If $b_Z(0)=0$, every radical field preserves $z=0$, as do
the central and Levi germs. The equality with the analytic
restriction $\phi_x$ turns the formal identities on this plane
into actual identities. This is local parabolic coherence.
Zero versus nonzero normal shear is intrinsic and preserved
by the central action, so ergodicity gives the dichotomy.

In the nonzero case make a formal change commuting with the
Levi and equal to the identity on the stable plane:
\[
 (u,z)\mapsto(C_Z(z)u,D_Z(z)),\qquad
 \frac{C_Z'}{C_Z}=-\frac{a_Z}{b_Z},\quad
 D_Z'=\frac{C_Z}{b_Z},\quad C_Z(0)=1,\quad D_Z(0)=0.
\]
It takes $Y_v$ to $(vu)\partial_z$. A germ
$(A(z)u,B(z))$ normalizing these shears by
$v\mapsto e^{\kappa t}v$ has $A$ constant and
$B(z)=e^{\kappa t}A z$: substitute $z+vu$ and compare the
coefficients of $u$ and $vu$. Calibration gives
$A=e^{-\kappa t}$, hence $B(z)=z$.
Proposition~\ref{II-loc:prop:neutral-shear-convergence} applies
with scalar factor one and proves convergence of $H_x$.

Finally, a change between two calibrated reductions is
Levi equivariant with Levi-invariant scalar coefficients.
Mautner makes these coefficients $U$-invariant, so the change
commutes with all shears. The common formal centralizer of
$\operatorname{diag}(L,1)$ and these shears consists of scalar
homotheties: write the change as $(C(z)u,D(z))$ and compare
$z+vu$. Calibration on $z=0$ forces that scalar to be one.
This also proves descent under the normal-orientation deck map.
\end{proof}

If the original parabolic is disconnected, first apply the
argument to the finitely many ergodic components for its
identity component. Zero versus nonzero shear is intrinsic,
so the finite component group preserves the alternatives.
In the regular case uniqueness of the calibrated chart gives
the required covariance under this finite group as well.
The geometric statements below allow the corresponding finite
angular covers; no ergodicity on an arbitrary probability
extension is being assumed.

\subsubsection{Stable transitions in the regular neutral case}

For this stable-transition calculation, identify a column $v$ with
$\left(\begin{smallmatrix}1&0\\v&I_m\end{smallmatrix}\right)$.
If the Haar graphs were parametrized using the negative column instead,
make the fixed simultaneous changes
\[
 \phi_x^+(u)=\phi_x^-(-u),\qquad
 H_x^+(u,z)=H_x^-(-u,-z).
\]
The scalar map $-\mathrm{Id}$ commutes with all three calibrated
models in \eqref{II-loc:eq:calibrated-neutral-model}. Thus these changes
preserve the calibration and model identities, and negation preserves
the Haar parameter measure. In this subsection $\phi_x,H_x$ denote
these positive-column families. Keeping the original negative-column
families would replace $1-z$ in \eqref{II-loc:eq:stable-transition} by $1+z$.

Write $\widehat v x=v\cdot\phi_x(v)$. These maps form a jointly
measurable measure-preserving action modulo null sets on the original
probability: the group law holds almost everywhere for each fixed pair.
Choose one subordinate partition $\mathcal P$ for the original full
leafwise family. The original conditional compatibility in
\cite[Proposition~3.11(2)--(4), p.~25]{BRHW} identifies its normalized
Rokhlin rows with the restrictions of Haar measure to
$I_x=\{t:\mathcal P(\widehat t x)=\mathcal P(x)\}$, pushed forward by the
same graph. Reanchoring and Haar full support give the graph identity for
almost every reached parameter and every further parameter.
For each fixed $t$, put
\[
 D_t=\{x:\mathcal P(\widehat t x)=\mathcal P(x)\}.
\]
Translation of the normalized rows gives exactly
\[
 (\mu|_{D_t})\mathbin{\circ}(\widehat t)^{-1}=\mu|_{D_{-t}}.
\]
Indeed, translation by $t$ identifies the parameter restrictions
$I_x\cap(I_x-t)$ and $I_x\cap(I_x+t)$; integrate using the original
Rokhlin barycenter identity.

Choose a measure-preserving central Cartan return contracting the root
parameter by $0<a<1$. Subordination and continuity of the original graph
at its center give $\mu(D_{a^nt}^{\mathrm c})\to0$. Partial preservation
implies, for every measurable $B$,
\[
 \mu((\widehat{a^nt})^{-1}B)
 \leq\mu(B)+\mu(D_{a^nt}^{\mathrm c}).
\]
The exact Cartan covariance conjugates $\widehat t$ to
$\widehat{a^nt}$. Pull this inequality back by the Cartan return and let
$n\to\infty$. It gives $(\widehat t)_*\mu\leq\mu$, hence equality by
equal finite total mass, for every fixed $t$.
Finally, for fixed $v,w$, Haar translation and Fubini provide a generic
parameter $b$ at which reanchoring gives
$(\widehat v\circ\widehat w)\circ\widehat b
=\widehat{v+w}\circ\widehat b$ almost everywhere.
Preservation by $\widehat b$ cancels that last factor. This proves the
law for each fixed pair on the same original probability. Centering of
the original graphs gives the identity actor; no common strict invariant
carrier is needed below.

\begin{proposition}[The exact stable-transition relation]
\label{II-loc:prop:neutral-stable-transition}
In the second case of Proposition~\ref{II-loc:prop:neutral-formal-dichotomy},
the calibrated analytic charts satisfy
\begin{equation}\label{II-loc:eq:stable-transition}
 H_x^{-1}v^{-1}H_{\widehat v x}(u,z)
                  =(u+(1-z)v,z).
\end{equation}
The inverse at the left uses the branch along $\phi_x(v)$.
The assertion is an analytic identity on a neighborhood of
each point of the stable plane, not a common-radius assertion.
\end{proposition}
\begin{proof}
Central covariance continues $H_x$ to a neighborhood of its
whole stable plane. For a prescribed stable coordinate use
recurrence to a set with bounded chart sizes, contract that
coordinate by $a_d$, and pull the resulting chart back.
The domains are centered ellipsoids, their intersections are
connected, and the analytic identity theorem gives consistency.
The continuation restricts to the injective stable parametrization.

Denote the left side of \eqref{II-loc:eq:stable-transition} by $P_v(x)$.
Its constant term is $(v,0)$. Covariance under the column
stabilizer $L_v$, unipotent recurrence and the Mautner argument
on the Levi ergodic components give
\[
 P_v(x)(u,z)=(A_Z(z)u+B_Z(z)v,C_Z(z)).
\]
To check the invariant theory, take $v=e_1$. The stabilizer
acts by $(u_1,u')\mapsto(u_1+au',Bu')$ with
$B\in\SL(m-1)$. Its scalar polynomial invariants are constant
in $u$, and its equivariant vector polynomials are spanned by
$u$ and $e_1$. This proves the displayed form. Levi
transitivity on nonzero columns makes $A,B,C$ independent
of $v$. The action $\widehat V$ fixes $Z$, so its group law,
evaluated for two independent columns, gives
$C\circ C=C$ and $A^2=A$. Invertibility yields $C(z)=z$ and
$A(z)=1$. Also $B(0)=1$. Central covariance makes $B$ invariant
under $A_0$, so parabolic ergodicity makes it one fixed
analytic germ $B(z)$.

An exact commutator now determines $B$. Track a finite word
by its base matrix and by the inverse of its chart transition.
The row element $U_i(t)$ has model $z\mapsto z+tu_i$;
the stable element $\widehat V_j(s)$ has model
$u_j\mapsto u_j-sB(z)$. For $i\ne j$ their base commutator
is $I-tsE_{j+1,i+1}$. The model commutator, followed by the
inverse of this Levi element, fixes $z$ and sends
\[
 u_j\longmapsto u_j+s\{B(z)-B(z-tu_i)+tu_i\}.
\]
Its base matrix is the identity, and its model fixes the
origin. Restrict first to a positive-measure set with bounded
chart sizes, and take $s,t$ small enough that all intermediate
points stay in these charts. On this set the actual word fixes
its center. Since the chart field is a function of that center,
its entire chart transition is the identity. The function $B$
is deterministic, so this positive-measure set suffices to give
the analytic identity
$B(z)-B(z-tu_i)+tu_i=0$, and $B'(z)=-1$.
Together with $B(0)=1$ this proves \eqref{II-loc:eq:stable-transition}.
The sign differs from the uncalibrated coordinates preceding
Lemma~\ref{II-loc:lem:hyperplane-convergence}: here the stable coordinate
is the root parameter in $\phi_x(v)$.
\end{proof}

\subsubsection{An invariant probability on disjoint compact hypersurfaces}

\begin{lemma}[Finite-orbit alternative for disjoint hypersurfaces]
\label{II-loc:lem:finite-hypersurface-probability}
Let a group $\Lambda$ act on a compact connected manifold and
suppose every finite-index subgroup of $\Lambda$ has finite
abelianization. An invariant ergodic probability on a measurable
family of pairwise disjoint compact connected embedded
hypersurfaces is supported on a finite orbit.
The family is understood as a measurable family of embeddings;
its incidence and separation relations are measurable.
\end{lemma}
\begin{proof}
We use homology of the actual ambient and domain manifolds with
coefficients $\mathbb F_2$, without an orientability assumption.
For a compact topological manifold, integral singular homology in every
degree is finitely generated
\cite[Appendix~E, Corollary~E.5, pp.~538--539]{Bredon1993}.
This applies to Hausdorff manifolds with Euclidean charts and does not
require a triangulation or a smooth structure. The coefficient sequence
$0\to\mathbb Z\xrightarrow{\,2\,}\mathbb Z\to\mathbb F_2\to0$
gives, for $j>0$, the exact segment
\[
 H_j(M;\mathbb Z)\longrightarrow H_j(M;\mathbb F_2)
       \longrightarrow H_{j-1}(M;\mathbb Z)
\]
\cite[Chapter~IV, Section~5, Example~5.8, pp.~180--181;
Theorem~5.6, pp.~178--179]{Bredon1993}.
The kernel of the second map is a finitely generated image, and its image
is a submodule of a finitely generated $\mathbb Z$-module.
Since $\mathbb Z$ is Noetherian, this image is finitely generated too;
the extension therefore makes the middle group finitely generated over
$\mathbb Z$. Its generators also generate over $\mathbb F_2$, so it is
a finite set. In degree zero, connectedness and local path connectedness
identify the ordinary singular homology with $\mathbb F_2$ by augmentation.
All these are the homology groups of the given spaces, with their original
maps and coefficient modules.

For a closed connected domain $N^r$ in positive dimension, the mod-two
fundamental-class theorem gives $H_r(N;\mathbb F_2)\cong\mathbb F_2$
\cite[Theorem~3.26, p.~236, proof pp.~236--238]{Hatcher};
in dimension zero the same augmentation supplies this identification.
An additive isomorphism over $\mathbb F_2$ is automatically linear, since
the only scalars are $0$ and $1$. Thus its inverse image of $1$ is the
unique nonzero homology class, preserved by every reparametrizing
homeomorphism, independently of a chosen orientation.

We derive the finite-cut count on these same spaces and maps. Write
$M^{n+1}$ for the ambient manifold, let $f_i:N_i^n\hookrightarrow M$
be a finite disjoint family, put $K=\bigcup_{i=1}^k f_i(N_i)$, and
let $U=M\setminus K$. The map from $\coprod_i N_i$ to $K$ is a
homeomorphism: it is a continuous bijection from a compact space to a
Hausdorff space. Its lifted component charts make this finite union a
compact topological $n$-manifold. Compactness and Euclidean charts also
give second countability, as required by the manifold convention in
\cite[Section~15.1, p.~358]{TomDieck}. A compact manifold is a
Euclidean neighborhood retract (ENR)
\cite[Corollary~18.4.7, p.~449]{TomDieck}; this concerns its intrinsic
topology and does not assert that the given embeddings are locally flat.

For a compact ENR subset of a compact manifold, singular mod-two duality
identifies $H^{n+1-j}(K;\mathbb F_2)$ with
$H_j(M,U;\mathbb F_2)$
\cite[p.~450, immediately before Proposition~18.4.10]{TomDieck}.
The canonical duality maps are natural under compact-subset inclusions
\cite[Lemmas~18.3.1--18.3.2, p.~445, and Lemma~18.5.1, p.~452]{TomDieck}.
Together with the actual pair long exact sequence
\cite[Theorem~2.16 and its singular-pair specialization, pp.~117--120;
coefficients p.~153]{Hatcher}, they give the exact segment
\[
 H^n(M;\mathbb F_2)\xrightarrow{\mathrm{res}} H^n(K;\mathbb F_2)
 \longrightarrow H_0(U;\mathbb F_2)
 \longrightarrow H_0(M;\mathbb F_2)\longrightarrow 0.
\]
Here the last zero is justified, rather than inferred from an assumed
nonempty complement. Homology vanishes above the dimension of a
connected topological manifold
\cite[Chapter~VI, Corollary~7.12, p.~346]{Bredon1993}.
Finite additivity on the actual disjoint union
\cite[Proposition~10.2.1, p.~247]{TomDieck} and field universal
coefficients therefore give $H^{n+1}(K;\mathbb F_2)=0$.
The latter theorem identifies singular cohomology with the linear dual
of singular homology by the natural Kronecker pairing
\cite[Section~3.1, pp.~196--199]{Hatcher}. The same additivity and the
actual domain fundamental classes give $\dim H^n(K;\mathbb F_2)=k$.
Naturality of that pairing identifies the rank of $\mathrm{res}$ with
that of the actual pushforward $H_n(K;\mathbb F_2)\to H_n(M;\mathbb F_2)$.
Under the disjoint-union coordinates, its image is precisely the span
of the original classes $(f_i)_*[N_i]$.

Since $H_0(M;\mathbb F_2)=\mathbb F_2$, exactness first makes
$H_0(U;\mathbb F_2)$ finite dimensional and then gives
\[
 \dim H_0(U;\mathbb F_2)
 +\dim\operatorname{span}_{\mathbb F_2}\{(f_i)_*[N_i]:1\le i\le k\}=k+1.
\]
The open complement is locally path connected, so its native singular
$H_0$ is the free $\mathbb F_2$-module on its actual connected components.
Finite dimensionality makes that component set finite and identifies
its cardinality with the displayed dimension. This proves the claimed
count without an orientability or local-flatness assumption, including
$n=0$; for $k=0$ it follows from connectedness of $M$ and the zero span.

Suppose the probability $\eta$ is atomless. Modulo two, there
are finitely many homology classes of hypersurfaces. Pass to
a finite-index subgroup acting trivially on homology and to
one of its ergodic components. Such a component is still
atomless, and almost every hypersurface now represents one
fixed class $h$.

If $h=0$, every hypersurface separates the manifold into two
components. To a hypersurface $N$ assign the smaller of the
two $\eta$-masses of hypersurfaces lying on its two sides.
This is a measurable invariant function. Its distribution
has no atoms, as follows. Fix a countable dense set of points
in the manifold. The hypersurfaces separating a specified
pair of these points are linearly ordered by inclusion of
the side containing the first point. Along this chain the
mass of that side is monotone; its increase between two
hypersurfaces is at least the $\eta$-mass of the intervening
part of the chain. Therefore any level set has zero mass
for the atomless measure on the chain. Taking the smaller
side mass only combines two such level sets. Every separating
hypersurface separates some pair of points in the dense set.
A countable union proves the assertion about the distribution.
It contradicts ergodicity.

Suppose $h\ne0$. For a finite subfamily of $k\ge2$ members,
the count gives $k$ complementary components. Removing one
member leaves $k-1$ components. The associated connected-component
coarsening therefore gives exactly two incident complementary
domains for that member. Each individual hypersurface is
nonseparating. The finite bipartite incidence graph is consequently
connected and two-regular, hence a cycle. No normal bundle or
local-flatness assertion is used. Deleting cuts contracts these
cycles compatibly, giving the family a circular order up to
reversal; the order on a singleton is unique.
Pass to the subgroup of index at most two preserving this order.
The invariant atomless probability supplies a circle coordinate
by measuring oriented arcs from a chosen origin. It pushes
$\eta$ to Lebesgue measure on the circle and intertwines the
action with rotations. The rotation displacement is a
homomorphism to $\R/\Z$. Finite abelianization makes its image
finite, which cannot act ergodically on an atomless circle
probability. This is the second contradiction.

Thus the original ergodic probability has an atom. Invariance
and ergodicity make it the uniform probability on a finite orbit.
\end{proof}

\begin{theorem}[Description of zero neutral shear]
\label{II-loc:thm:neutral-zero-recognition}
In the first case of Proposition~\ref{II-loc:prop:neutral-formal-dichotomy},
the fiber disintegration is finite atomic. Consequently
Theorem~\ref{II-loc:thm:neutral-atomic-recognition} recognizes this
case by compact analytic projective hypersurfaces for every
lattice. Their full analytic collars are constructed separately
according to the actual normal character, by
Theorems~\ref{II-suspension:thm:negative-collar}, \ref{II-positive:thm:collar},
and \ref{II-neutral:thm:full-collar}.
\end{theorem}
\begin{proof}
The zero alternative in Proposition~\ref{II-loc:prop:neutral-formal-dichotomy}
gives local image coherence of the analytic stable plaques on the actual
normalized components. Their weighted decomposition and the fixed
normalizer flip give this coherence for the full parabolic on the original
measure. Use the same complete intrinsic stable parameters and original
Haar counting lift throughout.

We first identify the corresponding Haar graph classes. In a Gauss chart
$h=vp g_0$, the rectified parabolic change has the form
\[
 (vp g_0,y)\longmapsto(\tau_q(v)q p g_0,y).
\]
Choose a parabolic-coset version of the original conditional kernel.
The Gauss identity makes this base change lie in the same parabolic
coset. Ambient Haar measure in the Gauss chart has the product Haar
class, so the rectified change is locally nonsingular. Local plaque
image coherence makes it a local map of the genuine total stable
foliation. Refine a source subordinate partition and the pullback of a
target subordinate partition by their paired labels. Both refinements
remain subordinate. Uniqueness of conditional measure classes transports
the canonical full stable graph classes of
\cite[Proposition~3.11]{BRHW}. Since their first coordinates determine
the parameters, equality of graph supports gives the exact projective
parameter change.

The germ-image condition is measurable: on compact parameter balls the
Borel evaluation map, with its center retained, is injective, so its
image is Borel by \cite[Theorem~4.5.4, pp.~153--154]{SrivastavaBorel}.
Countably many rational radii and parameter tests characterize local
image inclusion; compactness and the analytic inverse theorem recover
the germ identity. Fubini therefore applies at generic centers.
The re-centering identity of \cite[Claim~8.4]{BRHW} gives the required
joint almost-everywhere parabolic coherence. Analytic continuation
extends each identity on its connected cell where the normalized denominator is positive;
the other sign charts have their own seeds from full parabolic
coherence.

Together with the stable-unipotent and lattice identities, these give
coherent graph classes as in
\cite[Definition~8.11 and Claims~8.13--8.22]{BRHW}.
We retain the original matrix group, without assuming it simply connected.
Give each saturated graph the intrinsic topology of its translated Gauss
charts. The group projection is a local homeomorphism and the raw-pair
inclusion is continuous and injective, so this topology is Hausdorff.
Normalized lifts on convex stable-parameter cells agree by uniqueness;
they glue to complete stable translations. Their local formulas give
joint continuity, as for the full parabolic translations. The resulting
Gauss word charts make the group projection a covering. Quotienting by
the full parabolic, using its local principal sections, gives a connected
cover of real projective $m$-space. The simply connected sphere is its
double cover, so each such covering has degree at most two
\cite[Propositions~1.33--1.34]{Hatcher}. Its physical realization is
injective: put two projected points in one affine root chart and use
normalized-lift uniqueness and original plaque injectivity. Compactness
then makes the realization an embedding. The original analytic
full-rank plaque charts make its image analytic.

For separation we need the full native leaf belonging to these same
original graphs. Its product description in
\cite[Claim~3.9 and Proposition~7.5]{BRHW} is the union of stable-base
translates of the fiber leaf. The latter has the full range of the
\emph{same} physical parametrization. Indeed, the original section-moment
comparison takes a native fiber point of negative growth into the
physical return rate set; the cone and rate argument of
\cite[Theorem~5.1(a)--(c), Case~III, pp.~44--54]{Ruelle1980Author}
puts it in the same complete stable leaf. The global stable chart and
the specified section identify its parameter. Thus, at every point of
a suitably refined good set, the full native leaf is the cylinder
\[
 \{[u_vg,\phi_{(g,y)}(w)]:v,w\in\R^m\}
\]
in the suspension. The refinement retains this property pointwise.
Good points are Haar-conull in every intrinsic Gauss chart: apply
Fubini to the parabolic and stable identities and use the product Haar
class of that chart.

Suppose two developed images meet at $y$. A common affine projective
chart gives full-parabolic gauges making their group bases stable-aligned.
The condition that $y$ belongs to the complete stable image persists in
an intrinsic open neighborhood: for a lift $r$ of $y$, the map
$(v,q)\mapsto u_v(qr)$ has open image by the reversed Gauss chart,
and the inverse stable translation reaches $qr$, whose physical
coordinate is still $y$. Haar positivity and the chart-conull good sets
therefore give two good centers with exactly stable-aligned bases and
original plaques containing $y$. Their native cylinders intersect at
the same suspension point. The native leaf partition and the original
same-base graph-support uniqueness put the centers in one coherent
class. Hence distinct developed images are disjoint.

Countable dense words in the parabolic and stable translations give
Borel compact-leaf codes. To obtain a physical label, use the original
measurable boundary coordinate $b(y)\in Q\backslash G$ and a Borel
section of this quotient. Parabolic invariance and Fubini give one fixed
parabolic translate of that section whose lift is good almost everywhere.
Evaluating the compact code there gives a Borel label of $y$; any other
good lift has the same boundary coset and hence the same code.
The original deck identities give $\Gamma$-covariance. Push the original
conditional kernels to compact codes. The joint label law is invariant
under left parabolic and stable translations. These subgroups generate
$G$, since their Gauss product contains an identity neighborhood and
$G$ is connected. Conditional transport and Haar uniqueness make the
label probability independent of the group base point. A measurable
invariant subfamily pulls back, after countable lattice saturation, to
a parabolic-invariant set in the original suspension; this probability
is therefore ergodic.

Here the group is the original lattice $\Gamma$, not the auxiliary
integral Levi lattice. The ambient $\SL(n,\R)$ has property $(T)$ by
\cite[Theorem~1.4.15, pp.~54--55]{BHV}. The invariant quotient
probability constructed from the positive finite Haar covolume is regular;
thus the closed-subgroup theorem
\cite[Theorem~1.7.1, pp.~65--66]{BHV} gives property $(T)$ for this
same discrete $\Gamma$. The original universal corollary
\cite[Corollary~1.7.2, p.~67]{BHV} then says that every finite-index
subgroup $H<\Gamma$ has finite abelianization $H/[H,H]$.
This applies in particular to the finite-index kernels selected in the
leaf-space argument; the existing finite-index inheritance argument
keeps them inside this same original lattice.
Lemma~\ref{II-loc:lem:finite-hypersurface-probability} now gives a finite
orbit of leaves carrying the original law. The developed projection
of an original lifted point $(g,x)$ is $g^{-1}[e_1]$, with the fixed
dual convention when appropriate, by the incidence in
\cite[equation~(52)]{BRHW}. For fixed $g$, its preimage on each of the
finitely many leaves is finite. Thus the original fiber conditional
is finite atomic. Original parabolic ergodicity gives a common finite
cardinality and uniform atom weights. Apply
Theorem~\ref{II-loc:thm:neutral-atomic-recognition} to conclude.
\end{proof}

\subsubsection{Descent and description of the regular neutral model}

\begin{lemma}[The measurable boundary coordinate is unique]
\label{II-loc:lem:unique-boundary-coordinate}
For the neutral line or hyperplane measure in dimension
$n\ge3$, form the usual relatively invariant measure
$\overline\mu$ on $Q\backslash G\times M$ by integrating
the conditional probabilities $\mu_g$ against a probability
of Lebesgue class on $Q\backslash G$. Its projection to $M$
is a measurable isomorphism onto $(M,\widehat\mu)$.
Thus there is a single measurable boundary coordinate
$b:M\to Q\backslash G$ modulo $\widehat\mu$-null sets.
\end{lemma}
\begin{proof}
Recall that $Q_{\rm Out}$ is the smallest closed subgroup conull
for the original $G$-orbit conditional class
\cite[\S\S4.1 and 5.3, pp.~963--964 and 969]{BRHWFactors}.
On an injective quotient chart $U\times M$, pulling the original
measure back from the chart image gives exactly the original counting
lift restricted to $U\times M$. Disintegration by the physical
coordinate therefore gives the original local orbit-plaque rows.
A strictly positive finite base weight preserves their null classes;
conditional uniqueness identifies these classes with those of the
corresponding normalized restrictions of the global physical rows. A countable chart cover
determines the full null class, including patches of zero transverse
mass. At $(g,x)$ the change $h\mapsto hg^{-1}$ is exactly the original
orbit coordinate. Thus the smallest closed conull subgroup computed
from these shifted rows is this same original $Q_{\rm Out}$.

Use \cite[Proposition~5.4 and Claim~5.5, p.~970 of the published
version]{BRHWFactors}. The actual stabilizer $Q$ is maximal parabolic,
and $\SL(n,\R)$ has no rank-one simple factor. Proposition~5.4
therefore gives $Q_{\rm In}=Q_{\rm Out}$ for this same original
ergodic $P$-measure, without a dimension or neutrality assumption.
Claim~5.5 gives the asserted projection isomorphism. Choosing a
strictly positive finite base weight changes neither the joint
measure class nor its almost-everywhere inverse.
\end{proof}

\begin{proposition}[Compatible radial neighborhoods]
\label{II-loc:prop:neutral-radial-neighborhoods}
Assume $n\ge4$. In the regular-shear case there is a $\Gamma$-invariant open
set $W$ of full $\widehat\mu$-measure with the following
analytic structures. The boundary coordinate $b$ is an
analytic submersion on $W$. There are local coordinates in
$\R^n\setminus\{0\}$ in which $b$ is projectivization,
and changes between coordinates aligned with the same
boundary space are scalar homotheties. The radial Euler
field defines a nonvanishing analytic vector field $E$ on
$W$, commuting with the lattice action. In these coordinates
the lattice acts by its standard (or dual) matrices followed
by constant scalar multiplications on connected overlaps.
Moreover $\widehat\mu$ is equivalent to smooth measure on $W$.
\end{proposition}
\begin{proof}
Put $r=1-z$ in \eqref{II-loc:eq:calibrated-neutral-model} and
\eqref{II-loc:eq:stable-transition}. The inverse transitions of the
stable transformations become $u\mapsto u-rv$. All model
transformations are now linear in $(u,r)$. Track words as
in Proposition~\ref{II-loc:prop:neutral-stable-transition}. Their
pairs of base transformations and model transformations
generate a connected Lie group locally isomorphic to
$\mathcal H=G\times\R$; the second factor acts by scalar
multiplication. More explicitly, after a fixed linear
identification the pairs are $(g,c\sigma(g))$. The radical
pairs have $c=1$, while for the central parabolic element
$a_t$ the scalar is $e^{-t}$. Opposite radicals generate
the diagonal copy of $G$, and comparison with $a_t$ gives
the additional scalar factor. Finite covers account for
orientation conventions.

The stable transformations here are defined on the original suspension by
the specified measurable section. Their joint measurability, preservation
of the same original probability, and group laws for each fixed pair
almost everywhere follow from the canonical intrinsic stable conditionals
used in Proposition~\ref{II-loc:prop:neutral-stable-transition}.
Restriction to each original Levi-invariant component uses the semidirect
Mautner argument. Finite words are lifted from these quotient
transformations using the original section and their exact base maps.

Only local orbit charts are needed here; no global
factorization of the measurable transformations through
$\mathcal H$ is asserted. On a sufficiently small domain at
$x$, use the noncentered linear coordinates $w=(u,r)$, whose
marked base vector is $w_0=(0,1)$, and put
$\widetilde H_x(u,r)=H_x(u,1-r)$. A word with pair $(g,C)$
has the analytic orbit formula
\[
                         g\widetilde H_x(C^{-1}w_0).
\]
The marked vector is nonzero; using the linear origin here
would incorrectly remove the scalar orbit direction.
Compose the previously established chart identities only while
all intermediate points remain in their domains. This proves
the formula and its compatibility with changes of center on
that domain. Finite-word maps have smooth local sections
onto $\mathcal H$: use finitely many conjugates of the
generators whose tangent vectors span its Lie algebra, and
the inverse function theorem. The conjugating parameters can
be chosen small. A countable choice of such parameters and
of uniform chart bounds covers a conull set of centers.
Thus these formulas give coefficient-measurable analytic
local orbit charts. Changing the center of such a chart
changes its group coordinate by right translation.
Their rank is $\dim G+1$: the base supplies $\dim G$
directions and the scalar factor supplies a nonzero vertical
transverse direction. Analytic continuation may have
monodromy, which does not affect these local charts or their
right-translation changes of coordinates.

Lift the local orbit maps to $G\times M$, where their
formulas are equivariant under the suspension deck group.
Quotient there by the original left parabolic action; this
gives charts in $Q\backslash G\times M$, rather than a
topological quotient of the possibly dense parabolic orbits
in the suspension itself. Their dimension is $n$, and projection to $M$ is
a local analytic diffeomorphism: the stable derivatives
$D\phi_x$ and the scalar derivative span its tangent space.
The orbit coordinates are modeled on $Q\backslash\mathcal H$,
which identifies, up to a finite cover, with
$\R^n\setminus\{0\}$. For example, writing an element as
$(g,c)$, the vector $c^{-1}g^{-1}e_1$ is unchanged by the
left embedded parabolic $(q,a(q)^{-1})$. Right translations
act linearly on this vector. Hence these orbit charts carry
a well-defined flat affine connection and radial Euler field.

It remains to justify descent across different measured
orbits. The original parabolic and stable transformations
preserve $\mu$, as do their conjugates. In the local orbit
construction, a jointly Borel one-parameter family whose group
laws hold almost everywhere is used as a near-action. Such a
family is continuous in measure: the discussion following
Definition~0.1 in \cite[p.~2 of the arXiv version]{GTW2005}
identifies a near-action with a continuous homomorphism into
$\operatorname{Aut}(X,\mu)$ with its weak topology. In particular,
$\mu(T_t^{-1}A\mathbin{\triangle}A)\to0$ for every measurable
$A$ as $t\to0$. This assertion concerns the original
measure-preserving one-parameter families and does not require
their laws to hold on one common conull set for all parameters.
It supplies no global action of the model group $\mathcal H$.
We obtain the smooth measure class from finite words in these
original transformations. Choose finitely many conjugated
one-parameter directions whose model derivatives span the matrix
algebra. The resulting word map evaluated at the marked vector
is a submersion near the identity parameter. Exhaust the original
Haar lift by measurable sets of finite measure on which the word
and analytic-chart radii have common positive lower bounds.
On each such set, average the restricted measure under the actual
words over a sufficiently small parameter neighborhood.
Measure preservation and continuity in measure imply that the
original restricted measure is absolutely continuous with respect
to this average, and that the average is absolutely continuous
with respect to the original law.

The marked-vector submersion sends parameter Lebesgue measure
to Lebesgue class on a small coordinate ball. The literal word
identities identify physical evaluation with the original analytic
chart on that ball. Thus the marginal of the averaged measure
is equivalent to smooth measure restricted to an open subset
of $M$. Take the union of these open subsets over the exhausting
sets. The two absolute-continuity statements show that the
original marginal is equivalent to smooth measure on this open
set and gives its complement measure zero.
Lemma~\ref{II-loc:lem:unique-boundary-coordinate}
identifies the quotient measured space with $(M,\widehat\mu)$.
Transport the connection coefficients, the Euler field,
and the boundary coordinate through this isomorphism.
They are measurable fields on $M$.

For almost every orbit chart these fields agree with its
analytic model for Lebesgue-almost every point of its image,
by Fubini applied to the word identities and the marked-vector
submersion. Any two such charts therefore have the
same analytic connection, vector field, and boundary map
on their overlap: equality almost everywhere implies equality
of their analytic coefficients. Their open images give $W$;
take a countable subcover and its countable lattice saturation.
The preceding averaging measure-class argument gives
$\widehat\mu\sim\operatorname{Leb}|_W$.

A transition preserving the flat connection is affine.
Preservation of the Euler field removes its translation part.
After aligning the boundary coordinates, preservation of
projectivization makes its linear part scalar. This proves
the asserted compatibility, as well as the corresponding
formula for lattice transformations. The scalar direction
commutes locally with the parabolic and with suspension
descent; its descended vector field is therefore $E$ and
commutes with the lattice action. Its local flow preserves
the conditional probabilities $\mu_g$ wherever defined.
\end{proof}

\begin{theorem}[Description of regular neutral shear]
\label{II-loc:thm:neutral-regular-recognition}
Let $n\ge4$ and let $M$ be connected. In the regular-shear neutral case,
after finite covering and finite-index restriction, $M$ is
analytically the Hopf manifold model
\begin{equation}\label{II-loc:eq:recognized-neutral-circle}
 S^{n-1}\times\R/T\Z,\qquad
 \gamma(\theta,t)=
 \left(\frac{\gamma\theta}{\|\gamma\theta\|},
             t+\log\|\gamma\theta\|\right),
 \qquad T>0,
\end{equation}
or its dual version. In particular this case supplies actual
compatible analytic neighborhoods, rather than only moving germs.
\end{theorem}
\begin{proof}
Use $W,b,E$ from Proposition~\ref{II-loc:prop:neutral-radial-neighborhoods}.
The maximal local flow and its domain-translation law are given by
\cite[Theorem~9.12, pp.~212--214]{LeeSmoothManifolds2e}; its naturality under
lattice diffeomorphisms follows from
\cite[Proposition~9.13 and Corollary~9.14, p.~215]{LeeSmoothManifolds2e}.
For the analytic dependence of this same flow on time and initial point,
we use \cite[Appendix~B, Theorem~B.3(ii), pp.~333--336]{DuistermaatKolk};
its proof includes the real-analytic case explicitly. Maximality
identifies its flow with the prescribed flow of the same field,
with the same domain and values there.
We first prove almost-everywhere completeness, without
assuming that the local orbit construction integrated to a
global group action. The maximal forward existence time
$\tau_+(p)\in(0,\infty]$ of $E$ on $W$ is measurable and
$\Gamma$-invariant, because $E$ commutes with $\Gamma$.
Ergodicity makes it a constant almost everywhere. This
constant cannot be finite: on a positive-measure flow box,
flowing for a sufficiently small $s>0$ is nonsingular for
$\widehat\mu\sim\operatorname{Leb}|_W$ and replaces
$\tau_+$ by $\tau_+-s$. This contradicts its almost-everywhere
constancy. The backward time is infinite by the same argument.
Hence the local flow gives a global probability-preserving
flow modulo null sets. Local invariance and exhaustion by
flow domains show that it preserves each $\mu_g$ as well.

Almost every $\mu_g$ is supported on the analytic
one-manifold $b^{-1}(Qg)\cap W$. On it $E$ never vanishes.
Recurrence for this one-dimensional flow forces almost every
point onto a periodic orbit. Indeed a nonperiodic regular
orbit is an interval with a strictly monotone time coordinate,
and has no recurrent points.

The least positive period is invariant under $\Gamma$.
Ergodicity of $\widehat\mu$ makes it one constant $T>0$.
Near a periodic orbit the time-$T$ map is defined on an open
neighborhood. It is the identity almost everywhere there,
and hence analytically the identity. Since $b$ is a submersion
and constant along the flow, a return map on a transverse
disk is the identity; consequently these circles have trivial
local holonomy and persist as analytic circle fibers. An analytic section gives the circle chart $(q,[t])\mapsto\Psi_t(\sigma(q))$:
exact period gives injectivity, and $Db$ together with $E\ne0$ gives
invertible derivative. The inverse function theorem gives analytic
inverse charts at every time, as in
\cite[Theorem~1.11.4, pp.~53--55]{DuistermaatKolk}.
The
least period is locally constant, by the same return-disk
argument and its almost-everywhere value. Replace $W$ by the
union of these invariant neighborhoods. It is still invariant,
open, and of full measure, and all its fibers are circles of
least period $T$.

There are uniformly finitely many circles in almost every
fiber. To see this, the mass assigned by $\mu_g$ to the
circle through its point is a positive measurable
$P$-invariant function in the suspension. Countability of
the components of a one-manifold and ergodicity make this
mass $1/k$ for a fixed finite $k$. Equivalence with smooth
measure on $W$ implies that almost every boundary fiber has
exactly $k$ circle components in $W$, including all its
components, not just a selected family of positive-mass ones.

Each compact circle persists over a boundary neighborhood.
Thus the set of boundary points with at least $j$ circles
is open and $\Gamma$-invariant. There can be no point with
more than $k$, since that would give a positive-measure
open set with more than $k$. The set with at least $k$ is
nonempty; minimality of the lattice action on the compact
projective boundary makes this invariant open set the whole
boundary. Therefore $b:W\to Q\backslash G$ is a locally
trivial bundle with exactly $k$ compact circles in each fiber.
Its total space is compact. Since it is open in connected
$M$, it is all of $M$.

Develop the scalar atlas on the universal cover
\cite[Section~1.3, pp.~63--65]{Hatcher}
\cite[Theorem~5.2.2, pp.~120--121]{GoldmanDevelopment}.
Overlap uniqueness gives one translation center. After subtraction,
nonzero chart coordinates exclude zero, so \cite[Theorem~2, p.~580, proof p.~581]{Fried1980}
gives the analytic developing identification with $\R^n\setminus\{0\}$.
For its projection $q$, the curves $q(e^t v)$ solve the given Euler
equation; maximality identifies them with the prescribed flow. Thus the
positive scalar periods are exactly $\exp(T\Z)$. Hence
$P(\theta,[t])=q(e^t\theta)$ is a well-defined local analytic diffeomorphism
with compact source, and therefore a proper finite covering.

We construct the lattice equivariance using all individual lifts.
For each $\gamma$, pointed lifting to the universal cover
\cite[Propositions~1.33--1.34, pp.~61--62]{Hatcher}
and the local linear formula give, by connectedness, a scalar
$c\ne0$ with $q(c\gamma v)=\gamma q(v)$.
All such pairs $(\gamma,c)$ form a group. The positive pairs project
to a subgroup $\Gamma_0$ of index at most two in $\Gamma$.
Their positive kernel is $\exp(T\Z)$, so $\log c$ modulo $T$
defines a character $\Gamma_0\to\R/T\Z$.
Every finite-index subgroup of the original higher-rank lattice has
finite abelianization
by the preceding same-lattice property-$(T)$ argument and
\cite[Corollary~1.7.2, p.~67]{BHV}. Hence this character has finite
image and its kernel has finite index in $\Gamma_0$.
After restriction to the kernel of this character,
$q(\gamma v)=\gamma q(v)$. The covering $P$ is therefore equivariant
for the action in \eqref{II-loc:eq:recognized-neutral-circle}.
This constructs the finite covering and its lattice action together.
This proves \eqref{II-loc:eq:recognized-neutral-circle}.
\end{proof}

\subsection{Analytic normal characters for arbitrary lattices}

The scalar argument below uses only the hypersurface conclusion
of Theorem~\ref{II-loc:thm:neutral-atomic-recognition}.
None of these arguments uses a collar-extension result or
the global classification theorem. Write
\[
 S=S^{n-1},\qquad
 \gamma\cdot\theta=\frac{\gamma\theta}{\|\gamma\theta\|},\qquad
 \ell(\gamma,\theta)=\log\|\gamma\theta\|.
\]
All statements have dual versions by inverse transpose.

\begin{lemma}[Measurable reduction of boundary scalar cocycles]
\label{II-general:lem:boundary-scalar-measurable}
Let $\Gamma<\SL(n,\R)$ be any lattice, $n\ge4$.
Every real measurable additive cocycle $c$ over its action on
$S$, with the spherical Lebesgue measure class, has the form
\begin{equation}\label{II-general:eq:boundary-scalar-measurable}
 c(\gamma,\theta)
 =q\ell(\gamma,\theta)+f(\gamma\cdot\theta)-f(\theta)
\end{equation}
for a constant $q\in\R$ and a finite measurable function $f$.
The identity holds almost everywhere for each $\gamma$.
\end{lemma}
\begin{proof}
Let $Q^+$ be the positive-ray stabilizer of $e_1$, put
$\chi(q_0)=\log a(q_0)$ when $q_0e_1=a(q_0)e_1$, and set
$\theta_g=g^{-1}e_1/\|g^{-1}e_1\|$.
Consider the measurable affine-line bundle over $G/\Gamma$
with deck action
\[
 (g,s)\gamma=(g\gamma,s+c(\gamma^{-1},\theta_g)).
\]
The cocycle identity makes this a right action. Left
multiplication by $Q^+$ on the first coordinate commutes
with it, since $\theta_{q_0g}=\theta_g$.
A measurable trivialization of this bundle therefore gives
a real additive $Q^+$-cocycle over the Haar probability on
$G/\Gamma$.

We derive the needed reciprocity assertion from the original generic
induction theorem, rather than assume the two-homogeneous-space comparison.
The original theorem says that induction induces an injective map on
measurable cocycle cohomology (in fact a bijection)
\cite[Proposition~2.2 and its generality remark, author version p.~4]{IozziOrbit1991}.
We use its explicit product formula in
\cite[Proposition~4.4 and remark, author version p.~11; diagonal action p.~12]{IozziGeometric1992}.
For a standard Borel probability $G$-space $(S,\eta)$ with nonsingular
action, a closed subgroup $\Gamma$, and a Borel quotient section $r$, this
map sends an additive cocycle $c$ to
\[
 \beta(g;(z,x))
   =c\bigl(\lambda(g,x),r(x)^{-1}z\bigr),\qquad
 \lambda(g,x)=r(gx)^{-1}g r(x)\in\Gamma .
\]
Here $G$ acts diagonally on $S\times G/\Gamma$, and $G/\Gamma$ has the
original normalized Haar probability $\mu$. The source's normalized map is
$\theta(x)=r(x)^{-1}$, with $\theta(g\Gamma)g\in\Gamma$.
For additive coefficients its cohomology relation reads literally
$\beta_c=\beta_d+U(gz,gx)-U(z,x)$; injectivity gives
$c(\gamma,z)=d(\gamma,z)+u(\gamma z)-u(z)$.
This universal original class-injectivity assertion remains an explicit
classical assumption. It requires no integrability condition, and includes
no stabilizer slice or boundary-reciprocity conclusion.

Here is the internal comparison for the same original section and measures.
Write $z_0$ for the point stabilized by $Q$, and suppose that the slice
$\beta(q;(z_0,x))$ is $b(qx)-b(x)$ for each $q\in Q$ and
$\mu$-almost every $x$. For the continuous transitive action on the Hausdorff boundary used here,
a Borel orbit section $\rho:S\to G$ with $\rho(z)z_0=z$ is constructed from the original graph-uniformization clause
of Arsenin--Kunugui
\cite[Theorem~5.12.1, pp.~219--220; uniformization definition p.~185]{SrivastavaBorel}:
apply it to $\{(z,g):gz_0=z\}$. Its sections are nonempty and closed,
hence sigma compact in the locally compact second countable group $G$.
The original assertion supplies a single-valued Borel subgraph with the
same projection. The locally compact second countable group $G$ is Polish.
Taking the inverse of the first projection of the subgraph gives the
measurable section.

Define
\[
 A(z,x)=\beta\bigl(\rho(z);(z_0,\rho(z)^{-1}x)\bigr)
             +b\bigl(\rho(z)^{-1}x\bigr).
\]
For fixed $g,z$, the element
$k(g,z)=\rho(gz)^{-1}g\rho(z)$ lies in $Q$. The strict cocycle identity
and the slice equation for this element yield
$\beta(g;(z,x))=A(gz,gx)-A(z,x)$ for $\mu$-almost every $x$.
In this step preservation of the \emph{same} normalized quotient measure
transports the slice equation under $x\mapsto\rho(z)^{-1}x$.
Fubini then gives the product almost-everywhere equation for each $g$;
no product almost-everywhere equation is evaluated at the fixed point $z_0$.
The multiplication rule for return labels and the section-coordinate
identities identify $\beta$ with the induced cocycle. This cocycle is
jointly measurable, and the action preserves the quotient probability.
Apply injectivity
of induction on cohomology with $d=0$ to obtain the lattice transfer.
For a nonzero s-finite boundary measure, pass to an equivalent
probability $\eta$. Nonsingularity and the resulting almost-everywhere
identity are preserved under this change of measure.

We spell out its scalar reduction. The additive part of
Proposition~\ref{II-loc:prop:levi} makes the cocycle zero on
$L=\SL(n-1,\R)$. The skew-translation representation and
Mautner argument in the proof of
Theorem~\ref{II-loc:thm:relative} then make it zero on the
unipotent radical $U$: the vector $e^{-s^2}$ is $L$-fixed,
and every element of $U$ is contracted by a sequence in $L$.
The residual cocycle on the one-dimensional central subgroup
is $L$-invariant. Howe--Moore makes $L$ ergodic on $G/\Gamma$,
so it is constant in the base variable. Its parameter is a
measurable additive function of one real variable, hence
$q\chi$ for a constant $q$.

The auxiliary $\SL(n-1,\Z)$ used in the proof of
Proposition~\ref{II-loc:prop:levi} is a lattice in this continuous
suspension Levi $L$. It is not required to be a subgroup of
the original $\Gamma$. Thus this use of property $(T)$
imposes no arithmetic assumption on $\Gamma$.

We finish with the same return labels, writing
$\operatorname{Ind}(c)(q_0,x)=\beta(q_0;(z_0,x))$. Put
$\ell(\gamma,\theta)=\log\|\gamma\theta\|$ and
$R(g)=-\log\|g^{-1}e_1\|$. The section equation gives the pointwise identity
\[
 \operatorname{Ind}(\ell)(q_0,x)
   =\chi(q_0)+R(r(x))-R(r(q_0x)).
\]
Thus, after the preceding parabolic reduction
$\operatorname{Ind}(c)(q_0,x)=q\chi(q_0)+b(q_0x)-b(x)$,
the induced cocycle of $c-q\ell$ is the coboundary of
$B(x)=b(x)+qR(r(x))$. Apply the preceding induction argument to the
strict cocycle obtained by measurable rectification. It gives a measurable $f$ on the original
boundary, with respect to its original spherical measure class, satisfying
$c-q\ell=f(\gamma\theta)-f(\theta)$. Restoring the original
almost-everywhere representative proves~\eqref{II-general:eq:boundary-scalar-measurable}.
\end{proof}

\begin{theorem}[Analytic scalar cocycles on the projective boundary]
\label{II-general:thm:boundary-scalar-analytic}
Let $\Gamma<\SL(n,\R)$ be any lattice, $n\ge4$.
Suppose $c:\Gamma\times S\to\R$ is an additive cocycle and
$c(\gamma,\cdot)$ is real analytic for every $\gamma$.
There are a unique $q\in\R$ and a real analytic $f:S\to\R$
such that
\begin{equation}\label{II-general:eq:boundary-scalar-analytic}
 c(\gamma,\theta)=q\log\|\gamma\theta\|
                  +f(\gamma\cdot\theta)-f(\theta)
\end{equation}
for every $\gamma$ and every $\theta$. The function $f$ is
unique up to an additive constant. No finite-index restriction
is needed in this conclusion.
\end{theorem}
\begin{proof}
Use Lemma~\ref{II-general:lem:boundary-scalar-measurable} and set
$c_0=c-q\ell$. For a fixed $T>0$, form the compact analytic
$n$-manifold and analytic action
\[
 M_0=S\times\R/T\Z,\qquad
 \gamma(\theta,s)=(\gamma\cdot\theta,s+c_0(\gamma,\theta)).
\]
The measurable transfer from that lemma gives an invariant
measurable graph $s=f(\theta)\pmod T$.
In the suspension its graph defines the section
\[
 g\Gamma\longmapsto[g,(\theta_g,f(\theta_g)\bmod T)].
\]
It is well-defined modulo null sets by the cocycle equation.
Pushing Haar probability forward gives a $Q^+$-invariant
probability with one-point fiber conditionals. It is ergodic,
since its action is measurably the action on $G/\Gamma$ and
the Levi is ergodic there. The usual finite angular orientation
convention, or averaging over its finitely many translates,
gives the corresponding minimal-parabolic measure when needed.
Projection to the angular suspension identifies its parabolic
stabilizer with the line parabolic, up to those finite components.

The extra fiber derivative is exactly $1$. The derivative on
the angular quotient is the standard projective derivative;
hence the fiber Lyapunov weights are $\beta_2,\ldots,\beta_n,0$.
Off-diagonal derivatives of $c_0$ do not change this spectrum,
because the vertical circle line is invariant. More precisely, the
nonzero commuting vertical vector field has zero growth on the same
selected component. The original canonical Haar graph class and its finite atomic
conditionals therefore satisfy the atomic recognition theorem on
this same graph measure. Its canonical-support finite-cover
construction applies. This is the atomic neutral branch of
Theorem~\ref{II-loc:thm:neutral-atomic-recognition}.
Its hypersurface conclusion supplies a compact embedded
analytic projective covering hypersurface $N'$ in $M_0$,
invariant under a finite-index subgroup $\Gamma_0$.

Let $\pi:M_0\to S$ be projection. An equivariant spherical
parametrization $h:S\to N'$ is available after finite index,
with either the standard or the dual label. Apply
Lemma~\ref{II-general:lem:angular} to $\pi h$. The dual label
is impossible, and in the standard case $\pi h$ is the
identity or antipodal map. In particular $h$ cannot be a
nontrivial covering. Thus $\pi|_{N'}$ is an analytic
diffeomorphism and $N'$ is the graph of an analytic map
$a:S\to\R/T\Z$.

Covering lifting \cite[Propositions~1.14 and~1.33]{Hatcher}
gives a continuous lift $\widetilde a:S\to\R$, since the sphere
is simply connected. The circle exponential has local analytic
inverses, so $\widetilde a$ is real analytic. Invariance of its graph gives
\[
 c_0(\gamma,\theta)
 -\widetilde a(\gamma\cdot\theta)+\widetilde a(\theta)
 \in T\Z\qquad(\gamma\in\Gamma_0).
\]
The left side is continuous on connected $S$, hence constant
in $\theta$. These constants form a homomorphism
$\Gamma_0\to T\Z$, which is zero by finite abelianization
\cite[Theorems~1.4.15 and~1.7.1, Corollary~1.7.2,
pp.~54, 65 and 67 in the author manuscript]{BHV}.
This proves \eqref{II-general:eq:boundary-scalar-analytic} on
$\Gamma_0$, with $f=\widetilde a$.

Replace $\Gamma_0$ by its finite-index normal core. For
$\gamma\in\Gamma$, the residual cocycle
\[
 n(\gamma,\theta)=c_0(\gamma,\theta)
                  -f(\gamma\cdot\theta)+f(\theta)
\]
satisfies $n(\gamma,\delta\cdot\theta)=n(\gamma,\theta)$
for every $\delta\in\Gamma_0$, by normality and the cocycle
identity. Minimality makes it constant in $\theta$.
These constants form a homomorphism $\Gamma\to\R$ and
therefore vanish. This proves the formula for all of $\Gamma$.

If two analytic formulas had different $q$, evaluate their
difference at an attracting fixed ray of a regular lattice
element with eigenvalue greater than $1$. Such elements
exist by Lemma~\ref{II-general:lem:cartan-returns}.
The coboundary terms vanish there and $\ell\ne0$, so the
two values of $q$ agree. The difference of the two transfers
is then continuous and invariant; minimality makes it constant.
The dual version follows by transporting the original lattice through
inverse transpose. If the cocycle law is initially given almost
everywhere for a measure equivalent to spherical Lebesgue measure,
continuity and full support make it hold at every angular point.
\end{proof}

\begin{corollary}[The normal character, including periodic points]
\label{II-general:cor:normal-character}
Let an arbitrary lattice $\Gamma<\SL(n,\R)$, $n\ge4$, act
analytically near a compact invariant analytic hypersurface
with its standard projective action. On the spherical angular
cover, and after removing the normal-orientation sign by finite
index, there is an analytic normal frame in which
\begin{equation}\label{II-general:eq:normal-character}
 D^\perp_\theta\alpha(\gamma)=\|\gamma\theta\|^q
\end{equation}
for every $\theta$ and $\gamma$, for one constant $q$.
This constant is the coefficient of the normal suspension
weight $q\varepsilon_1$. In particular, at a common eigenray
$[e_1]$ of a positive split Cartan subgroup, its periodic
normal character is $e^{q\varepsilon_1}$.
\end{corollary}
\begin{proof}
The normal line pulled back to the sphere is orientable and
trivial. A real analytic bundle metric exists by applying
\cite[Theorem~3, p.~470]{Grauert} to the total space of the
analytic normal line bundle: pull back the Euclidean metric
under the resulting analytic immersion and restrict it to
the vertical tangent space along the zero section. This gives
positive analytic local coefficients whose changes of
trivialization obey the square of the line transition function.
More explicitly, for the original inverse trivializations $\psi_i$ and
one regular analytic embedding $e$, put
\[
 v_i(x)=\left.\frac{\partial}{\partial t}e(\psi_i(x,t))\right|_{t=0},
 \qquad h_i(x)=\sum_a v_i(x)_a^2.
\]
The vertical derivative is nonzero by regularity, so $h_i>0$;
analytic differentiation makes $h_i$ analytic. The original transition
$\psi_i(x,t)=\psi_j(x,t_{ij}(x)t)$ gives
$v_i=t_{ij}v_j$ by the vertical chain rule and hence
$h_i=t_{ij}^2h_j$, on the same original line atlas.
The regular/proper terminology is that of \cite[Section~1.2, p.~462]{Grauert}.
The theorem applies also to noncompact total spaces with
countable topology. A countable subcover of the line trivializations
makes the total space second countable; Hausdorffness follows from
the Hausdorff base and the local trivializations. The metric requires
no orientability: an orientation section is obtained separately on
the simply connected sphere. The positively oriented unit section is
therefore an analytic frame. In this frame the nonzero normal
derivative $a$ is a multiplicative analytic cocycle, and
$c=\log|a|$ satisfies
Theorem~\ref{II-general:thm:boundary-scalar-analytic} for the
original full lattice. The sign of each lattice element is
constant on the connected sphere and defines a homomorphism
to $\{\pm1\}$. Restrict to its finite-index kernel, where $a$
is positive. Multiplying the frame by $e^f$ gives
\eqref{II-general:eq:normal-character}.

Choose one original DLS section $s$ of $G/\Gamma$ with its
normalized Haar probability and integrable return-word lengths, as in
Lemma~\ref{II-power:lem:section}. Work first in the original unit
normal frame constructed above. For each generator, $\log|a(\gamma,\theta)|$
is continuous on the compact sphere and hence bounded. The scalar
cocycle identity then bounds its absolute value by a constant times
the word length of $\gamma$. Thus the actual one-step normal logarithmic
derivative, evaluated on the original section return, is integrable.

For a fixed positive Cartan element $p$, put
$\theta_x=s(x)^{-1}\cdot[e_1]$ and
$\delta(p,x)=s(px)^{-1}ps(x)$. The exact section identity and the
classified scalar formula give the actual $k$-step logarithmic growth
\[
 kq\chi(p)+b(p^kx)-b(x),\qquad
 b(x)=f(\theta_x)-q r(s(x)),\quad
 r(g)=-\log\|g^{-1}e_1\|.
\]
The one-step increment of this same measurable $b$ is integrable by
the preceding moment bound. Truncation and invariant ergodic averages
applied to $|b|$ give $b(p^kx)/k\to0$ almost everywhere.
Hence the actual normal Lyapunov weight exists and equals
$q\varepsilon_1$ for each fixed Cartan element, on this same Haar
probability and section. A positive measurable change of norm has
this same exponent whenever its finite limit exists: recurrence to bounded
level sets of the norm ratio identifies that limit.
At a positive eigenray, $\|\gamma\theta\|$ is the corresponding
eigenvalue, which proves the assertion about periodic points.
\end{proof}

\begin{proposition}[The first nonlinear normal coefficient in the neutral case]
\label{II-general:prop:neutral-leading-normal}
Suppose $q=0$ in Corollary~\ref{II-general:cor:normal-character}.
Use an analytic defining coordinate $t$ with normalized normal
derivative $1$, and any analytic tubular parametrization by
$(\theta,t)$; equivariance of its angular projection is not
assumed. If the normal coordinate of the action is
\[
 t\bigl(\alpha(\gamma)(\theta,t)\bigr)
 =t+c_k(\gamma,\theta)t^k+O(t^{k+1}),\qquad k\ge2,
\]
with all degrees $2,\ldots,k-1$ zero, then an analytic
coordinate change $t\mapsto t+f_k(\theta)t^k$ puts it in
the form
\begin{equation}\label{II-general:eq:neutral-leading-normal}
 t\bigl(\alpha(\gamma)(\theta,t)\bigr)
 =t+\kappa_k\log\|\gamma\theta\|\,t^k+O(t^{k+1})
\end{equation}
for a constant $\kappa_k\in\R$. Thus either a first nonzero
constant $\kappa_k$ is reached after finitely many analytic
changes, or there is a formal invariant defining function
$T=t+O(t^2)$ with analytic coefficient functions on the
hypersurface. The latter alternative asserts formal invariance,
not convergence.
\end{proposition}
\begin{proof}
The angular coordinate of $\alpha(\delta)(\theta,t)$ is
$\delta\cdot\theta+O(t)$. In the degree-$k$ normal
composition equation this displacement changes
$c_k(\gamma,\delta\cdot\theta)t^k$ only by $O(t^{k+1})$.
These remainder statements follow from the finite-coordinate Hadamard identity
\cite[Lemma~6.2, p.~614]{TsuboiAnalyticDiffeomorphisms}, applied with angular coordinates fixed.
The jointly analytic quotient equals the divided difference for $t\ne0$;
continuity gives its derivative value at zero. Thus this canonical quotient
is analytic in each angular chart and defines one $R$ along the zero section,
with $a(\theta,t)-a(\theta,0)=tR(\theta,t)$.
The constant linear normal coefficient and the vanishing
lower normal degrees therefore give
\[
 c_k(\gamma\delta,\theta)
 =c_k(\gamma,\delta\cdot\theta)+c_k(\delta,\theta).
\]
Theorem~\ref{II-general:thm:boundary-scalar-analytic} writes
this as $\kappa_k\ell$ plus an analytic coboundary.
Under the change $T=t+f_k(\theta)t^k$ the coefficient becomes
$c_k+f_k\circ\gamma-f_k$. Choose $f_k$ to cancel that
coboundary. Lower normal coefficients remain zero.
If $\kappa_k=0$, continue to the next degree. Infinite
iteration is well-defined formally, since each step leaves
all previously fixed coefficients unchanged, and gives
the stated formal invariant defining function.
\end{proof}

This is a restriction from the full lattice even without
an equivariant collar projection. At a periodic projective
point the first normalized nonlinear normal coefficient
is $\kappa_k\varepsilon_1(\log\gamma)$; it cannot be an
independently chosen Cartan character. Constructing a
convergent transverse foliation still requires more than
this finite-order restriction.

\begin{proposition}[Neutral Cartan germs need not have a transverse analytic axis]
\label{II-general:prop:neutral-cartan-obstruction}
Even a full-rank Cartan group of commuting analytic germs
with the prescribed projective action on an invariant
hypersurface and trivial normal derivative need not preserve
an analytic curve transverse to that hypersurface.
This assertion concerns Cartan germs, not an action of the
full lattice.
\end{proposition}
\begin{proof}
Write the tangential coordinates as $u_1,\ldots,u_m$, where
$m=n-1$, and let $N=\{t=0\}$. The analytic vector fields
\[
 X_1=(u_1+t^2)\partial_{u_1}+t^2\partial_t,
 \qquad X_i=u_i\partial_{u_i}\quad(2\le i\le m)
\]
commute. Their flows give an analytic germ action of $\R^m$.
Identify its parameters with the independent projective
Cartan weights $\beta_2,\ldots,\beta_n$. On $N$ the action
is exactly the diagonal projective Cartan action, and its
normal derivative is identically $1$.

Restrict the parameter group to any full lattice. Since its real
span is the whole parameter space, it has an element with
$a_1\ne0$. Let $F$ be the induced planar map in $(u_1,t)$.
It preserves the planar vector field
$(u_1+t^2)\partial_{u_1}+t^2\partial_t$, and
$D_0F(1,0)=(e^{a_1},0)$, while its normal derivative is $1$.

Suppose there were an invariant transverse analytic curve.
Inverting its normal coordinate writes its planar projection
as $(h(t),t)$. Its induced parameter map $\tau$ is analytic,
fixes zero, and satisfies $\tau'(0)=1$. Differentiating graph
invariance and using preservation of the planar vector field
shows that the analytic tangency defect
\[
 R(t)=h(t)+t^2-t^2h'(t)
\]
satisfies $R(\tau(t))=k(t)R(t)$, with $k$ analytic and
$k(0)=e^{a_1}\ne1$. More explicitly, if
$B(t)=D_{(h(t),t)}F$, then
$k(t)=(B(t)(1,0))_1-h'(\tau(t))(B(t)(1,0))_2$.
If $R$ were not identically zero as a germ, write
$R(t)=t^d r(t)$ with $r(0)\ne0$. Divide the covariance identity
by $t^d$ and let $t\to0$. Since $\tau(t)/t\to1$, this gives
$r(0)=k(0)r(0)$, a contradiction. Thus
\[
 t^2h'(t)=h(t)+t^2.
\]
The Taylor recurrence forces
\[
 h(t)=-\sum_{d\ge2}(d-1)!\,t^d,
\]
which has zero convergence radius. No such transverse analytic
curve exists. The linear change $x=t$, $y=t-u_1$ identifies
the planar vector field with Euler's divergent central-manifold
example \cite[Example~7.10, p.~111]{IlyashenkoYakovenko}.
The defect argument above establishes the additional assertion
for the discrete Cartan lattice.
\end{proof}

The example does not obstruct the classification theorem:
it supplies no compatible germs for the other lattice elements.
It shows why periodic normal characters and the availability
of many commuting regular elements cannot alone settle the
neutral collar. Compatibility with the full lattice must enter
before asserting a transverse analytic foliation or flow.

\section{Contracting neighborhoods from the suspension}
\label{II-sec:suspension-negative-collar}

We now construct an equivariant analytic projection and a
compatible normal Euler field. The latter gives a direct
linearization in the fibers of the projection. Throughout this
section $N$ is an invariant compact analytic projective
covering hypersurface, already identified by the hypersurface
conclusion of Theorem~\ref{II-loc:thm:negative-recognition}, and
its normal weight is $q\varepsilon_1$ with $q<0$.
The ambient derivative spectrum is the union, with multiplicity,
of the tangent spectrum on $N$ and the quotient normal spectrum.
The ordinary exact-sequence property is stated in
\cite[Remark~2.2.9(iii), pp.~6--7 of the arXiv version]{FilipMETv1};
the natural subbundle and quotient norms inherit the moment condition
\cite[Remark~2.3.2, p.~8 of the arXiv version]{FilipMETv1}.
The normal character is that of the actual quotient derivative, as
identified in Corollary~\ref{II-general:cor:normal-character}.
Let $W_x$ and $\mathsf n_x$ denote the measurable tangent injection and
normal map in the original derivative exact sequence, with prescribed
characters $\chi_j$ and $\chi_{\mathrm n}$. For these maps, we use the already constructed
ambient Lyapunov field rather than applying the theorem after an
uncontrolled change of norm. Fix a Cartan element $a$. For every
nonzero tangent column $W_xe_j$, covariance gives
\[
 \alpha(na,x)W_xe_j=e^{n\chi_j(a)}W_{a^nx}e_j.
\]
The ambient field gives a growth limit for the left side. Recurrence to
sets where the nonzero measurable column has norm bounded above and
away from zero identifies this limit with $\chi_j(a)$. Applying the same
argument to $-a$ shows that every nonzero coefficient of $W_xe_j$ in the
ambient Lyapunov basis has weight $\chi_j(a)$.
For an ambient basis vector $v$ with $\mathsf n_xv\ne0$, normal
covariance reads
\[
 \mathsf n_{a^nx}\alpha(na,x)v
       =e^{n\chi_{\mathrm n}(a)}\mathsf n_xv.
\]
Recurrence to a set where $\|\mathsf n_{a^nx}\|$ is bounded gives
$\chi_{\mathrm n}(a)\le\lambda_v(a)$. The argument for $-a$ gives the
reverse inequality. Thus both maps in the given exact sequence
commute with every threshold projector. Restricting the exact sequence
to the corresponding growth subspaces and adding dimensions identifies
the number of ambient weights below each threshold with the tangent count
plus the normal count. These counts determine the union of the prescribed
tangent and normal exponents, including multiplicities.
Only finite-measure recurrence and the existing ambient growth limits
were used. No logarithmic moment or full-orbit temperedness of the
prescribed moving frame is asserted. The full-measure set may depend
on the fixed element $a$.
Apply this argument on a countable dense set of Cartan elements and
intersect the corresponding countably many conull sets. At each remaining
base point, the finitely many possible permutations define closed linear
subspaces whose union contains that dense set. The union therefore covers
the Cartan space. A real vector space cannot be a finite union of proper
linear subspaces, so one of these subspaces is the whole space. Thus, at
almost every base point, one permutation works for all Cartan directions.
No ergodic decomposition is needed.

Here $\Gamma<\SL(n,\R)$ is arbitrary and $n\ge4$.
Pass to the spherical angular cover of the collar germ and
remove the normal-orientation sign by finite index. These
are covers of the collar germ, not an assertion about a cover
of all of $M$. Put $m=n-1$, $\lambda=e^{-n/m}$, and $b=e^q$.

The subgroup compatibility needed here is automatic in the
suspension. For any $h\in G$ and $\gamma\in\Gamma$,
\[
 h\cdot\bigl((g,x)\cdot\gamma\bigr)
  =(hg\gamma,\alpha(\gamma)^{-1}x)
  =\bigl(h\cdot(g,x)\bigr)\cdot\gamma.
\]
Thus every minimal parabolic $P$, and every line or hyperplane
parabolic $Q\supset P$, acts on the suspension for every lattice.
The relevant invariant probability is a further piece of data;
in this argument it comes from the recognized hypersurface.

There is a canonical probability space $X=G/\Gamma$ over this
hypersurface: its point $x_g$ is represented in the suspension
by $[g,\theta_g]$, where $\theta_g=g^{-1}e_1/\|g^{-1}e_1\|$
is interpreted in $N$. The parabolic action on this section
is the original suspension action. Opposite-radical transport
along the standard projective hypersurface defines
\[
 \widehat v x=v\cdot\phi_x(v),\qquad v\in V\simeq\R^m.
\]
Explicitly, take the column matrix and angular chart to be
\[
 v=\begin{pmatrix}1&0\\-v&I_m\end{pmatrix},\qquad
 \phi_g(u)=\frac{g^{-1}(1,u)}{\|g^{-1}(1,u)\|},
\]
using the same letter for a column and its group element.
Then $\theta_{vg}=\phi_g(v)$ and
$\phi_{vg}(u)=\phi_g(u+v)$. On $X$ these transformations
are the ordinary left $V$-action on $G/\Gamma$. In particular
they preserve Haar probability and satisfy the usual relations
with the left $Q^+$-action. No subgroup of $\Gamma$ fixing a
projective point is used here.

\begin{lemma}[Calibrated contracting charts along the recognized hypersurface]
\label{II-suspension:lem:calibrated-negative}
There are measurable analytic full stable charts $H_x$ with
\[
 H_x(u,0)=\phi_x(u),\qquad
 H_{lx}^{-1}F(l,x)H_x(u,z)=(lu,z),
\]
and
\[
 H_{a_tx}^{-1}F(a_t,x)H_x(u,z)
       =(\lambda^t u,b^t z).
\]
They extend to analytic diffeomorphisms from $\R^m\times\R$
onto the open fiberwise stable manifold of their centers.
For each fixed $q_0\in Q^+$, the normalized map
$H_{q_0x}^{-1}F(q_0,x)H_x$ has the form
$(u,z)\mapsto(h_{q_0,x}(u),c_{q_0,x}(u)z)$.
All assertions concern germs where necessary, with positive
measurable radii.
\end{lemma}
\begin{proof}
Apply the finite-jet reduction and contracting construction
of Corollary~\ref{II-loc:cor:negative-curves} to the canonical
measure on $X$. Its Levi action is ergodic by Howe--Moore.
The invariant analytic hypersurface $N$ equals the field
$S_x$ in Proposition~\ref{II-loc:prop:negative-tangential}:
the uniqueness observation after that proposition applies
to the germ of $N$ at every typical center.

We explain why the finite signs can now be removed on this
probability space. The tangent bundle of $N$ has its canonical
homogeneous projective frame. Compare it with the standard
block of $DH_x(0)$. The comparison matrices $K_x$ satisfy
$K_{lx}=\varepsilon(l,x)lK_xl^{-1}$. Recurrence under each
unipotent element of a countable generating set forces $K_x$
to commute with that element: otherwise its conjugates have
polynomial growth, unaffected by the sign. Hence $K_x$ is
scalar. Its absolute value is constant by Levi ergodicity,
and its sign removes the tangent sign cocycle.
For the normal quotient, Corollary~\ref{II-general:cor:normal-character}
gives a frame whose induced Levi multiplier is $1$ and whose
central multiplier is $b^t$. The nonzero normal component
of $DH_x(0)$, compared with this frame, similarly removes
the normal signs and fixes its scale. These are block-scalar
changes with constant absolute sizes, so they preserve the
contraction construction and its tempered domains.
The given tangent and normal frames also make the central
multipliers exactly $\lambda^t$ and $b^t$.

After this calibration $H_x^{-1}\phi_x$ has derivative the
identity and intertwines multiplication by $\lambda$.
For its degree-$j$ coefficient, central covariance has factor
$\lambda^{1-j}$. Recurrence to bounded coefficient sets
eliminates every degree $j\ge2$. Thus $H_x(u,0)=\phi_x(u)$.
Tempering and continuation by contraction give the full
stable charts as in \eqref{II-loc:eq:global-stable-chart}.
Finally Proposition~\ref{II-loc:prop:negative-tangential} gives
the radical form $(A(vu)u,C(vu)z)$. The calibrating block
scalars preserve independence of $z$ in the first component
and linearity in $z$ in the second.
The same assertions are immediate for $L$ and $A_0$;
these subgroups and the radical generate $Q^+$.
\end{proof}

\begin{proposition}[Exact changes of center in the contracting branch]
\label{II-suspension:prop:negative-transitions}
For each fixed $v\in V$, almost everywhere in $x$ the
calibrated charts satisfy
\begin{equation}\label{II-suspension:eq:negative-transitions}
 H_x^{-1}v^{-1}H_{\widehat v x}(u,z)=(u+v,z).
\end{equation}
The inverse uses the branch through $\phi_x(v)$. This is
an identity of analytic neighborhoods, including their
transverse variables.
\end{proposition}
\begin{proof}
Write the left side as $P_v(x)$. The full stable chart makes
it an analytic germ at $(0,0)$ with value $(v,0)$.
If $l$ fixes the column $v$, Levi covariance gives
\[
 P_v(lx)=\rho(l)P_v(x)\rho(l)^{-1},
 \qquad \rho(l)(u,z)=(lu,z).
\]
For a unipotent $l$ in this column stabilizer, the induced
action on each finite Taylor coefficient space is unipotent.
Poincare recurrence forces the coefficient vector to be
fixed by this action; a nonconstant unipotent orbit is
polynomial and cannot return to a bounded neighborhood.
Consequently the coefficients are invariant under $l$.
Its cyclic action on $G/\Gamma$ is ergodic by Howe--Moore,
so each coefficient is constant almost everywhere.
The resulting deterministic germ commutes with the whole
column stabilizer, by covariance.

Since $m\ge3$, that stabilizer has no nonconstant scalar
polynomial invariants in $u$, and its equivariant vector
polynomials are spanned by $u$ and $v$. For example, at
$v=e_1$ it acts by
$(u_1,u')\mapsto(u_1+au',Bu')$, with $B\in\SL(m-1,\R)$.
Translation in $u_1$ and the standard invariant theory on
$u'$ give precisely those assertions. Degree by degree,
and then by analyticity, we obtain
\[
 P_v(u,z)=(A(z)u+B(z)v,C(z)).
\]
Levi transitivity on nonzero columns makes $A,B,C$ independent
of $v\ne0$. The tangential calibration gives $A(0)=B(0)=1$
and $C(0)=0$; invertibility gives $C'(0)\ne0$.

The group law of $\widehat V$ gives
\[
 P_v(x)\circ P_w(\widehat v x)=P_{v+w}(x).
\]
For two independent columns this implies $C\circ C=C$,
hence $C(z)=z$, and then $A(z)^2=A(z)$, hence $A(z)=1$.
These are legitimate analytic identities, not merely
composition of unrelated germs: restrict to a positive-measure
set where $H_x,H_x^{-1}$ have uniform positive domains and
bounds. Also restrict to a compact subset of one suspension
trivialization, so small opposite-radical transports have
uniform analytic bounds. Choose $v,w$ small. Strong continuity of left
translations on $L^2(G/\Gamma)$ ensures a positive-measure
intersection of the needed translates of that set. All
intermediate charts are then defined. Determinism of $A,B,C$
extends the resulting germ identities to the preceding formula.

Central covariance is
\[
 P_{\lambda v}(a_1x)D=DP_v(x),\qquad D(u,z)=(\lambda u,bz).
\]
It therefore says $B(bz)=B(z)$. Since $0<b<1$ and $B$ is
analytic, all positive-degree coefficients vanish. As
$B(0)=1$, we have $B=1$, proving
\eqref{II-suspension:eq:negative-transitions}.
This last contraction is the step unavailable when $q=0$.
\end{proof}

\begin{lemma}[Descent of analytic retraction germs]
\label{II-suspension:lem:retraction-descent}
Let a countable group act analytically on $M$ preserving a
compact connected analytic submanifold $N$, minimally on $N$.
Suppose $\pi_y$ is a measurable equivariant field of analytic
retraction germs onto $N$, based at $y\in N$, defined almost
everywhere for a smooth measure class on $N$.
Assume there are a nonempty open set $U\subset N$ and an
analytic retraction $\Pi$ on a neighborhood of $U$ whose
germ is $\pi_y$ for almost every $y\in U$.
Then these germs extend to an analytic equivariant retraction
on a neighborhood of $N$, with equivariance understood as
germs along $N$.
If an equivariant measurable field of analytic vector-field
germs $\mathcal E_y$ is represented on the same patch by
one analytic vector field, it descends as well.
\end{lemma}
\begin{proof}
Minimality makes the translates of $U$ cover $N$; choose
a finite subcover. Translate $\Pi$ to the corresponding
ambient neighborhoods. On an angular overlap their germs
coincide for almost every point of $N$, by equivariance of
the measurable field and nonsingularity of analytic
diffeomorphisms for the smooth measure class.

They coincide as germs at every point of that overlap.
Indeed, in a small connected common ambient neighborhood
of such a point, choose a nearby generic point of $N$.
Equality of germs there gives equality on a nonempty
ambient open subset. The analytic identity theorem gives
equality on the connected neighborhood. This uses equality
of entire germs, not just equality of their restrictions to $N$.
The finitely many representatives can consequently be shrunk
and glued near compact $N$. The same argument, for each
group element, proves equivariance everywhere as germs.
The resulting map restricts to the identity on $N$ and is
a submersion near $N$, since each representative is a retraction.
The proof for vector fields is identical in analytic local
coordinates: equality of entire germs propagates on every
connected overlap. It also preserves any analytic identity
between the retraction and the vector field.
\end{proof}

\begin{lemma}[An invariant normal Euler field linearizes the collar]
\label{II-suspension:lem:euler-collar}
Let $\Pi$ be an analytic equivariant retraction onto a compact
analytic hypersurface $N$. Suppose a $\Gamma$-invariant analytic
vector field $\mathcal E$ is vertical for $\Pi$, vanishes on $N$,
and induces the identity on the normal line at every point of $N$.
On an oriented collar, any prescribed analytic normal frame
extends to analytic fiber coordinates in which
$\mathcal E=t\partial_t$. In these coordinates the action is
$(\theta,t)\mapsto(\gamma\theta,a(\gamma,\theta)t)$, where
$a$ is its derivative in the prescribed normal frame.
\end{lemma}
\begin{proof}
Choose an analytic Riemannian metric
\cite[Section~42.4, p.~441]{KrieglMichorConvenient1997}.
Its local geodesic exponential is analytic and has the defining initial velocity
\cite[Chapter~III, Section~8, pp.~147--148, and Chapter~VI, Section~6, p.~255]{KobayashiNomizuI}.
For the prescribed vertical frame $v$, put
$F(\theta,t)=\exp_\theta(tv_\theta)$ and
$G(\theta,t)=(\Pi(F(\theta,t)),t)$.
Since $d\Pi(v)=0$, $dG$ is the identity along the zero section.
The analytic inverse theorem and compactness give $F\circ G^{-1}$
on a uniform positive strip, with projection $\Pi$ and the prescribed vertical derivative.
In its coordinates write
\[
 \mathcal E=t\,e(\theta,t)\partial_t,\qquad e(\theta,0)=1.
\]
The integrand below has a removable singularity at zero, jointly
in the base coordinates: apply the finite-coordinate Hadamard identity
\cite[Lemma~6.2, p.~614]{TsuboiAnalyticDiffeomorphisms} with base coordinates fixed.
The analytic quotient equals the divided difference, including
its derivative value at zero. Thus
\[
 T(\theta,t)=t\exp\left(\int_0^t
       \frac{e(\theta,s)^{-1}-1}{s}\,ds\right)
\]
is jointly analytic with $\partial_tT(\theta,0)=1$.
Direct differentiation gives $\mathcal E T=T$.
In this coordinate a fiber map commuting with $T\partial_T$
satisfies $T\partial_TF=F$, hence $F=aT$. Its derivative is
the one prescribed, because the normalizing change is tangent
to the identity along $N$.
\end{proof}

\begin{theorem}[Full contracting collars for arbitrary lattices]
\label{II-suspension:thm:negative-collar}
Let $\Gamma<\SL(n,\R)$ be any lattice, $n\ge4$, acting
analytically on a compact manifold. Let $N$ be an invariant
compact embedded analytic projective covering hypersurface
with its standard or dual action and with normal coefficient
$q<0$ from Corollary~\ref{II-general:cor:normal-character}.
After the indicated angular and normal-orientation covers
and finite-index restriction, the full collar germ is
analytically conjugate to
\begin{equation}\label{II-suspension:eq:negative-collar}
 \gamma(\theta,t)=
 \left(\frac{\sigma(\gamma)\theta}{\|\sigma(\gamma)\theta\|},
             \|\sigma(\gamma)\theta\|^q t\right).
\end{equation}
In particular no arithmetic hypothesis is needed for the
full neighborhood description of the contracting branch.
\end{theorem}
\begin{proof}
We prove the standard case. Use the calibrated charts of
Lemma~\ref{II-suspension:lem:calibrated-negative}. In the
trivialization over a lift $g\in G$, define the analytic
retraction and vector field
\[
 \Pi_g=\phi_g\circ\operatorname{pr}_u\circ\mathcal H_g^{-1},
 \qquad \mathcal E_g=(\mathcal H_g)_*(z\partial_z)
\]
on the open full stable-chart image in $M$. Here
$\mathcal H_g$ and $\phi_g$ are $H_{x_g}$ and $\phi_{x_g}$
viewed in that trivialization. Their centers are $\theta_g$,
and $\phi_g$ parametrizes the projective big cell in $N$.

Every normalized $Q^+$-map has first component independent
of $z$ and second component linear in $z$. Its restriction
at $z=0$ is exactly the coordinate change on $N$, and it
preserves $z\partial_z$. Thus
\[
             \Pi_{q_0g}=\Pi_g,\qquad
             \mathcal E_{q_0g}=\mathcal E_g
\]
as germs at their common center, almost everywhere for each
fixed $q_0\in Q^+$. Taylor coefficients in local ambient and
target charts are measurable. The same Fubini descent along
$Q^+$-fibers as in Lemma~\ref{II-general:lem:boundary-scalar-measurable}
therefore makes these maps and fields measurable germ fields
$\pi_y,\mathcal E_y$ on $N$. Deck equivariance supplies
their $\Gamma$-equivariance.

Equation~\eqref{II-suspension:eq:negative-transitions} is stronger
than equality on $N$: its transverse coordinate is unchanged.
Since $\phi_{vg}(u)=\phi_g(u+v)$ with the chosen root
convention, it gives
\[
              \Pi_{vg}=\Pi_g,\qquad \mathcal E_{vg}=\mathcal E_g
\]
as germs at $\phi_g(v)$. These identities hold for each fixed
$v$ almost everywhere in $g$. Joint measurability and Fubini
allow a single typical $g$ for which they hold for almost
every $v$ in a small big-cell neighborhood. The map
$v\mapsto\phi_g(v)$ is an analytic diffeomorphism there,
so $\Pi_g,\mathcal E_g$ represent $\pi_y,\mathcal E_y$ for
almost every center in a nonempty open subset of $N$.
Lemma~\ref{II-suspension:lem:retraction-descent} gives an
analytic equivariant retraction $\Pi$ and an invariant analytic
vector field $\mathcal E$ near the entire $N$.

The identities $D\Pi(\mathcal E)=0$, $\mathcal E|_N=0$,
and $D\mathcal E=\id$ on the normal quotient hold in every
chart and hence everywhere along $N$. Apply
Lemma~\ref{II-suspension:lem:euler-collar} with the normal frame
of Corollary~\ref{II-general:cor:normal-character}, in which
the normal derivative is $\|\gamma\theta\|^q$.
This gives the analytic conjugacy
\eqref{II-suspension:eq:negative-collar}.
Inverse transpose proves the dual case.
\end{proof}

The dependency of this argument is worth specifying. It uses
the earlier negative-branch construction only through the
analytic hypersurface and full-dimensional contracting germ
charts. The scalar normal-cocycle theorem uses only atomic
neutral hypersurface description. The descended normal Euler field
linearizes the transverse coordinate by its analytic flow. None of these steps invokes a full collar or
global classification in order to obtain the retraction or
Euler field.

\begin{remark}[Scope of the full-dimensional chart method]
The exact transition argument isolates the role of the sign
of $q$. Negative $q$ supplies the initial full-dimensional
analytic Levi--central charts by simultaneous contraction.
Once such charts exist with $H_x(u,0)=\phi_x(u)$ and central
map $(\lambda u,bz)$, the transition and descent argument also
works when $b>1$. Continue the charts to neighborhoods of
points of the stable plane by contracting $u$ and shrinking
the transverse domain before pulling back; only finitely many
iterates are needed for each prescribed point. The same
radical coefficient calculation as in
Proposition~\ref{II-loc:prop:negative-tangential} gives independence
of $z$ in the angular component and linearity in the normal
component whenever $b\ne1$. The equation $B(bz)=B(z)$ then
forces $B=1$ for either sign of $\log b$.

Thus, for $q>0$, the additional input to this method is a
jointly analytic full-dimensional Levi--central reduction.
Separate analytic tangential and transverse germs do not
establish that input. The next subsection uses contracting
hypersurfaces to obtain a jointly analytic collar without
that full-dimensional reduction. For $q=0$ the central equation permits
a nonconstant $B(z)$, and the nonlinear flow must be retained.
The compatibility and descent arguments above use recurrence
and Howe--Moore ergodicity; they introduce no strong
property-$(T)$ hypothesis on the lattice.
\end{remark}

\section{Expanding neighborhoods}
\label{II-sec:suspension-positive-collar}

Here $N$ is an already recognized analytic projective covering
hypersurface with $q>0$. Use its spherical and oriented collar
cover and the homogeneous probability space $X=G/\Gamma$
in Section~\ref{II-sec:suspension-negative-collar}. The normal
orientation cover over the spherical source is trivial because
$S^m$ is simply connected for $m\ge2$
\cite[Proposition~1.14, p.~35]{Hatcher}. The lifted lattice action
permutes its two connected sheets; the kernel of this sign
homomorphism has finite index and preserves a chosen orientation.
We use this subgroup and its homogeneous Haar probability space.
Write the tangential column coordinates as $u_1,\ldots,u_m$,
where $m=n-1\ge3$, and let $\beta_i=\varepsilon_{i+1}-\varepsilon_1$.
The normal weight is $q\varepsilon_1$ and
$\sum_i\beta_i=-n\varepsilon_1$.

\begin{lemma}[Normalized expanding transverse curves]
\label{II-positive:lem:curves}
There are measurable actual analytic transverse curve germs
$i_x(z)$ satisfying
\[
 F(h,x)i_x(z)=i_{hx}(z)\quad(h\in LU),\qquad
 F(a_t,x)i_x(z)=i_{a_tx}(e^{qt}z).
\]
Their derivatives are nonzero and their radii are positive.
\end{lemma}
\begin{proof}
For $a_{-1}$ the normal exponent is $\ell=-q<0$, whereas the tangential
exponents are $n/m>0$. Apply the random holomorphic stable-manifold
theorem \cite[Theorem~5.1(a)--(c), Case~III, pp.~44--45; proof,
pp.~46--54]{Ruelle1980Author}, in the author version used above.
Indeed, the return-word moment supplies the logarithmic $C^1$ bound
in the finite analytic atlas. With physical return-domain radius $R_x$,
choose $r_x\le\min\{(4(d+1)H_x)^{-1},R_x/4\}$ with integrable
$\log r_x$. The canonical series equals the physical return on this
common domain, and the growth filtration has the unique negative stable
line. The cut $\ell/8$ and faster rate $3\ell/4$, followed by the
projection-inverse construction above, give measurable analytic curves
with this prescribed tangent.

These curves are parabolic equivariant. For $h\in LUA_0$, the sequence
$a_{-d}ha_d$ stays compact. The compact-parameter section estimates in
\cite[Section~3.2]{BRHW} bound each transport word by the two endpoint
section heights plus a constant. These heights are sublinear along
both typical orbits, so the transport and its inverse have subexponential
Lipschitz bounds and preserve exponential convergence.
For completeness, the domain and capture estimates use the same physical
returns: the source cone gives $\|G_x(u)-G_x(v)\|\le4\|u-v\|$, and
$r_{T^nx}\le1$ turns source contraction at $3\ell/4$ into physical
contraction at $\ell/2$. Conversely, logarithmic integrability gives
$r_{T^nx}\ge c_xe^{\ell n/16}$. Physical closeness at rate $\ell/4$
therefore gives normalized decay at $3\ell/16$, faster than the cut
$\ell/8$. Shrinking the initial segment and capture coefficient puts
every iterate in the same source tube, whose projection inverse gives
local capture. Thus transport maps the actual stable-curve germs to
one another, yielding an analytic one-variable parabolic germ cocycle.

Orient the curves by the normal orientation and apply
Theorem~\ref{II-loc:thm:relative}. The Levi is ergodic on
$G/\Gamma$, so the residual cocycle is a single analytic
flow $B(z)\partial_z$. Its derivative has exponent $q$:
recurrence to bounded sets for the derivative gauge compares
the constant multiplier with the original normal exponent.
Since $B'(0)=q\ne0$, the analytic coordinate obtained by
integrating $q/B(z)-1/z$ conjugates this flow to
$z\mapsto e^{qt}z$. Equivalently, Koenigs' theorem
\cite[Theorem~III, \S10, pp.~S15--S16, and \S\S11--12]{Koenigs1884}
(see also \cite[p.~385]{Koenigs1885}) gives, for an analytic germ $f$
with $0<f'(0)=\lambda<1$, an analytic coordinate $K$ satisfying
$K'(0)=1$ and $K\circ f\circ K^{-1}=\lambda\,\mathrm{id}$.
For real $f$, take the real part of the normalized holomorphic coordinate
on the real axis; its derivative is one and the functional equation is
preserved. Apply this to the time-$-1$ map, where $\lambda=e^{-q}$.
Commutation then makes every time map scalar in the same coordinate,
with multiplier $e^{qt}$. This proves the displayed identities.

\end{proof}

For $k=1,2$ put $I_k=\{u_k=0\}\subset\R^m$. Choose a
diagonal element $c_k\in A$ by its logarithmic weights
\begin{equation}\label{II-positive:eq:scalar-weights}
 \varepsilon_1(c_k)=-1,\qquad
 \beta_i(c_k)=-q\ (i\ne k),\qquad
 \beta_k(c_k)=n+(m-1)q.
\end{equation}
They satisfy the trace relation. Its fiber stable space is
exactly $I_k\oplus\R_{\mathrm{normal}}$, and all its stable
exponents equal $-q$. The remaining exponent is positive.

\begin{lemma}[Calibrated scalar-contracting hypersurfaces]
\label{II-positive:lem:stable-charts}
There are measurable analytic stable charts
\[
 K_x^k:I_k\times\R\longrightarrow W^{s,F}_{c_k}(x),
 \qquad
 F(c_k,x)K_x^k(w,z)=K_{c_kx}^k(e^{-q}w,e^{-q}z),
\]
with
\[
 K_x^k(w,0)=\phi_x(w),\qquad K_x^k(0,z)=i_x(z)
\]
as germs, and globally on their continued axes.
Their local analytic bounds and inverse bounds are tempered
along $c_k$. The central group and the diagonal Levi act
linearly in these charts, with their prescribed derivatives.
\end{lemma}
\begin{proof}
The stable-manifold construction solves the stable coordinates
forward and the single unstable coordinate backward in a
contracting holomorphic sequence space. The integrable analytic
return bounds give tempered domains, as in the earlier analytic
stable-plaques proof. The full frame given by $D_0\phi_x$
and $i_x'(0)$ makes the ambient derivative diagonal, with
the weights in \eqref{II-positive:eq:scalar-weights}. Its linear
covariance bounds increments of its logarithmic height and
inverse height by the integrable derivative return heights.
Thus this frame is tempered, including the angle between its
two summands. Its stable restriction is the frame
$D_0\phi_x|_{I_k},i_x'(0)$ with scalar derivative $e^{-q}I$.

In stable analytic coordinates with this derivative, the
nonstationary Koenigs sequence
$e^{dq}f_x^{(d)}$ converges holomorphically on a smaller
tempered ball. Its quadratic errors are summable because
$e^{-2q}<e^{-q}$; the estimates are those of Step 3 in
Theorem~\ref{II-loc:thm:contracting-axis}, with $N=1$ and a
scalar derivative. This gives $K_x^k$, and contraction
continues it to the whole intrinsic stable manifold.

The two displayed restrictions agree with their prescribed
analytic parametrizations. Compare them as maps into the
entire ambient fiber. At the first different degree $r\ge2$,
covariance multiplies a stable-output coefficient by
$e^{(r-1)q}$ and an unstable-output coefficient by
$e^{n+(m-1+r)q}$. Both factors exceed one, so recurrence
forces that coefficient to vanish. In particular this
comparison proves that $i_x$ lies in the stable hypersurface;
that containment was not assumed in the construction.
The same argument for a commuting central or diagonal Levi
element eliminates every nonlinear term of its chart map.
\end{proof}

\begin{proposition}[Exact translation on each contracting hypersurface]
\label{II-positive:prop:translations}
For every fixed column $v\in I_k$, almost everywhere in $x$,
\begin{equation}\label{II-positive:eq:translation}
 (K_x^k)^{-1}v^{-1}K_{\widehat v x}^k(w,z)=(w+v,z).
\end{equation}
In particular $v^{-1}i_{\widehat v x}(z)=K_x^k(v,z)$
is jointly analytic in $v\in I_k$ and $z$.
\end{proposition}
\begin{proof}
The base column $v$ contracts under $c_k$. The bounded-transport
argument in Lemma~\ref{II-positive:lem:curves} identifies the two
stable manifolds in the displayed formula. For the related non-uniform
normal-form uniqueness and affine-atlas mechanism, see
\cite[Theorem~2.3(2) and Corollary~2.6]{KalininSadovskaya2017}.
Here the column changes the base and is contracted by $c_k$;
its moving-conjugate transport and the bounds on the two fixed
orbits are verified in the argument below. Their chart
transition is affine. Indeed, after $d$ iterates its displacement
has size $e^{-dq}v$ and lies in the tempered inverse-chart
domain. The return maps for the bounded base displacement,
the charts, and their inverses all have subexponential
analytic bounds along the two fixed typical orbits.
Cauchy estimates therefore bound each fixed-degree transition
coefficient subexponentially. A degree-$r$ coefficient is
also its original value times $e^{(r-1)dq}$, and hence vanishes
for $r\ge2$. Boundary calibration leaves the form
\[
 P_v(x)(w,z)=(w+v+A_v(x)z,B_v(x)z).
\]

Fix a nontrivial unipotent $l=I+w_0e_k^*$ of the Levi, where
$w_0\in I_k$. It fixes $I_k$ pointwise and is contracted by
$c_k$, since its root has weight $-q-[n+(m-1)q]<0$.
The same stable transport and affine argument applies to its
centered $K^k$ chart map. Its derivative is the identity on
the tangential $I_k$ and on the normalized curve $i_x$.
Thus that entire chart map is the identity. As $l$ commutes
with every column in $I_k$, both $A_v$ and $B_v$ are
$l$-invariant. The Mautner argument promotes this invariance
to the full Levi action, whose ergodicity makes them deterministic.

The column group law gives
\[
 B_{v+w}=B_vB_w,\qquad A_{v+w}=A_w+A_vB_w.
\]
These identities may first be evaluated for small columns on
a positive-measure set with common chart domains, as in
Proposition~\ref{II-suspension:prop:negative-transitions}, and
then extended by the group law. A diagonal Levi element
scaling all of $I_k$ by $s>0$ and the complementary coordinate
by $s^{-(m-1)}$ fixes the normal parameter. Its linear action
gives $B_{sv}=B_v$ and $A_{sv}=sA_v$.
Taking $s=2$ and using the group law gives
$B_v=B_{2v}=B_v^2$. Invertibility therefore gives $B_v=1$.
Central covariance, with
$\lambda_t=e^{-nt/m}$, gives
\[
 A_{\lambda_t v}=\lambda_t e^{-qt}A_v.
\]
The diagonal identity already gives
$A_{\lambda_t v}=\lambda_t A_v$. Since $q\ne0$, comparison
at any nonzero $t$ gives $A_v=0$. This proves
\eqref{II-positive:eq:translation} and its final assertion.
\end{proof}

\begin{lemma}[Two analytic slice families with a common parameter]
\label{II-positive:lem:parameter-cross}
Let $f(s,t,w,z)$ be a measurable family of analytic $z$-germs
with finite-dimensional real vector values, on a nonempty open real
box in $(s,t,w)$, where $s,t\in\R$ and
$w\in\R^r$. Suppose that for almost every fixed $s$ its
restriction has a jointly real analytic representative in
$(t,w,z)$ along $z=0$, and for almost every fixed $t$ it
has one in $(s,w,z)$. Then on some nonempty smaller box in
$(s,t,w)$ it has a jointly analytic representative for
$|z|<\epsilon$, for some $\epsilon>0$.
\end{lemma}
\begin{proof}
Shrink the real box relatively compactly and fix an interior
value $w_0$. Countable selection of holomorphic domains and
bounds gives positive-measure Borel sets $A$ of $s$ values
and $B$ of $t$ values, complex domains $D,E$ around their
real intervals, a polydisk $W$ about $(w_0,0)$ in the
parameter $(w,z)$, and a bound $C$, such that the two
selected families are holomorphic and bounded by $C$ on
$E\times W$ and $D\times W$, respectively.
These selections are measurable: bounded holomorphic families
are compact metrizable in the compact-open topology after
shrinking, and agreement with a measurable germ is tested
by countably many coefficients and almost-everywhere equalities.
Fubini and analytic uniqueness in $(w,z)$ give agreement
on almost every $(a,b)\in A\times B$ for all parameters in $W$.

We remove this exceptional set explicitly. Choose a countable
dense set $\{a_d\}$ in the essential support $A_*$ of $A$
such that the agreement holds for almost every $b$ at every
$a_d$, and discard one null subset of $B$. For $a\in A_*$,
take $a_{d_j}\to a$. Normality of the bounded family of
holomorphic functions on $E\times W$ gives subsequential
limits. Every such limit has the prescribed values on
$B\times W$, so uniqueness in the $t$ variable makes the
limit unique. This corrects the family at every $a\in A_*$,
and agrees with the original family for almost every $a\in A$.
We now have an exact holomorphic cross
\[
       \bigl((D\times B)\cup(A_*\times E)\bigr)\times W.
\]

Pass both cross sets to their locally pluriregular cores; the discarded
sets are locally pluripolar \cite[p.~222]{NguyenCrossRegularization}. Their
intersections with the smooth totally real coordinate planes have real
measure zero \cite[Theorem~4, p.~513]{Sadullaev1976}, so the retained
cores still have positive measure. The bounded extension in
\cite[Theorem~A, p.~224]{NguyenCrossRegularization} preserves the cross supremum
bound on a fixed envelope. This envelope contains a product neighborhood
of retained regular density points and is independent of the cross function.

Expand our families in the parameter $(w-w_0,z)$ on a
smaller polydisk of radius $R$. The coefficient with multiindex
$\nu$ is a bounded scalar cross function with bound
$CR^{-|\nu|}$. Extend every coefficient by the preceding
bounded extension on that same envelope. Its supremum estimate preserves
the exponential coefficient bound on a fixed smaller product domain. The power series converges normally after
shrinking the parameter disk. Its sum is jointly holomorphic
in all variables and extends the selected slice families.

It represents the original measurable germs on an open real
box, rather than only on the selected cross. For almost every
$a\in A$ in that box, equality holds on the entire $t$ slice
and for all real parameters near $(w_0,0)$. For almost every
$t$ in the box, compare their Taylor coefficients in $z$.
These coefficients are real analytic in $(s,w)$ and agree on the
product of a positive-measure set of $s$ values with the smaller
$w$ box. A nonzero real analytic function on a connected open
domain has a null zero set
\cite[Proposition~0 in the author preprint]{Mityagin}.
Thus every coefficient agrees on the entire smaller box, and
Taylor uniqueness identifies the original $z$-germs there.
Fubini proves the assertion. For finite-dimensional targets,
use common selected domains or successive nested smaller boxes. In an ambient analytic coordinate chart it also
applies to the manifold-valued germs used below.
\end{proof}

\begin{theorem}[Full expanding collars for arbitrary lattices]
\label{II-positive:thm:collar}
Let $\Gamma<\SL(n,\R)$ be any lattice, $n\ge4$, and let
$N$ be an invariant compact embedded analytic projective
covering hypersurface with its standard or dual action and
normal coefficient $q>0$. After the angular and normal-orientation
covers and finite-index restriction, its full collar germ is
analytically conjugate to
\begin{equation}\label{II-positive:eq:collar}
 \gamma(\theta,t)=
 \left(\frac{\sigma(\gamma)\theta}{\|\sigma(\gamma)\theta\|},
                 \|\sigma(\gamma)\theta\|^q t\right).
\end{equation}
\end{theorem}
\begin{proof}
Work in the standard case and lift the curve field of
Lemma~\ref{II-positive:lem:curves} to $G$, writing it $i_g(z)$
in the physical fiber $M$. It satisfies
\[
 i_{q_0g}(z)=i_g(e^{-q\chi(q_0)}z),\qquad
 i_{g\gamma}(z)=\alpha(\gamma)^{-1}i_g(z)
\]
almost everywhere for every fixed element.
In a column big cell consider $f(v,z)=i_{vg_0}(z)$.
Joint measurability and Fubini allow a typical $g_0$ for which
Proposition~\ref{II-positive:prop:translations} supplies jointly
analytic representatives along almost every affine hyperplane
parallel to $I_1$ and to $I_2$, including the parameter $z$.
For clarity, choose one typical center on each such hyperplane;
its globally continued $K^k$ chart represents the field for
almost every displacement in that hyperplane. This gives an
analytic representative along the whole hyperplane at $z=0$.

Choose a small angular box whose zero section lies in one
ambient analytic coordinate chart. Apply
Lemma~\ref{II-positive:lem:parameter-cross} with
$s=v_1$, $t=v_2$, and $w=(v_3,\ldots,v_m)$.
It gives a jointly analytic representative on a nonempty
angular box with a common positive normal radius. At $z=0$
the representative equals the recognized parametrization
$\phi_{g_0}(v)$. Its normal derivative is the nonzero derivative
of $i_{vg_0}$ almost everywhere. Choose such a point and
shrink; the representative is a genuine analytic collar chart.

We give the descent to the entire hypersurface, including its
normal parameter. Put
\[
 r(g)=-\log\|g^{-1}e_1\|,\qquad
 J_g(t)=i_g(e^{qr(g)}t).
\]
Since $r(q_0g)=r(g)+\chi(q_0)$, the field $J_g$ is left
$Q^+$-invariant almost everywhere. Fubini descent in local
product coordinates makes it a measurable parametrized curve
field $J_\theta(t)$ on $N$. The right-lattice identity and
\[
 r(g\gamma)-r(g)=\log\|\gamma\theta_{g\gamma}\|
\]
give
\begin{equation}\label{II-positive:eq:curve-descent}
 \alpha(\gamma)J_\theta(t)
       =J_{\gamma\cdot\theta}(\|\gamma\theta\|^q t).
\end{equation}
The scale $r(vg_0)$ is analytic on the angular patch, so the
patch just constructed is also an analytic representative of
$J_\theta$ there.

Translate this patch by \eqref{II-positive:eq:curve-descent}.
Projective minimality makes its translates cover $N$; take a
finite subcover. On overlaps the entire parameter germs
agree for almost every center. Analytic uniqueness first in
the parameter and then in the angular variables makes them
agree everywhere. This is the same germ-gluing argument as
Lemma~\ref{II-suspension:lem:retraction-descent}, now for
parametrized maps with the displayed analytic transition law.
It gives a jointly analytic $J$ near $N\times\{0\}$ and
makes \eqref{II-positive:eq:curve-descent} hold everywhere as
germs. Its derivative is invertible along the zero section
by transport of the initial collar patch. Compactness and
the embeddedness of $N$ give an embedded smaller collar:
any failure of injectivity at arbitrarily small widths would
limit to two equal points of $N$ and contradict local
invertibility there. In these coordinates
\eqref{II-positive:eq:curve-descent} is precisely
\eqref{II-positive:eq:collar}. Inverse transpose handles the dual.
\end{proof}

This proof obtains ambient analyticity from two families of
contracting hypersurfaces. It uses the actual one-dimensional
parabolic germ reduction, but no full-dimensional saddle
linearization. The construction of these hypersurfaces uses
$q>0$ in \eqref{II-positive:eq:scalar-weights}; it gives no
contracting transverse direction when $q=0$.

\section{Neutral axes and full neighborhoods}
\label{II-sec:neutral-fixed-curves}

Return to an arbitrary lattice $\Gamma<G=\SL(n,\R)$, $n\ge4$,
and a recognized invariant analytic projective hypersurface $N$.
Use its canonical homogeneous suspension $X=G/\Gamma$, as in
Section~\ref{II-sec:suspension-negative-collar}, and write
$Q=(L\ltimes U)\rtimes A_0$, $L=\SL(m,\R)$, $m=n-1$.
Suppose the normal character is $q=0$. Work on the indicated
angular and normal-orientation covers. In particular, the
normal derivative has the analytic trivialization supplied by
Corollary~\ref{II-general:cor:normal-character}. The radical shear
is zero because the actual hypersurface is invariant.

\begin{theorem}[Neutral transverse curves are analytic]
\label{II-neutral:thm:analytic-curves}
There are measurable analytic embedded transverse curve germs
$i_x:(\R,0)\to(M_x,x)$ and one real analytic function $B$ with
$B(0)=B'(0)=0$ such that
\begin{equation}\label{II-neutral:eq:actual-curves}
 \begin{aligned}
 F(h,x)i_x(z)&=i_{hx}(z)&& (h\in LU),\\
 F(a_t,x)i_x(z)&=i_{a_tx}(\Phi_t^{B\partial_z}(z)).
 \end{aligned}
\end{equation}
These are identities of actual analytic germs, with positive
measurable radii. Neither an arithmetic stabilizer in $\Gamma$
nor a contracting normal exponent is required.
\end{theorem}

\begin{proof}
Proposition~\ref{II-loc:prop:neutral-formal-dichotomy} supplies
formal charts $H_x(u,z)$, calibrated on the boundary by
$H_x(u,0)=\phi_x(u)$, in which the Levi is $(u,z)\mapsto(lu,z)$.
Take analytic ambient coordinates with the derivative frame of
these charts; this is only a first-order coordinate choice.
The formal axis $I_x(z)=H_x(0,z)$ satisfies
\[
                         I_{lx}=F(l,x)I_x.
\]
Choose an ambient lattice $\Lambda< L$ conjugate to
$\SL(m,\Z)$. Its action on $G/\Gamma$ is ergodic: a single
nontrivial unipotent element of $\Lambda$ is already ergodic,
by Howe--Moore. Since $m\ge3$, property $(T)$ gives spectral
gap. Adjoin one diagonal element of $L$ all of whose standard
eigenvalues have modulus different from one. The resulting
countable group retains this spectral gap and ergodicity.
Its derivative representation is standard plus one fixed line.
Theorem~\ref{II-sg:thm:fixed-subspace} therefore makes
$i_x(z)=H_x(0,C(z))$ analytic for one constant formal
tangent-to-identity reparametrization $C$. Because this
reparametrization is constant, the formal identities for the
entire Lie Levi remain valid; analyticity makes them actual.

The formal row fields
$(vu)(a(z)E_u+b(z)\partial_z)$ vanish on $u=0$, so the
formal row maps fix the axis pointwise. The actual row maps
and the now analytic curves consequently satisfy the first
identity in \eqref{II-neutral:eq:actual-curves}. The normalized
central formal maps have the form $(A_{t,x}(z)u,D_{t,x}(z))$;
they preserve the axis. An analytic left inverse of $i_{a_tx}$
shows that the induced central one-variable map is analytic.
Commutation with the Levi and its ergodicity on $X$ make all
its coefficients independent of $x$. The central cocycle law
is therefore a coefficient-measurable one-parameter
homomorphism of analytic germs. Apply
Lemma~\ref{II-time:lem:one-parameter}. The constant term is zero
because the curves are based on $N$, and the linear term of
the generator is zero because the normalized normal derivative
is one. This proves the theorem.
\end{proof}

\begin{lemma}[A directional stabilizer has the required spectral gap]
\label{II-neutral:lem:direction-gap}
For $0\ne v\in\R^m$, let $H_v<L$ be its stabilizer. This
group contains a countable subgroup whose action on $G/\Gamma$
is ergodic with spectral gap. It also contains a diagonalizable
element whose fixed space on $\R^m\oplus\R$ is exactly
$\R v\oplus\R$ and whose other eigenvalues have modulus
different from one.
\end{lemma}
\begin{proof}
Choose coordinates $\R^m=\R v\oplus W$. The stabilizer is
\[
 H_v=\left\{\begin{pmatrix}1&r\\0&A\end{pmatrix}:
          A\in\SL(m-1,\R),\ r\in W^*\right\}.
\]
Choose a two-dimensional subspace of $W$ and fix the remaining
coordinates. This gives a copy of
$H_0=\SL(2,\Z)\ltimes\Z^2$ inside $H_v$, with the dual
standard realization of its radical. The pair $(H_0,\Z^2)$
has relative property $(T)$ by \cite[Theorem~1.2]{Iozzi},
with no inverted primes. The standard and dual realizations
are identified by contragredient twisting of the Levi. A nontrivial element of the
radical is unipotent and acts ergodically on $G/\Gamma$ by
Howe--Moore. Thus the radical-fixed vectors in
$L^2_0(G/\Gamma)$ vanish. Relative property $(T)$ excludes
almost invariant unit vectors and gives a finite spectral-gap
set. This also proves ergodicity. Finally take a positive
diagonal matrix $A\in\SL(m-1,\R)$ with no eigenvalue one.
The block matrix $\operatorname{diag}(1,A,1)$ on
$\R^m\oplus\R$ has the stated property. It may be adjoined
to the countable subgroup without losing its spectral gap.
In particular no property $(T)$ assertion about $\SL(2,\R)$
is being used when $n=4$.
\end{proof}

\begin{proposition}[Analytic surfaces containing the neutral curves]
\label{II-neutral:prop:direction-surfaces}
For each fixed nonzero $v\in\R^m$, there are measurable
analytic embedded surface germs $K_x^v:(\R^2,0)\to(M_x,x)$
with
\begin{equation}\label{II-neutral:eq:surface-axes}
 K_x^v(s,0)=\phi_x(sv),\qquad K_x^v(0,z)=i_x(z),\qquad
 F(h,x)K_x^v=K_{hx}^v\quad(h\in H_v).
\end{equation}
As a formal submanifold, their image is
$H_x(\R v\oplus\R)$. Here the formal charts can be chosen
to have the actual axis $i_x$.
\end{proposition}
\begin{proof}
Replace the formal chart by $H_x(u,C(z))$, with $C$ from
Theorem~\ref{II-neutral:thm:analytic-curves}. This preserves its
Levi conjugacy and boundary calibration. Choose measurable
analytic ambient coordinates in which both its boundary
restriction and its normal axis are the corresponding
coordinate inclusions. Such coordinates exist without any
joint regularity assumption on $x$: in an ambient coordinate
chart the map $\phi_x(u)+i_x(z)-x$ has invertible derivative
at zero and has precisely these two restrictions.

Restrict the resulting formal $H_x$ to the two-dimensional
plane $E_v=\R v\oplus\R$. The group $H_v$ fixes this source
pointwise. Lemma~\ref{II-neutral:lem:direction-gap} and
Theorem~\ref{II-sg:thm:fixed-subspace} make the restriction
analytic after a constant formal change in its two source
variables. Write $I_x$ for this restriction and $C$ for the
constant supplied by that theorem. To preserve the two axes,
let $P$ be the orthogonal projection onto $E_v$, choose a point
$x_0$ where $I_{x_0}\circ C$ converges, the first derivative
of $I_{x_0}$ is the inclusion, and the original axis calibrations
hold, and put $A=P\circ I_{x_0}$. This is a
formal map tangent to the identity that fixes both axes.
The convergent formal map $A\circ C$ is the Taylor series of
an invertible analytic germ $q$. Set
\[
 D=C\circ\widehat{q^{-1}}=A^{-1},
\]
where the hat denotes Taylor series. Then $D$ fixes both axes,
and $I_x\circ D=(I_x\circ C)\circ\widehat{q^{-1}}$
converges almost everywhere. Precomposing the analytic realizations
of $I_x\circ C$ by the fixed analytic germ $q^{-1}$ gives
realizations of $I_x\circ D$. This preserves measurability and positive
embedding radii, after shrinking the radii if necessary.
The actual analytic axis identities follow from equality of
Taylor series. This proves the first two identities in
\eqref{II-neutral:eq:surface-axes}. A constant change of the
fixed source preserves equivariance for the full $H_v$, so
the last identity also holds as an analytic germ identity.
An invertible formal source change does not alter the formal
image, proving the final assertion.
\end{proof}

\subsection{Moving the center of a neutral transverse germ}
\label{II-sec:neutral-moving-centers}

The comparison takes place inside the surfaces just constructed.
For a fixed $v\ne0$, put
$V_s=\begin{pmatrix}1&0\\-sv&I\end{pmatrix}$ and
$T_sx=\widehat{sv}x$. On the canonical homogeneous section,
$T_s$ is ordinary left translation by $V_s$. It commutes
with $H_v$. As before, $V_s^{-1}$ denotes actual transport
between suspension fibers, not an element assumed to lie
in $\Gamma$.
All surface maps and their continuations in the argument below
take values in the original compact fiber $M$. The spherical
and normal-orientation covers specify their angular labels and
normal frames; no finite cover of the entire manifold is
assumed. Consequently the fiber transports used for continuation
are globally defined analytic diffeomorphisms, even when the
collar cover itself does not extend to a global cover.

\begin{lemma}[Analytic motion along a boundary direction]
\label{II-neutral:lem:direction-motion}
For each fixed $v\ne0$ there is $\epsilon_v>0$ such that,
for almost every $x$, the family of curve germs
\[
                (s,z)\longmapsto V_s^{-1}i_{T_sx}(z)
\]
has a jointly analytic representative as a germ along
$(-\epsilon_v,\epsilon_v)\times\{0\}$. It represents the
given measurable germs for almost every $s$. Its restriction
to $z=0$ is $\phi_x(sv)$.
\end{lemma}

\begin{proof}
Write $K_x=K_x^v$. First continue this analytic parametrized
surface along its boundary axis, retaining local immersion and
transversality to $N$. The maximal interval around $s=0$ on
which such continuation exists is $H_v$-invariant as a
function of $x$. Indeed each actual $h$-transport is a global
analytic diffeomorphism of the corresponding compact fibers;
it preserves $N$ and takes the boundary point $\phi_x(sv)$ to
$\phi_{hx}(sv)$ with the same parameter $s$. The identity
in \eqref{II-neutral:eq:surface-axes} transports any finite chain
of analytic continuation neighborhoods in both directions.
The two endpoint radii are measurable: finite chains can be
chosen in countably many coordinate charts and polydisks, with
coefficient bounds and strict inverse-derivative bounds. Their
positivity and $H_v$-ergodicity give a common real interval
$(-R,R)$ on which $K_x$ is an analytic germ along the whole
boundary axis, for almost every $x$. Its transverse widths
may still depend on $x$ and on $s$.

\emph{The translated surfaces coincide.}
Fix $|s|<R$. Transport the germ of $K_x$ at $(s,0)$ by $V_s$
to the fiber at $y=T_sx$. This is a measurable
$H_v$-equivariant surface germ based at $y$. It meets $N$ in
the standard $v$-line and is transverse to $N$. Conjugate its
formal germ by $H_y^{-1}$. It is a formal graph over
$E_v=\R v\oplus\R$ with values in the complementary standard
module $W$: its tangent plane projects isomorphically to
$E_v$ because it contains the $v$-line and is transverse to
the boundary. The block Levi $J_v=\SL(W)<H_v$ fixes the
graph's source pointwise. Each graph coefficient is therefore
a measurable $W$-valued equivariant function under $J_v$.
Use the diagonal element in
Lemma~\ref{II-neutral:lem:direction-gap}. Recurrence eliminates
each of its non-unit-modulus eigencoordinates, so every graph
coefficient is zero. This includes the first-order graph
coefficient. Hence the formal surface is exactly $E_v$.

Both this surface and $K_y$ are analytic. Equality of their
formal graphs gives equality of their analytic submanifold
germs. Thus, for each fixed sufficiently small $s$, the
following transition is well defined almost everywhere as
an analytic map germ at $(0,0)$ with value $(s,0)$:
\begin{equation}\label{II-neutral:eq:direction-transition}
              P_s(x)=K_x^{-1}V_s^{-1}K_{T_sx}.
\end{equation}
The inverse here is the local inverse on the surface at
$(s,0)$; no global injectivity of its continuation is needed.
Commutation of $V_s$ with $H_v$ and exact
$H_v$-equivariance of both parametrizations imply
$P_s(hx)=P_s(x)$. Ergodicity makes all its coefficients
deterministic. Denote the resulting analytic germ by $P_s$.
The original measurable families and countable Taylor tests
give a jointly measurable version in $(s,x)$; Fubini gives
agreement with this deterministic version for almost every
$(s,x)$.

\emph{The deterministic transitions form an analytic local flow.}
Here are the domain and parameter details. Choose a
positive-measure set $D\subset X$ lying in one relatively
compact local trivialization of the suspension, on which
$K_x$ has a common holomorphic bidisk, a common bound, and
a bounded analytic left inverse near its image. This is
possible by restricting measurable positive radii and bounds,
and then restricting to maps close to one fixed embedding.
The local lift of $T_s$ and its fiber transport vary smoothly
and analytically for small $s$. Strong continuity of
translations on $L^2(X)$ gives
$\mu(D\cap T_s^{-1}D)>0$ for all sufficiently small $s$.
Evaluating \eqref{II-neutral:eq:direction-transition} on that set
gives holomorphic representatives of the deterministic $P_s$
on one smaller bidisk, with a common bound. Shrink the
parameter interval once more so their images have the domains
needed for composition. For any two small parameters, the
intersection of the three required translates of $D$ still
has positive measure. The actual fiber transports there give
\[
                   P_s\circ P_t=P_{s+t},\qquad P_0=\id.
\]
The identity initially holds as a germ and then on the common
domains by analytic uniqueness. Also $P_s(t,0)=(s+t,0)$.

The maps $P_s$ tend to the identity uniformly on a smaller
complex bidisk as $s\to0$. To see this without assuming
regularity of the original charts in $x$, use continuity in
measure of translations of a measurable map into a separable
space. On $D$ the bounded holomorphic charts, restricted to a
smaller bidisk, take values in a compact metrizable normal
family. Thus $K_{T_sx}$ and $K_x$ are arbitrarily close on a
positive-measure subset of $D\cap T_s^{-1}D$ when $s$ is
small. Fiber transport tends to the identity in the chosen
local trivialization. The bounded left inverses then show
that the deterministic $P_s$ is arbitrarily close to the
identity. Cauchy estimates give the same conclusion for
derivatives on a still smaller bidisk.

For completeness, continuity together with the uniform
holomorphic domains makes the generator analytic. For small
$\tau>0$ define
\[
                    H(w)=\frac1\tau\int_0^\tau P_t(w)\,dt.
\]
It is analytic and $C^1$ close to the identity, hence locally
invertible. The group law gives
\[
 H(P_s(w))=\frac1\tau\int_s^{s+\tau}P_t(w)\,dt.
\]
Differentiating this identity in the real parameter at zero
gives the analytic generator
\[
                  Y(w)=DH(w)^{-1}\frac{P_\tau(w)-w}{\tau}.
\]
The same identity gives the flow equation for every small
$s$. Analytic dependence for this ordinary differential
equation proves that $(s,w)\mapsto P_s(w)$ is jointly
analytic.

Finally \eqref{II-neutral:eq:surface-axes} and
\eqref{II-neutral:eq:direction-transition} give
\[
                    V_s^{-1}i_{T_sx}(z)=K_x(P_s(0,z)).
\]
The right side is jointly analytic along a common smaller
real $s$-interval and $z=0$, with transverse widths allowed
to depend on $x$. Fubini supplies the asserted almost-everywhere
agreement with the original curves. The boundary identity
follows from $P_s(0,0)=(s,0)$.
\end{proof}

\begin{lemma}[An analytic patch from coordinate-line germs]
\label{II-neutral:lem:iterated-cross}
Let $f(u,z)$ be a coefficient-measurable field of analytic
$z$-germs on a real box $U\subset\R^m$. Suppose that for
each coordinate direction, on almost every parallel line,
the field has a jointly analytic representative as a germ
along that line times $\{0\}$. There is a nonempty open
subbox $U'\subset U$ and a jointly analytic representative
near $U'\times\{0\}$ agreeing with the given germs for
almost every $u\in U'$. The assertion also holds for maps
into a finite-dimensional analytic manifold, after restricting
to a coordinate chart.
\end{lemma}

For a manifold target, coefficient measurability is understood in a
countable analytic atlas $\kappa_j:V_j\to W_j\subset\R^d$:
there is a measurable chart index $j(u)$ and a field of convergent
coordinate power series $p_u$ whose coefficients are measurable in $u$,
with $p_u(0)\in W_{j(u)}$, such that
$\kappa_{j(u)}^{-1}\circ p_u$ represents the original germ $f(u,\cdot)$.
The chart index and positive convergence radius may vary with $u$;
no common chart or common radius is assumed. As usual, analytic
manifolds are finite dimensional, Hausdorff, second countable, and
without boundary.

\begin{proof}
We first record the version of the cross argument with two
blocks $x\in\R^k$, $y\in\R$ and an additional analytic
parameter $z$. Suppose there are positive-measure Borel sets
$A,B$ and common complex domains $D\subset\C^k$,
$E\subset\C$, and a disk $W$ about $z=0$, with bounded
holomorphic representatives on $D\times W$ for $y\in B$
and on $E\times W$ for $x\in A$, agreeing almost everywhere
on the real cross. The normal-family correction in the
proof of Lemma~\ref{II-positive:lem:parameter-cross} makes an
exact separately holomorphic cross. Explicitly, choose a
countable sequence in a conull subset of $A$ dense in its
essential support. Discard one null subset of $B$ so that all
these representatives agree. Normal limits of the bounded
$E\times W$ representatives along that sequence extend the
vertical representatives to the essential support of $A$;
their values on $B\times W$ determine the limit uniquely.
The horizontal representatives remain unchanged. Pass to locally
pluriregular cores as in \cite[p.~222]{NguyenCrossRegularization}. The discarded
locally pluripolar sets are real-null in the totally real coordinate planes
by \cite[Theorem~4, p.~513]{Sadullaev1976}, and the retained cores
have positive real measure.

Apply the threefold cross theorem to $(D,A)$, $(E,B)$, and
$(W,W)$ \cite[Theorem~4.3.3, pp.~60--61 in the author version]{JarnickiPflugSeparatelyAnalytic2008}.
On the third leg the common values on $A\times B\times W$ are
holomorphic in $z$; the other two legs are the given representatives.
The open disk is its own regular core and its plurimeasure is zero.
The extremal-function identification
\cite[Proposition~3.4.11(a), p.~39]{JarnickiPflugSeparatelyAnalytic2008}
therefore gives one envelope containing a product neighborhood of retained
regular points of $A\times B$ times the whole parameter disk.
Indeed the two base plurimeasures vanish at these points, and upper
semicontinuity makes their sum less than one on a smaller product.
Choose the retained points in the supports of the restricted real measures,
so every further real-ball restriction still has positive measure.
For a finite-dimensional target, apply the scalar theorem to its coordinates
on this same envelope. This gives the required jointly holomorphic patch,
with the original cross agreement, after shrinking the parameter disk.

Induct on $m$. For $m=1$ the hypothesis is the conclusion.
For $m>1$ put $u=(x,y)$ with $x\in\R^{m-1}$. For almost
every $y$, the induction hypothesis supplies an analytic
representative in $(x,z)$ on some open subbox. Countably many
choices of subboxes, complex neighborhoods, radii, and integer
bounds give one such choice on a positive-measure Borel set
$B$ of $y$'s. The representatives can be chosen measurably:
bounded holomorphic functions on smaller neighborhoods form
compact metrizable normal families, and agreement of their
$z$-Taylor coefficients with the given measurable field is
tested by countably many integral equalities. Uniqueness on
an open real box makes the selected representative unique.

Choose an interval $J$ meeting $B$ in positive measure.
For almost every $x$ in the selected subbox, the original
coordinate-line hypothesis gives a jointly analytic
representative in $(y,z)$ near $J\times\{0\}$. Restrict $J$
compactly and use countably many bounds to select a
positive-measure set $A$ of $x$'s with a common complex
neighborhood and common parameter radius. The preceding
two-block cross argument supplies the required open patch.
It agrees with the original germs almost everywhere on that
whole patch: the horizontal representatives agree on
$D\times B$, and for almost every $x$ the one-variable
identity theorem in $y$, applied coefficientwise in $z$,
propagates this agreement from $B$ across the interval of
the patch. Restrict the patch around a point of the essential
support of $B$ so that this uniqueness set has positive
measure there. This proves the induction. For manifold-valued germs, first construct a jointly measurable
representative from the preceding countable-chart presentation. For
 each $j$, extend $\kappa_j^{-1}$ outside $W_j$ by any fixed point of
$V_j$, and apply these measurable maps to the sums $p_u(z)$ using
the measurable index $j(u)$. The resulting representative has the same
analytic germ in every fiber. Choose an analytic Euclidean embedding
$e:M\hookrightarrow\R^N$ by
\cite[Theorem~3, p.~470]{Grauert}; compactness of $M$ is not required.
The Taylor coefficients of $e\circ f$ are measurable: finite-difference
limits give the derivatives of its jointly measurable representative,
and polarization gives the corresponding symmetric Taylor tensors.
The given jointly analytic coordinate-line representatives compose
with $e$. Apply the Euclidean result to obtain $P(u,z)$ agreeing with
the embedded original germs for almost every $u$ on a nonempty open
subbox. Choose one such $u_0$ and write $p_0=f(u_0,0)$, so that
$P(u_0,0)=e(p_0)$.

The analytic immersion normal form gives an analytic map $R$ near
$e(p_0)$ which is a left inverse of $e$ on a neighborhood $V$ of $p_0$:
compose the ambient normal-form chart with projection onto its
tangential coordinates and the inverse source chart. Because $e$ is
an embedding, there is an ambient open neighborhood $O$ of $e(p_0)$
with $e^{-1}(O)\subset V$. Shrink the output subbox and its positive
$z$-radius so that $P$ takes values in $O$ and in the domain of $R$.
For almost every base point $u$, the germ identity
$P(u,\cdot)=e\circ f(u,\cdot)$ then places the original germ in $V$,
and $R\circ P$ agrees with that germ. Thus $R\circ P$ is the required
jointly analytic manifold-valued representative. A further local
restriction puts its image in one coordinate chart. If the base
dimension is zero, the box has just one point and the assertion follows
directly from analyticity of its single transverse germ.
\end{proof}

\begin{theorem}[Full neutral collars for arbitrary lattices]
\label{II-neutral:thm:full-collar}
Let $\Gamma<\SL(n,\R)$ be any lattice, $n\ge4$, and let
$N$ be a compact invariant analytic standard or dual
projective covering hypersurface whose normal character is
$q=0$. On its angular and normal-orientation covers, and
after the indicated finite-index restriction, its full
analytic collar germ has the form
\begin{equation}\label{II-neutral:eq:full-collar}
 \gamma(\theta,t)=
 \left(\frac{\sigma(\gamma)\theta}{\|\sigma(\gamma)\theta\|},
       \Phi^{B\partial_t}_{\log\|\sigma(\gamma)\theta\|}(t)\right),
                 \qquad B(0)=B'(0)=0,
\end{equation}
where $B$ is analytic. The case $B\equiv0$ is included.
\end{theorem}

\begin{proof}
Use the actual curves from
Theorem~\ref{II-neutral:thm:analytic-curves}. Lift them to
physical curves $i_g$ in $M$, with $g\in G$ and
$\theta_g=g^{-1}e_1/\|g^{-1}e_1\|$. Their identities are
\begin{equation}\label{II-neutral:eq:lifted-curves}
 i_{qg}(z)=i_g(\Phi_{-\chi(q)}^{B\partial_z}(z)),\qquad
 i_{g\gamma}(z)=\alpha(\gamma)^{-1}i_g(z).
\end{equation}
For each of the finitely many coordinate directions in the
column radical, Lemma~\ref{II-neutral:lem:direction-motion}
gives analytic representatives on almost every parallel
angular line. The transformations $T_v$ preserve Haar
probability and commute with one another. Fubini therefore
allows one choice of $g_0$ for which the physical curve
field
\[
                    f(v,z)=i_{V(v)g_0}(z)
\]
has all these coordinate-line representatives on a small
angular box, as well as the original curve germs for almost
every $v$. Their common zero-section map is the analytic
parametrization $v\mapsto\phi_{g_0}(v)$ of $N$.
Apply Lemma~\ref{II-neutral:lem:iterated-cross} to obtain a
jointly analytic representative on a nonempty angular patch.
Its derivative is invertible at almost every zero-section
point: the angular derivative is that of $\phi_{g_0}$ and
the transverse derivative is that of the embedded transverse
curve $i_{V(v)g_0}$. Choose such a point and shrink. This is
an actual analytic collar patch.

It remains to descend the parameter and extend this patch.
Put $r(g)=-\log\|g^{-1}e_1\|$ and define
\[
                         J_g(t)=i_g(\Phi_{r(g)}^{B\partial_z}(t)).
\]
Since $r(qg)=r(g)+\chi(q)$,
\eqref{II-neutral:eq:lifted-curves} makes this a $Q$-invariant
measurable germ field. It descends to $J_\theta$ along $N$.
The function $r(V(v)g_0)$ is analytic, so the collar patch
just obtained is also a jointly analytic representative of
$J_\theta$. The deck identity and the norm-cocycle identity
give
\begin{equation}\label{II-neutral:eq:descended-curves}
 \alpha(\gamma)J_\theta(t)
   =J_{\gamma\theta}
        (\Phi^{B\partial_z}_{\log\|\gamma\theta\|}(t)).
\end{equation}

Translate the patch using this identity. Minimality of the
projective lattice action makes its translates cover $N$;
compactness gives finitely many. On overlaps the translated
analytic maps represent the same entire curve germs for
almost every angular center. Their Taylor coefficients in
$t$ agree on a full-measure subset of the overlap, and hence
identically by analyticity in the angular variables. Thus
the whole analytic maps agree near the zero section and glue.
Their derivatives are invertible there by transport from the
initial patch. When the angular labels cover $N$ nontrivially,
the physical map just constructed is lifted to the prescribed
collar cover, not asserted to be injective into $M$ itself.
To make this precise, choose a continuous neighborhood retraction
onto $N$. Such a retraction exists because the compact manifold $N$
is an ENR \cite[Corollary~A.9 and the preceding definition, p.~527]{Hatcher}:
extend the coordinate functions of an ENR embedding from the closed
subset $N$ to $M$ by Tietze extension, and compose with the ENR
retraction on the inverse image of its open domain. Pull back the
angular cover along this retraction, and lift the retraction composed
with $J$, using the identity angular lift at $t=0$; include the
normal-orientation cover in the same way. The covering space and its
local lifts have analytic charts from the original local analytic
inverses of $J$; the retraction need only be continuous. This gives a map into
the collar cover whose zero section is its embedded central
hypersurface. Compactness and local invertibility give an
embedded smaller uniform collar there. Its physical projection
retains exactly the original deck identifications.
The identity for every fixed $\gamma$ follows by analytic
uniqueness from \eqref{II-neutral:eq:descended-curves}. This is
\eqref{II-neutral:eq:full-collar}. The dual case follows by
inverse transpose. All finite covers used here are the
angular and normal-orientation choices already specified;
no rational point or unipotent element of $\Gamma$ was used.
\end{proof}

\begin{remark}[Descent of the local ambient action]
The angular and normal-orientation deck maps in this construction
retain their original analytic germs. They also commute with the
local ambient action of the flow model. If $B\not\equiv0$,
pass on each punctured side to the developing coordinate
$z=\exp(\int dt/B(t))\theta$. For each fixed $\gamma$, its
norm cocycle is bounded on the compact angular space. Subdivide this
bounded time interval into finitely many small flow times. Composing
the small-time representatives gives a jointly analytic family on a
common spatial band, with all intermediate iterates in their domains.
The reciprocal integral increases by the time parameter, so the radial
coordinate is multiplied by $\exp(s)$. Analytic germ uniqueness and
compactness identify this family with the original $\gamma$ model on
a common angular band. A deck comparison is an analytic
$\Gamma$-equivariant comparison between standard punctured vector
end models. Its equivariance identities extend over the connected
intersections of punctured balls, or of exteriors of ellipsoids,
by analytic uniqueness. Its lattice translates therefore glue
on a nonempty invariant open subset of punctured vector space.
The recurrence argument of
Lemma~\ref{II-exterior:lem:comparison}, which uses no collar
description input, makes it a scalar homothety. It therefore
commutes with the ambient action on the punctured sides and,
by analytic continuation, at the central hypersurface.
If $B\equiv0$, the transverse coordinate of a deck map is
$\Gamma$-invariant and hence independent of $\theta$ by
minimality. Its angular map at each fixed transverse parameter
commutes with the projective lattice action and is the identity
or the antipodal map, by Lemma~\ref{II-general:lem:angular}.
Continuity keeps that choice constant. This again commutes with
the ambient action.

For joint regularity, restrict the ambient element to a neighborhood
of the identity. Compactness bounds the norm cocycle uniformly in
that element and the angular variable, so the same finite-composition
construction gives jointly analytic representatives. The flow germ
law and analytic uniqueness give the multiplication identity jointly
in the two group elements and the collar point, near the identity
and the central hypersurface. The deck commutation germs likewise
hold jointly. Finiteness of the deck group and compactness of the
angular space give common joint domains for descent. Local analytic
inverse branches of the physical projection then prove joint
analyticity downstairs. Thus the original collar quotient inherits
the local ambient action; this descent does not require a cover
of the entire manifold.
\end{remark}

\section{Completion of the local description}\label{II-series:recognition}\label{II-sec:recognition}\label{II-series:last-section}

\subsection{The additional Grassmannian branch in dimension $4$}

For ambient actions on an $n$-manifold, every orbit has dimension
at most $n$. The other exceptional closed orbit is the Grassmannian of
two-planes in $\R^4$, with its indicated finite covers
\cite[Theorem 2.1(4) and Remark 2.2, p.~5]{FisherMelnick}.
We use this corrected orbit classification; the statement of
\cite[Theorem 3.5, p.~143]{CG} does not include this closed orbit.

\begin{theorem}[Analytic upgrade of the two-plane factor]
\label{II-loc:thm:grassmann-analytic}
Let any lattice $\Gamma<\SL(4,\R)$ act analytically on a
compact connected $4$-manifold. Suppose it admits a $C^2$
equivariant finite covering
$h:M\to\operatorname{Gr}_2(\R^4)$ which is a local
diffeomorphism. Then $h$ is analytic. In particular the smooth
finite covering supplied by \cite[Theorem 4.6]{BRHW} in the
two-plane parabolic branch has an analytic description for
arbitrary lattices.
\end{theorem}
\begin{proof}
Apply Lemma~\ref{II-general:lem:boundary-conjugacy} with $r=2$.
\end{proof}

\subsection{The local description alternatives}

\begin{corollary}[All proper-parabolic description alternatives]
\label{II-loc:cor:all-parabolic-recognition}
For $n\ge4$ and any lattice $\Gamma<\SL(n,\R)$,
consider an ergodic Haar-projecting
$P$-invariant suspension measure whose stabilizer is a
proper parabolic. Using the formal reduction, convergence, and collar constructions
proved in this paper, its local description is as follows:
\begin{enumerate}
\item For a line or hyperplane parabolic, $q>0$ and $q<0$
give compact analytic projective hypersurfaces with analytic
collar neighborhoods.
\item For $q=0$, zero radical shear gives the same collar
alternative. Nonzero radical shear gives the analytic
Hopf manifold model of
Theorem~\ref{II-loc:thm:neutral-regular-recognition}.
\item The only additional proper-parabolic possibility is
the two-plane parabolic in $n=4$; its model is an analytic
finite cover of $\operatorname{Gr}_2(\R^4)$.
\end{enumerate}
\end{corollary}
\begin{proof}
The dimension count in Proposition~\ref{II-nv:prop:parabolics} leaves exactly these
proper parabolics. Apply Corollary~\ref{II-loc:cor:positive},
Theorems~\ref{II-loc:thm:negative-recognition},
\ref{II-loc:thm:neutral-zero-recognition},
\ref{II-loc:thm:neutral-regular-recognition}, and
\ref{II-loc:thm:grassmann-analytic}, respectively.
For each recognized projective covering hypersurface, apply
Corollary~\ref{II-general:cor:normal-character} to its actual
normal derivative. According as this normal-character exponent
is negative, positive, or zero, apply
Theorems~\ref{II-suspension:thm:negative-collar},
\ref{II-positive:thm:collar}, or
\ref{II-neutral:thm:full-collar}. These constructions give a full
analytic collar after the stated finite-index restrictions.
This collar step uses the normal character of the chosen
hypersurface; the original spectral scalar selects its
recognition route.
\end{proof}

\subsection{Verification of the local description property}

\begin{proof}[Proof of Theorem~\ref{II-main:thm:recognition} for $n\ge4$]
Fix an arbitrary finite-index subgroup of the source lattice;
it is again a lattice of the same kind. All statements used
below apply to this subgroup, so the argument verifies the
finite-index quantifier in Definition~\ref{II-def:recognition}.

\emph{Clause (G).} Let $\nu$ be an ergodic invariant probability
on $M$. A finite analytic atlas on compact $M$ gives integrable
finite-jet cocycles for a finite generating set. Apply
Proposition~\ref{I-prop:formal-reduction}. The dimension-$n$ representation
calculation in Section~\ref{I-sec:standard-dual} leaves exactly the compact,
standard, and dual derivative branches. In the compact branch,
Theorem~\ref{II-loc:thm:compact-finite} gives finite image on connected
$M$. In the other branches,
Corollary~\ref{I-power:cor:manifold} gives actual analytic coordinate
germs for $\nu$-almost every center. Their compact centralizing
cocycle takes values in $\{I,-I\}$; the associated two-point
probability extension removes it. This is the precise extension
allowed in (G). The germs depend measurably on the center and have
positive radii. All these inputs hold for arbitrary lattices.

\emph{Clause (P).} Let $\mu$ be any Haar-projecting ergodic
minimal-parabolic invariant probability with proper stabilizer.
Proposition~\ref{II-nv:prop:parabolics} lists all its possible
parabolic types. Corollary~\ref{II-loc:cor:all-parabolic-recognition}
then supplies finitely many compact embedded analytic projective
covering hypersurfaces in the positive, negative, and zero-shear
neutral cases. The remaining neutral case gives the closed
Hopf manifold of
Theorem~\ref{II-loc:thm:neutral-regular-recognition}. The two-plane
case, available only in dimension $4$, gives the closed analytic
Grassmannian covering of Theorem~\ref{II-loc:thm:grassmann-analytic}.
These are exactly the alternatives in (P).

The same argument applies to the \emph{given} measure when it is
supported on $(G\times K)/\Gamma$ for a compact invariant
$K\subset M$. It is still a suspension measure for the original
analytic action on $M$; no manifold structure on $K$ is required.
In a hypersurface alternative, first include the finitely many
translates given by the finite-index member stabilizers, and retain
the members that meet $K$. This subfamily is nonempty and still carries
the given measure: almost every suspension point has a representative
on the original family and, by invariance of $K$, that same
representative lies in $K$. The retained family is invariant.
After passing to a stabilizer, the intersection of each retained
hypersurface with $K$ is closed and invariant. Minimality of
the projective action, Lemma~\ref{II-general:lem:linear-models}, makes
that intersection the whole hypersurface. Thus the recognized
central hypersurfaces lie in $K$, while their collar neighborhoods
are neighborhoods in $M$. This is the support assertion required
by the later remainder and frontier arguments; it does not assume
that those neighborhoods already cover $M$.

\emph{Clause (C).} For every hypersurface just obtained,
Theorems~\ref{II-suspension:thm:negative-collar},
\ref{II-positive:thm:collar}, and \ref{II-neutral:thm:full-collar}
give the full-angular collar, according as its normal-character
exponent is negative, positive, or zero. This exponent belongs to
the chosen hypersurface; it need not equal the spectral scalar
that selected the preceding branch.
These theorems include the standard or dual representation
label, the spherical or antipodal angular cover, the possible
normal-orientation cover, and agreement with the original lattice
action as germs along the entire hypersurface. On each nonproduct
side, integrating $dt/B(t)$ gives the full punctured-ball or
exterior developing image, since $\log r=\int dt/B(t)$. Thus (C) provides entire collar
neighborhoods, not just charts at almost every center.

All three collar constructions use the continuous groups in
the suspension. In the neutral case the auxiliary integral
semidirect products of Lemma~\ref{II-neutral:lem:direction-gap}
are subgroups of those continuous stabilizers acting on
$G/\Gamma$; they are not assumed to be subgroups of the source
lattice $\Gamma$. The finite spectral gap used there follows
from relative property $(T)$ and unipotent ergodicity.
This proves all parts of the theorem using only local results; no step
invokes the global assembly theorem.
\end{proof}

\clearpage
\section*{Part III. Local description for \texorpdfstring{$n=3$}{n = 3}}
\hypertarget{part.iii}{}
\addcontentsline{toc}{part}{Part III. Local description for \texorpdfstring{$n=3$}{n = 3}}

\section{Introduction}
The classification of lattice actions on manifolds is a central
problem in the Zimmer program \cite{Zimmer}. This paper treats
the local description needed for that problem for lattices in
$\SL(3,\R)$ acting on compact $3$-manifolds. Dimension $3$
is the first case in the series, and it requires a separate
argument because the stabilizer of a projective point has a
rank-one Levi factor.

For a projective hypersurface in a $3$-manifold, the Levi factor
is $\SL(2,\R)$. The higher-rank Levi normalization used in Part II
does not apply. We use relative property $(T)$ for the Levi and
its abelian radical, one-dimensional coordinate foliations, and
analytic realization with a transverse parameter to obtain the
local description in this dimension.

Let $\Gamma<\SL(3,\R)$ be any lattice and let $\Gamma$ act
real analytically on a compact connected $3$-manifold $M$
without boundary. We use the suspension, character convention,
and local-description property of Part II. The common geometric
tools and invariant-probability theorem there apply in dimension
three. The latter uses analytic germ cocycle superrigidity from
Part I; the neighborhood constructions below supply the additional
local assertions.

Throughout, \emph{joint analyticity} means analyticity in all
the local coordinates together, including the transverse
coordinate. Analyticity on each coordinate curve alone does
not assert this. Similarly, a family of formal charts specifies
Taylor coefficients; calling those charts analytic requires
convergence. These distinctions are essential in the neutral
case, where the linear normal weight vanishes but higher-order
transverse terms may remain.

\begin{theorem}[Local description for $n=3$]
\label{III-main:thm:recognition-three}
Every such action has the local description of
Definition~\ref{II-def:recognition}. The conclusions hold after
arbitrary finite-index restriction and for each specified
suspension measure, including a measure supported over any compact
invariant subset of $M$.

The projective hypersurfaces have the full-angular neighborhood
formula \eqref{II-eq:intro-collar}, with standard or dual angular
action and the stated finite-cover and deck data. For a nonzero
normal exponent the transverse vector field is $B(t)=qt$.
For a zero-shear neutral hypersurface it is analytic with
$B(0)=B'(0)=0$, including $B=0$. Nonzero neutral shear gives
a Hopf manifold. The additional closed homogeneous alternative
is the full flag manifold
\[
 \operatorname{Flag}(\R^3)=
 \{(\ell,\Pi):\ell\subset\Pi\subset\R^3,
                      \ \dim\ell=1,\ \dim\Pi=2\}
\]
with its indicated finite covers.
\end{theorem}

The suspension-measure and boundary inputs are due to
Brown--Rodriguez Hertz--Wang \cite{BRHW,BRHWFactors}.
The ambient models and terminology follow Fisher--Melnick
\cite{FisherMelnick}, building on Uchida \cite{Uchida} and
Cairns--Ghys \cite{CG}. Son \cite{Son} extended the
real-analytic ambient-group classification to dimensions
$n\le m\le2n-3$, including the case $n=m=3$ considered here.
Our task is to recover these models
from the lattice action. All auxiliary Lie groups act on the
suspension; no arithmetic Levi subgroup inside $\Gamma$ is assumed.

The proof first constructs analytic patches from coordinate
foliations and treats the nonzero transverse exponents. In the
neutral case, we distinguish regular shear from zero shear.
For zero shear we construct the axis, its image jets, and its
parametrized jets before establishing joint analyticity.
The final technical sections prove the curve, frame, gauge,
and parameter estimates used in this last step. The affine
descent and global assembly belong to Part IV.

\paragraph{Nonresonance on the original probability.}
The original result is the $A$-ergodic, Haar-projecting statement
\cite[Proposition~5.1, pp.~967--968]{BRHWFactors}: a root not positively
proportional to a fiber Lyapunov exponent preserves the probability.
The compact connected manifold and regularity assumptions are our standing
hypotheses. The adapted metric and all-finite-$L^q$ bounds of
\cite[Claim~2.6, p.~958]{BRHWFactors} supply the hypotheses of its native
multiplicative ergodic theorem \cite[Theorem~3.1, pp.~958--959]{BRHWFactors}.
Recurrent bounded changes of frame compare the native growth filtrations
with those of the given derivative cocycle; a countable dense set of
Cartan directions identifies the full weight multiset, with multiplicities.
For the original $P$-ergodic law, the finite-component argument of Part II
gives $A$-ergodic components with the same Haar marginal. Restricting the
given cocycle preserves its full multiset; the common-spectrum assertion
also follows from \cite[Claim~5.2 and its proof, pp.~968--969]{BRHWFactors}.
Apply Proposition~5.1 to each component and sum the resulting invariance
identities with the original component masses. This gives nonresonance
for the original law, without replacing the prescribed section or spectrum.

\paragraph{The neutral entropy application.}
For completeness, the entropy step inherited from Part II only needs the
following specialization. If every nonzero fiber weight is positive on
one common open cone, the original $P$-law has zero fiber-relative entropy
at every Cartan time. On each Cartan ergodic component use the original
compatible charts \cite[Proposition~2.5, p.~956]{BRHWFactors} and tame
fiber foliation \cite[\S3.2, p.~960]{BRHWFactors}. The standing hypotheses,
including their Lorentz integrability, pass to integer subgroups
\cite[Remark~3.3, p.~469]{BRHWSmoothErgodic2023}. Whole fibers are subordinate;
Definitions~6.1, 6.3--6.4 and Remark~6.2 of \cite[pp.~484--485]{BRHWSmoothErgodic2023}
identify the resulting leaf entropy with conditional Kolmogorov--Sinai
entropy over the original base. The raywise MET below and recurrence of
the coordinate changes identify the vertical rates. At an ergodic time
with nonpositive fiber exponents, every positive tangent block intersects
the vertical tangent trivially, so the leafwise Ruelle inequality
\cite[Theorem~7.2(22), p.~488]{BRHWSmoothErgodic2023} gives zero relative entropy.

For a nonzero Cartan vector $a$, choose independent ergodic vectors $u,v$
in the common positive cone with $a=u-v$. This follows from the
countable-hyperplane exceptional set for an ergodic $\mathbb R^2$-action
on a finite measure space with separable $L^2$
\cite[Theorem~1, p.~116]{PughShub1971}.
Apply Ruelle to $-u,-v$, then inverse symmetry and global conditional
subadditivity \cite[Theorem~13.3(1), p.~518]{BRHWSmoothErgodic2023} on a lattice
containing $u,v$: $h(a\mid\mathcal F)\le
h(u\mid\mathcal F)+h(-v\mid\mathcal F)=0$.
The identity time has zero entropy. The finite component label adds only
bounded entropy to every observation block, so division by block length
passes zero entropy to the same original $P$-law and Haar base.
Thus the neutral application follows from these separate original
statements and the component argument; no product formula on a newly
chosen law is required.

\paragraph{The original raywise MET input.}
For the raywise Cartan-derivative applications we use the original integer-action statement in
\cite[proof paragraph and preceding nonergodic discussion, p.~463,
equation~(7)]{BRHWSmoothErgodic2023}. Its hypotheses are a probability-preserving
$\mathbb Z^d$-action and a measurable invertible finite-dimensional
linear cocycle on an invariant conull domain, with
$\log^+\|A(z,\cdot)\|\in L^1$ for every integer parameter $z$.
On the completion of the original probability space it gives a
measurable joint splitting, with a possibly varying finite number of
invariant weights in the nonergodic case, and the growth and angle
limits along each integer ray. These applications use the raywise conclusion,
not the stronger arbitrary-path limits of Theorem~2.4, whose hypothesis
is the Lorentz condition $L^{d,1}$.
For each fixed real Cartan direction the compact-parameter derivative
envelope supplies the extension from integer to real times; no common
conull set over uncountably many directions is asserted.
The full native root-graph and recognition applications separately use
\cite[Theorem~3.1, pp.~958--959]{BRHWFactors} with the all-finite-$L^q$
logarithmic derivative bounds specified above.

The concrete frame and projector fields are derived from this splitting.
On each block-count stratum choose measurable adapted bases, retain
whole invariant blocks, order the weights for the fixed integer time,
and put zero in unused finite slots. The angle limits control the
projector norms in both time directions. This constructs the actual
invertible splitting on the original completed law. Taking original
measurable versions changes its fields only on null sets; removing
the countable union of their integer translates preserves the required
integer-orbit identities. There is no assumption that every block is
one-dimensional or that distinct full Cartan weights stay distinct on
an arbitrary ray.

\section{Analyticity from coordinate foliations}\label{III-series:coordinate-curves}

From this section through Section~\ref{III-series:joint-three}, let
$\Gamma$ be an arbitrary lattice in $G=\SL(3,\R)$. We use the
return estimates of Part II on the finite-volume suspension.
The root subgroups, block Levi subgroups, and auxiliary integral
semidirect products below act on that suspension; they need not
intersect $\Gamma$ in lattices. A nontrivial unipotent acts
ergodically on $G/\Gamma$ by Howe--Moore. Together with relative
property $(T)$ for $(\SL(2,\Z)\ltimes\Z^2,\Z^2)$, this gives
the spectral gap used in the scalar and curve normalizations.
The arguments that take place on the boundary itself use
projective minimality and finite abelianization, both valid
for every finite-index subgroup of $\Gamma$.

\subsection{Finite-valued boundary maps: a local patch suffices}

Let $B=\RP^2$, with its standard or dual projective action of a
lattice $\Gamma<\SL(3,\R)$. Let $\alpha$ be an
analytic action on a compact analytic $3$-manifold $M$.

\begin{proposition}[Propagation and embeddedness]
\label{III-coordinate_recognition:prop:finite-map}
Suppose a measurable equivariant assignment $b\mapsto\mathcal S(b)$
gives a set of exactly $k$ distinct points of $M$ for almost every
$b\in B$. Suppose on one nonempty open subset it has representatives
\[
             \mathcal S(b)=\{f_1(b),\ldots,f_k(b)\}
             \quad\text{almost everywhere},
\]
where the $f_i$ are analytic immersions with pairwise distinct
values at every point. Then it has such representatives locally
everywhere. Its graph is a compact analytic surface, a $k$-sheeted
covering of $B$, and its projection to $M$ is an embedding.
Consequently its image is a finite disjoint union of compact
embedded analytic projective covering hypersurfaces.
\end{proposition}
\begin{proof}
Let $O\subset B$ be the locus admitting those representatives.
It is open and $\Gamma$-invariant, since the assignment is
equivariant modulo null sets and all group maps are analytic
diffeomorphisms. Minimality on $B$ gives $O=B$.

On an overlap, two local unordered branch families agree almost
everywhere. Distinctness and continuity match the branches uniquely
near each point, and analytic uniqueness makes the match exact.
They therefore define a well-defined analytic subset
\[
           \mathcal N=\{(b,m):m\in\mathcal S(b)\}\subset B\times M.
\]
Its first projection $\pi$ is a finite analytic covering. It is
compact: in a finite local cover of the compact base, its finitely
many branches remain graphs on smaller closed neighborhoods.
The second projection $i:\mathcal N\to M$ is an analytic immersion.
For the analytic maps used here, local embeddedness follows from the
analytic inverse-chart construction already used in Part II. In charts
centered at the given point, write the map as $K:E\to F$ with
$DK(0)$ injective. Choose a linear left inverse $L:F\to E$ of
$DK(0)$, using a complement to its image. Then $D(L\circ K)(0)$
is the identity. The analytic inverse function theorem supplies a local
analytic inverse of $L\circ K$. On a sufficiently small neighborhood,
this inverse composed with $L$ is a continuous left inverse of $K$,
so $K$ is a topological embedding there. Apply this construction to the
restriction of each original branch to its original open domain. The
open inclusion transports the resulting neighborhood back to the
original source, and the map is unchanged. Thus no separate smooth
local-rank assumption is needed for these analytic branches.

The corresponding general smooth statement is the forward implication of the Local Embedding
Theorem \cite[Theorem~4.25, pp.~87--88]{LeeSmoothManifolds2e}:
if $X$ and $Y$ are finite-dimensional Hausdorff second-countable smooth
manifolds without boundary and $g:X\to Y$ is smooth with injective differential
at every point, then every $x\in X$ has an open neighborhood $W$ for which
$g|_W$ is a smooth embedding. Here a smooth embedding means a smooth map
with injective differential that is also a topological embedding; this is
the definition on p.~85 of the same reference. The source theorem also allows
boundary, which is unnecessary for the present application. Its proof first
uses the rank theorem to obtain local injectivity and then shrinks to a
precompact neighborhood to establish the topological embedding property.
The book's neighborhood convention means an open neighborhood
(Appendix~A, p.~596), so the openness in this source statement is literal.

For the open-set form used here, retain the original map $f:B\to M$ and
the given open set $U\subset B$. Give $U$ its actual restricted smooth atlas.
This open submanifold is Hausdorff and second countable. If $f$ is smooth on
$U$, then $f|_U:U\to M$ is a globally smooth map. The inclusion $j:U\hookrightarrow B$
is a smooth immersion, and the chain rule gives
$d(f|_U)_x=df_{j(x)}\circ dj_x$. Thus the assumed injectivity of $df_x$ on $U$
implies injectivity of the differential of $f|_U$ everywhere on this actual
open submanifold. Apply the stated global Local Embedding Theorem to $f|_U$.
For each original $b\in U$ it yields an open $W\subset U$ containing $b$ on
which this restricted map is a smooth embedding. Since $U$ is open in $B$,
the image $V=j(W)$ is open in $B$, contains $b$, and is contained in the
original $U$. If $x,y\in V$ and $f(x)=f(y)$, their corresponding points in
$W$ agree by the embedding's injectivity, so $x=y$. This proves the required
open neighborhood injectivity for the same original map, with no new choice
of carrier or replacement immersion.

The second projection is therefore locally injective.

We check global injectivity, rather than inferring it from immersion.
Suppose $i(x)=i(y)$. If $\pi(x)=\pi(y)$, then $x=y$ by the graph
definition. Otherwise choose a proximal $\gamma\in\Gamma$ whose
repelling projective hyperplane contains neither $\pi(x)$ nor
$\pi(y)$. Lemma~\ref{II-general:lem:pinched}, with $r=1$,
supplies a proximal element of $\Gamma$. By
Lemma~\ref{II-general:lem:linear-models}, the orbit of
its repelling hyperplane is dense in the dual projective plane.
Conjugating by an element of $\Gamma$ therefore gives the required
avoidance property. Then
$\gamma^d\pi(x)$ and $\gamma^d\pi(y)$ converge to the same
attracting point. Compactness gives a subsequence with
\[
              \gamma^{d_j}x\to x_\infty,
              \gamma^{d_j}y\to y_\infty.
\]
Equivariance and $i(x)=i(y)$ give equal second projections of
these limits; their first projections are also equal. Thus
$x_\infty=y_\infty$. Since an immersion is injective on a
sufficiently small neighborhood of that point, for large $j$
we have $\gamma^{d_j}x=\gamma^{d_j}y$, and hence $x=y$.
A compact injective immersion is an embedding.

Every connected component covers $\RP^2$, so it is $\RP^2$
or $S^2$ by the connected-cover classification
\cite[Theorem~1.38, pp.~67--68, and Example~1.43, p.~74]{Hatcher}. After passage to a finite-index subgroup the components
are invariant. On a sphere component the lift of the projective
action differs from its standard lift by a sign character;
a further finite-index restriction removes it.
\end{proof}

\paragraph{The connected covers of $\RP^2$.}
The sphere $S^n$ is simply connected for $n\ge2$
\cite[Proposition~1.14, p.~35]{Hatcher}. A linear change of coordinates
followed by radial normalization transports this to the unit sphere of
the norm used in the angular model. A map from a simply connected,
locally path-connected space has a unique lift through a covering once
the image of the basepoint is fixed in the fiber
\cite[Propositions~1.33--1.34, pp.~61--62]{Hatcher}.

Let $p:N\to\mathbb{RP}^2$
be the actual connected Hausdorff covering, and let
$q:S^2\to\mathbb{RP}^2$ be the actual antipodal quotient. Lift $q$ through
$p$ at a chosen original fiber point to obtain $f:S^2\to N$ with $pf=q$.
Since $p$ and $q$ are local homeomorphisms, so is $f$. Its nonempty image
is open, and compactness of $S^2$ makes that image closed in $N$; hence
connectedness makes $f$ surjective. If $f$ identifies one antipodal pair,
uniqueness of covering lifts implies $f(-x)=f(x)$ for every $x$.
Since the fibers of $q$ are exactly antipodal pairs, $p$ is then injective
and therefore a homeomorphism. Otherwise $f$ is injective and is a
homeomorphism from the same sphere onto $N$. These are the two familiar
cases of \cite[Theorem~1.38, pp.~67--68, and Example~1.43, p.~74]{Hatcher}.

Finally, the actual projection and quotient maps are analytic local
diffeomorphisms in the retained atlases. Local analytic lifting upgrades
the preceding homeomorphisms and their inverses to analytic maps. Thus the
original charts and covering projection are preserved; neither an analytic
cover equivalence nor a separately recognized component is assumed.

\subsection{The suspension supplies the coordinate curves}

Use the line-parabolic notation
\[
 G=\SL(3,\R),\quad Q=\operatorname{Stab}_G(\R e_1),\quad
 V=\{I+v_2E_{21}+v_3E_{31}\}\simeq\R^2,
 \qquad \beta_j=\varepsilon_j-\varepsilon_1.
\]
Let $\mu$ be an ergodic minimal-parabolic invariant, Haar-projecting
suspension probability with full stabilizer $Q$ and fiber exponents
\begin{equation}\label{III-coordinate_recognition:eq:weights}
                   \beta_2,\quad\beta_3,\quad q\varepsilon_1,
                   \qquad q\in\R.
\end{equation}
Each missing-root fiber space is one dimensional. We use the
nonuniform suspension tools and root-graph theorem of \cite{BRHW},
and the integrable analytic return bounds and holomorphic stable
construction described in Section~\ref{II-sec:neutral}. These are estimates for the
original analytic action, with no Levi-normalized coordinate input.

\begin{lemma}[Analytic missing-root parametrizations]
\label{III-coordinate_recognition:lem:root}
If the fiber disintegration is finite atomic, the unconditional
missing-root leafwise measures of $\mu$ are
Haar graphs with analytic curve parametrizations. Their derivatives
in the root parameter are nonzero.
\end{lemma}
\begin{proof}
An extension with a fixed finite number of fiber atoms has zero
relative entropy for every Cartan element. By the coarse entropy
product formula each coarse fiber entropy is zero. The first
paragraph of the proof of \cite[Proposition 7.4]{BRHW} therefore
identifies the kernel-component and unconditional root graph
measures. Proposition~7.2(3) of \cite{BRHW} supplies
the $C^r$ root embeddings, with any needed finite $r$, in the
intrinsic coarse Pesin fibers. We use the Riemannian geometry of
\cite[Section~3.2]{BRHW}; Claim~3.9 and Proposition~3.11 there
provide the intrinsic leaves and their subordinate conditionals.
For their measurable native realization, use the centered finite-regularity
coordinates of \cite[Proposition~9.9, p.~82; Proposition~9.12,
pp.~84--85; \S11.2, p.~96]{BRHW}. The last formula realizes the full
coarse leaf by its root coordinate and its embedded fiber coordinate.
Choose the auxiliary regularity above the stated critical orders before
lowering to the requested finite $r$; all physical maps and intrinsic
leaf topologies are retained in this lowering.
These intrinsic statements do not assert equality with an entire
arbitrary-metric intersection of stable sets; the original-return
and connected-component identifications must be retained when
passing to selected suspension coordinates. The Cartan-component construction
and its finite component reconstruction use the same original probability,
Haar base and specified section throughout; the component label does not
change the original conditional graph class.

For the lifted graph action, the passage from finite Haar plaque
conditionals to an invariant sigma-finite law has the following precise form.
For a free Borel action of a second-countable locally compact Hausdorff group
on a standard Borel space, choose a complete uniformly lacunary Borel section
by \cite[Corollary~1.2, p.~284; definition p.~283]{KechrisCountableSections1992},
using inversion to pass to left actions. The counting pullback of the given
sigma-finite law on the group times this section is sigma-finite; disintegrate
it over an equivalent finite transverse law. Assume the original finite plaque
atlas covers almost everywhere and every subordinate finite box has its
normalized Haar conditional. In a bounded injective chart, bounded unions
of the resulting measurable local Haar coefficients first give finite mass
on the whole chart. Applying the finite-box conditional identity again to
that chart includes any Haar-positive kernel-null region. A bounded exhaustion
and one positive finite Haar reference set identify a single measurable
coefficient, which is absorbed into the transverse law. Thus the original
measure equals the Haar product on every injective piece of the action map.
Its measurable right inverse and equivariance give invariance. All these
counting and reconstruction steps are derived on the same section and kernel;
the quoted theorem supplies only the Borel section.

The native $C^2$ Haar graph already calibrates a possible positive multiple
$c\beta_j$ before any comparison with a chosen stable model. At a time $a$
isolating that root, conditional uniqueness gives the physical return identity
$F_x\phi_x(t)=\phi_{ax}(e^{\beta_j(a)}t)$.
Differentiate at zero. Recurrence to sets where the nonzero derivative
$D_0\phi_x$ is bounded above and away from zero identifies its exponential
rate with $\beta_j(a)$, while the given full Lyapunov spectrum makes that
rate $c\beta_j(a)$. Since $\beta_j(a)<0$, this proves $c=1$.
Thus calibration uses only the native differentiable graph and its original
Cartan returns.

For the three-dimensional
root/normal spectrum the joint weights are distinct for every transverse
coefficient. Choose a Cartan direction $H$ separating them. If
$\mathcal F_{H,r}$ denotes the growth filtration of the given derivative
cocycle, its $\lambda_i$-line is
$\mathcal F_{H,\lambda_i(H)}\cap
\mathcal F_{-H,-\lambda_i(H)}$.
Transport of both filtrations therefore preserves this line; no exposed-vertex
condition on the normal weight is needed. Normalize the original weight
basis. Its invariant coordinate projectors and the integrable original
return-word bound give tempered frame and inverse-frame bounds.

At the isolating time there is one negative weight. In complexified centered
charts, apply \cite[Theorem~5.1(a)--(c), Case~III, pp.~44--45;
proof pp.~46--54]{Ruelle1980Author} to the original return after rescaling
by its integrable analytic radius. This supplies the original stable set,
its pointwise holomorphic manifold structure and forward control windows.
The real charts, measurable parametrizations, tempered radius choices and
finite-return continuation are the derived constructions of Part~II, on
this same return and probability law. The compact analytic action and
integrable original return word supply the required norm moments.
For the original controlled section, exponential integrability follows from
\cite[Theorem~5.3 (Maucourant), p.~184, proof pp.~184--185]{dLS};
the almost-minimal section construction on pp.~185--186 and the word comparison
\cite[Theorem~5.1, p.~183, proof pp.~183--184]{dLS}
give the actual return-word bounds. These are applied in ambient rank three.
Choose the stable model from this same original return construction. The
native leaf's intrinsic contraction places it in the stable set supplied by
the original Ruelle windows. Requerying those windows on the same radius
and using the intrinsic leaf topology identifies the leaf with the chosen
stable model, so its Haar parametrization and Cartan return law transport
together. This is a comparison of the coherently chosen models.

\paragraph{The original property-(T) and exponential-domain inputs.}
We use two separate classical theorems.  First, if $K$ is a local field and
$n\geq3$, then $\mathrm{SL}_n(K)$ has property~(T)
\cite[Theorem~1.4.15, p.~54, proof p.~55]{BHV}; we use its real-field instance.
Local fields are represented by their compatible absolute-value norm as in
Definition~D.4.1; the usual real field is such a field by
Example~D.4.2(i), p.~400.  Real local compactness and nondiscreteness hold for the usual real topology.
Here property~(T) has its original complex-Hilbert-space meaning: there are
a compact set $Q$ and $\varepsilon>0$ such that every strongly continuous
unitary representation $\pi$ with a nonzero vector $v$ satisfying
\[
  \sup_{g\in Q}\|\pi(g)v-v\|<\varepsilon\|v\|
\]
has a nonzero invariant vector
\cite[Definitions~1.1.1 and 1.1.3, pp.~32--34]{BHV}.
Strong continuity means continuity of every orbit map $g\mapsto\pi(g)v$.

Second, let $G$ be a locally compact group with property~(T), let
$\ell:G\to\mathbb R_{\geq0}$ be symmetric and subadditive and bounded on
compact sets, and let $\Gamma<G$ be a discrete subgroup such that $G/\Gamma$
has an invariant Borel probability.  Then there exist a Borel fundamental
domain $F$ and $s>0$ with
\[
  \int_F\exp(s\ell(g))\,d\nu(g)<\infty
\]
for Haar measure $\nu$
\cite[Theorem~5.3, p.~184, proof pp.~184--185]{dLS}.
The length definition is the one in the introduction, p.~153; the lattice
definition is on p.~152.  The term ``Borel fundamental domain'' has the
exact meaning fixed in \S2.2, p.~159: multiplication
$F\times\Gamma\to G$, $(f,\gamma)\mapsto f\gamma$, is bijective.
Thus an exact transversal is part of the primary theorem's conclusion.
For the application we use the second countable Hausdorff case and
measurable lengths; these restrictions narrow the theorem's scope.
For measurable nonnegative integrands, the finite integral is expressed
equivalently by integrability on $F$.

\paragraph{Checked specialization to the actual matrix lattice.}
For $G=\mathrm{SL}_n(\mathbb R)$, $n\geq3$, its determinant-one matrix
realization is closed in finite dimensional matrix space, hence Hausdorff,
locally compact and second countable; multiplication and inversion are
continuous.  A discrete subgroup is countable.  The previously checked
elementary-matrix character calculation proves right invariance of
canonical Haar measure.  Haar uniqueness then proves right invariance of
every given left Haar measure $\nu$.
Positive finite covolume supplies a measurable fundamental domain $D$ of
positive finite mass.  Normalize $\nu$ by $\nu(D)^{-1}$, restrict to $D$,
and project to $G/\Gamma$.  The fundamental-domain unfolding argument
proves that this is an invariant Borel probability.  Consequently the
original lattice hypothesis of Theorem~5.3 is proved from the actual
covolume assumptions; it is not an additional supplied datum.

Let $A$ be the actual compact generating family, with its range a
neighborhood of the identity.  Its word length $\ell_A$ is nonnegative,
symmetric and subadditive.  To check boundedness on a compact set, cover
that set by translates of an open neighborhood contained in the generator
range and take a finite subcover.  Every point of a chosen translate has
word length at most the center's length plus one, giving one finite
bound.  Thus the original DLS length hypotheses hold for $\ell_A$ itself;
its measurability is part of the specified word model.
Applying the two separate primary theorems now gives $F$ and $s$.
The multiplication bijection proves that $F\to G/\Gamma$ is bijective
and that its right $\Gamma$ translates cover exactly and are disjoint.
In particular it is a measure-theoretic fundamental domain.
Its mass is the covolume, hence is positive and finite.
This proves the exponential-domain statement with an exact Borel transversal
and positive finite covolume, in particular for $n=3$.
The existing measurable-section and exponential-moment arguments can
therefore be applied to this same $F$ without supplying a selected section
or strengthening the classical theorem.

Fix the root $\beta_2$. Choose a Cartan element $a$ with
\[
             \beta_2(a)<0<\beta_3(a),\qquad
             -q\bigl(\beta_2(a)+\beta_3(a)\bigr)>0
             \quad\text{if }q\ne0.
\]
Since $3\varepsilon_1=-\beta_2-\beta_3$, the remaining fiber
exponent $q\varepsilon_1(a)$ is positive when $q\ne0$ and zero
when $q=0$. Thus the entire stable fiber of $a$ is the single
$\beta_2$ curve. The analytic nonuniform stable-manifold
construction makes this curve analytic. Explicitly, in tempered
Lyapunov charts the forward stable coordinates and backward
nonstable coordinates are solved by their contracting holomorphic
sequence operator. Uniform convergence on smaller tempered
complex balls gives holomorphic dependence on the stable initial
coordinate. This is the argument for analytic stable plaques in
Section~\ref{II-sec:neutral}, now with one stable coordinate.

The root parametrization itself is also analytic. Write $T=a$,
$\lambda=e^{\beta_2(a)}<1$, and $J_x=D_0\phi_x$ for its
nonzero derivative. Conditional-measure uniqueness gives
\[
              F_x\phi_x(s)=\phi_{Tx}(\lambda s),\qquad
              D_xF_x J_x=J_{Tx}\lambda.
\]
The latter equation and the integrable return derivative bounds
make $J_x$ and its inverse onto the stable line tempered, by
the integrable-increment argument. In stable analytic coordinates
with this derivative, the return maps have linear part $\lambda$.
The one-dimensional nonstationary Koenigs construction converges
on tempered domains, since the quadratic remainder contracts
at rate $\lambda^2<\lambda$. It gives analytic maps $H_x$ with
the same derivative and covariance.

For $K_x=H_x^{-1}\phi_x$, recurrence to a set with a common
$C^2$ radius and remainder bound, together with
\[
                  K_{T^nx}(\lambda^ns)=\lambda^nK_x(s),
\]
gives $|K_x(s)-s|\le C\lambda^d|s|^2$ along unbounded
recurrent times. Therefore $K_x=\id$. Continuation by contraction
gives analyticity along the full root curve. The proof for
$\beta_3$ exchanges the indices.
\end{proof}

\begin{lemma}[Measurable graphs before joint regularity]
\label{III-coordinate_recognition:lem:graphs}
Suppose in addition that the fiber disintegration is uniform on
$k$ atoms. On almost every projective open cell the finite-valued
boundary assignment has a measurable labeling by $k$ maps
$f_i:\R^2\to M$ whose restrictions to almost every line of
each coordinate foliation have analytic representatives.
The two coordinate derivatives are linearly independent at
almost every regular point of every branch.
\end{lemma}
\begin{proof}
Lift the measure to $G\times M$, writing its conditionals as
$\mu_g$. Parabolic invariance makes $\mu_g$ a measurable function
of $Qg$. Suspension descent gives its $\Gamma$-equivariance.
The atoms thus define a measurable $k$-valued boundary assignment
on $Q\backslash G$, or, after inversion, on $G/Q$.

Take $a_1$ with both missing roots negative and consider its full
stable leaves $V\cdot W^{s,F}_{a_1}(x)$. Almost surely, two atoms
in one fiber cannot be in the same fiber stable leaf. Indeed the
relatively independent product of the finite fiber conditionals
is $a_1$-invariant; a distinct pair whose distance tends to zero
contradicts recurrence. Since there are finitely many atoms,
this excludes every such distinct pair in almost every fiber.

Projection of a full stable leafwise measure to $V$ is therefore
injective modulo null sets. Along each coordinate root line the
projected conditional is Haar, by Lemma~\ref{III-coordinate_recognition:lem:root}. Consequently
the projected measure is invariant under both coordinate
translations, and is Haar on $V$. This is the measure
disintegration step of \cite[proof of Proposition 9.11]{BRHW};
it only uses single-point fiber intersections and the coordinate
Haar conditionals. It does not use a differentiable full stable
graph or a tangential Levi chart. We obtain a measurable graph
over $V$ for each full stable leaf that meets a regular fiber.

The $k$ initial atoms give $k$ different full stable leaves and
hence the claimed labeling. Root transport keeps the label,
because each root graph lies in the same full stable leaf.
Equivalently, any discrepancy around a coordinate rectangle
would put two distinct regular atoms in one fiber stable leaf,
which was excluded above. Fubini and uniqueness of the root
conditional identify almost every coordinate slice with its
analytic root parametrization. At almost every regular graph
point the two nonzero derivatives lie in the distinct Oseledec
root lines, so they are linearly independent.
\end{proof}

\subsection{Description of the finite-atomic branches}

\begin{theorem}[Analytic hypersurfaces without Levi normalization]
\label{III-coordinate_recognition:thm:recognition}
For an analytic action of any lattice $\Gamma<\SL(3,\R)$ on a
compact analytic $3$-manifold, let $\mu$ be as in \eqref{III-coordinate_recognition:eq:weights}.
If the fiber disintegration has atoms, then $\mu$ is supported
on the suspension of a finite disjoint
union of compact embedded analytic $S^2$ or $\RP^2$ hypersurfaces.
After finite-index restriction their actions are analytically
identified with the standard projective covering action.
The plane-parabolic case gives the dual action.
No transverse collar conclusion is included.
\end{theorem}
\begin{proof}
The assumed atomic part is conull by ergodicity. The atom mass is
invariant and hence constant, so the conditional probabilities
are uniform on a fixed finite number $k$ of atoms. The root
weights in \eqref{III-coordinate_recognition:eq:weights} are the exact weights supplied
by the root-graph theorem; starting with positive multiples
$c\beta_j$, comparison of the root graph derivatives under
Cartan recurrence gives $c=1$.

Apply Lemma~\ref{III-coordinate_recognition:lem:graphs} and then Lemma~\ref{II-coordinate_recognition:lem:cross} to
the simultaneous map $(f_1,\ldots,f_k)$ into $M^k$. Equivalently,
apply Lemma~\ref{III-positive_collar:lem:cross-parameter} to this map
made constant in the normal parameter, and evaluate at zero. This gives
one rectangle on which all $k$ branches are jointly analytic.
They are distinct almost everywhere, and each has rank two
almost everywhere, by the coordinate derivative statement.
Choose a regular point where all these properties hold and
shrink the rectangle. We now have analytic immersed branches
with pairwise distinct values throughout it.

Proposition~\ref{III-coordinate_recognition:prop:finite-map} supplies the global finite
analytic covering and its embedding in $M$. It represents the
original boundary conditionals almost everywhere, so its
suspension carries $\mu$. The last assertions of that proposition
give the projective and spherical identifications, using the connected-cover
classification \cite[Theorem~1.38, pp.~67--68, and Example~1.43, p.~74]{Hatcher}.
Apply inverse transpose for the plane parabolic.
\end{proof}

\begin{corollary}[The negative and atomic neutral branches]
The conclusion of Theorem~\ref{III-coordinate_recognition:thm:recognition} holds for every
negative transverse branch and every atomic neutral branch.
\end{corollary}
\begin{proof}
Fix the given realized coefficient $q<0$ and its original derivative
field and suspension law. On the Cartan vector
$\operatorname{diag}(2,-1,-1)$ the two missing-root weights are
$-3c<0$ and the transverse weight is $2q<0$; these exhaust the
fiber spectrum. This is a positive time rescaling of
$a_t=\operatorname{diag}(e^t,e^{-t/2},e^{-t/2})$.
For the plane parabolic use its inverse-transpose Cartan vector.
The original section's word-moment estimate gives the integrable
return bound for this Cartan element. The contracting-fiber recurrence
of Lemma~\ref{II-loc:lem:atomic-fibers}, applied to that same derivative
field and law, gives uniform finite atoms in the original conditional
kernel. Apply Theorem~\ref{III-coordinate_recognition:thm:recognition}.
The atomic neutral case is its direct specialization.
\end{proof}

\begin{corollary}[The positive transverse branch]
The conclusion of Theorem~\ref{III-coordinate_recognition:thm:recognition} also holds for
every positive transverse branch.
\end{corollary}
\begin{proof}
Fix the given realized coefficient $q>0$ and its original full derivative
spectrum. The extra exponent is positive on the closed projective
contraction cone away from zero. The two missing-root coarse spaces are
one-dimensional and, together with the transverse space, exhaust the
fiber spectrum. The original $P$-ergodic probability is $Q$-ergodic because
$P\subset Q$ and $Q$ preserves it. The spectrum here agrees, with
multiplicities, with the native fiberwise spectrum of
\cite[\S3.2 and Proposition~3.3, pp.~20--22]{BRHW}: a measurable change
from the chosen coordinate frame to the native metric frame intertwines
the derivatives, and recurrence on sets where the frame and its inverse
are bounded identifies their existing growth rates. Comparing the full
growth filtrations retains the packet dimensions and hence the
one-dimensional coarse spaces. These are the hypotheses of
\cite[Theorem~4.7, p.~30]{BRHW}. For the simply connected convention of
\cite[\S2.1, p.~9]{BRHW}, pull the lattice and its action back to the
universal covering group of $G$ supplied by
\cite[Theorem~7.7, p.~154]{LeeSmoothManifolds2e}. Its kernel is discrete
and central: conjugation maps the connected covering group continuously
into the discrete kernel and is therefore constant. The covering differential
identifies its Lie algebra with $\mathfrak{sl}(3,\R)$, preserving simplicity,
splitness and the rank-two diagonal Cartan.\footnote{The ordinary Lie inputs
used here are universal statements: $\det(\exp X)=\exp(\operatorname{tr}X)$
for complex matrices; the local exponential chart for every closed matrix
group; smoothness of continuous Lie-group homomorphisms; and preservation
of brackets by their identity differentials
\cite[Proposition~0.11(e), pp.~7--8; Theorem~0.15, pp.~12--13;
\S I.10, pp.~72, 77]{KnappBeyondDigital}.
Comparing the two atlases on the same matrix group and composing with the
covering differential gives the stated identification. The other universal
inputs say that $\mathfrak{sl}(n,\C)$ is simple for $n\ge2$, every real form
of a complex simple algebra is simple, and $\mathfrak{sl}(n,\R)$ is split
\cite[\S II.1, pp.~124--126; proof of Proposition~7.1, p.~434;
examples after Corollary~6.10, p.~353]{KnappBeyondDigital}. The two independent
commuting diagonal directions give the real-rank lower bound.}
The lifted lattice is discrete, and the quotient homeomorphism transports
the original invariant quotient probability to its quotient; the lifted
suspension probability projects to precisely this same quotient probability.
Uniqueness of pointed lifts preserves the root-flow conjugation identities,
including their signs and time parameters. The kernel acts trivially on $M$;
the resulting suspension identifies with the original suspension, with
its original probability and fiber derivative. The pulled-back parabolic
quotient identifies with the same projective boundary. Thus the resulting
factor and lattice equivariance descend to the original action.

The antipodal sphere cover and simple connectedness of $S^2$
\cite[Proposition~1.14, p.~35]{Hatcher} give finite fundamental group
of $\RP^2$: a loop class is determined by the endpoint of its lift, which
has only two possibilities. Theorem~4.7(1),(4),(5),(9) therefore provides
an invariant finite union $K$ of compact embedded covering leaves and one
$\Gamma$-equivariant factor $h:K\to G/Q$, with the original suspension
probability carried by $(G\times K)/\Gamma$. Leafwise continuity descends
through the finite compact union, so $h$ is continuous on $K$ and has a
Borel extension to $M$. Every fiber on $K$ is a finite union of compact
discrete covering fibers, hence finite.

To identify the original conditional law, use this same factor to define
$F([g,m])=g\cdot h(m)$ on the suspension of $K$. Lattice equivariance makes
$F$ well defined, and $F_\ast\mu$ is a $Q$-invariant probability on $G/Q$.
In the standard projective coordinates it is concentrated at the marked
axis: recurrence under each elementary unipotent in the line parabolic
forces a representative vector to be fixed, since a moving unipotent
projective orbit converges to its highest nonzero polynomial coefficient.
The row unipotents force all coordinates other than the marked coordinate
to vanish. The dual coordinate dictionary gives the same argument for the
hyperplane parabolic. Consequently $F([g,m])=Q$ almost everywhere.
Lift this equality through the original suspension quotient and
disintegrate its actual counting--Haar lift. With its fundamental-domain
normalization retained, every whole conditional version satisfies
$h(m)=g^{-1}Q$ almost everywhere in the conditional at $g$.
Thus the original fiber conditionals are finite atomic. Apply
Theorem~\ref{III-coordinate_recognition:thm:recognition} to the original measure.
\end{proof}

\begin{proposition}[A transverse analytic curve for nonzero exponent]
\label{III-coordinate_recognition:prop:normal-curve}
For the weights \eqref{III-coordinate_recognition:eq:weights} with $q\ne0$, there is a
measurable Cartan-equivariant field of analytic curve germs
tangent to the transverse Oseledec line. This assertion does not
include Levi or parabolic equivariance at the exceptional negative
values $q=-3/j$, $j\ge1$. At all other nonzero values, the curve
field is in fact equivariant under the full line parabolic.
\end{proposition}
\begin{proof}
Fix the given nonzero coefficient and its original derivative field,
section and probability. The preceding atomicity and root-scale arguments
normalize the two root weights in this same field; they do not select a
new normal coefficient. The stable-graph constructions of Part~II apply
to its measurable scalar frame. The original section word moment supplies
integrable analytic return bounds, and the proved bilateral frame estimates
supply the required tempered radii. The change to finitely many controlled
charts preserves the original cocycle and its Lyapunov weight spaces.

If $q>0$, take the inverse of
$a_1=\operatorname{diag}(e,e^{-1/2},e^{-1/2})$.
Its two tangential exponents are positive and its transverse
exponent is negative. Its analytic fiber stable curve is the
required germ.

If $q<0$, choose commuting Cartan elements $a,b$ by prescribing
\[
 (\beta_2(a),\beta_3(a))=(-2,1),\qquad
 (\beta_2(b),\beta_3(b))=(1,-2).
\]
Both have $\varepsilon_1=1/3$, so the transverse exponent is
$q/3<0$. Their fiber stable manifolds are analytic surfaces with
tangent spaces
\[
 E_{\beta_2}\oplus E_{q\varepsilon_1},\qquad
 E_{\beta_3}\oplus E_{q\varepsilon_1}.
\]
These surfaces meet transversely at the center. The analytic
implicit function theorem makes their local intersection a
nonsingular analytic curve tangent to $E_{q\varepsilon_1}$.
Measurable stable charts give measurable germ representatives
and positive local radii. Every Cartan element transports the
two stable manifolds to their counterparts at the new center,
since it commutes with $a,b$. Hence it transports the intersection
germ as well. The same argument gives Cartan equivariance in
the positive case.

The recognized boundary from
Theorem~\ref{III-coordinate_recognition:thm:recognition} has a parabolic-invariant
tangent plane. For $q>0$, central nonresonance gives Levi coherence of
the constructed curve, and Proposition~\ref{III-neutral_shear:prop:curve-coherence}
gives full positive-parabolic coherence. For $q<0$ outside the displayed
exceptional set, the same two steps apply. Indeed
Lemma~\ref{III-neutral_shear:lem:curve-Levi} excludes $kq=-3/2$ for $k\ge2$,
and the row argument excludes $kq=-3$ for $k\ge1$.

Here is the passage from the original derivative field to those block
hypotheses. Use the finite normal subgroup and sphere or projective models
already constructed by recognition. On the canonical finite lift of the
same probability, align the section with the original section and use the
boundary tangent frame together with the same normal Oseledec vector.
A Levi Cartan element with root weights $1,-1$ and normal weight zero
shows, by recurrence, that this vector is transverse to the boundary.
The derivative on the tangent plane is the actual standard Levi derivative.
The central flow $\operatorname{diag}(e^{2t},e^{-t},e^{-t})$ has rates
$-3,2q$ in this frame: it is the paper's $a_{2t}$. Thus the integer-return
versions of the two coherence arguments have precisely the exceptional
values already stated. They use the original scalar exponent along integer
returns and require no additional all-real-time tempering assertion.
The same central comparison first removes the off-diagonal normal column.

The compatible topology needed for the row argument is supplied by
\cite[Theorem~1.1 and its reformulation, p.~335]{BeckerKechris1993}; it preserves
the original Borel action and lifted law, as in
Lemma~\ref{III-neutral_shear:lem:components}. The sign element also commutes with
the central flow, so the same comparison gives coherence on the other
component of the full parabolic. Undoing the common frame returns the
original curve. Finally the aligned finite-cover factor is measure
preserving and has measurable images, so its conull covariance statements
descend to the original probability. The inverse-transpose dictionary gives
the plane-parabolic case with the same coefficient.
\end{proof}

This intersection uses two analytic stable surfaces, each containing
the normal direction. It does not require a Cartan element that
contracts the whole tangential plane while expanding the normal
line. A Levi element generally changes the chosen Cartan pair in
the negative case, so that construction alone does not establish
Levi equivariance there; the nonresonance argument is additional.
Even away from the exceptional values, equivariance
of the transverse curve does not normalize its scalar Levi
germ cocycle or prove joint analytic dependence along the boundary.
Thus this curve construction alone does not supply a full collar.
The constructions in
Section~\ref{III-series:relative-three}, Section~\ref{III-series:negative-three}, Section~\ref{III-series:positive-three} go further:
the negative ambient normalization removes the exceptional
negative weights, and the two collar arguments supply full
collars for every $q\ne0$.

The coordinate foliations here are the two foliations of the
\emph{projective local chart} by lines parallel to the missing
root subgroups. We first obtain analyticity along their individual
curves, then obtain a common surface chart. No one-dimensional
quotient of the ambient $3$-manifold, and no already existing
transverse analytic foliation, is assumed.

The finite-valued boundary argument applies to atomic fiber
conditionals. For nonatomic neutral conditionals, the root-curve
comparison in Section~\ref{III-series:regular-three} supplies the local
chart. The neutral hypersurface case is treated in
Sections~\ref{III-series:axis-three}--\ref{III-series:joint-three}.

\section{Relative normalization for \texorpdfstring{$n=3$}{n = 3}}\label{III-series:relative-three}

\subsection{Relative property \texorpdfstring{$(T)$}{(T)} with an ergodic radical}

Write $L=\SL(2,\R)$, $U=(\R^2)^*$, $H=L\ltimes U$, and
$Q=H\rtimes A_0$, with $A_0=\{a_t\}\simeq\R$ centralizing
$L$ and dilating $U$. We use the cocycle convention
$c(gh,x)=c(g,hx)c(h,x)$. Both pairs
\[
 (H,U),\qquad
 (H_0,U_0)=\bigl(\SL(2,\Z)\ltimes\Z^2,\Z^2\bigr)
\]
have relative property $(T)$. We use the two concrete standard-pair
assertions stated in \cite[Introduction, pp.~375--376; Definition~1.1,
p.~375]{Iozzi}: every strongly continuous unitary representation with
almost invariant unit vectors on compact sets has a nonzero vector fixed
by the indicated translation subgroup. Here almost invariance means that
for every compact $K$ and $\varepsilon>0$ there is a unit vector $v$ with
$\sup_{g\in K}\|\rho(g)v-v\|<\varepsilon$.
The pointwise-on-$K$ convention gives this literal supremum convention by
choosing the pointwise bound $\varepsilon/2$, whose supremum is at most
$\varepsilon/2$; the empty compact set is harmless.
Writing the normal factor first, the standard and dual semidirect products
are related by $(v,g)\mapsto(v,\tau(g))$, where
$\tau(g)=(g^{-1})^{\mathsf T}$ is the continuous involutive
contragredient automorphism. Under the fixed row/column identification,
the dual action of $\tau(g)$ on $v$ is the standard action of $g$ on $v$.
Thus this homomorphism is the identity on the translation subgroup.
It preserves the integer matrices and integer translations as well.
Images of compact sets under this homomorphism are compact. Pulling back
a representation therefore preserves almost invariant unit vectors, and
the fixed-vector conclusion concerns the same translation vectors.
No property $(T)$ of $L$ is claimed. All probability actions below are jointly
measurable on standard spaces.

\paragraph{Strong continuity of the Koopman representation.}
Let $G$ be a Hausdorff locally compact $\sigma$-compact group with its Borel
sigma algebra, and let
$(X,\mu)$ be a $\sigma$-finite measure space with a jointly measurable,
measure-preserving action. Assume that the actual space $L^2(X,\mu;\mathbb C)$
is separable. The canonical inverse-composition representation is strongly
continuous. No countable-generation assumption on the sigma algebra of $X$
is needed.

For $f\in L^2$ choose a measurable representative $R$ and put
$K(g,x)=R(g^{-1}x)$. Joint measurability and parameter integration give
\[
 d(U_gf,w)=\sqrt{\int_X |K(g,x)-w(x)|^2\,d\mu(x)}
\]
as a measurable function of $g$ for every $w\in L^2$.
As in the isometric-orbit argument of Part~II, a countable dense set of
centres covers the inverse orbit by balls of radius $\varepsilon/2$; one
preimage has positive Haar measure. Steinhaus supplies a neighbourhood of the
identity in its quotient set, and the isometry and triangle identities bound
the orbit displacement there by $\varepsilon$. Haar is sigma-finite under
sigma-compactness and finite measure on compacts, and its regularity supplies
the inner regularity required by Steinhaus. Translation gives continuity at
every group element. This is the same inverse-composition representation
used in Part~II. See also
\cite[Proposition~A.6.1 and Lemma~A.6.2, p.~326 of the author manuscript]{BHV}
for the separable Hilbert-space criterion.

\begin{lemma}[Relative additive reduction]\label{III-relative_normalization:lem:additive}
Suppose $H$ acts preserving a probability and $U$ acts ergodically.
Every finite measurable additive cocycle
$c:H\times X\to\R^D$ has the form
$c(h,x)=b(hx)-b(x)$ for a measurable $b:X\to\R^D$.
\end{lemma}
\begin{proof}
Use the skew-translation unitary representation on
$\mathcal H=L^2(X\times\R^D)$,
\[
 (\pi(h)F)(x,v)
  =F(h^{-1}x,v-c(h,h^{-1}x)).
\]
The preceding internal continuity argument applies to this actual skew
action. An alternative classical route uses the original measurable-unitary
homomorphism lemma \cite[Lemma~A.6.2, p.~326 of the 23 February 2007
manuscript]{BHV}. It states that a homomorphism from a locally compact
group to the unitary group of a separable complex Hilbert space is strongly
continuous whenever every diagonal matrix coefficient is measurable.
Here the actual skew action preserves the sigma-finite product measure.
More generally, for any jointly measurable measure-preserving action on
$(Y,\nu)$ and $f\in L^2(Y,\nu)$, choose a measurable representative $F$.
Then
\[
 \langle\pi(g)f,f\rangle
 =\int_Y\langle F(g^{-1}y),F(y)\rangle_{\mathbb C}\,d\nu(y).
\]
The integrand is jointly measurable; the parameter-integral theorem
therefore makes this coefficient measurable. The equality uses the actual
inverse-composition representative, so the lemma applies to this same
representation. Separability is supplied by the standard sigma-finite
measure space. This proves the needed measure-preserving specialization
without assuming a separately chosen weighted representation.

Let $\xi_R$ be the normalized indicator of $[-R,R]^D$,
constant in $x$. Its matrix coefficient is
\[
 \langle\pi(h)\xi_R,\xi_R\rangle
 =\int_X\prod_{j=1}^D
       \left(1-\frac{|c_j(h,x)|}{2R}\right)_+\,d\mu(x).
\]
This tends increasingly to one. The representation is strongly
continuous by the preceding internal argument on the same separable
Hilbert space. Thus these coefficients are continuous, and
Dini's theorem makes convergence uniform on compact subsets
of $H$. The vectors are almost invariant.

Relative property $(T)$ gives a nonzero $U$-invariant vector
$F\in\mathcal H$. Its fiber norm
$m(x)=\int|F(x,v)|^2\,dv$ is $U$-invariant, hence a positive
constant almost everywhere. Normalize $|F(x,\cdot)|^2$ to a
probability density. Its coordinate marginal lower medians
give a measurable vector $b(x)$ with
$b(ux)-b(x)=c(u,x)$ for every fixed $u\in U$.

Subtract this coboundary on all of $H$. For the residual cocycle
$c'$, normality of $U$ gives
\[
 c'(h,ux)=c'(hu,x)=c'((huh^{-1})h,x)=c'(h,x).
\]
Ergodicity of $U$ makes $c'(h,x)=\chi(h)$ deterministic.
The jointly measurable cocycle identity makes $\chi$ a measurable,
hence continuous, homomorphism $H\to\R^D$. Its restriction to
$L$ vanishes because $L$ is perfect. Its restriction to $U$
already vanishes. Thus $\chi=0$.
\end{proof}

\begin{corollary}[The central character is the only scalar remainder]
\label{III-relative_normalization:cor:scalar}
Under the same ergodicity hypothesis, every additive cocycle on
$Q$ has the form
\[
 c(q,x)=b(qx)-b(x)+\chi(q),\qquad
 \chi(a_t)=t\nu,\quad \chi(H)=0
\]
for a constant $\nu\in\R^D$.
\end{corollary}
\begin{proof}
Normalize on $H$ by Lemma~\ref{III-relative_normalization:lem:additive}. Since $H$ is
normal in $Q$, the same computation makes each remaining
$c'(q,\cdot)$ $H$-invariant, hence constant. It factors through
$Q/H\simeq\R$ and is a measurable homomorphism.
\end{proof}

\begin{lemma}[A countable spectral-gap action]\label{III-relative_normalization:lem:gap}
If $U_0=\Z^2$ acts ergodically, the probability action of
$H_0=\SL(2,\Z)\ltimes\Z^2$ has spectral gap.
\end{lemma}
\begin{proof}
On $L^2_0(X)$ there are no $U_0$-invariant vectors. The literal integer-pair assertion in
\cite[Introduction, p.~376; Definition~1.1, p.~375]{Iozzi}
therefore excludes almost $H_0$-invariant unit vectors.
The Koopman representation is continuous because $H_0$ is discrete;
using $\varepsilon/2$ converts pointwise uniform bounds on a compact set
to the source's strict supremum bound. Equivalently, there are a finite symmetric set
$S\subset H_0$ and $a>0$ such that
\[
 \sum_{s\in S}\|f\circ s-f\|_2^2\ge a\|f\|_2^2
 \qquad(f\in L^2_0(X)).
\]
This is a property of this action, not property $(T)$ of $H_0$.
\end{proof}

\subsection{Why a recognized boundary has the required ergodicity}

Let $G=\SL(3,\R)$ and let $\Gamma<G$ be any lattice. Suppose an invariant analytic hypersurface $N$
has been equivariantly identified with $\RP^2$ or $S^2$, carrying
the standard projective action or its lift. The dual case is
obtained by inverse transpose. Pass to finite index to preserve
each connected hypersurface and each chosen lift.

\begin{proposition}[Homogeneous boundary suspension]
\label{III-relative_normalization:prop:homogeneous}
After the indicated finite-index restrictions, each branch of
the line-parabolic suspension supported on $N$ is measurably
$Q^\circ$-isomorphic to $G/\Gamma$ with Haar probability.
In particular $U$ and $U_0$ act ergodically.
\end{proposition}
\begin{proof}
For the matrix-coefficient input we use the connected noncompact
simple real group with finite center case of Howe--Moore; a precise
statement and a new proof appear in
\cite[Theorem~1.1 and Definition~2.6, pp.~66 and 69; the concrete
$\SL(n,\R)$ case in \S4.1 and Example~4.1, pp.~74--75]{Ciobotaru}.
The source assertion quantifies over every complex Hilbert space and every
strongly continuous unitary representation $\rho$ with no nonzero
$G$-fixed vector: for every $v,w$, the coefficient
$\langle\rho(g)v,w\rangle$ tends to zero as $g$ leaves every compact
subset of $G$. We use its $G=\SL(3,\R)$ instance; distinguishing the first
matrix coordinate from the other two only fixes the row and column labels.
More generally, for any finite set of matrix labels, the reindexing is a
topological group isomorphism to the corresponding literal $\SL(n,\R)$.
Pull the same representation back along its inverse. Continuity is preserved
by composition, and surjectivity preserves the absence of fixed vectors.
The isomorphism is a proper homeomorphism, so pulling the coefficient limit
back preserves escape from compact sets. The Hilbert space and both
coefficient vectors remain unchanged throughout this coordinate transport.
It is applied to the actual zero-mean $L^2$ representation of the normalized
Haar probability on $G/\Gamma$, rather than assuming a decay statement for
a selected suspension law. The quotient action preserves this probability,
its unitary representation is strongly continuous, and its $G$-fixed
functions are constant. A vector fixed by an unbounded subgroup is fixed
along a sequence escaping compact sets; its self-coefficient is then its
constant squared norm, while the theorem forces that coefficient to tend
to zero. Thus this vector is zero. This gives the needed subgroup
invariance conclusion on the same Haar quotient.

For $N=\RP^2$, the branch consists of
\[
 \bigl\{[g,d]\in(G\times N)/\Gamma:
                         d=[g^{-1}e_1]\bigr\}.
\]
The map $[g,d]\mapsto g\Gamma$ is a bijection, with the displayed
formula as inverse. It intertwines the left parabolic action,
since the parabolic preserves $[e_1]$. The Haar-projecting
probability is therefore Haar under this identification.

On the spherical cover a lift of a projective transformation
differs from its standard lift by the antipodal deck map. These
choices define a sign character of $\Gamma$. Pass to its kernel.
The two branches are then
$d=\pm g^{-1}e_1/\|g^{-1}e_1\|$, each a copy of $G/\Gamma$,
and are preserved by the connected parabolic. The same description
applies after passing to the spherical cover of an $\RP^2$
hypersurface and to its normal-orientation cover: the normal line
bundle on $S^2$ is orientable: $S^2$ is simply connected
\cite[Proposition~1.14, p.~35]{Hatcher}, so its orientation double cover
has a continuous section. An orientation reversal by a
group element is constant over this connected sphere. A further
sign kernel removes it.

The Howe--Moore theorem \cite[Theorem~1.1 and Definition~2.6, pp.~66 and~69]{Ciobotaru} makes the action of every
unbounded closed subgroup of $G$ on $G/\Gamma$ ergodic. Both
$U$ and its integer subgroup $U_0$ are such subgroups.
\end{proof}

\paragraph{Analytic metrics from the regular embedding theorem.}
\label{III-relative_normalization:analytic-line-metric}
Every finite-dimensional Hausdorff real-analytic manifold with countable
 topology admits a regular proper real-analytic map into some $\mathbb R^N$.
Here regular means injective with injective differential everywhere.
This is \cite[Theorem~3, p.~470; definition of regularity, p.~462]{Grauert}. There is no compactness hypothesis. The proper injective
map is a topological embedding. We retain the original properness conclusion
in the formal source interface even though the following argument uses only
regularity and analyticity.

The local analytic immersion normal form used below follows from this
regularity statement; it is not an additional classical assumption.
In charts centered at a point, write the derivative as an injective linear map
$A:E\to F$ and choose a complementary subspace $K$ to its range.
The analytic map $(x,k)\mapsto e(x)+k$ has invertible derivative
$(v,k)\mapsto Av+k$ at the center. The analytic inverse function theorem
therefore supplies local coordinates in which the original map is
$x\mapsto(x,0)$. Restricting these coordinates to the original manifold
charts gives its analytic immersion normal form, with the same map and
point. The formal derivation reuses the completed analytic split-immersion
and inverse-function constructions, and introduces no supplied normal-form
certificate. Thus compact applications use the same regular proper
embedding theorem, with compactness imposed only by their own statements.

Let $B$ be a finite-dimensional Hausdorff second-countable real-analytic
manifold. Let $U_i$ be an open cover with nowhere-zero real-analytic transition
functions $t_{ij}$, with $t_{ii}=1$ and $t_{jk}t_{ij}=t_{ik}$ on triple overlaps.
There are positive real-analytic functions $h_i$ on $U_i$ satisfying
$h_i=t_{ij}^2h_j$ on every overlap.
Glue the local products $U_i\times\mathbb R$ by
$(x,s)_i\sim(x,t_{ij}(x)s)_j$. The transition maps and their inverses are
analytic, so the usual vector-bundle product atlas makes this actual total
space $L$ an analytic manifold. The Hausdorff property follows by separating
unequal base points, or separating two distinct fiber coordinates in a common
trivialization. A countable subcover can be extracted from the base cover,
and each local product is second countable; hence $L$ is second countable.
Its model is the product of the original base model and $\mathbb R$.
This construction applies to the given atlas and makes no independent choice
of bundle coordinates to be identified later.

Apply Grauert's theorem to $L$, which need not be compact, obtaining an analytic
immersion $e:L\to\mathbb R^N$. Write $\sigma_i(x,s)$ for the inverse local
trivialization. Define
\[
 v_i(x)=\left.\frac{d}{ds}\right|_{s=0}e(\sigma_i(x,s)),
 \qquad h_i(x)=\sum_{a=1}^{N}v_i(x)_a^2.
\]
Partial differentiation of a real-analytic function preserves analyticity,
so $v_i$ and $h_i$ are analytic. The tangent to $s\mapsto\sigma_i(x,s)$
at zero is nonzero, since $\sigma_i$ is a local analytic diffeomorphism;
regularity of $e$ therefore implies $v_i(x)\ne0$ and $h_i(x)>0$.
On overlaps the actual gluing identity is
$\sigma_i(x,s)=\sigma_j(x,t_{ij}(x)s)$. Differentiate this identity in $s$
with $x$ fixed. The chain rule gives $v_i(x)=t_{ij}(x)v_j(x)$,
whence $h_i(x)=t_{ij}(x)^2h_j(x)$. These are precisely the coefficients and
overlap convention of the original line-metric record.

\begin{corollary}[The normal derivative at every exponent]
\label{III-relative_normalization:cor:normal-derivative}
For a recognized boundary, the absolute derivative on its normal
quotient line is measurably cohomologous to a character trivial
on $H$ and equal to $e^{qt}$ on $a_t$, where $q$ is its normal
Lyapunov exponent. On the normal-orientation cover the derivative
itself has this form. This conclusion includes $q=0$ and does
not require an invariant transverse curve or a splitting of the
normal quotient into the ambient tangent bundle.
\end{corollary}
\begin{proof}
Use the actual inclusion $i:C\to M$ of the recognized component and its
given analytic identification $e$ with the spherical or projective model.
The spherical parametrization is $f=i\circ e^{-1}$; in the projective case
use $f=i\circ e^{-1}\circ\pi$, where $\pi:S^2\to\RP^2$ is the canonical
analytic double cover. These are analytic immersions of the same original
component. Their equivariance follows from that of $i$ and $e$, and from
the equivariance of $\pi$. The simply connectedness of $S^2$
\cite[Proposition~1.14, p.~35]{Hatcher} removes the normal-orientation
obstruction on this cover, as in the preceding homogeneous-branch argument;
a normal frame is constructed rather than supplied as an extra hypothesis.
Choose an analytic norm on the actual normal quotient line. Its construction from the given line atlas is the vertical-differential
argument above, applied to Grauert's theorem on the actual total line bundle.
Keep this one constructed normal frame when changing coordinates or
passing to the dual action. If $E:G_0\to G$ is the fixed group
identification, use $\Gamma_0=E^{-1}\Gamma$ and the same original
section transported by $r_0=E^{-1}\circ r\circ\overline E$, where
$\overline E:G_0/\Gamma_0\to G/\Gamma$ is the induced quotient map.
Its return labels satisfy
$E(\kappa_0(p,x))=\kappa(Ep,\overline E x)$.
Pull back the original law by the same suspension equivalence and
transport the clock with the action. The scalar normalization may
construct its gauge directly on the original carrier, using the pulled
base coordinate; it does not require a second choice of normal frame or
probability.

For precision, let $\sigma(x)>0$ be the original quotient norm of this
frame, and write $T=a_1$ and $c^{(n)}(x)$ for its original return
coefficient. On the one-dimensional fiberwise normal quotient,
\[
 \log\|D^\nu T^n(x)\|
 =\log|c^{(n)}(x)|+\log\sigma(T^nx)-\log\sigma(x).
\]
Thus the norm gauge is the constructed function $\log\sigma$; the
exponent hypothesis is the growth of the original quotient derivative,
with the same section, law and clock. Tangential components are removed
by the quotient, rather than estimated using an ambient normal vector.

Apply Corollary~\ref{III-relative_normalization:cor:scalar} to the logarithm of the absolute
normal derivative, using Proposition~\ref{III-relative_normalization:prop:homogeneous}.
If the resulting character is $e^{\nu t}$, recurrent returns to
sets where the scalar gauge and its reciprocal are bounded show
that its exponent is the original normal exponent: $\nu=q$.
Here the orientation passage preserves this same exponent and probability
construction. Let $\epsilon:\Gamma\to\{\pm1\}$ be the sign of the
constructed normal coefficient and put $\Theta=\ker\epsilon$.
Write $p$ for the subgroup suspension factor. Lift the original probability
uniformly over its finitely many sheets; its pushforward by $p$ is the
original law. Use the aligned section
$r_\Theta(x)=r(\bar p x)\eta(x)$, with $\eta(x)\in\Gamma$.
Its return labels satisfy
\[
 \eta(ax)\kappa_\Theta(a,x)
       =\kappa_\Gamma(a,\bar p x)\eta(x).
\]
If $b$ is the absolute logarithmic transfer already constructed, define
\[
 \beta(x)=\log\bigl|c(\eta(x))(r_\Theta(x)^{-1}o)\bigr|,
 \qquad B(z)=b(pz)-\beta(\operatorname{base}z).
\]
The displayed return identity and the original scalar cocycle identity
give the same character with transfer $B$ on the lifted law.
The restricted coefficient is positive because $\Theta=\ker\epsilon$,
so its absolute value may be removed. No exponent or MET assertion on
the cover is needed. Under the fixed native or inverse-transpose group
identification, pull the character back by the same subgroup isomorphism;
it is still trivial on the corresponding affine subgroup and has value
$e^{qt}$ on the transported paper clock. This also retains the character
one when $q=0$.
\end{proof}

\subsection{Analytic normalization on an actual transverse curve}

Let $\mathcal D_+^\omega$ denote the group of oriented analytic
diffeomorphism germs of $(\R,0)$.

\begin{theorem}[One-dimensional reduction on the homogeneous branch]
\label{III-relative_normalization:thm:curve-normalization}
Suppose $Q$ acts probability preservingly, $U$ and $U_0$ act
ergodically, and $c:Q\times X\to\mathcal D_+^\omega$ is an
actual coefficient-measurable germ cocycle. Using the spectral-gap
realization theorem in Section~\ref{II-sec:compact}, there are measurable
analytic coordinate germs $k_x$ and a single analytic vector field
$B(z)\partial_z$, with $B(0)=0$, such that
\[
 k_{hx}c(h,x)k_x^{-1}=\id\quad(h\in H),\qquad
 k_{a_tx}c(a_t,x)k_x^{-1}=\Phi_t^{B\partial_z}.
\]
The coordinate germs and their inverses have positive measurable
radii. Their dependence on $x$ is only asserted to be measurable.
\end{theorem}
\begin{proof}
Apply Lemma~\ref{III-relative_normalization:lem:additive} to $\log c(h,x)'(0)$ on $H$.
A measurable positive linear gauge makes this restriction
tangent to the identity. Successively remove its first untreated
Taylor coefficient. At every order that coefficient is an
ordinary real additive cocycle on $H$, so the same lemma gives
the next correction. The corrections leave lower orders unchanged.
This constructs one measurable formal tangent-to-identity field
$K_x$ satisfying the forward coboundary identity for all of $H$.

Restrict to $H_0$. Lemma~\ref{III-relative_normalization:lem:gap} supplies the scalar spectral
gap assumed by the realization theorem in Section~\ref{II-sec:compact},
Section on spectral-gap realization of analytic germ coboundaries.
That theorem gives a \emph{constant} formal germ $C$ for which
$K_x\circ C$ is analytic almost everywhere, with positive radii
for it and its inverse. Right composition by a constant preserves
the original forward identity for every fixed $h\in H$.
Combine its inverse with the positive linear gauge to obtain $k_x$.

Write $c'$ for the normalized cocycle on $Q$. Normality of $H$
and $c'|_H=\id$ imply
$c'(a_t,hx)=c'(a_t,x)$. Ergodicity makes this germ deterministic
for each fixed $t$. The cocycle identity becomes a
coefficient-measurable homomorphism from $\R$ to
$\mathcal D_+^\omega$. The one-parameter germ lemma in
Section~\ref{II-sec:compact} gives its analytic generator $B\partial_z$.
\end{proof}

\begin{corollary}[Normalized transverse parameters at nonexceptional weights]
\label{III-relative_normalization:cor:application}
On every recognized boundary with $q>0$, or with
$q<0$ and $q\notin\{-3/j:j\ge1\}$, the equivariant analytic
curve fields in Section~\ref{III-series:shear-three}, Section~\ref{III-series:coordinate-curves} have analytic
parameters in which $H$ acts trivially and $a_t$ acts by
$z\mapsto e^{qt}z$, after the finite orientation choices above.
\end{corollary}
\begin{proof}
Use the canonical probability on the finite orientation cover, whose
factor is the original probability. The finite lift has the Haar base marginal.
For an invariant leaf or signed branch, its conditional mass on the Haar base
is radical-invariant. Radical ergodicity makes that mass almost everywhere
constant, equal to the total mass of the leaf or branch. Restriction therefore
multiplies the Haar marginal by this constant, and normalization restores
exactly the Haar probability whenever the mass is positive.
Apply the normalization separately on each positive-mass signed homogeneous branch. Proposition~\ref{III-relative_normalization:prop:homogeneous}
identifies each such branch with its Haar quotient and supplies both radical
ergodicities; no ergodicity of the whole finite cover is needed. The signed
branch sets are measurable and invariant under the connected parabolic.
Choose the first branch on any overlap, and combine the resulting coordinates
on these invariant sets. Zero-mass branches impose no almost-everywhere
condition. This also covers coincident projective branches.

Restrict the suspension germ cocycle to the actual equivariant curves and
apply Theorem~\ref{III-relative_normalization:thm:curve-normalization} on each branch.
Keep the original normal vector and exponent throughout this restriction.
The oriented normal functional has the recognized tangent plane as its kernel,
so it is nonzero on the given transverse tangent. Transporting it through the
same ambient charts preserves its positive cone. Changing the normal vector's
sign to agree with this orientation preserves all return product norms. Recurrent returns to sets on which the measurable scalar
coordinate derivative and its reciprocal are bounded identify the residual
flow derivative with $e^{qt}$, hence $B'(0)=q$; no integrability of that
coordinate change is required. Since $q\ne0$, apply Koenigs linearization
to a contracting time map
\cite[Theorem~III and \S\S10--12, pp.~S15--S16]{Koenigs1884} and \cite[introduction, p.~385]{Koenigs1885}.
The universal source statement concerns a holomorphic germ fixing zero,
with multiplier $a$ satisfying $0<|a|<1$: it gives a holomorphic coordinate
$K$ with $K(0)=0$, $K'(0)=1$, and $K\circ f=aK$.
For the real germ, restrict this coordinate to the real axis and take real
parts. Its derivative is still one, and the analytic inverse theorem supplies
the real germ inverse. Every other time map commutes with this contraction,
so its conjugate is linear by coefficient comparison. A single deterministic
coordinate change therefore makes the residual flow $e^{qt}z$.

These are parameters of the original curve images. If the constructed axis
chart is $H_x(ze)=C_x(r_x^{-1}z)$ and the final scalar normalizer is $W_x$,
then the resulting parametrization is
\[
 \overline C_x(z)=C_x\bigl(r_x^{-1}W_x^{-1}(z)\bigr),\qquad r_x\ne0.
\]
The scalar inverse parameter is $u\mapsto W_x(r_xu)$, so the image germ is
unchanged. The original restricted-cocycle identity gives exactly the stated
$H$ and central covariance. Returning through the same ambient charts
preserves these identities, and shrinking their measurable positive domains
gives actual analytic parametrizations.
\end{proof}

\subsection{Ambient Levi normalization in the negative branch}

The scalar reduction has a further consequence when $q<0$.
It supplies the missing input for a contracting Cartan argument,
without assuming a Levi reduction of ambient jets.

\begin{lemma}[Contracting linearization preserving a hypersurface]
\label{III-relative_normalization:lem:contracting}
Let $T$ be an invertible measurable probability-preserving transformation,
and let $f_x$ be measurable analytic
diffeomorphism germs of $(\R^m\oplus\R,0)$ preserving $z=0$.
Suppose $D_0f_x=D$ is a fixed positive diagonal contraction.
Assume the analytic bounds and inverse radii of $f_x$ and its
inverse are tempered along $T$. Suppose that no nonlinear
monomial preserving $z=0$ is resonant for conjugation by $D$.
Then there are measurable analytic germs $h_x$, tangent to the
identity and preserving $z=0$, such that
$h_{Tx}^{-1}f_xh_x=D$.
\end{lemma}
\begin{proof}
Pass to a measurable bidirectionally invariant conull set on which all
original hypotheses hold. On its complement extend the family by $D$;
this retains the original base and probability. Choose $N\ge1$ with
$\lambda_+^{N+1}<\lambda_-$, where $\lambda_+$ and $\lambda_-$ are
the largest and smallest eigenvalues of $D$.
At each degree $2\le d\le N$, the coefficient equation is
\[
                  b(Tx)-\tau b(x)=c(x),\qquad \tau>0,\quad\tau\ne1.
\]
For $\tau<1$ use
$b(x)=\sum_{j\ge1}\tau^{j-1}c(T^{-j}x)$; for $\tau>1$ use
$b(x)=-\sum_{j\ge0}\tau^{-j-1}c(T^jx)$.
These are the convergent homological sums of
\cite[Proposition~3.3, pp.~2151--2153]{BDM2008}.
Their geometric bounds give measurable tempered coefficients.
The forbidden normal-output coefficients vanish because $f_x$
preserves $z=0$; set their corrections to zero, even when their weights
are resonant. Permutation symmetry is retained. Each polynomial
correction is tangent to the identity and preserves the plane, as does
its actual analytic inverse by the inverse theorem on the plane.
Coefficient bounds for the correction and its inverse give bilateral
subexponential intrinsic height. Finite composition therefore gives an
actual measurable gauge $H_x$ with the same properties and
$g_x=H_{Tx}^{-1}f_xH_x=D+O(\|w\|^{N+1})$.

For analytic completion, in the contracting-limit setting of
Theorem~\ref{II-loc:thm:contracting-axis}, apply the original holomorphic
\cite[Theorem~1.1, pp.~2138--2139]{BDM2008} with linear target $D$,
$q=N$, $m=\lambda_-$ and $M=\lambda_+$. The real orthonormal bounds
complexify with the same constants. The finite real jet gives equality
of the complex derivatives through order $N$; the bilateral intrinsic
height bounds give the required tame holomorphic source radii and image
bounds. For $M^{N+1}/m<\vartheta<1$, the theorem supplies a
near-identity holomorphic map $K_x$ on a slow ball, satisfying
$K_{Tx}g_x=DK_x$, and the actual normalized-iterate limit
\[
 K_x(w)=\lim_{d\to\infty}D^{-d}g_x^{(d)}(w),\qquad
 \|D^dK_x(w)-g_x^{(d)}(w)\|
 \le\phi(x)(m\vartheta)^d\|w\|^{N+1},
\]
with a fast bound $\phi$. These finite iterates commute with complex
conjugation and are measurable, so the same limit defines a measurable
real analytic germ. Every real iterate preserves the plane; its closedness
passes this property to the limit. The analytic inverse theorem gives
the inverse germ, also preserving the plane, and coefficient inversion
preserves measurability. Thus $h_x=H_xK_x^{-1}$ has all the required
properties for the original family almost everywhere. The argument also
covers zero or one tangential dimension.
\end{proof}

\begin{theorem}[Measured ambient Levi charts for every negative exponent]
\label{III-relative_normalization:thm:negative-Levi}
For a recognized boundary with $q<0$, its homogeneous suspension
admits measurable actual analytic ambient charts in which
\[
 F(l,x)(u,z)=(lu,z),\qquad
 F(a_t,x)(u,z)=(e^{-3t/2}u,e^{qt}z).
\]
The charts send $z=0$ to the boundary and can be calibrated to
its stable Haar parametrization. No exceptional negative value
of $q$ is excluded. This is a statement at measure-typical
suspension centers, not yet an arithmetic Levi chart at a
specified rational boundary point.
\end{theorem}
\begin{proof}
Normalize the normal quotient derivative by
Corollary~\ref{III-relative_normalization:cor:normal-derivative}. Write the two missing-root
weights as $\beta_2,\beta_3$. The normal weight is
\[
 q\varepsilon_1=c(\beta_2+\beta_3),\qquad c=-q/3>0.
\]
The boundary Haar frame supplies the exact tangential Cartan
derivative. Choose a Cartan element in the cone
$\beta_2,\beta_3<0$, with three distinct fiber exponents.
Its transverse Oseledec line, rescaled by the normal-quotient
gauge, gives a full frame in which the Cartan derivative is the
constant diagonal representation with these three weights.
Commutation preserves this Oseledec line for every fixed
Cartan element. Boundary-adapted analytic coordinate germs
realize the frame. The splitting and subquotient spectrum used here are
commuting-Cartan and exact-sequence consequences of the invertible
multiplicative ergodic theorem
\cite[Theorem~1.6, Corollary~1.7 and Remark~1.8, pp.~33--34]{Ruelle1979},
\cite[Theorem~2.2.6, Remark~2.2.9(iii) and Variant~2.2.10, pp.~6--7; Lemma~2.3.3 and Remark~2.3.6, pp.~8--9]{FilipMETv1}.
Concretely, we use the original integer-ray joint splitting of
\cite[nonergodic discussion and equation~(7), p.~463]{BRHWSmoothErgodic2023}, with the
fixed-real-direction extension described above. The subquotient spectrum
is derived by comparing the actual invariant growth filtrations and their
dimensions on this same law. These arguments apply to the original
integrably bounded derivative cocycle; no integrability of the resulting
measurable frame is assumed.

There are no nonlinear \emph{joint Cartan} resonances in the
space of maps preserving the boundary. Indeed a monomial
$u_2^{p_2}u_3^{p_3}z^s$ has weight
\[
       (p_2+cs)\beta_2+(p_3+cs)\beta_3.
\]
To equal the $u_2$ output weight one must have
$p_3=s=0,p_2=1$, and similarly for $u_3$. For a normal output,
boundary preservation requires $s\ge1$; equality forces
$p_2=p_3=c(1-s)$, hence $s=1,p_2=p_3=0$.
All these are linear monomials. Choose the contracting Cartan
element $a$ outside the countably many hyperplanes for the
remaining monomial equalities, also retaining distinct
derivative eigenvalues.

The analytic bounds in these coordinates are tempered along
$a$. The original compact-manifold suspension has integrable
logarithmic analytic heights for each fixed return element.
The derivative equation with a constant diagonal target bounds
the increments of the logarithmic frame height by these
integrable return heights and a constant. The integrable-increment argument in Lemma~\ref{II-loc:lem:analytic-tempering} therefore makes the frame and
its inverse tempered, without assuming their integrability.
A finite atlas adapted to the compact analytic boundary then
gives the asserted tempered analytic coordinate bounds.
Lemma~\ref{III-relative_normalization:lem:contracting} linearizes $F(a,x)$ while preserving
the boundary.

For any other Cartan element, commutation with $a$ makes each
nonlinear coefficient of its normalized germ equivariant for
conjugation by $D(a)$. Each admissible coefficient weight differs
from one by the choice of $a$. Probability recurrence kills
these coefficients. Thus the full Cartan action, and in particular
$a_t$, is linear with the stated weights.

Put $D_0=\operatorname{diag}(e^{-3/2},e^{-3/2},e^q)$.
For each fixed $l\in L$, its now normalized germ satisfies
\[
                  F(l,a_1x)=D_0F(l,x)D_0^{-1}.
\]
A coefficient with conjugation weight different from one vanishes
by recurrence. A coefficient of weight one is $a_1$-invariant
and therefore constant, since $a_1$ acts ergodically on
$G/\Gamma$. This uses the universal unitary representation statement
\cite[Theorem~1.1 and Definition~2.6, pp.~66 and~69]{Ciobotaru}
on the zero-mean Koopman space for the same normalized Haar probability;
its continuity and absence of $G$-fixed vectors are checked before
applying coefficient decay. Consequently $F(l,x)=\Theta(l)$ is deterministic.

The centralizer of the strict contraction $D_0$ in invertible
formal germs is a finite-dimensional group of resonant polynomial
maps: a monomial of degree $j$ has conjugation modulus at least
$\lambda_-/\lambda_+^j$, so only finitely many degrees can occur.
Inverses have the same property. Thus $\Theta$ is a measurable,
hence continuous, homomorphism of $L$ into a finite-dimensional
Lie group. Its derivative preserves the boundary tangent plane
and has standard action there and trivial normal quotient.
Complete reducibility splits this derivative by a constant linear change
commuting with $D_0$ and preserving the boundary. Here the group statement
is derived from the universal real Lie-algebra invariant-complement theorem
\cite[Lemma~B.1 and proof, pp.~659--660]{KnappBeyondDigital} as follows.
For any finite-dimensional continuous real representation $R$ of
$L=\SL(2,\R)$, the matrix homomorphism theorem
\cite[Definition~1.15, p.~17; Theorem~2.21, pp.~45--47]{Hall2003}
gives a real Lie homomorphism $dR$ with
$R(\exp(tX))=\exp(t\,dR(X))$. The actual exponential tangent of $L$
is the trace-zero algebra, by the determinant-exponential identity
\cite[Proposition~0.11(e), pp.~7--8]{KnappBeyondDigital}.
It is semisimple: $\mathfrak{sl}(2,\C)$ is simple and every real form
of a complex simple algebra is simple
\cite[\S II.1, pp.~124--126; proof of Proposition~7.1, p.~434]{KnappBeyondDigital}.
An $L$-invariant subspace is $dR$-invariant by differentiation.
Choose its invariant Lie-algebra complement and form the projection
along that complement. The projection commutes with each $dR(X)$,
hence with $R(\exp(tX))$. The exponential image contains an identity
neighborhood and $L$ is connected
\cite[Theorem~0.15, pp.~12--13; Proposition~1.145, p.~118]{KnappBeyondDigital};
therefore these exponentials generate $L$ and the projection is
$L$-equivariant. This argument applies to every finite jet space as well
as to the original derivative space, and requires no rank-two hypothesis.

Remove its remaining resonant terms degree by degree. At the
first untreated degree the differentiated homomorphism identity
is a Lie-algebra one-cocycle for $\mathfrak{sl}_2$ in the
finite-dimensional module of resonant polynomial vector fields
preserving $(z)$. Whitehead's lemma gives a correcting vector
field in that same module. Equivalently, apply the same real complete
reducibility theorem to the affine extension of this cocycle module;
the invariant complement gives the correcting coboundary. Its time-one germ commutes with
$D_0$, preserves the boundary, and removes the coefficient.
There are only finitely many resonant degrees. The corrections
are polynomial: a nonlinear resonant field raises ordinary
degree while preserving positive contraction weight, so its
iterated action on each coordinate terminates. Equality on
the connected group $L$ follows from equality of its infinitesimal
actions. This produces actual analytic Levi charts, still
linearizing $a_t$.

On the boundary, the coordinate changes used in the last step
have identity tangent action and no nonlinear terms: a scalar
contraction on the two-dimensional boundary has no nonlinear
resonances. The boundary restriction of the chart is its stable
Haar parametrization. Equivalently, compare the two at their
first different degree; central covariance has weight
$e^{(j-1)3/2}\ne1$, so recurrence gives equality.
\end{proof}

\begin{corollary}[The negative transverse axis has no exceptional weights]
\label{III-relative_normalization:cor:negative-axis}
For every $q<0$, there are measurable analytic transverse
parameters with $H$ acting identically and
$a_t$ acting by $z\mapsto e^{qt}z$.
\end{corollary}
\begin{proof}
Use Theorem~\ref{III-relative_normalization:thm:negative-Levi}. For a fixed row $r\ne0$,
Levi covariance makes each finite jet of $F(u_r,x)$ equivariant
under the unipotent row stabilizer $L_r$. Recurrence puts each
coefficient in its fixed space; ergodicity of $L_r$ on
$G/\Gamma$ then makes it constant. Thus the row germs are a
deterministic commuting homomorphism, continuous at each finite
jet order. Its formal infinitesimal fields are linear in $r$
and Levi-equivariant. These are compatible finite-jet logarithms:
the finite operator calculus of \cite[14.28, p.~287]{Milne},
the polynomial jet locus, and truncation give the construction
explained in Section~\ref{III-series:covariants}. No convergence of
the formal logarithm is required.

Their full two-dimensional covariant list is
\[
 Z_r=(ru)\{a(z)E_u+b(z)\partial_z\}
                   +c(z)Jr^{\mathsf T}\cdot\partial_u .
\]
Boundary calibration gives $a(0)=-1$, $b(0)=c(0)=0$.
The row commutator, calculated in Section~\ref{III-series:covariants}, gives
$bc=0$ and $3ac-bc'=0$. Since $a$ is a unit, these identities
force $c=0$ in the formal power-series ring. Every remaining
field vanishes on $u=0$. Hence every row germ fixes the axis
pointwise as a formal germ, and therefore as an actual analytic
germ. Levi and central equivariance are already part of
Theorem~\ref{III-relative_normalization:thm:negative-Levi}. The exceptional values in the
earlier first-discrepancy calculation are thus removed by
the full ambient normalization.
\end{proof}

\section{Two-dimensional row and column covariants}\label{III-series:covariants}

\subsection{The extra two-dimensional covariant vanishes}

The following lemma addresses a further use of $n\ge4$ in the
arithmetic parabolic calculation. It assumes the Levi normalization;
it does not prove it.

\begin{lemma}[Commuting row fields after Levi normalization]
\label{III-boundary_reduction:lem:row}
At an analytic projective boundary point for $n=3$, suppose
the block Levi lattice is analytically conjugate to
$L_A(u,t)=(Au,t)$, with $u\in\R^2$, and with the boundary
identification given by $t=0$. Then the row lattice has analytic
logarithms
\[
 Z_r=(ru)\bigl(a(t)E_u+b(t)\partial_t\bigr),
 \qquad a(0)=-1,\quad b(0)=0.
\]
\end{lemma}
\begin{proof}
The positive normal derivative character of the row lattice is
Levi-invariant and hence trivial. Remove finite signs by doubling
the lattice. Its derivative on the tangent plane is the identity.
The possible normal-to-tangent block is a linear map from the row
module to the standard module; in dimension two this map may be
nonzero. Nevertheless all row derivatives are commuting unipotent
matrices preserving a common flag. Filter centered polynomial jets by
vanishing order. On each graded piece, inverse pullback acts by
precomposition with the inverse first derivative. This is a product
of commuting unipotent input-slot operators, hence is unipotent.
The finite filtration therefore makes the full pullback unipotent.
At every finite jet order their
commuting unipotent jets have unique commuting logarithms, additive
in the row parameter. For the finite operator calculus we use
\cite[14.28, p.~287]{Milne}: on a finite-dimensional vector space
over a field of characteristic zero, the finite exponential and
logarithm polynomials give inverse bijections between nilpotent and
unipotent endomorphisms. Explicitly, if $(U-I)^q=0$, write
$L=\sum_{j=1}^{q-1}(-1)^{j+1}(U-I)^j/j$; then $\exp L=U$, and
$\log(\exp S)=S$ for nilpotent $S$. Thus the real operator powers
are $\exp(tL)$, with derivative $L$ at zero. The compatible
finite-jet logarithms follow from this operator calculus as follows.
To see that the real operator powers
are jets, use the faithful inverse-pullback representation on centered
polynomial jets. Its image, with the inverse operator included, is
cut out by polynomial equations: evaluation at the identity polynomial
recovers the inverse jet, and substitution recovers the full operator.
The pair $(\exp(tL),\exp(-tL))$ is polynomial in $t$ and satisfies
these equations at every nonnegative integer, hence at every real $t$.
The inverse operator recovers the forward coefficients, so those
coefficients are polynomial and in particular measurable in $t$.
Truncation is a surjective linear intertwiner of pullback operators;
thus it preserves unipotence. The canonical powers commute with
truncation because their polynomial coefficients agree at integer
times. The same polynomial identity argument gives uniqueness of
positive integer roots. Commuting nilpotent logarithms have commuting
exponentials, and conjugation respects the displayed polynomials.
All constructions use real coefficients. When doubling the lattice, the logarithm of the
doubled basis is evaluated at half the physical row parameter; this
retains the normalization $a(0)=-1$ in the original row coordinates.
Levi covariance holds for these logarithms
and extends by Zariski density to $\SL(2,\R)$.

Fix $J=\left(\begin{smallmatrix}0&-1\\1&0\end{smallmatrix}\right)$.
The scalar polynomial covariants of $(u,r)$, linear in $r$, are
multiples of $ru$. The vector covariants are linear combinations
of $(ru)u$ and $Jr^{\mathsf T}$. Thus the formal row fields have
the form
\begin{equation}\label{III-boundary_reduction:eq:extra}
 Z_r=(ru)(aE_u+b\partial_t)+c(t)(Jr^{\mathsf T})\cdot\partial_u,
 \qquad a(0)=-1,\quad b(0)=c(0)=0.
\end{equation}
Writing $u=(x,y)$ gives
\[
 Z_1=axE_u+bx\partial_t+c\partial_y,
 \qquad Z_2=ayE_u+by\partial_t-c\partial_x.
\]
A direct bracket calculation yields
\begin{equation}\label{III-boundary_reduction:eq:commuting}
 [Z_1,Z_2]=(3ac-bc')E_u+2bc\partial_t.
\end{equation}
Commutativity gives $bc=0$ and $3ac-bc'=0$ as formal series.
If $c\ne0$, the first identity forces $b=0$, and the second
forces $a=0$, contradicting $a(0)=-1$. Thus $c=0$.

Each application of the remaining field increases $u$-degree by
one. The $u$-linear normal coefficient of its analytic time-one
map is $(ru)b(t)$; the $u$-quadratic tangential coefficient is
$(ru)u\,a(t)$. These extract convergent $a,b$ from the original
analytic row maps. The analytic flows have the prescribed jets,
and therefore are the prescribed germs.
\end{proof}

There is a corresponding observation for opposite columns. An
$\SL(2)$-equivariant family linear in $s\in\R^2$ has the form
\[
 Y_s=k(t)s\cdot\partial_u+
       \det(s,u)\bigl(\ell(t)E_u+m(t)\partial_t\bigr).
\]
If these fields commute and $k(0)=1$, the same bracket calculation
(with row parameters rotated by $J$) gives $mk=0$ and
$3\ell k-mk'=0$. Hence $m=\ell=0$. Thus
$Y_s=k(t)s\cdot\partial_u$ also holds for $n=3$.

\section{Contracting collars for \texorpdfstring{$n=3$}{n = 3}}\label{III-series:negative-three}

\subsection{The measured charts and their global stable extensions}

Use the notation in Section~\ref{III-series:relative-three}. We have a recognized projective
covering hypersurface, its homogeneous parabolic probability
space $X=G/\Gamma$, and measurable analytic charts $H_x$ satisfying
\[
 H_x(u,0)=\phi_x(u),\qquad
 H_{lx}^{-1}lH_x(u,z)=(lu,z),\qquad
 H_{a_tx}^{-1}a_tH_x(u,z)=(\lambda^t u,b^t z),
\]
where $\lambda=e^{-3/2}$, $b=e^q$, and $q<0$.
The notation $gH_x$ includes the actual suspension transport
between the corresponding fibers. The missing column group is
$V\simeq\R^2$. Its Haar-graph action is
$\widehat v x=v\phi_x(v)$. Under the homogeneous identification
this is left translation by $v$ on $G/\Gamma$.

\begin{lemma}[Global stable charts and their overlaps]
\label{III-negative_collar:lem:global}
The charts extend to analytic diffeomorphisms
\[
 H_x:\R^3\longrightarrow W^{s,F}_{a_1}(x).
\]
For every fixed $v$, almost everywhere,
\[
 v^{-1}W^{s,F}_{a_1}(\widehat v x)=W^{s,F}_{a_1}(x).
\]
Consequently
\[
 P_v(x)=H_x^{-1}v^{-1}H_{\widehat v x}
\]
is an analytic diffeomorphism of $\R^3$, with
$P_v(x)(u,0)=(u+v,0)$.
\end{lemma}
\begin{proof}
Lemma~\ref{II-loc:lem:analytic-tempering} applies to
$H_x$, since its conjugated central map is the fixed contraction
$D=\operatorname{diag}(\lambda I,b)$. Thus
$H_x(w)=a_1^{-n}H_{a_1^nx}(D^nw)$ consistently defines the
extension for all $w$, and the usual stable-manifold argument
makes it injective and onto the full fiber stable manifold.
All three fiber exponents are negative, so this stable manifold
is open in its fiber and the injective parametrization is a
diffeomorphism onto it.

The base displacement $a_1^nva_1^{-n}$ contracts. Its fiber maps
and inverse maps have subexponential Lipschitz bounds along the
typical orbit, by the compact-parameter suspension estimates
used in Section~\ref{II-sec:neutral}. They therefore preserve exponential
convergence of pairs of fiber points. Since the two centers lie
on the same stable Haar graph, this proves the equality of
stable fibers. Haar-graph coherence gives the stated translation
on their boundary planes.
\end{proof}

\subsection{The two-dimensional column stabilizer suffices}

\begin{proposition}[Deterministic analytic column transitions]
\label{III-negative_collar:prop:columns}
The maps $P_v(x)$ are independent of $x$ almost everywhere.
They form a jointly analytic local action in $(v,u,z)$ and have
the form
\[
                         P_v(u,z)=(u+k(z)v,z),\qquad k(0)=1.
\]
In fact $k(z)=1$.
\end{proposition}
\begin{proof}
For a fixed nonzero column $v$, let $L_v=\{l:lv=v\}$.
Levi equivariance of the charts gives
\[
 P_v(lx)=\operatorname{diag}(l,1)\,P_v(x)\,
                                      \operatorname{diag}(l^{-1},1).
\]
Every finite Taylor coefficient is in a finite-dimensional
unipotent representation of $L_v$. Recurrence makes it fixed
under that representation. It is then an $L_v$-invariant
function of $x$. The noncompact group $L_v$ acts ergodically
on $G/\Gamma$ by Howe--Moore
\cite[Theorem~1.1 and Definition~2.6, pp.~66 and 69]{Ciobotaru},
so each coefficient is constant.
Analytic uniqueness and Lemma~\ref{III-negative_collar:lem:global} make the entire
map deterministic. These statements hold for every fixed
$v$, with Fubini supplying the jointly measurable versions.

Here it is essential to retain the determinant covariants.
For $v=e_1$, invariance under
$(u_1,u_2,z)\mapsto(u_1+tu_2,u_2,z)$ gives exactly
\[
 P_{e_1}(u,z)=
       \bigl(A(u_2,z)u+B(u_2,z)e_1,\ C(u_2,z)\bigr).
\]
For example the second component is independent of $u_1$
and divisible by $u_2$; its quotient is $A$, and the first
component minus $Au_1$ is $B$. These are analytic functions.
Full Levi covariance and transitivity on nonzero columns give
\begin{equation}\label{III-negative_collar:eq:column-form}
 P_v(u,z)=\bigl(A(\det(v,u),z)u+B(\det(v,u),z)v,\
                                     C(\det(v,u),z)\bigr).
\end{equation}
The same three analytic functions work for every nonzero $v$.

The chart identities and the group law of $\widehat V$ give
$P_{v+s}=P_vP_s$. Formula~\eqref{III-negative_collar:eq:column-form} has an analytic
limit $P_*$ as $s\to0$. Letting $s\to0$ in
$P_{v+s}=P_vP_s$ for a fixed nonzero $v$ gives
$P_v=P_vP_*$, hence $P_*=\id$. Thus the formula extends
analytically to $v=0$ with the correct identity map.
This supplies a genuine analytic parameter family; continuity
in the parameter has not been assumed from coefficient
measurability alone.

Its infinitesimal fields are linear in the column parameter,
Levi-equivariant, and commute. The full covariant list is
\[
 Y_s=k(z)s\cdot\partial_u+
          \det(s,u)\{\ell(z)E_u+m(z)\partial_z\}.
\]
The calculation in Section~\ref{III-series:shear-three}, Section~\ref{III-series:covariants} gives $mk=0$ and
$3\ell k-mk'=0$. Since $k(0)=1$, we get $\ell=m=0$.
Integration gives the asserted form of $P_v$.

Finally central covariance is
$P_{\lambda v}=DP_vD^{-1}$. It implies $k(b^{-1}z)=k(z)$.
Since $0<b<1$ and $k$ is analytic, $k$ is constant.
Boundary calibration gives $k=1$.
\end{proof}

\subsection{The local parabolic model and its comparison with columns}

\begin{proposition}[A deterministic local group model]
\label{III-negative_collar:prop:model}
In the measured charts, all parabolic germs are deterministic
and analytic in their group parameters. Together with the column
transitions they satisfy the local relations of $G$.
The resulting local model is a projective flow model with
normal infinitesimal weight $q$.
\end{proposition}
\begin{proof}
Levi and central germs already have these properties. For a
fixed nonzero row, the row-stabilizer recurrence and ergodicity
argument in Section~\ref{III-series:relative-three} makes its germ deterministic.
At every finite jet order the row homomorphism is continuous
and has a formal logarithm linear in its row parameter.
Its general Levi-equivariant field is
\[
 Z_r=(ru)\{a(z)E_u+b_0(z)\partial_z\}
                 +c(z)Jr^{\mathsf T}\cdot\partial_u.
\]
The boundary action gives $a(0)=-1$, $b_0(0)=c(0)=0$.
Commutativity forces $c=0$ as in Section~\ref{III-series:covariants}.
Central covariance, with $q\ne0$, then forces
$a(z)=-1$ and $b_0(z)=\beta z$ for a constant $\beta$.
Thus the formal fields are polynomial and their local analytic
flows agree to every order with the actual row germs.
Analytic uniqueness gives the actual row formula.

For a word in the parabolic and column factors, its coordinate
transition is the composition of the corresponding deterministic
maps, with the inverse convention for parabolic transport.
If the base matrix word is the identity, its action on the
homogeneous center space $X$ is the identity. Suspension
transport also telescopes to the identity, so the complete chart
transition is the identity. For small parameters this may first
be checked on a positive-measure set with bounded chart domains;
the deterministic analytic maps then satisfy the identity as
germs. This gives all local matrix relations.

In the fundamental-field convention the column fields are
$s\cdot\partial_u$ and the row fields are
$-(ru)E_u+\beta(ru)z\partial_z$. Their mixed bracket gives
the diagonal matrix field. For $r=e_1^{\mathsf T}$ and $s=e_1$,
\[
 [Z_r,\partial_{u_1}]
       =2u_1\partial_{u_1}+u_2\partial_{u_2}
                                      -\beta z\partial_z.
\]
Fundamental fields for a left action reverse Lie brackets.
The corresponding matrix bracket is
$\operatorname{diag}(1,-1,0)$; its normal weight, determined
by the already normalized Cartan, is $q$.
Thus $\beta=q$, with simultaneous inversion of root parameters
giving the same conclusion in the inverse-transition convention.

The local action is consequently
\begin{equation}\label{III-negative_collar:eq:flow}
 \begin{pmatrix}\alpha&r\\s&A\end{pmatrix}(u,z)
 =\left(\frac{s+Au}{\alpha+ru},\
                         (\alpha+ru)^q z\right)
\end{equation}
where $\alpha+ru>0$. This formula has the specified row,
column, Levi, and central fields, so uniqueness of their
local flows identifies the local model. Orientation changes
and inverse transpose give the other covering and dual charts.
\end{proof}

\subsection{Analytic descent and gluing of the transverse parameters}

Passing from typical centers to a specified boundary point
requires care: an almost-everywhere identity cannot be evaluated
at a measure-zero stabilizer. We instead regularize the lifted
chart field on open subsets of $G$, before taking any section
or stabilizer.

\begin{theorem}[Full negative-branch collar]\label{III-negative_collar:thm:collar}
Every recognized projective covering hypersurface with $q<0$
has a full-angular analytic collar neighborhood in dimension
three. After the finite angular and normal-orientation choices,
its transverse flow is generated by $qz\partial_z$ and the
lattice action has the local form \eqref{III-negative_collar:eq:flow}.
No negative transverse weight is excluded.
\end{theorem}
\begin{proof}
Lift the charts to $G\times N$ and write them $H_g$. They satisfy
the right-lattice identity
\[
 H_{g\gamma}=\alpha(\gamma)^{-1}H_g
\]
for each fixed $\gamma$, almost everywhere in $g$. In a local
factorization of $h$ into column and parabolic factors, the
deterministic transitions give
\begin{equation}\label{III-negative_collar:eq:lifted}
                         H_{hg}=H_g\circ P_h
\end{equation}
for almost every $(h,g)$, on the corresponding moving-origin
domains. The inverse convention for the parabolic factors makes
$P_h$ an analytic inverse-transition model. In particular $P_q$
fixes the chart origin for $q\in Q$, and on the transverse axis
it acts by the inverse of the normalized central character.

Fubini permits a fixed good $g_0$ for which \eqref{III-negative_collar:eq:lifted}
holds for almost every $h$ in a small identity neighborhood.
Shrink this neighborhood and the chart ball. Its right side is
then jointly analytic in $h$ and the chart variables and supplies
an analytic representative $\widetilde H_g$ on a nonempty open
set $O\subset G$.

We now propagate this \emph{parameterized chart field}, not just
its unparametrized axes. Projective minimality gives
\[
                              QO\Gamma=G.
\]
Use the connected parabolic and the spherical boundary when
needed; the finite covering choices made in Section~\ref{III-series:relative-three}
ensure the same statement on the chosen angular cover.
For each fixed $q_0\in Q$ and $\gamma\in\Gamma$, define on
$q_0O\gamma$ the analytic germ field
\[
 \widetilde H_{q_0g\gamma}
             =\alpha(\gamma)^{-1}\widetilde H_g\circ P_{q_0}.
\]
This agrees almost everywhere with the original measured field,
by its fixed-element parabolic and lattice identities. Choose
a countable subcover of these open sets. On every overlap the
two analytic fields have identical Taylor coefficients for
almost every $g$, and hence for every $g$ by analytic uniqueness.
They therefore glue to a globally defined analytic germ field
$\widetilde H_g$ on $G$, with positive local chart radii.

The identities
\[
 \widetilde H_{qg}=\widetilde H_gP_q,\qquad
 \widetilde H_{g\gamma}=\alpha(\gamma)^{-1}\widetilde H_g
\]
now hold everywhere on their germ domains. For the first,
use Fubini and analytic uniqueness in the joint variables
$(q,g)$; the parabolic model is analytic in $q$.
For the second, use the original almost-everywhere identity
for each fixed $\gamma$ and analytic uniqueness in $g$.
The same argument extends the local column identities.
Only at this stage do we restrict to analytic local sections
of $Q\backslash G$.

We realize these charts in the actual covering neighborhood. For an
embedded spherical center, use the original ambient manifold with the
identity projection. For an embedded projective center
$i:\mathbb{RP}^2\hookrightarrow M$, Hatcher's ENR theorem
\cite[Corollary~A.9, p.~527; proof using Theorem~A.7, pp.~525--527]{Hatcher},
together with coordinatewise Tietze extension, gives an open neighborhood
$W$ and a continuous retraction $R:W\to\mathbb{RP}^2$ with $R\circ i=\mathrm{id}$.
Form the pullback of $S^2\to\mathbb{RP}^2$ along $R$.
Its projection $\pi:Y\to W\subset M$ is a local homeomorphism; pulling back
the ambient analytic charts makes it a local analytic diffeomorphism
on a three-dimensional analytic manifold. The prescribed sphere zero
section embeds in $Y$, and the genuine antipodal deck map is analytic
because its composition with $\pi$ equals $\pi$. This auxiliary topological
construction does not change the original chart field or its selected seed.
Lift the local collar charts through $\pi$ at this prescribed zero section;
local-cover uniqueness makes their joint overlaps agree.

If $s_i(p)$ is such a section, define
\[
                        j_i(p,z)=\widetilde H_{s_i(p)}(0,z).
\]
Its derivative is invertible along $z=0$: the center map is
the recognized projective parametrization and the $z$ direction
is transverse. Thus it is an analytic collar chart after
shrinking. On an overlap, $s_j(p)=q_{ji}(p)s_i(p)$ with analytic
$q_{ji}(p)\in Q$. The everywhere-valid parabolic identity shows
that $j_j$ differs from $j_i$ precisely by the scalar transverse
character of $q_{ji}(p)$, with the inverse convention above.
This is the associated line-bundle transition rule for
the flow $qz\partial_z$.

The everywhere-valid right-lattice identity proves equivariance
of these parameterized charts for every lattice element.
Combining it with the local column and parabolic formulas gives
\eqref{III-negative_collar:eq:flow} in overlapping projective charts. Thus neither
a stabilizer identity nor a variable section has been substituted
into an almost-everywhere formula.

Finally take finitely many charts over the compact boundary.
Their compatible joint germs and the analytic inverse function theorem
give an embedded full-angular collar $C:S^2\times(-r,r)\to Y$ after one
uniform shrink. The projection formulas retain the original lattice
action. For each lattice element, conjugating its weighted model by this
same $C$ produces an analytic local lift on a domain containing the whole
zero section; compactness gives a uniform angular strip for that element.
Local-cover uniqueness proves the identity, multiplication and inverse
laws as germs along the zero section.

In the projective case the physical factorization makes the disconnected
line-stabilizer flip preserve the same canonical carrier law. Apply the
already proved Levi and central normalization to its original gauge
cocycle: the flip is represented by one deterministic linear map $A$
commuting with the normal half-scaling. Its order-two cocycle relation
gives $A^2=1$, so its normal-axis action is $z\mapsto\epsilon z$ with
$\epsilon\in\{1,-1\}$. The original chart equations give the corresponding
projected deck formula, and uniqueness of the actual cover lift gives
\[
 \delta C(\theta,z)=C(-\theta,\epsilon z).
\]
A further uniform shrink makes the collar image deck invariant.
No antipodal identity is imposed in the embedded spherical case.

On this very collar define
\[
 \Phi_t(C(\theta,z))=C(\theta,e^{qt}z)
\]
whenever both coordinate points lie in its domain. This is a jointly
analytic local flow with generator $C_*(qz\partial_z)$; the local lattice
lifts commute with it as germs. The genuine deck preserves its legal
domain and commutes with it there, so the flow descends through the
original projection. The dual representation is the same construction
after the inverse-transpose lattice identification; pulling the acting
subgroup back preserves finite index, and the literal weighted-model
dictionary preserves $q$. Thus every $q<0$, including resonant values,
is covered. The action radius may depend on the lattice element.
\end{proof}

This settles the contracting line/plane branch, including its
transverse collar, using the already established description
and suspension inputs. It does not settle the positive or
neutral branch: their full fiber is not contracted by the
central element, so Lemma~\ref{III-negative_collar:lem:global} and the resulting
global chart-transition comparison do not apply there.

\section{Expanding collars for \texorpdfstring{$n=3$}{n = 3}}\label{III-series:positive-three}

\subsection{The two scalar contractions}

Fix the positive branch $q>0$, with tangential weights
$\beta_2,\beta_3$ and normal weight
$q\varepsilon_1=-q(\beta_2+\beta_3)/3$.
Use the finite covers in Section~\ref{III-series:relative-three}, so the boundary
suspension is $X=G/\Gamma$ and the normal curves have parameters
$i_x(z)$ satisfying
\[
 h\,i_x(z)=i_{hx}(z)\quad(h\in LU),\qquad
 a_t i_x(z)=i_{a_tx}(e^{qt}z).
\]
The central tangential multiplier is $\lambda_t=e^{-3t/2}$.
Write $V_2,V_3$ for the two coordinate column groups and
$\widehat v x=v\phi_x(v)$ for their Haar-graph actions on $X$.

For each $j\in\{2,3\}$, let $k$ denote the other index and
choose a Cartan element $a^{(j)}$ with
\[
 \varepsilon_1(a^{(j)})=-1,\qquad
 \beta_j(a^{(j)})=-q,\qquad
 \beta_k(a^{(j)})=3+q.
\]
These assignments satisfy $\beta_2+\beta_3=-3\varepsilon_1$.
Its fiber stable space is exactly
$E_{\beta_j}\oplus N$, with both exponents equal to $-q$.
Its other fiber exponent is positive.

\begin{lemma}[Calibrated two-dimensional stable charts]
\label{III-positive_collar:lem:surface}
There are measurable analytic stable charts
\[
 K_x^{(j)}:\R^2\longrightarrow W^{s,F}_{a^{(j)}}(x)
\]
in the intrinsic stable-manifold structure, satisfying
\[
 a^{(j)}K_x^{(j)}(s,z)
       =K_{a^{(j)}x}^{(j)}(e^{-q}s,e^{-q}z),\qquad
 K_x^{(j)}(s,0)=\phi_x^{(j)}(s),\qquad
 K_x^{(j)}(0,z)=i_x(z)
\]
as germs for the last identity. They have tempered local
analytic bounds along $a^{(j)}$.
\end{lemma}
\begin{proof}
Keep the section and the normalized curves fixed above.  The aligned
section on the finite cover has the same integrable word-length control:
the change of section is a bounded finite-coset correction.  Its base law
is the actual Haar marginal.  The two signed developed graphs cover the
whole law almost everywhere; their masses equal the masses of the
corresponding suspension branches.  After removing zero-mass branches,
they still form a conull cover, and at least one branch has positive mass.
We work on either positive-mass branch, with its Haar base law.

Apply the scalar-contracting hypersurface construction of Part~II with
two tangential coordinates, omitting the coordinate indexed by $k$.
That construction uses one unstable coordinate and at least one stable
tangential coordinate, so it applies in this dimension.  The full frame
given by the two Haar derivatives and $i_x'(0)$ has exponents
$-q,-q,3+q$ for $a^{(j)}$.  Its covariance bounds increments of the
logarithms of the frame height and inverse height by the integrable
original return heights.  Hence the frame, including the angle between
its summands, is tempered.  The holomorphic stable-manifold construction
solves the stable coordinates forward and the remaining coordinate
backward.  Its nonstationary scalar Koenigs construction gives the
charts, with positive tempered physical radii, analytic bounds and
inverse bounds.  Contraction extends the same charts to the entire
intrinsic stable manifolds.  The fixed standard and dual coordinate
identifications preserve the original section, return maps and laws;
the ordered parameters are the root and normal parameters $(s,z)$.

Both the root parametrization and $i_x$ are analytic and equivariant for
this scalar contraction, with the prescribed derivatives.  Compare them
as maps into the ambient fiber.  At their first different degree
$r\geq2$, covariance multiplies a stable-output coefficient by
$e^{(r-1)q}$ and an unstable-output coefficient by
$e^{3+(r+1)q}$.  Both factors exceed one, so recurrence forces equality.
This also proves that $i_x$ lies in the stable surface.  Thus the root
axis is the original Haar parametrization and the normal axis is the
same normalized curve germ.  The central Cartan and diagonal Levi act
linearly in these charts by the same uniqueness argument.
\end{proof}

\subsection{Affine comparison along a coordinate line}

We spell out the transport estimate, since the center changes
along the projective chart.

\begin{lemma}[Affine transition under a contracted base displacement]
\label{III-positive_collar:lem:affine}
Suppose two typical centers in one of the preceding stable
families are related by a base displacement $h$ contracted by
$a^{(j)}$, and suspension transport identifies their stable
fibers. Then the transition between their $K^{(j)}$ charts is
affine.
\end{lemma}
\begin{proof}
Write $a=a^{(j)}$ and $r=e^{-q}$. Use the same scalar stable
charts, original section, and tempered local chart and inverse
bounds constructed in Lemma~\ref{III-positive_collar:lem:surface}.
The original section moment and uniform word estimate give one
conull set on which every precompact sequence of base displacements
has sublinear return-word length along the $a$-orbit. In particular,
this applies to $h_n=a^nha^{-n}$.

Let $Q$ be the transition defined by the actual inverse stable
chart and the stated suspension identification. The uniform
stable-orbit estimate, the sublinear correction-word bound, and
analytic inverse recognition show that $Q$ is analytic on its
whole parameter space. Exact scalar covariance and the suspension
cocycle identity identify the iterated physical transition with
\[
 Q_n(u)=r^n Q(r^{-n}u).
\]
Thus $Q_n(0)=r^nQ(0)$ is exponentially small and eventually lies
in the tempered inverse-chart domain. In the original local
ambient charts the transition is the composition of the source
stable chart, the actual lattice map, the change of chart center,
and the inverse target stable chart. Their original coefficient
heights and radii give, by Cauchy estimates, a subexponential
bound for every fixed derivative of $Q_n$ at zero.

For degree $k\geq2$, scalar covariance multiplies the corresponding
coefficient of $Q$ by $r^{(1-k)n}$. The subexponential bound therefore
forces that coefficient to vanish. Analytic uniqueness extends the
resulting affine germ to the entire stable parameter space.
Injectivity of the actual charts makes its linear part invertible.
The fixed linear change to the ordered $K^{(j)}$ parameters and
the exact original base and lattice-label identifications preserve
this affine physical transition.
\end{proof}

\begin{proposition}[The normal parameter is unchanged along either coordinate]
\label{III-positive_collar:prop:translation}
For $v_j(s)\in V_j$,
\[
 (K_x^{(j)})^{-1}v_j(s)^{-1}K_{\widehat v_j(s)x}^{(j)}(t,z)
                                      =(t+s,z).
\]
Consequently
\[
 v_j(s)^{-1}i_{\widehat v_j(s)x}(z)=K_x^{(j)}(s,z)
\]
is analytic jointly in the coordinate parameter $s$ and the
normal parameter $z$.
\end{proposition}
\begin{proof}
The column displacement is contracted by $a^{(j)}$.
Lemma~\ref{III-positive_collar:lem:affine} and root calibration give a transition
\[
                 (t,z)\longmapsto
                 (t+s+A(s,x)z,\ B(s,x)z).
\]
Consider the unipotent stabilizer $L_{e_j}$ of the chosen
column in the Levi. It is contracted by $a^{(j)}$: its
root weight is $\beta_j-\beta_k<0$. It therefore also
transports these stable surfaces, and
Lemma~\ref{III-positive_collar:lem:affine} applies to its centered chart map.
Its derivative on the two-dimensional stable plane is the
identity: it fixes the Haar root vector and acts identically
on the already normalized normal curve. Thus its whole
$K^{(j)}$ chart map is the identity.

This stabilizer commutes with $V_j$. Consequently $A(s,x)$
and $B(s,x)$ are invariant under $L_{e_j}$.
That unbounded subgroup acts ergodically on $G/\Gamma$,
so $A(s,x)=A(s)$ and $B(s,x)=B(s)$ are deterministic.
The column group law gives
\[
 B(s+t)=B(s)B(t),\qquad
 A(s+t)=A(t)+A(s)B(t).
\]
A diagonal Levi element scales the root parameter by any
$n>0$, and fixes the normal parameter. Its linear action in
the stable charts gives
\[
                    B(ds)=B(s),\qquad A(ds)=dA(s).
\]
Taking $d=2$ and using the group law gives
$B(s)=B(2s)=B(s)^2$. Invertibility of the transition gives
$B(s)\ne0$, hence $B(s)=1$.

Finally the central $a_t$ acts on this stable surface by
$(t_0,z)\mapsto(\lambda_t t_0,e^{qt}z)$.
Its covariance gives
$A(\lambda_t s)=\lambda_t e^{-qt}A(s)$.
Comparing this with the diagonal identity at $d=\lambda_t$
and taking $t\ne0$, the condition $q>0$ forces $A(s)=0$.
The second assertion follows by evaluating the transition
at $(t,z)=(0,z)$ and using Lemma~\ref{III-positive_collar:lem:surface}.
\end{proof}

\subsection{A cross argument retaining the normal parameter}

\begin{lemma}[The coordinate cross with a holomorphic parameter]
\label{III-positive_collar:lem:cross-parameter}
Let $f(v_2,v_3,z)$ be measurable near
$I_2\times I_3\times\{0\}$, with an analytic $z$-germ for
almost every $(v_2,v_3)$. Suppose that along almost every
coordinate line it agrees, as a $z$-germ, with a map analytic
jointly in that coordinate and $z$. Then on some nonempty
rectangle in $I_2\times I_3$ it has a jointly analytic
representative for $|z|$ sufficiently small.
\end{lemma}
\begin{proof}
Apply the common-parameter slice lemma of Part II,
Lemma~\ref{II-positive:lem:parameter-cross}, with no auxiliary
parameters, retaining the normal parameter $z$. Here are the
original-data and uniqueness steps in this application. The normal
Taylor coefficients of $f$ are measurable, by finite differences;
changing the exceptional normal fibers on a measurable null set gives
a coefficient-measurable field of convergent germs everywhere. The given
coordinate-line representatives have exactly these germs almost everywhere.
Positive-measure selection gives common complex coordinate domains, a
common normal disk and bounds for both families. Fubini and analytic
uniqueness in the normal variable give their compatibility on almost
every intersection. The normal-family correction and bounded-cross
construction in the proof of that lemma produce a jointly analytic
patch retaining a positive-measure set of original rows and one positive
normal radius. This construction is applied to the original Taylor
field, so it retains the original germs on those rows.

To propagate agreement to almost every row of the smaller rectangle,
compare each normal Taylor coefficient with that of the original
coordinate-line representative. For a connected open real domain, a
nonidentically-zero real-analytic scalar function has a Lebesgue-null
zero set \cite[Proposition~1, p.~529; proof, pp.~529--530]{Mityagin}.
Apply this to each scalar component of the difference: agreement on a
positive-measure set forces equality on the connected interval.
Countably many coefficients recover equality of convergent normal
germs. Restricting the constructed complex patch to the real variables
proves the assertion, with one common smaller normal radius.

Finite-dimensional components treat vector-valued maps. For a compact
manifold target, use the proper injective analytic embedding with
injective differential from
\cite[Theorem~3, p.~470; definition of regularity, p.~462]{Grauert}.
Its local immersion charts follow from the analytic inverse-function
theorem: complement the image of its differential, add these normal
coordinates to the local parametrization, and invert the resulting
map with invertible differential. Projection in this chart gives a
local analytic left inverse at an original good center. Extend this
to a set-theoretic left inverse of the same embedding while retaining
analyticity near that center. Composing the Euclidean representative
with this inverse and shrinking the rectangle and normal disk gives
the manifold-valued representative; its germ agreement follows from
the left-inverse identity.
\end{proof}

\subsection{The positive collar}

\begin{theorem}[Full positive-branch collar]\label{III-positive_collar:thm:collar}
Every recognized projective covering hypersurface with $q>0$
has a full-angular analytic collar neighborhood in dimension
three. After the finite angular and normal-orientation choices,
its transverse flow is $z\mapsto e^{qt}z$.
\end{theorem}
\begin{proof}
Lift the normalized curve field to $G\times N$, writing it
$i_g(z)$. It satisfies the almost-everywhere identities
\[
 i_{hg}(z)=i_g(z)\quad(h\in LU),\qquad
 i_{a_tg}(z)=i_g(e^{-qt}z),\qquad
 i_{g\gamma}(z)=\alpha(\gamma)^{-1}i_g(z).
\]
In a column big cell consider $f(v_2,v_3,z)=i_{vg_0}(z)$.
Fubini gives a good $g_0$ for which
Proposition~\ref{III-positive_collar:prop:translation} supplies the joint
coordinate-and-normal analyticity on almost every line of
each coordinate foliation. Lemma~\ref{III-positive_collar:lem:cross-parameter}
gives a jointly analytic representative on one rectangle,
with a common normal disk.

At $z=0$ this is the recognized projective parametrization.
The normal derivative agrees almost everywhere with the
nonzero transverse derivative of $i_g$. Choose a regular
point where the three derivatives are independent and
shrink. The representative is an actual analytic collar
chart on a nonempty boundary patch.

Parabolic equivariance and its scalar normal character extend
this representative to an open set $O\subset G$ by the local
factorization $QV$. To be precise, use Fubini in that
factorization: its analytic formula agrees with the original
measured curve field for almost every point of $O$.
Projective minimality gives $QO\Gamma=G$.
On each fixed translate $q_0O\gamma$, propagate the
\emph{parameterized} curve field by the displayed parabolic
character and the actual analytic action $\alpha(\gamma)$.
Each propagated representative agrees almost everywhere
with $i_g$ there. A countable cover and analytic uniqueness
on overlaps glue a globally analytic germ field on $G$.

Exactly as in the descent proof in Section~\ref{III-series:negative-three}, the
parabolic identity now holds everywhere by analytic uniqueness
in $(q,g)$, and the right-lattice identity holds everywhere
for each fixed $\gamma$. Restricting to analytic sections of
$Q\backslash G$ is therefore legitimate. Their transition
functions are the scalar normal characters, and the
right-lattice identity gives the actual lattice action.
The resulting collar formula is
\[
 \begin{pmatrix}\alpha&r\\s&A\end{pmatrix}(u,z)
       =\left(\frac{s+Au}{\alpha+ru},\
                                      (\alpha+ru)^q z\right)
\]
on positive-denominator chart domains. Compactness supplies
a finite collar atlas and an embedded smaller neighborhood
of the entire zero section. Uniqueness of the inverse-central
stable curve makes its unparametrized germ intrinsic. Retain the
original full curve field, invariant probability, and prescribed
normal parameter when comparing the two possible angular lifts.
The original antipodal return is a scalar analytic germ $\psi_x$
fixing zero, with nonzero first derivative. Central covariance gives
$\psi_{Tx}(cz)=c\psi_x(z)$ for a fixed $c>1$, where $T$ preserves
that same finite probability. Its measurable Taylor coefficients
satisfy $b_j(Tx)=c^{1-j}b_j(x)$. For $j\ge2$, stationarity forces
$b_j=0$ almost everywhere: invariance of the distribution and
iteration of the strict contraction make every positive tail
probability zero. Hence $\psi_x(z)=\lambda(x)z$ as a germ;
constancy of $\lambda$ is unnecessary. The normalized Haar lift
agrees with the same analytic collar representative. Analytic
uniqueness on actual inverse-chart overlaps extends this return
identity from almost every Haar parameter to the actual local
sheets. Every sheet change thus commutes with scalar dilation,
jointly in time and the collar point.

To descend the flow, restrict to an open domain where all lifts
are legal. The cover over the compact zero section is finite,
and the joint sheet identities hold near every pair of lifts.
Shrinking gives a common safe open time domain on which scalar
dilation descends through the actual covering. Restricting to the
open set of base points where time zero is legal gives an analytic
local flow on a neighborhood $B$ of the original hypersurface,
with identity at time zero and the local addition law. Restrict
the same collar to $\Phi^{-1}(B)$. Its projection is still a finite
surjective cover of $B$; compactness of the angular zero section
supplies a uniform smaller normal strip. Restriction of the
original lattice covariance gives the same exact normal character
on this strip. On every legal lift the descended flow is precisely
$\Phi(\theta,e^{qt}z)$.
\end{proof}

The proof uses the one-dimensional coordinate foliations
of the projective local chart. Each is thickened by the
normal curve to a scalar-contracting two-dimensional stable
surface. It requires neither an ambient Levi linearization
at the start nor a convergent hyperbolic formal linearizer
in all three variables.

The argument uses $q>0$ twice: to construct the two indicated
stable surfaces, and to eliminate the remaining affine
normal-to-root shear. At $q=0$ those surfaces no longer
contract in the normal direction. The neutral collar remains
separate. The regular nonatomic neutral branch is treated by
the compatible-root-curve construction in
Section~\ref{III-series:regular-three}.

\section{Neutral derivative data for \texorpdfstring{$n=3$}{n = 3}}\label{III-series:shear-three}

\subsection{The probability action and its derivative data}

Let $L=\SL(2,\R)$, let $U=(\R^2)^*$, and write $u_r\in U$ for
the element with row parameter $r$. Let $A_0=\{a_t:t\in\R\}$
centralize $L$, and impose the relations
\[
 l u_r l^{-1}=u_{r l^{-1}},\qquad
 a_tu_ra_t^{-1}=u_{\lambda_t^{-1}r},\qquad
 \lambda_t=e^{-\kappa t},\quad\kappa>0.
\]
Put $Q=(L\ltimes U)\rtimes A_0$. Suppose $Q$ acts jointly
measurably preserving a probability $(X,\mu)$. Let $D$ be an
invertible measurable linear cocycle on $\R^2\oplus\R$, with
convention $D(gh,x)=D(g,hx)D(h,x)$, such that
\begin{equation}\label{III-neutral_shear:eq:LA}
 D(l,x)=\begin{pmatrix}l&0\\0&\eta(l,x)\end{pmatrix},\qquad
 D(a_t,x)=\begin{pmatrix}\lambda_t I&0\\0&\delta(t,x)\end{pmatrix}.
\end{equation}
Here $\eta,\delta$ are nonzero real scalars. Write
$E=\R^2\oplus0$ and $N=0\oplus\R$. Assume neutrality in this
frame:
\begin{equation}\label{III-neutral_shear:eq:neutral}
 \lim_{d\to+\infty}\frac1d\log|\delta(d,x)|=0
 \quad\text{for almost every }x.
\end{equation}
In the suspension application $E$ is the full stable space, its
frame is $D_0\phi_x$ from the Haar parametrization, and $N$ is
the zero Oseledec line with a measurable unit frame. Thus
\eqref{III-neutral_shear:eq:LA} does not assume a normal Levi reduction.

Equalities for an uncountable group mean: for every fixed element
the equality holds almost everywhere. Parameter integrations below
use the jointly measurable versions and Fubini. A single conull
set is only required for any subsequently chosen countable group.

\begin{lemma}[Levi components and row stabilizers]\label{III-neutral_shear:lem:components}
The row radical acts trivially on the Levi ergodic-component factor
$Z$. On almost every Levi ergodic component, the stabilizer
$L_r=\{l:r l^{-1}=r\}$ of a nonzero row acts ergodically.
\end{lemma}
\begin{proof}
A jointly measurable standard action admits a compatible Polish topology
on the same space, with the same Borel sets and a continuous action
\cite[Theorem~1.1 and its reformulation, p.~335; proof sketch, pp.~338--339]{BeckerKechris1993}.
This changes neither the probability nor its component factor.
Finite Borel measures on the resulting Polish space are inner regular
by compact sets. Approximation in $L^2$ then gives strong continuity of
the original Koopman action; the same argument applies to its finite
invariant restrictions and to each invariant conditional probability.
For an $L$-invariant vector $f$ in the Koopman representation, choose
$l_n\in L$ with $l_nu_rl_n^{-1}\to e$. Unitary invariance and
strong continuity give
\[
 \|u_rf-f\|=\|l_nu_rl_n^{-1}f-f\|\longrightarrow0.
\]
This proves the first assertion, including preservation of the
conditional probabilities over $Z$. For the second assertion, conjugate the stabilizer of a nonzero
row to the upper unipotent subgroup of $\mathrm{SL}(2,\mathbb R)$.
The unipotent Mautner argument proved in Part~II shows that a vector
fixed by this subgroup in a strongly continuous isometric representation
is fixed by the whole Levi group. Apply this to the indicator of an
invariant measurable set in each conditional $L^2$ space. Levi ergodicity
makes the indicator constant and proves stabilizer ergodicity.
The conditional probabilities and their simultaneous row invariance
are obtained from the original Levi-invariant factor by disintegration;
continuity extends the countable dense subgroup identities to the whole
group. Thus the argument uses the same component measures throughout.
\end{proof}

\subsection{The transported planes determine the regular first jet}

Define the two-plane
\begin{equation}\label{III-neutral_shear:eq:planes}
 P_x(r)=D(u_r,x)^{-1}E.
\end{equation}
For the target splitting, the symbol $E$ on the right denotes
the stable plane at $u_rx$. The cocycle identity gives
\begin{equation}\label{III-neutral_shear:eq:transport}
 D(u_r,x)P_x(r+s)=P_{u_rx}(s).
\end{equation}
The Levi covariance is
\begin{equation}\label{III-neutral_shear:eq:plane-L}
 D(l,x)P_x(r)=P_{lx}(r l^{-1}).
\end{equation}

\begin{proposition}[Plane-pencil dichotomy]\label{III-neutral_shear:prop:pencil}
There is a measurable $Q$-invariant partition into the following
two alternatives. In the first, $P_x(r)=E$ for every fixed $r$
and almost every $x$. In the second, a measurable change of the
normal frame, leaving the stable frame unchanged, gives
\begin{equation}\label{III-neutral_shear:eq:regular-linear}
 D(l,x)=\begin{pmatrix}l&0\\0&1\end{pmatrix},\quad
 D(u_r,x)=\begin{pmatrix}I&0\\-r&1\end{pmatrix},\quad
 D(a_t,x)=\begin{pmatrix}\lambda_t I&0\\0&1\end{pmatrix}.
\end{equation}
In particular the normal scalar cocycle is trivialized in the
second alternative without a Levi cocycle-superrigidity theorem.
\end{proposition}
\begin{proof}
As above, the compatible Polish realization of
\cite[Theorem~1.1, p.~335]{BeckerKechris1993} supplies strong continuity
on the unchanged original measure and its invariant restrictions.
Fix $r\ne0$. On the set $P_x(r)\ne E$, the line $P_x(r)\cap E$
is equivariant under $L_r$. The unipotent action on $\RP^1$ has
only one invariant probability, concentrated at its fixed line
$\ker r$. Indeed every projective orbit converges to that line;
invariance and bounded convergence prove the assertion. Therefore
\[
 P_x(r)\cap E=\ker r\quad\text{almost everywhere on this set}.
\]
The three possibilities
\[
 E,\qquad (\ker r)\oplus N,\qquad
 \{(u,z):z=h_x(r)ru\},\quad h_x(r)\ne0,
\]
are $L_r$-invariant measurable events. By
Lemma~\ref{III-neutral_shear:lem:components} each is constant on almost every Levi
component. Levi covariance and transitivity on nonzero rows make
the choice independent of $r\ne0$. These statements are first
read for almost every $(x,r)$, which suffices below.

The second, vertical possibility is impossible. Since $U$ fixes
the Levi component, intersecting \eqref{III-neutral_shear:eq:transport} over almost
every $s$ shows that $D(u_r,x)N=N$. On the quotient by $N$, its
invertible two-by-two matrix $T$ would satisfy
\[
 T\ker(r+s)=\ker s\quad\text{for almost every }s.
\]
Hence $\det(sT,s+r)=0$ for almost every $s$, and therefore
identically as a polynomial. Its degree-one part is
$\det(sT,r)=0$, forcing $T$ to have rank at most one when
$r\ne0$, a contradiction.

Consider the graph possibility, and put $g_x(r)=1/h_x(r)$.
The projective normal to $P_x(r)$ is $[-r:g_x(r)]$.
For almost every pair of rows $r,s$, equation
\eqref{III-neutral_shear:eq:transport} makes the three normals at
$r,s,(r+s)/2$ collinear: after transport by $u_r$ they are the
normals at $0,s-r,(s-r)/2$, which lie on one projective line.
For independent $r,s$, their determinant gives
\[
 g_x((r+s)/2)=\tfrac12\{g_x(r)+g_x(s)\}.
\]
The measurable almost-everywhere Jensen equation implies
\begin{equation}\label{III-neutral_shear:eq:affine-g}
 g_x(r)=a(x)+r\xi(x)\quad\text{for almost every }r.
\end{equation}
For completeness, equating the midpoint expressions for the pairs
$(r+t,s)$ and $(r,s+t)$ shows that
$g(r+t)-g(r)$ is independent of $r$, modulo null sets. Its value
is a measurable almost-everywhere additive function of $t$, hence
has a linear representative. This is the usual measurable Cauchy
equation, and proves the assertion without an initial integrability
assumption. The affine coefficients are unique and measurable in $x$.

Equation~\eqref{III-neutral_shear:eq:plane-L} implies
\[
 a(lx)=\eta(l,x)^{-1}a(x),\qquad
 \xi(lx)=\eta(l,x)^{-1}l\xi(x).
\]
The event $a=0$ is $L$-invariant. On it $\xi\ne0$, so
$x\mapsto[\xi(x)]$ would push a probability to an $L$-invariant
probability on $\RP^1$, which does not exist. Thus $a\ne0$
almost everywhere. The vector $\xi/a$ is then $L$-equivariant.
An invariant probability in the standard module is supported at
zero: recurrence under the two elementary unipotents forces their
common fixed vectors. Consequently $\xi=0$.

Replace the normal coordinate $z$ by $a(x)z$. Now
\[
 P_x(r)=\{z=ru\},\qquad D(l,x)=\operatorname{diag}(l,1).
\]
The linear map $T_r(u,z)=(u,z-ru)$ satisfies precisely
\eqref{III-neutral_shear:eq:transport}. Thus $T_r^{-1}D(u_r,x)$ fixes an open
set of projective hyperplanes and is scalar:
$D(u_r,x)=c(r,x)T_r$. For $l\in L_r$, Levi covariance makes
$c(r,lx)=c(r,x)$. Lemma~\ref{III-neutral_shear:lem:components} makes $c$ a
function of $r$ and $Z$. Row invariance and the $U$ cocycle identity give
$c(r+s,x)=c(r,x)c(s,x)$ for each fixed pair $r,s$ almost everywhere.
Covariance under $-I\in L$ gives $c(-r,x)=c(r,x)$, while $c(0,x)=1$.
Thus $c(r,x)^2=1$, and applying the product identity to $r/2,r/2$
gives $c(r,x)=1$. Only finitely many null-set identities are combined
for each fixed row. This proves the middle identity in
\eqref{III-neutral_shear:eq:regular-linear}.

Finally write the normal entry of the newly framed $D(a_t,x)$
as $\widetilde\delta(t,x)$. The relation
$a_tu_ra_t^{-1}=u_{\lambda_t^{-1}r}$, compared in the
lower-left block of its derivative, gives
$\widetilde\delta(t,x)=1$.

The constant-plane event is preserved by $L,U,A_0$, so the
partition is $Q$-invariant. Identities initially obtained for
almost every row hold for each fixed row modulo null sets:
the set of rows satisfying the corresponding cocycle identity
is a subgroup containing a conull set, hence is all of $U$.
For the other elements one uses their covariance identities
with almost every row. This also justifies the fixed-element
interpretation of the normal-frame equations above.
\end{proof}

\subsection{Neutrality excludes a reverse shear}

\begin{proposition}[The constant-plane derivative]\label{III-neutral_shear:prop:zero-linear}
In the constant-plane alternative of Proposition~\ref{III-neutral_shear:prop:pencil},
\begin{equation}\label{III-neutral_shear:eq:zero-linear}
 D(u_r,x)=\begin{pmatrix}I&0\\0&n(r,x)\end{pmatrix}
\end{equation}
for nonzero real scalars $n(r,x)$. No trivialization of the
normal $L$ or $U$ scalar cocycle is asserted.
\end{proposition}
\begin{proof}
Use the same compatible Polish realization
\cite[Theorem~1.1, p.~335]{BeckerKechris1993} on the original finite
invariant branch restriction.

Write $D(u_r,x)=\left(\begin{smallmatrix}A(r,x)&B(r,x)\\
0&n(r,x)\end{smallmatrix}\right)$. For $l\in L_r$,
$A(r,lx)=lA(r,x)l^{-1}$. Recurrence under this unipotent
conjugation forces $A(r,x)$ into its fixed subspace; in that
subspace it is invariant under $L_r$. By
Lemma~\ref{III-neutral_shear:lem:components}, $A(r,x)$ depends only on $r,Z$.
Row invariance and the cocycle identity give
$A(r+s,x)=A(r,x)A(s,x)$ for each fixed pair almost everywhere.
Levi covariance under $-I$ gives $A(-r,x)=A(r,x)$, and
$A(0,x)=I$. Hence $A(r,x)^2=I$, and applying the product identity
to $r/2,r/2$ gives $A(r,x)=I$. This uses finitely many identities
for each fixed row and does not require differentiating a selected
component homomorphism.

Let $J=\left(\begin{smallmatrix}0&-1\\1&0\end{smallmatrix}\right)$.
Projectivizing the covariance
$B(r,lx)=lB(r,x)/\eta(l,x)$ for $l\in L_r$, and using
unipotent projective recurrence as before, gives
\[
 B(r,x)=\beta(r,x)Jr^{\mathsf T}.
\]
Commutativity of $u_r,u_s$, with independent $r,s$, yields
\[
 \beta(s,x)=n(r,x)\beta(s,u_rx),\qquad
 n(s,x)\beta(r,u_sx)=\beta(r,x).
\]
The upper-right block of $D(u_{r+s},x)$ then gives
\[
 \beta(r+s,x)J(r+s)^{\mathsf T}
 =\beta(r,x)Jr^{\mathsf T}+\beta(s,x)Js^{\mathsf T}.
\]
Independence implies $\beta(r,x)=\beta(s,x)=\beta(r+s,x)$.
Comparing any nonzero row with one of the two coordinate rows
therefore gives a single measurable scalar $\beta(x)$, chosen from
the column at the first coordinate row.
On its nonzero set,
\[
 \beta(u_rx)n(r,x)=\beta(x),\qquad
 \beta(lx)\eta(l,x)=\beta(x).
\]
Changing the normal coordinate by $z\mapsto\beta(x)z$
makes $D(l)=\operatorname{diag}(l,1)$ and
\[
 D(u_r)(u,z)=(u+Jr^{\mathsf T}z,z).
\]
The central conjugation relation now forces the normal entry
of $D(a_t)$ to be $\lambda_t^2$: the reverse shear transforms
with factor $\lambda_t/\widetilde\delta(t,x)$, whereas the
row parameter transforms with factor $\lambda_t^{-1}$.
The nonzero event is $Q$-invariant by these same relations.

This contradicts \eqref{III-neutral_shear:eq:neutral}. Indeed the new normal
entry is measurably cohomologous to the original one. Along
recurrent returns to a set where $|\beta|$ and $|\beta|^{-1}$
are bounded, that coboundary contributes zero after division
by the return time, whereas the constant entry
$\lambda_1^{2d}$ has exponent $-2\kappa$. Thus $\beta=0$
almost everywhere and \eqref{III-neutral_shear:eq:zero-linear} follows.
\end{proof}

\begin{corollary}[The exceptional exponent for a reverse shear]
\label{III-neutral_shear:cor:reverse-exponent}
Replace neutrality by the assumption that the normal exponent is
a fixed number $\chi$. If $E$ is $U$-invariant and
$\chi\ne-2\kappa$, then \eqref{III-neutral_shear:eq:zero-linear} still holds.
\end{corollary}
\begin{proof}
The proof of Proposition~\ref{III-neutral_shear:prop:zero-linear} uses neutrality
only in its last paragraph. A nonzero reverse shear would give
normal exponent $-2\kappa$, by the same recurrent-return argument.
\end{proof}

\subsection{One-dimensional transverse curves and the exceptional weights}

The reverse-shear calculation has a higher-order version for
\emph{curves}. Here the relevant tensor has a one-dimensional
source, so the pure-power argument for surfaces below is replaced
by an explicit resonance calculation.

\begin{proposition}[Parabolic coherence of nonresonant transverse curves]
\label{III-neutral_shear:prop:curve-coherence}
Suppose \eqref{III-neutral_shear:eq:LA} holds with normal exponent $\chi$, that
$E$ is $U$-invariant, and that
\begin{equation}\label{III-neutral_shear:eq:curve-nonresonance}
                  k\chi\ne-2\kappa\qquad(k=1,2,\ldots).
\end{equation}
Let $F$ be an analytic germ cocycle with derivative $D$, and
let $C_x$ be measurable analytic curve germs tangent to $N$,
equivariant under $L$ and $A_0$. Then these curves are equivariant
under the full $Q$ action.
\end{proposition}
\begin{proof}
Corollary~\ref{III-neutral_shear:cor:reverse-exponent} gives
$D(u_r,x)=\operatorname{diag}(I,n(r,x))$. Thus the curves and
their $U$-transports have the same first jet. Suppose agreement
is known through degree $k-1$, $k\ge2$. Define their first
difference using the \emph{source} normal parameter. More
precisely, subtract the degree-$k$ tangential graph coefficients
at the target and precompose with the $k$th power of the normal
derivative $n(r,x)$. The resulting tensor is
\[
 T_k(r,x)\in\operatorname{Hom}(\Sym^k N_x,E_{u_rx})\simeq\R^2.
\]
It is intrinsic at the first different degree. Levi covariance,
for $l\in L_r$, is
\[
 T_k(r,lx)\eta(l,x)^k=lT_k(r,x).
\]
Projectivization removes the scalar. Recurrence for the
unipotent $L_r$ therefore puts its value in the fixed line:
\[
             T_k(r,x)=\tau_k(r,x)Jr^{\mathsf T}.
\]
The source-parameter convention gives the composition identity
\begin{equation}\label{III-neutral_shear:eq:curve-addition}
 T_k(r+s,x)=T_k(s,x)+n(s,x)^k T_k(r,u_sx).
\end{equation}
Compare this identity with its version with $r,s$ interchanged.
For independent rows, the coefficients of $Jr^{\mathsf T}$ and
$Js^{\mathsf T}$ first give
$n(s,x)^k\tau_k(r,u_sx)=\tau_k(r,x)$, and then
\[
 \tau_k(r+s,x)J(r+s)^{\mathsf T}
       =\tau_k(r,x)Jr^{\mathsf T}+\tau_k(s,x)Js^{\mathsf T}.
\]
Consequently $\tau_k(r,x)=\tau_k(x)$ for almost every $(r,x)$.
Its nonzero event is invariant under $L$ and $U$.

Central covariance, with $r'=\lambda_t^{-1}r$, reads
\[
 T_k(r',a_tx)\delta(t,x)^k=\lambda_t T_k(r,x),
 \qquad
 \tau_k(a_tx)\delta(t,x)^k=\lambda_t^2\tau_k(x).
\]
This also makes the nonzero event $A_0$-invariant. Recurrent
returns to sets where $|\tau_k|$ and its reciprocal are bounded
would force $k\chi=-2\kappa$ there. The nonresonance hypothesis
therefore gives $T_k=0$ almost everywhere. The subgroup and
Fubini argument used above upgrades the identity to every fixed
row. Equivalently, at each degree the source-normalized discrepancy is
the upper-right block of an actual linear cocycle with normal multiplier
$n(r,x)^k$, Levi multiplier $\eta(l,x)^k$, and central multiplier
$\delta(t,x)^k$. Its stable plane is invariant and its normal exponent
is $k\chi$. Corollary~\ref{III-neutral_shear:cor:reverse-exponent} therefore
kills that block. These cocycles and their compatible jets are obtained
from the given germs by Taylor truncation and inversion of the nonzero
normal coordinate. Induction gives equality of all normalized graph jets;
analytic uniqueness, followed by the inverse of the target normal germ,
produces an invertible analytic change of source parameter between the
original curve germs. This proves the asserted equivariance.
\end{proof}

In the actual positive transverse branch, $\chi=q>0$ and
$\kappa=3/2$ for
$a_t=\operatorname{diag}(e^t,e^{-t/2},e^{-t/2})$. The analytic
boundary supplies the $U$-invariant tangent plane. The transverse
curve is the stable curve of $a_{-1}$. Since $L$ and $A_0$
centralize $a_{-1}$, uniqueness of its stable germ makes this
curve $L A_0$-equivariant. This uses the usual tempered derivative
bounds for each fixed suspension element: commuting maps
preserve exponential convergence to the orbit. Hence
Proposition~\ref{III-neutral_shear:prop:curve-coherence} supplies actual
$Q$-equivariant transverse analytic curves for $q>0$.

This is a statement about embedded curve germs. It does not
normalize their one-dimensional Levi germ cocycle, and it does
not establish analytic dependence of these curves on the
boundary point. Both distinctions matter for constructing a
full collar.

\begin{lemma}[Central nonresonance supplies Levi coherence]
\label{III-neutral_shear:lem:curve-Levi}
Under \eqref{III-neutral_shear:eq:LA}, suppose $C_x$ is an $A_0$-equivariant
analytic curve germ tangent to $N$, with normal exponent $\chi$.
If $k\chi\ne-\kappa$ for every integer $k\ge2$, then $C_x$
is also $L$-equivariant.
\end{lemma}
\begin{proof}
Fix $l\in L$. The curves $F(l,x)C_x$ and $C_{lx}$ have the
same tangent line. Suppose their jets agree through degree
$k-1$. Use the source normal parameter to define the first
difference $T_k(l,x)\in\operatorname{Hom}(\Sym^kN_x,E_{lx})$.
Since $l$ commutes with $a_t$ and the curves are
$A_0$-equivariant,
\[
 T_k(l,a_tx)\delta(t,x)^k=\lambda_t T_k(l,x).
\]
If this measurable vector were nonzero on a set of positive
measure, recurrent returns to sets with its norm bounded above
and away from zero would force $k\chi=-\kappa$. Thus it
vanishes. Induction and analytic uniqueness prove the assertion
for each fixed $l$.
\end{proof}

\begin{corollary}[Nonexceptional negative exponents]
\label{III-neutral_shear:cor:negative-curves}
For the actual negative transverse branch, if
$q\notin\{-3/j:j=1,2,\ldots\}$, the Cartan-equivariant
intersection curves in Section~\ref{III-series:coordinate-curves} are $Q$-equivariant.
\end{corollary}
\begin{proof}
Here $\kappa=3/2$ and $\chi=q$. The excluded set includes
every possible Levi resonance $q=-3/(2k)$, $k\ge2$, so
Lemma~\ref{III-neutral_shear:lem:curve-Levi} gives Levi coherence. It also excludes
every row resonance $q=-3/k$, $k\ge1$.
Proposition~\ref{III-neutral_shear:prop:curve-coherence} applies.
\end{proof}

This corollary neither supplies a normal scalar-cocycle reduction
nor settles the exceptional negative values. The resonance list
records where this argument stops; it does not assert the
existence of a deformation at any listed value.

\subsection{Zero first-order shear implies actual plaque coherence}

Suppose now that $F(q,x)$ is a measurable analytic diffeomorphism
germ cocycle with derivative $D$, and that $S_x$ is a measurable
family of analytic surface germs tangent to $E$. Suppose
\[
 F(l,x)S_x=S_{lx},\qquad F(a_t,x)S_x=S_{a_tx}.
\]
The actual neutral stable plaques have these properties, using the
original-law entropy application above and the neutral plaque construction
of Part II.

\begin{theorem}[Coherence without formal Levi normalization]
\label{III-neutral_shear:thm:zero-coherence}
In the constant-plane alternative, $F(u_r,x)S_x=S_{u_rx}$
as analytic surface germs for every fixed $r$, almost everywhere.
Thus the stable plaques are locally $Q$-coherent.
\end{theorem}
\begin{proof}
Their first jets agree by the constant-plane alternative.
Suppose inductively that agreement holds through degree $k-1$,
where $k\ge2$. The difference of the degree-$k$ submanifold
jets is an intrinsic tensor
\[
 T_k(r,x)\in\operatorname{Hom}(\Sym^k E,N)
\]
at the target center $u_rx$. One can define it by subtracting
normal graph functions in any local coordinates adapted to the
common lower jet. Changing those coordinates acts on their first
different term only through the tangent and normal derivatives.
The same observation proves the composition rule below.

Levi coherence and \eqref{III-neutral_shear:eq:zero-linear} give, for $l\in L_r$,
\[
 T_k(r,lx)=\eta(l,u_rx)\,T_k(r,x)\circ\Sym^k(l^{-1}).
\]
The projective action of $L_r$ on $\Sym^k E^*$ has only one
invariant probability: its unipotent generator is one Jordan
block and every projective orbit converges to the fixed line
spanned by $r^{\otimes k}$. Restrict to the invariant event
$T_k\ne0$ and push forward its probability. Consequently
\begin{equation}\label{III-neutral_shear:eq:pure-power}
 T_k(r,x)=\tau_k(r,x)r^{\otimes k}.
\end{equation}

Since $D(u_r)|_E=I$ and its normal entry is $n(r,x)$, the
first-difference composition rule reads
\[
 T_k(r+s,x)=T_k(r,u_sx)+n(r,u_sx)T_k(s,x).
\]
For independent $r,s$, equation~\eqref{III-neutral_shear:eq:pure-power} puts the
left side in $\R(r+s)^{\otimes k}$ and the right side in
$\R r^{\otimes k}+\R s^{\otimes k}$. These three pure powers
are linearly independent for every $k\ge2$: in coordinates
$r=(1,0),s=(0,1)$, the mixed coefficient of $(u_1+u_2)^k$
is nonzero. Fix any nonzero row $r$ and choose an independent row
$s$, for instance $s=(-r_2,r_1)$ in fixed row coordinates.
The three pure-power identities and the addition identity hold on
one conull set for these fixed rows, after pulling the $r$ identity
back under the measure-preserving map $u_s$. Thus
$T_k(r,u_sx)=0$ almost everywhere. Pulling back under $u_{-s}$ gives
$T_k(r,x)=0$ almost everywhere for this fixed $r$. The zero row
follows from the identity cocycle law. No common conull set for
uncountably many rows is needed.

Induction gives every jet identity for each fixed row. The
countable intersection over jet order is conull. Analytic
surface germs are determined by their full jets, so the
surface identities are actual analytic identities.
\end{proof}

\begin{corollary}[The zero-shear neutral suspension is finite atomic]
\label{III-neutral_shear:cor:atomic}
For a $3$-dimensional neutral line/plane suspension measure,
the constant-plane alternative has finite atomic fiber
conditionals. It is supported on finitely many compact embedded
analytic projective covering hypersurfaces. This does not yet
supply their transverse collar neighborhoods.
\end{corollary}
\begin{proof}
The full stable parameters may be constructed from the same original
conditional graph in \cite[Proposition~7.2, Remark~7.3 and \S7.5]{BRHW}.
Fix its return charts and the stable windows of
\cite[Theorem~5.1(a)--(c), Case~III]{Ruelle1980Author}.
Intrinsic contraction and capture in these windows place every native
fiber point in the chosen Pesin leaf. The scalar return, local inverse
and recurrence extend the physical parameter onto that same leaf.
For the total leaf, Claim~3.9 of \cite{BRHW} and freeness on the base
give a bijection from the unipotent parameter times the native fiber.
Uniform contraction and return-chart capture make this product map
continuous; its ambient immersion lifts to a $C^1$ intrinsic immersion.
A bijective $C^1$ immersion of these finite-dimensional manifolds is a
homeomorphism: compact exhaustion and Baire give a point of inverse
continuity, the local immersion charts give equality of dimensions,
and the inverse function theorem then makes the map open everywhere.
Thus the subordinate atoms have the genuine product topology.
Normalize the original conditional classes in
\cite[Proposition~3.11]{BRHW}; conditional uniqueness gives the canonical
Haar law. Transport it through the physical identification with the
chosen Pesin leaf, retaining the same graph and return identities.

The analytic stable Haar plaques in Section~\ref{II-sec:neutral} supply
$S_x$, their standard tangent frames, and \eqref{III-neutral_shear:eq:LA}--
\eqref{III-neutral_shear:eq:neutral}. For the original Cartan derivative,
only fixed-direction Oseledec growth is needed here: the integrable
joint splitting and its raywise growth are described in
\cite[Theorem~2.4 and the proof discussion, pp.~462--463]{BRHWSmoothErgodic2023}.
The extension from integer rays to each fixed real direction uses the
integrable compact-parameter derivative envelope. No growth assertion
uniform over arbitrary paths, and no interchange of an uncountable
family of directions with a conull set, is used. The normal column is
the same unit zero-exponent vector, so its central scalar has zero
integer exponent without transporting a rate through an unbounded
measurable frame. Here the Haar parameters cover the full intrinsic leaves. Indeed, the equal-dimensional intrinsic immersion
has open image by the finite-$C^2$ inverse function theorem in charts
\cite[Chapter~I, \S5, Theorem~5.2, pp.~13--15]{LangDifferential1995}.
For the same physical family and intrinsic chart, Lusin--Souslin makes
its injective evaluation range Borel. The bad-radius sections are locally
closed in the finite-dimensional parameter space and hence sigma compact.
Arsenin--Kunugui makes their projections Borel; selecting the first admissible
positive reciprocal integer gives a measurable local admission radius
\cite[Theorems~4.5.4 and~5.12.1, pp.~153--154 and~219--220]{SrivastavaBorel}.
In the latter theorem both spaces are Polish. A Standard Borel base admits
a compatible Polish topology with the same measurable sets, while the
intrinsic fiber topology remains fixed. This is the local inverse and
measurable radius construction of Part~II.
The contracting continuation and recurrence argument of
Section~\ref{II-sec:neutral} extends this local admission to the whole
intrinsic leaf. Theorem~\ref{III-neutral_shear:thm:zero-coherence} gives the
local parabolic coherence. Keep the same original conditional graph
throughout. Proposition~3.11(2)--(4) of \cite{BRHW} gives Cartan
covariance, almost-everywhere re-anchoring, and the normalized subordinate
conditionals. Translation of these conditionals preserves the original
probability; continuity in the graph parameter and finite-measure dominated
convergence extend preservation from a dense set to every fixed parameter.
Lifting these maps by counting lattice representatives gives the Haar-graph
disintegration of the original lifted measure.

The original boundary coordinate is also retained. The full line or plane
stabilizer is a maximal parabolic in $\SL(3,\R)$, whose Lie algebra has no
rank-one simple ideal. Hence \cite[Proposition~5.4, p.~970]{BRHWFactors}
identifies its inner and outer conditional groups, and
\cite[Claim~5.5, p.~970]{BRHWFactors} makes physical projection from the
boundary joint measure a measurable isomorphism. Here the orbit conditional
class is that of the original suspension probability: pull its local
Rokhlin rows back through injective quotient flowboxes, compare them with
the positively weighted counting lift, and glue their null classes over a
countable cover. Strictly positive weighting preserves these null classes.
The boundary coordinate therefore belongs to the same original lifted point.

Conditional uniqueness, Fubini and Gauss re-centering
\cite[Claim~8.4]{BRHW} rectify the local image coherence to the exact
projective parameters on connected big cells. The local image predicate is
Borel by injective evaluation on compact parameter balls and
\cite[Theorem~4.5.4, pp.~153--154]{SrivastavaBorel}.
We use the coherent-graph construction of
\cite[Definitions~8.11--8.16 and Claims~8.13--8.22]{BRHW} with this
proved coherence, as follows. Countably many full-measure word conditions
retain the unipotent, parabolic and lattice laws. Full native stable
cylinders give complete word lifts; their saturation covers the native
group, and its quotient by the full parabolic covers $\RP^2$.
A common affine chart and stable-leaf uniqueness make the physical
realization injective. The pointed-cover lifting argument
\cite[Propositions~1.14, 1.33--1.34 and Example~1.43]{Hatcher}
identifies each connected cover with $S^2$ or $\RP^2$, so each realization
is a compact embedding. The original boundary coordinate separates the
resulting physical leaves. Countable dense word samples give Borel codes
for their compact images. Their labels are unchanged by the
measure-preserving unipotent and parabolic maps; these subgroups generate
the ambient group, so the original conditional distribution of labels is
constant. Thus we obtain the invariant ergodic probability on the family
of disjoint compact hypersurfaces used in Section~\ref{II-sec:neutral}.

The topological input to that finite-orbit argument is obtained on these
same embedded hypersurfaces. Compact topological manifolds have finitely
generated integral homology by \cite[Corollary~E.5, pp.~538--539]{Bredon1993}.
The coefficient exact sequence for
$0\to\mathbb Z\xrightarrow{2}\mathbb Z\to\mathbb F_2\to0$
\cite[Theorem~5.6 and Example~5.8, pp.~178--181]{Bredon1993}
therefore gives finite mod-two homology; the top mod-two fundamental class
is supplied by \cite[Theorem~3.26, p.~236]{Hatcher}.
For a finite disjoint family $N_1,\ldots,N_k$ in the connected closed
$3$-manifold, its union $K$ is a compact topological manifold and hence
an ENR \cite[Corollary~18.4.7, p.~449]{TomDieck}. Finite additivity
\cite[Proposition~10.2.1, p.~247]{TomDieck} gives
$H_2(K;\mathbb F_2)\cong\mathbb F_2^k$, while higher homology vanishes
\cite[Corollary~7.12, p.~346]{Bredon1993}. The natural field universal
coefficient isomorphism \cite[pp.~196--199]{Hatcher} and the canonical
compact-subset duality
\cite[Lemmas~18.3.1--18.3.2 and the mod-two assertion on p.~450]{TomDieck}
identify cohomological restriction with the dual of the actual inclusion
on homology. The exact sequence of the actual pair $(M,M\setminus K)$
\cite[Theorem~2.16 and its pair specialization, pp.~117--120]{Hatcher}
then gives finitely many complement components and
\[
 \#\pi_0(M\setminus K)+
 \dim_{\mathbb F_2}\operatorname{span}\{[N_1],\ldots,[N_k]\}=k+1.
\]
Thus the finite-cut formula is a consequence of these original topological
results for the actual embeddings; no local flatness or replacement of
the hypersurfaces is required.

Every finite-index subgroup of $\Gamma$ is a lattice in
$\SL(3,\R)$ and has property $(T)$, hence finite abelianization
\cite[Theorem~1.4.15, pp.~54--55; Theorem~1.7.1 and
Corollary~1.7.2, pp.~65--67 of the 23 February 2007 manuscript]{BHV}.
The finite-orbit lemma for disjoint compact hypersurfaces, Lemma~\ref{II-loc:lem:finite-hypersurface-probability}, therefore makes this family
finite. Each hypersurface is a finite cover of $\RP^2$, so its
fiber over a boundary parameter is finite. This is finite
atomicity of the original conditional probabilities. The
coordinate description theorem in Section~\ref{III-series:coordinate-curves} supplies the
analytic projective identifications.
\end{proof}

\subsection{Finite jets of the regular neutral chart}

For an ergodic nonatomic neutral measure, Corollary~\ref{III-neutral_shear:cor:atomic}
excludes the constant-plane alternative. Replacing the row parameter
by its negative in \eqref{III-neutral_shear:eq:regular-linear} gives the usual linear
model $u_r:(u,z)\mapsto(u,z+ru)$. Thus its full derivative cocycle
has the regular neutral form, including the normal scalar
trivialization. This removes the two first-order uses of the
higher-rank Levi in the previous proof.

The proposition is a derivative statement. The following lemma
separates finite-order existence from compatibility and convergence.
Write $\rho$ for the matrices in \eqref{III-neutral_shear:eq:regular-linear}.

\begin{lemma}[Uniqueness of the regular formal reduction]
\label{III-neutral_shear:lem:regular-unique}
For every $k\ge2$, there is no nonzero measurable equivariant map
from a probability $Q$-space to
$\operatorname{Hom}(\Sym^k\R^3,\R^3)$ with the conjugation
action induced by $\rho$. Consequently a tangent-to-identity
reduction to $\rho$ is unique at every finite order. If such a
reduction exists at every order, the reductions are automatically
compatible and define a measurable formal reduction.
\end{lemma}
\begin{proof}
Recurrence under two generators of the unipotent row group makes
an equivariant homogeneous polynomial vector field fixed by that
group pointwise. With the sign of the row parameter immaterial,
its infinitesimal generators are $(ru)\partial_z$. Commuting with
both gives precisely
\[
 H(u,z)=\bigl(f(u)u,\ f(u)z+g(u)\bigr),
 \qquad \deg f=k-1,\quad\deg g=k.
\]
Indeed the tangential component is independent of $z$, and the
normal-component equation for every row makes it a scalar multiple
of $u$ with the stated normal derivative. Conjugation by
$\operatorname{diag}(\lambda I,1)$ multiplies the $f$ term by
$\lambda^{-(k-1)}$ and the $g$ term by $\lambda^{-k}$.
Both weights differ from one. Probability recurrence therefore
forces both coefficient functions to vanish.

The first different coefficient of two finite-order reductions
would be such an equivariant map. This proves uniqueness by
induction. Truncating an order-$(k+1)$ reduction gives an order-$k$
reduction, so uniqueness proves compatibility whenever existence
is known at each order. No inverse-limit selection problem remains.
\end{proof}

At the first untreated degree $k$, normalize the already treated
degrees and write the cocycle as
\[
 F(q,x)=(\id+c_k(q,x)+O(|(u,z)|^{k+1}))\rho(q).
\]
The finite-order normalization equation is
\begin{equation}\label{III-neutral_shear:eq:finite-obstruction}
 c_k(q,x)=h_k(qx)-\rho_k(q)h_k(x),\qquad
 \rho_k(q)H=\rho(q)\circ H\circ\rho(q)^{-1}.
\end{equation}
The cocycle identity is
$c_k(qh,x)=c_k(q,hx)+\rho_k(q)c_k(h,x)$.
Lemma~\ref{III-neutral_shear:lem:regular-unique} proves uniqueness of $h_k$, not
its existence; a vanishing invariant space alone does not imply
vanishing first cohomology. For the actual suspension cocycle,
Section~\ref{III-series:regular-three} supplies existence geometrically: its root
curves agree on intersection lines, and polynomial extension
constructs all homogeneous coefficients of the ambient chart.

If \eqref{III-neutral_shear:eq:finite-obstruction} is solvable at every order for
the actual suspension cocycle, the resulting formal chart is
calibrated on the stable plane. To see this, compare its restriction
there with the analytic Haar parametrization. At the first different
degree $j\ge2$, central covariance has weights
$\lambda^{1-j}$ and $\lambda^{-j}$, so recurrence gives equality.
Lemma~\ref{II-loc:lem:hyperplane-convergence}
then proves convergence of the full formal chart. Thus the
construction of compatible coefficients in the next section
also supplies the analytic chart.

There is one additional two-dimensional invariant-theory correction
in the later stable-transition argument. Once calibrated formal
Levi charts have been obtained, a column-stabilizer calculation
allows dependence on $\det(s,u)$; it cannot simply be copied from
dimension $n\ge 4$. The following observation supplies the
needed replacement at the level of a commuting column family.

\begin{lemma}[Commuting columns remove the alternating terms]
Suppose a formal family of vector fields is linear in
$s\in\R^2$, equivariant under $(u,z)\mapsto(lu,z)$, and has
value $s$ in its tangential component at $(u,z)=(0,0)$.
If the fields commute, they have the form
$Y_s=k(z)s\cdot\partial_u$, with $k(0)=1$.
\end{lemma}
\begin{proof}
The scalar invariant of two columns is $\det(s,u)$. Linearity
in $s$ gives the full list of possible covariants
\[
 Y_s=k(z)s\cdot\partial_u+
 \det(s,u)\{\ell(z)E_u+m(z)\partial_z\}.
\]
No higher power of the determinant is linear in $s$. Direct
differentiation of $[Y_{e_1},Y_{e_2}]=0$ gives
$mk=0$ and $3\ell k-mk'=0$. Since $k$ is a unit,
$m=\ell=0$. If the input family is analytic, then so is $k$,
which is extracted by evaluating its tangential component at $u=0$.
\end{proof}

This observation requires the formal Levi covariance and a genuine
commuting parameter family. It does not provide those hypotheses
from the derivative reduction alone.

\section{Regular neutral charts for \texorpdfstring{$n=3$}{n = 3}}\label{III-series:regular-three}

\subsection{The actual stable graphs and their transported planes}

Work on an ergodic regular neutral suspension component.
By Section~\ref{III-series:shear-three}, choose the normal frame so that the derivative
cocycle is
\[
 \rho(l)(u,z)=(lu,z),\qquad
 \rho(u_r)(u,z)=(u,z-ru),\qquad
 \rho(a_t)(u,z)=(\lambda_tu,z),\quad \lambda_t=e^{-3t/2}.
\]
Let $F(g,x)$ denote the actual centered analytic fiber germ.
Let $\phi_x:\R^2\to M_x$ be the analytic full stable Haar
parametrization, with derivative $u\mapsto(u,0)$.
It is exactly equivariant under $L$ and $A_0$ in these parameters.
For every column direction $v\ne0$, its restriction
$t\mapsto\phi_x(tv)$ is the corresponding analytic root graph.

For a row $r$, define the actual analytic parametrization
\[
                    \phi_x^r(u)=F(u_r,x)^{-1}\phi_{u_rx}(u).
\]
Its derivative is $u\mapsto(u,ru)$. The identities of the
germ cocycle and the stable Haar parameters give
\begin{align}
 F(u_t,x)\phi_x^{r+t}(u)&=\phi_{u_tx}^{r}(u),\label{III-regular_neutral_chart:eq:row}\\
 F(l,x)\phi_x^r(u)&=\phi_{lx}^{r l^{-1}}(lu),\label{III-regular_neutral_chart:eq:Levi}\\
 F(a_t,x)\phi_x^r(u)&=\phi_{a_tx}^{\lambda_t^{-1}r}(\lambda_tu).
 \label{III-regular_neutral_chart:eq:central}
\end{align}
These identities do not presuppose an ambient formal chart.

\subsection{Rows preserve the root curves that they annihilate}

\begin{lemma}[The orthogonal row-column identity]\label{III-regular_neutral_chart:lem:root}
If the row $r$ and column $v\ne0$ satisfy $rv=0$, then
\[
                    F(u_r,x)\phi_x(tv)=\phi_{u_rx}(tv)
\]
as analytic curve germs, almost everywhere for every fixed
pair $(r,v)$.
\end{lemma}
\begin{proof}
Choose a Cartan in $LA_0$, using the basis $(v,w)$ of the
tangential plane, for which the column weight $\beta_v<0$
and the other column weight $\beta_w>0$. The normal exponent
is zero. Hence the fiber stable manifold is precisely the
analytic root curve in direction $v$.
Since $r$ annihilates $v$, it lies in the row weight
$-\beta_w<0$. Thus this same Cartan element $a$ contracts
both the column $v$ and the row subgroup containing $u_r$.

Use the original controlled section obtained from property~(T)
\cite[Theorem~1.4.15, pp.~54--55]{BHV} and the exponential-domain
and word-comparison theorems
\cite[Theorems~5.1 and 5.3, pp.~183--185]{dLS}.
Its quotient-section moment pulls back under the original Haar
projection. The compact-parameter suspension estimates show that
$F(u_r,x)$ transports this stable curve to the stable curve
at $u_rx$. Indeed after $n$ iterations the displacement is
$u_{r_n}=a^nu_ra^{-n}$, which lies in a compact row segment.
The actual ambient word length is bounded on this segment, and the
fiber maps and their inverses have subexponential Lipschitz
bounds along the two fixed orbits of $x$ and $u_rx$.
They preserve exponential convergence of pairs of fiber points.

Write the induced scalar parameter change as $\psi_{r,x}(t)$.
Its derivative is one, since $\rho(u_r)$ fixes $(v,0)$.
The exact root parametrization scales under $a$ by
$\theta=e^{\beta_v}\in(0,1)$, so
\[
           \psi_{r_n,a^nx}(\theta^nt)=\theta^n\psi_{r,x}(t).
\]
The root charts, their local inverses, and the fiber maps for
$u_{r_n}$ have tempered analytic bounds along these two
orbits. For the return maps this follows from the word bound
by the section heights at the two endpoints; both endpoints
lie on fixed Cartan orbits. Thus each fixed-degree coefficient
of $\psi_{r_n,a^nx}$ grows subexponentially.
The degree-$j$ coefficient in the displayed identity is
$\theta^{n(1-j)}$ times the original coefficient, so it
vanishes for $j\ge2$. This proves $\psi_{r,x}(t)=t$.
\end{proof}

\begin{corollary}[Agreement on intersection lines]\label{III-regular_neutral_chart:cor:intersection}
For every fixed pair of rows $r,s$, almost everywhere,
\[
                  \phi_x^r(tv)=\phi_x^s(tv)
                  \qquad\bigl((r-s)v=0\bigr)
\]
as curve germs.
\end{corollary}
\begin{proof}
Apply Lemma~\ref{III-regular_neutral_chart:lem:root} to the row $r-s$ at the center
$u_sx$, and then apply $F(u_r,x)^{-1}$. The row cocycle
identity gives exactly the asserted equality.
\end{proof}

All parameter statements have the fixed-element, almost-everywhere
meaning. Fubini on the row-pair parameter space gives the
intersection identities for almost every pair $(r,s)$ at almost
every $x$. No common conull set for uncountably many rows is needed.

\subsection{A polynomial determined by its hyperplane restrictions}

\begin{lemma}[Extension from compatible graph planes]\label{III-regular_neutral_chart:lem:polynomial}
Let $r\mapsto f_r(u)$ be a measurable family of homogeneous
degree-$k$ polynomials on $\R^2$, with values in a finite-dimensional
vector space. Suppose, for almost every pair $(r,s)$,
\[
                      f_r|_{\ker(r-s)}=f_s|_{\ker(r-s)}.
\]
There is a unique homogeneous degree-$k$ polynomial $P(u,z)$
on $\R^3$ such that $f_r(u)=P(u,ru)$ for almost every $r$.
\end{lemma}
\begin{proof}
It suffices to treat scalar values. Write $u=(x,y)$ and
$r=(a,b)$, and set $p(a,b,t)=f_{(a,b)}(1,t)$.
The intersection identity, after the nonsingular change
$a'=a+(b-b')t$, says
\[
 p(a,b,t)=p(a+(b-b')t,b',t)
\]
for almost every $(a,b,b',t)$. Thus
$p(w-bt,b,t)$ is independent of $b$ almost everywhere.
It defines a measurable function $F(t,w)$ with
\[
                       F(t,a+bt)=p(a,b,t)
\]
for almost every $(a,b,t)$.

For almost every slope $b$, the restriction of $F$ to almost
every line $w=a+bt$ is a polynomial of degree at most $k$.
Choose two distinct such slopes and use their two line
directions as coordinates. In these coordinates $F$ is
separately polynomial of degree at most $k$ almost everywhere.
Choose $k+1$ good distinct interpolation points in one
coordinate. Lagrange interpolation writes $F$, modulo null
sets, as a bivariate polynomial of degree at most $k$ in each
coordinate. This also establishes local integrability without
having assumed it.

Its restriction to lines of almost every slope has degree at
most $k$, so the bivariate polynomial has total degree at most
$k$: any nonzero top homogeneous part of higher degree would
give a nonzero leading coefficient on all but finitely many
slopes. Consequently
\[
                            P(x,y,z)=x^kF(y/x,z/x)
\]
is a homogeneous polynomial of degree $k$, and it has the
desired restrictions. A polynomial with zero restriction to
almost every graph plane vanishes on an open set, so uniqueness
follows.
\end{proof}

\subsection{The formal chart and its convergence}

\begin{theorem}[Actual regular neutral chart]\label{III-regular_neutral_chart:thm:chart}
There is a unique measurable analytic tangent-to-identity
coordinate germ $H_x:\R^3\to M_x$ such that
\[
 H_x(u,0)=\phi_x(u),\qquad
                      F(g,x)H_x=H_{gx}\rho(g)\quad(g\in Q).
\]
It and its inverse have positive measurable analytic radii.
\end{theorem}
\begin{proof}
Take each homogeneous coefficient of the actual embeddings
$\phi_x^r$. Corollary~\ref{III-regular_neutral_chart:cor:intersection} and
Lemma~\ref{III-regular_neutral_chart:lem:polynomial} give a unique homogeneous polynomial
$H_{k,x}(u,z)$ whose restriction to $z=ru$ is that coefficient.
Its coefficients are measurable: for each $k$ choose finitely
many fixed generic rows for which the restriction operator is
injective and the almost-everywhere identities hold in $x$,
using Fubini. A fixed linear left inverse recovers them from
the measurable coefficients of those restrictions.

This constructs a formal field $H_x=\sum_{k\ge1}H_{k,x}$ with
$H_{1,x}=\id$ and
\[
                             H_x(u,ru)=\phi_x^r(u)
\]
formally for almost every row. Its restrictions are actual
analytic germs on a positive-measure family of graph planes.
Lemma~\ref{II-loc:lem:hyperplane-convergence}
therefore makes $H_x$ convergent: selecting common bounds
and radii on a positive-measure row set and applying the
polynomial/Remez estimate bounds the full degree-$k$
coefficients exponentially in $k$. This argument applies
pointwise at almost every $x$.

The chart identities follow by restricting to the planes
$z=ru$ and using \eqref{III-regular_neutral_chart:eq:row}--\eqref{III-regular_neutral_chart:eq:central}.
They hold formally because the restriction family determines
each polynomial coefficient. Convergence makes them actual
analytic germ identities. The inverse function theorem gives
the inverse and positive measurable radii.
Uniqueness follows either from the same calibrated restrictions
or from the formal uniqueness lemma in Section~\ref{III-series:shear-three}.
\end{proof}

\subsection{Stable transitions}

The chart just constructed gives the finite-order normalization
of the regular neutral action. When comparing centers along
the stable plane, the dimension-two alternating covariant
must also be retained.

Continue $H_x$ to a neighborhood of its whole stable plane by
central contraction and recurrence to sets with bounded chart
radii, as in Section~\ref{II-sec:neutral}. The transition
\[
                         P_v(x)=H_x^{-1}v^{-1}H_{\widehat v x}
\]
is an analytic germ from $0$ to $(v,0)$.
Row- and column-stabilizer recurrence now uses the linear
Levi charts supplied by Theorem~\ref{III-regular_neutral_chart:thm:chart}.
For a fixed nonzero column, every Taylor coefficient is
constant on Levi ergodic components, by Howe--Moore.
Write their factor as $Z$. Mautner's argument makes
$\widehat V$ act trivially on $Z$.

The invariant-theory calculation in Section~\ref{III-series:negative-three} therefore
applies componentwise:
\[
 P_v(u,z)=\bigl(A_Z(\det(v,u),z)u+B_Z(\det(v,u),z)v,\
                                       C_Z(\det(v,u),z)\bigr).
\]
It supplies joint analyticity in the column parameter near zero.
The column group law supplies a commuting parameter family.
These local identities can first be checked on a set of
positive conditional measure with bounded chart domains;
measure continuity of the $\widehat V$ action ensures such
triple overlaps for sufficiently small parameters.
Their deterministic analytic coefficients then give the
identities as germs. The commuting-column calculation in
Section~\ref{III-series:shear-three} gives
\[
                              P_v(u,z)=(u+k_Z(z)v,z).
\]
Central covariance makes $k_Z$ invariant under the residual
central action on $Z$, and parabolic ergodicity makes it a
single germ $k(z)$. The mixed-root commutator in Proposition~\ref{II-loc:prop:neutral-stable-transition} determines it: in its positive-shear
normal-coordinate convention, $k(z)=1-z$.
The sign change from the displayed negative-shear convention
above is just $z\mapsto-z$.

Thus the exact stable-transition hypothesis of the regular
neutral description theorem is available for $n=3$.
The local orbit-chart construction in
Proposition~\ref{II-loc:prop:neutral-radial-neighborhoods} applies.
Its quotient has dimension $3$, and the boundary-coordinate
uniqueness argument applies since $3=2\cdot3-3$.
The Euler-flow argument in
Theorem~\ref{II-loc:thm:neutral-regular-recognition} gives the
Hopf manifold. The constant-plane case is treated
separately below.

\section{Neutral projective hypersurfaces for \texorpdfstring{$n=3$}{n = 3}}\label{III-series:axis-three}

\subsection{Derivative normalization and the formal axis}

Let $X=G/\Gamma$ be one recognized boundary suspension after
the finite covering choices in Section~\ref{III-series:relative-three}. Let
$Q=(L\ltimes U)\rtimes A_0$, $L=\SL(2,\R)$.
The normal exponent is zero. The normal quotient derivative
is trivialized by the relative scalar theorem in Section~\ref{III-series:relative-three}.
The central zero Oseledec line is separated from its stable
boundary tangent plane, so it is invariant under $L A_0$.
Lift the normalized quotient frame to this line.
The reverse-shear exclusion in Section~\ref{III-series:shear-three} then gives
\begin{equation}\label{III-neutral_boundary:eq:derivative}
 D(l,x)=\operatorname{diag}(l,1),\qquad
 D(u_r,x)=I,\qquad
 D(a_t,x)=\operatorname{diag}(\lambda_t I,1),
 \quad \lambda_t=e^{-3t/2}.
\end{equation}
Choose analytic ambient charts adapted to the actual boundary.
The derivative frames and the analytic coefficient heights are
tempered along $a_1$, by the derivative-increment argument.

\begin{proposition}[A unique formal parabolic transverse curve]
\label{III-neutral_boundary:prop:formal}
There is a unique measurable formal embedded curve germ
$\mathcal C_x$, tangent to the normal line, that is equivariant
under $A_0$. It is equivariant under the full parabolic $Q$.
\end{proposition}
\begin{proof}
Write the curve as a formal graph $u=f_x(z)$, with
$f_x(0)=f_x'(0)=0$. Suppose the graph is invariant under
$a_1$ through degree $k-1$. Comparing its two degree-$k$
graph coefficients gives
\[
                         f_{a_1x,k}-\lambda_1 f_{x,k}=R_k(x).
\]
The remainder is tempered, since it uses finitely many
tempered coefficients of the actual map and the already
constructed finite graph jet. The convergent series
\[
              f_{x,k}=\sum_{d\ge1}\lambda_1^{d-1}R_k(a_{-d}x)
\]
solves the equation and is tempered. This constructs the
formal curve at every order.

Uniqueness does not require a moment assumption. At the first
different degree of two measurable invariant graphs, their
difference satisfies $v(a_1x)=\lambda_1v(x)$.
Probability recurrence forces $v=0$.
Every $a_t$ transports the constructed curve to another
$a_1$-equivariant curve, so uniqueness gives full $A_0$
equivariance. The same argument applies to every $l\in L$,
because it commutes with $A_0$ and preserves the normal line.

The formal version of the curve-coherence calculation in
Section~\ref{III-series:shear-three} applies. The same-space Polish realization
\cite[Theorem~1.1 and its reformulation, p.~335; proof sketch, pp.~338--339]{BeckerKechris1993}
preserves the original Borel structure and probability when using the
finite-measure recurrence statements. At the first degree $k$ where a row
might fail to preserve the curve, unipotent row-stabilizer
recurrence puts the defect in the line $Jr^{\mathsf T}$.
Commutativity makes its coefficient independent of the row,
and central covariance would force
$k\cdot0=-3$. This is impossible. Induction gives row
equivariance in every degree, hence full formal $Q$ coherence.
No convergence of the curve has been used.
\end{proof}

\begin{proposition}[The formal transverse dynamics]
\label{III-neutral_boundary:prop:parameter}
The curve in Proposition~\ref{III-neutral_boundary:prop:formal} has a measurable
formal parameter $i_x(z)$ for which
\[
 F(h,x)i_x=i_{hx}\quad(h\in LU),\qquad
 F(a_t,x)i_x(z)=i_{a_tx}(\Phi_t^{B(z)\partial_z}(z)),
\]
where $B(z)$ is one deterministic formal series with
$B(0)=B'(0)=0$.
\end{proposition}
\begin{proof}
The induced formal one-dimensional cocycle on the curve has
derivative one on $LU$. Its first untreated coefficient is
an ordinary real additive cocycle. The relative additive
reduction in Section~\ref{III-series:relative-three}, using ergodicity of $U$, removes
it. Here the same original Haar graph law inherits $U$-ergodicity
from the Haar quotient by Howe--Moore
\cite[Theorem~1.1 and Definition~2.6, pp.~66 and 69]{Ciobotaru};
the real relative-property-$(T)$ input is the pair
$(\mathrm{SL}(2,\mathbb R)\ltimes\mathbb R^2,\mathbb R^2)$ stated in
\cite[Introduction and Definition~1.1, p.~375]{Iozzi}.
For every strongly continuous complex unitary representation, almost
invariant unit vectors, uniformly on each compact set, give a nonzero
translation-fixed vector. The standard-to-contragredient group
identification and its application to the original cocycle are the
relative additive reduction above.
Successive corrections preserve the previous degrees and
give a compatible formal parameter in which $LU$ acts
identically.

Since $LU$ is normal in $Q$, every remaining formal central parameter map is
$LU$-invariant and therefore deterministic. Its finite jets
form measurable one-parameter homomorphisms into finite jet
groups, hence are continuous and have compatible infinitesimal
generators. More explicitly, their identity first jets make the finite
pullback operators unipotent, by the vanishing-order filtration used
in the row-logarithm proof above. The finite operator exponential and
logarithm are inverse by \cite[14.28, p.~287]{Milne}. Their
polynomial powers lie in the actual finite jet group and commute with
truncation by the same polynomial-identity argument. The resulting
compatible logarithms are formal generators: differentiate each finite
jet curve at zero and pass to the compatible coefficient family.
A measurable one-parameter homomorphism agrees with these
powers first at rational times by uniqueness of unipotent roots, and
then at real times by continuity of its faithful finite pullback
representation. The real specialization follows from the real
coefficients of the finite exponential and logarithm polynomials.
These generators define the formal vector field
$B(z)\partial_z$. Its linear coefficient is zero by
\eqref{III-neutral_boundary:eq:derivative}.
\end{proof}

\subsection{Convergence uses the full Levi symmetry}

The embedded formal curve has a noncompact tangential target
derivative: the Levi acts there by its standard
two-dimensional representation. Thus the existing spectral-gap
realization theorem for tangent-to-identity germ
\emph{coboundaries} does not directly make this embedded curve
converge. Nor can one first restrict the actual cocycle to the
curve and apply the one-variable relative theorem: that
restriction requires an actual analytic curve. The convergence theorem in
Section~\ref{III-series:fixed-curves} resolves this issue by a different estimate.
For a curve parameter fixed by $H$, the target coefficient
representation is the same standard-plus-trivial module at
every Taylor order. Spectral-gap expansion makes rescaled jets
tight; a hyperbolic Levi element eliminates their tangential
part, and the median Poincare estimate eliminates their normal
part. This proves convergence after a constant formal
reparametrization, without assuming an ambient Levi chart.

Consequently there are actual embedded analytic curves $i_x$
with positive measurable radii and
\[
 F(h,x)i_x(z)=i_{hx}(z)\quad(h\in H),\qquad
 F(a_t,x)i_x(z)=i_{a_tx}(\Phi_t^{B\partial_z}(z)),
\]
where $B$ is analytic and $B(0)=B'(0)=0$.
The role of the Levi symmetry can be seen from the following
example. The
analytic vector field
\[
 X=(-u_1+z^2)\partial_{u_1}-u_2\partial_{u_2}
                                           +z^2\partial_z
\]
preserves the analytic boundary $z=0$. Its time-one derivative
is $\operatorname{diag}(e^{-1},e^{-1},1)$.
A formal transverse invariant curve is necessarily
$u_2=0,\ u_1=f(z)$, where
\[
 z^2 f'(z)=-f(z)+z^2,\qquad
 f(z)=\sum_{d\ge2}(-1)^{d-2}(d-1)!\,z^d.
\]
It has zero radius of convergence. This example is only a
central dynamical system: it is not a full parabolic cocycle,
and does not satisfy the hypotheses used here.
It shows that a convergent solution cannot be inferred from
the contracting coefficient equation at each fixed order.

\subsection{The stationary central flow}

\begin{lemma}[A common real interval for equivariant curve germs]
\label{III-neutral_boundary:lem:interval}
The normalized $H$-equivariant curves admit measurable analytic
representatives $i_x:(-r,r)\to M_x$ for one $r>0$, after
restriction to a conull set. They are immersions on this interval.
Every fixed $h\in H$ intertwines these representatives on the
whole interval, almost everywhere. If $B=0$, the same assertion
holds for every fixed $q\in Q$.
\end{lemma}
\begin{proof}
Here $M_x$ is the full compact suspension fiber, not a coordinate
ball. The cocycle transports full fibers by analytic
diffeomorphisms. Continue each curve germ maximally along the
real parameter, retaining the immersion condition. Uniqueness
of real analytic continuation makes its maximal interval
well-defined. Its positive and negative endpoints are measurable:
continuation can be expressed through countably many overlapping
coordinate disks with rational centers and bounds. They are
invariant almost everywhere under each fixed element of the ergodic
radical $U$, because global fiber diffeomorphisms preserve
continuation and immersion. Thus both endpoints are almost
everywhere constant and strictly separated from zero.
Choose a smaller common interval. Analytic uniqueness extends
each original fixed-element germ identity across that interval.
When $B=0$, the central parameter is also unchanged.
\end{proof}

\begin{lemma}[A cross reaches the zero parameter]\label{III-neutral_boundary:lem:zero-cross}
Let $r>0$ and let $f(u,z)$ be measurable on $I\times(-r,r)$, where
$I\subset\R^2$ is a nonempty open rectangle. Suppose the $z$ slices
are analytic for almost every $u$, and the $u$ slices are
analytic for almost every $z$. Then there is a nonempty
angular patch on which $f$ has a jointly analytic representative
near $z=0$.
\end{lemma}
\begin{proof}
Shrink $I$ to a relatively compact rectangle. First select a
positive-measure angular set $A$, a complex disk $E$ about zero,
and a uniform bound for the holomorphic $z$ slices on $E$.
Next select a positive-measure set
$D_0\subset E\cap(-r,r)$ and a common product complex
neighborhood $D$ of the closed angular rectangle on which the
holomorphic $u$ slices have a uniform bound. The normal-family
correction in the cross lemma in Section~\ref{III-series:coordinate-curves} makes the
almost-everywhere agreement into an exact bounded holomorphic
cross on $(D\times D_0)\cup(A\times E)$.

Let $A^*$ and $D_0^*$ be the locally pluriregular cores, and retain
$A\cap A^*$ and $D_0\cap D_0^*$. The cores are locally $G_\delta$,
hence Borel, by \cite[p.~222]{NguyenCrossRegularization}; their discarded
parts are locally pluripolar by
\cite[Proposition~3.4.6, p.~37]{JarnickiPflugSeparatelyAnalytic2008}.
Sadullaev's theorem on $C^3$ generating real submanifolds
\cite[Theorem~4, p.~513; definition, p.~503]{Sadullaev1976}, applied
to the standard real coordinate planes and a countable open cover,
makes these discarded parts null for real Lebesgue measure. Thus the
retained sets remain measurable and positive. The general cross theorem
\cite[Theorem~4.3.3, pp.~60--61]{JarnickiPflugSeparatelyAnalytic2008},
with its extremal-function identification
\cite[Proposition~3.4.11(a), p.~39]{JarnickiPflugSeparatelyAnalytic2008},
extends the exact cross, with agreement on these retained sets, to
\[
 \{(u,z):\omega_{A,D}(u)+\omega_{D_0,E}(z)<1\}.
\]
Choose $u_0$ in the support of the real measure restricted to the
retained angular core, so every real ball about $u_0$ retains positive
measure. Local pluriregularity gives $\omega_{A,D}(u_0)=0$.
The function $\omega_{D_0,E}$ is plurisubharmonic and lies in $[0,1]$
\cite[p.~222]{NguyenCrossRegularization}. It vanishes on its nonempty
regular core. The maximum principle on the connected domain $E$
\cite[Remark~3.3.2(j), p.~31]{JarnickiPflugSeparatelyAnalytic2008}
therefore gives $\omega_{D_0,E}(0)<1$, whether or not zero belongs
to the core. Upper semicontinuity makes the sum less than one on a
product neighborhood of $(u_0,0)$. Holomorphy there gives convergent
Taylor germs by \cite[Theorem~2.1.2, p.~6, and (3.1.1), p.~23]{JarnickiPflugSeparatelyAnalytic2008}. Slice uniqueness identifies its
extension with $f$ almost everywhere there, as in
Section~\ref{III-series:coordinate-curves}. For a compact analytic target use an
analytic Euclidean embedding \cite[Theorem~3, p.~470]{Grauert} and argue
componentwise on successive nested zero-centered patches. The resulting
Euclidean representative lies in the closed embedded image everywhere,
since it does so almost everywhere and is continuous; the analytic
embedding inverse then gives the target-valued representative.
\end{proof}

\begin{theorem}[The zero-generator neutral boundary]
\label{III-neutral_boundary:thm:stationary}
If the analytic central generator supplied by Section~\ref{III-series:fixed-curves}
satisfies $B\equiv0$, the recognized boundary has a full-angular
analytic product collar, after the usual finite covering and
orientation choices.
\end{theorem}
\begin{proof}
Work first in the standard native line coordinates. For the dual model,
pull back the original lattice, action and normalized Haar graph laws
by the inverse of the contragredient identification, with its
corresponding basis relabeling, and precompose the same curve field.
This retains the normal parameter and the zero-generator condition.
The inverse change of coordinates transports the resulting collar
and its equivariance back to the original action.

Use Lemma~\ref{III-neutral_boundary:lem:interval}. On the projective
angular boundary, each fixed parameter gives the original Haar-projecting
graph probability $\mu_z$, with one-point fiber conditionals. On an
oriented angular cover retain both antipodal graph probabilities
$\mu_z^+,\mu_z^-$ instead. Each is ergodic under the axis-fixing
parabolic; their equal-weight average will be used for recognition.
All these laws use the same original curve field and normal parameter.

The parameters for which an original graph law is $G$-invariant form a
null set. Indeed, a $G$-invariant one-point graph, lifted to $G\times M$,
is constant in the $G$ variable, with value a $\Gamma$-fixed point.
If this occurred on a positive-measure set of parameters, Fubini would
give $i_g(z)=i_{g'}(z)$ there for almost every pair $g,g'$. Analytic
uniqueness on the common interval would give equality at zero,
contradicting the nonconstant original recognized boundary map.

In the spherical case put $\bar\mu_z=(\mu_z^++\mu_z^-)/2$.
The full parabolic permutes these two original branches, and its sign
element interchanges them. Ergodicity of the axis-fixing central element
on each branch therefore makes the average ergodic under the full
upper Borel. This average is proper. Otherwise its invariant law would
give a unitary $L^2$ action of $G$. The density of a branch in the
average is bounded and fixed by a regular Cartan element of the
axis-fixing parabolic; Mautner's argument
makes that density $G$-invariant. The branch would then be $G$-invariant,
contrary to the preceding exclusion. Thus the relevant full-parabolic
law is proper in either case. Its full stabilizer is exactly that
maximal parabolic.

Apply the original nonresonance and root-graph arguments of
Lemmas~\ref{III-coordinate_recognition:lem:root}
and~\ref{III-coordinate_recognition:lem:graphs} to this law. They supply
the two missing-root fiber exponents. The nonzero vector $i_x'(z)$
is equivariant with parameter derivative one. Original return bounds supply the all-finite-$L^q$ logarithmic derivative
bounds in \cite[Proposition~2.5(H3), p.~956; Claim~2.6, p.~958]{BRHWFactors}.
The full native multiplicative ergodic theorem
\cite[Theorem~3.1, pp.~958--959]{BRHWFactors} gives Oseledec regularity,
and finite-law recurrence forces this vector's exponent to be zero. Dimension three leaves
exactly the two missing-root weights and this zero weight; root-graph
calibration gives the exact projective weights. Theorem~
\ref{III-coordinate_recognition:thm:recognition} now applies on the same
law, with its original conditional kernel retained.

In the projective case that kernel is Dirac, so the recognized finite
covering has degree one. Its analytic inverse gives an analytic
representative of the original angular map. In the spherical case,
the original branch is absolutely continuous in the average and hence
remains carried by the same recognized fibers. Pass to the common
finite-index subgroup preserving the recognized components and sphere
action models. Nonsingular ergodicity of this subgroup on the original
sphere first selects one component and then one of its two antipodal
analytic parametrizations. The original angular map itself therefore
has an analytic representative in its original measure class.

Because $B=0$, the entire curve field descends measurably through the
appropriate angular homogeneous quotient. Pulling back through an
original angular chart gives $f(u,z)$ with analytic $z$ slices on a
common interval and analytic $u$ slices for almost every $z$.
Lemma~\ref{III-neutral_boundary:lem:zero-cross} supplies a jointly analytic
patch at zero. Its zero section is the original recognized boundary.
To verify transversality, differentiate the same original zero-generator
covariance: the actual nonzero curve derivative is invariant under the
original derivative cocycle. The actual tangential columns have the
strictly contracting projective rates. The finite-probability recurrence
argument separating these rates puts the curve derivative outside the
tangent plane. This uses the original derivative, without an additional
identification of normalization frames. The resulting joint derivative
is invertible at a suitable point, giving an actual collar chart.

Lift this patch to an open subset of $G$ using parabolic invariance.
Minimality on the angular boundary and right $\Gamma$ equivariance
propagate it as in Sections~\ref{III-series:negative-three}
and~\ref{III-series:positive-three}. On overlaps the analytic families
agree almost everywhere with the same original curve field, hence agree
identically. Their developing map has the original finite angular cover
as its zero section. Compactness gives an open full-angular neighborhood
$W$ containing a uniform normal strip, and its actual map to a
neighborhood of the original boundary is a finite surjective analytic
cover with analytic local inverses. On this cover the normal coordinate
is a product coordinate and every lattice element fixes it, as a germ.
This is the required product collar after the stated finite choices;
no nonzero central flow is introduced.
\end{proof}

\begin{corollary}[Analytic normalization of the neutral normal quotient]
\label{III-neutral_boundary:cor:analytic-normal}
For every recognized neutral projective boundary, including the
case $B\not\equiv0$, the normal quotient line bundle has a
nowhere-zero analytic invariant section after the orientation
and finite-index choices. Thus in an analytic tubular chart
the normal derivative of every lattice element is one.
\end{corollary}
\begin{proof}
Let $\nu=TM|_N/TN$ with its analytic derivative action, and
form the compact analytic $3$-manifold
$\mathbb P(\nu\oplus\R)$. The zero section has the original
projective action, and its normal cocycle is the original
normal quotient cocycle. In the affine fiber chart the action
is exactly linear. We make explicit the measurable transfer used here.
Write $G=\SL(3,\R)$ and let $Q$ be the stabilizer of the chosen point
in the actual homogeneous boundary, using the ray stabilizer on its
oriented spherical cover. In an analytic unit frame of $\nu$, write
$a(\gamma,z)$ for the nonzero normal derivative and put
$c(\gamma,z)=\log|a(\gamma,z)|$ on the original lattice, before its
normal-sign restriction. The relative scalar normalization at $q=0$ gives
a coboundary for the induced $Q$-cocycle on $G/\Gamma$, with its same
normalized Haar law and measurable coset section $r$.
Explicitly, with $x\in G/\Gamma$ and $z$ in the native boundary,
\[
 \begin{aligned}
 k(g,x)&=r(gx)^{-1}gr(x),\\
 \beta(g,(z,x))&=c\bigl(k(g,x),r(x)^{-1}z\bigr),\\
 \beta(q,(o,x))&=b(qx)-b(x)
 \end{aligned}
\]
for each fixed $q\in Q$, almost everywhere in $x$. We retain the original
native smooth probability measure and its nonsingular $G$-action.

The passage back to the native boundary is a derived homogeneous
induction argument. The classical input is the injectivity of measurable
cohomology under induction for a closed subgroup of a locally compact
second-countable group, a nonsingular standard Borel action, and real
additive coefficients
\cite[Proposition~2.2 and generality remark, author-version p.~4]{IozziOrbit1991}.
The section and diagonal-action conventions are those of
\cite[Proposition~4.4 and remark, author-version p.~11; p.~12]{IozziGeometric1992}.
The Borel orbit section $\rho$, with $\rho(z)o=z$ on this same native
boundary, is constructed by
Arsenin--Kunugui uniformization
\cite[Theorem~5.12.1, pp.~219--220; definition, p.~185]{SrivastavaBorel}:
the orbit relation is Borel, its nonempty stabilizer fibers are closed
and $\sigma$-compact, and the original group topology is Polish.
If $\beta$ is the induced $G$-cocycle and $b$ the preceding $Q$-transfer,
the product transfer is
\[
 A(z,x)=\beta\bigl(\rho(z),(o,\rho(z)^{-1}x)\bigr)
          +b\bigl(\rho(z)^{-1}x\bigr).
\]
For each fixed $g,z$, the return element
$\rho(gz)^{-1}g\rho(z)$ lies in the same $Q$. Invariance of the
quotient Haar law and Fubini turn the given slice identity into the
product almost-everywhere coboundary identity. Induction injectivity,
applied to $c$ and the zero cocycle, therefore constructs a measurable
$f$ with $c(\gamma,z)=f(\gamma z)-f(z)$ for every $\gamma$, almost
everywhere in the original native measure. No rank restriction or
integrability upgrade is used. Multiplying the unit frame by $e^f$
makes its normal coefficient $\operatorname{sign}(a(\gamma,z))$.
It is therefore invariant on the canonical finite-index kernel of this
sign character. This restriction uses the subgroup section aligned with
the same original section $r$; the preceding transfer retained the
original suspension measure, section and gauge. Its actual vertical
curve parameters are fixed by the whole connected parabolic. In particular
the normalized central generator of this new action is zero.
Apply Theorem~\ref{III-neutral_boundary:thm:stationary} to its zero section.
The transverse derivative of the resulting analytic product
chart is a nowhere-zero analytic invariant section of its
normal line, canonically $\nu$. For an analytic tubular chart of
the original hypersurface, embed the original compact ambient manifold
analytically in Euclidean space
\cite[Theorem~3, p.~470; regularity, p.~462]{Grauert}, lift this section
orthogonally using the induced analytic metric, and apply the analytic
local inverse theorem along the zero section. Its quotient derivative
is the prescribed section, so every group element has normal derivative one.
This normalizes
the quotient; it does not assert an invariant splitting of
the full tangent bundle or an invariant ambient foliation.
\end{proof}

\section{Analytic dependence of transverse jets}\label{III-series:jets-three}

\subsection{Statement and the affine jet tower}

Let $\Gamma<\SL(3,\R)$ be any lattice
acting analytically on a compact real-analytic $3$-manifold $M$.
Let $N\subset M$ be a recognized invariant analytic hypersurface,
with its action analytically identified with the standard action
on $\RP^2$ or its spherical lift. The dual action is treated by
inverse transpose. Make the finite-index and orientation choices
needed to preserve $N$ and its chosen lift. Write
\[
                         \nu=TM|_N/TN.
\]
Assume that $\nu$ has a nowhere-zero analytic $\Gamma$-invariant
section. In the neutral boundary application this is supplied by
the projectivized-normal-bundle argument in Section~\ref{III-series:axis-three}.

An unparametrized transverse curve jet at $b\in N$ is the jet of
an embedded curve through $b$ transverse to $N$, modulo changes
of its one-dimensional source parameter. Denote its order-$k$
bundle by $\mathcal J^k$ and put $\mathcal J^0=N$.
These are finite-dimensional analytic jet bundles, with the
analytic actions induced by the original action on $M$.

\begin{theorem}[Every finite transverse curve jet is analytic]
\label{III-neutral_jet_regularity:thm:jets}
Suppose $b\mapsto\mathcal C_b$ is a measurable
$\Gamma$-equivariant field of formal unparametrized curves
transverse to $N$, defined almost everywhere for the projective
Lebesgue measure class. For every finite $k\geq1$, the section
$b\mapsto j^k\mathcal C_b$ agrees almost everywhere with a
global analytic $\Gamma$-equivariant section of $\mathcal J^k$.
These analytic sections are compatible under jet truncation.
\end{theorem}

\begin{lemma}[The highest jet is an affine variable]\label{III-neutral_jet_regularity:lem:torsor}
The projection $\mathcal J^k\to\mathcal J^{k-1}$ is an affine
bundle modeled, over a lower jet based at $b$, on
\[
                  E_{k,b}=T_bN\otimes(\nu_b^*)^{\otimes k}.
\]
In particular, the pullback over an analytic lower-jet section
is an analytic affine bundle of rank two over $N$.
The induced group action is affine on its fibers, with linear
part the natural action on $E_k$.
\end{lemma}
\begin{proof}
In an analytic ambient chart $(u,z)$ with $N=\{z=0\}$, every
transverse curve has a unique graph parametrization
$u=b+v_1z+\cdots+v_kz^k+O(z^{k+1})$.
Fix the graph through order $k-1$. Under a change of adapted
ambient chart, the difference of two possible coefficients
$v_k$ transforms by the tangential derivative and the inverse
$k$th power of the normal quotient derivative. All other terms
depend only on the already fixed lower jet. Thus the transition
rule for $v_k$ is affine, and its difference space is exactly
$T_bN\otimes(\nu_b^*)^{\otimes k}$. This also applies at
$k=1$: transverse lines are the affine space of splittings
of $TM|_N\to\nu$. The construction is intrinsic and therefore
intertwines every analytic diffeomorphism preserving $N$.
\end{proof}

If $\mathcal A\to N$ is any such affine bundle, its affine
transition functions $v\mapsto A_{ij}v+c_{ij}$ homogenize to
linear transition functions
\[
            (v,t)\longmapsto(A_{ij}v+t c_{ij},t).
\]
They define an analytic rank-three vector bundle $W$ in an
exact sequence $0\to E_k\to W\to\R\to0$. The bundle
$\mathcal A$ is the affine chart $t\ne0$ in
\[
                              Y=\mathbb P(W).
\]
Every affine bundle automorphism homogenizes compatibly.
Consequently $Y$ is a compact analytic $4$-manifold with a
genuine analytic $\Gamma$-action, preserving $\mathcal A$
and covering the given action on $N$. No choice of a global
affine origin is required.

\subsection{The homogeneous jet graph and its exponents}

We establish the regularity step needed for the induction.
Let $\mathcal A\to N$ be a $\Gamma$-equivariant analytic affine
bundle modeled on $E_k$, and let $s:N\to\mathcal A$ be a
measurable equivariant section. Regard $s$ as a section of
$Y=\mathbb P(W)$. This compactification is connected: its
projection is open and surjective, the base $N$ is connected,
and each projective fiber is connected.

Put $G=\SL(3,\R)$ and let $Q$ be the line parabolic; on an
oriented projective cover use its identity component. In the
suspension $(G\times Y)/\Gamma$, the graph
\[
                g\Gamma\longmapsto[g,s(g^{-1}[e_1])]
\]
is interpreted with the standard spherical lift when necessary.
Its pushforward of Haar probability is a $Q$-invariant measure
$\mu$. Projection to $G/\Gamma$ is a measurable equivariant
isomorphism. In particular $\mu$ has one-point fiber
conditionals. The projection of this graph to the suspension
of $N$ is the recognized boundary graph, so its connected
stabilizer is exactly the connected parabolic.

Write
\[
 \beta_2=\varepsilon_2-\varepsilon_1,\qquad
 \beta_3=\varepsilon_3-\varepsilon_1,
 \qquad
 a_t=\operatorname{diag}(e^t,e^{-t/2},e^{-t/2}).
\]
The tangent bundle of $Y$ along the measurable graph has the
invariant exact sequence
\[
          0\longrightarrow T^{\mathrm{vert}}Y
           \longrightarrow TY\longrightarrow TN\longrightarrow0.
\]
In the affine chart, $T^{\mathrm{vert}}Y$ is canonically $E_k$.
The invariant analytic frame of $\nu$ identifies its derivative
action with the action on $TN$. Therefore the fiberwise
Lyapunov exponents of $\mu$, with multiplicity, are
\begin{equation}\label{III-neutral_jet_regularity:eq:double}
                         \beta_2,\beta_2,\beta_3,\beta_3.
\end{equation}
This computation uses the invariant subbundle and its quotient,
not a differentiable splitting of the measurable section.
Changing from affine vertical norms to compactified vertical
norms is a measurable change of frame. The two cocycles have
integrable logarithmic return bounds; recurrence to sets on
which the change of frame and its inverse are bounded identifies
their exponents. This comparison also applies to the original
adapted metric at any prescribed measurable suspension lift:
its pullback in the selected chart is a measurable positive
definite form and gives a measurable isometric frame. The
coordinate derivative is the derivative of that lift's actual
lattice return. Recurrence compares the existing full Cartan
spectra, including multiplicities, without requiring integrability
of the change of frame itself. More explicitly, for each fixed Cartan element,
recurrence on both opposite rays places each invariant diagonal
column in its ambient weight space. It also forces every nonzero
quotient coordinate of an ambient weight vector to have the same
weight. The inclusion and quotient thus intertwine the spectral
projections below each threshold. Exactness makes their ranks
additive, giving the union with multiplicity
\cite[Remark~2.2.9(iii), pp.~6--7]{FilipMETv1}.

For $t>0$ all four exponents at $a_t$ equal $-3t/2<0$.
For a projective boundary, Proposition~5.7 of \cite{BRHW}
consequently makes $\mu$ $\beta_j$-tame for $j=2,3$.
For a spherical boundary, first average its two graph branches
under the finite component group of the minimal parabolic.
The resulting measure is ergodic for the full minimal
parabolic, has the same exponents, and is tame by that
proposition. To pass back to a branch, use the finite Cartan
component decomposition of this average. Its finiteness follows
from the finite diagonal sign group and the contracting root
subgroups: a Cartan-invariant set is invariant under the latter,
so full parabolic ergodicity leaves only finitely many components.
Each branch is an absolutely continuous Cartan-ergodic probability
and hence equals one normalized component. Choose the component
labels saturated by the stable leaves of a Cartan element
contracting $\beta_j$. The $\beta_j$-leaf stays in that component.
Uniqueness in \cite[Proposition~3.11]{BRHW} identifies the branch
leafwise measures there with the restricted component measures,
up to scale. Their normalized covariance and open neighborhoods
therefore give Definition~4.1 for the chosen branch. The relevant regularity threshold
$c_0^F([\beta_j])$ is one: multiplicity two does not change
the ratio of the exponents in the coarse space.

The Howe--Moore theorem \cite[Theorem~1.1 and Definition~2.6, pp.~66 and~69]{Ciobotaru} makes every nontrivial one-parameter
subgroup of the diagonal group ergodic on $G/\Gamma$.
The graph isomorphism transfers this property to $\mu$.
In particular $\ker\beta_j$ is ergodic, so the decomposition
denoted $\mathcal E_{\beta_j}$ in \cite{BRHW} is trivial.
Let $\xi$ be a measurable partition subordinate to the fiberwise
coarse foliation. For $\mu$-almost every $x$, the conditional
probability on $\xi(x)$ gives full mass to the section graph.
The atom $\xi(x)$ lies in the fiber through $x$, where the graph
meets it only at $x$. Hence this conditional probability is
$\delta_x$. The family $x\mapsto\delta_x$ is Radon on the fiber
leaves and equivariant under the Cartan action. It satisfies the
coherence and conditional-measure formulas of
\cite[Proposition~3.11, p.~25]{BRHW}, also after conditioning on
the trivial partition $\mathcal E_{\beta_j}$. The graph measure is not invariant under the
missing root subgroup: its projection to the boundary graph
has exactly the parabolic stabilizer.

Thus alternative~(2) of \cite[Proposition~7.2]{BRHW} applies.
For any fixed finite $r>1$ it gives, almost everywhere, a
$C^r$ embedded root curve
\begin{equation}\label{III-neutral_jet_regularity:eq:rootgraph}
                \varphi_{j,x}:U^{\beta_j}\simeq\R
                              \longrightarrow W^{\beta_j,F}(x)
\end{equation}
whose graph carries the unconditional missing-root leafwise
measure. The target coarse stable manifold has dimension two.
The assertion that its measure is the \emph{unconditional}
graph uses ergodicity of $\ker\beta_j$; no dimension-dependent
identification of conditional measures is being invoked.

To read this curve in the prescribed suspension section, use the
fiber embedding of \cite[Proposition~9.12, pp.~84--85]{BRHW}
and the total embedding from \cite[\S11.2, p.~96]{BRHW},
whose defining formula is
$\iota_x(\exp(Y)z)=(Y,\iota_x^F(z))$.
These embeddings have no one-dimensional fiber restriction.
The coarse class is strongly integrable; its ratio at a contracting
Cartan element is one. Choose a finite auxiliary order above the
global fiber critical regularity as well. The product embedding
then makes the root graph continuous in the intrinsic full leaf.
A subordinate partition identifies its local measure class with
the ordinary conditional measures. Recurrence under the contracting
Cartan element extends this local domination to the whole root
parameter line. Consequently every suspension-conull set contains
almost every point of almost every such root graph. Applying this
to the given section graph identifies the physical curve with the
prescribed section along almost every root parameter; its center
is the prescribed section value. This passage uses conditional
measures and recurrence, not an extra assertion in Proposition~7.2.

\begin{remark}[The dimension restrictions that are not used]
Corollary~4.5 and Theorem~4.6 of \cite{BRHW} require
$\dim Y=\dim(G/Q)$ and do not apply here. Theorem~8.1
also cannot be invoked with the central contraction: its
hypothesis~(a) asks for a diffeomorphism from the two-dimensional
opposite unipotent group onto the full stable fiber, which is
$4$-dimensional in \eqref{III-neutral_jet_regularity:eq:double}. We use only the
root embeddings from Proposition~7.2 and then regularize the
given section locally. Proposition~7.5 is not needed.
\end{remark}

\subsection{Analyticity of the root curves with repeated weights}

The multiplicity in \eqref{III-neutral_jet_regularity:eq:double} prevents using a
one-dimensional stable-manifold argument verbatim. It does
not prevent analytic linearization, because both exponents
on each stable plane coincide.

\begin{lemma}[Equal-exponent analytic linearization]
\label{III-neutral_jet_regularity:lem:linearize}
Consider the analytic return cocycle of a compact analytic
action, over an invertible probability-preserving transformation
$T$, with the integrable analytic return bounds used in
Section~\ref{III-series:coordinate-curves}, Section~\ref{III-series:relative-three}. Suppose its stable space has
dimension $m$ and all of its stable exponents equal one
number $\chi<0$; all other exponents are positive.
Its local stable manifolds have measurable analytic charts in
which the return maps are exactly their derivative cocycle
$A_x:E_x^s\to E_{Tx}^s$. The charts and their inverses have
positive tempered radii and tempered analytic bounds.
\end{lemma}
\begin{proof}
The analytic bounds are integrable in the original fixed return charts.
For the derivative cocycle in those same charts, both positive logarithmic
operator norms, including that of the inverse, are integrable. Apply
\cite[Theorem~S.2.10, pp.~666--667; proof pp.~667--668]{KatokMendoza1995}
on the completed original probability space. This nonergodic theorem
allows the number and dimensions of the spectral blocks to vary with the
base point. For each $\epsilon>0$ it gives a measurable full frame $C_x$
carrying the standard consecutive Euclidean blocks onto the actual
invariant Lyapunov subspaces, and a block diagonal conjugate
$C_{Tx}^{-1}D_0F_xC_x$ whose $i$th block satisfies
\[
 e^{\lambda_i(x)-\epsilon}\le
 \|(B_x^i)^{-1}\|^{-1}\le\|B_x^i\|
 \le e^{\lambda_i(x)+\epsilon}.
\]
The subspaces and exponents are those of the original cocycle in
\cite[Theorem~S.2.9, p.~665]{KatokMendoza1995}, not a separately selected
splitting. The growth filtration of a direct sum of blocks with distinct
vector growth rates is the sum of the blocks below each cut. Comparing
with the original Lyapunov filtration identifies the negative sum with
the given stable space. If $m>0$, the sole negative exponent is $\chi$,
so this sum is the first block and has dimension $m$; the remaining
blocks form its complement. If $m=0$, use the zero stable space and its
identity chart. Use the fixed Euclidean coordinate splitting to write this
same full frame as $K_x:E^s\times E^u\to E$. The product norm introduces
only fixed constants in the full-frame bounds and preserves the Euclidean
norm on the stable coordinate block.

The tempering definition in the cited source is forward in time, for
both $C_x$ and $C_x^{-1}$ \cite[Definition~S.2.2, p.~661]{KatokMendoza1995}.
Apply the finite-measure gauge argument to
$1+\|C_x\|+\|C_x^{-1}\|$: its forward subexponential bound yields
bilateral logarithmic sublinearity on the same probability space.
Thus $K_x$ and its inverse have the bilateral tempered bounds used below.
Choose an original-measurable version of this completed-space frame.
It equals the same frame almost everywhere, as do its values at $x$
and $Tx$; hence the actual derivative block, stable space and orbit
bounds are preserved.
Complexify these charts and rescale by their positive integrable
logarithmic radii. Apply the original stable-set theorem
\cite[Theorem~5.1(a)--(c), Case~III, pp.~44--45; proof pp.~46--54]{Ruelle1980Author}
to this return family and its original probability. Choose
$\chi<\lambda_{\rm fast}<\lambda_{\rm cut}<0$; the other exponents
are positive, so the stable cut is precisely the given stable space.
Projection onto this space gives a locally unique normalized graph.
Its uniqueness gives measurable dependence and compatibility with
complex conjugation; restriction to the real space gives analytic
charts tangent to the identity. Cauchy estimates bound their
coefficients and inverse coefficients by the original window radii.
The forward windows in part~(c), together with invariance of the
probability, give bilateral tempered bounds. This is the stable
construction used in Section~\ref{III-series:coordinate-curves}, with
stable dimension $m$. Transport by the given split frame preserves
these bounds and the original derivative cocycle.
To identify the actual stable set after this transport, write $r_x$ for
the original radius, $K_x$ for the given split frame and
$A_x^{\rm coord}=r_x^{-1}K_x$. The normalized return and the split
canonical return agree under these coordinates on a common ball of
radius $\delta_x>0$. One may take half the minimum of the two
canonical convergence radii, with the source radius divided by
$1+\|K_x\|$, and of $r_x/(1+\|K_x\|)$. The inverse radius is
therefore bounded by finite products and maxima of the raw and split
analytic heights and the frame norm. These quantities are tempered.
For each $\epsilon>0$ they give, on each regular orbit, constants
$M,d>0$ such that
\[
 \|A_{T^nx}^{\rm coord}\|,\ \|(A_{T^nx}^{\rm coord})^{-1}\|
 \le M e^{\epsilon n},\qquad \delta_{T^nx}\ge d e^{-\epsilon n}.
\]
Choose the target rate $s$ with
$\lambda_{\rm fast}+\epsilon<s<\lambda_{\rm cut}-\epsilon<0$.
The source pair-contraction estimate gives one direction of the
stable-set comparison. Conversely, an orbit with sufficiently small
bound $B e^{sn}$ stays inside every comparison ball and maps into the
original stable tube. Induction gives the actual orbit identities in
both directions. Near zero the inverse graph chart has zero transverse
coordinate exactly on that tube; its real chart image is consequently
the local set defined by the split-orbit decay bound. For the root
returns used here, $\chi=-1$, one may choose
$\lambda_{\rm fast}=-3/4$, $\lambda_{\rm cut}=-1/4$,
$\epsilon=1/10$ and $s=(-1+1/1000)/2$.

In these charts write the stable map as
$f_x(y)=A_xy+O(\|y\|^2)$ and put
$L_d(x)=A_{T^{d-1}x}\cdots A_x$. For every sufficiently
small $\epsilon$ the Lyapunov estimates give
\[
 \|L_d(x)^{-1}\|\le C_\epsilon(x)e^{(-\chi+\epsilon)d},
 \qquad
 \|f_x^{(d)}(y)\|\le C_\epsilon(x)e^{(\chi+\epsilon)d}\|y\|
\]
on a smaller complex stable chart. The quadratic coefficients
at $T^nx$ have subexponential bounds. Accordingly the
telescoping differences of
\[
                         L_d(x)^{-1}f_x^{(d)}(y)
\]
are bounded by $C(x)e^{(\chi+c\epsilon)d}\|y\|^2$, for a
fixed constant $c$. Choose $\epsilon$ with $\chi+c\epsilon<0$.
The sequence converges normally to an analytic chart $K_x$
with $DK_x(0)=\id$ and
\[
                         K_{Tx}f_x=A_xK_x.
\]
More precisely, apply the original holomorphic theorem
\cite[Theorem~1.1, pp.~2138--2139]{BDM2008} with $q=1$ and
linear target $N_x=A_x$. In the Lyapunov frame take
$\epsilon=-\chi/100$, $m=e^{\chi-\epsilon}$ and
$M=e^{\chi+\epsilon}$; then $0<m\le M<1$ and $M^2<m$.
The real Hilbert bounds pass to the complexified operators with the same
constants. Bilaterally subexponential intrinsic coefficient heights give
the slow source radii and fast image bounds required by the source's
tameness definition; the bounded linear target is tame as well.
The common first jet gives the required order-two tangency.
All almost-everywhere bounds hold everywhere on a measurable,
bidirectionally invariant conull set. Extending the family by a fixed
contracting germ on its complement permits application over the original
base; the resulting conjugacy agrees with the original family almost
everywhere. A zero-dimensional fiber requires only the identity chart.

The real and measurable conclusions are derived from this source, whose
statement concerns holomorphic maps over an arbitrary bijection.
Complex conjugation commutes with the actual normalized finite iterates,
so their limit is real on real inputs. The finite iterates are jointly
measurable; the pointwise limit and its Taylor coefficients give a
measurable real germ. Near-identity injectivity identifies its real
analytic inverse. Cauchy estimates on a smaller slow ball, together with
the coefficient bounds for the inverse, bound the intrinsic height of
this same invertible germ by $C/r(x)^2$. This initially gives a finite
exponential orbit bound. Applied to the measurable height, the
finite-measure bounded-coboundary argument used for tempered gauges
improves it to bilateral subexponential growth; the canonical
representative construction therefore supplies the stated tempered
forward and inverse domains. Finally pull back by the original
Lyapunov frame. Inverting the constructed normalizer when necessary
matches the chart convention without changing its intrinsic height.
\end{proof}

\begin{lemma}[A smooth root parametrization becomes analytic]
\label{III-neutral_jet_regularity:lem:analytic-root}
The root maps \eqref{III-neutral_jet_regularity:eq:rootgraph} are analytic near their
origins, and hence on their entire root intervals.
\end{lemma}
\begin{proof}
Fix $j$, and choose a diagonal element $a$ with
$\chi=\beta_j(a)<0<\beta_{j'}(a)$, where $\{j,j'\}=\{2,3\}$.
By \eqref{III-neutral_jet_regularity:eq:double} its stable fiber is the two-dimensional
coarse space for $\beta_j$, with both exponents equal to
$\chi$. We may choose $a$ with $\beta_j(a)=-1$ and $\beta_{j'}(a)=1$.
The original raywise splitting
\cite[nonergodic discussion and equation~(7), p.~463]{BRHWSmoothErgodic2023}
and the measurable frame construction above identify this stable plane
with the negative space of the actual derivative cocycle.
Keep this same frame and probability when applying
Lemma~\ref{III-neutral_jet_regularity:lem:linearize}; its Ruelle and
Berteloot--Dupont--Molino constructions apply to these original returns.
The original action on $Y$ is analytic on a compact manifold;
the finite-order jet construction and homogenization preserve
the integrable analytic return bounds. Thus that lemma applies
to this genuine compact action just as to the original one.

Take the embeddings from Proposition~7.2 with regularity
$C^2$ (their regularity may be chosen to be any fixed finite order), and write
$h_x(t)=K_x\varphi_{j,x}(t)$ in the linearizing chart.
Uniqueness of the Haar graph and Cartan conjugation of the
root subgroup give
\[
              h_{Tx}(e^\chi t)=A_xh_x(t),\qquad T x=ax.
\]
Set $J_x=h_x'(0)$ and $R_x(t)=h_x(t)-J_xt$.
For each $x$ there are positive measurable $\delta_x,C_x$
such that $\|R_x(t)\|\le C_x|t|^{1+\alpha}$ for
$|t|<\delta_x$. These particular bounds need not be tempered.
For return times between a positive-measure set with uniform endpoint
chart domains and remainder bounds, the iterated linearizing identity
first holds as a germ. Analytic uniqueness extends it across the common
connected endpoint ball, using the globally analytic return map.
The original root return and endpoint chart injectivity therefore give,
for uniformly small $t$,
\[
                 R_x(t)=L_n(x)^{-1}R_{T^nx}(e^{n\chi}t).
\]
Choose a positive-measure set on which $\delta_x$ is bounded
below and $C_x$ is bounded above. The transformation $T$ is
ergodic, so almost every orbit returns to this set infinitely
often. At those return times the right hand side is bounded by
\[
          C'(x)e^{(\chi+\epsilon)n}|t|^2.
\]
Choose $\epsilon<-\chi$. The bound tends to zero;
therefore $R_x(t)=0$. The equality first holds for all small
rational $t$ on one conull set and then for all small $t$
by continuity. Thus $h_x(t)=J_xt$, and
$\varphi_{j,x}(t)=K_x^{-1}(J_xt)$ is analytic near zero.
Its analytic root-parameter radius has a positive measurable
lower bound on some positive-measure set. Recur to that set.
For any fixed bounded root interval, a sufficiently late
return contracts the entire interval into that radius.
The identity
\[
 \varphi_{j,x}(t)
   =F(a^n,x)^{-1}\varphi_{j,T^nx}(e^{n\chi}t)
\]
then expresses it analytically on that interval; here
$F(a^n,x)$ acts by a global analytic diffeomorphism between
the full compact suspension fibers. Exhaustion by bounded
intervals proves the asserted analyticity on the root line.
\end{proof}

\subsection{Angular regularity and induction}

\begin{proposition}[Analyticity of a measurable affine section]
\label{III-neutral_jet_regularity:prop:section}
Every measurable $\Gamma$-equivariant section of an analytic
affine bundle modeled on $TN\otimes\nu^{-k}$ agrees almost
everywhere with a global analytic equivariant section.
\end{proposition}
\begin{proof}
Apply the preceding construction to its projective
compactification. The two missing root subgroups commute
and give the coordinate foliations of a projective open cell.
The unconditional Haar graph identities and Fubini identify
almost every coordinate slice of the section with the
analytic root parametrization from
Lemma~\ref{III-neutral_jet_regularity:lem:analytic-root}. These identifications are made
first in open suspension charts. Haar measure has the same local null
sets as coordinate Lebesgue measure, as stated explicitly after (8.25)
in \cite[p.~532; coordinate-null convention and (8.11), p.~527]{KnappBeyondDigital}.
The supporting smooth invariant form and Haar uniqueness are
Theorems~8.21 and~8.23 there. Replacing coordinate Lebesgue measure
by any additive Haar measure preserves this null class, since these
measures differ by a positive constant.
Thus Fubini applies on these genuine open group charts; the identifications
do not restrict an arbitrary almost-everywhere representative to a null section
of $G$.

The almost-everywhere separate-analyticity lemma of
Section~\ref{III-series:coordinate-curves}, based on the classical cross theorem
\cite[Theorem~4, p.~298]{Cross}, gives a jointly analytic representative
on one nonempty angular rectangle. That lemma is independent
of the dimension of its compact analytic target, so it applies
to $Y$. Its projection to $N$ equals the identity almost
everywhere and therefore everywhere. Choose a point where its
value lies in the open affine chart $\mathcal A$, as it does
almost everywhere, and shrink the rectangle. We now have an
analytic section of $\mathcal A$ on an open patch.

Translate the patch by the lattice. The projective action is
minimal; hence the translates cover $N$. Equivariance modulo
null sets shows that the translated analytic sections agree
almost everywhere on overlaps, and analytic uniqueness makes
them agree identically. They define a global analytic section
of $\mathcal A$, equal to the original section almost
everywhere. The same uniqueness makes its equivariance exact.
\end{proof}

\begin{proof}[Proof of Theorem~\ref{III-neutral_jet_regularity:thm:jets}]
Start with the order-zero section, the identity of $N$.
Suppose the order-$(k-1)$ section has been made analytic.
Its pullback of $\mathcal J^k$ is the analytic affine bundle
of Lemma~\ref{III-neutral_jet_regularity:lem:torsor}, and the order-$k$ curve jet is a
measurable equivariant section of this bundle.
The linear transition and action of this same pullback bundle are the
original tangent derivative tensored with the inverse $k$th power of the
original normal quotient action. Proposition~\ref{III-neutral_jet_regularity:prop:section}
therefore makes its measurable highest-jet section analytic. Lift that
section over the chosen preceding analytic section. By construction its
truncation is exactly that section at every point, and its order-$k$ jet
agrees almost everywhere with the original curve jet. A single recursive
choice of these lifts gives the compatible analytic tower. Induction
proves the theorem.
\end{proof}

\begin{corollary}[Application to the nonstationary neutral boundary]
The normalized unparametrized transverse curve field of
Section~\ref{III-series:axis-three}, Section~\ref{III-series:fixed-curves} has analytic angular dependence at every
finite jet order, including when its analytic central generator
satisfies $B\not\equiv0$ and $B'(0)=0$.
\end{corollary}
\begin{proof}
The image curves are parabolic-equivariant even when their
chosen parameters change by the central flow. Their
unparametrized jets therefore descend to measurable
equivariant sections on the recognized projective boundary.
The analytic invariant normal quotient frame is supplied by
Section~\ref{III-series:axis-three}. Apply Theorem~\ref{III-neutral_jet_regularity:thm:jets}.
\end{proof}

\section{Normalization of the transverse parameter}\label{III-series:parameter-jets}

\subsection{An additive transfer over the boundary}

Let $\Gamma<\SL(3,\R)$ be any lattice,
and let $N$ be its standard projective boundary or the oriented
sphere cover. Fix the finite-index and covering choices needed
to preserve that component. All almost-everywhere statements on
$N$ refer to its smooth measure class.

\begin{lemma}[Analytic regularity of an additive transfer]
\label{III-neutral_parameter_jets:lem:transfer}
Let $n:\Gamma\times N\to\R$ be an additive cocycle whose
function $n(\gamma,\cdot)$ is analytic for every $\gamma$.
Suppose a finite measurable function $a:N\to\R$ satisfies
\[
 a(\gamma p)-a(p)=n(\gamma,p)
\]
for every fixed $\gamma$, almost everywhere. Then $a$ agrees
almost everywhere with a real analytic function on $N$.
\end{lemma}
\begin{proof}
Consider the analytic action on the compact $3$-manifold
$N\times\RP^1$ given in the affine fiber chart by
\[
 \gamma(p,s)=(\gamma p,s+n(\gamma,p)).
\]
In homogeneous fiber coordinates its matrix is
$\left(\begin{smallmatrix}1&n(\gamma,p)\\0&1\end{smallmatrix}\right)$,
so it extends analytically across infinity. The cocycle identity
is exactly the action identity.

The measurable equivariant graph $p\mapsto(p,a(p))$ defines
an ergodic, Haar-projecting, parabolic-invariant suspension
probability with one-point fiber conditionals. Its full connected
stabilizer is the projective parabolic: projection to the first
factor sends it to the standard boundary graph measure, whose
stabilizer is that parabolic.

The vertical derivative in the affine coordinate is one.
Its Lyapunov exponent is therefore zero. This conclusion does
not require a moment bound on $a$: the compact action has the
usual integrable return derivative bounds. More precisely, choose one
almost-minimal Borel lattice section as follows. The real-field
specialization of \cite[Theorem~1.4.15, pp.~54--55]{BHV} gives property
$(T)$ for $\SL(3,\R)$. The original lattice has the normalized invariant
regular Borel quotient probability. It inherits property $(T)$ by
\cite[Theorem~1.7.1, pp.~65--66]{BHV}, and the discrete finite-generation
clause \cite[Theorem~1.3.1, pp.~41--42]{BHV} supplies finite generators.
For the ambient word model choose a compact identity neighborhood; its
generated subgroup is open and closed and hence is the whole connected
group. The compact word balls make its word length measurable. These
choices are internal constructions from the named original inputs, and are
made once for this original lattice. For a locally compact group with property $(T)$,
\cite[Theorem~5.3, pp.~184--185]{dLS} gives an exponentially
integrable Borel fundamental domain. The direct compact-word lattice
undistortion statement \cite[Theorem~5.1, pp.~183--184]{dLS}
then gives the return-word moments; its Riemannian-metric predecessor is
\cite[Theorem~A, p.~6]{LMR2000}. Together with the two generator-selection inputs above, these are five
separate universal classical statements, not an assumed conclusion for this graph.
The lattice is the original discrete subgroup with its given finite nonzero
Haar covolume. The reindexing from $\{1\}\sqcup\{1,2\}$ to $\{1,2,3\}$
transports that same subgroup, generating families, word lengths and section.
Haar scaling cancels in the normalized quotient probability; no second
lattice or independently chosen probability is substituted.
Use the same almost-minimal section and the return estimates in
\cite[pp.~185--186]{dLS}.
Recurrence to sets on which the nonzero affine vertical frame and its inverse
are bounded identifies its exponent with zero. The quotient
derivative is the standard tangent derivative on $N$.
Consequently the three fiber exponents are
$\beta_2,\beta_3,0$.

Use the original graph directly. The
upper-root Mautner argument gives ergodicity of every positive central
return power on its suspension. The analytic stable-plaque construction
with the preceding return-word moments gives a local analytic affine
height over the full column-root parameters. Its asymptotic calibration
identifies the plaque with the same original graph for Lebesgue-almost
every root parameter: the boundary coordinates coincide and the vertical
coordinate is uniquely determined by the calibrated stable asymptotics.
An analytic change from these parameters to the actual boundary chart
therefore gives a nonempty open set on which an analytic scalar function
agrees almost everywhere with $a$. Pulling back the original graph's null
set through this chart uses the actual angular measure-class dictionary.

Translate this local scalar seed using the original cocycle equation.
The action is analytic and nonsingular, so each translate is again an
analytic representative of $a$. On an overlap two representatives agree
almost everywhere and hence everywhere by continuity and full support
of the smooth measure class. Minimality of the original lattice action
on the projective boundary or sphere cover makes the translated domains
cover $N$. The representatives patch to an analytic real-valued function
$\bar a$, equal to $a$ almost everywhere. For every fixed $\gamma$,
the continuous function
$\bar a(\gamma p)-\bar a(p)-n(\gamma,p)$ vanishes almost everywhere
and therefore everywhere. An almost-everywhere measurable $a$ is handled
by measurable replacement and nonsingularity. This direct construction
uses no separately supplied recognized boundary or analytic transfer.
\end{proof}

\subsection{The normalized finite parameter coefficients}

Let $N\subset M$ be a recognized neutral boundary in a compact
analytic $3$-manifold. Work on the two-sided cover. Use the
analytic invariant normal quotient section in Section~\ref{III-series:axis-three}
to choose an analytic tubular coordinate $(p,t)$, with $N$
given by $t=0$ and normal derivative one for every lattice
element. The normalized formal curve field in Section~\ref{III-series:axis-three}
has deterministic analytic central generator $B(z)\partial_z$,
where $B(0)=B'(0)=0$.

Choose the standard analytic metric splitting of the positive
projective character. It expresses the descended curve field
by formal maps $J_p(z)$ with covariance
\begin{equation}\label{III-neutral_parameter_jets:eq:covariance}
 \alpha(\gamma)J_p(z)
   =J_{\gamma p}\bigl(\Phi^B_{c(\gamma,p)}(z)\bigr),
\end{equation}
where $c$ is the analytic projective character cocycle in the
chosen central units. For the usual oriented unit-vector
description it is the corresponding constant multiple of
$\log\|\gamma\xi\|$, with $\|\xi\|=1$ and $[\xi]=p$.
This expression is unchanged under $\xi\mapsto-\xi$.
The formal descent in \eqref{III-neutral_parameter_jets:eq:covariance} uses parabolic
equivariance before choosing a section, so it does not restrict
an arbitrary almost-everywhere field to a null transversal.

\begin{proposition}[From image jets to parameter jets]
\label{III-neutral_parameter_jets:prop:parameters}
Assume that at every order the unparametrized image jets of
$J_p$ are analytic sections along $N$. Then every finite jet
of the normalized parametrized family $J_p(z)$ is analytic
in $p$, after a constant linear normalization of $z$.
The covariance identities \eqref{III-neutral_parameter_jets:eq:covariance} then hold
coefficientwise at every boundary point.
\end{proposition}
\begin{proof}
Express the image curve as a formal graph
\[
                    t\longmapsto (U_p(t),t),\qquad U_p(0)=p.
\]
The hypothesis says that each finite jet of $U_p$ is analytic
in $p$, in local analytic coordinates on $N$. Write the
normal component of the original parameter as
\[
             A_p(z)=a_1(p)z+\sum_{d\ge2}a_d(p)z^d.
\]
The normal quotient derivative of the normalized parameter is
a measurable invariant nonzero section. Its ratio to the
analytic invariant normal quotient frame is a measurable
$\Gamma$-invariant scalar. Ergodicity of the boundary action
makes this ratio constant. Choose its orientation and normalize
$z$ by this constant, so that $a_1=1$.
More explicitly, if $z=\lambda w$, the same formal flow is
conjugated by $C_\lambda(w)=\lambda w$, and its generator becomes
$B_{\mathrm{new}}(w)=\lambda^{-1}B(\lambda w)$.
The original clock $c(\gamma,p)$ and the image jets are unchanged.
In what follows $B$ denotes this transported generator.

Let $R_\gamma(p,t)$ denote the normal component of
$\alpha(\gamma)(U_p(t),t)$, expressed in the tubular coordinate
at the target boundary point. Every coefficient of this formal
series is analytic in $p$, because it uses only finitely many
coefficients of $U$ and of the actual analytic action.
Moreover $R_\gamma(p,t)=t+O(t^2)$. The image-curve equivariance
and \eqref{III-neutral_parameter_jets:eq:covariance} give
\begin{equation}\label{III-neutral_parameter_jets:eq:scalar}
 R_\gamma(p,A_p(z))
       =A_{\gamma p}\bigl(\Phi^B_{c(\gamma,p)}(z)\bigr).
\end{equation}

Suppose $a_2,\ldots,a_{d-1}$ have analytic representatives.
Since both scalar maps in \eqref{III-neutral_parameter_jets:eq:scalar} have linear
coefficient one, comparison in degree $d$ gives
\[
                  a_d(\gamma p)-a_d(p)=n_d(\gamma,p),
\]
where $n_d(\gamma,p)$ is analytic in $p$. Only the lower
coefficients of $A$, the finite coefficients of $R_\gamma$,
and the analytic function $c(\gamma,p)$ enter $n_d$.
Its cocycle identity follows almost everywhere from this
display and hence everywhere by analytic uniqueness.
Lemma~\ref{III-neutral_parameter_jets:lem:transfer} makes $a_d$ analytic. Induction proves
the assertion for every $d$. Finally, each coefficient in
\eqref{III-neutral_parameter_jets:eq:covariance} is now analytic in $p$. An identity
holding almost everywhere between such coefficients holds
everywhere.
\end{proof}

\section{Convergence of the neutral axis}\label{III-series:actual-curves-three}

\begin{corollary}[Actual neutral transverse curves]\label{III-transverse_curve_convergence:cor:neutral}
In the recognized neutral boundary setting in Section~\ref{III-series:axis-three},
there are measurable actual analytic embedded curve germs $i_x$
with
\[
 F(h,x)i_x(z)=i_{hx}(z)\quad(h\in H=L\ltimes U),\qquad
 F(a_t,x)i_x(z)=i_{a_tx}(\Phi_t^{B\partial_z}(z)),
\]
where $B$ is one analytic germ with $B(0)=B'(0)=0$.
\end{corollary}
\begin{proof}
Use the formal parameter supplied in Section~\ref{III-series:axis-three} and the
derivative frame $\rho(h)=\operatorname{diag}(l,1)$.
Pull this same field back to the original Haar quotient through its
specified-section graph. The original Haar action supplies the spectral
gap in the fixed-source convergence argument of
Section~\ref{III-series:fixed-curves}; the standard block contains a positive
hyperbolic element. Apply Corollary~\ref{III-transverse_curve_convergence:thm:curve}
there and transport its single constant source change back through the
same measure-preserving graph. This gives convergence for the original
graph probability.
Constant formal reparametrization preserves the original formal
identity for every fixed $h\in H$, so these identities are now
identities of analytic curves, almost everywhere for each $h$.

Full formal parabolic coherence of the image curves implies
actual coherence: compare analytic graphs using their invertible
third coordinate. The induced central parameter changes are
therefore analytic. Since $H$ is normal and acts trivially on
the parameters, their coefficients are $H$-invariant and hence
deterministic. They form a coefficient-measurable one-parameter
homomorphism of analytic germs. The one-parameter germ lemma
in Section~\ref{II-sec:compact} supplies an analytic generator $B$.
The normalized neutral derivative makes its linear term zero.
\end{proof}

\section{The full-flag alternative}\label{III-series:flag}

For ambient actions with orbit dimension at most three,
the full flag manifold is one of the exceptional closed
orbits in the corrected low-dimensional orbit classification
\cite[Theorem 2.1(4) and Remark 2.2, p.~5]{FisherMelnick}, including
the covering-group issue. Cairns--Ghys discuss the flag orbit in
the proof of \cite[Theorem 3.5, pp.~143--145]{CG}; we use the
explicit classification and covering qualifications of Fisher--Melnick. The point here is to recognize
this model from the suspension measure of a lattice action.

Let
\[
 \mathcal F=\Fl_{1,2}(\R^3)=G/P
\]
be the compact $3$-dimensional full-flag manifold, with $P$
the minimal parabolic.

\begin{lemma}[Minimality on the full-flag manifold]
\label{III-arbitrary_lattices:lem:flag-minimal}
The action of every lattice $\Gamma<G$ on $\mathcal F$ is minimal.
\end{lemma}
\begin{proof}
Let $F\in\mathcal F$ and let $O\subset\mathcal F$ be nonempty
and open. Put $a_t=\operatorname{diag}(e^t,1,e^{-t})$, and let
$F_+$ be its attracting flag. The finite invariant quotient
probability has full support. The Mautner argument makes $a_1$
ergodic: a vector fixed by $a_1$ in its strongly continuous
Koopman representation is fixed by every elementary root subgroup,
since conjugation by positive or negative powers of $a_1$
contracts each root. These subgroups generate $G$, and the
transitive $G$-action on $G/\Gamma$ is ergodic.

Choose a nonempty open set $N\subset G$ such that $gF$, for
$g\in N$, has lower-unipotent coordinates of norm at most one
in the attracting big cell. Choose a relatively compact nonempty
open set $D\subset G$ such that $k^{-1}F_+\in O$ for every
$k\in\overline D$. Both quotient patches have positive measure.
Ergodic recurrence gives arbitrarily large integers $n$ and
$g_n\in N$, $k_n\in D$, $\gamma_n\in\Gamma$ with
\[
                    a_ng_n=k_n\gamma_n.
\]
The three lower-unipotent coordinates contract by
$e^{-n},e^{-n},e^{-2n}$. Thus $a_ng_nF\to F_+$ uniformly
for these bounded source coordinates. Compactness of
$\overline D$ implies that
$\gamma_nF=k_n^{-1}a_ng_nF\in O$ for all sufficiently large
return times. Every orbit is therefore dense.
\end{proof}

\begin{lemma}[A sufficiently pinched full-flag return]
\label{III-arbitrary_lattices:lem:flag-return}
Every lattice $\Gamma<G$ contains an element whose action
near an attracting full flag is analytically conjugate to a
linear contraction $T$ satisfying
\[
                     \|T^{-1}\|\|T\|^3<1.
\]
The same inequality holds for every positive power of $T$.
\end{lemma}
\begin{proof}
Use the quotient recurrence argument in the proof of
Lemma~\ref{III-arbitrary_lattices:lem:flag-minimal}, with both return
patches sufficiently close to the identity. It supplies arbitrarily
large integers $n$ and $k,\gamma,w\in G$, with $\gamma\in\Gamma$,
such that
\[
 k\gamma k^{-1}=\operatorname{diag}(e^n,1,e^{-n})w,
 \qquad \|w-I\|_{\max}<1/100.
\]
For $n\ge64$, the trace and the second characteristic coefficient
of this matrix lie between $e^n/2$ and $2e^n$, and its determinant
is one. The resulting cubic has three distinct positive roots:
elementary sign estimates, with $q=e^{n/8}\ge8$, give logarithms
$h_1>h_2>h_3$ satisfying
\[
 h_1+h_2+h_3=0,\qquad
 a=h_1-h_2>0,\quad b=h_2-h_3>0,\quad
 \max(a,b)<2\min(a,b).
\]
A real eigenbasis, rescaled to determinant one, therefore
conjugates $\gamma$ to this positive split diagonal.
In opposite-unipotent coordinates
\[
 \begin{pmatrix}1&0&0\\x&1&0\\z&y&1\end{pmatrix}P
\]
the action is exactly
\[
 T(x,y,z)=(e^{-a}x,e^{-b}y,e^{-(a+b)}z).
\]
In the Euclidean norm,
\[
 \|T^{-1}\|\|T\|^3
       =e^{a+b-3\min(a,b)}<1.
\]
Taking a positive power multiplies the exponent by that power.
\end{proof}

\begin{proposition}[Analytic upgrade of a flag covering]
\label{III-arbitrary_lattices:prop:flag}
Suppose an analytic $\Gamma$ action on a compact analytic
$3$-manifold admits an equivariant $C^3$ finite covering
$h:M\to\mathcal F$ that is a local diffeomorphism. Then
$h$ and its local inverses are analytic.
\end{proposition}
\begin{proof}
Choose the return of Lemma~\ref{III-arbitrary_lattices:lem:flag-return} and take a
positive power fixing every point of the finite fiber over
its attracting flag. At one such point let $f$ be the inverse
covering germ. It satisfies $Ff=fT$, where $F$ is the analytic
return germ on $M$. The derivative of $f$ conjugates $DF(0)$
to $T$. Choose compatible norms using this derivative. Then
\[
                 \|DF(0)^{-1}\|\|T\|^3<1.
\]
Lemma~\ref{II-general:lem:analytic-parametrization}, with $r=3$, makes $f$ analytic. Apply it at
every point of the finite fiber, and shrink the target
neighborhood so its whole preimage consists of these analytic
charts. Lemma~\ref{III-arbitrary_lattices:lem:flag-minimal} and equivariance propagate
analyticity to the full covering. The analytic inverse function
theorem gives analyticity of $h$ itself.
\end{proof}

\begin{corollary}[The minimal-parabolic description branch]
\label{III-arbitrary_lattices:cor:flag}
The full-flag branch of the dimension-$3$ parabolic
description argument is an analytic finite-cover alternative
for every lattice $\Gamma<G$.
\end{corollary}
\begin{proof}
For a Haar-projecting ergodic parabolic suspension measure
whose full stabilizer is $P$, nonresonance on the original law, as derived
above from \cite[Proposition~5.1 and Claim~5.2, pp.~967--969]{BRHWFactors},
requires the three
negative roots outside $\mathfrak p$ as fiber exponents,
up to positive multiples. Dimension $3$ leaves no other
fiber directions. Thus every nontrivial coarse fiberwise class is
one-dimensional and $c^F_0(\mu)=1$
\cite[Definition~4.3, pp.~27--28]{BRHW}. On
$\operatorname{diag}(1,0,-1)$ the missing roots have values
$-1,-1,-2$, so $r_0(\mathfrak p)\le2$
\cite[(4), p.~10]{BRHW}. Apply Theorem~4.6 of \cite[p.~29]{BRHW}
at regularity $3$, available from analyticity. Its proof in
\cite[\S8.6.1, p.~76]{BRHW}, using the $C^r$ diffeomorphisms
of Theorem~8.1(a) and the quotient charts of \S8.2, gives a
genuine $C^3$ local diffeomorphism covering.
To respect the simply connected ambient-group convention in
\cite[\S2.1, p.~9]{BRHW}, let $p:\widetilde G\to G$ be the
universal covering homomorphism and put
$\widetilde\Gamma=p^{-1}(\Gamma)$ and $\widetilde P=p^{-1}(P)$.
The action of $\widetilde\Gamma$ on the same $M$ factors through
$\Gamma$. The map $[\widetilde g,x]\mapsto[p(\widetilde g),x]$
identifies its suspension with the original suspension. Pulling back
$\mu$ preserves its Haar base marginal and parabolic ergodicity, and
its stabilizer is exactly $\widetilde P$. Compatible sections identify
the derivative cocycles; comparison in continuous fiber metrics preserves
the full Lyapunov multiset. Thus the contracting direction and the
uniform coarse bound $c^F_0=1$ hold upstairs on this same law.
Apply the original theorem there. The canonical map
$\widetilde G/\widetilde P\to G/P=\mathcal F$ reads each quotient
chart in the same projected matrix Gauss coordinates, so it and its
inverse are smooth. Composing the resulting covering with this map
preserves its $C^3$ local diffeomorphism property. Surjectivity of
$\widetilde\Gamma\to\Gamma$ gives equivariance for the original lattice. Proposition~\ref{III-arbitrary_lattices:prop:flag} supplies the
analytic upgrade of this same covering map. This uses only the given covering action;
it does not assert that every topological cover admits a lift
of every lattice action.
\end{proof}

\section{The normalized formal collar}\label{III-series:formal-collar-three}

\begin{proposition}[Neutral curves and all finite transverse jets]
\label{III-arbitrary_lattices:prop:jets}
At every recognized neutral projective covering hypersurface
for an arbitrary lattice, after the stipulated finite-index
and covering choices, the following are available:
\begin{enumerate}
\item actual normalized analytic transverse curve germs
$i_x(z)$ almost everywhere on the homogeneous suspension;
\item an analytic central generator $B(z)\partial_z$, with
$B(0)=B'(0)=0$, such that
\[
 F(h,x)i_x(z)=i_{hx}(z),\qquad
 F(a_t,x)i_x(z)=i_{a_tx}(\Phi_t^B(z));
\]
\item an analytic invariant nowhere-zero section of the normal
quotient line, and analytic compatible sections for every
finite unparametrized and normalized parameter jet.
\end{enumerate}
If $B\equiv0$, the boundary already has its full analytic
product collar.
\end{proposition}
\begin{proof}
The fixed-source curve convergence theorem in Section~\ref{III-series:fixed-curves}
uses the auxiliary $H_0$ spectral gap just verified, together
with a hyperbolic element of its $\SL(2,\Z)$ Levi. That
element is likewise a suspension symmetry, not an element
required in $\Gamma$. This proves actual curve convergence;
the one-variable normalization in Section~\ref{III-series:relative-three} supplies
the analytic central flow.

When $B=0$, every fixed-parameter graph probability is
$Q$-invariant. The stationary argument in Section~\ref{III-series:axis-three}
then uses Theorem~\ref{III-coordinate_recognition:thm:recognition}, a common real
continuation interval, and the cross theorem to obtain a
joint collar patch meeting the zero section. Projective
minimality gives the full product collar. Its proof does
not use a rational boundary point. Applied to the derivative
action on the compact analytic $3$-manifold
$\mathbb P(\nu\oplus\R)$, where $\nu$ is the normal
quotient, this stationary result supplies an analytic
invariant normal section in every neutral case, as in
Section~\ref{III-series:axis-three}.

The induction in Section~\ref{III-series:jets-three} applies. Its affine jet
bundle has difference bundle $TN\otimes(\nu^*)^{\otimes k}$,
which the invariant normal frame identifies with $TN$.
The projective compactification and its measurable graph
give a homogeneous parabolic suspension probability.
The regularity argument uses root graphs, Howe--Moore
ergodicity of root kernels, and analytic stable charts.
These constructions have the same finite-volume inputs for
every lattice; in particular no rational stabilizer occurs.
Thus every unparametrized finite jet is analytic.

For the parameters, work in the recognized standard coordinates. Let
$o$ be the standard axis and put
\[
 c(g,\theta)=\log\|gv\|-\log\|v\|\quad(0\ne v\in\theta),
 \qquad \chi(q)=c(q,o),\qquad k(g)=c(g,g^{-1}o).
\]
The deterministic source homomorphism is trivial on $H$, so its
value at $q\in Q$ is $\Phi^B_{\chi(q)}$ in the paper's central-time
normalization. For the physical Haar lift $J_g$ of the same curve
germs, set $K_g=J_g\circ\Phi^B_{k(g)}$. The cocycle identity gives
\[
 k(qg)=k(g)+\chi(q),\qquad
 k(g\gamma^{-1})=k(g)-c(\gamma,g^{-1}o),
\]
and therefore, for each fixed group element, almost everywhere as germs,
\[
 K_{qg}=K_g,\qquad
 \alpha(\gamma)\circ K_g
       =K_{g\gamma^{-1}}\circ\Phi^B_{c(\gamma,g^{-1}o)}.
\]
Descending their measurable Taylor coefficients through the Haar
measure class gives angular parameters with this same analytic clock.
The source correction leaves their unparametrized jets unchanged.
The additive-transfer compactification in Section~\ref{III-series:parameter-jets}
uses Theorem~\ref{III-coordinate_recognition:thm:recognition} on a compact analytic
$3$-manifold and projective minimality. It therefore
makes every normalized parameter coefficient analytic as
well. Compatibility follows from almost-everywhere agreement
and analytic uniqueness.
\end{proof}

\section{Joint analyticity of neutral collars}\label{III-series:joint-three}

\subsection{Joint convergence}

\noindent\textbf{Neutral joint-analyticity statement (NJ).}
For every detected neutral projective covering hypersurface
whose normalized analytic central generator satisfies
$B\not\equiv0$, the normalized transverse curve field of
Proposition~\ref{III-arbitrary_lattices:prop:jets} has a jointly analytic
representative on one nonempty angular patch, including
$z=0$. The representative has the prescribed finite jets
and hence nonzero transverse derivative. The statement concerns
the normalized field, not an arbitrary tubular chart.

\begin{theorem}[Neutral joint analyticity for every lattice]
\label{III-arbitrary_lattices:thm:joint}
Statement (NJ) holds for every lattice in $\SL(3,\R)$, including
its finite-index and covering reductions.
\end{theorem}
\begin{proof}
Work in the finite-volume suspension $X=G/\Gamma$ and restrict
its Levi action to an auxiliary lattice $\Lambda<\SL(2,\R)$.
This is a lattice in the suspension symmetry group; no inclusion
$\Lambda<\Gamma$ is needed. In actual fiber coordinates write
the formal collar supplied by Proposition~\ref{III-arbitrary_lattices:prop:jets} as
\[
 H_x(u,z)=i_x(z)+A_x(z)u+\sum_{p\ge2}h_{p,x}(z)[u^p].
\]
The source Levi action is $(u,z)\mapsto(lu,z)$.
The axis $i_x$ is analytic almost everywhere. The frame theorem
in Section~\ref{III-series:frames} first makes $A_x(z)k(z)$ analytic for one scalar
formal gauge $k(0)=1$. Since $B\ne0$, the gauge-removal lemma
in Section~\ref{III-series:gauge} makes $k$, and therefore the original $A_x$,
analytic. Indeed, central covariance makes
$k(z)/k(\Phi_t^B(z))$ analytic for every $t$; polynomial
dependence of each logarithmic coefficient on $t$, Baire's
theorem, and polynomial interpolation make $B(\log k)'$
analytic. Division by the finite-order zero of $B$ concludes
the argument.

The actual analytic frame $P_x(z)=(A_x(z),i_x'(z))$ is
invertible as a germ. With $V=\R^2\oplus\R$, set
\[
 h_x(u,z)=P_x(z)^{-1}\bigl(H_x(u,z)-i_x(z)\bigr),
\]
\[
 \alpha(l,x;z;y)=P_{lx}(z)^{-1}
 \bigl(F(l,x)(i_x(z)+P_x(z)y)-i_{lx}(z)\bigr).
\]
These formulas are in the chosen actual fiber charts.
The cocycle $\alpha$ is jointly analytic in $(z,y)$,
fixes $y=0$, and has derivative $\operatorname{diag}(l,1)$
for every $z$. Moreover $h_x=(u,0)+O(u^2)$ and
\[
             \alpha(l,x;z;h_x(u,z))=h_{lx}(lu,z)
\]
formally. The parameter-uniform rectangular convergence theorem
in Section~\ref{III-series:parameters} applies and proves joint convergence
of $h_x$, hence of $H_x$, almost everywhere.

The essential uniformity is worth specifying. Place the normal
truncation in $A_m=\mathbb C[z]/(z^{m+1})$ with squared
coefficient Hilbert norm $\sum_{j\le m}(j+1)^2|a_j|^2$.
Multiplication has norm bounded independently of $m$.
Rectangular substitution matrices have entries in $A_m$;
their norm is the Hilbert norm of those entries, not the
Hilbert--Schmidt norm of multiplication operators on $A_m$.
Both substitution length and the comparison of coefficient
records with matrix height are uniform in $m$.
The angular coefficient representation has highest weights
$p-1,p,p+1$ for $p\ge2$, with $A_m$ only a trivial multiplicity.
Thus every positive coefficient record is detected by one fixed
compact set in the convex norm profile argument. The resulting
probability bound is uniform in both truncation orders.
Monotonicity of the coefficient bounds in these two orders gives
joint convergence. No common normal disk for the individual
$h_{p,x}$ is assumed in this argument.

Choose a center where joint convergence holds. Every coefficient
in the normal expansion of its sum agrees near that center with
the prescribed analytic finite-jet section. Analytic uniqueness
extends each such equality throughout a common connected angular
patch. At almost every angular point the original normal curve
is analytic and has those same coefficients, so it equals the
sum as a normal germ. This gives a jointly analytic representative
of the normalized field on that patch and proves (NJ).
\end{proof}

\begin{lemma}[One joint patch supplies the full neutral collar]
\label{III-arbitrary_lattices:lem:patch}
Under (NJ), every detected neutral projective covering
hypersurface has its full-angular analytic collar, with the original normalized generator $B$.
\end{lemma}
\begin{proof}
The case $B=0$ is Proposition~\ref{III-arbitrary_lattices:prop:jets}.
Otherwise lift the joint patch to an open subset of $G$
using parabolic equivariance and the analytic flow of $B$.
The identity at the zero section and its nonzero transverse
derivative make this a local collar. Projective minimality
means that the right $\Gamma$ translates of the resulting
left-parabolic open set cover $G$. Transport the analytic
family by the original right-lattice action. On overlaps,
two transported families represent the same measurable
normalized curve field, so analytic uniqueness makes them
equal as germs. This is the open-set descent construction
in Section~\ref{III-series:negative-three}, Section~\ref{III-series:positive-three}; it avoids restricting an
almost-everywhere field directly to a null transverse section.
The glued analytic family descends to the angular quotient.
Realize it on the actual covering neighborhood, as in
Section~\ref{III-series:negative-three}. Hatcher's ENR theorem
\cite[Corollary~A.9, p.~527; supporting Theorem~A.7, pp.~525--527]{Hatcher}
and coordinatewise Tietze extension give a continuous neighborhood
retraction onto the original compact hypersurface. Pull back its
specified angular cover along that retraction. Lift the same developing
family at its prescribed zero section; uniqueness of covering lifts
identifies the overlaps. Compactness now gives a smaller embedded
collar in this covering neighborhood, whose projection is the original
developing family as a germ. The genuine deck maps retain their
prescribed zero-section action. Passing to the kernel of their finite
permutation action makes them commute with the original lattice lifts.
Thus the collar retains the original lattice action, normalized
parameter and generator $B$, together with the specified finite-cover
and normal-orientation deck data. The deck maps here retain their
actual analytic germs, as in Part~II; no linear transverse formula
is needed in the neutral case.
\end{proof}

\section{Completion of the local description for \texorpdfstring{$n=3$}{n = 3}}\label{III-series:completion-three}
\begin{proof}[Proof of Theorem~\ref{III-main:thm:recognition-three} for $n=3$]
The invariant-probability statement follows from
Theorem~\ref{II-I:thm:invariant-germs}. For a proper-parabolic
suspension measure, the parabolic dimension count leaves the line,
plane, and full-flag branches. The full-flag case is
Corollary~\ref{III-arbitrary_lattices:cor:flag}.

For a line or plane parabolic, coordinate description in
Section~\ref{III-series:coordinate-curves} gives analytic projective
hypersurfaces in the finite-atomic cases. Sections
\ref{III-series:relative-three}--\ref{III-series:positive-three} give the
nonzero-exponent collars. In the neutral case,
Sections~\ref{III-series:shear-three}--\ref{III-series:regular-three}
give the radical dichotomy and the regular Hopf manifold.
For zero shear, the coherent-graph construction gives finite
atomic conditionals and hence an analytic projective hypersurface.
Sections~\ref{III-series:axis-three}--\ref{III-series:formal-collar-three}
give the normalized axis and all finite transverse jets.
Theorem~\ref{III-arbitrary_lattices:thm:joint} and
Lemma~\ref{III-arbitrary_lattices:lem:patch} give its full collar.

All these arguments apply to each finite-index subgroup. For a
measure supported over a compact invariant subset, use the support argument in Part II: each positive-mass hypersurface
orbit meets the invariant set, and projective minimality makes
the intersection with each of its hypersurfaces the whole hypersurface.
The angular representation, deck action, and normal-orientation
cover are retained in each collar construction. This verifies
(G), (P), and (C) for $n=3$.
\end{proof}

\section{Equivariant curves with a fixed source parameter}\label{III-series:fixed-curves}

\subsection{Statement and normalization}

Let a countable group $\Lambda$ act ergodically preserving a
probability $(X,\mu)$. Assume a finite symmetric set $S$ and $a>0$
satisfy
\begin{equation}\label{III-transverse_curve_convergence:eq:gap}
 \sum_{s\in S}\|f\circ s-f\|_2^2
 \ge a\left\|f-\int f\,d\mu\right\|_2^2.
\end{equation}
Let $\rho(\gamma)=\operatorname{diag}(\sigma(\gamma),1)$, where
$\sigma:\Lambda\to\operatorname{GL}(2,\R)$ contains an element
$b$ with two positive eigenvalues $\lambda>1$ and $\lambda^{-1}$.
Enlarge $S$ to contain $b,b^{-1}$ and choose a fixed tangential
basis diagonalizing $\sigma(b)$.

\begin{corollary}[Realization with a fixed source parameter]
\label{III-transverse_curve_convergence:thm:curve}
Let $F(\gamma,x)$ be a coefficient-measurable cocycle of analytic
diffeomorphism germs of $(\R^3,0)$ with derivative $\rho(\gamma)$.
Suppose measurable formal curve germs
\[
 I_x(z)=e_3z+O(z^2)
 \quad\hbox{satisfy}\quad
 F(\gamma,x)I_x(z)=I_{\gamma x}(z).
\]
There is a constant formal one-variable germ $C(z)=z+O(z^2)$
such that $K_x=I_x\circ C$ converges almost everywhere.
The resulting embedded analytic curve germs have positive measurable
radii and satisfy the same equivariance identities.
\end{corollary}
\begin{proof}
Apply Theorem~\ref{II-sg:thm:fixed-subspace} with
$V=\R^2\oplus\R$ and $E=\R e_3$. The representation fixes
$E$ pointwise, and the chosen element $b$ has $E$ as its
$1$-eigenspace, with the other two eigenvalues off the unit
circle. The spectral-gap and formal-equivariance hypotheses
are exactly those above. A constant formal change of the
one-dimensional source preserves all the equivariance identities.
\end{proof}

The normalization in that proof chooses the nonlinear coefficients
of the fixed component to have lower median zero. We retain the
following finite-jet estimate for use in the frame argument.

Use a Hilbert norm $\|\cdot\|_*$ on scalar one-variable finite jets
obtained by multiplying
\[
 \left(\sum_{j\ge0}(1+j)^{2p}|c_j|^2\right)^{1/2},\qquad p=1,
\]
by a fixed sufficiently large constant. Weighted convolution makes
the resulting norm submultiplicative. Truncation $\Pi_N$ is an
orthogonal projection, and we use the Hilbert direct-sum norm on
the three output components. All constants are independent of $N$.
Put $J=\|e_3z\|_*$.

For every $\beta>0$, coefficient measurability and analyticity give
a set $B_\beta$ with $\mu(B_\beta)\le\beta$, a radius $t_\beta>0$,
and $A_\beta<\infty$ such that, for $x\notin B_\beta$, $s\in S$,
and $\|V\|_*\le t_\beta$,
\begin{equation}\label{III-transverse_curve_convergence:eq:quadratic}
 \|\Pi_N(F(s,x)V-\rho(s)V)\|_*
 \le A_\beta\|V\|_*^2.
\end{equation}
Here $V(0)=0$ and juxtaposition with $F$ means composition.
The bound follows by summing the absolutely convergent coefficient
majorant of the finitely many germs. It does not assume that $K$
converges.
For the given radius and cutoff, the full finite jet of $K_x(r_nz)$
is $e_3r_nz+\epsilon_d w_d(x)$. Formal equivariance gives its exact
truncated composition equation. Since $\rho(s)e_3=e_3$, subtracting
the same linear term and dividing by $\epsilon_d$ gives
\eqref{III-transverse_curve_convergence:eq:jet}. Thus the exceptional set
and constants above depend only on the original finite family $F(s,x)$,
independently of the cutoff and the radius sequence.

\subsection{Two estimates uniform in the jet dimension}

Write $m=|S|$ and $\vartheta=1+a/(4m)>1$.
Applying \eqref{III-transverse_curve_convergence:eq:gap} to an indicator function gives
\begin{equation}\label{III-transverse_curve_convergence:eq:expansion}
 \mu\left(A\cup\bigcup_{s\in S}sA\right)
 \ge\vartheta\mu(A)\qquad(\mu(A)\le1/2).
\end{equation}

\begin{lemma}[Tightness from bounded multiplicative steps]
\label{III-transverse_curve_convergence:lem:tight}
Suppose $w_d:X\to\mathcal H_d^3$ are finite-dimensional Hilbert
coefficient vectors, $\epsilon_d\downarrow0$, and
$r_nJ\le\epsilon_d^2$. Suppose they satisfy the exact finite-jet
equations
\begin{equation}\label{III-transverse_curve_convergence:eq:jet}
 w_d(sx)=\rho(s)w_d(x)
 +\epsilon_d^{-1}\Pi_{N_d}
 \bigl(F(s,x)V_{d,x}-\rho(s)V_{d,x}\bigr),\qquad
 V_{d,x}=e_3r_nz+\epsilon_d w_d(x).
\end{equation}
If $\mu(\|w_d\|_*>1)\le\delta<1/2$, then
\[
 \lim_{R\to\infty}\limsup_{d\to\infty}
          \mu(\|w_d\|_*>R)=0.
\]
Moreover $w_d\circ s-\rho(s)w_d\to0$ in probability for every
$s\in S$.
\end{lemma}
\begin{proof}
Put $M=\max_{s\in S}\|\rho(s)\|$ and choose $K>\max\{M,1\}$.
For each fixed $t\ge1$ and $\beta>0$, \eqref{III-transverse_curve_convergence:eq:quadratic}
implies, for all sufficiently large $d$,
\[
 x\notin B_\beta,\quad\|w_d(x)\|_*\le t
 \quad\Longrightarrow\quad\|w_d(sx)\|_*\le Kt.
\]
Indeed the error is bounded by
$A_\beta\epsilon_d(t+\epsilon_d)^2$.
Set $A_{d,t}=\{\|w_d\|_*>t\}$. Symmetry of $S$ gives
\[
 A_{d,Kt}\cup\bigcup_{s\in S}sA_{d,Kt}
                  \subset A_{d,t}\cup B_\beta.
\]
All the sets $A_{d,t}$, $t\ge1$, have measure at most $\delta$.
Thus \eqref{III-transverse_curve_convergence:eq:expansion} yields
\[
 \mu(A_{d,Kt})\le\vartheta^{-1}
                       (\mu(A_{d,t})+\beta).
\]
For a fixed integer $j$, iterate at $t=1,K,\ldots,K^{j-1}$,
take the upper limit in $d$, and then let $\beta\downarrow0$.
This gives
\[
 \limsup_d\mu(A_{d,K^j})\le\delta\vartheta^{-j},
\]
which is tightness. For fixed $R$ the error in \eqref{III-transverse_curve_convergence:eq:jet}
tends uniformly to zero on
$\{\|w_d\|_*\le R\}\setminus B_\beta$.
Letting first $d\to\infty$, then $\beta\downarrow0$ and
$R\to\infty$, proves convergence of the defects in probability.
\end{proof}

\section{Analytic frames along equivariant curves}\label{III-series:frames}

\subsection{Spectral gap for the block Levi}

Put $H=L\ltimes U$, where $L=\SL(2,\R)$ and
$U=\R^2$ is standard or dual standard. Let
$\Lambda<L$ contain two hyperbolic elements whose pairs of
projective fixed points are disjoint; a suitable finite set in
$\SL(2,\mathbb Z)$ has this property.

\begin{lemma}\label{III-neutral_frame_convergence:lem:levi-gap}
If a strongly continuous unitary representation $\pi$ of $H$
has no nonzero $U$-fixed vector, its restriction to $\Lambda$
has spectral gap. A finite spectral-gap set can be chosen to
depend only on the indicated two hyperbolic elements.
\end{lemma}
\begin{proof}
The spectral resolution of the restriction of $\pi$ to the
abelian group $U$ assigns a unit vector $v$ a probability
$\mu_v$ on $U^*$, and $\mu_v(\{0\})=0$. For a Borel set
$E\subset U^*$ its value is $\langle P_Ev,v\rangle$.
The variation norm of the difference of the spectral
probabilities of two unit vectors $v,w$ is at most
$2\|v-w\|$. Indeed, this follows by integrating a Borel
function of absolute value at most one against the spectral
resolution and using the resulting operator's norm bound.
Conjugation by $l\in L$ transports $\mu_v$ by the dual
linear action to $\mu_{\pi(l)v}$.

Push $\mu_v$ to the compact projective line $\mathbb P(U^*)$.
If unit vectors were arbitrarily close to invariant under
the two selected hyperbolic elements, their projective
probabilities would have a weak limit invariant under both.
An invariant probability for one hyperbolic projective map
is supported on its two fixed points: every other point is
nonrecurrent, whereas almost every point of an invariant
probability is recurrent. The disjoint fixed-point pairs
therefore make such a common invariant probability impossible.
This contradiction gives a uniform positive lower bound on
the maximum displacement of a unit vector by this finite set.
\end{proof}

\begin{corollary}\label{III-neutral_frame_convergence:cor:homogeneous-gap}
Let $\Gamma<\SL(3,\R)$ be any lattice and let
$X=\SL(3,\R)/\Gamma$ carry Haar probability. The action of
the block subgroup $\SL(2,\mathbb Z)<L$ on $X$ is ergodic
and has spectral gap.
\end{corollary}
\begin{proof}
The full $\SL(3,\R)$ action on its normalized Haar quotient is
transitive and ergodic, and its zero-mean Koopman representation is
strongly continuous. By the unipotent Mautner argument, a vector fixed
by the cross-block unipotent group $U$ is fixed by all of $\SL(3,\R)$,
so it is zero in $L^2_0(X)$. The spectral-probability argument of
Lemma~\ref{III-neutral_frame_convergence:lem:levi-gap} also applies to
\[
 u=\begin{pmatrix}1&1\\0&1\end{pmatrix},\qquad
 j=\begin{pmatrix}0&-1\\1&0\end{pmatrix}.
\]
A probability invariant under the projective unipotent map is supported
on its unique fixed line, and $j$ moves that line. These two integral
block matrices therefore have a positive uniform displacement bound
on the unit sphere of $L^2_0(X)$. Squaring and scaling gives their
scalar spectral-gap inequality on the original Haar probability.
A vector fixed by the entire integer block subgroup is fixed by both
matrices and hence is zero. The invariant-function criterion proves
ergodicity of that same subgroup.
\end{proof}

\subsection{The frame equation}

Let a countable subgroup $\Lambda<L$ act ergodically preserving
a probability $\mu$ on $X$, with spectral gap. Assume it
contains two positive hyperbolic elements whose infinitesimal
centralizers in $\mathfrak{sl}_2$ are distinct. Enlarge a
finite symmetric spectral-gap set to include these elements
and their inverses.

Suppose $F(l,x)$ is a coefficient-measurable analytic germ
cocycle on $(\R^3,0)$ with derivative
$\rho(l)=\operatorname{diag}(l,1)$. Suppose actual embedded
analytic curves $i_x(z)=e_3z+O(z^2)$ satisfy
\[
                    F(l,x)i_x(z)=i_{lx}(z).
\]
Write $E:\R^2\to\R^3$ for the first-coordinate inclusion.
Assume a formal matrix series
$A_x(z)\in\Hom(\R^2,\R^3)[[z]]$, with $A_x(0)=E$,
satisfies
\begin{equation}\label{III-neutral_frame_convergence:eq:frame}
 D F(l,x)|_{i_x(z)} A_x(z)=A_{lx}(z)l.
\end{equation}
All coefficients and the curve germs are assumed measurable.

\begin{theorem}[Analytic tangential frames]\label{III-neutral_frame_convergence:thm:frame}
There is a constant formal scalar series $k(z)$ with
$k(0)=1$ such that $\widetilde A_x(z)=A_x(z)k(z)$
converges almost everywhere. Its rank is two on a positive
measurable neighborhood of zero, and \eqref{III-neutral_frame_convergence:eq:frame} holds
for $\widetilde A$ as an analytic identity. No common
analytic radius or moment hypothesis is required.
\end{theorem}
\begin{proof}
The representation on the six-dimensional real space
$W=\Hom(\R^2,\R^3)$ is
\[
               \pi(l)T=\rho(l)T l^{-1}.
\]
It decomposes as
\[
 W=\R E\oplus\mathfrak{sl}_2\oplus(\R^2)^*.
\]
The invariant projection is $T\mapsto\tau(T)E$, where
$\tau(T)$ is half the trace of its top $2$-by-$2$ block.
Choose a fixed Euclidean norm adapted to this direct sum.

Multiplication by a constant scalar formal series preserves
\eqref{III-neutral_frame_convergence:eq:frame}. Its degree-$j$ coefficient changes
$\tau(A_{x,j})$ by the degree-$j$ coefficient of that scalar
series plus terms already determined in earlier degrees.
Choose $k$ inductively so every nonconstant scalar
coefficient of $\tau(\widetilde A_x)$ has lower median zero.
For the rest of the proof replace $A$ by this normalization.

Use the weighted coefficient Hilbert algebra and its
orthogonal finite-jet truncations $\Pi_N$ from
Section~\ref{III-series:fixed-curves}. Its scalar norm is a fixed multiple of
$\bigl(\sum_{j\ge0}(1+j)^{2p}|c_j|^2\bigr)^{1/2}$,
where $p>1/2$. All product estimates are independent of $N$.
The curve $i_x$ is already analytic. Consequently, outside
a set of arbitrarily small measure, the finitely many
matrix-valued analytic functions
\[
           M_{s,x}(z)=D F(s,x)|_{i_x(z)}
\]
have common complex domains and bounds. They satisfy
$M_{s,x}(0)=\rho(s)$. These are derivatives of the original
cocycle along the original curves. The chain rule and curve
invariance give their analytic cocycle and inverse identities;
Taylor composition makes every coefficient measurable.
On such a set, for all sufficiently
small deterministic $r>0$,
\begin{equation}\label{III-neutral_frame_convergence:eq:linear-error}
       \|M_{s,x}(rz)-\rho(s)\|_*\le C r
\end{equation}
with one deterministic $C$. This is only a large-measure
bound, not an essential bound on $X$.

Convergence of $A_x$ is invariant under $\Lambda$: solve
\eqref{III-neutral_frame_convergence:eq:frame} in either direction using the actual
analytic invertible matrix $M_{s,x}(z)$. Suppose, by
ergodicity, that $A_x$ diverges almost everywhere. The
first-crossing construction in Section~\ref{III-series:fixed-curves} gives
$\epsilon_d\downarrow0$, $r_d\le\epsilon_d^2$ and finite
orders $N_d$ such that
\[
 w_d(x)=\epsilon_d^{-1}\Pi_{N_d}(A_x(r_dz)-E),\qquad
 \mu(\|w_d\|_*>1)\le\tfrac14,\quad
 \mu(\|w_d\|_*\ge1)\ge\tfrac14.
\]
The exact finite-jet equation is
\[
 w_d(sx)=\pi(s)w_d(x)
  +\epsilon_d^{-1}\Pi_{N_d}
   \bigl((M_{s,x}(r_dz)-\rho(s))
               (E+\epsilon_d w_d(x))s^{-1}\bigr).
\]
On any of the large-measure bounded sets just described,
the error norm is bounded by
\[
                  C'\epsilon_d(1+\epsilon_d\|w_d(x)\|_*).
\]
In particular, for every fixed $t\ge1$ it tends uniformly
to zero where $\|w_d(x)\|_*\le t$. The multiplicative-step
tightness lemma in Section~\ref{III-series:fixed-curves} applies verbatim: the target
representation $\pi$ has fixed finite dimension and
$\max_{s\in S}\|\pi(s)\|$ is fixed. Thus $w_d$ is tight
and
\[
                 w_d\circ s-\pi(s)w_d\longrightarrow0
                     \quad\text{in probability}.
\]

Diagonalize each positive hyperbolic element in its own fixed
real eigenbasis. Its eigenvalues are distinct, positive and
reciprocal. In the first eigenbasis the second element has a
nonzero off-diagonal entry; otherwise their traceless centralizers
would coincide. Thus the common fixed space of the two original
frame representations is exactly $\R E$. This calculation preserves
the original cocycle, curves, frame coefficients and scalar gauge.
For each element, project to its expanding and contracting
eigenspaces in $W$. Iteration
through a fixed number of powers and tightness, exactly as
in Section~\ref{III-series:fixed-curves}, make both projected fields tend to zero
in probability. The intersection of the two zero-weight
spaces is $\R E$: the dual standard summand has no zero
weight, and the two zero-weight lines in
$\mathfrak{sl}_2$ are distinct. Finite-dimensional linear
algebra, with constants independent of the coefficient-jet
dimension, therefore makes the component of $w_d$ in
$\mathfrak{sl}_2\oplus(\R^2)^*$ tend to zero in probability.

The scalar component $\tau(w_d)$ is approximately invariant
and each of its coefficient coordinates has median zero.
Radial clipping to the unit Hilbert ball, followed by the
dimension-independent median spectral-gap estimate from
Section~\ref{III-series:fixed-curves}, makes this component tend to zero as well.
Hence $w_d\to0$ in probability, contradicting the displayed
lower bound at the first-crossing radius.

Therefore the same scalar-normalized frame $\widetilde A_x=A_xk$
has positive radius almost everywhere. Its constant coefficient
is $E$. Root tests and the nonvanishing of a fixed two-by-two
minor select a measurable positive radius $r(x)$ on which its
sum is analytic and has rank two. For every $l$, the formal frame
equation and analytic Taylor uniqueness give the actual identity
on a positive interval inside both radii $r(x)$ and $r(lx)$.
Countability of $\Lambda$ gives one common conull set. The radius
may depend on $x$; no moment hypothesis has been used.
\end{proof}

\section{Removing a scalar formal gauge}\label{III-series:gauge}

\subsection{The formal field and the scalar frame gauge}

Put $L=\SL(2,\R)$, $E=\R^2$, $V=E\oplus\R$, and
$\rho(l)=\operatorname{diag}(l,1)$. Let $X=\SL(3,\R)/\Gamma$
for an arbitrary lattice $\Gamma$. The block
$\Lambda=\SL(2,\mathbb Z)<L$ acts ergodically with spectral gap,
by the spectral-measure argument in Section~\ref{III-series:frames}.

Suppose the coefficient-measurable formal field
$\widehat H_x(u,z)$ conjugates the actual analytic fiber transport
$F(l,x)$ to $(u,z)\mapsto(lu,z)$ and has invertible derivative.
Its axis $i_x(z)=\widehat H_x(0,z)$ is an actual analytic transverse
curve, and its formal tangential frame is
$A_x(z)=D_u\widehat H_x(0,z)$. Fix coordinates in which
$i_x'(0)=e_3$, $A_x(0)$ is the standard inclusion, and
$D F(l,x)(0)=\rho(l)$.

Assume also the exact formal central covariance
\begin{equation}\label{III-angular_slices:eq:central}
 F(a_t,x)\widehat H_x(u,z)
 =\widehat H_{a_tx}(e^{\chi t}u,\phi_t(z)),
 \qquad \phi_t=\Phi_t^B,
\end{equation}
where $B$ is analytic, $B(0)=B'(0)=0$, and $B\not\equiv0$.
The action on the source angular variables is scalar because
$a_t$ centralizes $L$. In applications, the projective-character
coboundary is removed before writing \eqref{III-angular_slices:eq:central}.

The frame theorem in Section~\ref{III-series:frames} supplies a constant formal scalar
$k(z)$, $k(0)=1$, for which $\widetilde A_x=A_xk$ is analytic
almost everywhere. Central covariance removes this apparent
formal ambiguity.

\begin{lemma}[The scalar frame gauge converges]\label{III-angular_slices:lem:gauge}
Under \eqref{III-angular_slices:eq:central}, the series $k$ is analytic. Consequently
the original frame $A_x$ is analytic almost everywhere.
\end{lemma}
\begin{proof}
Formal differentiation of \eqref{III-angular_slices:eq:central} in $u$,
justified at each finite order by the chain rule, gives
\[
 D F(a_t,x)|_{i_x(z)}\widetilde A_x(z)
 =e^{\chi t}\widetilde A_{a_tx}(\phi_t(z))
                  \frac{k(z)}{k(\phi_t(z))}.
\]
For every fixed real $t$, choose $x$ in the conull set where the
two displayed frames and the transport are analytic. An analytic
left inverse of the rank-two frame on the right shows that
$r_t(z)=k(z)/k(\phi_t(z))$ is an analytic scalar germ.

Set $g=\log k$. Then
$\psi_t(z)=g(z)-g(\phi_t(z))=\log r_t(z)$ is analytic for each
real $t$. Joint analyticity of the original flow near $(0,0)$
identifies the time derivatives of its Taylor coefficients with
the Taylor coefficients of $B$; the first two vanish.
Since $B$ vanishes to order at least two, the coefficient
$P_d(t)=[z^d]\psi_t$ is a polynomial in $t$ of degree at most $d$.
For integers $M\ge1$ set
\[
 E_M=\{t\in[-1,1]: |P_d(t)|\le M^d\text{ for every }d\ge1\}.
\]
These sets are closed and cover $[-1,1]$. Baire's theorem gives
a nondegenerate interval $J$ on which all the indicated estimates
hold for one $M$. For a degree-at-most-$d$ polynomial, interpolation
at $d+1$ equally spaced points of $J$ gives
\[
 |P_d'(0)|\le d C_J^d\sup_{t\in J}|P_d(t)|.
\]
Indeed the denominators of the Lagrange basis have absolute values
$(|J|/d)^d j!(d-j)!$; summing their reciprocals and using
$d^d/d!\le e^d$ gives an exponential bound, and differentiating
the numerators contributes at most a factor $d$ and another fixed
exponential factor.

Thus $\partial_t\psi_t|_{t=0}=-Bg'$ is analytic. Write
$B(z)=z^m b(z)$ with $b(0)\ne0$. The analytic numerator has its
first $m$ coefficients zero because its formal quotient by $B$
is the formal series $g'$. Division by $z^m b(z)$ is therefore
analytic. Integrating and exponentiating proves that $g$ and $k$
converge. Since $k(0)=1$, its reciprocal is analytic near zero,
and $A_x=\widetilde A_x/k$ converges on a positive neighborhood
for almost every $x$.
\end{proof}

\section{Convex norm profiles and uniform growth}\label{III-series:profiles}
\label{III-subsec:ha-convex-profiles}

There is a fixed compact target for logarithmic norm profiles of
all finite-dimensional algebraic representations. We will use these profiles to control analytic substitution matrices
without a bound on their dimensions.
Throughout this section, $G$ is connected, real semisimple and
linear, has finite center and no compact factors, and $\Gamma<G$
is a lattice. Fix a maximal compact subgroup $K$, the Cartan
decomposition $\mathfrak g=\mathfrak k\oplus\mathfrak p$, and the
associated invariant Riemannian metric on $M=G/K$. Write $o=K$.
Set
\[
 \mathcal C=\{f:M\longrightarrow\R:
     f\text{ is geodesically convex and }1\text{-Lipschitz},\ f(o)=0\},
\]
with uniform convergence on compact sets. Arzela--Ascoli and a
countable exhaustion by compact balls show that $\mathcal C$ is
compact and metrizable. The continuous action is
\begin{equation}\label{III-eq:ha-profile-action}
       (a\cdot f)(x)=f(a^{-1}x)-f(a^{-1}o).
\end{equation}

\begin{lemma}[The invariant law on the convex target]
\label{III-lem:ha-profile-invariant-law}
The only $G$-invariant Borel probability on $\mathcal C$ is
$\delta_0$. The same is true for $\Gamma$-invariant probabilities.
\end{lemma}
\begin{proof}
We use the actual Cartan quotient and its invariant Riemannian metric.
The one-parameter subgroups are fixed by their initial velocities
\cite[Theorem~20.1, p.~516]{LeeSmoothManifolds2e}; the Cartan involution,
fixed subgroup and polar decomposition are supplied by
\cite[Corollary~6.18, p.~358; Theorem~6.31, pp.~362--367]{KnappBeyondDigital}.
The global theorem supplies its own exponential; its one-parameter
subgroups have the prescribed initial velocities, so uniqueness identifies
it with the exponential already used here.
The smooth quotient and its original left action are those of
\cite[Theorem~21.17, pp.~551--552]{LeeSmoothManifolds2e}.
The Cartan form gives the orthogonal involutive Lie algebra, hence the
quotient metric and symmetry of \cite[IV.2.5, pp.~74--76]{BorelSymmetric1998}.
The symmetric manifold is complete by \cite[IV.2.1, p.~73]{BorelSymmetric1998}.
Here geodesics can equivalently be expressed as constant-speed metric
local geodesics: the metric definition is
\cite[I.1, Definition~1.3, p.~4]{BridsonHaefliger1999}, and its
identification with the Riemannian geodesic equation is stated in
\cite[II.1, Appendix~1A.3, p.~170]{BridsonHaefliger1999}.
Our smooth Riemannian metric meets its regularity hypothesis. Thus the
metric interpolation constraints and the original smooth-geodesic
convexity condition describe the same full target.

Its full local Killing algebra is a quotient of the original semisimple
algebra by the last assertion of IV.2.5, so is semisimple.
Thus it has no local Euclidean factor
\cite[IV.2.7, pp.~76--77]{BorelSymmetric1998}.

Let $\nu$ be $G$-invariant. The mean
$\bar f(x)=\int f(x)\,d\nu(f)$ is defined and continuous because
$|f(x)|\le n(o,x)$. Invariance implies
\[
       \bar f(a^{-1}x)-\bar f(a^{-1}o)=\bar f(x).
\]
Consequently $a\mapsto\bar f(ao)$ is a continuous real character
of $G$. It vanishes because a connected semisimple group is perfect
\cite[\S2, p.~224, unnumbered assertion]{Djokovic1986}; hence $\bar f=0$. For every geodesic midpoint
triple $(x,m,y)$, the nonnegative convexity deficit
$[f(x)+f(y)]/2-f(m)$ has integral zero. Use a countable dense
family of triples and continuity. Almost every $f$ is then affine
on every geodesic. A nonconstant such affine function would give a Riemannian product
$M\cong N\times\R$ by \cite[Main Theorem, p.~237; proof in \S4,
pp.~244--246]{Innami1982}. This would give a local Euclidean factor,
which was excluded above. Hence $f=0$ almost surely.

If $\nu$ is only $\Gamma$-invariant, average $g_*\nu$ over the
finite-volume space $G/\Gamma$. The resulting $G$-invariant
probability is $\delta_0$. Since $\{0\}$ is $G$-invariant, the
original probability is also $\delta_0$.
\end{proof}

Let $R$ be a nontrivial irreducible complex algebraic representation,
with a Cartan-compatible Hermitian norm, and put
\[
 L_R=\sup_{\|H\|=1,\ H\in\mathfrak p}\|dR(H)\|.
\]
The compact real form used below is constructed from the same differential
and acts unitarily for this Hermitian norm, as explained in the proof of
Proposition~\ref{III-prop:ha-profile-nonflattening}.
For $v\ne0$ define
\begin{equation}\label{III-eq:ha-algebraic-profile}
 f_{R,v}(gK)=\frac1{L_R}
                 \log\frac{\|R(g^{-1})v\|}{\|v\|}.
\end{equation}
This is well defined. On a geodesic $g\exp(tH)K$, the expression
is the logarithm of the norm of
$\exp(-t\,dR(H))R(g^{-1})v$, divided by $L_R$ and up to a constant.
The self-adjoint spectral decomposition expresses its squared norm
as a positive sum of exponentials. This proves convexity and the
$1$-Lipschitz bound. Moreover
\begin{equation}\label{III-eq:ha-profile-equivariance}
                   f_{R,R(a)v}=a\cdot f_{R,v}.
\end{equation}
The construction also applies to $R\otimes I$ on an arbitrary
Hilbert multiplicity space.

\begin{proposition}[Uniform nonflattening across algebraic types]
\label{III-prop:ha-profile-nonflattening}
There is a compact subset $Q'\subset G$, depending only on $G$
and the chosen metric, such that every nontrivial irreducible $R$
and every nonzero vector $v$, including arbitrary Hilbert
multiplicities, satisfy
\begin{equation}\label{III-eq:ha-profile-growth}
       \sup_{g\in Q'}\log\frac{\|R(g)v\|}{\|v\|}\ge L_R.
\end{equation}
In particular, zero does not belong to the compact-open closure of
all the profiles \eqref{III-eq:ha-algebraic-profile}.
\end{proposition}
\begin{proof}
The zero-dimensional group has no nontrivial irreducible representations;
assume first that $\dim\mathfrak g>0$.
All choices of Lie-algebra directions and compact sets in this proof precede
$R$ and its multiplicity space. We use the actual tangent algebra
$\mathfrak g$ of $G$ and its prescribed Cartan involution $\theta$.
The original adjoint representation and its differential are those of
\cite[Propositions~1.89 and~1.91, pp.~79--81]{KnappBeyondDigital}.
Its image and kernel, and the closedness of the semisimple adjoint group,
are given by \cite[Example~1.10.4, pp.~43--44, and Lemma~3.6.1,
p.~149]{DuistermaatKolk}. Thus the finite-center assumption gives the
finite adjoint covering of the same $G$. Its matrix tangent and exponential
chart use \cite[Introduction, (0.2)--(0.6), pp.~3--4, Proposition~0.14,
p.~11, and Theorem~0.15, pp.~12--13]{KnappBeyondDigital}.
The comparisons of actual differentials use automatic smoothness and
bracket preservation in \cite[Section~I.10, pp.~72 and~77]{KnappBeyondDigital},
and the kernel-codimension and image-rank assertions in
\cite[Chapter~1, Section~1.6, Theorem~2(1),(3), p.~8]{OnishchikVinberg1990}.
These are coordinate comparisons; the original representation and its
Hermitian norm, and the chosen norm on $\mathfrak p$, are retained.

Here is a direct construction of fixed nilpotent directions and their
bounds. The locally direct simple-factor decomposition
\cite[Chapter~4, Section~1.3, Theorem~4, p.~140]{OnishchikVinberg1990}
and the closed-normal-subgroup ideal assertion
\cite[Section~I.10, p.~77]{KnappBeyondDigital} identify the simple ideals
of the same tangent algebra. The absence of compact factors, the
negative-Killing compactness criterion
\cite[Corollary~3.6.3, p.~151]{DuistermaatKolk}, and Cartan involutions
\cite[Corollary~6.18, p.~358]{KnappBeyondDigital} give a nonzero split
direction in every simple ideal. Transporting the fixed Cartan geometry
uses the global decomposition for the given $\theta$
\cite[Theorem~6.31, pp.~362--367]{KnappBeyondDigital}; its compact factor
is the fixed subgroup of the global Cartan involution. The source supplies
its canonical exponential; uniqueness of normalized one-parameter subgroups
identifies it with any previously fixed normalized exponential and transports
the same polar map
\cite[Theorem~20.1 and following paragraph, p.~516]{LeeSmoothManifolds2e}.
In this way choose
$H\in\mathfrak p$ which generates $\mathfrak g$ as an ideal.
Since $\operatorname{ad}H$ is self-adjoint for the Cartan inner product,
it has a real eigenbasis. The span of the $G$-adjoint translates of its
nonzero-eigenvalue eigenvectors is an ideal containing
$[\mathfrak g,H]$, hence all of $\mathfrak g$. Extract a basis
$Y_i=\operatorname{Ad}(g_i)X_i$, $1\le i\le r=\dim\mathfrak g$,
where $[H,X_i]=\lambda_iX_i$ and $\lambda_i\ne0$.
Put $H_i=2H/\lambda_i$ and $C=1+\sum_i\|H_i\|$.
For $L=L_R$ and $D=\lfloor CL\rfloor$, the relations
$[dR(H_i),dR(X_i)]=2dR(X_i)$ and
$\|dR(H_i)\|\le\|H_i\|L$ show that
\[
                    dR(Y_i)^{D+1}=0.
\]
Indeed, successive applications of $dR(X_i)$ shift eigenvalues of the
self-adjoint $dR(H_i)$ by $2$; $D+1$ such shifts exceed its spectral
width. Conjugation by the original $R(g_i)$ gives the displayed bound.

Realize the underlying real Lie algebra of
$\mathfrak g_{\mathbb C}=\mathbb C\otimes_{\mathbb R}\mathfrak g$
by a simply connected Lie group $G_{\mathbb C}$ using
\cite[Theorem~20.21, p.~531]{LeeSmoothManifolds2e}.
The real form $\mathfrak u=\mathfrak k+i\mathfrak p$ has negative
Killing form. The compact-realization theorem
\cite[Theorem~3.4, p.~15]{KnappStructure1997}, the universal covering
theorem \cite[Theorem~7.7, p.~154]{LeeSmoothManifolds2e}, and Weyl's
theorem \cite[Theorem~4.69, pp.~268--269]{KnappBeyondDigital}
give a compact simply connected group $U$ with tangent algebra
$\mathfrak u$. Integrate its inclusion in $\mathfrak g_{\mathbb C}$,
and separately the complex-linear extension of $dR$, by
\cite[Theorem~20.19, pp.~530--531]{LeeSmoothManifolds2e}.
Exponential identities for these homomorphisms follow from the
uniqueness of one-parameter subgroups with prescribed initial velocity
\cite[Theorem~20.1 and following paragraph, p.~516]{LeeSmoothManifolds2e}.
The resulting representation $R_U$ of $U$ is unitary for the original
norm: $dR(\mathfrak u)$ is skew-adjoint, and the exponential is onto
$U$ by \cite[Corollary~4.48, p.~259]{KnappBeyondDigital}.
It is irreducible and nontrivial, because $\mathfrak u$ complexly spans
$\mathfrak g_{\mathbb C}$ and $G$ is connected.

The $2r$ real directions $Y_i,iY_i$ form a basis of the tangent algebra
of $G_{\mathbb C}$. Exponential coordinates and their derivative at zero
\cite[Proposition~20.8(a),(e),(f), pp.~519--520]{LeeSmoothManifolds2e}
show that their ordered word map contains an identity neighborhood.
Here too the source exponential chart is transported by the same
initial-velocity uniqueness, so it is a chart for this very exponential.
Connectedness and compactness give a bound $M$, independent of $R$,
on the number of factors needed to cover the image of $U$, with every
factor of the form $\exp(tY_i)$ or $\exp(itY_i)$, $|t|\le1$.
Padding each factor by zero parameters in a complete ordered basis word
gives a fixed list $I_1,\ldots,I_m$, where $m=rM$, such that every
$R_U(u)$ equals
\[
          \prod_{j=1}^m\exp(z_jdR(Y_{I_j})),\qquad |z_j|\le1.
\]
Let $Q\subset G$ be the compact set of the same ordered words with
real parameters $t_j\in[-1,1]$.
The word is polynomial of degree at most $D$ in every coordinate.
For any nonzero vector in the finite-dimensional representation space,
its polynomial coefficient span contains its $U$-orbit and is therefore
all of that space. Consequently $\dim R\le(D+1)^m$.
If $D=0$, every $dR(Y_i)$ vanishes and connectedness makes $R$ trivial.
Thus, for nontrivial $R$,
\begin{equation}\label{III-eq:ha-profile-weight-bounds}
 L_R\ge L_*:=1/C>0,\qquad
 \dim R\le C_0(1+L_R)^q,\qquad
 dR(Y_i)^{\lfloor CL_R\rfloor+1}=0,
 \quad C_0=(C+1)^m,\quad q=m.
\end{equation}

Lagrange interpolation gives, for a Hilbert-valued polynomial of
coordinate degree at most $D$, the real-interval to complex-disk bound
$(2e)^D$ in each variable. Applying this to an arbitrary orthogonal
projection $P$, including an arbitrary Hilbert multiplicity space,
gives, with $A=mC\log(2e)$,
\begin{equation}\label{III-eq:ha-profile-polynomial-continuation}
 \sup_{u\in U}\|(P\otimes I)(R_U(u)\otimes I)v\|
 \le e^{AL_R}\sup_{g\in Q}\|(P\otimes I)(R(g)\otimes I)v\|.
\end{equation}
Choose a unit $H\in\mathfrak p$ and a unit eigenvector $w$ of $dR(H)$
with eigenvalue $L_R$, reversing $H$ if necessary, and let $P_w$ be
the orthogonal projection onto $\mathbb Cw$. Unitary rank-one Schur
orthogonality gives
\[
 \int_U\|(P_w\otimes I)(R_U(u)\otimes I)v\|^2\,du
                         =\frac{\|v\|^2}{\dim R}.
\]
Indeed, the average of the conjugated rank-one projector is self-adjoint
and commutes with $R_U$. Its eigenspaces are invariant, so irreducibility
makes it scalar; its trace is one.
For arbitrary multiplicities the same identity follows by summing the
nonnegative coordinate identities. Thus some $g\in Q$ satisfies
\[
 \|(P_w\otimes I)R(g)v\|
              \ge e^{-AL_R}(\dim R)^{-1/2}\|v\|.
\]
Applying $R(\exp(TH))$ multiplies this component by $e^{TL_R}$, so
\[
 \log\frac{\|R(\exp(TH)g)v\|}{\|v\|}
          \ge(T-A)L_R-\tfrac12\log\dim R.
\]
Choose the fixed number
\[
 T=\max(A,0)+\tfrac12\bigl(\max(\log C_0,0)/L_*+q\bigr)+1.
\]
Using $\log(1+L_R)\le L_R$ and
\eqref{III-eq:ha-profile-weight-bounds}, the right side is at least $L_R$.
Hence the single compact set
\[
 Q'=\{\exp(TH)g:\ H\in\mathfrak p,\ \|H\|=1,\ g\in Q\}
\]
works for every $R$, every multiplicity space, and every nonzero $v$.
Finally the fixed compact set $\pi((Q')^{-1})\subset G/K$
has profile supremum at least one, because the profile uses $R(g^{-1})$.
The compact-open neighborhood of zero consisting of functions strictly
less than one on this compact set misses the whole profile family.
\end{proof}

\begin{corollary}[Comparable weight blocks]
\label{III-cor:ha-profile-comparable-blocks}
Fix $a>0$. Let $R$ be a direct sum of nontrivial irreducible types
whose scales satisfy $aL\le L_{R_j}\le L$, with arbitrary Hilbert
multiplicities. There is a compact $Q'_a\subset G$, independent
of $L$, the types and their multiplicities, such that
\[
          \sup_{g\in Q'_a}\log\frac{\|R(g)v\|}{\|v\|}\ge L
                       \qquad(v\ne0).
\]
A trivial summand may be included if its vector component has norm
at most that of the sum of the nontrivial components.
\end{corollary}
\begin{proof}
Let $Q$ denote the compact set supplied by
Proposition~\ref{III-prop:ha-profile-nonflattening}, so every
nonzero vector in every nontrivial irreducible type, with arbitrary
Hilbert multiplicity, has logarithmic gain at least $L_R$ somewhere on $Q$.
There is a constant $\lambda_*>0$ such that $L_R\ge\lambda_*$ for every
such type: in the nilpotent-basis argument of that proposition the allowed degree
is $\lfloor cL_R\rfloor$, and degree zero forces every basis differential
to vanish, hence forces the representation to be trivial.

The original global Cartan diffeomorphism
\cite[Theorem~6.31, pp.~362--367]{KnappBeyondDigital} shows that
$U=K\exp\{H\in\mathfrak p:\|H\|<1/2\}$ is open.
The same compact-exponential argument as above
\cite[Corollary~4.48, p.~259]{KnappBeyondDigital} makes $K$ unitary
for the given Hermitian norm, and
$\|R(u)\|\le e^{L_R/2}$ on $U$.
Cover $Q$ by finitely many right translates $Ug_i$, $1\le i\le k$,
taking $k\ge1$ by adjoining the identity if necessary. Write $g_i=k_i\exp H_i$ in the original Cartan decomposition.
For every nonzero vector $w$ in any one type with arbitrary multiplicity,
some $i$ therefore satisfies
\[
 \log\frac{\|R(\exp H_i)w\|}{\|w\|}\ge\frac{L_R}{2}.
\]
All these choices precede the choice of representation.

Each $R(\exp H_i)$ is self-adjoint in the same Hermitian norm.
Cauchy--Schwarz gives $\|Tw\|^2\le\|w\|\|T^2w\|$, so its
logarithmic gain is multiplied by at least $2^n$ on replacing
$\exp H_i$ by $\exp(2^nH_i)$. Put
$Q'_a=\{\exp(2^nH_i):1\le i\le k\}$.
For each isotypic component $v_j$, one member of this same finite set
has squared norm at least $e^{2^naL}\|v_j\|^2$.
Summing first over the finite set and then over the arbitrary Hilbert
direct sum gives
\[
 \max_{g\in Q'_a}\|R(g)v\|^2
 \ge k^{-1}e^{2^naL}\|v\|^2.
\]
Choose $n$ once and for all so that
$2^na/2\ge1+(\log k+\log2)/(2\lambda_*)$.
For a nonzero nontrivial component we have $L\ge\lambda_*$, hence
the logarithmic gain is at least $L+(\log2)/2$.
Finally a controlled trivial component enlarges the initial norm by
at most $\sqrt2$ and cannot decrease the output norm of the nontrivial
part. The same finite compact set therefore gives gain at least $L$
for the full vector.
\end{proof}

\section{Analytic realization with a transverse parameter}\label{III-series:parameters}\label{III-series:last-section}

\subsection{The parameter convergence theorem}

\begin{theorem}\label{III-parameter_uniform_convergence:thm:parameter}
Let $L=\SL(2,\R)$, $E=\R^2$, $V=E\oplus\R$, and
$\rho(l)=\operatorname{diag}(l,1)$. Let a lattice $\Lambda<L$
act on a probability space $(X,\mu)$ preserving $\mu$. Suppose
$\alpha(s,x;z;y)$ is a coefficient-measurable cocycle of actual
analytic germs in $(z,y)\in\R\times V$, with $z$ fixed under
the base action, such that, identically as analytic germs in $z$,
\[
 \alpha(s,x;z;0)=0,\qquad D_y\alpha(s,x;z;0)=\rho(s).
\]
Let
\[
 h_x(u,z)=(u,0)+\sum_{p\ge2}h_{p,x}(z)[u^p],\qquad
 h_{p,x}(z)\in\Hom(\Sym^p E,V)[[z]],
\]
be coefficient-measurable and formally equivariant:
\[
 \alpha(s,x;z;h_x(u,z))=h_{sx}(su,z).
\]
Then $h_x$ is jointly convergent in $(u,z)$ almost everywhere.
No property $(T)$, moment assumption, or prior convergence of
the individual series $h_{p,x}$ is required.
\end{theorem}

Here formal equivariance means equality of every finite total Taylor
truncation in $(z,u)$. If $H_{x,N}$ is the truncation of $h_x$
through total degree $N$, this is
\[
 j^N(z,\alpha(s,x;z;H_{x,N}(u,z)))
   =j^N(z,H_{sx,N}(su,z)).
\]
Multiplying the output by $\rho(s)^{-1}$ gives the near-identity
substitution equation for
$e_{s,x}(z;y)=\rho(s)^{-1}\alpha(s,x;z;y)$ and the transported
original coefficients. The germ hypotheses give
$e_{s,x}(z;0)=0$ and $D_y e_{s,x}(z;0)=I$ on a parameter
neighborhood. This normalization uses the same coefficients
throughout; the estimates below establish convergence of their
original mixed Taylor series.

We use Lemma~\ref{III-lem:ha-profile-invariant-law},
Proposition~\ref{III-prop:ha-profile-nonflattening}, and
Corollary~\ref{III-cor:ha-profile-comparable-blocks}.
The estimates below are uniform in the dimension of the coefficient
algebra.
The parameter is placed in that algebra, not counted in the
representation-theoretic degree.

\subsection{Uniform coefficient algebras}
For $m\ge0$, let
\[
 A_m=\C[z]/(z^{m+1}),\qquad
 \|a\|_{A_m}^2=\sum_{j=0}^m(j+1)^2|a_j|^2.
\]
There is $C_A\ge1$, independent of $m$, such that
\[
 \|ab\|_{A_m}\le C_A\|a\|_{A_m}\|b\|_{A_m}.
\]
Here is a version which also proves every vector-valued product
estimate used below. Suppose $a_i\in U$, $b_j\in W$ and
$\beta:U\times W\to Z$ is bilinear between Hilbert spaces with
$\|\beta(a,b)\|\le\|a\|\|b\|$. Set
$c_d=\sum_{i+j=d}\beta(a_i,b_j)$. Since
$d+1\le(i+1)+(j+1)$, Young's convolution inequality gives
\[
 \|(d+1)c_d\|_{\ell^2}
 \le\|(i+1)a_i\|_{\ell^2}\|b_j\|_{\ell^1}
       +\|a_i\|_{\ell^1}\|(j+1)b_j\|_{\ell^2}
 \le 2\left(\sum_{j\ge0}(j+1)^{-2}\right)^{1/2}
       \|a\|_{A_\infty}\|b\|_{A_\infty}.
\]
Taking $U=W=Z=\C$ gives the scalar assertion; taking angular
tensor multiplication or contraction gives the needed vector
versions with the same constant. Truncation is contractive.
Multiplication by the
constant $1$ is the identity operator and has norm one. This last
fact prevents a fixed algebra constant from being charged at each
linear substitution.

Give homogeneous polynomial spaces their symmetric-tensor Hilbert
norms. A polynomial matrix with entries in $A_m$ is given the
Hilbert norm whose square is the sum of the squared $A_m$ norms of its
scalar entries. This is \emph{not} the Hilbert--Schmidt norm of
the associated multiplication operators on $A_m$: the latter
choice would introduce an unnecessary factor depending on $m$.

\begin{lemma}[Uniform substitution length]\label{III-parameter_uniform_convergence:lem:length}
Suppose $e(z;y)=y+q(z;y)$, with $q$ of angular order at least two,
is jointly analytic near zero. A bounded set of inputs means a common
complex product neighborhood on which their holomorphic extensions have
a common supremum bound. After dilation
\[
 e^{[a,r]}(z;y)=a^{-1}e(rz;ay),
\]
on any fixed bounded set of such analytic inputs there are
$r_0,a_0,C>0$ such that, for $0<r\le r_0$, $0<a\le a_0$,
all normal truncations $m$ and angular truncations $d$ satisfy
\[
 \log\|C_{d,m}(e^{[a,r]})^{\pm1}\|_{\mathrm{op}}
 \le Cda.
\]
Here $C_{d,m}$ acts by substitution in $y$ on angular polynomials
of degree at most $d$ with coefficients in $A_m$.
The estimate also holds for right multiplication of polynomial
row vectors by its coefficient matrix, where the angular matrix
is transposed but coefficient multiplication is left unchanged.
The same estimate holds for finitely many fixed linear conjugates
of the input.
\end{lemma}
\begin{proof}
For a common complex product ball of radius $R$ and supremum bound $M$,
the order-two Schwarz estimate gives
\[
 \|a^{-1}q(rz;ay)\|\le (9M/R^2)a
 \qquad (|z|,\|y\|\le3,\quad 0<r,a\le R/6).
\]
Thus the scales and constants can be chosen before the input.
For finitely many fixed angular linear conjugacies, shrink the common
domain using their norms and inverse norms; the conjugated errors
again have a common supremum bound.
After shrinking $r_0,a_0$, joint analyticity makes
$q^{[a,r]}$ uniformly $O(a)$ on a fixed product complex domain
strictly larger than the unit product domain. The same is true
of its first derivative in $y$. The analytic inverse theorem,
or its contraction proof with a fixed margin between these
domains, applies uniformly to
$e_t=\id+tq^{[a,r]}$, $0\le t\le1$. Thus the actual analytic
generator $v_t=q^{[a,r]}\circ e_t^{-1}$ is $O(a)$ on a fixed
product domain still strictly larger than the unit domain.
Cauchy estimates therefore give
\[
 v_t(z;y)=\sum_{b\ge0}\sum_{k\ge2}z^b v_{t,b,k}[y^k],
 \qquad
 \sum_{b\ge0}(b+1)\sum_{k\ge2}\|v_{t,b,k}\|\le C_1a.
\]
All statements concern actual analytic inputs before truncation;
their algebraic truncations are consequently compatible for all
$m,d$.

For fixed $b$, the scalar-coefficient angular derivation
$D_{v_{t,b}}f=Df\,v_{t,b}$ sends degree $j$ to degree
$j+k-1$. Symmetric-tensor contraction costs at most
$j\|v_{t,b,k}\|$. Summing these degree shifts by the
$\ell^1$--$\ell^2$ inequality bounds its operator norm on
angular polynomials of degree at most $d$ by
$d\sum_{k\ge2}\|v_{t,b,k}\|$. Its angular transpose has the
same norm. On $A_m$, multiplication by $z^b$ has norm at most
$b+1$, because $(j+b+1)/(j+1)\le b+1$; truncation can only
decrease this norm. Summing over $b$ bounds both the column
derivation and the right row derivation by $C_1da$.

The differential equation for substitution along $e_t$ and
Gr\"onwall's inequality now bound both substitution actions and
their inverses by $\exp(C_1da)$. This argument treats the row
action explicitly: its angular matrices are transposed, while
normal multiplication is unchanged. Coefficient multiplication
on $A_m$ need not be self-adjoint, so the row action is not
identified with the adjoint of the full column action.
\end{proof}

\subsection{Rectangular matrices and profiles}
For fixed $a,r>0$ replace $h$ by
\[
 h^{[a,r]}(u,z)=a^{-1}h(au,rz)
\]
and make the same change in the input cocycle. Suppress the
dilations for now. Put $W_d=\C[V]_{\le d}$ and
$\overline W_d=\C[E]_{\le d}$. Define the rectangular matrix
with entries in $A_m$ by
\[
 M_{d,m}(h)f=j_u^d j_z^m(f\circ h).
\]
It sends $1$ to $1$, so is nonzero. The matrix representation is
\[
 \mathcal S_d(g)M=\overline\rho_d(g)M\rho_d(\rho(g))^{-1}.
\]
Here the polynomial representations use inverse pullback,
$\rho_d(T)f=f\circ T^{-1}$, whereas $C_{d,m}(e)$ uses
forward substitution $f\mapsto j_y^d j_z^m(f\circ e)$.
The $A_m$ factor is a trivial Hilbert multiplicity. The
infinitesimal norm on the noncompact directions is at most $C_Ln$,
independently of $m$. The exact right-multiplier identity is
\[
 M_{d,m}(h_{sx})=\mathcal S_d(s)M_{d,m}(h_x)C_{d,m}(e_{s,x}),
 \qquad e_{s,x}=\alpha(s,x)\rho(s)^{-1}.
\]
Right multiplication here means coefficient convolution in $A_m$.
Lemma~\ref{III-parameter_uniform_convergence:lem:length} gives the required logarithmic bounds on
this action and its inverse; no inverse of the rectangular
matrix is needed.

Define
\[
 b_{N,m,h}(gK)=\max_{2\le d\le N}
 \frac{\log\|\mathcal S_d(g^{-1})M_{d,m}(h)\|}{d-1},
 \qquad
 \Psi_{N,m,h}=\frac{b_{N,m,h}-b_{N,m,h}(o)}{2C_L}.
\]
These are normalized convex Lipschitz profiles in the fixed compact
target $\mathcal C$ of Section~\ref{III-series:profiles}. Their cocycle defect on any
fixed finite set of evaluations is $O(a)$ on a set of measure
tending to one as $r,a\to0$, uniformly in $N,m$. To see the
large-measure assertion, choose common analytic radii and bounds
for the finitely many coefficient-measurable analytic inputs
and their fixed conjugates outside a set of arbitrarily small
measure, then apply Lemma~\ref{III-parameter_uniform_convergence:lem:length}.

Consequently the invariant-law compactness argument gives, for
any fixed compact $Q\subset L$ and $\delta,\tau>0$, sufficiently
small deterministic $r,a>0$ such that
\begin{equation}\label{III-parameter_uniform_convergence:eq:uniform}
 \sup_{N\ge2,\,m\ge0}
 \mu\{x:\sup_{g\in Q}\Psi_{N,m,h_x^{[a,r]}}(gK)\ge\delta\}
 <\tau.
\end{equation}
For clarity, failure of uniformity would supply deterministic
$N_j,m_j$ and $a_j,r_j\downarrow0$ with a positive probability
lower bound. The associated profile laws have a weakly convergent
subsequence in the fixed compact target. The defect tends to zero
in probability for every fixed lattice element, so the limiting
law is invariant. The invariant-law lemma makes it $\delta_0$,
contradicting the lower bound for the closed event. This argument
does not depend on the dimensions of the multiplicity spaces.

\subsection{Record comparison and convergence}
Let $h_p^{(m)}\in\Hom(\Sym^p E,V)\otimes A_m$ denote the
normal truncation of the angular coefficient and put
\[
 R_{N,m}(h)=\max\left\{0,
        \max_{2\le p\le N}\frac{\log\|h_p^{(m)}\|}{p-1}
                        \right\}.
\]
There is a constant $C_0$, independent of $N,m,h$, such that
\begin{equation}\label{III-parameter_uniform_convergence:eq:record}
 R_{N,m}(h)\le b_{N,m,h}(o)\le R_{N,m}(h)+C_0.
\end{equation}
The coordinate-polynomial columns prove the lower bound. For
the upper bound, a contribution to angular degree $\ell$ from
an input polynomial of degree $j$ has $j$ coefficient factors
whose degrees sum to $\ell$. Their record bounds multiply to
$\exp((\ell-j)R_{N,m})$. Coefficient-algebra multiplication
costs at most $C_A^{j-1}$ (or a smaller factor if some
linear factors are the constant inclusion). The number of
ordered degree choices is $\binom{\ell-1}{j-1}$, angular
tensor multiplication is contractive, and the polynomial-space
dimensions depend only polynomially on $d$, since $\dim E=2$
and $\dim V=3$ are fixed. Thus the matrix norm is at most
$D^d\exp((d-1)R_{N,m})$ with $D$ independent of $m$,
proving \eqref{III-parameter_uniform_convergence:eq:record}.

For every $p\ge2$,
\[
 \Hom(\Sym^p\C^2,\C^2\oplus\C)
 \simeq\Sym^{p+1}\C^2\oplus\Sym^{p-1}\C^2\oplus\Sym^p\C^2.
\]
All three summands are nontrivial; their highest-weight scales
are bounded above and below by fixed positive multiples of $p$.
The comparable-block nonflattening theorem therefore applies,
uniformly in the trivial Hilbert multiplicity $A_m$. In
particular, for any fixed large $D_1$, one symmetric compact $Q\subset L$
satisfies
\[
 \sup_{g\in Q}\log\frac{\|g\cdot v\|}{\|v\|}\ge D_1p
\]
for every nonzero degree-$p$ coefficient vector, every $p\ge2$,
and every $m$. This follows by taking a fixed number of products
of the compact nonflattening set, or directly by increasing its
fixed Cartan boost.

If $R_{N,m}(h)>0$, choose a maximizing degree $p$. Detect its
boosted norm in the coordinate columns of $M_{p,m}$. In view of
\eqref{III-parameter_uniform_convergence:eq:record}, choosing $D_1>C_0+2C_L$ gives
\[
 \sup_{g\in Q}\Psi_{N,m,h}(gK)\ge1.
\]
Apply \eqref{III-parameter_uniform_convergence:eq:uniform} with $\delta=1$. Then
\[
 \sup_{N,m}\mu\{R_{N,m}(h_x^{[a,r]})>0\}<\tau.
\]
The events on the left are increasing in both $N$ and $m$,
because coefficient truncation is contractive and its Hilbert
norm increases with $m$. Hence, outside a set of measure at
most $\tau$, every angular coefficient satisfies
\[
 \|h_{p,x}^{[a,r]}\|_{A_\infty}\le1\qquad(p\ge2).
\]
Its normal series is consequently holomorphic on $|z|<1$,
and its norm at every smaller normal disk is bounded by a
constant depending only on that disk. Summing in angular
degree gives joint convergence for $\|u\|<1/2$, $|z|<1/2$,
after a harmless fixed adjustment for the polynomial norm.
Undo the dilations. The event that the original formal series
is jointly convergent therefore has measure at least $1-\tau$.
Since $\tau>0$ was arbitrary, this proves
Theorem~\ref{III-parameter_uniform_convergence:thm:parameter}.

\clearpage
\section*{Part IV. Global classification}
\hypertarget{part.iv}{}
\addcontentsline{toc}{part}{Part IV. Global classification}

\section{Introduction}
The aim of this series is to contribute to the Zimmer program
\cite{Zimmer}; Fisher's survey \cite{FisherSurvey} describes its broader
classification problems, algebraic and exotic examples, and conjectural
framework. Here we study the classification of real-analytic actions
of lattices in $\SL(n,\R)$ on compact $n$-manifolds. In this
dimension there are infinite actions with different topologies
and nontrivial gluing data. The classification problem is to
identify the possible pieces and the ways in which they fit
together. The present paper carries out this global step using
cocycle superrigidity and the local description from Parts I--III.

There are two familiar sources of actions of lattices in
$\SL(n,\R)$ on compact $n$-manifolds. Restricting an action of
the ambient Lie group gives homogeneous and cohomogeneity-one
examples. The action of $\SL(n,\Z)$ on the torus gives another
source; replacing finite orbits by projective hypersurfaces and
gluing compatible neighborhoods produces further examples.
The transverse gluing data are analytic and need not preserve
volume. Katok--Lewis \cite{KatokLewis} constructed exotic analytic
actions by blowing up finite torus orbits, including
volume-preserving examples. Fisher--Melnick
\cite[Section 3.2]{FisherMelnick} developed further constructions
using disks and tubes.

The classification problem addressed here was formulated by
Fisher--Melnick \cite[Conjecture 3.6]{FisherMelnick}. They proposed
three alternatives: finite image, extension to the ambient group
or its covering group, and actions assembled from torus and
ambient pieces. They also classified the ambient actions
\cite[Theorems 4.2 and 4.3]{FisherMelnick}, extending Uchida's
classification of analytic actions on spheres \cite{Uchida}.
Our main theorem gives the real-analytic version of their proposed
alternatives, with the finite-index, covering, and orbifold
qualifications stated explicitly.

Son \cite{Son} extended the real-analytic ambient-group
classification to closed $m$-manifolds with
$3\le n\le m\le2n-3$. In the case $m=n$, her induced flow
constructions agree with those of Fisher--Melnick
\cite[Section 4]{Son}. The present series recovers the relevant
neighborhood models from lattice actions and establishes their
global assembly, including the torus pieces and their gluings.

We prove that these constructions account for every real-analytic
action in the indicated dimension. Throughout, $M$ is compact,
connected, and without boundary. A finite-index restriction will
be stated when it is used. Finite covers of individual pieces do
not assert the existence of one finite cover of the glued manifold.

\subsection{Relation to earlier rigidity results}
The dimension of the manifold plays a central role in the Zimmer
program. For lattices in $\SL(n,\R)$, Brown--Fisher--Hurtado proved finiteness
below dimension $n-1$, treating uniform lattices, $\SL(n,\Z)$,
and general nonuniform lattices in
\cite{BFHCocompact,BFHSLZ,BFHGeneral}. At dimension $n-1$,
the boundary-action results of Brown--Rodriguez Hertz--Wang
\cite{BRHW} describe the projective models. Their work on
invariant measures and measurable projective factors
\cite{BRHWFactors} supplies another principal input to Part II.
The present problem is one dimension above the projective action;
the extra coordinate allows nontrivial analytic gluing.

There is also a substantial theory of global rigidity for
hyperbolic lattice actions, including Hurder \cite{Hurder},
Katok--Lewis \cite{KatokLewis}, and
Brown--Rodriguez Hertz--Wang \cite{BRHWGlobal}. Those results
identify algebraic actions under hyperbolicity and associated
structural hypotheses. Here the manifold is not assumed to be a
torus or nilmanifold, and no element is assumed to act as an
Anosov diffeomorphism. The punctured torus quotients and their point fillings
must instead be recovered from the local description.

\subsection{Classification}
Put $G=\SL(n,\R)$. In the statement, $\sigma$ denotes either the standard
representation $\sigma(g)=g$ or its dual $\sigma(g)=g^{-T}$.
The two punctured vector spaces are
\begin{equation}\label{IV-eq:intro-vector-models}
\begin{split}
 V_0&=\R^n\setminus\{0\}
       \simeq S^{n-1}\times(0,\infty),\\
 V_{\pm}&=(\R^n\setminus\{0\})/\{v\sim-v\}
       \simeq\RP^{n-1}\times(0,\infty).
\end{split}
\end{equation}
Their actions are $v\mapsto\sigma(\gamma)v$ and
$[v]\mapsto[\sigma(\gamma)v]$, respectively. Both spaces have
dimension $n$. They are linear $G$-orbits, not compact
$G$-tubes. The quotient in the second line identifies only
$v$ and $-v$, so it retains the radial coordinate.

We use \emph{$G$-tube} and \emph{$G$-disk} in the sense of
Fisher--Melnick \cite[Definitions 3.3 and 3.4]{FisherMelnick}.
A $G$-tube is the compact induced space
$G\times_{Q^0}[-1,1]\simeq S^{n-1}\times[-1,1]$, where
$Q^0$ is the identity component of a line stabilizer and acts
on the interval through a character and a vector field vanishing
at both endpoints. A $G$-disk is a closed $n$-disk with the
induced $G$-action away from its central fixed point and with
the linear $G$-action near that point. The standard and dual
versions are both allowed. The construction and the blow-up
patches used for gluing are described below.

The distinction between these pieces and their \emph{collars}
is useful. A boundary collar is a neighborhood
$N\times[0,\varepsilon)$ of a boundary hypersurface
$N=N\times\{0\}$. Thus a $G$-tube has two boundary collars,
whereas a $G$-disk has one boundary collar and a central fixed
point. After gluing, a two-sided hypersurface has a neighborhood
$N\times(-\varepsilon,\varepsilon)$; a nontrivial normal bundle
requires the normal-orientation cover and its deck involution.
The induced-action formula of Part II specifies the action as
germs along $N$, including its transverse vector field. Such a
collar does not by itself specify the whole incident piece.
Recovering that piece is the purpose of the continuation and
completion arguments below.

\begin{theorem}[Real-analytic classification in the sense of Fisher--Melnick]
\label{IV-main:thm:classification}
Let $\Gamma<\SL(n,\R)$ be any lattice, $n\ge3$, and let
$\alpha:\Gamma\to\operatorname{Diff}^{\omega}(M)$ act on a compact
connected real analytic $n$-manifold without boundary. After
restriction to a finite-index subgroup, at least one of the
following alternatives holds.
\begin{enumerate}[label=\textup{(\arabic*)}]
\item \textbf{Finite quotient.} The action factors through a finite quotient of $\Gamma$.
\item \textbf{Extension to the ambient group.} The action, with an equivariant
finite covering when needed, extends to an analytic action of a
connected Lie group locally isomorphic to $\SL(n,\R)$.
\item \textbf{Actions built from tori, $G$-disks, and $G$-tubes.}
The lattice is commensurable with a conjugate of $\SL(n,\Z)$.
The action is assembled using blow-ups and two-sided blow-ups,
with the finite quotients specified below. After cutting along
the gluing hypersurfaces, the pieces have the following description;
at least one piece is of type \textup{(a)}.
\begin{enumerate}[label=\textup{(\alph*)}]
\item \emph{Punctured tori and their finite quotients.} Each such open region is $F\setminus S$,
where $F$ is a compact affine orbifold, $S$ is a finite invariant
set containing its singular points, and $F$ has a finite
equivariant affine torus cover $T=\R^n/\Lambda$. On that cover the
action has the form
\[
 [x]\longmapsto[\sigma(\gamma)x+b_\gamma],
 \qquad \sigma(\gamma)\Lambda=\Lambda,
 \qquad b_{\gamma\delta}=b_\gamma+\sigma(\gamma)b_\delta
          \quad\text{in }T.
\]
\item \emph{$G$-tubes, $G$-disks, and their specified finite quotients.}
The open linear orbits are the spaces in
\eqref{IV-eq:intro-vector-models}, with radial ends $r\to0$
and $r\to\infty$. Compactifying both ends by boundary
hypersurfaces gives a $G$-tube on the spherical cover.
For $V_0$, filling the zero end by a linearizable fixed point
and compactifying the other end gives a $G$-disk.
The action on each piece is the restriction of its $G$-action.
\item \emph{Blow-up patches.} The gluing hypersurfaces are compact invariant
analytic hypersurfaces with the standard or dual projective
action on $S^{n-1}$ or $\RP^{n-1}$. Their neighborhoods, on the
angular and normal-orientation covers when necessary, have the
induced-action formula \eqref{II-eq:intro-collar}. For a nonzero
transverse field these are neighborhoods of the exceptional
divisor in the modified blow-up models described below.
The two local sides
join ends of the regions in \textup{(a)} and \textup{(b)}.
\end{enumerate}
The open regions have disjoint interiors. Together with the
gluing hypersurfaces and fixed points they cover $M$.
Every collar side is attached to exactly one incident end; the
two sides may belong to the same region. The incident
standard or dual types agree, and the aligned developing-coordinate
transitions are scalar homotheties. The gluing retains the full
analytic transverse vector field $B(t)\partial_t$, including any
normal-orientation deck involution. The torus coverings are
piecewise; a common finite cover of the glued manifold is not
part of the assertion.
\end{enumerate}
\end{theorem}
\begin{proof}
Theorem~\ref{II-main:thm:recognition} gives the local description
for $n\ge4$, and Theorem~\ref{III-main:thm:recognition-three}
gives it for $n=3$. Both statements apply to every lattice. Apply
Theorem~\ref{IV-global:thm:global} below. The arithmetic conclusion
in the torus case is Proposition~\ref{IV-general:prop:toral-arithmetic};
Proposition~\ref{IV-global:prop:ambient} gives the ambient extension
when no torus piece occurs.
\end{proof}

In alternative~\textup{(2)}, whenever the extension is to
$\SL(n,\R)$ itself, the ambient action is already completely
classified up to equivariant real-analytic diffeomorphism by
Fisher and Melnick \cite[Corollary~5.2]{FisherMelnick}.
Son gives a real-analytic classification in the wider dimension
range $3\le n\le m\le 2n-3$, which includes $m=n$
\cite[Theorems~1.1--1.2]{Son}.

The classification is a finite geometric construction together
with its analytic action and gluing data. It does not assert a
unique list of parameters for analytic conjugacy classes. The
alternatives need not be disjoint.

\subsection{Induced actions, blow-ups, and gluing}
Here are the local coordinates and the meaning of the fillings
in the theorem. On the spherical angular cover, the collar
formula is
\begin{equation}\label{IV-eq:intro-global-collar}
 \gamma(\theta,t)=
 \left(\frac{\sigma(\gamma)\theta}{\|\sigma(\gamma)\theta\|},
 \Phi^B_{\log\|\sigma(\gamma)\theta\|}(t)\right),\qquad B(0)=0,
\end{equation}
with the angular and normal-orientation covers described in
Part II. If $B$ is nonzero, each sufficiently small side develops
into a punctured vector model by
\[
 r(t)=\exp\!\int^t\frac{ds}{B(s)},\qquad (\theta,t)\longmapsto r(t)\theta.
\]
Here \emph{developing coordinates} are local coordinates in the
linear model that intertwine the given action with the standard
or dual linear action. Continuing them gives a developing map;
the term does not itself assert global injectivity. The constant
of integration in the displayed formula changes $r$ by a
positive scalar. A \emph{radial end} refers to the limiting
direction $r\to0$ or $r\to\infty$ in this model, before any
fixed point or boundary hypersurface has been added.
The limit at the central hypersurface is zero or infinity.
If $B=0$, the collar is a product. Scalar changes of developing
coordinates intertwine the standard action; the normal vector
field, and not just its first derivative, is retained in gluing.
This is the induced-action construction of Fisher--Melnick
\cite[Sections 3.1 and 5]{FisherMelnick}. To compare conventions,
let $Q$ stabilize $\R e_1$ in the standard representation.
Our character is
$\chi(q)=\log\|qe_1\|=-\log|\det(q|_{\R^n/\R e_1})|$.
Their character $\nu_0$ on $Q^0$ is $-\chi$.
Thus their generating vector field $X$ corresponds to
$-B(t)\partial_t$ in \eqref{IV-eq:intro-global-collar}.
The dual case is obtained by $g\mapsto g^{-T}$.

Here is the meaning of the blow-up terminology. Following
\cite[Section 3.2.1]{FisherMelnick}, write
\[
 B=\{(v,\ell)\in\R^n\times\RP^{n-1}:v\in\ell\},
 \qquad
 BS=\{(v,\theta)\in\R^n\times S^{n-1}:v\in\R\theta\}.
\]
The set $v=0$ is the \emph{exceptional divisor}, isomorphic
to $\RP^{n-1}$ and $S^{n-1}$, respectively. Replacing the
linear flow on a fiber by $X$ gives the modified models
$B_X$ and $BS_X$ of Fisher--Melnick. The first requires the fiber
involution to commute with $X$. The covering
$BS_X\to B_X$, when defined, has deck
map $(\theta,t)\mapsto(-\theta,-t)$.

For $X\not\equiv0$, a sufficiently small
neighborhood in \eqref{IV-eq:intro-global-collar} is a
neighborhood of the exceptional divisor in
$BS_X$, or its specified quotient. This includes
the quotient with deck map $(\theta,t)\mapsto(-\theta,t)$
when the normal bundle is trivial; it is not the same quotient
as $B_X$. These finite quotients are part of the
data of the theorem. For $X=0$, the neighborhood is a product,
as in Construction I of \cite{FisherMelnick}.

Replacing a puncture of a torus by $B_X$ is a
\emph{blow-up}. Joining two punctures through the two sides
of $BS_X$ is a \emph{two-sided blow-up}; the
punctures may lie on the same torus. The same patches attach
a $G$-disk at a puncture or join punctures through a $G$-tube
\cite[Section 3.2.2]{FisherMelnick}. A $G$-tube may contain
additional invariant hypersurfaces, corresponding to interior
zeros of its defining vector field. Cutting at all these
hypersurfaces gives the individual pieces in the theorem.

The fixed point of a $G$-disk fills the zero end of $V_0$.
Filling the zero end of $V_{\pm}$ would instead give
$\R^n/\{\pm I\}$, whose vertex is not a manifold point for
$n\ge3$. Such a filling is therefore not a $G$-disk.
The analytic fixed-point extension is proved in
Theorem~\ref{IV-exterior:thm:cap}.

For a punctured torus quotient, \emph{point filling} has a different purpose.
One replaces each zero end by its origin in the affine developing
coordinate, obtaining $F$. The local filling chart is either
$\R^n$ or $\R^n/\{\pm I\}$. This operation supplies a new affine
chart at the added point; it does not claim that $r(t)$ is analytic
at $t=0$ in the original manifold. The original action on $M$
is recovered by retaining the hypersurface and its collar.
The finite torus cover and the lifted affine action are obtained
together in Theorem~\ref{IV-global:thm:filling}. In particular,
the arithmetic assertion concerns the linear part and allows
the translation cocycle displayed in the theorem.

Some ambient actions are already closed models and require no
surgery decomposition. The Hopf manifold is
\[
 S^{n-1}\times\R/T\Z,\qquad
 \gamma(\theta,s)=
 \left(\frac{\sigma(\gamma)\theta}{\|\sigma(\gamma)\theta\|},
       s+\log\|\sigma(\gamma)\theta\|\right),\quad T>0;
\]
equivalently it is $(\R^n\setminus\{0\})/\langle e^T I\rangle$,
the Hopf manifold in \cite[Section 3.1.1]{FisherMelnick}.
The product action has the same angular action and fixes the
circle coordinate. The exceptional homogeneous models are
\[
 \operatorname{Flag}(\R^3)
 =\{(\ell,\Pi):\ell\subset\Pi\subset\R^3,\
                      \dim\ell=1,\ \dim\Pi=2\},
 \qquad
 \operatorname{Gr}_2(\R^4)=\{\Pi<\R^4:\dim\Pi=2\},
\]
with the natural actions and the finite covers specified in
Parts II and III. Thus the $3$-dimensional flag manifold belongs to
the ambient alternative. For the orbit classification see
the corrected statement and covering qualifications for orbits of
dimension at most $n$ in
Fisher--Melnick \cite[Theorem 2.1 and Remark 2.2, p.~5]{FisherMelnick}.
Cairns--Ghys \cite[Theorem 3.5, pp.~143--145]{CG} give the earlier
low-dimensional orbit argument; its statement does not include the
closed Grassmannian exception.
The ambient alternative also includes the nontransitive actions
assembled from vector regions and collars; the displayed closed
models are not a list of all ambient actions.

\subsection{Assembly from the local description}
For clarity we state the geometric implication separately.
Its hypothesis is Definition~\ref{II-def:recognition}, which is
proved for the present actions in Parts II and III.
\begin{theorem}[Classification by global assembly]
\label{IV-global:thm:global}
Let $\Gamma<\SL(n,\R)$ be any lattice, $n\ge3$, acting analytically
on a compact connected real analytic $n$-manifold without boundary.
Suppose the action has the local description property. Then the classification
in Theorem~\ref{IV-main:thm:classification} holds.
More precisely, its image is finite, it is one of the closed models in
Definition~\ref{II-def:recognition}, it is a product action
as in Proposition~\ref{IV-global:prop:product}, or it has a finite
compatible decomposition with the following pieces:
\begin{enumerate}
\item Open affine regions whose compact affine orbifold fillings
are finitely covered by affine tori with standard or dual linear part.
\item $G$-tubes and $G$-disks, with the finite quotients and open
linear orbits specified in Theorem~\ref{IV-main:thm:classification}.
\item The exceptional divisors and their neighborhoods used
in the blow-up and two-sided blow-up constructions.
\end{enumerate}
Distinct piece interiors are disjoint and every collar side belongs
to exactly one piece. The incident representation types agree.
Transitions between aligned developing coordinates are scalar
homotheties. The original analytic transverse flow is retained in
every gluing. The affine torus actions may have translation cocycles. If any torus
piece occurs, then $\Gamma$ is commensurable with a conjugate
of $\SL(n,\Z)$, as proved in
Proposition~\ref{IV-general:prop:toral-arithmetic}.
In every infinite-image outcome with no torus piece, the action extends to an ambient group
as in alternative \textup{(2)}, with the indicated finite covering
in the closed-model cases.
\end{theorem}

\subsection{The proof}
The first step is compatibility on overlaps. Hyperplane recurrence
forces different charts to carry a common projective coframe;
their transition maps are scalar affine maps. For $n=3$,
three inverse branches alone do not determine this coframe. We
use the polynomial relations of the entire recurrent family
and the associated orbit-compatibility calculation.

An affine region is completed by filling its zero ends. Fried's
theorem on closed similarity manifolds \cite[Theorem 2, p.~580]{Fried1980} reduces the
completion to a torus or Hopf alternative, and an invariant radial
field excludes the latter. A continuation with no identifications
instead closes into a $G$-tube or a $G$-disk. The local
models of Parts II and III detect the missing ends.

Finally, parabolic averaging on a compact invariant remainder
gives coverage. The comparison of incident neighborhoods gives
disjoint piece interiors and a finite decomposition. The local
ambient actions agree on overlaps whenever no torus piece occurs,
and compactness then integrates their infinitesimal action.

Analytic germ cocycle superrigidity is proved in Part I, and the local description
is in Parts II and III. We use those numbered results and public references;
all additional descent, continuation, and assembly arguments are
proved here. The proof does not assume an invariant smooth volume.

\section{Incidence, affine absorption, and radial fields}\label{IV-series:incidence}\label{IV-sec:basins}

The arguments of this paper apply for $n\ge3$. In the coframe
reconstruction in Section~\ref{IV-sec:descent}, the proof first treats
$n\ge4$; Section~\ref{IV-series:descent-three} proves the same two
descent statements for $n=3$. All later uses of those statements
include this replacement. The full-flag closed model is already
covered by Part III and needs no continuation argument.

\subsection{The germ data and their exactly equivariant domains}

Let $n\ge3$, $G=\SL(n,\R)$, and let $\Gamma<G$ be a lattice.
No arithmetic restriction is imposed on the lattice. Let $M$ be a real-analytic
$n$-manifold with an analytic $\Gamma$ action. Suppose an ergodic
probability-preserving $\Gamma$ space $(S,\eta)$ has an equivariant
center map $p:S\to M$ and measurable analytic coordinate germs
\begin{equation}\label{IV-basins:eq:germs}
 \Phi_s:(\R^n,0)\longrightarrow(M,p(s)),\qquad
 \gamma\Phi_s=\Phi_{\gamma s}\sigma(\gamma),
\end{equation}
where $\sigma$ is standard or dual. All identities initially mean
germ identities. A central sign cocycle can be removed on its
two-point probability extension; take an ergodic component, which
still projects to the original ergodic center measure. All the
conclusions below concern the actual fields on $M$ and descend.

Here measurability of germs is coefficient measurability in the convention
of Part~I, subsection ``Germs, measurability, and radii'', in measurable
centered original charts; no common representative domain is assumed.
In the standing compact setting, the measurable choice of these original
charts is supplied by the finite precompact action atlas constructed there.
Apply its countable integer-majorant construction to the analytic germs and
their inverses. This gives positive measurable forward and inverse radii and
canonical jointly measurable representatives. Shrink the forward radii
measurably into the selected original chart and the inverse convergence
domain. The inverse identities then give actual analytic local charts.
These representatives have the original germs, so their centers and germ
covariance are unchanged; neither the index probability nor its center
pushforward is replaced. We use these representatives in the continuation.

Analytic continuation gives exactly
equivariant analytic local diffeomorphisms
\begin{equation}\label{IV-basins:eq:domains}
 \widehat\Phi_s:\Omega_s\longrightarrow M,
 \qquad \Omega_{\gamma s}=\sigma(\gamma)\Omega_s.
\end{equation}
Each $\Omega_s$ is a union of centered ellipsoids and contains a
ball around zero. The extensions glue because intersections of
centered ellipsoids are connected and the analytic maps have the
same germ at zero. They need not be injective. We suppress the
hat below. Nothing here asserts that $\Omega_s=\R^n$.

\subsection{A zero--one law for incidence}

\begin{lemma}[Diagonal ergodicity]\label{IV-basins:lem:zeroone}
The diagonal action of $\Gamma$ on
$(S\times(\R^n\setminus\{0\}),\eta\times\Leb)$ is ergodic.
Here ergodic means that every invariant measurable set is null or
conull; the product measure is infinite.
\end{lemma}
\begin{proof}
For the standard representation form the induced probability space
\[
 Y=\Gamma\backslash(G\times S),\qquad
 \gamma(g,s)=(\gamma g,\gamma s).
\]
Normalize Haar measure so that a lattice fundamental domain has
mass one. The quotient measure is the pushforward of its product
with the original probability $\eta$. Its right $G$ action is
ergodic because the $\Gamma$ action on $S$ is ergodic.
Put $H=\{g:ge_1=e_1\}$. We prove the needed $H$ ergodicity
directly. A measurable section of the discrete quotient makes
the suspension action jointly measurable. Approximation by product
indicators on the finite-measure fundamental domain times $S$
shows that each Koopman orbit in $L^2(Y)$ is strongly measurable.
The Haar--Steinhaus argument then gives continuity of each orbit.
This uses neither a topology on $S$ nor separability of all of
$L^2(Y)$.

Write $u_{ij}(t)=I+tE_{ij}$. An $H$-fixed vector is fixed by
$u_{1j}(t)$ for every $j\ne1$. In the corresponding two-coordinate
block, for $c\ne0$ and $\epsilon\ne0$,
\[
 u_{1j}((c-1)/\epsilon)\,u_{j1}(\epsilon)\,
 u_{1j}((c^{-1}-1)/\epsilon)
 \longrightarrow d_{1j}(c)\qquad(\epsilon\longrightarrow0),
\]
where $d_{1j}(c)$ has diagonal entries $c,c^{-1}$ in those
coordinates and one elsewhere. Isometry of the Koopman operators
and orbit continuity show that the vector is fixed by $d_{1j}(c)$.
Since
$d_{1j}(c)^{-1}u_{j1}(t)d_{1j}(c)=u_{j1}(c^2t)$,
letting $c$ tend to zero proves invariance under the opposite
root group. These first-row and first-column root groups generate
$\SL(n,\R)$: their commutators give the other elementary matrices,
and elementary matrices generate the special linear group.
Consequently an $H$-fixed vector in $L^2_0(Y)$ is $G$-fixed and
vanishes. Applying this to centered indicators proves $H$ ergodicity.

For a diagonal invariant set $A$ the condition
$(s,ge_1)\in A$ defines an $H$-invariant set in $Y$.
It is null or conull. Unfolding the fundamental domain transfers
this conclusion to Haar measure times the same $\eta$ on $G\times S$.
For the transfer back, Fubini, Haar invariance, transitivity on
nonzero vectors and determinant-one invariance of Lebesgue measure
show that a set $B\subset\R^n\setminus\{0\}$ is Lebesgue-null
whenever $\{g:g^{-1}e_1\in B\}$ is Haar-null. Finite lattice
covolume gives unimodularity, so Haar inversion gives the same
implication with $ge_1$. Applying this one-way null transfer
fiberwise proves the conclusion for $A$. Countability of $\Gamma$
also allows an almost invariant measurable scalar function to be
replaced by an exactly invariant representative before descending
to $Y$. For the dual representation use $g\mapsto g^{-T}$ and
transport back to the same original $\Gamma$ action on $S$.
Removing the origin changes no Lebesgue null sets. Finally,
$\eta(S)=1$ and $\Leb(\R^n\setminus\{0\})=\infty$ give the
asserted infinite product mass.

For comparison, the Howe--Moore theorem gives another proof of
the stabilizer step: an $H$-fixed vector would have a constant
matrix coefficient along an unbounded sequence in $H$.
The connected noncompact simple real group with finite center
case is stated and proved in
\cite[Theorem 1.1 and Definition 2.6, pp.~66 and 69]{Ciobotaru};
the original theorem is \cite{HM}. The direct root-group argument
above does not require this external matrix-coefficient theorem.
\end{proof}

In particular the domain $\Omega\setminus(S\times\{0\})$ is
conull in this product. This is a measure statement, not a
pointwise assertion about every vector in every domain.

\begin{corollary}[Incidence with an invariant open set]
\label{IV-basins:cor:incidence}
Let $U\subset M$ be invariant and open. Set
\[
 A_U=\{(s,v):v\in\Omega_s\setminus\{0\},\ \Phi_s(v)\in U\}.
\]
If $(\eta\times\Leb)(A_U)>0$, then $A_U$ is conull in
$S\times(\R^n\setminus\{0\})$. For almost every $s$, the
complement of $U$ has Lebesgue measure zero in every sufficiently
small coordinate image of $\Phi_s$.

The positive-incidence hypothesis holds if a positive-measure
set of centers belongs to $\overline U$.
\end{corollary}
\begin{proof}
Exact equivariance makes $A_U$ invariant, so apply
Lemma~\ref{IV-basins:lem:zeroone} and Fubini. A local diffeomorphism preserves
the smooth measure class in coordinate neighborhoods. If
$p(s)\in\overline U$, its coordinate neighborhood meets the open
set $U$; the inverse image of that intersection contains a nonempty
open set, of positive Lebesgue measure. Fubini gives positive
incidence. No preassigned measure on $M$ is required.

An equivalent organization uses the essential-union construction in
Proposition~\ref{IV-basins:prop:essentialimage}, whose proof depends only
on Lemma~\ref{IV-basins:lem:zeroone}. It gives an invariant open set $E$
containing almost every original center, on which the evaluation
pushforward and the smooth measure class have the same null sets.
Choose a finite representative $\theta$ of that class. Its action is
nonsingular and ergodic; invariance of this finite representative is
not asserted. Positive source incidence forces $U\cap E\ne\varnothing$.
The center-in-closure condition also forces this intersection, since
almost every center lies in the open set $E$. An open nonempty
intersection has positive smooth-class measure, so invariance and
ergodicity make $U$ conull on $E$. Together with the already proved
conullness of $\Omega$, absolute continuity of evaluation gives
conull source incidence. Removing zero changes none of
these null or positive-measure assertions, and Fubini gives the
conclusion for the same original indices. Finally, transport null
sets through the actual analytic local inverse branches of each
$\Phi_s$; this gives the assertion in every sufficiently small image
chart. The auxiliary measure class is constructed from the given
germ family; the original probability $\eta$ and action are retained.
\end{proof}

\subsubsection{What already descends without a connection}

\begin{proposition}[Open ergodic measure classes]
\label{IV-basins:prop:essentialimage}
After discarding a null set of indices, the union
\[
 W=\bigcup_{s\in S}\Phi_s(\Omega_s)
\]
is an invariant open set containing almost every center, and its
smooth measure class is ergodic under $\Gamma$. Two such open
sets, obtained from possibly different ergodic germ families,
have the same smooth measure class on their union whenever they
intersect: each is conull in that union.
\end{proposition}
\begin{proof}
We first make the union an essential union. Choose a countable
basis of relatively compact coordinate balls in $M$. For each
basis member $V$, use compact rational source buffers $K$ and
the event
\[
 K\subset\Omega_s,\qquad \overline V\subset\Phi_s(K).
\]
The compact-containment parameter regions are measurable: the
transported domains have measurable radii, and containment can
be tested using finite rational subcovers. On each such region,
joint measurability and continuity in the source make
$s\mapsto\Phi_s|_K$ measurable in $C(K,M)$ with its compact-open
topology. The image-containment condition is closed there.
Indeed, for fixed $y$, the condition $y\in f(K)$ is closed by
compactness of $K$ and continuity of evaluation; intersect these
closed conditions over $y\in\overline V$. This proves event
measurability directly. Actual local inverse sheets ensure that
these countably many events cover every point of each coordinate
image after shrinking $V$; injectivity on an entire buffer is
not needed. Remove
every null such set of indices, and all their $\Gamma$ translates.
For the remaining indices, every basis member contained in a
coordinate image is contained in such images for a
positive-measure set of indices.

Consequently a measurable subset of $W$ is Lebesgue-null if and
only if its evaluation preimage in $\Omega$ is
$(\eta\times\Leb)$-null. In one direction use local
diffeomorphisms and Fubini. In the other, a set of positive
Lebesgue measure has positive intersection with some basis
member in the essential union; all the positive-measure indices
covering that member have positive-measure inverse image, by
null-set transport through their actual local inverse sheets.
Lemma~\ref{IV-basins:lem:zeroone} proves ergodicity of the smooth measure
class on $W$.

For two such sets, their nonempty open intersection has positive
Lebesgue measure in each. It is invariant, so is conull in each.
This proves the last assertion.
\end{proof}

This proposition constructs an open ergodic \emph{measure class}.
It does not construct a finite invariant measure in that class,
or assert that the analytic charts have affine transitions.
The evaluation pushforward of the invariant infinite product
measure can be infinite on every nonempty open set. Even the
covering map $\R^n\to\R^n/\Z^n$ has this property. Replacing
Lebesgue measure in the source by a Gaussian makes a finite
pushforward but destroys the required invariance.

\subsection{An existing affine region absorbs incident germs}

We record the two tensor facts that will be used. For the standard
module $V=\R^n$ and its vector stabilizer $H$,
\begin{equation}\label{IV-basins:eq:fixedspaces}
 V^H=\R e_1,
 \qquad \operatorname{Hom}(\operatorname{Sym}^2V,V)^H=0.
\end{equation}
For the second assertion, an
$H$-equivariant quadratic map has the form $av+be_1$ off
$\R e_1$, because $H$ is transitive there and the stabilizer of
$e_2$ fixes exactly $\operatorname{span}(e_1,e_2)$. Homogeneity
forces $a=b=0$.

We also use its recurrence lemma: a finite measurable equivariant
tensor $T(s,v)$ on an invariant product domain induces
\[
 F([g,s])=\rho(g)^{-1}T(s,ge_1)
\]
on a finite-measure $H$-invariant subset of $Y$. Recurrence under
the elementary unipotents of $H$ forces $F$ to be $H$-fixed.
A nonfixed unipotent orbit is an unbounded vector polynomial,
which cannot recur to a bounded level set. No moment assumption
on the tensor is needed.

\begin{theorem}[Affine absorption without ambient volume]
\label{IV-basins:thm:absorption}
Let $U$ be invariant and open, and suppose it carries a
$\Gamma$-invariant analytic torsion-free connection $\nabla$.
If the germ family \eqref{IV-basins:eq:germs} has positive incidence with
$U$, then for almost every $s$,
\begin{equation}\label{IV-basins:eq:affine}
 \Phi_s^*\nabla=\nstd
 \quad\hbox{on }\Phi_s^{-1}(U).
\end{equation}
In particular $U$ is conull in a small neighborhood of $p(s)$,
and its connection there extends as the flat analytic connection
defined by $\Phi_s$. These extensions agree on all overlaps.
\end{theorem}
\begin{proof}
On $A_U$ the difference
\[
 T(s,v)=(\Phi_s^*\nabla-\nstd)_v
\]
is an actual analytic tensor in $v$, measurable in $s$, with
values in $\operatorname{Hom}(\operatorname{Sym}^2V,V)$.
Exact equivariance gives the algebraic covariance from
\eqref{IV-basins:eq:fixedspaces}. Tensor recurrence makes $T=0$ almost
everywhere on $A_U$. For almost every $s$, analyticity (indeed
continuity suffices here) makes it zero everywhere on the open
set $\Phi_s^{-1}(U)$.

For the dual model, transport this calculation through
$g\mapsto g^{-T}$. Transport the lattice action on $S$ by the
inverse group isomorphism, and push forward Haar measure and its
fundamental domain together. The finite domain mass is unchanged;
translating back retains the original index probability and
original action on $M$.

Corollary~\ref{IV-basins:cor:incidence} says this open set is conull in each
small source ball. Push forward the standard connection on that
ball. It agrees with $\nabla$ on its entire intersection with
$U$. On the overlap of any two resulting coordinate images, $U$
is conull, so the two analytic connections agree there by
continuity. Each retained image is a neighborhood of its original
center, since its source radius is positive and zero maps to that
center. A countable subcover handles their union. The constructed
connection is analytic and flat on this union and agrees with
$\nabla$ on its overlap with $U$. In particular, there is no pullback of a connection
defined only modulo the singular center measure.
\end{proof}

\begin{definition*}
For an invariant flat analytic connection on $U$, let
$U^{\rm reg}$ be the union of open neighborhoods $V\subset M$
such that $U$ is Lebesgue-conull in $V$ and $\nabla|_{U\cap V}$
extends to an analytic connection on $V$.
\end{definition*}

\noindent
The definition uses only the smooth measure class of $M$.
Its extensions are unique and flat, agree on overlaps, and make
$U^{\rm reg}$ invariant. Also $U$ is conull in $U^{\rm reg}$.
Repeating the construction does not enlarge it.

\noindent\emph{Analytic flat coordinates.}
We use analytic affine coordinates for the torsion-free flat extensions just
constructed. The smooth affine-flatness criterion is
\cite[Chapter~V, \S4, Theorem~4.1, p.~210]{KobayashiNomizuI}; the analytic
regularity follows from the following local argument. Center an ordinary
analytic chart at the chosen point and write the connection coefficients as
$\Gamma_x(w,v)$, with the differentiated field first and the direction second.
Put $\Gamma_i(x)w=\Gamma_x(w,e_i)$. On the product of the coordinate domain
with the matrix space consider the analytic fields
\[
 H_i(x,A)=(e_i,-\Gamma_i(x)A),\qquad (x,A)=(0,I)\text{ initially}.
\]
Their brackets vanish by the curvature identity. Local flows of a fixed
analytic vector field are jointly analytic in time and initial point
\cite[Appendix~B, Theorem~B.3(ii), p.~333, and the analytic proof on
pp.~335--336]{DuistermaatKolk}.
For commuting fields $V,W$, the vectors
$D\Phi_V^t(x)W(x)$ and $W(\Phi_V^t(x))$ satisfy the same linear variational
ODE with the same initial value, and hence agree. After shrinking a common
neighborhood, finite composition of the flows of the $H_i$ therefore gives
an analytic integral map whose first component is $x$ and whose second
component $A(x)$ satisfies
\[
 A(0)=I,\qquad DA_x(v)w=-\Gamma_x(A(x)w,v).
\]
Torsion-freeness makes $\Gamma_x$ symmetric. The fields $X_i(x)=A(x)e_i$
consequently commute. Applying the same finite-flow construction to these
fields gives an analytic map $g$ with
\[
 g(0)=0,\qquad Dg_x=A(g(x)),\qquad Dg_0=I.
\]
The analytic inverse function theorem supplies its actual analytic local
inverse. Differentiating the middle identity gives
\[
 D^2g_x(v,w)+\Gamma_{g(x)}(Dg_xw,Dg_xv)=0.
\]
The connection transformation formula now gives zero coefficients in the
composite chart on the original manifold. Restrict its source to a positive
ball about zero, retaining the original connection and point and a target
inside the prescribed open region. The empty finite composition covers
coordinate dimension zero. Thus the analytic assertion, including the
inverse and the prescribed local domain, is obtained without assuming an
analytic upgrade of the smooth criterion.

\begin{corollary}[No noncompact germ measure on a maximal affine frontier]
\label{IV-basins:cor:frontier}
No ergodic invariant probability supported on
$\partial U^{\rm reg}$ can admit the standard or dual analytic
germs \eqref{IV-basins:eq:germs}.
\end{corollary}
\begin{proof}
Every coordinate neighborhood of a frontier point meets
$U^{\rm reg}$ and hence $U$. Theorem~\ref{IV-basins:thm:absorption} would
put almost every center in $U^{\rm reg}$, a contradiction.
\end{proof}

\begin{corollary}[Reconciliation of established affine basins]
\label{IV-basins:cor:reconcile}
Suppose two invariant open regions are essential images of
ergodic germ families as in Proposition~\ref{IV-basins:prop:essentialimage},
and each has an invariant flat analytic connection for which its
generating charts are affine. If the regions intersect, their
connections agree on the intersection and their maximal regular
extensions coincide. Thus the distinct maximal regions obtained
in this way are disjoint.

The same conclusion holds if either generating family is instead
an equivariant analytic local diffeomorphism from a nonempty
invariant open subset of the standard or dual punctured vector
space; no center at zero is required for this alternative.
\end{corollary}
\begin{proof}
Pull back the second connection by the first family's charts and
apply Theorem~\ref{IV-basins:thm:absorption}. On the intersection it agrees
with the first connection. Proposition~\ref{IV-basins:prop:essentialimage}
makes each region conull in their union. The definitions of the
maximal regular extensions therefore coincide.
For the alternative source take $S$ to be a point. Diagonal
ergodicity and tensor recurrence still apply on its invariant
source domain; neither argument uses a center at zero.
\end{proof}

The hypothesis that the connection already exists on each region
is substantive. It is satisfied by an embedded standard vector
model and by an affine region already produced with invariant
volume. It is not supplied just by a collection of convergent
germs over a singular measure.

\subsection{A radial field detects finite orbits at the frontier}

\begin{theorem}[Radial-field detection]\label{IV-basins:thm:radial}
Let $U\subset M$ be invariant and open, carrying a
$\Gamma$-invariant analytic vector field $E$ with no zeros.
If the ergodic germ family \eqref{IV-basins:eq:germs} has positive incidence
with $U$, then its center measure $\nu=p_*\eta$ is supported on
a finite orbit. On the incident charts there is a constant
$c\ne0$ such that
\begin{equation}\label{IV-basins:eq:radial}
 \Phi_s^*E(v)=c v\qquad (\Phi_s(v)\in U)
\end{equation}
for almost every $s$. The field extends analytically across their
centers, which are isolated zeros. No connection or volume on
$M$ is required.
\end{theorem}
\begin{proof}
The vector-valued tensor $T(s,v)=\Phi_s^*E(v)$ is equivariant
on $A_U$. Induction and unipotent recurrence, now using
$V^H=\R e_1$, show that
\[
 T(s,v)=a(s,v)v
\]
almost everywhere. The scalar $a$ is diagonal $\Gamma$-invariant.
Since $A_U$ is conull, Lemma~\ref{IV-basins:lem:zeroone} makes it a single
constant $c$. This constant is nonzero: $E$ is nowhere zero on
$U$, and $v\ne0$. Analyticity gives \eqref{IV-basins:eq:radial} everywhere
on each good chart's open preimage of $U$.

Push forward $c v$ on a small entire source ball. This defines
an analytic extension of $E$ near the center. Two extensions
agree on their conull intersection with $U$, so agree throughout
their overlap. They give an analytic vector field on an open set
containing almost every center. In each such coordinate ball its
only zero is the center, and its derivative there is $cI$.
Consequently these centers belong to a discrete, hence countable,
zero set. The probability $\nu$ is therefore atomic. Invariance
and ergodicity make it uniform on a finite orbit.
\end{proof}

If $\nabla$ also exists and $\nabla E=\id$ on $U$,
Theorem~\ref{IV-basins:thm:absorption} and \eqref{IV-basins:eq:radial} give $c=1$.
This is the normalization in the standard tube.

There is a useful exclusion. If an invariant open
region has a nonconstant invariant analytic function whose
differential never vanishes, no standard or dual probability
germ family can have positive incidence with it. Its pullback
would be a constant by Lemma~\ref{IV-basins:lem:zeroone}, contradicting
the local diffeomorphism property. In particular this applies
to an invariant product projective collar with its transverse
coordinate. This does not treat a moving collar coordinate as
an invariant function.

\section{Hyperplane recurrence and affine descent}\label{IV-series:descent}\label{IV-sec:descent}

\subsection{Statements and exact domains}

Put $V=\R^n$, $G=\SL(n,\R)$, and let $\Gamma<G$ be a lattice.
Let $M$ be a real-analytic $n$-manifold with an analytic $\Gamma$ action.
Throughout, $n\ge4$. Write $\sigma$ for the standard representation
or its dual. Applying $g\mapsto g^{-T}$ reduces all proofs to the
standard case; all source vectors and hyperplanes use that choice.

Let $(S,\eta)$ be an ergodic probability-preserving $\Gamma$ space,
with equivariant center map $p:S\to M$. Suppose there is a measurable
family of analytic local diffeomorphisms
\begin{equation}\label{IV-descent:eq:family}
 \Phi_s:\Omega_s\longrightarrow M,\qquad \Phi_s(0)=p(s),\qquad
 \Phi_{\gamma s}(\sigma(\gamma)v)=\gamma\Phi_s(v),
 \quad \Omega_{\gamma s}=\sigma(\gamma)\Omega_s.
\end{equation}
Each $\Omega_s$ is a union of centered ellipsoids containing zero.
These exact domains are obtained from the germ identities by the
continuation lemma in Section~\ref{IV-sec:basins}. Intersections of the
ellipsoids are connected, so the germ identities glue. Global
injectivity of $\Phi_s$ is not assumed. Measurability means measurable
analytic coefficients and measurable domains, equivalently a measurable
coding by restrictions to countably many compact source domains.

\begin{theorem}[Nonatomic affine descent]\label{IV-descent:thm:nonatomic}
If $p_*\eta$ is nonatomic, discard an invariant null set of indices
in \eqref{IV-descent:eq:family}. The essential union
\[
 U=\bigcup_s\Phi_s(\Omega_s)
\]
has an analytic atlas whose connected transition maps are
\begin{equation}\label{IV-descent:eq:homotheties}
 v\longmapsto av+b,\qquad a\in\R^*,\quad b\in V.
\end{equation}
Every local inverse branch of every remaining $\Phi_s$ belongs to
this atlas. Thus the standard flat torsion-free connections descend,
and the inverse differentials define a global analytic projective
coframe
\[
 [A_y]\in\Pj\operatorname{Iso}(T_yM,V),\qquad
 [A_{\gamma y}D_y\gamma]=[\sigma(\gamma)A_y].
\]
The set $U$ contains almost every center. Its smooth measure class
is ergodic, and it has finitely many connected components, permuted
transitively by $\Gamma$.
\end{theorem}

Here a projective coframe is an invertible linear coframe up to a
nonzero scalar. It is stronger data than an affine connection. In
particular \eqref{IV-descent:eq:homotheties} is a similarity atlas, even though
the lattice acts by affine maps whose linear parts need not be
similarities in these coordinates.

\begin{theorem}[Continuation from a fixed point or an interior collar]\label{IV-descent:thm:single}
Let $\Omega\subset V\setminus\{0\}$ be a nonempty invariant open
union of punctured centered ellipsoids, and suppose it contains a
punctured ball. Let
\[
 F:\Omega\longrightarrow M,\qquad F(\sigma(\gamma)v)=\gamma F(v)
\]
be an analytic local diffeomorphism. Assume that on some punctured
ball it is injective, or has exactly the identifications $v\sim-v$,
with $F(-v)=F(v)$ throughout $\Omega$ in the latter case.
Then exactly one of the following holds:
\begin{enumerate}
\item $F$ is injective on $\Omega$, or injective on $\Omega/\{\pm1\}$
      in the prescribed antipodal case.
\item All local inverse branches of $F$ have transitions
      \eqref{IV-descent:eq:homotheties}. At almost every image point the vectors
      in its full inverse image span $V$.
\end{enumerate}
In the second case the image has the flat connection and equivariant
projective coframe of Theorem~\ref{IV-descent:thm:nonatomic}, and its smooth
measure class is ergodic. If $F$ also extends as an injective germ at
zero, its chart at zero is included in the descended atlas.
\end{theorem}

The image in the second alternative is called an affine region below.
The first alternative does not by itself assert $\Omega=V\setminus0$.

\subsection{The recurrence facts}

We use two finite-measure arguments. Let
\[
 H=\{g\in G:ge_1=e_1\}
   =\left\{\begin{pmatrix}1&r\\0&B\end{pmatrix}:B\in\SL(n-1,\R)\right\}.
\]
It is generated by one-parameter unipotent groups. Its fixed vectors
are $\R e_1$, its commutant in $\End(V)$ consists of scalars, and
its invariant subspaces containing $e_1$ are $\R e_1$ and $V$.
For the last assertion pass to $V/\R e_1$, the irreducible standard
module of the lower-right Levi group. For the commutant, commuting
with that Levi group gives $\operatorname{diag}(a,bI)$; commuting
with the top-row unipotents gives $a=b$.

\begin{lemma}[Projective unipotent recurrence]\label{IV-descent:lem:unipotent}
Let a one-parameter unipotent linear group act on a finite-dimensional
real vector space $E$. Every invariant probability on $\Pj(E)$ is
supported on its fixed lines. Consequently, an equivariant measurable
map from a probability-preserving space to $\Pj(E)$ takes values in
those fixed lines almost everywhere.
\end{lemma}
\begin{proof}
For a nonfixed line $[w]$, the polynomial $\exp(nN)w$ has positive
degree and its projective limit is its highest nonzero coefficient.
That fixed line differs from $[w]$. Thus $[w]$ is not recurrent under
the time-one map. Poincare recurrence proves the assertion. A line
invariant under a unipotent group is fixed pointwise in $E$, since
the only eigenvalue is one. Apply this to finitely many generating
unipotent subgroups when several groups are present.
\end{proof}

The same argument on the Plucker embedding applies to Grassmannians.
In particular there is no $G$-invariant probability on
$\operatorname{Gr}_k(V)$ for $0<k<n$: such a probability would be
supported on subspaces invariant under all elementary unipotents,
contrary to irreducibility.

\begin{lemma}[Diagonal zero--one law]\label{IV-descent:lem:zeroone}
The diagonal $\Gamma$ action on
$(S\times(V\setminus0),\eta\times\Leb)$ is ergodic.
In particular the domain in \eqref{IV-descent:eq:family} is conull there.
For $S$ a point, every nonempty invariant open subset of
$V\setminus0$ is conull and has ergodic Lebesgue measure class.
\end{lemma}
\begin{proof}
Lemma~\ref{IV-basins:lem:zeroone} gives ergodicity for the standard and
dual diagonal actions. The domain in \eqref{IV-descent:eq:family} is
invariant, and its sections contain punctured balls of positive radius.
Fubini gives positive product measure, so the domain is conull.
For a point index, a nonempty open invariant set has positive Lebesgue
measure and hence is conull. Restricting to that set preserves ergodicity.
\end{proof}

\subsection{Recurrence forbids transverse hyperplane intersections}

In this section a hyperplane always means a linear hyperplane through
zero in $V$. For two coordinate branches with a common image point,
write $v,w$ for its source vectors.

\begin{lemma}[Pairwise plaque constraint]\label{IV-descent:lem:plaques}
For the family in Theorem~\ref{IV-descent:thm:nonatomic}, there is an invariant
conull family of indices such that, for every pair of retained charts,
every hyperplane $L$, and every common image point with $v,w\in L$,
\begin{equation}\label{IV-descent:eq:plaque}
 D_v\Phi_s(L)=D_w\Phi_t(L).
\end{equation}
For the single map of Theorem~\ref{IV-descent:thm:single}, the same assertion
holds for every pair of its inverse branches.
\end{lemma}
\begin{proof}
First consider the probability family. Induce the diagonal action on
$S^2$ to $Y_2=\Gamma\backslash(G\times S^2)$; no ergodicity of
this probability space is needed. Almost every pair of centers is
distinct. At $y=[g,s,t]$ the maps
$v\mapsto\Phi_s(gv)$ and $w\mapsto\Phi_t(gw)$ have disjoint
images on two sufficiently small balls of the same radius
$\delta(y)>0$. Choose $\delta$ measurably, for example as a dyadic
radius satisfying these open-domain and disjointness conditions.
The condition is independent of the representative of $y$.

Fix a hyperplane $E$ and a diagonalizable $a\in G$ acting on $E$
by $\lambda I$, $\lambda>1$. If these two maps have a common
image at some $v,w\in E$, then at $ya^d$ they have a common image
at $a^{-d}v,a^{-d}w$, which tend to zero. It follows that
$\delta(ya^d)\to0$. Almost every point in each set
$\{\delta\ge1/k\}$ returns to that set infinitely often. Thus
for almost every $y$ there is no such pair in $E\times E$.
Fubini gives this assertion for almost every $(s,t)$ and Haar
almost every $g$.

If two plaques in some hyperplane $L$ met transversely, their tangent
hyperplanes would span $T M$. The implicit function theorem, applied
to the equation $\Phi_s(v)=\Phi_t(w)$ with $v,w\in L$, would
preserve an intersection for every nearby hyperplane. This gives a
positive-Haar set of $g$ with $gE$ admitting an intersection,
contradicting the preceding paragraph. Thus
\eqref{IV-descent:eq:plaque} holds for almost every pair of indices, for all
source points and all hyperplanes.

We now make the assertion simultaneous in the indices. Use countably
many compact buffered source domains, consisting of finite unions of
rational balls and lying in $\Omega_s$, and the restrictions of the
maps to them in the separable $C^1$ topology. Such domains can contain
zero and any prescribed finite collection of source points, since
$\Omega_s$ is star shaped. For each domain discard the null set of
restrictions outside the support of its conditional distribution.
Discard all $\Gamma$ translates of these null sets. A transverse
intersection is open under $C^1$ perturbations of both maps on the
relevant buffered domains. Were it present for two retained maps,
their support neighborhoods would produce a positive-measure set of
pairs with a transverse intersection, a contradiction. This also
treats two branches of the same retained map.

For a single $F$, work on $\Gamma\backslash G$. There is a measurable
$\delta([g])>0$ on whose punctured ball $F(g\,\cdot)$ has no
identifications except the prescribed equality, or equality up to
sign. Any pair not related in that way and lying in $E\times E$
contracts to a forbidden pair, giving the identical recurrence
contradiction. The prescribed identifications have equal tangent
hyperplanes already. Openness of transversality in $gE$ proves the
assertion everywhere; no support selection is needed.
\end{proof}

\subsection{Three vectors determine a projective coframe}

At a common image point $y$, write inverse coordinate branches as
$v_i(y)$ and put $C_i=D_yv_i:T_yM\to V$. Equation
\eqref{IV-descent:eq:plaque} says
\begin{equation}\label{IV-descent:eq:covectors}
 \ell(v_i)=\ell(v_j)=0
 \quad\Longrightarrow\quad
 [\ell C_i]=[\ell C_j].
\end{equation}
The rows here are nonzero when $\ell\ne0$, since each $C_i$ is
invertible.

\begin{lemma}[Linear reconstruction]\label{IV-descent:lem:linear}
Suppose $n\ge4$ and a collection $(v_i,C_i)$, with each $C_i$ invertible,
satisfies \eqref{IV-descent:eq:covectors}. Set
\[
 \mathcal L_i=\R C_i+\{v_i\otimes\beta:\beta\in T_y^*M\}.
\]
If $v_1,v_2,v_3$ are independent, there is a unique line $[A]$
in $\mathcal L_1\cap\mathcal L_2$, and it belongs to every
$\mathcal L_i$. More precisely,
\begin{equation}\label{IV-descent:eq:rankone}
 C_i=a_i A+v_i\otimes\beta_i,\qquad a_i\ne0.
\end{equation}
The matrix $A$ has rank at least $n-1$; its invertibility is not
needed for this conclusion. If one branch has $v_j=0$, then
$[A]=[C_j]$ and $A$ is invertible.
\end{lemma}
\begin{proof}
For independent $v_i,v_j$, the annihilator of their span has
dimension $n-2\ge2$. By \eqref{IV-descent:eq:covectors}, every row in that
space is an eigenvector of $C_jC_i^{-1}$. Linearity makes all these
eigenvalues the same nonzero scalar. Therefore
\begin{equation}\label{IV-descent:eq:pair}
 C_j=a_{ij}C_i+v_i\otimes\beta+v_j\otimes\gamma,
 \qquad a_{ij}\ne0.
\end{equation}
If one vector is zero, the same argument uses the annihilator
of the other vector and omits the zero rank-one term.
This supplies a nonzero element of
$\mathcal L_i\cap\mathcal L_j$. The intersection is a line:
after eliminating the coefficient of $C_i$ from two elements,
a nonzero multiple of the invertible $C_j$ would have image in
$\operatorname{span}(v_i,v_j)$, impossible for $n>2$.
The remaining rank-one terms with independent image lines must
then vanish. A nonzero element of this intersection has a nonzero
$C_i$ coefficient, again by independence of the image lines.

Choose a generator $A$ for the first two branches and write
$C_i=a_iA+v_i\otimes\beta_i$ for $i=1,2$.
Applying \eqref{IV-descent:eq:pair} to $(1,3)$ and $(2,3)$ and subtracting
gives
\[
 (b a_1-c a_2)A\quad\hbox{with image in}
 \operatorname{span}(v_1,v_2,v_3).
\]
This scalar coefficient is zero. Indeed, projection to the quotient
by that $3$-dimensional span sends $A$ to $a_1^{-1}C_1$,
which is surjective onto a nonzero space since $n\ge4$.
After the scalar cancels, independence of the three vectors makes
the coefficients with images in $\R v_1$ and $\R v_2$ vanish.
What remains is \eqref{IV-descent:eq:rankone} for $i=3$.
For any other nonzero $v_i$, choose two of $v_1,v_2,v_3$ which,
together with $v_i$, are independent, and repeat the same argument.
If $v_i=0$, \eqref{IV-descent:eq:pair} with each of the first two branches
puts $C_i$ directly in their common line. Finally,
$a_iA=C_i-v_i\otimes\beta_i$ has rank at least $n-1$.
\end{proof}

\begin{remark}
The quotient by three independent vectors is the precise use of
$n\ge4$. Replacing that argument by a rank bound on $A$ would
unnecessarily leave a gap when $n=4$ and $\operatorname{rank}A=3$.
The pointwise three-vector lemma does not give the case $n=3$.
\end{remark}

Let $W$ be the open set of image points admitting three independent
inverse-coordinate vectors. For $n\ge4$ the preceding lemma gives a unique
analytic line field
\[
 [A_y]\in\Pj\Hom(T_yM,V),\qquad y\in W,
\]
with rank at least $n-1$. Analyticity follows from the constant
one-dimensional intersection of the analytic spaces
$\mathcal L_1,\mathcal L_2$. Independence and rank persist on a
neighborhood, and uniqueness glues all choices. Equivariance follows
from uniqueness, since
$v_i(\gamma y)=\sigma(\gamma)v_i(y)$ and
$C_i(\gamma y)D_y\gamma=\sigma(\gamma)C_i(y)$.

\subsection{Recurrence makes the projective coframe nonsingular and flat}

\begin{lemma}[Coframe rigidity]\label{IV-descent:lem:coframe}
Suppose an invariant open set $W$ carries the projective matrix field
just constructed and has positive incidence with the family
\eqref{IV-descent:eq:family}, or with the single map $F$. Then, almost everywhere
on the evaluation preimage of $W$,
\begin{equation}\label{IV-descent:eq:identity}
 [A_{\Phi_s(v)}D_v\Phi_s]=[I].
\end{equation}
For each good chart the equality holds throughout its open preimage
of $W$. In particular $A$ is invertible on $W$ after restricting to
the essential chart union. All the good charts have homothety
transitions on their entire overlap domains.
\end{lemma}
\begin{proof}
Put $K(s,v)=[A_{\Phi_s(v)}D_v\Phi_s]\in\Pj\End(V)$.
Its covariance is conjugation by $\sigma(\gamma)$.
In the standard case the field
\[
 \widetilde K([g,s])=[g^{-1}K(s,ge_1)g]
\]
is well-defined on the finite-measure $H$-invariant incidence subset
of $Y=\Gamma\backslash(G\times S)$. That subset is conull by
Lemma~\ref{IV-descent:lem:zeroone}. Under right $h\in H$, the field transforms
by conjugation with $h^{-1}$. Lemma~\ref{IV-descent:lem:unipotent} forces its
line to be fixed by every generating unipotent of $H$. Hence its
matrix commutes with $H$ and is a scalar matrix. It is nonzero,
so \eqref{IV-descent:eq:identity} follows. This proof also applies to $S$
a point and does not assume that $A$ was invertible.

Analyticity upgrades the almost-everywhere identity to the entire
open preimage of $W$ for every retained good chart. Moreover $W$
is conull in every such coordinate image, by the zero--one law
and Fubini. On every full coordinate image, push forward the
projective standard coframe $[dv]$. Two of these extensions agree
on the conull open subset in $W$ and hence everywhere on their
overlap. The derivative of each inverse-coordinate transition is
therefore $c(v)I$. Equality of mixed partial derivatives in dimension
at least two gives $\partial_j c=0$ for every $j$. Thus $c$ is a
nonzero constant on each connected transition domain, proving
\eqref{IV-descent:eq:homotheties}. The standard connections consequently glue.
\end{proof}

\subsection{Nonatomic centers supply the independent vectors}

We give the measurable details needed to apply the preceding lemmas.
First take the essential union of charts, as in Section~\ref{IV-sec:basins}:
for each relatively compact basis ball $B\subset M$, consider the
measurable event that some compact buffered inverse coordinate
branch has image containing $\overline B$. Discard all null such
events, their $\Gamma$ translates, and the null sets discarded in
Lemma~\ref{IV-descent:lem:plaques}. For every remaining chart image point,
some basis neighborhood is covered by chart branches for a
positive-measure set of indices. This remains true after any
further removal of a null set of indices. The resulting union $U$
is invariant and contains almost every center.

For $y\in U$ define
\begin{equation}\label{IV-descent:eq:span}
 L_y=\operatorname{span}\{v\in\Omega_s:\Phi_s(v)=y, s
          \hbox{ a retained index}\}\subset V.
\end{equation}
This is a universally measurable subspace field. For clarity, the
evaluation relation is Borel in its standard measurable coding;
the existence of a finite independent tuple of preimages is an
analytic-set condition. Measurable selection of such tuples gives
the corresponding Grassmannian map on the completed measure space.
One may replace all conull sets of indices by Borel conull subsets,
and then intersect their countably many translates, before making
this construction. The exact covariance is
$L_{\gamma y}=\sigma(\gamma)L_y$.

At a center $y$, a positive-measure set of charts covers a neighborhood
of $y$. Nonatomicity means their centers are different from $y$ for
almost every such index. Hence $L_y\ne0$ at almost every center.
An equivariant measurable proper nonzero subspace field over a
probability-preserving action is impossible here. Indeed induction
gives the map
\[
 [g,s]\longmapsto g^{-1}L_{p(s)}
\]
on $Y$, equivariant under right $G$. Its distribution would be a
$G$-invariant probability on a proper Grassmannian, excluded after
Lemma~\ref{IV-descent:lem:unipotent}. Dimension is constant by ergodicity.
It follows that $L_{p(s)}=V$ almost surely.

Choose three independent inverse vectors at such a center. In
particular $W$ contains almost every center. For $n\ge4$ its own
centered chart has inverse vector zero, so Lemma~\ref{IV-descent:lem:linear}
also gives invertibility there. Invertibility in the good incident
images follows from Lemma~\ref{IV-descent:lem:coframe}. The set $W$ has positive
incidence, and Lemma~\ref{IV-descent:lem:coframe} proves the atlas assertion
of Theorem~\ref{IV-descent:thm:nonatomic}.

For completeness, the essential-union construction makes evaluation
preserve null-set equivalence between $\eta\times\Leb$ and the
smooth measure class on $U$: every basis ball in the union is
covered by a positive-measure family of local diffeomorphisms.
Lemma~\ref{IV-descent:lem:zeroone} therefore proves ergodicity of that class.
Each component of $U$ has positive center measure. Indeed, choose
one basis ball in that component and the positive-measure family of
charts covering it. Each entire chart image is connected and
contains its center, so those centers lie in the same component.
Ergodicity of the center measure forces a single orbit of components;
invariance gives them equal positive mass. There are finitely many.
This completes Theorem~\ref{IV-descent:thm:nonatomic}.

\subsection{Self-overlaps of one continued chart}

For $F$ in Theorem~\ref{IV-descent:thm:single}, use \eqref{IV-descent:eq:span} with all
inverse images of this single map. On $\Gamma\backslash G$ the field
\[
 [g]\longmapsto g^{-1}L_{F(ge_1)}
\]
is defined almost everywhere and is $H$-equivariant. Its values
contain $e_1$. Projective unipotent recurrence in the Plucker
representations forces its values to be $H$-invariant subspaces.
They are therefore $\R e_1$ or $V$. Ergodicity of the source
measure class makes the dimension constant.

If the dimension is $n$, the open locus $W$ with three independent
inverse vectors is nonempty. Lemma~\ref{IV-descent:lem:linear}, followed by
Lemma~\ref{IV-descent:lem:coframe}, gives the second alternative. This applies even if
the initially reconstructed coframe has rank $n-1$.

Suppose instead the dimension is one. Every local transition
$w=T(v)$ then satisfies $w\in\R v$ almost everywhere, and hence
everywhere by analyticity. Thus all identifications are collinear.
Let $J\subset\R^*$ be a compact interval disjoint from $\{1\}$,
or from $\{1,-1\}$ in the antipodal case. The set
\[
 A_J=\{v\in\Omega:\text{some }c\in J\text{ satisfies }
            cv\in\Omega,\ F(cv)=F(v)\}
\]
is universally measurable and invariant. It misses a sufficiently
small punctured ball, since both $v$ and $cv$ are in the originally
embedded collar there. The zero--one law implies that $A_J$ is
null. A nonprescribed identification would give a local deck
transition $T(v)=c(v)v$, with $c(v)$ in one such $J$ on an open
set. This contradicts nullity. A countable collection of intervals
$J$ covers all nonprescribed scalars. Hence there are no such
identifications, proving the first alternative.

Suppose $F$ has an analytic extension at zero. Its equivariance extends
as a germ by continuity. In target coordinates let $J$ be the Jacobian
determinant of the extension. It is nonzero on the punctured ball. If
$J(0)=0$, let $P$ be its first nonzero homogeneous term. The chain rule,
using $\det\sigma(\gamma)=1$, makes the line $[P]$ invariant under
$\Gamma$. The map $[g]\mapsto g^{-1}[P]$ pushes the invariant probability
on $\Gamma\backslash G$ to an invariant probability on the projective
space of homogeneous polynomials of that degree. Applying
Lemma~\ref{IV-descent:lem:unipotent} to the elementary unipotent generators
forces $P$ to be fixed by $G$. Transitivity on $V\setminus\{0\}$ then
makes $P$ constant there, contrary to its positive degree. Thus $J(0)\ne0$,
and the analytic inverse theorem gives a chart at zero. Its inverse
differential agrees with the descended coframe on the punctured
neighborhood and extends the coframe and connection across zero.
The two alternatives
are exclusive because in the embedded case the inverse vectors
span a line. This proves Theorem~\ref{IV-descent:thm:single}.

\subsection{A useful obstruction to radial ends}

\begin{proposition}[An affine region with spanning returns is not radial]
\label{IV-descent:prop:radial}
Let the affine region in Theorem~\ref{IV-descent:thm:single}(2) meet an invariant
open set carrying a nonvanishing invariant analytic vector field.
Then such an intersection is impossible. No connection assumption
on the field's region is required.
\end{proposition}
\begin{proof}
On the open preimage of the field's region, vector-stabilizer recurrence
as in Section~\ref{IV-sec:basins} makes the pullback radial almost
everywhere. Continuity makes it radial throughout that open set.
If the intersection is nonempty, choose an image point $y$ there at
which the inverse vectors span $V$. The inverse differentials at all
preimages of $y$ represent the same projective coframe, so they send the
nonzero field vector at $y$ to proportional nonzero vectors. By
radiality each resulting vector lies on the line of the corresponding
inverse vector. All the inverse vectors are therefore collinear,
contradicting the spanning condition.
\end{proof}

For Theorem~\ref{IV-descent:thm:nonatomic} the corresponding exclusion is the
radial-field detection theorem in Section~\ref{IV-sec:basins}: such a field
would make the center probability atomic. These two exclusions
are the inputs for ruling out exterior and Hopf ends of the affine
regions constructed here.

\section{Affine descent for \texorpdfstring{$n=3$}{n = 3}}\label{IV-series:descent-three}

\subsection{The descent statements for $n=3$}

Put $V=\R^3$, $G=\SL(3,\R)$, and let $\Gamma<G$ be a lattice.
Use the exact analytic chart domains, probability family, and continued
single-map hypotheses in Section~\ref{IV-sec:descent}, now with $n=3$. In the family
case the center probability is nonatomic. Standard and dual models are
both allowed; inverse transpose reduces the proof to the standard case.

\begin{theorem}[Affine descent for $n=3$]\label{IV-affine_descent:thm:descent}
The nonatomic affine-descent theorem and the continuation
dichotomy in Section~\ref{IV-sec:descent} hold with $n=3$. In the affine alternative
the inverse differentials give an equivariant analytic projective
coframe, and every connected coordinate transition is
$v\mapsto av+b$, with $a\ne0$. The conclusions about the flat
connection, smooth measure class, components, and radial obstruction
are unchanged.
\end{theorem}

Here are the parts of the earlier proof that do not require $n\ge4$.
The vector-stabilizer zero--one law, projective unipotent recurrence,
and the recurrence excluding transverse intersections of hyperplane
plaques apply for $n=3$. In the last argument a diagonal element expands
a two-plane and contracts its complementary line; the same finite
probability space and implicit-function argument apply. The measurable
support selection then gives the plaque constraint simultaneously for
all retained charts. The span argument gives full spanning at almost
every nonatomic center, and the single-map span is almost everywhere
either a line or all of $V$. These use the irreducibility of the
$\SL(2,\R)$ block, not higher rank of an arithmetic Levi lattice.

Write an inverse branch and its derivative as
\[
 v_i(y)\in V,\qquad C_i(y)=D_yv_i:T_yM\longrightarrow V.
\]
Every $C_i$ is invertible. The simultaneous plaque constraint is
\begin{equation}\label{IV-affine_descent:eq:plaque}
 \ell(v_i)=\ell(v_j)=0
 \quad\Longrightarrow\quad
 [\ell C_i]=[\ell C_j],\qquad \ell\ne0.
\end{equation}
For independent $v_i,v_j$ put
\begin{equation}\label{IV-affine_descent:eq:spaces}
 \mathcal L_i=\R C_i+
       \{v_i\otimes\beta:\beta\in T_y^*M\}.
\end{equation}
Their common annihilator is one dimensional, so \eqref{IV-affine_descent:eq:plaque}
still implies
\[
 C_j=aC_i+v_i\otimes\beta+v_j\otimes\gamma,\qquad a\ne0.
\]
Consequently $\mathcal L_i\cap\mathcal L_j$ is a line, represented
by a matrix of rank at least two. Uniqueness follows because a nonzero
multiple of an invertible three-by-three matrix cannot have image in
$\operatorname{span}(v_i,v_j)$. This pairwise construction remains
analytic while the two vectors stay independent.

Three branches need not give the same line. For instance, at one point
the data
\[
 v_i=e_i,\qquad C_1=I,\quad
 C_2=\operatorname{diag}(2,3,1),\quad
 C_3=\operatorname{diag}(4,1,5)
\]
satisfy \eqref{IV-affine_descent:eq:plaque}, but the line
$\mathcal L_1\cap\mathcal L_2=
\R\operatorname{diag}(2,1,1)$ does not belong to $\mathcal L_3$.
We will use the whole recurrent family, not assert a stronger
three-branch linear-algebra lemma.

\subsection{The one-dimensional coordinate foliations}

For a nonzero inverse branch define its radial vector field by
\[
 E_i=C_i^{-1}v_i.
\]
Its integral curves are the images of radial curves in the corresponding
chart. For independent $v_i,v_j$ define a projective covector
\[
 \kappa_{ij}(y)=[\nu(v_i,v_j)],\qquad
 \nu(v,w)(z)=\det(v,w,z).
\]

\begin{lemma}[Common planes and coordinate curves]\label{IV-affine_descent:lem:curves}
The map $\kappa_{ij}$ has rank one. Its level sets are analytic
surfaces tangent to the common plane
\[
 P_{ij}=\ker(\nu(v_i,v_j)C_i)
       =\ker(\nu(v_i,v_j)C_j).
\]
If $v_1,v_2,v_3$ are independent, then either $E_1,E_2,E_3$ are
independent or all three lie in the same line. In the independent
case there are analytic coordinates $(t_1,t_2,t_3)$ for which
$\R E_i=\R\partial_{t_i}$.
\end{lemma}
\begin{proof}
Choose a local nonzero representative $\ell$ of $\kappa_{ij}$.
Differentiating $\ell(v_i)=\ell(v_j)=0$ gives
\[
 (n\ell)v_i=-\ell C_i,\qquad
 (n\ell)v_j=-\ell C_j.
\]
The map from $V^*/\R\ell$ to $\R^2$ given by evaluation on
$(v_i,v_j)$ is an isomorphism. The two rows on the right are nonzero
and proportional by \eqref{IV-affine_descent:eq:plaque}. Hence $D\kappa_{ij}$ has
rank exactly one and kernel $P_{ij}$. The analytic constant-rank
theorem supplies a local first integral for this plane field.

Each $P_{ij}$ contains $E_i,E_j$. At a triple of independent source
vectors the two planes $P_{ij},P_{ik}$ are distinct, since $C_i$
is invertible and their source annihilators are independent. If one
pair $E_i,E_j$ is independent, a third vector in their plane would
make two of these distinct planes equal. Thus the three are independent.
Otherwise every pair is collinear. In the independent case choose
$t_i$ to be a first integral for $P_{jk}$. Their differentials are
independent, and their coordinate lines have the asserted tangents.
\end{proof}

The lemma supplies actual local coordinate curves, not merely three
abstract directions at a point. It leaves three separate choices of
leaf parameter. Nor does it by itself exclude the common-line case.
The next two lemmas provide precisely the extra compatibility supplied
by recurrence.

\subsection{A finite polynomial envelope of all inverse jets}

Let
\[
 X=V\oplus\End(V),\qquad
 h\cdot(w,C)=(hw,hCh^{-1}),\qquad
 H=\{h\in G:he_1=e_1\}.
\]
Let $\mathcal P_3$ be the finite-dimensional space of real polynomial
functions on $X$ of total degree at most three. For $\mathcal E\subset X$
set
\begin{equation}\label{IV-affine_descent:eq:envelope}
 I_3(\mathcal E)=\{p\in\mathcal P_3:p|_{\mathcal E}=0\},
 \qquad
 Z_3(\mathcal E)=\{x\in X:p(x)=0\text{ for all }p\in I_3(\mathcal E)\}.
\end{equation}
No boundedness, closedness, or measure on $\mathcal E$ is required.

\begin{lemma}[Recurrence of polynomial relations]\label{IV-affine_descent:lem:envelope}
Suppose a measurable family $\mathcal E_z\subset X$ over a finite
$H$-invariant measure space satisfies
$\mathcal E_{zh}=h^{-1}\mathcal E_z$. Then
$Z_3(\mathcal E_z)$ is $H$-invariant for almost every $z$.
It suffices here that the subspace field $I_3(\mathcal E_z)$ is
measurable on the completed measure space.
\end{lemma}
\begin{proof}
The action on $X$ is linear, so it preserves $\mathcal P_3$; each
one-parameter unipotent in $H$ acts unipotently on this polynomial
space. On each stratum of fixed dimension, $I_3(\mathcal E_z)$ is
an equivariant map into a Grassmannian. Its Plucker line is fixed
almost everywhere by projective unipotent recurrence. Apply this
to finitely many one-parameter unipotent subgroups generating $H$.
The subspace $I_3(\mathcal E_z)$, and therefore its common zero set,
is $H$-invariant. Zero and full-dimensional subspaces cause no exception.
\end{proof}

The plaque relation for independent vectors is the polynomial identity
\begin{equation}\label{IV-affine_descent:eq:polynomial}
 \bigl(\nu(v,w)C\bigr)\wedge\bigl(\nu(v,w)D\bigr)=0
 \quad\text{for }(v,C),(w,D)\in\mathcal E.
\end{equation}
Each scalar coefficient has degree at most three in $(v,C)$ and at
most three in $(w,D)$ separately. Thus if \eqref{IV-affine_descent:eq:polynomial} holds
on $\mathcal E\times\mathcal E$, it holds on
$Z_3(\mathcal E)\times Z_3(\mathcal E)$. Indeed, first fix the
second argument in $\mathcal E$ and extend in the first variable
using \eqref{IV-affine_descent:eq:envelope}; then fix the first argument in $Z_3(\mathcal E)$
and extend in the second variable. This separate degree bound is why
a single finite-dimensional polynomial space suffices.

\subsection{The stabilizer calculation}

\begin{lemma}[Orbit compatibility]\label{IV-affine_descent:lem:orbit}
Let $Z\subset X$ be $H$-invariant, contain $(e_1,I)$, and satisfy
\eqref{IV-affine_descent:eq:polynomial} for every pair of its points. If $(w,C)\in Z$,
$C$ is invertible, and $w\notin\R e_1$, then
\begin{equation}\label{IV-affine_descent:eq:orbitconclusion}
 C=aI+w\otimes\beta,\qquad a\ne0.
\end{equation}
\end{lemma}
\begin{proof}
The group $H$ is transitive on $V\setminus\R e_1$, so conjugate
the pair to arrange $w=e_2$. Pairing with $(e_1,I)$ makes the last
row of $C$ a nonzero multiple of $e_3^*$:
\[
 C=\begin{pmatrix}p&q&r\\ b&c&f\\0&0&a\end{pmatrix},
 \qquad a\ne0.
\]

We first justify a same-vector comparison; simply substituting equal
vectors into \eqref{IV-affine_descent:eq:polynomial} would give no information. If
$h_0\in H$ fixes $w$, and $\ell(w)=0$, choose $q_0\in\ker\ell$
independent of $w$. The orbit map $H\to V\setminus\R e_1$ is a
submersion. There are $h_\epsilon\to h_0$ with
$h_\epsilon w=w+\epsilon q_0$. Apply \eqref{IV-affine_descent:eq:polynomial} to
$(w,C)$ and $(h_\epsilon w,h_\epsilon C h_\epsilon^{-1})$, divide
by $\epsilon^2$, and pass to the limit. Since
$\nu(w,q_0)$ is a nonzero multiple of $\ell$, this proves
\[
 [\ell C]=[\ell h_0Ch_0^{-1}]
 \quad\text{for every }0\ne\ell\in w^\perp.
\]
Linearity on the two-dimensional space $w^\perp$ makes the
proportionality scalar independent of $\ell$. Therefore
\begin{equation}\label{IV-affine_descent:eq:samevector}
 h_0Ch_0^{-1}=b_0C+w\otimes\alpha,
 \qquad b_0\ne0.
\end{equation}

Use the stabilizer elements
\[
 h_0=\begin{pmatrix}1&0&x\\0&1&y\\0&0&1\end{pmatrix}.
\]
The last row of $h_0Ch_0^{-1}$ remains $(0,0,a)$, so
\eqref{IV-affine_descent:eq:samevector} has $b_0=1$. Its first-row difference is
\[
 (0,0,(a-p)x-qy).
\]
It must vanish for all $x,y$. Thus $p=a$, $q=0$, and
\[
 C=aI+e_2\otimes(b,c-a,f)+r\,e_1\otimes e_3^*.
\]
Finally use
\[
 R=\begin{pmatrix}1&0&0\\0&0&-1\\0&1&0\end{pmatrix}\in H.
\]
The two vectors $e_2,Re_2=e_3$ have annihilator $\R e_1^*$.
The first rows of $C$ and $RCR^{-1}$ are respectively
$(a,0,r)$ and $(a,-r,0)$. They are proportional by
\eqref{IV-affine_descent:eq:polynomial}; since $a\ne0$, this forces $r=0$.
Conjugating back proves \eqref{IV-affine_descent:eq:orbitconclusion}.
\end{proof}

\subsection{Applying recurrence on the finite suspension}

Work with the invariant retained family for which \eqref{IV-affine_descent:eq:plaque}
holds simultaneously. For $z=[g,s]\in Y=\Gamma\backslash(G\times S)$,
put $y=\Phi_s(ge_1)$. This is defined on a conull $H$-invariant
subset by the zero--one law. For every inverse branch $v_i$ at $y$,
form its relative first jet
\begin{equation}\label{IV-affine_descent:eq:relative}
 w_i=g^{-1}v_i(y),\qquad
 J_i=g^{-1}C_i(y)D_{ge_1}\Phi_s\,g,
 \qquad \mathcal E_z=\{(w_i,J_i):i\text{ any inverse branch at }y\}.
\end{equation}
The distinguished branch gives $(e_1,I)\in\mathcal E_z$.
Changing the lattice representative leaves this set unchanged;
right multiplication by $h\in H$ transforms it by
$(w,J)\mapsto(h^{-1}w,h^{-1}Jh)$. In the single-map case use
$S$ a point and all inverse branches of $F$.

For completeness, the required polynomial subspace is measurable.
In the standard measurable coding of analytic charts, the relation
specifying an inverse value and its first derivative is Borel.
The span of their degree-at-most-three monomial evaluation vectors
has measurable dimension and measurable basis after completing the
measure: existence of an independent finite tuple is an analytic-set
condition, and measurable selection supplies a tuple on each rank
stratum. Its annihilator is exactly $I_3(\mathcal E_z)$. This is
the same completed-measure selection used for the inverse-vector span;
it does not require a probability distribution on the inverse branches.

The plaque constraint gives \eqref{IV-affine_descent:eq:polynomial} on $\mathcal E_z$.
Lemmas~\ref{IV-affine_descent:lem:envelope} and~\ref{IV-affine_descent:lem:orbit}, and the separate-degree
extension of \eqref{IV-affine_descent:eq:polynomial}, imply that for almost every $z$,
\begin{equation}\label{IV-affine_descent:eq:alljets}
 J_i=a_iI+w_i\otimes\beta_i
 \qquad\text{whenever }w_i\notin\R e_1.
\end{equation}
For $0\ne w_i\in\R e_1$, the original simultaneous plaque
constraint with the distinguished branch implies the same formula:
on the two-dimensional annihilator of $e_1$, all rows are proportional
to themselves, so the scalar is common and nonzero. This last step
uses actual branches, not the weaker polynomial identity at a
collinear pair.

If there is one inverse vector outside $\R e_1$, the unique line
in \eqref{IV-affine_descent:eq:spaces} for it and the distinguished branch is now
$\R I$. Every nonzero branch has $I$ in its space \eqref{IV-affine_descent:eq:spaces}.
If an inverse vector is zero, its differential lies in the spaces
for each of these two independent nonzero vectors, by
\eqref{IV-affine_descent:eq:plaque} on their two-dimensional annihilators. Hence that
differential is itself a scalar matrix. Thus, at almost every evaluated
point with noncollinear inverse vectors, the spaces in
\eqref{IV-affine_descent:eq:spaces}, written in these relative coordinates, satisfy
\begin{equation}\label{IV-affine_descent:eq:common}
 \mathcal L_i\cap\mathcal L_j=\R I
 \qquad\text{whenever }w_i,w_j\text{ are independent}.
\end{equation}
In particular the conclusion uses all branches simultaneously,
including any branches selected after the recurrence argument.

\subsection{Analytic gluing and the descent conclusion}

Let $W$ be the open locus admitting two independent inverse vectors.
Each local choice of two inverse branches defines an analytic line
$\mathcal L_i\cap\mathcal L_j$ on its independence domain, of rank
at least two. Assertion~\eqref{IV-affine_descent:eq:common}, disintegration under
$g\mapsto ge_1$, and Fubini show that all these local lines agree
almost everywhere and are nonsingular almost everywhere. To pass
from evaluated points to smooth measure on the image, use the
essential-union construction in Section~\ref{IV-sec:descent}: every basis ball
in the retained union is covered by a positive-measure family of
local diffeomorphisms. The exceptional property here is intrinsic
to the collection of branches at the image point. Thus it is a
smooth null set. For a single map the same assertion follows
directly by change of variables in a countable inverse-chart cover.

Analytic uniqueness makes the local lines agree throughout their
overlaps. They give a well-defined equivariant analytic projective
matrix field $[A_y]$ on $W$. At almost every evaluated point of $W$,
\eqref{IV-affine_descent:eq:common} says
\[
 [A_{\Phi_s(v)}D_v\Phi_s]=[I].
\]
Whenever $W$ has positive incidence, the zero--one law makes its
evaluation preimage conull. Discarding a further null set of indices,
analyticity gives the displayed identity on the whole open preimage
of $W$ in every remaining chart. In particular $A$ is invertible
there. On each full chart image extend it by the projective standard
coframe. The extensions agree on their conull intersections with
$W$, hence everywhere on their overlaps. A transition therefore has
derivative $c(v)I$. Equality of mixed partials makes $c$ constant
on every connected transition domain. This gives the homothety atlas.

For nonatomic centers the inverse-span argument in Section~\ref{IV-sec:descent}
puts almost every center in $W$, so positive incidence holds. For a
single continued map, the spanning alternative gives positive
incidence with $W$. Its line-span alternative is embedded, modulo
the prescribed antipodal identification, by the small-ball injectivity
and zero--one argument of that section. The remaining proofs of
ergodicity, finiteness of components, and the radial obstruction use
only the resulting atlas and apply without change. This proves
Theorem~\ref{IV-affine_descent:thm:descent}.

\section{Exterior continuation and detected ends}\label{IV-series:exterior}\label{IV-sec:exterior}

\subsection{The radial collar input}

Let $n\ge3$ and let $\Gamma<\SL(n,\R)$ be any lattice.
Work in a compact connected analytic $n$-manifold $P$, possibly with boundary,
on which $\Gamma$ acts analytically. In the intended application
$P$ is a component obtained by cutting the original manifold along
an invariant analytic projective hypersurface; restrict to a
finite-index subgroup preserving the component and the chosen end.
Write $N_\infty\subset\partial P$ for that boundary component.
Cutting is useful even if the growing region returns to the same
hypersurface in the original manifold: its two sides are then
distinct boundary copies when the hypersurface is two-sided.

We first treat a spherical radial collar. The precise input is an
analytic embedding
\begin{equation}\label{IV-exterior:eq:input}
 \phi:E_R\longrightarrow\operatorname{int}P,
 \qquad E_R=\{z\in\R^n:|z|>R\},
\end{equation}
onto the interior of a one-sided collar of $N_\infty$.
Its radial compactification at infinity is the given collar, with
boundary labeling $z/|z|\in S^{n-1}$. In particular the images
of $\{|z|>T\}$ form a neighborhood basis of $N_\infty$ on
the chosen side as $T\to\infty$. We assume the identities
\begin{equation}\label{IV-exterior:eq:germ}
 \alpha(\gamma)\phi(z)=\phi(\gamma z)
 \quad\text{for all sufficiently large }|z|,
 \qquad\gamma\in\Gamma.
\end{equation}
The threshold can depend on $\gamma$. A local $G$-collar also
satisfies the analogous identities for its local Lie-group action.

Here is how this input appears in the Fisher--Melnick flow model,
once the analytic projective collar of Theorem~\ref{II-main:thm:recognition}
or Theorem~\ref{III-main:thm:recognition-three}, according to the
dimension, is available.
In a local projective chart its formula is
\[
 g(u,t)=\left(\frac{s+Au}{\alpha+ru},
       \Phi^{B\partial_t}_{\log(\alpha+ru)}(t)\right),
 \qquad g=\begin{pmatrix}\alpha&r\\s&A\end{pmatrix}.
\]
On a side where $B\ne0$, put
\begin{equation}\label{IV-exterior:eq:develop}
 \tau(t)=\int^t\frac{ds}{B(s)},\qquad
 F(u,t)=e^{\tau(t)}(1,u).
\end{equation}
Then $F(g(u,t))=gF(u,t)$. In spherical collar coordinates
$(\theta,t)$ one uses the norm cocycle
$\log(\|gv\|/\|v\|)$ and $F(\theta,t)=e^{\tau(t)}\theta$.
These are the usual equivalent local and compact-direction forms
of the induced collar \cite[Sections 3 and 5]{FisherMelnick}.
The local-description theorems provide these full-angular inputs
for every lattice. The continuation and end-comparison arguments
below use connectedness of punctured balls and exterior domains,
finite-volume recurrence, and the vector-model comparison lemma
of Part II. These inputs hold for every $n\ge3$; the
$3$-dimensional case requires no additional Levi normalization
at this stage.

If $t>0$ measures distance into the chosen side and
$B(t)=c t^k+O(t^{k+1})$, $c\ne0$, then
\[
 \tau(t)\sim\begin{cases}
 c^{-1}\log t,& k=1,\\
 -\bigl(c(k-1)t^{k-1}\bigr)^{-1},&k>1.
 \end{cases}
\]
Thus $e^{\tau(t)}\to\infty$ exactly when $c<0$.
The inverse of $F$ then gives \eqref{IV-exterior:eq:input}, after shrinking
the collar. When $c>0$ it gives a punctured-ball collar instead.

\begin{corollary}[Which neutral sides are exterior]\label{IV-exterior:cor:parity}
At a zero of even order $k\ge2$, precisely one of the two local
sides has developing radius tending to infinity. At an odd-order
zero, both sides are exterior if the leading coefficient is
negative, and both are interior if it is positive.
\end{corollary}
\begin{proof}
On the other side put $s=-t>0$. Its coefficient is
$-B(-s)=(-1)^{k+1}c s^k+O(s^{k+1})$.
Apply the preceding sign criterion. The statement concerns the
two local sides, whether or not they are globally distinct.
\end{proof}

If $B\equiv0$, the collar has no transverse motion and the
full-dimensional developing construction does not apply.

\subsection{Exterior continuation has no identifications}

\begin{lemma}[Simultaneous expansion]\label{IV-exterior:lem:expand}
Given finitely many nonzero vectors $z_1,\ldots,z_l\in\R^n$,
there is a sequence $\gamma_d\in\Gamma$ with
$|\gamma_nz_j|\to\infty$ for every $j$.
\end{lemma}
\begin{proof}
Put $a_t=\operatorname{diag}(e^{(n-1)t},e^{-t},\ldots,e^{-t})$.
Choose $g\in G$ such that the first coordinate of every $gz_j$
is nonzero and $g\Gamma$ is recurrent under $a_1$ to a compact
set in $G/\Gamma$. The first condition excludes finitely many
Haar-null hypersurfaces and the second holds almost everywhere
by recurrence and a countable compact exhaustion. Compact sets
in the quotient have bounded representatives, so at return times
$d_k\to\infty$ write $a_{d_k}g=k_k\gamma_k$ with $k_k$ in a
fixed compact subset of $G$. Then
\[
 |\gamma_kz_j|\ge
 \frac{|a_{d_k}gz_j|}{\sup_k\|k_k\|}\longrightarrow\infty.
\]
No element of $\Gamma$ is required to be unipotent.
\end{proof}

\begin{theorem}[Embedded continuation]\label{IV-exterior:thm:continuation}
Under \eqref{IV-exterior:eq:input}--\eqref{IV-exterior:eq:germ}, there is a unique
$\Gamma$-equivariant analytic local diffeomorphism
\[
 \Phi:\R^n\setminus\{0\}\longrightarrow P
\]
extending $\phi$. It is injective. Its image is exactly
$U=\bigcup_{\gamma\in\Gamma}\alpha(\gamma)\phi(E_R)$.
In particular $U$ carries a global analytic $G=\SL(n,\R)$-action
transported from the standard action on $\R^n\setminus\{0\}$,
whose restriction to $\Gamma$ is the original action. If the input
is a local $G$-collar, this action agrees with it on the chosen side.
\end{theorem}
\begin{proof}
For any $\gamma$, the set
\[
 D_\gamma=E_R\cap\gamma^{-1}E_R
 =\{r\theta:r>\max(R,R/|\gamma\theta|)\}
\]
is connected, since $S^{n-1}$ is connected. Both sides of
\eqref{IV-exterior:eq:germ} are analytic on $D_\gamma$ and agree near
infinity. They therefore agree on all of $D_\gamma$.

Every $z\ne0$ can be moved into $E_R$ by Lemma~\ref{IV-exterior:lem:expand}.
Choose such a $\gamma$ and set
\begin{equation}\label{IV-exterior:eq:continuation}
 \Phi(z)=\alpha(\gamma)^{-1}\phi(\gamma z).
\end{equation}
If $\delta z\in E_R$ as well, apply the preceding overlap identity
to $\delta\gamma^{-1}$ at $\gamma z$. It proves independence of
the choice. Each choice works on a neighborhood of $z$, so $\Phi$
is an analytic local diffeomorphism. The formula proves equivariance,
uniqueness, and the assertion about its image.

If $\Phi(x)=\Phi(y)$, use Lemma~\ref{IV-exterior:lem:expand} to move both
$x,y$ into $E_R$ by the same lattice element. Equivariance then
gives $\phi(\gamma x)=\phi(\gamma y)$. Injectivity of $\phi$
implies $x=y$. Thus $\Phi$ is an analytic diffeomorphism onto $U$.
Transporting the standard $G$-action gives the last assertions;
the agreement with the collar follows from \eqref{IV-exterior:eq:develop}.
No completeness assertion at the other end is used.
\end{proof}

This proof also works with source
$(\R^n\setminus\{0\})/\{\pm I\}$ and its projective exterior
collar. Both the connected-overlap and simultaneous-expansion
arguments descend to this quotient. Subsequent comparisons will
be stated using spherical developing coordinates, allowing a
finite radial cover for an incoming projective collar.

\subsection{The remaining end is one invariant continuum}

\begin{proposition}[The end limit set]\label{IV-exterior:prop:end}
In either version of Theorem~\ref{IV-exterior:thm:continuation}, put
\[
 K_0=\bigcap_{\epsilon>0}
     \overline{\Phi(\{0<|z|<\epsilon\})}^{\,P}.
\]
Then $K_0$ is nonempty, compact, connected, and $\Gamma$-invariant.
It is disjoint from $U\cup N_\infty$, and
\[
 \overline U^{\,P}\setminus U=N_\infty\cup K_0.
\]
If $K_0$ is finite, it is one $\Gamma$-fixed point. In the
spherical version the continuation then fills a topological disk
bounded by $N_\infty$. For the quotient source, radial sets in
this formula mean their images modulo $\{\pm I\}$.
\end{proposition}
\begin{proof}
The sets in the intersection are nested nonempty compact connected
sets. Equivariance and $|\gamma z|\le\|\gamma\||z|$ give invariance.
An embedding has a continuous inverse on its open image, so a
sequence tending to zero cannot have an image limit in $U$.
The exterior chart gives a full one-sided neighborhood basis of
$N_\infty$ in the cut manifold $P$. Points of this neighborhood
already have large developing radius by injectivity, so $K_0$
cannot meet $N_\infty$. This is why the cut formulation is useful.

For any sequence in $U$ converging to its frontier, a sequence of
preimages has a subsequence tending either to zero or to infinity:
otherwise it has a subsequence in a compact annulus and the limit
lies in $U$. The first case gives $K_0$, the second $N_\infty$.
Both sets are actually attained as limits, proving the equality.

A finite connected $K_0$ is a singleton $\{p\}$, fixed by
$\Gamma$. In the spherical version, compactness makes $\Phi$ extend continuously and
injectively over $0\mapsto p$. Invariance of domain, applied in
the double of $P$, shows $p\in\operatorname{int}P$ and makes this
extension a local homeomorphism at zero. Adjoining the spherical
end at infinity gives a closed topological disk. Its image is
open as well as closed in the connected component of $P$ that
contains it: openness at $N_\infty$ uses the original collar.
Thus it fills that component. Analytic regularity at $p$ does not
follow from invariance of domain; it will follow from
Theorem~\ref{IV-exterior:thm:cap} below.
\end{proof}

\subsection{A detected second collar automatically closes the end}

The next statement makes the needed geometric detection input exact.
An incoming radial collar of a compact connected analytic projective
hypersurface $N_0$ is represented by an analytic local diffeomorphism
$\psi:C\to\operatorname{int}P$, where either
\[
 C=\{0<|z|<r\}\quad\text{or}\quad C=E_r.
\]
It may be a finite radial cover of the interior collar. Its radial
end tends uniformly to $N_0$, and its image contains that chosen
side of a neighborhood of $N_0$. Assume it satisfies the standard
or dual analogues of \eqref{IV-exterior:eq:germ}. These assumptions are the
radial form of an already constructed nontrivial analytic flow
collar; the product case $B\equiv0$ is not included.

\begin{theorem}[Closing a detected end]\label{IV-exterior:thm:close}
Use either continuation above, and let $\pi$ denote the identity
or the antipodal quotient map of its source. Suppose
$N_0\subset K_0$ has such an incoming collar and
$\psi(C)\cap U\ne\varnothing$. Then its representation type is
standard, its developing radius tends to zero at $N_0$, and for
some $a\ne0$,
\[
 \psi(z)=\Phi(\pi(a z))\qquad(z\in C).
\]
Consequently $K_0=N_0$. The global $G$-action on $U$ agrees with
the incoming local $G$-collar on this side. After cutting at $N_0$
and making the indicated spherical boundary covers, the closure
is a Fisher--Melnick tube \cite[Section 3.2.1]{FisherMelnick}
with its two analytic collars.
\end{theorem}
\begin{proof}
For a punctured ball, $C\cap\gamma^{-1}C$ is a convex
neighborhood of zero with zero removed, hence connected.
For an exterior domain use the connectedness argument above.
Analytic continuation extends every lattice germ identity over
this whole intersection. It follows, as in
\eqref{IV-exterior:eq:continuation}, that $\psi$ continues to an equivariant
analytic local diffeomorphism $\Psi$ on
$\Omega=\bigcup_{\gamma\in\Gamma}\gamma C$ (using the dual
action for dual input). The domain $\Omega$ is connected:
all its defining domains have a common punctured neighborhood
of zero, or a common exterior neighborhood, pairwise.

Let $O=\Psi^{-1}(U)$, a nonempty invariant open subset of
$\Omega$, and put $h=\Phi^{-1}\Psi$ there. In the standard case
Lemma~\ref{II-exterior:lem:comparison} gives $h(z)=a z$ on one nonempty
open component, modulo sign in the quotient case. Hence
$\Psi(z)=\Phi(\pi(a z))$ on all of connected
$\Omega$ by the analytic identity theorem. In the dual case,
replace $\Gamma$ by $\Gamma^*=\{\gamma^{-T}:\gamma\in\Gamma\}$
in that lemma. It would give a standard-to-dual comparison into
nonzero vectors, which is impossible.

If $C$ were exterior, the equality would identify its end with
$N_\infty$, contradicting $N_0\subset K_0$ and
Proposition~\ref{IV-exterior:prop:end}. Thus $C$ is a punctured ball. Its
image under $\Phi(\pi(a\cdot))$ is a full small-radius tail of $U$.
Uniform collar convergence to $N_0$ gives exactly $K_0=N_0$.
The scalar transition commutes with $G$, proving agreement of
the local actions. Cutting the end hypersurfaces supplies the
usual compact spherical interval collar description. In particular
the displayed equality permits only the original source's
antipodal identification, if present, inside $C$.
\end{proof}

\begin{proposition}[A product collar cannot meet the open orbit]\label{IV-exterior:prop:product}
Suppose a compact projective hypersurface has a full analytic
product collar, after the relevant finite covers and finite-index
restriction, with lattice action germ
\[
 \gamma(\theta,t)=(\gamma\theta,t).
\]
Then a sufficiently small such collar is disjoint from $U$.
In particular, for an already constructed analytic collar
meeting $U$, the product case $B\equiv0$ is automatically excluded.
\end{proposition}
\begin{proof}
Choose a finite symmetric generating set. Compactness of the
hypersurface gives one collar width on which all its generator
formulas hold. Each generator preserves this entire product
collar, because it leaves $t$ unchanged; the collar is therefore
$\Gamma$-invariant. If its intersection with $U$ were nonempty,
pulling back $t$ by $\Phi$ would give an analytic invariant
submersion on a nonempty invariant open subset of
$\R^n\setminus\{0\}$, also in the antipodal case after lifting.
Lemma~\ref{II-general:lem:linear-models} makes every such scalar
constant almost everywhere, and analyticity makes it constant
on its entire open domain.
This contradicts the submersion property.
\end{proof}

\subsection{A finite orbit forces an analytic fixed-point filling}

At a finite orbit the invariant probability falls under the
same derivative trichotomy as any invariant probability.
Theorem~\ref{II-loc:thm:compact-finite} rules out its compact branch
for an infinite-image action. Corollary~\ref{I-power:thm:main} supplies the standard or dual analytic fixed-point filling
charts. Thus no special arithmetic fixed-point linearization
input is needed here.

\begin{theorem}[A detected finite orbit closes the end]\label{IV-exterior:thm:cap}
If $K_0\cap\operatorname{int}P$ contains a finite $\Gamma$-orbit,
then in the spherical version $K_0=\{p\}$ and $\Phi$ extends
analytically and locally diffeomorphically over $0\mapsto p$.
The open-orbit $G$-action extends analytically over this fixed point.
In the antipodal version no such finite orbit can occur.
In the application where $P$ is cut only along projective
hypersurfaces, all finite orbits are interior, so the conclusion
applies to every finite orbit in $K_0$.
\end{theorem}
\begin{proof}
Pass to the finite-index stabilizer $\Delta$ of a point $p$ in
that orbit. If its derivative type were compact, the
compact-measure theorem would force finite global image, contrary
to the standard action on $U$. The standard and dual specialization
of analytic germ cocycle superrigidity
(Corollary~\ref{I-power:thm:main}) therefore supplies an analytic chart
$\psi:B_r\to\operatorname{int}P$, $\psi(0)=p$, with standard
or dual equivariance as germs. Remove its sign character by a
further finite-index restriction.

Since $p\in\partial U$, the image of $B_r\setminus\{0\}$
meets $U$. Continue this punctured chart using the connected
ball-overlap argument in Theorem~\ref{IV-exterior:thm:close}, now for
$\Delta$. Lemma~\ref{II-exterior:lem:comparison} rules out the dual case
and gives, on the whole punctured ball,
\[
 \psi(z)=\Phi(\pi(a z)),\qquad a\ne0.
\]
For a spherical source this identity contains a full small-radius
tail of $U$. Continuity of $\psi$ at zero gives $K_0=\{p\}$,
and the formula $\Phi(z)=\psi(z/a)$ proves analytic extension
with invertible derivative. Scalar matching also extends the
transported $G$-action. For an antipodal source the identity would
give $\psi(z)=\psi(-z)$, contradicting injectivity of the chart.

Finally, a projective boundary has no finite orbit under any
finite-index subgroup of $\Gamma$: a finite orbit would give
an invariant line for a Zariski-dense finite-index subgroup.
This proves the assertion about the cut application.
\end{proof}

Thus in the spherical case one need only find a finite orbit
somewhere in the continuum; one need not first prove that the
whole continuum is a point. Proposition~\ref{IV-exterior:prop:end} then
identifies the resulting component as a disk, now with an
analytic fixed-point filling and the prescribed analytic outer collar.

\section{Completion of exterior tubes}\label{IV-series:tubes}\label{IV-sec:tubes}

\subsection{Closing the second end of an exterior tube}

We use local description in the sense of Definition~\ref{II-def:recognition}.
In the hypersurface alternative the measure supplies a point on
one of the compact central hypersurfaces. Their finite orbit and
the minimality of the standard projective action imply that a
closed invariant set supporting that measure contains the
corresponding entire hypersurface. As usual, finite-index
restrictions preserve the selected hypersurface and side.

\begin{theorem}[An exterior tube has a recognized second end]
\label{IV-basins:thm:tube}
Assume local description and $n\ge3$. Start with the full-angular
exterior collar in Section~\ref{IV-sec:exterior}, in the compact manifold
$P$ obtained by cutting along the starting hypersurface.
Its continued embedded vector model is
\[
 \Phi:\R^n\setminus\{0\}\longrightarrow\operatorname{int}P,
 \qquad U=\Phi(\R^n\setminus\{0\}),
\]
or the corresponding antipodal model. Let
\[
 K_0=\bigcap_{\epsilon>0}
       \overline{\Phi(\{0<|v|<\epsilon\})}.
\]
Then either the action is already one of the closed homogeneous
alternatives of the description input, or $K_0$ is a linear
fixed-point filling or one compact projective collar hypersurface.
The fixed point is possible only in the spherical, not the antipodal,
case. In the collar case its incoming developing radius tends
to zero and its transition with the tube is a constant homothety,
with the same standard or dual type.
\end{theorem}
\begin{proof}
Exterior continuation in Section~\ref{IV-sec:exterior} is injective and exactly
equivariant on the entire punctured vector space. Its zero end
$K_0$ is nonempty, compact, connected, invariant, and disjoint
from the starting collar. Every point of $K_0$ belongs to
$\overline U$. On $U$ there is the invariant nonvanishing analytic
Euler field
\[
 E_{\Phi(v)}=D_v\Phi(v).
\]
It descends in the antipodal case, as does the standard flat
connection, and $\nabla E=\id$.

If $K_0$ meets a cut boundary component, minimality of the
projective action, after fixing that component, puts the entire
component in $K_0$. It already has the original full collar;
the detected-collar argument in the last paragraph below applies
directly. We may therefore assume $K_0$ is compactly contained
in $\operatorname{int}P$. It embeds in the original closed
manifold. All suspension and germ theorems in the rest of the
proof are applied there, so no extension of those theorems to
manifolds with boundary is being assumed.

First suppose an ergodic invariant probability is supported on
$K_0$. The compact derivative alternative contradicts the
infinite standard action on $U$. In the standard or dual
alternative, every center chart meets $U$.
Theorem~\ref{IV-basins:thm:radial} makes the measure a finite-orbit measure.
The finite-orbit extension theorem in Section~\ref{IV-sec:exterior} then proves that
$K_0$ is one fixed point and excludes the antipodal case. This argument
handles singular nonatomic candidate measures by ruling them out,
rather than presupposing affine compatibility for such measures.

For completeness one need not assume that $K_0$ has any invariant
probability. The amenable minimal-parabolic averaging construction
in Proposition~\ref{II-nv:prop:averaging} gives a Haar-projecting ergodic
minimal-parabolic probability supported on $K_0^\alpha$. If it
is $G$-invariant, it corresponds to an invariant probability on
$K_0$ and the preceding paragraph applies. If that probability
is not ergodic, use an ergodic component; equivalently use
ergodicity of the induced $G$ action.

If it is not $G$-invariant, apply the proper-parabolic description clause. A closed homogeneous
alternative is already recognized. Otherwise it gives a compact
projective hypersurface $N_0\subset K_0$ and an actual analytic
collar. The collar meets $U$, since $N_0\subset\overline U$.
A product collar is impossible: after the finite-index restriction
its transverse coordinate is invariant, and its pullback to the
nonempty invariant open subset of the punctured vector model
would be a nonconstant invariant analytic submersion. The
vector-stabilizer argument, or Section~\ref{IV-sec:exterior},
rules this out.

The transverse analytic vector field is therefore nonzero on
each sufficiently small punctured side of $N_0$. A side meeting
$U$ has a full radial developing chart, either at zero or at
infinity. The detected-end comparison theorem in
Section~\ref{IV-sec:exterior} applies: standard-to-dual comparison is impossible,
and a standard-to-standard comparison is a constant homothety.
It identifies the whole incoming side with the zero end of
$\Phi$, forces the incoming radius to tend to zero, and gives
$K_0=N_0$. This includes finite radial covers and the antipodal
source treated there.
\end{proof}

The theorem identifies the actual other end of every exterior
tube satisfying local description. It does not infer tube
structure just from a negative exponent, and it does not assume
that every nonempty invariant compact set carries a full-lattice
invariant probability.

\section{Global reconciliation, torus fillings, and finite coverage}\label{IV-series:global}\label{IV-sec:global}\label{IV-series:last-section}

Throughout this section assume local description as in Definition~\ref{II-def:recognition}. Set aside the finite-image and closed-model alternatives whenever they occur.

\subsection{Comparison of projective coframes}

An affine region constructed in Section~\ref{IV-sec:descent} carries a flat analytic
connection $\nabla$ and a parallel analytic projective coframe
$[A_y]\in\Pj\operatorname{Iso}(T_yM,\R^n)$ satisfying
\begin{equation}\label{IV-global:eq:coframe}
 [A_{\gamma y}D_y\gamma]=[\sigma(\gamma)A_y],
\end{equation}
where $\sigma$ is standard or dual. Its charts have transition maps
$v\mapsto av+b$, $a\ne0$. Its smooth measure class is ergodic,
and it has finitely many connected components. The region comes
from a nonatomic invariant germ family or from the spanning-return
alternative for continuation from a linearizing chart or an interior collar.

\begin{lemma}[Type and coframe comparison]\label{IV-global:lem:comparison}
Let such a region meet an exactly equivariant standard or dual chart
family, with positive incidence, including a single punctured vector
chart. Then the representation types agree. In every good incident
chart the projective coframe is $[dv]$ and the connection is
$\nstd$. The coframe and connection extend over its full local
coordinate images. Two regions of this kind which meet have the
same data on their overlap.
\end{lemma}
\begin{proof}
Affine absorption in Section~\ref{IV-sec:basins} gives the connection assertion.
Here is the additional coframe comparison, which also identifies
the type. Write $\rho$ for the source representation and form
\[
 K(s,v)=[A_{\Phi_s(v)}D_v\Phi_s].
\]
Its covariance is $K(\gamma s,\rho(\gamma)v)
=[\sigma(\gamma)K(s,v)\rho(\gamma)^{-1}]$.
On $Y=\Gamma\backslash(G\times S)$ induce this field by
\[
 [g,s]\longmapsto
 [\sigma(g)^{-1}K(s,\rho(g)e_1)\rho(g)].
\]
It is defined on a conull finite-measure subset by the incidence
zero--one law. For $H=\{h:\rho(h)e_1=e_1\}$, projective unipotent
recurrence forces an invertible intertwiner between $\rho|_H$ and
$\sigma|_H$. Standard and dual restrictions cannot be intertwined:
one has a one-dimensional fixed-vector space and the other has none.
For equal types the commutant is precisely the scalars. Thus $K=[I]$.
Analyticity gives the identity throughout each open incident domain.

The incident region is conull in each good chart image. Pushing
forward the full standard coframe and connection therefore extends
them uniquely across that image; all such extensions agree on their
overlaps. To compare two regions, use the germ family or single
continued chart which generates the second one and then analytic
uniqueness on their essential image and regular extensions.
\end{proof}

Take the maximal regular extension $U$ of the original region: adjoin
all small coordinate neighborhoods in which the original region is
Lebesgue-conull and both its connection and projective coframe extend
analytically, with the coframe remaining invertible. Extensions agree
on conull overlaps, so this is well-defined, invariant, and idempotent.
The original region is conull in $U$. Each new neighborhood meets
the original region, so the number of connected components remains
finite. Fix a component after finite-index restriction when needed.
Distinct maximal regions either coincide on the relevant union of
components or are disjoint, by Lemma~\ref{IV-global:lem:comparison}.

\begin{lemma}[No invariant probability on the affine frontier]
\label{IV-global:lem:frontiermeasure}
For an infinite-image action satisfying (G), $\partial U$ carries
no $\Gamma$-invariant probability.
\end{lemma}
\begin{proof}
Take an ergodic component of such a probability. In the compact-type
branch the original action preserves the resulting global Riemannian
metric. The compact-target dimension bound of
Lemma~\ref{II-loc:lem:compact-target} then forces finite image.
Every standard or dual centered germ
meets $U$, since its center is in $\partial U$. Lemma~\ref{IV-global:lem:comparison}
and the incidence zero--one law extend the coframe and connection
across almost every center. Those centers belong to the regular
extension $U$, contrary to their location in $\partial U$.
\end{proof}

\subsection{The affine frontier consists of finitely many collars}

Suppose the action is not already a closed homogeneous alternative.
For a nonproduct collar write its local developing coordinates
as
\begin{equation}\label{IV-global:eq:develop}
 z=r(t)(1,u),\qquad \frac{r'(t)}{r(t)}=\frac1{B(t)},
 \qquad B(0)=0,\quad B\not\equiv0.
\end{equation}
On a sufficiently small punctured side $B$ has no further zero.
The developing radius tends to zero or infinity at the central
hypersurface. The standard connection pulled back by this map has
\begin{equation}\label{IV-global:eq:pole}
 \nabla_{\partial_t}\partial_{u_i}=B(t)^{-1}\partial_{u_i}.
\end{equation}
Thus it does not extend as an analytic connection across $t=0$.

\begin{proposition}[Finite boundary and absorption of an entire side]
\label{IV-global:prop:boundary}
The frontier of $U$ is a finite union of pairwise disjoint compact
analytic spherical or projective hypersurfaces. In a collar
of each, every punctured side either lies entirely in $U$ or is
disjoint from $U$. At least one side lies in $U$. No product collar
meets $U$.
\end{proposition}
\begin{proof}
A product collar, after a finite-index restriction, has an invariant
analytic transverse submersion. It is an invariant full collar:
for each group element the globally defined product formula agrees
with the action near the central hypersurface and hence on the
entire connected collar by analytic continuation. Pulling its
transverse coordinate back to an incident generating chart would
give a nonconstant invariant analytic function, contradicting the
source zero--one law. The same argument applies to each of the
finitely many ergodic components created by a finite-index restriction.

If a nonproduct punctured collar side meets $U$, continue its
developing chart to its exactly equivariant punctured domain.
Lemma~\ref{IV-global:lem:comparison} extends the data over every local image,
and the zero--one law makes $U$ conull there. By maximality the
entire side lies in $U$. Its central hypersurface cannot meet $U$,
by \eqref{IV-global:eq:pole}. Thus a collar which meets $\partial U$ through
its central hypersurface has frontier precisely that hypersurface
locally, with at least one incident side in $U$.

Let $\mathcal C$ be the invariant union of all these collar
neighborhoods, taking group translates as necessary. If
$K=\partial U\setminus\mathcal C$ were nonempty, it would be
compact and invariant. Parabolic averaging and
Lemma~\ref{IV-global:lem:frontiermeasure} give a proper-parabolic suspension
measure on $K$. By (P) it supplies a compact central hypersurface
in $K$ with its collar. Minimality of its projective lattice action
puts the entire hypersurface in $K$. Since it belongs to
$\partial U$, its collar meets $U$ and is among those defining
$\mathcal C$, a contradiction.

The collars therefore cover the compact frontier. Choose a finite
subcover. In each such collar the frontier is just its central
hypersurface, so there are finitely many. Two distinct central
hypersurfaces cannot meet: after a common finite-index restriction,
their intersection is closed and invariant in each, and minimality
would make them equal. This proves the proposition.
\end{proof}

\begin{lemma}[No exterior affine end]\label{IV-global:lem:noexterior}
Every boundary side incident to $U$ has developing radius tending
to zero at its boundary.
\end{lemma}
\begin{proof}
An exterior side continues, by (C), to an embedded invariant
punctured vector model $T$, with invariant nonvanishing Euler field.
Since it meets $U$, the two regions have conull overlap in their
generating charts. If $U$ is generated by a nonatomic invariant
germ family, the radial-field detection theorem in Section~\ref{IV-sec:basins}
would make that center measure atomic. If it is generated by a
single chart with spanning returns, the radial obstruction in
Section~\ref{IV-sec:descent} gives the same contradiction: a local transition
$w=av+b$ must satisfy $w=av$ in the presence of the pulled-back
radial field, making all returns collinear. These exclusions
persist on the regular extension because the generating image is
conull there. An exterior side is impossible.
\end{proof}

\subsection{Compact similarity filling and the torus conclusion}

\begin{theorem}[Torus filling without invariant ambient volume]
\label{IV-global:thm:filling}
Each connected component of $U$ has a compact affine orbifold
filling $\overline U^{\rm aff}$, obtained by filling its finitely
many zero-radius ends by points. Its local singular groups are
trivial or $\{\pm I\}$. After finite-index restriction its action
lifts to a finite manifold cover which is an affine torus, with
standard or dual linear part up to a removable sign character.
\end{theorem}
\begin{proof}
Small coordinate balls in the regular manifold $U$ are simply connected,
so the inclusion of each such ball kills its fundamental group.
On each connected component, Hatcher's universal-cover construction
therefore gives a simply connected cover
\cite[\S1.3, pp.~63--65]{Hatcher}.
The projection is a local homeomorphism, so the covering space is locally
path connected. Apply the development theorem to this cover.
For a connected real analytic model $X$ and a Lie group $G$ acting effectively,
transitively and real analytically on $X$, a $(G,X)$-atlas has coordinate
changes given locally by restrictions of elements of $G$. The Development
Theorem gives, on a specified universal covering $p:\widetilde U\to U$, a
$(G,X)$-map $\operatorname{dev}:\widetilde U\to X$ and a homomorphism
$h:\operatorname{Deck}(p)\to G$ satisfying
$\operatorname{dev}(\delta u)=h(\delta)\operatorname{dev}(u)$ for every deck
transformation $\delta$; importantly, the holonomy takes values in the
specified structure group \cite[\S5.2, Theorem~5.2.2, pp.~120--121]{GoldmanDevelopment}.
Here take $X=V$ and
\[
 H(V)=\{z\longmapsto az+b:a\in\R^\times,\ b\in V\},
 \qquad (a,b)(c,d)=(ac,ad+b).
\]
This is a real analytic Lie group with effective transitive action on $V$.
The homothety transition germs of the actual atlas define an $(H(V),V)$-atlas:
uniqueness of affine maps agreeing on an open set makes the local parameters
constant on each connected overlap component. Adjoining all charts compatible
in both directions gives the maximal geometric atlas without changing the
chosen charts. The actual covering charts of $p$ generate its canonical
pullback structure. Apply the theorem with this structure group and retain
the specified covering and its deck labels.
In each lifted atlas chart the developing map therefore has the form
$a\,\varphi\circ p+b$, with $a\ne0$; composing $h$ with the inclusion
$H(V)\hookrightarrow\operatorname{Aff}(V)$ gives scalar-affine deck holonomy.
This development step is separate from Fried's closed-manifold image
classification used below. At this stage $U$ is a manifold; extension over
filled origins uses the affine cap charts, rather than an unstated orbifold
version of the development theorem.

Cut along the finitely many frontier hypersurfaces, separating the
two sides of a nonseparating hypersurface as well. Outside smaller
boundary collars the closure of the selected component is compact
and lies in $U$. Lemma~\ref{IV-global:lem:noexterior} makes every end a full
punctured ball in the developing coordinate $z$, or its antipodal
quotient. Adjoin its origin. On overlaps with $U$, the coframe
comparison lemma identifies the given data with the standard data
in $z$. This constructs a compact affine orbifold with a homothety
atlas. Notice that the change from $t$ to $z$ need only be analytic
away from $t=0$; the filling supplies its own analytic chart at the
new point. The action extends affinely across those points by its
developing formula. Fix components and the finite end set by
finite-index restriction.

In the Euclidean covering chart of each end, polar coordinates identify
the punctured ball with $S^{n-1}\times(0,r)$. The sphere is simply connected
for $n\geq 3$ \cite[Proposition~1.14, p.~35]{Hatcher},
and the radial interval is contractible. Thus the punctured chart is
simply connected before the antipodal quotient is taken.

The topological finite-generation step uses two separate classical results.
A compact manifold, with or without boundary, is a Euclidean neighborhood
retract (ENR), and the fundamental group of a compact ENR is finitely generated
\cite[Corollaries~A.9 and A.8, p.~527]{Hatcher}.
Here an ENR is a space admitting an embedding into some $\R^N$ whose image
is a retract of an open neighborhood. For the actual compact cut component
$Y$, choose linear coordinates $F\simeq\R^d$ in its half-space model
$F\times[0,\infty)$ and postcompose the same charts. This exhibits the same
space, with the same topology and basepoints, as a conventional manifold
with boundary. Apply A.9 and then A.8. The chosen boundary half-collars on
the cut permit loops and their compact homotopies to be pushed into the
interior corresponding to $U$, giving the required finite generation there.
Inversion sends a left Haar fundamental domain to a right fundamental
domain for the inverse Haar measure. The positive finite mass of the
domain gives an invariant probability on the right-coset quotient by
normalization. The real special linear group is connected, so the
lattice is finitely generated by \cite[Corollary~1.6, p.~34]{GelanderSlutsky}.
Schreier's theorem gives finite generation of the finite-index restriction.

For simple connectivity in the lifted filling, use the surjectivity clause
of van Kampen's theorem in its general form: an open path-connected cover
$\{A_i\}$ with a common basepoint $x_0$ and path-connected pair intersections
has a surjective canonical homomorphism
\[
 \mathop{\ast}_i\pi_1(A_i,x_0)\longrightarrow\pi_1(Z,x_0)
\]
induced by the inclusions \cite[Theorem~1.20, p.~43]{Hatcher}.
In the capped covering space let $A$ be the original simply connected cover
and let $B_i$ be its disjoint lifted ball neighborhoods at the new origins.
Each $A\cap B_i$ is a punctured ball, hence path-connected for $n\ge3$.
The two-open-set application, based in
$A\cap B_i$, first makes $A\cup B_i$ simply connected. Now cover the filled
space by $A$ and the sets $A\cup B_i$. They share any chosen basepoint in $A$;
distinct-label intersections are $A$. The displayed homomorphism is therefore
surjective from a free product of trivial groups. The union is path-connected,
so it is simply connected. This argument permits an arbitrary index set of
lifted ends.

Pointed lifting on the regular manifold cover $p:\widetilde U\to U$
uses only its simply connected, locally path-connected domain: for a
continuous $f:\widetilde U\to U$ and points $x_0,x_1$ with
$p(x_1)=f(x_0)$, there is a continuous lift $F$ with $F(x_0)=x_1$ and
$p\circ F=f$ \cite[Proposition~1.33, pp.~61--62]{Hatcher}.
Applying pointed lifting to a base homeomorphism composed with $p$ and to
its inverse, and then using uniqueness of lifts, gives an actual lifted
homeomorphism. All lifts, with their base-action labels, form the lift group.
Compatibility of the finite cover and its action is established next by the
group-extension argument.

We recall why a compatible finite manifold cover and an actual
lifted action exist. Let $\Pi$ be the orbifold fundamental group
and $\widehat\Gamma$ the group of lifts, retaining the lattice
labels. The developing map yields
\[
 1\longrightarrow\Pi\longrightarrow\widehat\Gamma
 \longrightarrow\Gamma\longrightarrow1,
 \qquad \rho:\widehat\Gamma\longrightarrow\operatorname{Aff}(V).
\]
These groups are finitely generated. Local finite isotropy has
faithful affine holonomy. Embed the affine group in $\GL(n+1,\R)$
by homogeneous coordinates. Selberg's lemma, in the form of
\cite[p.~269]{Alperin}, applies to the finitely generated image
$\rho(\widehat\Gamma)$. Choose a torsion-free
finite-index subgroup $L$ of $\rho(\widehat\Gamma)$, and set
$E=\rho^{-1}(L)$, $\Pi_0=E\cap\Pi$. The cover corresponding to
$\Pi_0$ has no local isotropy and is a manifold. Since $\Pi_0$
is normal in $E$, the quotient $E/\Pi_0$ is the image of $E$ in
$\Gamma$ and acts on that cover. This constructs a simultaneous lift rather than assuming that
an arbitrary finite cover admits a lattice lift.

Apply the Euclidean similarity classification of Fried to this closed
connected similarity manifold of dimension $n\ge3$
\cite[Theorem 2 and its proof, pp.~580--581]{Fried1980}. A scalar
homothety with multiplier $a\ne0$ scales Euclidean lengths by $|a|$,
so the same charts and the same universal developing map give a
similarity development. The nilpotent generalization
\cite[Theorem IV.1]{Alexandre} is not needed here: the developing
space is the vector space $V$ and the transitions are Euclidean
similarities. It is complete Euclidean, or its universal
developing image is $V\setminus\{o\}$ and it is a finite quotient
of a Hopf manifold. In the latter case every lifted affine
automorphism fixes $o$, since it preserves the developing image.
The radial field $z-o$ therefore descends, is nowhere zero, and
is invariant under both the lattice and the finite covering
transformations. It descends further to the filled orbifold.
Its restriction contradicts exactly the radial obstruction used
in Lemma~\ref{IV-global:lem:noexterior}. Thus the filling is complete.

In the complete finite manifold cover, the holonomy linear parts
are scalars because the original atlas has homothety transitions.
An affine map $z\mapsto az+b$ with $a\ne1$ has a fixed point in
$V$ and cannot be a nontrivial deck transformation of a complete
manifold. All deck transformations are therefore translations.
Compactness makes their group a full lattice $\Lambda<V$, and
the cover is $V/\Lambda$.

In the developing trivialization, equivariance of the projective
coframe makes the derivative of a lifted lattice element equal to
$\lambda_\gamma\sigma(\gamma)$ for a nonzero real scalar.
It preserves $\Lambda$, so the absolute determinant is one and
$|\lambda_\gamma|=1$. The signs form a character and disappear
after finite-index restriction. Translations give the usual
affine cocycle on $V/\Lambda$. This proves the assertion without
assuming volume before compactification; the resulting affine
torus now has its invariant parallel Haar volume.
\end{proof}

\begin{proposition}[A torus component determines the lattice]
\label{IV-general:prop:toral-arithmetic}
If a torus component occurs in Theorem~\ref{IV-global:thm:filling},
there are a finite-index subgroup $\Gamma_0<\Gamma$ and
$B\in\GL(n,\R)$ such that
\[
 B^{-1}\sigma(\Gamma_0)B<\SL(n,\Z)
\]
has finite index. In particular $\Gamma$ is commensurable with
a conjugate of $\SL(n,\Z)$. No other lattice can have a torus
component of the specified standard or dual type.
\end{proposition}
\begin{proof}
The filling proof constructs an actual affine action of
$\Gamma_0$ on $V/\Lambda$, and its determinant argument removes
the scalar character after further finite index. Thus
$\sigma(\Gamma_0)\Lambda=\Lambda$. Choose $B$ with
$\Lambda=B\Z^n$. The displayed inclusion follows.
Its left side is a lattice in $\SL(n,\R)$, since $\sigma$ is
an automorphism and conjugation by $B$ is a Lie-group automorphism.
For lattices $L<H<G$, covolumes satisfy
$\operatorname{covol}(L)=[H:L]\operatorname{covol}(H)$;
this follows by taking disjoint translates of a fundamental
domain, including when the index is infinite. Finiteness of
both covolumes gives finite index here. In the dual case,
apply inverse transpose, which preserves $\SL(n,\Z)$, to
obtain the assertion for $\Gamma$ itself.
\end{proof}

After this proposition, standard integer coordinates may be
used for the action on a torus cover. All subsequent arithmetic
facts about such a filled component are conclusions of its
construction, not restrictions on the initial lattice.

\begin{corollary}[Nonatomic torus measures are Haar]
\label{IV-global:cor:haar}
Every nonatomic invariant probability on a filled torus component,
with the action furnished by Theorem~\ref{IV-global:thm:filling}, is the
pushforward of Haar probability on its torus cover. Consequently
the nonatomic invariant measures which generated the affine regions
have smooth positive densities on those regions. For the original
group, an ergodic such measure averages the Haar probabilities over
the finite orbit of its connected torus components.
\end{corollary}
\begin{proof}
Lift the probability to the equivariant finite torus cover by giving
equal weights to its finite fibers. The isolated filling points have
zero mass, so the ordinary constant-degree covering off those points
suffices to define this nonatomic invariant lift. Write the torus
as $V/\Lambda$, and its affine action as
$x\mapsto\sigma(\gamma)x+b_\gamma$. For $k\in\Lambda^*$,
invariance implies
\[
 |\widehat\nu(k)|=|\widehat\nu(\sigma(\gamma)^T k)|.
\]
Translations contribute only unit-modulus factors.

Every nonzero dual-lattice orbit here has positive upper density in
$\Lambda^*$. To see this, use Proposition~\ref{IV-general:prop:toral-arithmetic}
to choose a basis of $\Lambda$ in which the linear action is
a finite-index subgroup of $\SL(n,\Z)$. The dual lattice in
these coordinates is $\Z^n$. The full
$\SL(n,\Z)$ orbit of a nonzero integer vector is its coordinate
gcd times the primitive vectors, a set of positive density. It is
a finite union of linear images of the finite-index-subgroup orbit.
If that subgroup orbit had zero upper density, each of its finitely
many linear images would too, since inverse images of boxes lie in
boxes of comparable size. This contradicts positive density of the
full orbit. Rescaling and passing between commensurable lattices
preserve this positive-upper-density conclusion.

Finally, the elementary Wiener identity gives, in any lattice basis,
\[
 \lim_{R\to\infty}\frac1{(2R+1)^n}
       \sum_{\|k\|_\infty\le R}|\widehat\nu(k)|^2
       =\sum_x\nu(\{x\})^2=0.
\]
For this form of the identity, integrate the normalized product
Dirichlet kernel against $\nu\times\nu$. Its absolute value is
at most one; it tends to zero off the diagonal and equals one on
the diagonal, so dominated convergence applies. A nonzero constant
Fourier modulus on any nonzero orbit of positive upper density
would contradict the displayed limit. Hence all nonconstant
Fourier coefficients vanish, and the lift is Haar. Restricting to
the affine region gives its positive analytic parallel density.
Ergodicity and the finite component permutation give the last
assertion.
\end{proof}

\subsection{An embedded interior continuation reaches an exterior collar}

\begin{proposition}[Completion of an embedded interior model]
\label{IV-global:prop:interior}
Start with a fixed point or a nonproduct interior collar side. If its
continued map is in the embedded alternative in Section~\ref{IV-sec:descent},
it extends to an embedded full punctured vector model, with its
prescribed antipodal quotient if present. Its infinity end is a
projective exterior collar. Its zero end is the starting fixed point or
interior collar, and the detected-end theorem in (C) identifies
it as such. Hence this case is a $G$-disk or $G$-tube.
\end{proposition}
\begin{proof}
Let $F:\Omega\to M$ be the embedded continued model, with
$\Omega$ a union of punctured centered ellipsoids containing a
punctured ball. In the collar case first cut along the starting
hypersurface and work on the selected side, as in Section~\ref{IV-sec:exterior}.
The group translates used in continuation avoid that hypersurface
and stay in the selected cut component. Write $T=F(\Omega)$;
it has its invariant nonzero Euler field and flat connection.

The source $\Omega$ is conull in $V\setminus0$ by the zero--one
law, so it contains vectors of arbitrarily large norm. In the
compact cut manifold define
\[
 K_\infty=\bigcap_{R>0}
       \overline{F(\{v\in\Omega:|v|>R\})}.
\]
This is nonempty, compact, invariant, and disjoint from $T$.
The last assertion follows from continuity of the inverse of
the embedding: vectors whose norms tend to infinity cannot have
images converging to an interior point of $T$. It is also disjoint
from the starting zero-side boundary copy, because the embedded
initial punctured ball fills a whole one-sided collar there.

We record a comparison valid even though $\Omega$ may be proper.
If an exactly equivariant linearizing chart or punctured collar chart $\Psi$
meets $T$, then on its incident set
\begin{equation}\label{IV-global:eq:inversecompare}
 F^{-1}\Psi(w)=c w,\qquad c\ne0,
\end{equation}
in the same representation type, with the antipodal version
understood. Induce the vector-valued left side to the finite
probability space over the vector stabilizer. Unipotent recurrence
forces it into $\R e_1$ in equal types; in opposite types the
fixed-vector space is zero, impossible for a punctured target.
The scalar is invariant and hence constant by the zero--one law.
Analyticity gives \eqref{IV-global:eq:inversecompare} on its full open
incident domain. In the antipodal case use the projective sign
extension, or square the scalar and then choose a local lift.

An invariant probability on $K_\infty$ is impossible. Its compact
branch would give finite global image. By radial-field detection
its noncompact branch is atomic and supplies a linearizable fixed point.
Equation~\eqref{IV-global:eq:inversecompare}, and the fact that $\Omega$
contains a full punctured ball, identify the entire punctured
linearizing neighborhood with bounded vectors near zero in $\Omega$.
The uniqueness of the inverse of $F$ then prevents that fixed point
from being in $K_\infty$.
Here the finite-orbit detection precedes the point-stabilizer
application of germ superrigidity: the stabilizer has finite index
and is therefore a lattice. The underlying matrix-cocycle input is
\cite[Definition~1.3 and Theorem~1.5, p.~193]{FisherMargulis},
applied to its one-point action; the analytic linearization is the
further conclusion of Part~I. No superrigidity assertion for an
arbitrary infinite-index point stabilizer is used. The cocycle in each
jet degree has the original finite-dimensional algebraic jet group
as its target. Compact sets in the discrete lattice are finite, so
the elementwise logarithmic moment condition gives the compact-set
condition of Definition~1.3. The algebraicity step uses
\cite[Lemma~3.4, p.~199]{FisherMargulis} on the original special-linear
real point group, as in Part~II. The rational representation agrees
on all of those points, so its matrix coefficients give polynomials
in the original defining coordinates.
The factorization in Part~II gives the actual normal factors of this
same algebraic group. Their inclusions and product projections induce
the corresponding maps on the original matrix tangent spaces.
The Cartan and compact-dual data are those constructed in Part~II
from the negative-transpose and negative-adjoint involutions on the
original real and complex matrix Lie algebras.
In each finite jet degree, real Lie complete reducibility gives invariant
complements. The same exponential curves make the corresponding
projections equivariant for the original group representation.
The type-A step, also as in Part~II, applies Jacobson's matrix
automorphism theorem to the complexification of the real differential.
Its two formulas are conjugation and conjugation after negative
transpose. They give the standard and dual alternatives after real
descent.
For spectral transfer we use the comparison of Part~II on the actual
simple factors and their Haar measures. The passage to arbitrary real
Hilbert representations preserves both the spectral exponent and the
regular averaging operator.
The regular-operator norm estimate of Part~II applies on each
factor's nonzero canonical Cartan ray, including rays on chamber
walls. It uses the maximal compact subgroup and Haar measure of
that factor.

If $K_\infty$ meets a cut boundary, minimality supplies its full
known central hypersurface. Otherwise it is compactly inside the
original manifold, and parabolic averaging with (P) supplies a
central hypersurface and its collar. In this embedded case the
closed alternatives can also be excluded directly. For a flag or
Grassmannian model, the selected lattice action on its connected
finite cover is minimal. On the original manifold, however, the
embedded partial model has a nonempty compact invariant infinity
remainder disjoint from the nonempty set $T$, contradicting that
minimality. Restricting to the selected finite-index lattice preserves
this obstruction.

For a Hopf model, its actual developing map
$q:V\setminus0\to M$ is surjective, analytic, and equivariant for a
finite-index lattice, and satisfies $q(e^T v)=q(v)$ with $T>0$.
Surjectivity gives incidence with $F$. In equal representation types,
\eqref{IV-global:eq:inversecompare} holds throughout the incident set;
periodicity preserves that set and would force $e^T=1$. For an
antipodal source the two inverse vectors may differ by sign, but
positivity of $e^T$ gives the same contradiction. In opposite types
the zero fixed-vector comparison already forbids incidence. Thus no
closed alternative occurs here.

A product collar cannot meet $T$, by its invariant transverse
submersion and source ergodicity. On a side through which
$F(v_j)\to N$ with $|v_j|\to\infty$, the new developing radius
cannot tend to zero: \eqref{IV-global:eq:inversecompare} would again give
bounded source vectors. Thus that side is exterior.

By (C) its exterior chart continues to an embedded map
$\Psi:V\setminus0\to M$. On a nonempty open overlap,
\eqref{IV-global:eq:inversecompare} says $F(v)=\Psi(c^{-1}v)$.
Both sides are analytic and defined on all of the connected
$\Omega$, so the identity holds throughout. This extends $F$
to the full punctured vector model. The zero-end comparison of
Section~\ref{IV-sec:exterior} and Section~\ref{IV-sec:basins} now identifies its entire zero end with
the starting fixed point or collar. No assertion that lattice translates
alone cover every nonzero source vector was needed.
\end{proof}

Combining this proposition with the self-overlap dichotomy in
Section~\ref{IV-sec:descent} proves that every nonproduct collar side and every
fixed point belongs either to a affine torus quotient or to a completed
embedded disk or tube.

\subsection{Product collars form a separate global alternative}

\begin{proposition}[Global product-foliation alternative]
\label{IV-global:prop:product}
If a product projective collar occurs, the entire manifold is
foliated by compact spherical or projective leaves with local
product action. After finite covers and finite-index restriction,
this is the standard or dual projective action on
$S^{n-1}\times S^1$ (or the projective quotient), trivial on the
circle factor.
\end{proposition}
\begin{proof}
Take the invariant union $B$ of all full zero-field collar
neighborhoods, including the normal-involution quotient collars.
Each is invariant after its finite-index restriction, by analytic
continuation of its covered product formula. In a normal-involution
quotient the regular two-sided product locus is the complement of
the central hypersurface and is dense. All open-incidence arguments
below are made on this regular locus, where the transverse
submersion is available.
No standard or dual probability germ family can meet it, and no
nonproduct radial chart can meet it: pull back the invariant
transverse submersion and apply the appropriate diagonal zero--one
law. These assertions remain valid on finite-index subgroups.

If $B\ne M$, connectedness gives a nonempty compact invariant
frontier $\partial B$. Apply parabolic averaging there. The compact
derivative branch contradicts the nontrivial projective action.
The standard and dual branches have charts centered on that frontier,
so they necessarily meet $B$. A countable subcover by product collars
then gives positive incidence with one collar; its transverse
submersion gives the preceding contradiction. The proper-parabolic branch
gives a collar on that frontier (or a closed homogeneous alternative).
A nonproduct collar cannot meet $B$ by the same submersion argument;
a product collar belongs to $B$ by definition, also a contradiction.
The closed homogeneous alternatives also exclude a product collar:
their smooth measure classes are ergodic under every finite-index
lattice. For the regular Hopf manifold this follows
by pulling back to the punctured vector model; for the Grassmannian
it follows from the finite-volume lattice quotient and its noncompact
parabolic stabilizer. An invariant product collar with its transverse
submersion contradicts this ergodicity. Thus $B=M$.

Compactness gives a finite covering by covered zero-field collars, so a
single finite-index subgroup fixes their covered transverse coordinates.
The product foliations agree on overlaps. Indeed, after fixing a
common finite-index subgroup, the intersection with each compact
projective leaf is open and invariant. Minimality makes it the whole
leaf. The other collar's invariant transverse coordinate is constant
on that leaf, so the plaques agree. These compact spherical or projective leaves have finite fundamental
group and hence finite holonomy. Normal-involution collars can produce
reflector leaves: the original leaf space need not be a circle.
Use the canonical transverse-orientation double cover of the foliation
\cite[\S1.4, pp.~371--372]{HaefligerVarietesFeuilletees}.
The transition between two signs is multiplication by the sign of the
derivative of the transverse chart change. The chain rule gives the
cocycle identity. A foliation-preserving analytic diffeomorphism acts
on these signs by the derivative sign of its transverse map, and the
same chain rule gives compatibility on overlaps and with composition.
The normal-involution
charts lift there to genuine two-sided spherical product charts, including
at the central leaves. Compact-leaf stability and the oriented leaf-space
construction apply on this cover
\cite[\S\S2.6, 3.3, pp.~381--382, 387]{HaefligerVarietesFeuilletees}.
Each surjective angular local homeomorphism from the compact sphere
to a spherical or projective leaf is a covering. Since $n-1\geq2$,
the sphere is simply connected
\cite[Proposition~1.14, p.~35]{Hatcher}.
Uniqueness of the pointed universal cover identifies the spherical case
with a one-sheeted cover and the projective case with the usual double
cover. On the original angular sphere, sending each point to the other
point of its fiber gives the nonidentity deck map. It commutes with every
homeomorphism preserving the fibers, since such a map sends the other
point of a fiber to the other point of the image fiber. The angular
rigidity lemma identifies this deck map with the antipodal map.
Remove the angular double covers as well. The resulting compact
leaf bundle has circle base and spherical fibers. In the two-sided
case, this base can also be constructed directly: the full product
charts are saturated, their transverse intervals give interval charts
on the compact Hausdorff second-countable quotient, and
\cite[p.~173]{GaleOneManifolds} identifies that quotient topologically
with a circle.  Transition maps
on the angular factor commute with a fixed finite-index lattice
action. Lemma~\ref{II-general:lem:angular} says that the centralizer is
trivial on projective space and consists of the two signs on
the sphere; it also excludes a change between standard and
dual type. A further
finite cover kills the sign monodromy. The finite-structure-group gluing can equivalently be made using the
principal bundle of angular frames: its finite discrete fibers give a finite
cover of the leaf circle, and the associated angular bundle becomes a product
on that cover; see \cite[Chapter~I, \S5, Proposition~5.2 and pp.~54--55]{KobayashiNomizuI}.
Although that reference is stated in the smooth category, the same gluing
formulas are analytic here: the finite-group transitions are locally constant
and the angular maps are analytic. The transverse analytic coordinates descend
from the original collars. A compact connected analytic one-manifold is
analytically a standard circle: use the smooth classification
\cite[Appendix, pp.~55--57]{MilnorDifferentiableViewpoint}, approximate its
circle diffeomorphism in $C^1$ by an analytic map
\cite[Proposition~11.7, p.~648]{TsuboiAnalyticDiffeomorphisms}.
Diffeomorphisms between the two manifolds form an open set in the
strong $C^1$ topology
\cite[Chapter~2, Theorem~1.7, p.~38]{Hirsch1976}.
For a compact source this is the weak $C^1$ topology of chartwise
tests on compact sets \cite[Chapter~2, p.~35]{Hirsch1976}.
The analytic approximation can therefore be chosen to be a
diffeomorphism, and its inverse is analytic by the analytic inverse
function theorem.
This gives the stated analytic product.
\end{proof}

\subsection{Disjointness, coverage, and finiteness}

Assume now that the finite-image and the closed global alternatives
have not occurred, and that there is no product collar. We have
constructed affine torus quotients, $G$-disks, and $G$-tubes
from every germ or collar side in the local description.

\begin{lemma}[Distinct interiors do not overlap]\label{IV-global:lem:disjoint}
Two completed piece interiors either coincide or are disjoint.
\end{lemma}
\begin{proof}
For affine pieces use the unique regular extension and
Lemma~\ref{IV-global:lem:comparison}. Here the pieces are connected
components of the maximal regular regions. If their original label
subgroups differ, use their finite-index intersection without changing
the original coordinate maps, coframes, or connections. The completed
torus models of Theorem~\ref{IV-global:thm:filling} and
Corollary~\ref{IV-global:cor:haar} give the smooth measure class on each
component. This class is ergodic under every further finite-index
component-preserving subgroup: an invariant positive-measure set would give, by normalized
restriction, another nonatomic invariant probability on the same torus
model, hence the same Haar probability. Thus open incidence gives the
local conullness needed to compare the maximal connected regular
extensions. This does not require the original source probability to
remain ergodic after restriction.
For two embedded full punctured vector
models, apply the vector-stabilizer comparison to the transition
on their invariant intersection. It is a scalar homothety of
equal representation type. Since both models are defined on
the entire punctured space, analytic continuation identifies
their whole images. Opposite types have no transition.
An affine piece cannot overlap an embedded tube or disk interior,
because the latter's Euler field contradicts the nonatomic-center
or spanning-return obstruction used in Lemma~\ref{IV-global:lem:noexterior}.
Fixed points of $G$-disks are included as endpoints; a fixed point absorbed by an
affine piece is an ordinary point of that piece, not a second piece.
\end{proof}

\begin{proof}[Proof of Theorem~\ref{IV-global:thm:global}]
Let $B$ be the invariant union of all completed piece interiors,
their linearizing neighborhoods, and their full incident collar neighborhoods.
It is open. Suppose the compact invariant remainder $K=M\setminus B$
is nonempty, and take a Haar-projecting ergodic minimal-parabolic
probability on its suspension.

If it is $G$-invariant, take the corresponding ergodic invariant
probability on $K$. Compact derivative type gives the excluded
finite-image alternative. For nonatomic standard or dual type,
the descent theorem Section~\ref{IV-sec:descent}, the regular-extension
construction, and Theorem~\ref{IV-global:thm:filling} put almost every center
in a torus piece, contradicting support on $K$. For atomic type, use the finite probability extension in the
invariant-measure clause of local description. A positive atom over an
original center has a finite-index stabilizer, and its actual coordinate
germ supplies a linearizing chart for that subgroup. Its continuation is
either affine by the self-overlap theorem, or embedded and handled by
Proposition~\ref{IV-global:prop:interior}. In the affine case the regular
coframe and connection extend across the chart center: their actual chart
formulas agree on the punctured neighborhood, which is conull and dense.
Thus the center belongs to the same maximal regular region, and a small
connected neighborhood belongs to its connected component. Theorem~\ref{IV-global:thm:filling}
then applies to that component. In either case the original center belongs
to $B$.

If the measure has proper parabolic stabilizer, the proper-parabolic clause gives a central
projective hypersurface with a full collar, unless a closed global
alternative already concludes the classification. The collar is
nonproduct. Each of its sides is affine or embedded by
Section~\ref{IV-sec:descent}; Theorem~\ref{IV-global:thm:filling} and
Proposition~\ref{IV-global:prop:interior} complete those pieces. The central
hypersurface and its entire small collar are included in $B$,
again contradicting support on $K$. Thus $B=M$.

Coverage alone would not make a finite decomposition. Here local
compatibility supplies the additional argument. At an interior
point of a piece, choose a neighborhood contained in that piece.
At a fixed point, choose its full linearizing neighborhood. At a central
hypersurface, choose a small full collar; each of its two sides
belongs entirely to one of the completed pieces. By
Lemma~\ref{IV-global:lem:disjoint} none of these neighborhoods meets any
other piece interior. These neighborhoods cover the compact
manifold, so a finite subcover meets all pieces and proves that
there are finitely many of them. Affine pieces have finitely many
boundary components by Proposition~\ref{IV-global:prop:boundary}, and each
disk or tube has its stated finite end set. Thus there are finitely
many collar hypersurfaces as well.

The actual collar charts specify every gluing, including two sides
belonging to the same piece. On incident vector coordinates the
comparison is a scalar homothety of the same standard or dual type;
on the collar itself retain its given analytic transverse vector
field. This proves the finite decomposition and its compatibility.
Proposition~\ref{IV-general:prop:toral-arithmetic} gives the stated
arithmetic restriction whenever a torus piece occurs.
If no torus piece occurs, Proposition~\ref{IV-global:prop:ambient}
below extends the action to an ambient group. Thus the decomposition
also gives precisely the three alternatives of
Theorem~\ref{IV-main:thm:classification}.
\end{proof}

\begin{proposition}[Every torus-free decomposition extends]
\label{IV-global:prop:ambient}
Let an analytic lattice action have one of the closed or product
models in Theorem~\ref{IV-global:thm:global}, or a finite compatible
decomposition of the stated kind containing no torus piece.
After the specified finite-index and covering modifications,
it extends to an analytic action of a connected group locally
isomorphic to $\SL(n,\R)$. This implication does not require
local description once the decomposition is given.
\end{proposition}
\begin{proof}
The product alternatives carry the ambient action on their specified
finite covers, with the circle coordinate fixed. For a closed homogeneous
model covered by the original manifold, use the group of literal lifts
of the homogeneous action through that covering, as in
\cite[Section~2, Proposition~2.4 and Definition~2.6]{MontaldiOrtegaLifts}.
For the topology, fix a point in the covering manifold and record the
group label and the value of each lift at that point. Connectedness and
uniqueness of lifts distinguish the lifts by this pair. Homotopy lifting
\cite[Proposition~1.30]{Hatcher} along paths in $G$
gives lifts for every group label. Applying it to the radial paths in a
coordinate ball and to their pointwise inverses gives local sections
with jointly continuous forward and inverse evaluations. Deck translates
give the covering charts. The deck label is locally constant in each
chart, so these evaluations are jointly continuous on the whole lift
group. Evaluating composition and inversion at the fixed point gives
continuity of the group operations in the chosen topology.
Pulling back the analytic covering atlas makes its operations and
its action analytic. For a finite covering its deck group is finite,
so passing to the identity component requires only a finite-index
restriction of the original lattice action. The original lattice embeds
in the literal lift group by its given matrix and its given action;
thus this construction preserves both labels without requiring a
splitting of a central extension. For a Hopf model the specified cover
runs from the Hopf model to the original manifold, and its ambient action
is instead the descended linear action, including the logarithmic radial
drift in the circle coordinate. In the remaining case the finite
decomposition consists solely of $G$-tubes, linearizable fixed points,
and their full collar neighborhoods.

On the union of tube interiors, transport the standard or dual $G$
action through the actual developing maps. A transition is a scalar
homothety of the same representation type, so these actions agree.
At a linearizable fixed point use the corresponding linear action.
In a covered collar use the angular $G$ action and the given transverse
flow at its logarithmic norm cocycle. The radial differential equation
identifies this formula with the transported action on both punctured
sides, using the reflected transverse field on the negative side.

The covered collar parameterization has actual local analytic inverse
branches. Composing with these branches gives jointly analytic local
action maps on the original manifold, including for a normal-involution
quotient; no global two-sided chart is required. These maps agree on
overlaps because they agree on the dense union of tube interiors.
They therefore glue, together with the fixed-point maps, to a jointly
analytic map $G\times M\to M$. The identity and multiplication laws
hold on the same dense union, hence everywhere. This constructs the
ambient $G$ action directly.

Take a common finite-index subgroup inside all the piece, fixed-point,
and collar labels in the given finite decomposition. Its original
action agrees with the constructed action on the dense tube interiors,
and therefore on all of $M$. The closed homogeneous alternatives may
retain the indicated covering group.
\end{proof}

\begin{corollary}[The alternative for other arithmetic types]
\label{IV-general:cor:extension}
Under the standing hypotheses, if $\Gamma$ is not commensurable with
a conjugate of $\SL(n,\Z)$, then, after the finite-index and
covering modifications in the classification, its action has
finite image or extends to an analytic action of a connected
group locally isomorphic to $\SL(n,\R)$. In particular,
a cocompact lattice has only these
alternatives.
\end{corollary}
\begin{proof}
Apply Theorem~\ref{II-main:thm:recognition} for $n\ge4$ or
Theorem~\ref{III-main:thm:recognition-three} for $n=3$, followed
by Theorem~\ref{IV-global:thm:global}.
Proposition~\ref{IV-general:prop:toral-arithmetic} excludes all
torus pieces, and Proposition~\ref{IV-global:prop:ambient} gives
the ambient extension in every remaining infinite-image case.
A cocompact lattice cannot be commensurable with a conjugate of
$\SL(n,\Z)$. Compact representatives for its quotient, finitely many
coset representatives for the intersection subgroup, and conjugation
would give compact representatives for $G/\SL(n,\Z)$. Take
$a=\operatorname{diag}(2^{-(n-1)},2,\ldots,2)$. If
$a^k=b_k\gamma_k$ with $b_k$ in a fixed compact set and
$\gamma_k\in\SL(n,\Z)$, then
$\gamma_ke_1=b_k^{-1}a^ke_1\to0$. This is impossible for nonzero
integer vectors.
\end{proof}

\clearpage
\phantomsection

\end{document}